\documentclass[11pt]{amsart}
\usepackage[foot]{amsaddr}
\usepackage[a4paper,margin=2.65cm]{geometry}
\usepackage[english]{babel}
\usepackage[T1]{fontenc}
\usepackage[utf8]{inputenc}
\usepackage{textcomp}
\usepackage{amsmath,amssymb,amsthm,mathtools}
\usepackage{newtxtext,newtxmath}
\usepackage[bb=boondox,cal=boondoxo,scr=boondoxo]{mathalfa}
\allowdisplaybreaks
\usepackage[renew-dots,renew-matrix]{nicematrix}
\mathtoolsset{showonlyrefs}
\usepackage{enumitem}
\usepackage{needspace}
\usepackage{etoolbox}
\makeatletter
\patchcmd{\@sect}{\@tocwrite{#1}{#8}}{\@tocwrite{#1}{#7}}{}%
  {\PackageError{BTP}{Cannot set short contents headings}{}}
\patchcmd{\@sect}{\@tocwrite{#1}{#8}}{\@tocwrite{#1}{#7}}{}%
  {\PackageError{BTP}{Cannot set short run-in contents headings}{}}
\makeatother

\usepackage{microtype,booktabs}
\usepackage{tikz}
\usepackage{graphicx}
\usetikzlibrary{arrows.meta}
\usepackage[hidelinks]{hyperref}
\usepackage{orcidlink}
\hypersetup{
 pdftitle={Banded Totally Positive Matrices and Positive Bidiagonal Factorizations: Existence, Stability, and Toeplitz Orthogonality},
 pdfsubject={Positive bidiagonal factorizations and Toeplitz orthogonality},
 pdfauthor={Manuel Ma\~nas}
}
\numberwithin{equation}{section}
\numberwithin{figure}{section}
\newtheorem{theorem}{Theorem}[section]
\newtheorem{proposition}[theorem]{Proposition}
\newtheorem{lemma}[theorem]{Lemma}
\newtheorem*{uniformchoicelemma}{Lemma}
\newtheorem{corollary}[theorem]{Corollary}
\newtheoremstyle{remarkstyle}
{6pt}{12pt}{\normalfont}{}{\itshape}{.}{0.5em}{}
\theoremstyle{remarkstyle}
\newtheorem{remark}[theorem]{Remark}
\newtheorem{example}[theorem]{Example}
\newcommand{\dd}{\mathrm d}
\newcommand{\NN}{\mathbb N_0}
\newcommand{\RR}{\mathbb R}
\newcommand{\CC}{\mathbb C}
\newcommand{\im}{\operatorname{Im}}
\newcommand{\diag}{\operatorname{diag}}
\newcommand{\supp}{\operatorname{supp}}
\newcommand{\BTP}{\mathrm{BTP}}
\newcommand{\PBF}{\mathrm{PBF}}
\newcommand{\hc}{h}
\newcommand{\gbinom}[2]{\genfrac{[}{]}{0pt}{}{#1}{#2}_{\!\tau}}
\DeclareMathOperator{\sgn}{sgn}
\title[Banded total positivity and bidiagonal factorizations]
{Banded Totally Positive Matrices and\\
Positive Bidiagonal Factorizations:\\
Existence, Stability, and Toeplitz Orthogonality}
\author{Manuel Ma\~nas\,\orcidlink{0000-0003-3764-5737}}
\address{Department of Theoretical Physics, Faculty of Physical Sciences,
Complutense University of Madrid, 28040 Madrid, Spain}
\email{manuel.manas@ucm.es}
\date{}
\newcommand{\LL}{\mathcal L}

\begin{document}
\begin{abstract}
Positive bidiagonal factorizations of banded totally positive
matrices of arbitrary finite bandwidth are products of unit lower and upper
bidiagonal factors with positive off-diagonal entries and a positive
diagonal matrix. In finite dimension these factorizations always exist,
with parameter domains described by polynomial inequalities.
For semi-infinite matrices, existence is equivalent to convergence of
explicit positive series. Every bounded banded totally positive matrix
is an operator-norm limit of matrices admitting these factorizations,
through perturbations confined to a fixed initial block.

Necessary and sufficient criteria characterize analytic dependence of
prescribed factors on matrix perturbations. With the initial entries fixed, a relative
entrywise criterion allows unbounded factors with different growth rates. For
uniformly positive bounded factors with variable initial entries, the
criterion involves every pair of factors within each triangular side
and the signed triangular inverses. The dimension of the local
factorization families is determined, and normalized LDU factors depend
analytically at every interior point of the factorization class.

For banded totally positive Toeplitz matrices, a stable positive
bidiagonal factorization exists if and only if the polynomial symbol
has distinct roots: every sufficiently small operator-norm perturbation
within the same band admits nearby positive factors depending
analytically on the perturbed matrix. Root formulas give two vector
polynomial families satisfying mixed-type multiple orthogonality,
paired by a rectangular matrix of positive measures. The initial
normalization depends on the chosen factorization. The smallest
attainable first off-diagonal entries of the two normalization matrices
are determined explicitly, and one factorization attains all these minima.
Although polynomial vectors and matrix measures depend on the
factorization, the associated weighted linear forms are independent of it. Their expansions in classical Chebyshev polynomials yield asymptotic
zero spacing inside the support and convergence of normalized kernels
to the sine kernel after rescaling by the local zero density.
\end{abstract}

\maketitle
\begingroup
\small
\setlength{\baselineskip}{9.6pt}\vspace{-22pt}
\tableofcontents
\endgroup

\section{Introduction}

A banded matrix has only finitely many nonzero diagonals. Banded
total positivity (BTP) requires every minor not forced to vanish
by the band pattern to be positive. A positive bidiagonal
factorization (PBF) expresses the matrix as a product of unit lower
bidiagonal factors, a positive diagonal matrix, and unit upper
bidiagonal factors, with positive entries on every allowed
off-diagonal. The questions addressed here are when these factorizations
exist, how they depend on the matrix, and how they determine spectral
data. The three main contributions are a global existence criterion, complete
criteria for analytic dependence through prescribed factors, and
an explicit spectral construction for BTP Toeplitz matrices and
finite boundary changes. Different admissible choices of initial factor entries give different
factorizations and different normalizations of polynomial vectors and
measures. In the weighted linear forms identified below, those changes
cancel, recovering spectral quantities determined by the matrix itself.

Already for bounded semi-infinite matrices, PBF existence is not
an open condition: arbitrarily small perturbations far from the
initial rows can destroy it. Existence, membership in the interior
of the PBF class, and analytic continuation through a prescribed
factorization are distinct properties. Here interiority and bounded stability
refer to the operator norm on the space of all real matrices with
the same band. The choice of initial factor entries is essential:
an interior matrix can have a PBF whose initial entries must vary
under perturbation.

\medskip\noindent
\emph{Global existence.}
In finite dimension, BTP and PBF are equivalent. Starting from
Metelmann's factorization, the construction replaces its initial zeros
by positive entries and solves triangular linear systems row by row,
preserving the product. This gives rational coordinates for all positive
factorizations and a nonempty, open, bounded domain of initial
parameters. The classical foundations are
\cite{Metelmann,Cryer1973,BFM2024}; related minor conditions and
finite factorizations appear in
\cite{deBoorPinkus,GascaMicchelliPena1992,GascaPena1995}.
Positive exchanges and word maps provide another route
\cite{FominZelevinsky2000,Lusztig1994,DavisHershMiller}.

In the semi-infinite setting, one choice of initial entries must
work in every row. For a BTP matrix with lower and upper bandwidths $p,q\ge1$, this is
possible if and only if $p+q-2$ explicit series of positive terms
converge to finite, strictly positive sums
(Theorem~\ref{thm:global-lower-criterion} and
Corollary~\ref{cor:global-two-sided-criterion}). The series compare adjacent Metelmann factors. The reciprocals of
their sums are strictly positive limits of decreasing inverse-minor
ratios. These limits give bounds uniform in the section size and
hence admissible global choices.
The criterion allows unbounded matrices and factors. Certified
tail bounds make the construction effective. Moreover, every
bounded BTP matrix is an operator-norm limit of bounded matrices
admitting PBFs, through perturbations confined to a fixed initial
block (Theorem~\ref{thm:boundary-approximation}). The constructed
factors are bounded when the original Metelmann factors are bounded.

In lower bandwidth two, the first subdiagonal normalizes the row
equations to a chain sequence in the sense of Wall and Chihara
\cite{Wall1948,Chihara1962}. Its maximal initial parameter gives the
endpoint in Proposition~\ref{rem:bandwidth-two-endpoint}. In higher
bandwidth, simultaneous convergence for all adjacent pairs gives one
choice of initial entries valid in every row.

The inverse-minor ratios are convergents of scalar Jacobi
continued fractions, extending the tetradiagonal criterion of
\cite{BFMTetra}. Branched fractions instead compute individual
factors after the initial entries have been chosen; their algebraic
form is related to \cite{PetreolleSokalZhu,Sokal2024}.
Positive factorizations of differential and difference operators
are studied in
\cite{Trench1974,Hartman1978,HarrisKruegerTrench1998,ErbeMathsenPeterson2001}.
The existence series incorporate the initial rows, which can
decide global existence even for a Toeplitz tail.

\medskip\noindent
\emph{Stability through prescribed factors.}
First suppose that all positive factor entries, including the diagonal
pivots, have common positive lower and finite upper bounds. The desired
continuation assigns positive factors to every nearby matrix, with that
matrix as their product, and recovers the prescribed PBF at the original
matrix. Real analytic dependence means that these factor sequences have
convergent power-series expansions in the perturbation, in the supremum
norm over all rows. In particular, their changes are bounded by a
constant times the operator norm of the perturbation.
With initial entries fixed, analytic dependence holds precisely
when products of adjacent factor ratios decay uniformly
geometrically over every block of rows, together with a condition
coupling the last lower
and upper factors
(Theorem~\ref{stab:thm:general-full-stability}). These products
measure how changes accumulate through successive rows.

Allowing the initial entries to vary permits a growing forward
recurrence to be solved backward. The necessary and sufficient
conditions then require uniform geometric decay of ratio products
or their reciprocals for every pair of factors, including
nonadjacent pairs, on each triangular side
(Theorem~\ref{stab:thm:free-separated-limits}). For the full matrix,
one must also control the uniquely determined normalized
decomposition $T=LDU$: a uniform bound on joint row sums of the
signed inverses of $L$ and $U^{\mathsf T}$ characterizes its
analytic dependence. Every operator-norm interior point of the PBF
class satisfies this LDU bound
(Proposition~\ref{stab:prop:interior-ldu}). Combining it with the pair
conditions gives
the complete criterion
(Proposition~\ref{stab:prop:free-full-transfer}). Under these
conditions, nearby factorizations of a fixed matrix form a
finite-dimensional analytic family.

With fixed initial entries, an entrywise relative norm also gives
a necessary and sufficient criterion for factors with arbitrary
positive growth, including different unbounded growth rates
(Theorem~\ref{stab:thm:relative-criterion}). Here smallness means
$|\widetilde T_{i,j}-T_{i,j}|\le\varepsilon T_{i,j}$ within the band,
with corresponding relative bounds for the factor changes. The
factorizations and perturbations in this relative criterion are defined
entry by entry.
The proofs use classical connections between bounded solutions of
first-order difference equations and geometric decay
\cite{BravermanKarabash,HuyMinh2001}, followed by the analytic inverse
function theorem \cite[Chapter~4]{Deimling}.
They develop the componentwise perturbation viewpoint of
\cite{Higham1994,BuenoDopico2004} for semi-infinite products.

\medskip\noindent
\emph{Toeplitz factorizations and spectral normalization.}
For a Toeplitz matrix $T=[t_{i-j}]_{i,j\in\NN}$ with positive
extreme coefficients, put $P(z)\coloneq\sum_{j=-q}^{p}t_jz^{j+q}$.
BTP is equivalent to all roots of this polynomial symbol being
negative, as follows from \cite{ASW}. Dividing these roots into two ordered groups constructs
a bounded PBF, including when roots repeat. The further condition that these roots be distinct is equivalent
to interiority and to the existence of a PBF whose factors depend
analytically on every sufficiently small operator-norm perturbation
within the same band
(Theorem~\ref{thm:toeplitz-stable-simple-roots}). Every ordered
partition of distinct roots permits variable initial entries;
a suitable ordering permits fixed ones. The construction uses
the positive upper--lower exchange of \cite{ManasPBF2026}.
Limits of triangular factors have antecedents in
\cite{Bauer1956,YoulaKazanjian1978,BottcherSilbermann2006}.

For each bounded Toeplitz PBF, the spectral Favard theorem
\cite{BFM2023,BFMUnbounded} supplies left and right vector
polynomials and a rectangular matrix of compactly supported
positive scalar measures. Their orthogonality is of mixed type:
polynomial vectors enter on both sides. The chosen PBF fixes their
initial conditions and normalization. The global existence series
bound the first off-diagonal entries of these initial-condition
matrices; inverse-minor limits express the bounds in the symbol
roots, and an ordered-root PBF attains them simultaneously
(Theorem~\ref{thm:toeplitz-normalization-minima}). For distinct roots, the minimizing factorization is locally isolated
among uniformly close PBFs of the same matrix, and it has a unique
nearby analytic continuation when its initial entries may vary.

Rational generating functions now give explicit polynomial vectors,
and root formulas remain valid at repeated roots by taking
derivatives. The matrix measure has interval support and an explicit
rank-one density: a scalar weight times the outer product of two vector
weights. Combining each polynomial vector
with its corresponding weight gives a scalar linear form.
Changes of PBF normalization cancel in these forms, so they and
the scalar weight depend only on the Toeplitz matrix
(Theorem~\ref{prop:spectral-functions}). This identifies the
common spectral data behind the different factorizations. The weights
depend on the spectral variable, so the linear forms are generally
algebraic functions, although their polynomial components are genuine
polynomials. Finite Chebyshev convolutions give explicit oscillatory
amplitudes; an identity relating these amplitudes to the density yields
zero spacing and the normalized sine-kernel limit, also for repeated
symbol roots.

Finite-section eigenvectors give another construction of the same
forms through ratios of determinantal polynomials. Schur formulas
\cite{MaximenkoMoctezuma2017,Alexandersson2012} connect these
polynomials with the vector equilibrium problem of
\cite{DuitsKuijlaars2008}. Together with
\cite{SchmidtSpitzer,Ullman}, this describes the limiting zero
measures; BTP places the equilibrium supports on alternating real
semiaxes. For bounded BTP matrices, convergence of the diagonals to those of a
BTP Toeplitz matrix already preserves the global limiting eigenvalue
measure. The reality and common compact bound of the finite spectra
allow the classical moment argument to apply. Local asymptotics are established below for finite boundary changes.

Changes confined to the initial $p\times q$ block are also accessible
whenever the modified matrix has a bounded PBF. A finite matrix inverse
gives their polynomial vectors, continuous density, and possible point
masses. For reordered constant-factor products, an exact interior
criterion involves distinct root parameters and the
separation of the limiting lower and upper parameters. On compact intervals inside the continuous support where an explicit
boundary determinant is nonzero, these changes modify the
oscillation phase and amplitude while preserving the interior zero
density and normalized sine-kernel limit. Scalar Jacobi
antecedents are \cite{GeronimoIliev2017,WebbOlver2021}.

Binomial and geometrically spaced roots provide explicit families
throughout the construction, including the systems of
\cite{Lima2023,Uniform} and Gaussian-binomial matrices
\cite{ArikanKilic2017}. Related higher-order recurrences, Nikishin
systems, symmetric factorizations, and moment-matrix constructions
are studied in
\cite{AptekarevKalyaginSaff2009,BarriosGeronimoLopez,DelvauxLopez,LimaSymmetric2026,BFMBidiagonal}.

Lean proofs cover the semi-infinite existence criterion for compatible
Metelmann factors, the Toeplitz BTP--PBF equivalence, and the bounded
and relative stability criteria specified in the verification table
\cite{BTPPBFVerification}. The spectral proofs in
Section~\ref{sec:mixed-systems} lie outside this formalization.

Sections~\ref{sec:finite}, \ref{chap:stability}, and
\ref{sec:mixed-systems} develop existence, stability, and the
spectral construction in this order.

\section[Existence of positive bidiagonal factorizations]{Banded total positivity and bidiagonal factorizations}\label{sec:finite}

The banded total positivity criterion and its Metelmann factorization
provide the starting point. Solving the row equations then replaces
the initial zeros by positive entries while preserving the finite
matrix. In the
semi-infinite setting, inverse-minor limits give the criterion for
one choice to work in every row. Initial-block perturbations provide
operator-norm approximation by matrices with PBF.
Section~\ref{chap:stability} studies the dependence of the factors on the matrix,
with fixed or variable initial entries.
Section~\ref{sec:mixed-systems} applies both criteria to Toeplitz
matrices and computes their polynomial families and measures.

\subsection{Positivity and Metelmann factors}\label{subsec:conventions}

The positivity and bandwidth conventions distinguish Metelmann's
initial zero entries from the positive entries inside its nonidentity
bidiagonal blocks.

Throughout, positive and nonnegative scalar quantities mean $>0$
and $\ge0$, respectively. Put $\mathbb N\coloneq \{1,2,\dots\}$
and $\NN\coloneq \{0,1,2,\dots\}$. Matrix entries are indexed from zero.
An integer index set $\{a,\dots,b\}$ is empty when $b<a$.
A matrix of lower bandwidth $p$ and upper
bandwidth $q$ satisfies $T_{i,j}=0$ outside $-q\le i-j\le p$.
For an increasing row list $(i_1,\ldots,i_r)$ and column list
$(j_1,\ldots,j_r)$, its minor is structurally admissible if
$-q\le i_k-j_k\le p$ for every $k\in\{1,\dots,r\}$.
Otherwise the band structure forces the determinant to vanish.
The matrix is \emph{banded totally positive} (BTP) when every
admissible minor is positive; it is \emph{totally
nonnegative} (TN) when all minors are nonnegative.
For a semi-infinite matrix these conditions concern all finite minors,
or equivalently all finite leading sections, with the inherited band.

A \emph{positive bidiagonal factorization} (PBF) means
\begin{equation}\label{eq:strict-definition}
 T=L_1\cdots L_p\,D\,U_q\cdots U_1,
\end{equation}
where every $L_i$ is unit lower bidiagonal with positive
subdiagonal, every $U_j$ is unit upper bidiagonal with positive superdiagonal, and $D$ is positive diagonal.
Here ``unit'' means that every diagonal entry of the triangular
factor equals one. The expression \emph{positive bidiagonal} always refers to its allowed
off-diagonal entries. Products in the semi-infinite
setting are well defined entrywise because there are finitely many
banded factors. The notation $\PBF$ also denotes the class of matrices
admitting such a factorization. The positive-minor criterion below
proves the implication $\PBF\Rightarrow\BTP$.

Write $T=LDU$ for the normalized triangular decomposition.
In \eqref{eq:strict-definition}, its unit triangular factors are
$L\coloneq L_1\cdots L_p$ and $U\coloneq U_q\cdots U_1$.
Thus $L$ and $U$ denote the complete triangular products, whereas
$L_i$ and $U_j$ denote the individual bidiagonal factors.
The same convention applies to the Metelmann products below.

For finite matrices of order $d$, assume $p,q\in\{1,\dots,d-1\}$.
Let $\BTP_{p,q}^{(d)}$ denote the class of such BTP matrices
and $\PBF_{p,q}^{(d)}$ the class admitting a PBF with $p$
lower and $q$ upper bidiagonal factors.

\begin{proposition}[PBF implies BTP, recalled from \cite{BFM2024}]
\label{prop:positive-product}
Let $L_1,\ldots,L_p$ be unit lower bidiagonal matrices of order $d$
with positive subdiagonals, let $U_1,\ldots,U_q$ be unit upper
bidiagonal matrices with positive superdiagonals, and let $D$ be
positive diagonal. For
$T\coloneq L_1\cdots L_pDU_q\cdots U_1$, every minor is nonnegative and
\[
 \det T[I,J]>0
 \quad\text{if and only if}\quad
 -q\le i_a-j_a\le p\quad(a\in\{1,\dots,r\}),
\]
where $T[I,J]$ denotes the submatrix with row indices
$I\coloneq (i_1<\cdots<i_r)$ and column indices $J\coloneq (j_1<\cdots<j_r)$,
in the indicated order. In particular, $T$ is BTP. The same implication holds for an
entrywise semi-infinite PBF.
\end{proposition}
\begin{proof}
Consider first the minors of a single lower bidiagonal factor.
A nonzero term in its determinant expansion must match every selected
row to a selected column differing from it by zero or one.
This matching preserves order. Indeed, if two rows $i<i'$ were
matched to columns $j>j'$, the bidiagonal restrictions would give
$j\le i<i'\le j'+1\le j$, which is impossible.
Thus at most one permutation contributes. When it contributes, it
is the identity permutation in the increasing row and column lists,
and its contribution is the product of the selected entries.
The minors are therefore nonnegative, and are positive whenever
all the matched entries are positive. Transposition proves the
same assertions for upper bidiagonal factors. The positive diagonal
matrix $D$ is also TN. Cauchy--Binet, applied successively to these
finitely many factors, now shows that every minor of $T$ is
nonnegative.

Suppose that the displayed band inequalities hold. Put
$k_a\coloneq \min(i_a,j_a)$ and $K\coloneq (k_1,\dots,k_r)$.
Both original lists increase by at least one at every step, so their
componentwise minimum does too. Hence $K$ is an increasing list.
The band inequalities imply
$0\le i_a-k_a\le p$ and $0\le j_a-k_a\le q$.
For the lower product, choose the intermediate lists
$I^{(t)}\coloneq [\max(k_a,i_a-t)]_{a\in\{1,\dots,r\}}$ for $t\in\{0,\dots,p\}$.
The componentwise maximum of these two strictly increasing integer
lists is strictly increasing. Also, $I^{(0)}=I$, $I^{(p)}=K$,
and the $a$th index decreases by either zero or one at each step.
Consequently, the minor of $L_t$ with row list $I^{(t-1)}$ and
column list $I^{(t)}$ is the positive product identified above.
These minors form one positive summand in the iterated
Cauchy--Binet expansion of $\det L[I,K]$. All its other summands
are nonnegative, so $\det L[I,K]>0$.

For the upper product, use
$J^{(t)}\coloneq [\min(j_a,k_a+t)]_{a\in\{1,\dots,r\}}$ for $t\in\{0,\dots,q\}$.
These lists are strictly increasing, start at $K$, end at $J$,
and increase each coordinate by zero or one at each step.
In the order $U_q,\dots,U_1$, they give a positive summand in
$\det U[K,J]$. Finally, the Cauchy--Binet expansion for $LDU$
contains the summand
$\det L[I,K](\prod_{a=1}^{r}D_{k_a,k_a})\det U[K,J]>0$.
This proves positivity for every structurally admissible minor.

For completeness, suppose that a band inequality fails at position
$a$. If $i_a>j_a+p$, then rows $i_a,\dots,i_r$ cannot be matched
to any of the first $a$ selected columns: all those entries lie
below the band. There are $r-a+1$ such rows but only $r-a$
remaining columns. If $i_a<j_a-q$, the first $a$ selected rows
cannot be matched to any of the last $r-a+1$ columns; only $a-1$
columns remain for them. In either case, every permutation in the
determinant expansion contains a zero entry. The minor vanishes.

This proves the finite statement of
\cite[Proposition~2.21]{BFM2024} directly. For the semi-infinite
assertion, choose a leading section containing the given minor
and of order greater than $\max(p,q)$. Leading sections commute
with multiplication of lower triangular matrices, and also with
multiplication of upper triangular matrices. In an entry of $LDU$,
the intermediate diagonal index is at most both the row and the
column index. Thus its leading section is exactly the product of
the corresponding leading sections of all the factors. Applying
the finite argument proves the assertion for the given minor.
\end{proof}

With $i\in\{1,\dots,p\}$ and $j\in\{1,\dots,q\}$,
Metelmann's factorization takes the form
\begin{equation}\label{eq:metelmann}
\begin{aligned}
 T&=\widehat L_1\cdots\widehat L_p\,D\,
          \widehat U_q\cdots\widehat U_1,\\
 \widehat L_i&=
 \begin{bmatrix}I_{p-i}&0\\0&L_i^{\rm M}\end{bmatrix},
 &L_i^{\rm M}&\in\RR^{(d-p+i)\times(d-p+i)},\\
 \widehat U_j&=
 \begin{bmatrix}I_{q-j}&0\\0&U_j^{\rm M}\end{bmatrix},
 &U_j^{\rm M}&\in\RR^{(d-q+j)\times(d-q+j)}.
\end{aligned}
\end{equation}
Here $I_r$ is the $r\times r$ identity matrix, $I_0$ is the empty
block, and every displayed full factor has size $d\times d$.
The off-diagonal zero blocks have the corresponding rectangular
sizes. The blocks $L_i^{\rm M},U_j^{\rm M}$ are unit bidiagonal with
positive allowed off-diagonal entries; $D$ is $d\times d$
positive diagonal. In the semi-infinite notation below, $I$ denotes
the identity on $\ell^2(\NN)$.

The following characterization is \cite[Theorem~2.16(iv)]{BFM2024},
based on Metelmann's factorization \cite{Metelmann}.
\begin{proposition}[Metelmann characterization, recalled from \cite{BFM2024}]
\label{prop:metelmann-input}
For a matrix of order $d$ with $p,q\in\{1,\dots,d-1\}$, BTP is equivalent to
the existence of \eqref{eq:metelmann}, with $D$ positive diagonal
and every allowed off-diagonal entry of each bidiagonal block
$L_i^{\rm M}$ and $U_j^{\rm M}$ positive.
\end{proposition}
The initial identity blocks impose the zeros in
\eqref{eq:metelmann}; positivity is required only inside the nonidentity
bidiagonal blocks. The refactorization below removes these initial
zeros while preserving $T$.

A finite matrix is called \emph{oscillatory} if it is totally
nonnegative and some positive integer power has all its minors
strictly positive. This condition gives the spectral information
needed below.

\begin{corollary}[Classical oscillation and spectrum]
\label{cor:btp-oscillatory}
Every finite BTP matrix is oscillatory and has distinct positive
eigenvalues.
In particular, this conclusion applies to every finite leading
section of a BTP matrix of bandwidth $(p,q)$ with $p,q\in\mathbb N$.
\end{corollary}
\begin{proof}
The matrix is totally nonnegative, its determinant is positive,
and every entry on its first upper and lower subdiagonals is positive.
The Gantmacher--Krein criterion therefore makes it oscillatory
\cite{GantmacherKrein}. The positive simple spectrum follows from the classical spectral
theorem for oscillatory matrices; see \cite{GantmacherKrein} or
\cite[Theorem~1.4]{BFM2023}.
\end{proof}

\subsection{Refactorization with positive parameters}\label{subsec:refactorization}

The initial identity blocks in \eqref{eq:metelmann} leave some subdiagonal entries
equal to zero. Replacing these zeros by small positive parameters
and adjusting the other entries will keep the entire lower
triangular product fixed. The construction proceeds row by
row: at each row the new parameters are free, and a triangular
linear system determines the remaining entries. Applying the same
construction to the transpose of the upper product will then give
a PBF of the original matrix. After the existence theorem, an
argument determining one factor at a time will identify the associated linear
recurrences and their branched continued-fraction representation.

\begin{lemma}[Refactorization of a fixed lower product, row by row]
\label{lem:rowwise}
Let $p\in\{1,\dots,d-1\}$ and
$L=\widehat L_1\cdots\widehat L_p$, where the lower factors have the
pattern in \eqref{eq:metelmann}. There are analytic families of unit
lower bidiagonal matrices $L_i(\boldsymbol u)$, near $\boldsymbol u=0$,
such that $L=L_1(\boldsymbol u)\cdots L_p(\boldsymbol u)$ and
$L_i(0)=\widehat L_i$.
The parameter vector $\boldsymbol u$ has $p(p-1)/2$ entries, one for each
initial zero subdiagonal parameter. For all sufficiently small
entrywise positive $\boldsymbol u$, every factor is positive
bidiagonal.
\end{lemma}
\begin{proof}
For $i\in\{1,\dots,p\}$ and $n\in\{1,\dots,d-1\}$, denote the
subdiagonal entries by $a_i(n;\boldsymbol u)\coloneq [L_i(\boldsymbol u)]_{n,n-1}$
and $a_i^0(n)\coloneq a_i(n;0)=[\widehat L_i]_{n,n-1}$.
The superscript $0$ denotes evaluation at the initial parameter
$\boldsymbol u=0$. Write $a_i(n)$ for $a_i(n;\boldsymbol u)$ when
the dependence on $\boldsymbol u$ is understood. At $\boldsymbol u=0$, the initial
identity blocks in \eqref{eq:metelmann} give
\begin{equation}\label{eq:metelmann-parameters}
 a_i^0(n)=0\quad(n\le p-i),\qquad
 a_i^0(n)>0\quad(n>p-i).
\end{equation}
Proceed in increasing row order $n$ and put $k\coloneq \min(p,n)$.
By \eqref{eq:metelmann-parameters}, an initial entry $a_i^0(n)$
vanishes exactly when $n\le p-i$, equivalently $i\le p-n$.
There are therefore $\max(p-n,0)=p-k$ such entries. Choose
$a_1(n),\ldots,a_{p-k}(n)$ freely near zero; if $n\ge p$, this
list is empty. The remaining $k$ entries are initially positive.
They are determined by requiring row $n$ of $L_1\cdots L_p$ to
agree with row $n$ of the fixed matrix $L$.

\emph{Derivation of the row equations.}
Interpreting an empty product as the identity gives
$L_1\cdots L_p-I=\sum_{i=1}^p (L_i-I)L_{i+1}\cdots L_p$.
Indeed, the $i$th summand is
$L_i\cdots L_p-L_{i+1}\cdots L_p$, so the intermediate products
cancel in the sum. Since $L_i$ is unit lower bidiagonal, row $n$
of $L_i-I$ has just one possible nonzero entry: $a_i(n)$ in
column $n-1$. Hence the entry $(n,n-r)$ of its product with
$L_{i+1}\cdots L_p$ is
$a_i(n)[L_{i+1}\cdots L_p]_{n-1,n-r}$.
For $r\ge1$, the identity matrix has zero in position $(n,n-r)$.
Matching that entry of $L_1\cdots L_p$ with the fixed entry
$L_{n,n-r}$ therefore gives
\begin{equation}\label{eq:row-system}
 L_{n,n-r}=\sum_{i=1}^p
 a_i(n)[L_{i+1}\cdots L_p]_{n-1,n-r},
 \qquad r\in\{1,\dots,k\}.
\end{equation}
The construction proceeds in increasing row order. At row $n$, the
entries $a_j(t)$ with $1\le t<n$ have already been determined, and
\eqref{eq:row-system} determines the remaining entries of row $n$. The left-hand sides are known,
because they are entries of the given matrix $L$.

The coefficient multiplying $a_i(n)$ is the entry $(n-1,n-r)$
of $L_{i+1}\cdots L_p$. To calculate it, multiply the leading
$n\times n$ blocks of these factors, consisting of rows and
columns $0,\dots,n-1$. Lower triangularity ensures that this gives
the same entry as the full product. These blocks contain only
entries from earlier rows, so the coefficient is known before
row $n$ is determined.

This computation is illustrated for three factors in
Example~\ref{ex:three-factor-row-equation}.

This reasoning applies to every equation in
\eqref{eq:row-system}: each current entry $a_i(n)$ is multiplied
by a coefficient determined entirely by the preceding rows.

\emph{The triangular system.}
Recall that $k=\min(p,n)$. In row $n$, the entries
$a_i(n)$ with $i\in\{1,\dots,p-k\}$ have been chosen freely. The remaining $k$
entries, with $i\in\{p-k+1,\dots,p\}$, are the unknowns to be
determined by \eqref{eq:row-system}. First move the terms
containing the chosen entries to the right-hand side and denote
the resulting known quantity by
\[
 y_r\coloneq L_{n,n-r}-\sum_{i=1}^{p-k}
 a_i(n)[L_{i+1}\cdots L_p]_{n-1,n-r},
 \qquad r\in\{1,\dots,k\}.
\]
If $p=k$, the sum is empty and $y_r=L_{n,n-r}$.

Number the remaining unknowns in reverse factor order:
$x_j\coloneq a_{p-j+1}(n)$ for $j\in\{1,\dots,k\}$.
Thus $x_1=a_p(n)$ and $x_k=a_{p-k+1}(n)$.
Substituting $i=p-j+1$ in \eqref{eq:row-system} gives
\[
 \sum_{j=1}^k C_{r,j}x_j=y_r,
 \qquad
 C_{r,j}\coloneq[L_{p-j+2}\cdots L_p]_{n-1,n-r},
 \qquad r\in\{1,\dots,k\}.
\]
Writing $x$ and $y$ as column vectors, these equations are $Cx=y$.

To see why this ordering gives an upper triangular matrix, consider
$C_{r,j}$. It is an entry of a product of $j-1$ lower bidiagonal
factors, and its position $(n-1,n-r)$ lies $r-1$ places below
the diagonal. To obtain an entry this far below the diagonal,
at least $r-1$ factors must contribute a subdiagonal entry.
There are only $j-1$ factors available. Consequently, if $j<r$,
the coefficient $C_{r,j}$ is zero. The $r$th equation therefore
contains only $x_r,\dots,x_k$, which is precisely the upper
triangular form.

For a diagonal coefficient $C_{r,r}$, there are exactly $r-1$
factors. Each must contribute its subdiagonal entry, so there is
just one term in the product. The first factor contributes its
entry in row $n-1$, the second in row $n-2$, and so on, down to
row $n-r+1$. Hence the diagonal coefficient, or pivot, is
\begin{equation}\label{eq:pivot}
 C_{r,r}=[L_{p-r+2}\cdots L_p]_{n-1,n-r}
 =\prod_{s=1}^{r-1}a_{p-r+1+s}(n-s).
\end{equation}
For $r=1$, the empty product gives $C_{1,1}=1$.
For $r\ge2$, positivity at $\boldsymbol u=0$ follows from
\eqref{eq:metelmann-parameters}. The entry
$a_{p-r+1+s}(n-s)$ is positive there whenever
$n-s>p-(p-r+1+s)$. Cancelling $s$ reduces this condition to
$n>r-1$, which holds because $r\le k\le n$.
Thus $C$ is upper triangular and all its diagonal entries are
positive at the initial factorization.

\emph{Back substitution and positivity.}
The unique solution is obtained by taking
$r\in\{k,k-1,\dots,1\}$ in decreasing order:
\[
 x_r=\frac{y_r-\sum_{j=r+1}^{k}C_{r,j}x_j}{C_{r,r}}.
\]
Inductively, the coefficients and the right-hand sides are analytic
functions of the free parameters chosen in this and earlier rows.
Since the pivots are nonzero at $\boldsymbol u=0$, back substitution is
analytic in a neighborhood of that point. At $\boldsymbol u=0$ the
prescribed entries solve the same system, so uniqueness recovers
$L_i(0)=\widehat L_i$ at every row. The diagonal is automatically
one; entries outside the lower bandwidth vanish automatically.
Thus matching these $k$ entries in every row matches the whole
product $L$.

There are finitely many rows and hence finitely many pivots and
originally positive entries. Intersect their neighborhoods of
positivity. Within this common neighborhood all solved entries
remain positive; choosing every free entry positive makes every
subdiagonal entry of every factor positive. The free entries are
exactly the $p-n$ initial zeros in each row $n\in\{1,\dots,p-1\}$, so their
total number is
$\sum_{n=1}^{p-1}(p-n)=p(p-1)/2$.
\end{proof}

\begin{example}[The row equations for three factors]
\label{ex:three-factor-row-equation}
Take $p=3$ and write $a_i(n)\coloneq [L_i]_{n,n-1}$ for the entries
of the factors in Lemma~\ref{lem:rowwise}. For $n\ge3$,
$[L_2L_3]_{n-1,n-2}=a_2(n-1)+a_3(n-1)$, so
\eqref{eq:row-system} gives the second-subdiagonal equation
\[
 L_{n,n-2}
 =\bigl(a_2(n-1)+a_3(n-1)\bigr)a_1(n)
  +a_3(n-1)a_2(n).
\]
The quantities $a_2(n-1)+a_3(n-1)$ and $a_3(n-1)$ have already
been calculated; they are the coefficients of the current entries
$a_1(n)$ and $a_2(n)$. Thus this is a linear equation for row $n$.
\end{example}

\begin{theorem}[Finite BTP equals PBF]\label{thm:finite}
For every $d$ and $p,q\in\{1,\dots,d-1\}$,
\begin{equation}\label{eq:finite-equivalence}
 \BTP_{p,q}^{(d)}=\PBF_{p,q}^{(d)}.
\end{equation}
More precisely, for every $T\in\BTP_{p,q}^{(d)}$, its Metelmann
factorization \eqref{eq:metelmann} is a limit of PBFs of
\emph{the same matrix $T$}, with convergence
of each individual bidiagonal factor.
Define its unit triangular factors by $L\coloneq \widehat L_1\cdots\widehat L_p$
and $U\coloneq \widehat U_q\cdots\widehat U_1$, so that $T=LDU$.
Along the approximating PBFs, $L_1\cdots L_p=L$ and
$U_q\cdots U_1=U$, and $D$ is unchanged. Thus it is the normalized
LDU decomposition of $T$ that remains fixed; only the individual
bidiagonal factors of $L$ and $U$ vary.
\end{theorem}
\begin{proof}
For $T\in\BTP$, Proposition~\ref{prop:metelmann-input} gives
\eqref{eq:metelmann}. Apply Lemma~\ref{lem:rowwise} to its lower
product. Apply the same lemma, with $p$ replaced by $q$, to the
transpose of the upper product, namely
$\widehat U_1^{\mathsf T}\cdots\widehat U_q^{\mathsf T}$.
The first application gives $L_1(\boldsymbol u)\cdots L_p(\boldsymbol u)=L$.
For the second application, write
$U^{\mathsf T}=M_1(\boldsymbol{\theta})\cdots M_q(\boldsymbol{\theta})$ and set
$U_j(\boldsymbol{\theta})\coloneq M_j(\boldsymbol{\theta})^{\mathsf T}$. Transposing the product reverses
its order, so $U_q(\boldsymbol{\theta})\cdots U_1(\boldsymbol{\theta})=U$. Thus
\[
 T=L_1(\boldsymbol u)\cdots L_p(\boldsymbol u)D
 U_q(\boldsymbol{\theta})\cdots U_1(\boldsymbol{\theta})
\]
for all parameters in the neighborhoods supplied by the lemma.
Choose all components of $\boldsymbol u$ and $\boldsymbol{\theta}$ equal to
$\varepsilon>0$. Since both parameter spaces are finite-dimensional,
one sufficiently small $\varepsilon$ belongs to both neighborhoods.
Every allowed bidiagonal entry is then positive, while $L,U,D$ remain
unchanged. As $\varepsilon\downarrow0$, the analytic factors converge
entry by entry to the prescribed Metelmann factors. If one parameter
space has dimension zero, its factors are already strictly positive
and remain fixed.

Conversely, Proposition~\ref{prop:positive-product} applies to any PBF and makes every
structurally admissible minor positive. This proves the reverse
inclusion and hence the equality of the two classes.
\end{proof}

\begin{proof}[Alternative proof of \eqref{eq:finite-equivalence}]
The existence of a positive factorization also follows by rearranging
the elementary factors in a Metelmann product. The following explicit
rearrangement includes the changes that replace its initial zeros by
positive entries.

Put $m\coloneq d-1$ and, using the coordinate vectors $e_0,\dots,e_m$,
define $\mathsf X_i(t)\coloneq I+t e_i e_{i-1}^{\mathsf T}$ for
$i\in\{1,\dots,m\}$ and $t>0$.
The product $\mathsf X_1(t_1)\cdots \mathsf X_m(t_m)$ is the unit lower
bidiagonal matrix with subdiagonal entries $t_1,\dots,t_m$.
For $a,b,c>0$, the following identities rearrange such elementary factors:
\begin{subequations}\label{eq:elementary-refactorization}
\begin{align}
 \mathsf X_i(a)\mathsf X_j(b)&=\mathsf X_j(b)\mathsf X_i(a),
 &&|i-j|>1,\label{eq:elementary-refactorization-commutation}\\
 \mathsf X_i(a)\mathsf X_j(b)\mathsf X_i(c)
 &=\mathsf X_j\!\left(\frac{bc}{a+c}\right)
   \mathsf X_i(a+c)
   \mathsf X_j\!\left(\frac{ab}{a+c}\right),
 &&|i-j|=1,\label{eq:elementary-refactorization-braid}\\
 \mathsf X_i(a)\mathsf X_i(b)&=\mathsf X_i(a+b).
 \label{eq:elementary-refactorization-addition}
\end{align}
\end{subequations}
Identities \eqref{eq:elementary-refactorization-commutation} and
\eqref{eq:elementary-refactorization-braid} are the elementary
relations used in \cite[proof of Theorem~23, p.~30]{FominZelevinsky2000};
they can also be checked directly by multiplication.
The three-parameter identity \eqref{eq:elementary-refactorization-braid}
is equivalent to the scale-invariant
form of the Yang--Baxter relation in
\cite[Appendix~A, Proposition~A.1 and Eq.~(A.5), p.~146]{BerensteinFominZelevinsky1996}.
Set $u\coloneq a$, $v\coloneq c$, and $\lambda\coloneq b/(a+c)$.
Then it reads
\[
 \mathsf X_i(u)\mathsf X_j\!\bigl(\lambda(u+v)\bigr)\mathsf X_i(v)
 =\mathsf X_j(\lambda v)\mathsf X_i(u+v)\mathsf X_j(\lambda u),
 \qquad |i-j|=1.
\]
The three positive parameters $u,v,\lambda$ are independent:
the inverse change is $a=u$, $b=\lambda(u+v)$, and $c=v$.
The choice $\lambda=1$ gives the usual additive Yang--Baxter relation.
In the present proof, the identity with arbitrary positive $a,b,c$
reorganizes the factors while preserving their product and the
positivity of their parameters.
The identities \eqref{eq:elementary-refactorization} can be used in
either direction: \eqref{eq:elementary-refactorization-braid} has the
same form with $i$ and $j$ interchanged, and
\eqref{eq:elementary-refactorization-addition} is reversed by splitting
a positive parameter into two positive summands.

\emph{Initial and final index lists.}
Record a product of elementary matrices by the list of their
indices. Each occurrence of an index carries its own positive
parameter, which changes according to \eqref{eq:elementary-refactorization}.
Put $\mathsf C_a\coloneq(a,a+1,\dots,m)$ for
$a\in\{1,\dots,m\}$, and let $\mathsf C_{m+1}$ be the empty list.
A product with list $\mathsf C_a$ is a unit lower bidiagonal
matrix with zero subdiagonal entries in rows $1,\dots,a-1$ and
positive entries in rows $a,\dots,m$. In particular, the full
list $\mathsf C_1=(1,\dots,m)$ gives a bidiagonal matrix with
all its subdiagonal entries positive.

Juxtaposition denotes concatenation of lists, and a single index
$i$ in a concatenation denotes the one-entry list $(i)$.
Thus $\mathsf C_a i=(a,a+1,\dots,m,i)$, whereas
$(i+1)\mathsf C_a=(i+1,a,a+1,\dots,m)$.
Write $A\sim B$ when the identities \eqref{eq:elementary-refactorization}
transform a product with index list $A$ into a product with index
list $B$, for every choice of positive parameters, and conversely.
The parameters may change, but they remain positive and the resulting matrix is unchanged.
In particular, $(i,i)\sim(i)$: two consecutive factors combine by
\eqref{eq:elementary-refactorization-addition} as
$\mathsf X_i(a)\mathsf X_i(b)=\mathsf X_i(a+b)$.
Conversely, a factor $\mathsf X_i(t)$ can be split into
$\mathsf X_i(s)\mathsf X_i(t-s)$ for any $s\in(0,t)$.
Empty portions of the displayed lists are omitted.

For $T\in\BTP_{p,q}^{(d)}$, take the Metelmann factorization
from Proposition~\ref{prop:metelmann-input}.
The first positive subdiagonal entry of $\widehat L_j$ occurs
in row $p-j+1$, so its elementary factors have list
$\mathsf C_{p-j+1}$. Consequently the lower product
$L=\widehat L_1\cdots\widehat L_p$ has list
$\mathsf M_p$, where
$\mathsf M_r\coloneq\mathsf C_r\mathsf C_{r-1}\cdots\mathsf C_1$
for $r\in\{1,\dots,m\}$.
The desired positive factorization has list $\mathsf C_1^p$:
each of its $p$ copies of $\mathsf C_1$ represents one bidiagonal
factor with all its subdiagonal entries positive.
It therefore suffices to transform $\mathsf M_p$ into
$\mathsf C_1^p$ while preserving the matrix product.
The equivalence follows by reducing $\mathsf C_1^p$ to
$\mathsf M_p$ and then reversing the changes.

\emph{Combining two copies of the same list.}
Example~\ref{ex:three-index-list-rearrangement} illustrates the
rearrangement for $m=3$.

For arbitrary $a\in\{1,\dots,m\}$, the statement to prove is
$\mathsf C_a^2\sim\mathsf C_{a+1}\mathsf C_a$.
To establish it, move the indices
$a,\dots,m-1$ of the second copy across the first copy,
then combine the two final factors with index $m$.
For $a<m$, the rearrangement needed for each
$i\in\{a,\dots,m-1\}$ is
$\mathsf C_a i\sim(i+1)\mathsf C_a$.

Here are the individual moves. Start with
$\mathsf C_a i=(a,\dots,m,i)$ and move the final $i$ to the left
past $m,m-1,\dots,i+2$. The corresponding factors commute by
\eqref{eq:elementary-refactorization-commutation}, because their
indices differ from $i$ by at least two.
The resulting list is
$(a,\dots,i-1,i,i+1,i,i+2,\dots,m)$.
Apply \eqref{eq:elementary-refactorization-braid} to the consecutive indices
$(i,i+1,i)$, replacing them by $(i+1,i,i+1)$.
The first of these new indices, $i+1$, now commutes to the left
past $i-1,\dots,a$. The resulting list is
$(i+1,a,\dots,m)=(i+1)\mathsf C_a$.

Apply this rearrangement to the indices
$a,a+1,\dots,m-1$ of the second copy of $\mathsf C_a$ in $\mathsf C_a^2=\mathsf C_a\mathsf C_a$.
At each step, one index moves into the list on the left,
while a copy of $\mathsf C_a$ remains in the middle.

For the first step, separate the initial index $a$ of the
second copy. The block $\mathsf C_a\,a$ can be replaced by
$(a+1)\mathsf C_a$, so
\[
 \mathsf C_a^2
 =\underbrace{\mathsf C_a\,a}_{\sim\,(a+1)\mathsf C_a}
   (a+1,\dots,m)\sim(a+1)\mathsf C_a(a+1,\dots,m).
\]
The next index immediately after the middle copy of
$\mathsf C_a$ is now $a+1$.
When $a\le m-2$, the second step replaces the block
$\mathsf C_a(a+1)$ by $(a+2)\mathsf C_a$:
\[
(a+1)\underbrace{\mathsf C_a(a+1)}_{\sim\,(a+2)\mathsf C_a}
   (a+2,\dots,m)\sim(a+1,a+2)\mathsf C_a(a+2,\dots,m).
\]
The list on the left has gained the index $a+2$, and the
remaining list on the right now starts with $a+2$.

More generally, just before the step with
$i\in\{a,\dots,m-1\}$, the list on the left is
$(a+1,\dots,i)$ and the remaining list on the right is
$(i,\dots,m)$. The left list is empty at the first step, $i=a$.
Applying the rearrangement to the block $\mathsf C_a\,i$ gives
\[
 (a+1,\dots,i)
   \underbrace{\mathsf C_a\,i}_{\sim\,(i+1)\mathsf C_a}
   (i+1,\dots,m)\sim(a+1,\dots,i+1)\mathsf C_a(i+1,\dots,m).
\]
The final $i+1$ in the left list is produced by this replacement.
The first $i+1$ in the right list is the next index of the
original second copy, still to the right of $\mathsf C_a$.
At $i=m-1$ the resulting list is $\mathsf C_{a+1}\mathsf C_a m$.
Since $\mathsf C_a$ ends in $m$, the last two indices are $(m,m)$.
Combining their matrices by \eqref{eq:elementary-refactorization-addition},
$\mathsf X_m(s)\mathsf X_m(t)=\mathsf X_m(s+t)$, gives
$\mathsf C_{a+1}\mathsf C_a$, as required.
For $a=m$, this combination alone proves the result:
$\mathsf C_m^2=(m,m)\sim(m)=\mathsf C_{m+1}\mathsf C_m$.

\emph{Reducing the full list to the Metelmann list.}
Appending one full copy of $\mathsf C_1$ to
$\mathsf M_r$ produces a list equivalent to $\mathsf M_{r+1}$, for
$r\in\{1,\dots,m-1\}$.
In $\mathsf M_r\mathsf C_1$, the two rightmost copies form
$\mathsf C_1^2$. Replace them by $\mathsf C_2\mathsf C_1$.
If $r>1$, the new $\mathsf C_2$ is adjacent to the existing
copy of $\mathsf C_2$, so replace
$\mathsf C_2^2$ by $\mathsf C_3\mathsf C_2$.
More generally, before the step with index $a$, the list is
$\mathsf C_r\cdots\mathsf C_{a+1}\mathsf C_a^2\mathsf C_{a-1}\cdots\mathsf C_1$.
For $a<r$, replacing $\mathsf C_a^2$ by
$\mathsf C_{a+1}\mathsf C_a$ creates the next pair
$\mathsf C_{a+1}^2$. At $a=r$, the replacement gives
$\mathsf C_{r+1}\mathsf C_r\cdots\mathsf C_1=\mathsf M_{r+1}$.
Thus $\mathsf M_r\mathsf C_1\sim\mathsf M_{r+1}$.
Starting from $\mathsf M_1=\mathsf C_1$ and adding a full copy
at each step proves $\mathsf C_1^p\sim\mathsf M_p$.

Return to the given matrix $L$, whose Metelmann
factors have list $\mathsf M_p$.
Reverse the finite sequence of changes just established, splitting
each positive sum into two positive summands when a combination
is undone. Every step preserves the matrix product and keeps all
parameters positive. The final list is $\mathsf C_1^p$.
Grouping its $p$ full copies of $\mathsf C_1$ gives
$L=L_1\cdots L_p$, with every subdiagonal entry of every factor
strictly positive.

Apply the same argument with $q$ in place of $p$ to
$U^{\mathsf T}=\widehat U_1^{\,\mathsf T}\cdots
\widehat U_q^{\,\mathsf T}$. Transposing the resulting product
gives $U=U_q\cdots U_1$ with positive upper bidiagonal factors.
Together with the unchanged positive diagonal $D$, these products
give a PBF of $T$. The converse follows from
Proposition~\ref{prop:positive-product}, completing the proof of
\eqref{eq:finite-equivalence}.

The triangular-system proof above also provides convergence to a
prescribed Metelmann factorization and permits every sufficiently
small positive choice of the initial parameters.
\end{proof}

\begin{example}[Combining two lists with three indices]
\label{ex:three-index-list-rearrangement}
Use the index lists and positive refactorization identities from
the alternative proof of Theorem~\ref{thm:finite}. For $m=3$, the lists are
$\mathsf C_1=(1,2,3)$ and $\mathsf C_2=(2,3)$.
The desired transformation takes
$\mathsf C_1^2=(1,2,3,1,2,3)$ into
$\mathsf C_2\mathsf C_1=(2,3,1,2,3)$.

First commute the adjacent indices $3,1$ using
\eqref{eq:elementary-refactorization-commutation}, then apply
\eqref{eq:elementary-refactorization-braid} to the initial
triple $(1,2,1)$:
\[
 (1,2,3,1,2,3)
 \sim (1,2,1,3,2,3)
 \sim (2,1,2,3,2,3).
\]
The resulting list is $2\mathsf C_1(2,3)$.
Thus the first index of the second copy has passed across
the first copy, changing from $1$ to $2$.

Next apply \eqref{eq:elementary-refactorization-braid}
to the triple $(2,3,2)$ in positions three through five,
and commute the adjacent indices $1,3$:
\[
 (2,1,2,3,2,3)
 \sim (2,1,3,2,3,3)
 \sim (2,3,1,2,3,3).
\]
The resulting list is $\mathsf C_2\mathsf C_1\,3$.
The last two indices are both $3$, so
\eqref{eq:elementary-refactorization-addition} combines
their factors:
\[
 (2,3,1,2,3,3)\sim(2,3,1,2,3)
 =\mathsf C_2\mathsf C_1.
\]
Every step preserves the matrix product and keeps the
parameters positive.
\end{example}

\begin{remark}[Free parameters and earlier factorization results]
\label{rem:dhm}
The reduction from $\mathsf C_1^p$ to $\mathsf M_p$ combines
$r$ pairs of equal indices at the step
$\mathsf M_r\mathsf C_1\sim\mathsf M_{r+1}$,
for $r\in\{1,\dots,p-1\}$. It therefore makes exactly
$s\coloneq p(p-1)/2$ such combinations.
Each combination in \eqref{eq:elementary-refactorization-addition}
is reversed by splitting its positive parameter $c$ into
$\nu c$ and $(1-\nu)c$, with $\nu\in(0,1)$.
These $s$ choices explain the number of free parameters in the
positive factorizations of a fixed lower product $L$.
This count agrees with the general description of finite positive
factorization families by Davis, Hersh, and Miller
\cite[Lemma~2.24 and Theorem~4.70]{DavisHershMiller}.

The initial subdiagonal entries in Lemma~\ref{lem:rowwise}
provide coordinates that can be kept fixed across leading sections.
The uniform bounds in Lemma~\ref{lem:global-inverse-minors}, obtained
from the adjacent series \eqref{eq:global-adjacent-series}, provide
positive choices valid in every row; Theorem~\ref{thm:global-lower-criterion}
characterizes their existence.
\end{remark}

Lemma~\ref{lem:rowwise} gives positive factorizations for sufficiently
small positive initial entries. The full set of admissible choices
of these entries describes every ordered positive factorization
of the same lower product $L$.

Fix $p\in\{2,\dots,d-1\}$ and
$L=\widehat L_1\cdots\widehat L_p$ as in
Lemma~\ref{lem:rowwise}. For the unit lower bidiagonal factors
$L_i$ to be constructed, write $a_i(n)\coloneq[L_i]_{n,n-1}$.
Their freely chosen initial entries occupy the positions that
are zero in the Metelmann factors, listed in
\eqref{eq:metelmann-parameters}. Set $s\coloneq p(p-1)/2$ and
collect these entries in the vector
$\boldsymbol u\coloneq[a_i(n)]_{n\in\{1,\dots,p-1\},\,i\in\{1,\dots,p-n\}}$.
Let $\mathcal D_L\subset(0,\infty)^s$ consist of the initial
vectors for which the row-by-row construction
\eqref{eq:row-system} can be carried out through row $d-1$ and produces
strictly positive subdiagonal entries in every factor.

\begin{proposition}[Initial parameters for all positive factorizations]
\label{prop:positive-parameter-domain}
The construction using \eqref{eq:row-system} gives a bijection
\[
 \boldsymbol u\in\mathcal D_L\longmapsto(L_1(\boldsymbol u),\ldots,L_p(\boldsymbol u))
\]
onto all ordered factorizations $L=L_1\cdots L_p$ by unit positive
lower bidiagonal matrices. The admissible set $\mathcal D_L$ is
nonempty, bounded, and open, and is described by finitely many
strict polynomial inequalities.
The factor entries are rational functions of $\boldsymbol u$.
The inverse map reads the chosen initial entries from the factors.
Thus the factorizations depend analytically on exactly $s$
independent parameters: the derivative of the displayed map
has rank $s$ at every point of $\mathcal D_L$.
\end{proposition}
\begin{proof}
Calculate the factor entries in increasing row order.
In row $n$, put $k\coloneq \min(p,n)$ and order the unknown entries as
$x_j\coloneq a_{p-j+1}(n)$, $j\in\{1,\dots,k\}$.
As in the proof of Lemma~\ref{lem:rowwise}, moving the chosen
initial entries to the right-hand side of \eqref{eq:row-system}
gives an upper triangular system $Cx=y$.
Solve it in the order $x_k,x_{k-1},\dots,x_1$.

Suppose every subdiagonal entry already constructed has been
written as a quotient of polynomials with positive numerator and
denominator on the parameters retained so far. Initially this
holds because each free entry is a positive coordinate of
$\boldsymbol u$, with denominator $1$.
For the next unknown $x_r$, take a common denominator and write
\[
 y_r-\sum_{j=r+1}^{k}C_{r,j}x_j
 =\frac{N(\boldsymbol u)}{Q(\boldsymbol u)},
 \qquad Q(\boldsymbol u)>0.
\]
Here $Q$ is a product of the positive denominators of entries
already constructed. The pivot \eqref{eq:pivot} is a product of
$r-1$ earlier subdiagonal entries. Write these entries as
$P_\ell/R_\ell$, with $P_\ell,R_\ell>0$ on the retained parameters.
Back substitution now gives the explicit formula
\[
 x_r
 =\frac{y_r-\sum_{j=r+1}^{k}C_{r,j}x_j}{C_{r,r}}
 =\frac{N\prod_{\ell=1}^{r-1}R_\ell}
        {Q\prod_{\ell=1}^{r-1}P_\ell}.
\]
The products are $1$ when $r=1$.
Since $Q$, $P_\ell$, and $R_\ell$ are positive,
$x_r>0$ is equivalent to the polynomial inequality
$N(\boldsymbol u)>0$.
After imposing this inequality, the last fraction has positive polynomial
numerator and denominator, so the same calculation applies to
the next unknown.

The $p$ bidiagonal factors have $p(d-1)$ subdiagonal entries.
Of these, $s=p(p-1)/2$ are freely chosen initial entries.
The triangular systems determine the remaining
\mbox{$M\coloneq p(d-1)-s$} entries.
Let $N_1,\dots,N_M$ be their polynomial numerators, with the
positive denominators constructed above. Consequently,
\[
 \mathcal D_L
 =\{\boldsymbol u\in(0,\infty)^s:
       N_j(\boldsymbol u)>0\text{ for }j\in\{1,\dots,M\}\}.
\]
Thus $\mathcal D_L$ is described by finitely many strict polynomial
inequalities and is open. Lemma~\ref{lem:rowwise} supplies positive
initial vectors in this set, so it is nonempty.

Every positive factorization of $L$ satisfies the same row equations
\eqref{eq:row-system}. Its positive pivots make back substitution
unique, so its initial entries recover the entire factorization.
Conversely, the construction matches every subdiagonal entry
of $L$ within its bandwidth; the diagonal and the entries outside
the band also agree. Hence the resulting product is $L$.
This proves the claimed bijection.

Write $\Phi(\boldsymbol u)$ for the constructed factors, and let
$\pi$ be the linear map that selects their $s$ initial entries.
The rational formulas have nonzero denominators on $\mathcal D_L$,
so $\Phi$ is analytic there. Since the chosen entries are preserved,
\(
 \pi\Phi(\boldsymbol u)=\boldsymbol u\), \(\pi\,\mathrm D\Phi(\boldsymbol u)=I_s\).
The derivative therefore has rank $s$.

Finally, \eqref{eq:row-system} with $r=1$ gives
$L_{n,n-1}=\sum_{i=1}^{p}a_i(n)$.
All $p\ge2$ summands are positive, so every initial entry satisfies
$0<a_i(n)<L_{n,n-1}$. Consequently
\(\mathcal D_L\subset
 \prod_{n=1}^{p-1}(0,L_{n,n-1})^{p-n}\),
which proves boundedness.
\end{proof}

\begin{example}[A three-parameter domain in dimension five]
\label{ex:five-dimensional-domain}
Take $p=3$, $d=5$, and
\[
 L\coloneq \begin{bNiceMatrix}
 1&0&0&0&0\\
 1&1&0&0&0\\
 1&2&1&0&0\\
 1&3&3&1&0\\
 0&1&3&3&1
 \end{bNiceMatrix}.
\]
Its Metelmann factors have the three initial zeros
$\widehat a_1(1)=\widehat a_1(2)=\widehat a_2(1)=0$;
all their other subdiagonal entries equal one.
For a positive factorization $L=L_1L_2L_3$, choose these initial
entries as $u\coloneq a_1(1)$, $v\coloneq a_2(1)$, and $w\coloneq a_1(2)$.
The row systems \eqref{eq:row-system} determine the remaining entries.
Row $1$ gives $a_3(1)=1-u-v$, and row $2$ gives
\[
 a_2(2)=\frac{1-w+uw}{1-u-v},\qquad
 a_3(2)=\frac{1-2u-2v+vw}{1-u-v}.
\]
Once the preceding rows are positive, the denominators in each
new row are positive, so the numerators give its polynomial
conditions. For instance,
\[
 a_2(3)=\frac{(1-u-v)(1-2w+3uw)}
 {(1-w+uw)(1-2u-2v+vw)}
\]
adds $1-2w+3uw>0$. Continuing through row $4$ shows that
$\mathcal D_L$ consists exactly of the triples $u,v,w>0$ satisfying
\begin{equation}\label{eq:five-dimensional-domain}
\begin{aligned}
 \text{row }1:\quad &1-u-v>0,\\
 \text{row }2:\quad &1-w+uw>0,
                    &&1-2u-2v+vw>0,\\
 \text{row }3:\quad &1-2w+3uw>0,
                    &&1-3u-3v+3vw>0,\\
 \text{row }4:\quad &1-3w+6uw>0,
                    &&1-4u-4v+6vw>0.
\end{aligned}
\end{equation}
For example, $u=v=w=1/10$ satisfies all seven inequalities.
\end{example}

For a full matrix $T=LDU$, the normalized triangular factors and $D$ are
unique, and the lower and upper refactorizations are independent.
For $p,q\ge2$, the domain is consequently
$\mathcal D_L\times\mathcal D_{U^{\mathsf T}}$, of dimension
$p(p-1)/2+q(q-1)/2$. A bandwidth-one side contributes a single
fixed factor.

The triangular systems establish the existence result in
Theorem~\ref{thm:finite}. Determining the factors one at a time
also gives linear recurrences for auxiliary sequences and rational
recurrences for the bidiagonal entries themselves. The data for
Proposition~\ref{prop:successive-recurrences} and
Corollary~\ref{cor:branched-factorization} are specified next.

To describe a single step, let $L^{[r]}$ be a unit lower triangular
matrix of order $d$ and lower bandwidth $r\in\{2,\dots,d-1\}$.
Write $g_j^{[r]}(n)\coloneq[L^{[r]}]_{n,n-j}$ for
$j\in\{0,\dots,\min(r,n)\}$, with $g_0^{[r]}(n)=1$,
and assume $g_r^{[r]}(n)\ne0$ for $n\in\{r,\dots,d-1\}$.
Choose nonzero values $\theta_h$, $h\in\{1,\dots,r-1\}$, for
the first $r-1$ subdiagonal entries of the factor to be determined.
Define the initial values of an auxiliary sequence by
\begin{equation}\label{eq:extraction-initial-data}
 \varphi_{r-1}^{[r]}\coloneq 1,\qquad
 \varphi_j^{[r]}\coloneq \prod_{h=j+1}^{r-1}\theta_h
 \quad(j\in\{0,\dots,r-2\}).
\end{equation}
Continue this sequence by the order-$r$ linear recurrence
\begin{equation}\label{eq:extraction-linear-recurrence}
 g_r^{[r]}(n)\varphi_n^{[r]}
 =\sum_{j=1}^{r}(-1)^{j-1}
       g_{r-j}^{[r]}(n-j)\varphi_{n-j}^{[r]},
 \qquad n\in\{r,\dots,d-1\}.
\end{equation}
Since $g_r^{[r]}(n)\ne0$, the linear recurrence
\eqref{eq:extraction-linear-recurrence} determines the entire
sequence from the chosen initial values. First compute this
sequence and check whether all its terms are nonzero.
For initial data satisfying this condition, put
$a_n^{[r]}\coloneq\varphi_{n-1}^{[r]}/\varphi_n^{[r]}$ for
$n\in\{1,\dots,d-1\}$,
and denote by $\mathcal L(a^{[r]})$ the unit lower bidiagonal
matrix with these subdiagonal entries.

\begin{proposition}[Determining successive factors by linear recurrences]
\label{prop:successive-recurrences}
With the preceding data, suppose
$\varphi_n^{[r]}\ne0$ for every $n\in\{0,\dots,d-1\}$.
Then
\begin{equation}\label{eq:extraction-product}
 L^{[r]}=\mathcal L(a^{[r]})L^{[r-1]},\qquad
 L^{[r-1]}\coloneq \mathcal L(a^{[r]})^{-1}L^{[r]},
\end{equation}
where $L^{[r-1]}$ has lower bandwidth at most $r-1$.
Conversely, every factorization of this form with nonzero bidiagonal
entries is obtained in this way, with $\theta_h=a_h^{[r]}$.
For a positive bidiagonal factorization of a lower product of
bandwidth $p$, determining one factor at a time therefore gives $p-1$
recurrences, of orders $p,p-1,\ldots,2$, carrying respectively
$p-1,p-2,\ldots,1$ initial parameters. The final bandwidth-one
matrix is the last factor.
\end{proposition}
\begin{proof}
For a single extraction, omit the superscript $[r]$ and write
$L_{\mathrm{rem}}\coloneq \mathcal L(a)^{-1}L$. The unit triangular inverse
exists in finite dimension. Comparing entry $(n,k)$ in
$L=\mathcal L(a)L_{\mathrm{rem}}$ gives
$[L_{\mathrm{rem}}]_{n,k}=L_{n,k}-a_n[L_{\mathrm{rem}}]_{n-1,k}$ for $k<n$.
Fix $n\ge r$ and put $k\coloneq n-r$. Iterate this equation through
rows $n,n-1,\dots,n-r+1$. At the last step the remaining entry
is $[L_{\mathrm{rem}}]_{k,k}=1=g_0(k)$, because both factors
are unit lower triangular. Therefore
\[
 [L_{\mathrm{rem}}]_{n,n-r}
 =\sum_{j=0}^{r}(-1)^j
 \left(\prod_{h=0}^{j-1}a_{n-h}\right)g_{r-j}(n-j).
\]
The product for $j=0$ equals one. All terms of $\varphi$ are
nonzero by hypothesis, so the substitutions
$a_t=\varphi_{t-1}/\varphi_t$ are legitimate and give
$\prod_{h=0}^{j-1}a_{n-h}=\varphi_{n-j}/\varphi_n$.
Multiplying the preceding sum by $\varphi_n$ shows that its
vanishing is precisely \eqref{eq:extraction-linear-recurrence},
with the term $j=0$ moved to the other side. Thus the entire
$r$th subdiagonal of $L_{\mathrm{rem}}$ vanishes.

To check entries farther below the diagonal, fix a column $k$
with $k+r\le d-1$. The preceding calculation gives
$[L_{\mathrm{rem}}]_{k+r,k}=0$. For $n>k+r$, the bandwidth of
$L$ gives $L_{n,k}=0$, so
$[L_{\mathrm{rem}}]_{n,k}=-a_n[L_{\mathrm{rem}}]_{n-1,k}$.
Induction on $n$ makes every such entry zero. Columns with
$k+r>d-1$ have no entries to check this far below the diagonal.
This proves the claimed bandwidth reduction.

The initial values satisfy
$\varphi_{h-1}/\varphi_h=\theta_h$ for $1\le h\le r-1$.
Conversely, start with a factorization of the stated form whose
bidiagonal entries $a_n$ are all nonzero. Normalize
$\varphi_{r-1}\coloneq 1$, define the earlier values by
$\varphi_{h-1}\coloneq a_h\varphi_h$, and the later ones by
$\varphi_n\coloneq \varphi_{n-1}/a_n$. This gives a nonzero sequence and
recovers \eqref{eq:extraction-initial-data} with $\theta_h=a_h$.
The bandwidth of the remaining factor makes
$[L_{\mathrm{rem}}]_{n,n-r}=0$, so the identity just proved
forces \eqref{eq:extraction-linear-recurrence}. Its coefficient
$g_r(n)$ is nonzero, hence these initial values determine every
later term uniquely. This proves the converse.

If the given matrix is a product of $r$ positive bidiagonal factors,
removing its first factor leaves a product of $r-1$ positive
bidiagonal factors. On its extreme subdiagonal there is exactly
one contributing term: each factor supplies its subdiagonal entry.
That term is positive, so the next recurrence has the required
nonzero leading coefficients. Repeat until the remaining matrix
has bandwidth one and is itself the last factor.
At order $r$, an order-$r$ recurrence needs $r$ initial values.
Multiplying all of them by one common nonzero scalar leaves every
ratio $\varphi_{n-1}/\varphi_n$ unchanged; fixing
$\varphi_{r-1}=1$ removes exactly this freedom. The remaining
$r-1$ initial ratios are the chosen entries $\theta_1,\dots,\theta_{r-1}$.
Summing over $r\in\{p,p-1,\dots,2\}$ gives $p(p-1)/2$ parameters.
At each later extraction, the recurrence coefficients are the
entries of the remaining product already determined.
\end{proof}

With the normalization $\varphi_{r-1}^{[r]}=1$ in
\eqref{eq:extraction-initial-data}, positivity of the extracted
factor is equivalent to $\varphi_n^{[r]}>0$ for every
$n\in\{0,\dots,d-1\}$.
In the finite problem, Lemma~\ref{lem:rowwise} supplies positive
factorizations. Taking the initial entries of their successive
factors gives positive sequences by the converse just proved.
Proposition~\ref{prop:positive-parameter-domain} describes all
initial vectors that produce positive factors at every stage.

In the semi-infinite BTP problem, positivity must hold throughout
the infinite sequence at each extraction.
Lemma~\ref{lem:global-extraction-determinant} expresses each later
term as a determinant divided by a product of positive coefficients.
Theorem~\ref{thm:global-lower-criterion} then shows that convergence
of the adjacent Metelmann series provides initial parameters making
these determinants positive in every row and at every stage.

Taking ratios in \eqref{eq:extraction-linear-recurrence} gives
a recurrence for the bidiagonal entries themselves. Repeated
substitution expresses each entry in terms of the prescribed initial
values. For $r=2$ this gives an ordinary finite continued fraction;
for higher orders, products of earlier entries produce the branching.

For the subdiagonal entries
$a_n^{[r]}=\varphi_{n-1}^{[r]}/\varphi_n^{[r]}$ of the factor
$\mathcal L(a^{[r]})$ in \eqref{eq:extraction-product}, write
\[
 \gamma_{n,j}^{[r]}
 \coloneq\frac{(-1)^{j-1}g_{r-j}^{[r]}(n-j)}{g_r^{[r]}(n)},
 \qquad x_n\coloneq a_n^{[r]},
\]
where $\gamma_{n,j}^{[r]}$ is defined for
$n\in\{r,\dots,d-1\}$ and $j\in\{1,\dots,r\}$, and $x_n$
for $n\in\{1,\dots,d-1\}$. To express a fixed entry $x_N$, with $N\in\{r,\dots,d-1\}$,
number the entries backwards by setting $F_k\coloneq x_{N-k}$ for
$k\in\{0,\dots,N-1\}$.
For $k\in\{0,\dots,N-r\}$ and $\ell\in\{1,\dots,r-1\}$, set
$c_{k,0}\coloneq\gamma_{N-k,1}^{[r]}$ and
$c_{k,\ell}\coloneq-\gamma_{N-k,\ell+1}^{[r]}$.

\begin{corollary}[Finite branched continued fractions]
\label{cor:branched-factorization}
With the preceding data and notation, suppose that the sequence
defined by \eqref{eq:extraction-initial-data} and
\eqref{eq:extraction-linear-recurrence} satisfies
$\varphi_n^{[r]}\ne0$ for every $n\in\{0,\dots,d-1\}$.
Then the bidiagonal entries satisfy
\begin{equation}\label{eq:extraction-ratio-recurrence}
 x_n=\frac{1}{
 \gamma_{n,1}^{[r]}+
 \sum_{j=2}^{r}\gamma_{n,j}^{[r]}
                   \prod_{h=1}^{j-1}x_{n-h}},
 \qquad x_h=\theta_h\quad(h\in\{1,\dots,r-1\}).
\end{equation}
For each $N\in\{r,\dots,d-1\}$, the backward notation gives
\begin{equation}\label{eq:extraction-branched-fraction}
 F_k=\frac{1}{ c_{k,0}
                   -\sum_{\ell=1}^{r-1}c_{k,\ell}
                                 \prod_{h=1}^{\ell}F_{k+h}},
 \quad k\in\{0,\dots,N-r\},
\end{equation}
with terminal values $F_{N-h}=\theta_h$ for $h\in\{1,\dots,r-1\}$.
Repeated substitution therefore gives a finite branched continued
fraction for $x_N=F_0$, ending at the prescribed initial values.
This is a finite $(r-1)$-branched Jacobi-type fraction:
its denominator contains products of up to $r-1$ later fractions.
The coefficients may have either sign.
\end{corollary}
\begin{proof}
The numbers $g_r^{[r]}(n)$ and $\varphi_{n-1}^{[r]}$ are
nonzero, so divide
\eqref{eq:extraction-linear-recurrence} by their product.
For $j\ge2$, the adjacent ratios telescope:
\[
 \frac{\varphi_{n-j}^{[r]}}{\varphi_{n-1}^{[r]}}
 =\frac{\varphi_{n-j}^{[r]}}{\varphi_{n-j+1}^{[r]}}
   \cdots\frac{\varphi_{n-2}^{[r]}}{\varphi_{n-1}^{[r]}}
 =\prod_{h=1}^{j-1}x_{n-h}.
\]
The divided equation is consequently
\[
 \frac1{x_n}=\gamma_{n,1}^{[r]}+
 \sum_{j=2}^{r}\gamma_{n,j}^{[r]}\prod_{h=1}^{j-1}x_{n-h}.
\]
Its right-hand side equals the nonzero number
$\varphi_n^{[r]}/\varphi_{n-1}^{[r]}$, so taking its reciprocal
proves \eqref{eq:extraction-ratio-recurrence}.

Now put $n=N-k$ and $j=\ell+1$. Then
$x_{n-h}=F_{k+h}$ and
$\gamma_{N-k,\ell+1}^{[r]}=-c_{k,\ell}$.
These substitutions give \eqref{eq:extraction-branched-fraction}.
The indices $n\in\{1,\dots,r-1\}$ correspond exactly to
$F_{N-h}=\theta_h$. Every substitution of a nonterminal fraction
increases $k$ by at least one; the indices never exceed $N-1$.
Thus every branch of the substitution ends after finitely many
steps at one of those initial values. Each denominator equals
$1/x_{N-k}$, so the nonvanishing hypothesis also justifies all
intermediate divisions.
\end{proof}

For comparison with the usual generating-function convention for
branched Jacobi fractions \cite[Section~4]{PetreolleSokalZhu},
write $\beta_k^{(0)}\coloneq 1-c_{k,0}$ and
$\beta_{k+\ell}^{(\ell)}\coloneq c_{k,\ell}$ for
$\ell\in\{1,\dots,r-1\}$. The denominator in that convention is
\[
 1-\beta_k^{(0)}t-
 \sum_{\ell=1}^{r-1}\beta_{k+\ell}^{(\ell)}t^{\ell+1}
                   \prod_{h=1}^{\ell}F_{k+h}(t).
\]
At $t=1$ it gives \eqref{eq:extraction-branched-fraction}, with
the initial factor entries as terminal values.
For positive bidiagonal products, $c_{k,\ell}$ has sign
$(-1)^{\ell+1}$ for $\ell\in\{1,\dots,r-1\}$; the positive
values furnished by these finite fractions are the factor entries.

For $p=3$, the successive factors are calculated using a two-branched
fraction with two initial values and an ordinary continued fraction
with one initial value. For $p=4$, the successive orders $4,3,2$
require $3+2+1=6$ initial parameters. At $p=2$ this reduces to the
ordinary continued fraction of the tetradiagonal analysis
in \cite{BFMTetra}.

Positivity of the factors restricts the admissible initial entries.
The next proposition gives upper bounds for $[L_1]_{1,0}$,
the first subdiagonal entry of the leftmost factor, using only
the inverse of the fixed product $L$. These bounds decrease
strictly as the row index increases. In the semi-infinite setting,
a global positive factorization must satisfy all the bounds
with the same initial entry. If they tend to zero, no positive
choice is possible. This observation will provide an obstruction
to global factorization.

In dimension $d$, let
$\Sigma\coloneq\diag[(-1)^n]_{n\in\{0,\dots,d-1\}}$ be the diagonal
sign matrix. Multiplying $L^{-1}$ by $\Sigma$ on both sides removes
its alternating signs, so $\mathcal Q\coloneq \Sigma L^{-1}\Sigma$ has
nonnegative entries.
For each unit lower bidiagonal factor, the entries of its inverse,
after the same sign change, are products of consecutive subdiagonal
entries, as in \eqref{eq:signed-bidiagonal-inverse-entries} below.
These inverse formulas and the use of the diagonal sign matrix
are classical; see \cite[Lemma~2.1 and Theorem~3.2]{Higham2024}.
Here they bound initial factor entries; Lemma~\ref{lem:global-inverse-minors}
uses them to calculate inverse minors and bound their normalized values
uniformly when the adjacent series converge.

\begin{proposition}[Strictly decreasing inverse-entry bounds]
\label{prop:inverse-parameter-bound}
Let $p\in\{2,\dots,d-1\}$ and
$L=\widehat L_1\cdots\widehat L_p$, where the lower Metelmann
factors are as in \eqref{eq:metelmann}, and define
$\mathcal Q_{n,j}\coloneq (-1)^{n-j}[L^{-1}]_{n,j}$ and
$\kappa_n\coloneq \mathcal Q_{n,0}/\mathcal Q_{n,1}$ for
$n\in\{1,\dots,d-1\}$.
Both entries in the ratio are positive, and
\begin{equation}\label{eq:inverse-bound-monotonicity}
 \kappa_n>\kappa_{n+1}>0
 \qquad(n\in\{1,\dots,d-2\}).
\end{equation}
Every positive factorization $L=L_1\cdots L_p$ satisfies
\begin{equation}\label{eq:inverse-bound-finite}
 0<[L_1]_{1,0}<\kappa_{d-1}.
\end{equation}
\end{proposition}
\begin{proof}
	For each $i\in\{1,\dots,p\}$, write $\widehat L_i=I+A_i$ and
put $\ell_{i,h}\coloneq [A_i]_{h,h-1}\ge0$ for $h\in\{1,\dots,d-1\}$,
	where all other entries of $A_i$ vanish.
	Since $\Sigma^2=I$ and each subdiagonal entry changes sign under
multiplication by $\Sigma$ on both sides, $\Sigma\widehat L_i\Sigma=I-A_i$.
	Taking inverses therefore gives
$\mathcal Q_i\coloneq \Sigma\widehat L_i^{-1}\Sigma=(I-A_i)^{-1}$.
	Moreover, $L=\widehat L_1\cdots\widehat L_p$, so reversing the
order of the factors when taking the inverse yields
$\mathcal Q=\Sigma L^{-1}\Sigma=\mathcal Q_p\cdots\mathcal Q_1$.
	
	Because $A_i$ is strictly lower triangular, $A_i^d=0$. Thus
$\mathcal Q_i=I+A_i+\cdots+A_i^{d-1}$.
	For $n>j$, only the power $A_i^{n-j}$ contributes to entry
	$(n,j)$, and that entry is the product of the consecutive
	subdiagonal entries between columns $j$ and $n$. Consequently,
	\begin{equation}\label{eq:signed-bidiagonal-inverse-entries}
	[\mathcal Q_i]_{n,j}=
	\begin{cases}
		0, & n<j,\\
		1, & n=j,\\
		\displaystyle\prod_{h=j+1}^{n}\ell_{i,h}, & n>j.
	\end{cases}
	\end{equation}
	In particular, $\mathcal Q_i$ is entrywise nonnegative and
	unit lower triangular.
	
	Nonnegativity of its minors follows from an elementary factorization. Define
$E_{i,h}\coloneq I+\ell_{i,h}e_he_{h-1}^{\mathsf T}$.
	Formula~\eqref{eq:signed-bidiagonal-inverse-entries} gives
$\mathcal Q_i=E_{i,d-1}\cdots E_{i,1}$.
	Each $E_{i,h}$ is TN: its minors are either $0$, $1$, or
	$\ell_{i,h}$. Products of TN matrices are TN by
	Cauchy--Binet. Hence every $\mathcal Q_i$, and every product
	of these matrices, is TN.
	
	The initial zero pattern in \eqref{eq:metelmann-parameters} determines
	the first two columns.
	The last factor satisfies $\ell_{p,h}>0$ for
$h\in\{1,\dots,d-1\}$,
	so \eqref{eq:signed-bidiagonal-inverse-entries} implies
$[\mathcal Q_p]_{n,j}>0$ whenever $n\ge j$.
	For every $i\in\{1,\dots,p-1\}$, however, $\ell_{i,1}=0$.
	Every entry of column zero of $\mathcal Q_i$ below the diagonal
	therefore vanishes. Equivalently, $\mathcal Q_i e_0=e_0$.
	
	Set $R\coloneq \mathcal Q_{p-1}\cdots\mathcal Q_1$.
	It follows that $Re_0=e_0$, or explicitly, $R_{0,0}=1$ and
$R_{n,0}=0$ for $n\in\{1,\dots,d-1\}$.
	To examine column one, write $R=\mathcal Q_{p-1}S$, where
	$S\coloneq \mathcal Q_{p-2}\cdots\mathcal Q_1$ if $p>2$ and $S\coloneq I$
	if $p=2$. The matrix $S$ is entrywise nonnegative and unit
	lower triangular. Since $\ell_{p-1,h}>0$ for $h\in\{2,\dots,d-1\}$,
	formula~\eqref{eq:signed-bidiagonal-inverse-entries} gives
	$[\mathcal Q_{p-1}]_{n,1}>0$ for $n\ge1$.
	The term with intermediate index $1$ in the matrix product
then gives $R_{n,1}\ge[\mathcal Q_{p-1}]_{n,1}S_{1,1}
=[\mathcal Q_{p-1}]_{n,1}>0$.
	
	Since $\mathcal Q=\mathcal Q_pR$ and $Re_0=e_0$,
$\mathcal Q_{n,0}=[\mathcal Q_p]_{n,0}>0$.
	Similarly, using $R_{1,1}=1$ and nonnegativity of all summands,
the product gives $\mathcal Q_{n,1}\ge[\mathcal Q_p]_{n,1}R_{1,1}
=[\mathcal Q_p]_{n,1}>0$ for $n\in\{1,\dots,d-1\}$.
	Thus all the ratios $\kappa_n$ are well defined and positive.
	
	Fix $n\in\{1,\dots,d-2\}$. Apply Cauchy--Binet to the minor
	of $\mathcal Q_pR$ with rows $\{n,n+1\}$ and columns $\{0,1\}$.
	Every summand is nonnegative because both matrices are TN.
	For the intermediate indices $\{0,n+1\}$, the two determinants
	are
	\[
	\det\begin{bNiceMatrix}
		[\mathcal Q_p]_{n,0}&0\\
		[\mathcal Q_p]_{n+1,0}&1
	\end{bNiceMatrix}
	=[\mathcal Q_p]_{n,0}>0,
	\qquad
	\det\begin{bNiceMatrix}
		1&0\\
		0&R_{n+1,1}
	\end{bNiceMatrix}
	=R_{n+1,1}>0.
	\]
	Hence
	\(\mathcal Q_{n,0}\mathcal Q_{n+1,1}
	-\mathcal Q_{n,1}\mathcal Q_{n+1,0}>0\).
	Dividing by the positive product
	$\mathcal Q_{n,1}\mathcal Q_{n+1,1}$ gives
	\[
	\kappa_n-\kappa_{n+1}
	=
	\frac{\mathcal Q_{n,0}\mathcal Q_{n+1,1}
		-\mathcal Q_{n,1}\mathcal Q_{n+1,0}}
	{\mathcal Q_{n,1}\mathcal Q_{n+1,1}}>0.
	\]
	Together with positivity of the ratios, this proves
	\eqref{eq:inverse-bound-monotonicity}.
	
	Finally, consider any positive factorization $L=L_1\cdots L_p$,
and put $L_{\mathrm{rem}}\coloneq L_2\cdots L_p$ and $u\coloneq [L_1]_{1,0}>0$.
	The identity $L=L_1L_{\mathrm{rem}}$ gives
$L_{\mathrm{rem}}^{-1}=L^{-1}L_1$ and
$\Sigma L_{\mathrm{rem}}^{-1}\Sigma=\mathcal Q(\Sigma L_1\Sigma)$.
	Column zero of $\Sigma L_1\Sigma$ is $e_0-ue_1$. Taking entry $(n,0)$
therefore yields $[\Sigma L_{\mathrm{rem}}^{-1}\Sigma]_{n,0}
=\mathcal Q_{n,0}-u\mathcal Q_{n,1}$.
	
	There is at least one factor in $L_{\mathrm{rem}}$, since $p\ge2$.
	For each of its factors $L_i$, formula
	\eqref{eq:signed-bidiagonal-inverse-entries} shows that $\Sigma L_i^{-1}\Sigma$ has strictly positive
	entries on and below the diagonal. Their product
$\Sigma L_{\mathrm{rem}}^{-1}\Sigma=(\Sigma L_p^{-1}\Sigma)\cdots(\Sigma L_2^{-1}\Sigma)$
	also has strictly positive entries on and below the diagonal.
	Thus $0<\mathcal Q_{n,0}-u\mathcal Q_{n,1}$, and hence
$0<u<\kappa_n$, for $n\in\{1,\dots,d-1\}$.
	Taking $n=d-1$ proves \eqref{eq:inverse-bound-finite}.
\end{proof}

\subsection{Semi-infinite factorizations and their obstructions}\label{sec:limit}

The initial parameters of a global positive factorization must belong
to the domains of Proposition~\ref{prop:positive-parameter-domain}
for all finite leading sections. For the following bounded tetradiagonal
matrix, these domains are explicit intervals with empty intersection.
Compatibility of the Metelmann factors across matrix sizes provides
the starting point for the complete criterion in
Section~\ref{subsec:global-criterion} and the norm approximation in
Section~\ref{subsec:boundary-approximation}.

Let $\Lambda$ be the upper shift on $\ell^2(\NN)$, defined by
$\Lambda e_0\coloneq 0$ and $\Lambda e_n\coloneq e_{n-1}$ for $n\in\mathbb N$.
Its transpose is the lower shift, $\Lambda^{\mathsf T}e_n=e_{n+1}$.
For a bounded operator $A$ on $\ell^2(\NN)$, the induced norm is
$\|A\|_2\coloneq \sup_{\|u\|_{\ell^2}=1}\|Au\|_{\ell^2}$, where
$\|u\|_{\ell^2}^2\coloneq \sum_{n\in\NN}|u_n|^2$.
For a real sequence $a=(a_n)_{n\in\mathbb N}$ write $\mathcal L(a)$ for
the unit lower bidiagonal matrix with subdiagonal $a_n$. It defines a
bounded operator on $\ell^2$ exactly when the sequence is bounded.
Consider
\begingroup
\setlength{\arraycolsep}{3pt}
\begin{equation}\label{eq:counter-matrix}
\begin{aligned}
T_*&\coloneq \begin{bNiceMatrix}[parallelize-diags=false, margin=5pt,small]
 1&1&0&0&0&0&\Cdots[shorten-end=-10pt]&\phantom{0}\\
 1&2&1&0&0&0&\Cdots[shorten-end=-10pt]&\phantom{0}\\
 1&3&3&1&0&0&\Cdots[shorten-end=-10pt]&\phantom{0}\\
 0&1&3&3&1&0&\Cdots[shorten-end=-10pt]&\phantom{0}\\
 0&0&1&3&3&1&0&\phantom{0}\\
 \Vdots[shorten-end=-7pt]&\Vdots[shorten-end=-7pt]&\Ddots[shorten-end=-10pt]&\Ddots[shorten-end=-10pt]&\Ddots[shorten-end=-10pt]&\Ddots[shorten-end=-10pt]&\Ddots[shorten-end=-10pt]&\\
 \phantom{0}&\phantom{0}&\phantom{0}&\phantom{0}&
 \phantom{0}&\phantom{0}&\phantom{0}&\phantom{0}
 \end{bNiceMatrix}=\mathcal L(0,1,1,\ldots)\mathcal L(1,1,1,\ldots)(I+\Lambda).
\end{aligned}
\end{equation}
\endgroup
The matrix $T_*$ is bounded and has bandwidth $(2,1)$.
Its first three diagonal entries are $1,2,3$.
Deleting its first two rows and columns gives a Toeplitz matrix:
\[
 [[T_*]_{i+2,j+2}]_{i,j\in\NN}=[t_{i-j}]_{i,j\in\NN},
 \qquad (t_{-1},t_0,t_1,t_2)\coloneq (1,3,3,1),
\]
with $t_k\coloneq 0$ for $k\notin\{-1,\dots,2\}$.
This property defines the \emph{eventually Toeplitz} structure used here.
Its normalized lower--upper decomposition is $T_*=L_*U_*$
with $U_*=I+\Lambda$, and $L_*$ has diagonal $1$, first subdiagonal
$m_1=1$, $m_n=2$ for $n\in\mathbb N\setminus\{1\}$, and second subdiagonal $1$.

\begin{example}[Finite-section PBFs and a global obstruction]
\label{ex:counterexample}
In a factorization $L_*=\mathcal L(a)\mathcal L(b)$, set $a_1=t>0$.
Matching the two subdiagonals gives $b_1=1-t$ and, for $n\ge2$,
$a_n=1/b_{n-1}$ and $b_n=2-1/b_{n-1}$. Hence
$b_n=(1-nt)/(1-(n-1)t)$.
Starting with $0<t<1$, positivity through row $N$ is therefore
equivalent to $0<t<1/N$. Every leading section of order $N+1$
thus has a PBF and is BTP, with bandwidth $(\min(2,N),1)$ for
$N\ge1$; the section of order one is $[1]$.
Uniqueness of the normalized LDU decomposition forces the lower
product of a global PBF to be $L_*$. Its first entry would have
to satisfy $0<t<1/N$ for every $N$, which is impossible.
The bandwidth here is $(2,1)$, within the tetradiagonal setting of
\cite{BFMTetra}; its continued-fraction criterion also decides this
example. The global criterion and the initial-block approximation below
treat arbitrary lower and upper bandwidths.
\end{example}

The admissible intervals $(0,1/N)$ have empty intersection, and their
common limiting endpoint $t=0$ recovers the Metelmann factorization
with $a_1=0$. This identifies the parameter responsible for the global
obstruction.

The finite characterization in \cite[Theorem~2.16(iv)]{BFM2024},
recalled in Proposition~\ref{prop:metelmann-input}, extends to
semi-infinite matrices by taking leading sections of order
$d>\max(p,q)$. The ordered product \eqref{eq:metelmann} respects
truncation. Its factors are unique: the normalized LDU decomposition
is unique, and the prescribed initial zeros determine the individual
bidiagonal factors through the triangular systems
\eqref{eq:row-system}, as in the proof of Lemma~\ref{lem:rowwise}.
Thus the factorizations of nested leading sections agree on their
common entries and define one semi-infinite Metelmann factorization.
Conversely, every sufficiently large leading section of such a
product is BTP by the same finite characterization.

\subsection{A complete criterion for global positive factorization}
\label{subsec:global-criterion}

For a semi-infinite BTP matrix $T$ of bandwidth $(p,q)$, the
global criterion starts from its normalized decomposition $T=LDU$
and the Metelmann products for $L$ and $U^{\mathsf T}$ in
\eqref{eq:metelmann}, whose compatibility was established in
Section~\ref{sec:limit}.
These factors determine the $p+q-2$ nonnegative numbers
\[
 \underbrace{\kappa^-_{1,\infty},\dots,\kappa^-_{p-1,\infty}}
 _{\text{from }L},
 \qquad
 \underbrace{\kappa^+_{1,\infty},\dots,\kappa^+_{q-1,\infty}}
 _{\text{from }U^{\mathsf T}}.
\]
A side of bandwidth one contributes no numbers.
Corollary~\ref{cor:global-two-sided-criterion} will show that $T$
admits a global PBF exactly when all these numbers are strictly
positive. Each is the limit of a decreasing sequence of ratios of
inverse minors, the value of a scalar Jacobi continued fraction,
and the reciprocal of a positive series determined by two adjacent
Metelmann factors. The computation of these factors from $T$ is
described first.

Compute $T=LDU$ by Gaussian elimination without row exchanges
on finite leading sections. All leading principal minors of $T$
are positive, so the pivots are positive. Uniqueness of the normalized LDU decomposition makes the
triangular factors agree on overlapping leading sections.

For the lower product $L=\widehat L_1\cdots\widehat L_p$, write
$\widehat a_i(n)\coloneq[\widehat L_i]_{n,n-1}$. Proceed in
increasing row order $n\ge1$. With $k\coloneq\min(p,n)$, set
$\widehat a_i(n)=0$ for $i\in\{1,\dots,p-k\}$, and determine
the remaining $k$ entries from \eqref{eq:row-system}, using the
Metelmann factors in place of $L_i$. Solve the equations for
$r\in\{1,\dots,k\}$ in decreasing order. Their coefficients use
only earlier rows, and their divisors are the positive products
in \eqref{eq:pivot}. Applying the same procedure to $U^{\mathsf T}$
and transposing the resulting product gives the upper Metelmann factors.
Thus the factors are computed uniquely, row by row, from $T$.

\begin{example}[Computing three Metelmann factors]
\label{ex:three-metelmann-factors}
For example, when $p=3$, denote the subdiagonal entries in row $n$
of $\widehat L_1,\widehat L_2,\widehat L_3$ by $a(n),b(n),c(n)$,
respectively. The initial zeros and the first two rows give
\[
 \begin{aligned}
 a(1)=a(2)=b(1)&=0, & c(1)&=L_{1,0},&
 b(2)&=\frac{L_{2,0}}{c(1)}, & c(2)&=L_{2,1}-b(2).
 \end{aligned}
\]
For $n\ge3$, compute in the displayed order:
\[
 \begin{aligned}
 a(n)&=\frac{L_{n,n-3}}{b(n-1)c(n-2)},\\[1mm]
 b(n)&=\frac{L_{n,n-2}-a(n)\bigl(b(n-1)+c(n-1)\bigr)}{c(n-1)},\\[1mm]
 c(n)&=L_{n,n-1}-a(n)-b(n).
 \end{aligned}
\]
These formulas follow by matching the third, second, and first
subdiagonals of $L=\widehat L_1\widehat L_2\widehat L_3$, in that
order. Proposition~\ref{prop:metelmann-input} gives positivity of every
permitted entry on each leading section; the triangular systems
\eqref{eq:row-system} determine the factor entries uniquely.
\end{example}

The factors determine the following series and inverse-minor ratios
for the lower side. The same definitions
applied to $U^{\mathsf T}$ will give the upper conditions.

Fix $p\in\mathbb N$. Throughout this subsection, let
$L=\widehat L_1\cdots\widehat L_p$ be a compatible unit lower
Metelmann product, with $\widehat L_i=\mathcal L(\widehat a_i)$,
$\widehat a_i(n)=0$ for $n\le p-i$, and $\widehat a_i(n)>0$ for
$n\ge p-i+1$. Put
$\Sigma\coloneq \diag[(-1)^n]_{n\in\NN}$ and $\mathcal Q\coloneq \Sigma L^{-1}\Sigma$.
The inverse is defined entrywise by the triangular inverses of the
finite leading sections. Starting with the rightmost factor, define
\begin{equation}\label{eq:global-cumulative-products}
	\pi_j(n)\coloneq
	\prod_{t=j+1}^{n}\widehat a_{p-j}(t),
	\qquad n\ge j,\quad j\in\{0,\dots,p-1\}.
\end{equation}
The empty product equals $1$, so $\pi_j(j)=1$.
For $k\in\{1,\dots,p-1\}$ define
\begin{equation}\label{eq:global-adjacent-series}
 r_k(n)\coloneq\frac{\pi_k(n)}{\pi_{k-1}(n)}\quad(n\ge k),\qquad
 \mathfrak S_k\coloneq\sum_{n=k}^{\infty}r_k(n)\in(0,\infty].
\end{equation}
Each $r_k(n)$ compares the accumulated products of the two consecutive
Metelmann factors $\widehat L_{p-k}$ and $\widehat L_{p-k+1}$.
In terms of their subdiagonal entries, \eqref{eq:global-cumulative-products}
gives
\[
 r_k(n)=\frac{1}{\widehat a_{p-k+1}(k)}
 \prod_{t=k+1}^{n}
 \frac{\widehat a_{p-k}(t)}{\widehat a_{p-k+1}(t)},\qquad n\ge k.
\]
The prefactor accounts for the extra entry at $t=k$ in the denominator
product. The series $\mathfrak S_k$ sums these positive ratios as $n$
increases.

For $k\ge1$ and $n\ge0$, define
$\mathcal I_{k,n}\coloneq\{n,\dots,n+k-1\}$,
$\mathcal C_k^{(0)}\coloneq\mathcal I_{k,0}$, and
$\mathcal C_k^{(1)}\coloneq(\mathcal C_k^{(0)}\setminus\{k-1\})\cup\{k\}$.
Thus $\mathcal C_k^{(0)}$ consists of the first $k$ column indices, while
$\mathcal C_k^{(1)}$ is obtained by replacing the last index $k-1$ with $k$.
In particular, $\mathcal C_1^{(0)}=\{0\}$ and $\mathcal C_1^{(1)}=\{1\}$.
If $I,K$ are increasing finite index sets,
$\mathcal Q[I,K]$ denotes the corresponding submatrix. Define the
\emph{initial contiguous minors} $\Delta_{k,n}$ for
$k\in\{1,\dots,p\}$ and $n\ge0$, and
the ratios $\kappa_{k,n}$ for $k\in\{1,\dots,p-1\}$ and $n\ge1$, by
\begin{equation}\label{eq:global-minor-ratios}
 \Delta_{k,n}\coloneq\det\mathcal Q[\mathcal I_{k,n},\mathcal C_k^{(0)}],\qquad
 \kappa_{k,n}\coloneq\frac{\Delta_{k,n}}
 {\det\mathcal Q[\mathcal I_{k,n},\mathcal C_k^{(1)}]}.
\end{equation}
Lemma~\ref{lem:global-inverse-minors} shows that the denominator is
positive for $n\ge1$.
The minors $\Delta_{k,n}$ use consecutive rows and the first $k$
columns. Such contiguous minors are called initial in
\cite[p.~75]{FallatJohnson2011}.

For $k=1$, this is the ratio $\kappa_n$ in
Proposition~\ref{prop:inverse-parameter-bound}.

The next lemma computes the inverse-minor ratios from the adjacent
series and bounds the normalized minors uniformly in the row index.
This bound is what allows one choice of positive initial entries to
work in every row in Theorem~\ref{thm:global-lower-criterion}.

When all the adjacent series converge, define
\begin{equation}\label{eq:global-minor-bound-constant}
 \mathfrak B\coloneq
 \max\left(\{1\}\cup
 \left\{\prod_{h=\ell+1}^{j}\mathfrak S_h:
                 0\le\ell<j\le p-1\right\}\right).
\end{equation}
This is the maximum of $1$ and all products of consecutive terms
among $\mathfrak S_1,\dots,\mathfrak S_{p-1}$; in particular,
$\mathfrak B=1$ when $p=1$.

\begin{lemma}[Inverse minors and a uniform bound]
\label{lem:global-inverse-minors}
The minors and ratios defined in
\eqref{eq:global-minor-ratios} satisfy:
\begin{enumerate}
\item For $k\in\{1,\dots,p\}$ and $n\ge0$,
\begin{equation}\label{eq:global-leading-minor}
 \Delta_{k,n}=\prod_{j=0}^{k-1}\pi_j(n+j)>0.
\end{equation}
For $k\in\{1,\dots,p-1\}$ and $n\ge0$,
\begin{equation}\label{eq:global-ratio-series}
 \det\mathcal Q[\mathcal I_{k,n},\mathcal C_k^{(1)}]
 =\Delta_{k,n}\sum_{s=k}^{n+k-1}r_k(s).
\end{equation}
The determinant in \eqref{eq:global-ratio-series} is strictly
positive for $n\ge1$.

\item For each $k\in\{1,\dots,p-1\}$, the sequence
$(\kappa_{k,n})_{n\ge1}$ is positive and strictly decreasing, with
\[
 \kappa_{k,n}^{-1}=\sum_{s=k}^{n+k-1}r_k(s),
 \qquad
 \kappa_{k,\infty}\coloneq\lim_{n\to\infty}\kappa_{k,n}
 =\frac1{\mathfrak S_k},
\]
where $1/\infty=0$.

\item If all the adjacent series converge, the constant
$\mathfrak B$ in \eqref{eq:global-minor-bound-constant} gives the
uniform bound
\begin{equation}\label{eq:global-uniform-minors}
 \left|
 \frac{\det\mathcal Q[\mathcal I_{k,n},K]}{\Delta_{k,n}}
 \right|
 \le k!\mathfrak B^k
\end{equation}
for every $k\in\{1,\dots,p\}$, $n\ge0$, and $k$-element subset
$K\subseteq\{0,\dots,p-1\}$.
\end{enumerate}
\end{lemma}
\begin{proof}
\emph{(1) Identities for the inverse minors.}
All inverse entries used below are entries of triangular inverses
of finite leading sections. If the largest selected row is $N$,
all intermediate indices in a lower triangular product lie in
$\{0,\dots,N\}$. Thus the identities can be proved on that
finite section and are unchanged when a larger section is taken.
As in the proof of Proposition~\ref{prop:inverse-parameter-bound}, put
$\mathcal Q_i\coloneq \Sigma\widehat L_i^{-1}\Sigma$ for $i\in\{1,\dots,p\}$,
so that $\mathcal Q=\mathcal Q_p\cdots\mathcal Q_1$.
	Recall the cumulative products $\pi_j$ from
    \eqref{eq:global-cumulative-products}. For $j\in\{0,\dots,p-1\}$,
    formula~\eqref{eq:signed-bidiagonal-inverse-entries}, applied to
	$\widehat L_{p-j}$ in any finite leading section containing row $t$,
	gives
	\begin{equation}\label{eq:global-factor-inverse-entries}
	[\mathcal Q_{p-j}]_{t,s}
	=\prod_{u=s+1}^{t}\widehat a_{p-j}(u)
	=\frac{\pi_j(t)}{\pi_j(s)},\qquad t\ge s\ge j.
	\end{equation}
	Moreover, $\mathcal Q_{p-j}$ fixes every coordinate vector $e_s$ with $s<j$.
	Indeed, $s+1\le j$, so $\widehat a_{p-j}(s+1)=0$.
	Every inverse entry below position $(s,s)$ in that column
	contains this zero as a factor, while the diagonal entry is one.
	
	The relevant columns of $\mathcal Q$ admit expressions as nested
	sums. For $\ell,j\in\{0,\dots,p-1\}$ and $t\ge\ell$, define
	$\chi_{\ell,\ell}(t)\coloneq1$ and
	$\chi_{\ell,j}(t)\coloneq0$ when $j<\ell$. For $j>\ell$, define
	\begin{equation}\label{eq:global-nested-sums}
		\chi_{\ell,j}(t)\coloneq
		\sum_{s=j}^{t}r_{\ell+1}(s)\chi_{\ell+1,j}(s).
	\end{equation}
	An empty sum equals zero, so $\chi_{\ell,j}(t)=0$ whenever $t<j$.
	
	Fix $j\in\{0,\dots,p-1\}$. For
	$\ell\in\{0,\dots,j\}$ and $t\ge j$,
	\[
	[\mathcal Q_{p-\ell}\cdots \mathcal Q_{p-j}]_{t,j}
	=\pi_\ell(t)\chi_{\ell,j}(t).
	\]
	This follows by descending induction on $\ell$. For $\ell=j$,
	\eqref{eq:global-factor-inverse-entries} gives
	$[\mathcal Q_{p-j}]_{t,j}=\pi_j(t)/\pi_j(j)=\pi_j(t)$,
	since $\pi_j(j)=1$; this is the claimed identity because
	$\chi_{j,j}(t)=1$.
	For $\ell<j$, multiply $\mathcal Q_{p-\ell}$ by the remaining
	product. Both matrices are lower triangular, so the intermediate
	index $s$ ranges from $j$ to $t$. In particular, $s\ge j>\ell$,
	so \eqref{eq:global-factor-inverse-entries}, with $\ell$ in place
	of $j$, and the induction hypothesis apply:
	\[
	\begin{aligned}
		[\mathcal Q_{p-\ell}\cdots \mathcal Q_{p-j}]_{t,j}
		&=\sum_{s=j}^{t}
		\frac{\pi_\ell(t)}{\pi_\ell(s)}
		\pi_{\ell+1}(s)\chi_{\ell+1,j}(s)=\pi_\ell(t)\sum_{s=j}^{t}
		r_{\ell+1}(s)\chi_{\ell+1,j}(s)\\
		&=\pi_\ell(t)\chi_{\ell,j}(t).
	\end{aligned}
	\]
	If $i<p-j$, then $j+1\le p-i$, so $\widehat a_i(j+1)=0$ and
	\eqref{eq:signed-bidiagonal-inverse-entries} gives
	$\mathcal Q_i e_j=e_j$.
	Hence column $j$ of $\mathcal Q$ equals column $j$ of
	$\mathcal Q_p\cdots\mathcal Q_{p-j}$, and the claim with
	$\ell=0$ gives $\mathcal Q_{t,j}=\pi_0(t)\chi_{0,j}(t)$ for
	$t\ge j$. For $t<j$, both sides vanish, so this identity holds
	for every $t\ge0$ and $j\in\{0,\dots,p-1\}$.
	
	Subtracting consecutive partial sums in
	\eqref{eq:global-nested-sums} gives
	$\chi_{\ell,j}(t)-\chi_{\ell,j}(t-1)=r_{\ell+1}(t)\chi_{\ell+1,j}(t)$ for $t\ge\ell+1$,
	for $\ell\in\{0,\dots,p-2\}$. For $j>\ell$ and $t\ge j$,
	the difference leaves just the final summand; at $t=j$, the
	preceding sum is empty. For $t<j$, the inequalities
	$\ell+1\le t<j$ force $j\ge\ell+2$, so
	$\chi_{\ell,j}(t)$, $\chi_{\ell,j}(t-1)$ and
	$\chi_{\ell+1,j}(t)$ are all empty sums.
	For $j\le\ell$, both sides are zero by the definitions.
	
	Fix $k\in\{1,\dots,p\}$ and $n\ge0$. Consider the matrix formed
	by rows $n,\dots,n+k-1$ and columns $0,\dots,p-1$ of $\mathcal Q$.
	Number its rows by $i\in\{0,\dots,k-1\}$.
	First divide row $i$ by the positive number $\pi_0(n+i)$.
	Its entries become $\chi_{0,j}(n+i)$.
	
	Next perform $k-1$ stages of row operations; no such stage is
	needed when $k=1$. Before stage $\ell+1$, where
	$\ell\in\{0,\dots,k-2\}$, the rows with $i\ge\ell$ have entries
	$\chi_{\ell,j}(n+i)$, while the rows with $i<\ell$ already have
	their final form.
	
	For $i\in\{\ell+1,\dots,k-1\}$, working in decreasing order,
	subtract row $i-1$ from row $i$ and divide the result by
	$r_{\ell+1}(n+i)$. Since $n+i\ge\ell+1$, this divisor is
	defined and positive. The difference
	identity shows that the new entries of row $i$ are
	\[
	\frac{\chi_{\ell,j}(n+i)-\chi_{\ell,j}(n+i-1)}
	{r_{\ell+1}(n+i)}
	=\chi_{\ell+1,j}(n+i).
	\]
	The decreasing order ensures that row $i-1$ still has its old
	value when it is subtracted. Row $\ell$ is left unchanged.
	
	After the last stage, the resulting matrix has entries
$\mathcal E(n,k)_{i,j}\coloneq\chi_{i,j}(n+i)$ for
$i\in\{0,\dots,k-1\}$ and $j\in\{0,\dots,p-1\}$.
	The row subtractions preserve every determinant obtained by
	selecting $k$ columns. After its initial division by $\pi_0(n+i)$,
	row $i$ is divided once in each of stages $1,\dots,i$ and is
	left unchanged in all later stages. Thus its divisions contribute
	the total factor
$\pi_0(n+i)\prod_{\ell=1}^{i}r_\ell(n+i)=\pi_i(n+i)$,
	where the equality follows by cancellation of the consecutive
	ratios defining $r_\ell$. Thus, for every $k$-element subset
	$K\subseteq\{0,\dots,p-1\}$,
	\begin{equation}\label{eq:global-normalized-minor-identity}
		\frac{\det\mathcal Q[\mathcal I_{k,n},K]}
		{\prod_{i=0}^{k-1}\pi_i(n+i)}
		=\det\mathcal E(n,k)[\{0,\dots,k-1\},K].
	\end{equation}
	
	For $K=\mathcal C_k^{(0)}=\{0,\dots,k-1\}$, the first column set
    in \eqref{eq:global-minor-ratios}, the matrix on the right is upper triangular:
	$\chi_{i,j}(n+i)=0$ for $j<i$, and $\chi_{i,i}(n+i)=1$.
	Its determinant is therefore one. Consequently,
$\Delta_{k,n}=\prod_{i=0}^{k-1}\pi_i(n+i)>0$,
	which proves \eqref{eq:global-leading-minor}.
	
	Now take $k\in\{1,\dots,p-1\}$ and $n\ge0$.
	For the columns $\mathcal C_k^{(1)}$, obtained from $\mathcal C_k^{(0)}$
    by replacing index $k-1$ with $k$, the last row has zero entries
    in the first $k-1$ selected columns. The submatrix formed by the first $k-1$
	rows and those columns is upper triangular with unit diagonal.
	Expanding along the last row therefore gives
	$\det\mathcal E(n,k)[\{0,\dots,k-1\},\mathcal C_k^{(1)}]=\chi_{k-1,k}(n+k-1)$.
	For $k=1$, this is simply the determinant of a one-entry matrix.
	Since $\chi_{k,k}(s)=1$ for $s\ge k$,
\eqref{eq:global-nested-sums} yields
$\chi_{k-1,k}(n+k-1)=\sum_{s=k}^{n+k-1}r_k(s)$.
	Multiplying by $\Delta_{k,n}$ in
	\eqref{eq:global-normalized-minor-identity} proves
	\eqref{eq:global-ratio-series}. For $n=0$, the sum is empty.
	For $n\ge1$, it contains at least one positive term, so
	$\det\mathcal Q[\mathcal I_{k,n},\mathcal C_k^{(1)}]>0$.

\emph{(2) Positivity, monotonicity, and limits.}
Fix $k\in\{1,\dots,p-1\}$. Part~(1) and the definition
\eqref{eq:global-minor-ratios} give
$\kappa_{k,n}^{-1}=\sum_{s=k}^{n+k-1}r_k(s)$ for $n\ge1$.
These partial sums are positive and increase strictly with $n$:
passing from $n$ to $n+1$ adds the positive term $r_k(n+k)$.
Their reciprocals $\kappa_{k,n}$ are therefore positive and
strictly decreasing. The sums tend to $\mathfrak S_k$, so
$\kappa_{k,n}$ tends to $1/\mathfrak S_k$, with limit zero when
the series diverges.
	
\emph{(3) The uniform bound.}
Suppose that all the adjacent series converge.
	The defining sums show that every $\chi_{\ell,j}(t)$ is
	nonnegative. Induction on $j-\ell$ gives
    $0\le\chi_{\ell,j}(t)\le\prod_{h=\ell+1}^{j}\mathfrak S_h$
    for $\ell<j$.
	Indeed, using the bound for $\chi_{\ell+1,j}$ in
	\eqref{eq:global-nested-sums} gives
	\[
		\chi_{\ell,j}(t)
		\le
		\left(\prod_{h=\ell+2}^{j}\mathfrak S_h\right)
		\sum_{s=j}^{t}r_{\ell+1}(s)\le
		\left(\prod_{h=\ell+2}^{j}\mathfrak S_h\right)
		\mathfrak S_{\ell+1}.
	\]
	The last inequality holds because $j\ge\ell+1$, so the partial
	sum is bounded by the full series $\mathfrak S_{\ell+1}$.
	When $j=\ell+1$, the product in this calculation is empty and
	equals one, while $\chi_{\ell+1,j}=1$ on the summation range.
	
	By the definition \eqref{eq:global-minor-bound-constant},
	every entry of the matrix $\mathcal E(n,k)$ constructed in
	part~(1) lies between zero and $\mathfrak B$,
	independently of $n$. The determinant expansion has $k!$ terms,
	each a product of $k$ entries bounded in absolute value by
	$\mathfrak B$. Hence
    $|\det\mathcal E(n,k)[\{0,\dots,k-1\},K]|\le k!\mathfrak B^k$.
	Using \eqref{eq:global-normalized-minor-identity} and
	\eqref{eq:global-leading-minor} proves
	\eqref{eq:global-uniform-minors}.
\end{proof}

The inverse-minor ratios are also convergents of scalar Jacobi
continued fractions, one for each adjacent Metelmann pair. This
extends the tetradiagonal description of \cite{BFMTetra} to
arbitrary bandwidth. Convergence of the positive series to finite,
strictly positive sums is equivalent to strictly positive limits
of these fractions.

Fix $k\in\{1,\dots,p-1\}$. Recall that
$\widehat a_i(n)=[\widehat L_i]_{n,n-1}$ is the subdiagonal
entry of the $i$th Metelmann factor, and define
	\begin{equation}\label{eq:global-jacobi-parameters}
		\alpha_j\coloneq\widehat a_{p-k}(k+j)\quad(j\ge1),\qquad
		\beta_j\coloneq\widehat a_{p-k+1}(k+j)\quad(j\ge0).
	\end{equation}
	These parameters are positive. Consider the scalar Jacobi
	continued fraction
	\begin{equation}\label{eq:global-jacobi-fraction}
		\beta_0-\cfrac{\beta_0\alpha_1}{
			\alpha_1+\beta_1-\cfrac{\beta_1\alpha_2}{
				\alpha_2+\beta_2-\cfrac{\beta_2\alpha_3}{
					\alpha_3+\beta_3-\cdots}}}.
	\end{equation}
	Set $\mathcal K_{k,1}\coloneq\beta_0$. For $N\ge2$, denote by
	$\mathcal K_{k,N}$ its finite truncation ending at the
	denominator $\alpha_{N-1}+\beta_{N-1}$.

\begin{proposition}[Jacobi fractions for the inverse-minor ratios]
	\label{prop:global-jacobi-fractions}
	For each $k\in\{1,\dots,p-1\}$, the inverse-minor ratios
	$\kappa_{k,N}$ are exactly the convergents of a scalar Jacobi
	continued fraction whose coefficients depend only on the
	adjacent Metelmann factors $\widehat L_{p-k}$ and
	$\widehat L_{p-k+1}$.

\begin{enumerate}
\item Every denominator obtained when evaluating these truncations
	from the bottom upward is strictly positive, and
	\begin{equation}\label{eq:global-jacobi-identity}
		\mathcal K_{k,N}
		=\kappa_{k,N}
		=\frac{\beta_0}{
			1+\sum_{j=1}^{N-1}
			\prod_{i=1}^{j}\frac{\alpha_i}{\beta_i}}
		\qquad(N\ge1),
	\end{equation}
	where $\kappa_{k,N}$ is the ratio of minors defined in
	\eqref{eq:global-minor-ratios}. For $N=1$, the sum is empty
	and the common value is $\beta_0$.
	
\item The convergents are positive and strictly
	decreasing, and the infinite fraction, defined as their limit,
	satisfies
	\begin{equation}\label{eq:global-jacobi-limit}
		\mathcal K_k\coloneq\lim_{N\to\infty}\mathcal K_{k,N}
		=\kappa_{k,\infty}
		=\frac1{\mathfrak S_k}\ge0,
	\end{equation}
	with $1/\infty=0$. In particular, this limit is strictly
	positive if and only if the corresponding Metelmann series
	$\mathfrak S_k$ in \eqref{eq:global-adjacent-series} converges.
\end{enumerate}
\end{proposition}
\begin{proof}
For $N\ge2$, denote the successive denominators by $d_j$,
where $1\le j\le N-1$. Evaluation from the bottom upward gives
\[
 d_{N-1}=\alpha_{N-1}+\beta_{N-1},\qquad
 d_j=\alpha_j+\beta_j
       \left(1-\frac{\alpha_{j+1}}{d_{j+1}}\right),
 \quad 1\le j\le N-2.
\]
The first denominator satisfies $d_{N-1}>\alpha_{N-1}>0$.
If $d_{j+1}>\alpha_{j+1}>0$, then
$0<\alpha_{j+1}/d_{j+1}<1$, so $d_j>\alpha_j>0$.
Backward induction proves that every division is valid and that
$\mathcal K_{k,N}=\beta_0(1-\alpha_1/d_1)>0$.

Consider the tridiagonal matrix $\mathcal J_N$ of order $N$,
with indices $0,\dots,N-1$, whose nonzero entries are
\[
 \begin{aligned}
 [\mathcal J_N]_{0,0}&\coloneq \beta_0,&
 [\mathcal J_N]_{j,j}&\coloneq \alpha_j+\beta_j,\\
 [\mathcal J_N]_{j-1,j}&\coloneq 1,&
 [\mathcal J_N]_{j,j-1}&\coloneq \alpha_j\beta_{j-1},
 \end{aligned}
 \qquad 1\le j\le N-1.
\]
Direct multiplication gives $\mathcal J_N=\mathcal L_N\mathcal U_N$,
where $\mathcal L_N$ is unit lower bidiagonal with subdiagonal
$\alpha_j$, and $\mathcal U_N$ is upper bidiagonal with diagonal
$\beta_j$ and superdiagonal $1$. Both matrices are invertible:
their diagonal entries are nonzero. Solving the triangular equations
for their inverses yields
\[
 [\mathcal U_N^{-1}]_{0,j}
 =\frac{(-1)^j}{\beta_0\cdots\beta_j},\qquad
 [\mathcal L_N^{-1}]_{j,0}=(-1)^j\prod_{i=1}^{j}\alpha_i.
\]
Since $\mathcal J_N^{-1}=\mathcal U_N^{-1}\mathcal L_N^{-1}$,
the signs cancel in each summand of its entry $(0,0)$, giving
\[
 [\mathcal J_N^{-1}]_{0,0}
 =\frac1{\beta_0}\left(
 1+\sum_{j=1}^{N-1}\prod_{i=1}^{j}\frac{\alpha_i}{\beta_i}
 \right).
\]

Its reciprocal equals the stated fraction, with positive denominators
throughout elimination, as follows. Let $\delta_j$ be the determinant
of the trailing principal block with indices $j,\dots,N-1$,
and set $\delta_N\coloneq 1$. Then
$\delta_{N-1}=\alpha_{N-1}+\beta_{N-1}$ and expansion of the
first row of each tridiagonal block gives
$\delta_j=(\alpha_j+\beta_j)\delta_{j+1}-\beta_j\alpha_{j+1}\delta_{j+2}$ for $1\le j\le N-2$.
Starting with $\delta_{N-1}/\delta_N=d_{N-1}$, this recurrence
shows inductively that $\delta_j/\delta_{j+1}=d_j$ and
$\delta_j>0$ for $1\le j\le N-1$.
The first row of the full matrix gives
\[
 \frac{\det\mathcal J_N}{\delta_1}
 =\beta_0-\beta_0\alpha_1\frac{\delta_2}{\delta_1}
 =\beta_0\left(1-\frac{\alpha_1}{d_1}\right)
 =\mathcal K_{k,N}.
\]
The cofactor formula gives
$[\mathcal J_N^{-1}]_{0,0}=\delta_1/\det\mathcal J_N$.
Combining the two calculations proves the rational expression in
\eqref{eq:global-jacobi-identity}. For $N=1$, the matrix is
$[\beta_0]$, its complementary minor is the empty determinant $1$,
and the same expression is simply $\beta_0$.

Finally, the definitions of the adjacent products give
\[
 r_k(k+j)=\frac1{\beta_0}
             \prod_{i=1}^{j}\frac{\alpha_i}{\beta_i},\qquad j\ge0.
\]
The extra factor $\beta_0^{-1}$ is the factor at index $k$ in
the denominator product; the later factors are exactly the displayed
ratios. With the change of index $s=k+j$,
\eqref{eq:global-ratio-series} says that $\kappa_{k,N}$ is the
reciprocal of the sum over $j\in\{0,\dots,N-1\}$. This is the value
already found for $\mathcal K_{k,N}$. Adding the next positive
summand increases the denominator strictly, so the convergents
decrease strictly. Their denominators tend to $\mathfrak S_k$,
which proves \eqref{eq:global-jacobi-limit}, including the case
of an infinite sum and limit zero.
\end{proof}

For $p=2$, the first and second subdiagonals of $L$ are
$m_1=\beta_0$, $m_{j+1}=\alpha_j+\beta_j$, and
$\ell_{j+1}=\alpha_j\beta_{j-1}$ for $j\ge1$.
Thus \eqref{eq:global-jacobi-fraction} is exactly
$m_1-\ell_2/(m_2-\ell_3/(m_3-\cdots))$, the fraction used in the
tetradiagonal analysis of \cite{BFMTetra}.
For arbitrary $p$, each adjacent Metelmann pair gives one such
fraction. Its coefficients are fixed by $L$, independently of any
positive factorization to be constructed.

A global positive factorization can be constructed by removing one
positive bidiagonal factor at a time. Proposition~\ref{prop:successive-recurrences}
determines such a factor from its initial entries, provided the
recurrence produces positive entries in every later row. The next
lemma gives a determinant condition that guarantees this positivity
and reduces the lower bandwidth by one.

For the data of this construction, let $F$ be semi-infinite, unit lower
triangular, of lower bandwidth $m\ge2$. Write
$f_s(n)\coloneq F_{n,n-s}$ and assume $f_m(n)>0$ for $n\ge m$.
Choose positive numbers $a(1),\dots,a(m-1)$ as the first $m-1$
subdiagonal entries of the unit lower bidiagonal factor $\mathcal L(a)$
in the desired factorization $F=\mathcal L(a)F_{\mathrm{rem}}$, where
$F_{\mathrm{rem}}$ denotes the remaining unit lower triangular factor.
Define $\varphi$ by \eqref{eq:extraction-initial-data} and
\eqref{eq:extraction-linear-recurrence}, with $r=m$,
$\theta_h=a(h)$, and coefficients $f_s$, continuing the recurrence
for every $n\ge m$.

Put $R\coloneq \Sigma F^{-1}\Sigma$, where
$\Sigma=\diag[(-1)^n]_{n\in\NN}$ is the diagonal sign matrix
used in Lemma~\ref{lem:global-inverse-minors}; the triangular inverse
is computed entrywise. Let $M$ be the
$m\times(m-1)$ matrix with diagonal entries $1$, subdiagonal entries
$-a(1),\dots,-a(m-1)$, and all other entries zero. For $n\ge0$, set
\begin{equation}\label{eq:extraction-test-determinants}
 \mathcal I_{m-1,n}\coloneq\{n,\dots,n+m-2\},\qquad
 \mathcal W_n\coloneq\det\bigl(R[\mathcal I_{m-1,n},\{0,\dots,m-1\}]M\bigr).
\end{equation}
Thus each determinant uses $m-1$ consecutive rows and the first $m$
columns of $R$, in increasing order, with this submatrix multiplied
by $M$. Each $\mathcal W_n$ depends only on $F$ and the
chosen initial entries.

\begin{lemma}[A determinant test for a positive bidiagonal factor]
\label{lem:global-extraction-determinant}
With the preceding data, if the determinants $\mathcal W_n$ in
\eqref{eq:extraction-test-determinants} are positive for every $n\ge1$, then
\[
 F=\mathcal L(a)F_{\mathrm{rem}},\qquad
 a(n)=\frac{\varphi_{n-1}}{\varphi_n}>0\quad(n\ge1).
\]
Here $\mathcal L(a)$ is unit lower bidiagonal with the prescribed
initial entries, and $F_{\mathrm{rem}}$ is unit lower triangular of lower
bandwidth $m-1$, with positive entries on its $(m-1)$st subdiagonal.

For every choice of the positive initial entries, the exact identity is
\begin{equation}\label{eq:global-extraction-determinant}
 \mathcal W_n=\varphi_{n+m-1}\prod_{t=m}^{n+m-1}f_m(t),
 \qquad n\ge0.
\end{equation}
This identity holds even when some later term of $\varphi$ vanishes.
The product is empty when $n=0$, and $\mathcal W_0=\varphi_{m-1}=1$.
\end{lemma}

The identity is a discrete form of Abel's identity, related to
determinant methods for factorizations of difference operators; see
\cite{Hartman1978,HarrisKruegerTrench1998,ErbeMathsenPeterson2001}.
\begin{proof}
	The determinant identity will first be established and then used to
	prove positivity of the extracted factor. To include cases where
	some later value of $\varphi$ vanishes, first treat the initial parameters
	$a(1),\dots,a(m-1)$ as independent variables.
	
	The initial values are $\varphi_{m-1}=1$ and
    $\varphi_j=\prod_{h=j+1}^{m-1}a(h)$ for $j\in\{0,\dots,m-2\}$.
	These are distinct monomials and hence are linearly independent
	as polynomials in the initial parameters. Every later
	$\varphi_n$ is also a polynomial: the recurrence only divides
	by the fixed nonzero numbers $f_m(n)$.
	
	Moreover, no $\varphi_n$ is the zero polynomial. Indeed, each
	recurrence step sends
    $(\varphi_{n-1},\dots,\varphi_{n-m})$ to
    $(\varphi_n,\dots,\varphi_{n-m+1})$
	by a linear transformation. The coefficient of
	$\varphi_{n-m}$ in the formula for $\varphi_n$ is
	$(-1)^{m-1}/f_m(n)\ne0$.
	The preceding $m-1$ values are retained in the new group, and
	the recurrence recovers the remaining old value $\varphi_{n-m}$
	from them and $\varphi_n$, because its coefficient is nonzero.
	Thus the transformation is invertible and preserves linear
	independence. Every consecutive group of $m$ values remains
	linearly independent as polynomials, so each individual value
	is a nonzero polynomial.
	
	The ratios $a(n)\coloneq\varphi_{n-1}/\varphi_n$, for $n\ge1$,
	are therefore defined in the field $\mathbb R(a(1),\dots,a(m-1))$ of rational
	functions in the initial parameters. Dividing in this field
	requires a nonzero polynomial, which was just established.
	For each determinant index $n$, only finitely many of these
	ratios enter the calculation. Specialization of the parameters
	follows after obtaining a polynomial identity for that index.
	For $n\in\{1,\dots,m-1\}$, these ratios recover the chosen
	initial entries. The algebraic construction in
	Proposition~\ref{prop:successive-recurrences} gives
    $F=\mathcal L(a)F_{\mathrm{rem}}$, with
    $F_{\mathrm{rem}}=\mathcal L(a)^{-1}F$,
	where $F_{\mathrm{rem}}$ is unit lower triangular of bandwidth at most
	$m-1$.
	
	Use the determinants $\mathcal W_n$ in \eqref{eq:extraction-test-determinants}.
	To explain the appearance of the matrix $M$ of initial entries,
    invert the factorization to obtain
    $F_{\mathrm{rem}}^{-1}=F^{-1}\mathcal L(a)$.
	Since $\Sigma^2=I$, this gives
    $S\coloneq \Sigma F_{\mathrm{rem}}^{-1}\Sigma=R(\Sigma\mathcal L(a)\Sigma)$.
	The matrix $\Sigma\mathcal L(a)\Sigma$ has diagonal entries $1$ and
	subdiagonal entries $-a(n)$. Its first $m-1$ columns consist
	of $M$ followed by zero rows. Thus the first $m-1$ columns
	of $S$ are obtained by multiplying the first $m$ columns
	of $R$ by $M$.
	
	Let $v_t$ denote row $t$ of these first $m-1$ columns of $S$,
	and write
    $f_s^{\mathrm{rem}}(t)\coloneq[F_{\mathrm{rem}}]_{t,t-s}$.
	Then
	\[
	\mathcal W_n=
	\det\begin{bNiceMatrix}
		v_n\\
		v_{n+1}\\
		\Vdots\\
		v_{n+m-2}
	\end{bNiceMatrix}.
	\]
	A recurrence for these determinants follows from the inverse identity.
	The identity $(\Sigma F_{\mathrm{rem}}\Sigma)S=I$
	implies
    $v_t+\sum_{s=1}^{m-1}(-1)^sf_s^{\mathrm{rem}}(t)v_{t-s}=0$
    for $t\ge m-1$.
	Indeed, the entries of $\Sigma F_{\mathrm{rem}}\Sigma$ on its $s$th subdiagonal are
	$(-1)^s f_s^{\mathrm{rem}}(t)$, and the first $m-1$ entries of row $t\ge m-1$
	of the identity matrix are zero. Hence
	$v_t=\sum_{s=1}^{m-1}(-1)^{s+1}f_s^{\mathrm{rem}}(t)v_{t-s}$.
	
	In the determinant $\mathcal W_{n+1}$, substitute this expression for
	its last row $v_{n+m-1}$. For
	$s\in\{1,\dots,m-2\}$, the row $v_{n+m-1-s}$ already occurs
	among the preceding rows, so the corresponding determinant
	vanishes. Only the term containing $v_n$ remains, with
	coefficient $(-1)^mf_{m-1}^{\mathrm{rem}}(n+m-1)$.
	Moving $v_n$ from the last row to the first requires $m-2$
	row exchanges and contributes the sign $(-1)^{m-2}$.
	The two signs cancel, giving
	$\mathcal W_{n+1}=f_{m-1}^{\mathrm{rem}}(n+m-1)\mathcal W_n$ for $n\ge0$.
	Also, $\mathcal W_0=1$, because the leading $(m-1)\times(m-1)$ block of
	the unit lower triangular matrix $S$ has determinant $1$.
	
	It remains to express the coefficient $f_{m-1}^{\mathrm{rem}}$ in terms of
	the original matrix. Comparing the entries on the $m$th
	subdiagonal in $F=\mathcal L(a)F_{\mathrm{rem}}$ gives
    $f_m(t)=a(t)f_{m-1}^{\mathrm{rem}}(t-1)$ for $t\ge m$.
	There is no contribution from the $m$th subdiagonal of
	$F_{\mathrm{rem}}$, since its bandwidth is at most $m-1$. Iterating the
	recurrence for $\mathcal W_n$ therefore yields
    $\mathcal W_n=\prod_{t=m}^{n+m-1}(f_m(t)/a(t))$.
	The factors involving $a(t)$ telescope:
	\[
	\prod_{t=m}^{n+m-1}\frac1{a(t)}
	=
	\prod_{t=m}^{n+m-1}\frac{\varphi_t}{\varphi_{t-1}}
	=
	\frac{\varphi_{n+m-1}}{\varphi_{m-1}}
	=
	\varphi_{n+m-1}.
	\]
	Thus
	\[
	\mathcal W_n=
	\left(\prod_{t=m}^{n+m-1}f_m(t)\right)
	\varphi_{n+m-1},
	\]
	including $n=0$, when the products are empty.
	
	Although the calculation used ratios of the polynomials
	$\varphi_n$, both sides of the resulting identity are
	polynomials in the initial parameters. On the left, $R$
	is fixed and the entries of $M$ are polynomials in those
	parameters. On the right, the recurrence introduces only
	fixed nonzero denominators depending on the coefficients
	$f_m(t)$. Equality as rational functions therefore gives
	a polynomial identity. It holds for every choice of the
	initial parameters, including choices for which some
	$\varphi_n$ vanishes. This proves
	\eqref{eq:global-extraction-determinant}.
	
	Finally, take positive initial parameters and suppose
	$\mathcal W_n>0$ for every $n\ge1$. The initial values
	$\varphi_0,\dots,\varphi_{m-1}$ are positive by construction.
	Since every $f_m(t)$ in the displayed product is positive,
	the determinant identity gives
	$\varphi_{n+m-1}>0$ for every $n\ge1$.
	Hence all the ratios $a(n)=\varphi_{n-1}/\varphi_n$
	are well defined and strictly positive.
	Proposition~\ref{prop:successive-recurrences} now applies
	in every row and gives the required positive factor
	$\mathcal L(a)$. The remaining matrix has bandwidth at
	most $m-1$, and
    $f_{m-1}^{\mathrm{rem}}(t)=f_m(t+1)/a(t+1)>0$
    for $t\ge m-1$.
	Its bandwidth is therefore exactly $m-1$, with positive
	extreme subdiagonal.
\end{proof}

The next lemma makes the choice of initial parameters uniform in the
row index. When the adjacent Metelmann series converge, the normalized
inverse minors have bounds independent of the row. Cauchy--Binet then
controls all the determinants needed in the extraction.

\begin{uniformchoicelemma}[Uniform choice of initial parameters]
Let $L=\widehat L_1\cdots\widehat L_p$ be a compatible unit lower
Metelmann product as in \eqref{eq:metelmann}, with $p\ge2$.
Assume that all its adjacent series
$\mathfrak S_1,\dots,\mathfrak S_{p-1}$ converge. Let $L^f$ be the
leading section of $L$ of order $p+1$, with row and column indices
$\{0,\dots,p\}$, and let
$L^f=L_1^f(\boldsymbol u)\cdots L_p^f(\boldsymbol u)$ be the family
in Lemma~\ref{lem:rowwise}, with the $p(p-1)/2$ free initial
parameters $\boldsymbol u$.
	Write $\Sigma^f\coloneq\diag[1,-1,\dots,(-1)^p]$.
	For $r\in\{0,\dots,p-1\}$, let $X_r$ be the
	$p\times(p-r)$ matrix formed by the first $p$ rows and
	the first $p-r$ columns of $\Sigma^f(L_1^f\cdots L_r^f)\Sigma^f$.

	For $m\in\{1,\dots,p-1\}$, define
	\begin{equation}\label{eq:global-W-determinants}
		\mathcal W_{m,n}\coloneq
		\det\bigl(
		\mathcal Q[\mathcal I_{m,n},\{0,\dots,p-1\}]X_{p-m}
		\bigr),
		\qquad n\ge0.
	\end{equation}

For the same values of $m$, define
	\begin{equation}\label{eq:global-initial-minor-error}
	\epsilon_m\coloneq
	\sum_{\substack{K\subseteq\{0,\dots,p-1\}\\|K|=m}}
	\left|
	\det X_{p-m}[K,\{0,\dots,m-1\}]
	-\mathbf1_{K=\mathcal C_m^{(0)}}
	\right|.
	\end{equation}
	Here $\mathbf1_{K=\mathcal C_m^{(0)}}$ equals $1$ if $K=\mathcal C_m^{(0)}$ and
	$0$ otherwise.
Then, with $\mathfrak B$ as in
\eqref{eq:global-minor-bound-constant},
\[
 \left|\frac{\mathcal W_{m,n}}{\Delta_{m,n}}-1\right|
 \le m!\mathfrak B^m\epsilon_m
 \qquad\bigl(m\in\{1,\dots,p-1\},\ n\ge0\bigr).
\]
There is one neighborhood of $\boldsymbol u=0$ such that every
positive parameter vector in it gives positive finite factors and
$\mathcal W_{m,n}>\frac12\Delta_{m,n}>0$ for every
$m\in\{1,\dots,p-1\}$ and $n\ge0$.
\end{uniformchoicelemma}
\begin{proof}
Lemma~\ref{lem:rowwise} gives continuous dependence of the finite
factors on $\boldsymbol u$. At $\boldsymbol u=0$, they are the leading
sections of the prescribed Metelmann factors. For every sufficiently
small positive parameter vector, all the finite factors are positive.

	For $r=0$, the product is the identity and $X_0=I_p$.
	A product of $r$ lower bidiagonal matrices has bandwidth
	at most $r$. Thus its column $j$ vanishes below row
	$j+r$. For the columns selected here, $j+r\le(p-r-1)+r=p-1$.
	Consequently, $X_r$ contains all the possibly nonzero
	entries of these columns.
	
	At the Metelmann parameters, each of the first $r$ factors
	fixes $e_0,\dots,e_{p-r-1}$. To check the indices, let
	$i\in\{1,\dots,r\}$ and $j\in\{0,\dots,p-r-1\}$. Then
	$j+1\le p-r\le p-i$, so
	$\widehat a_i(j+1)=0$ by
	\eqref{eq:metelmann-parameters}; this is the only possible
	off-diagonal entry in column $j$ of $\widehat L_i$.
	Conjugating by the sign matrix also fixes that coordinate
	column, since its diagonal entry remains one.
	Hence
	\[
	X_r\longrightarrow
	[e_0\ \cdots\ e_{p-r-1}]
	\qquad\text{as }\boldsymbol u\longrightarrow0,
	\]
	where the coordinate vectors in this display belong
	to $\mathbb R^p$.
	
	At $\boldsymbol u=0$, the matrix $X_{p-m}$ consists
	of the first $m$ coordinate columns. Its minor with
	row set $K$ is therefore $1$ when $K=\mathcal C_m^{(0)}$ and $0$
	otherwise. The changes in these minors are measured by
	\eqref{eq:global-initial-minor-error}.
	There are finitely many terms, each
	continuous in $\boldsymbol u$, so
	$\epsilon_m\longrightarrow0$ as
	$\boldsymbol u\longrightarrow0$.
	
	Cauchy--Binet gives
	\[
		\frac{\mathcal W_{m,n}}{\Delta_{m,n}}-1=
		\sum_{\substack{K\subseteq\{0,\dots,p-1\}\\|K|=m}}
		\frac{\det\mathcal Q[\mathcal I_{m,n},K]}{\Delta_{m,n}}
		\left(
		\det X_{p-m}[K,\{0,\dots,m-1\}]
		-\mathbf1_{K=\mathcal C_m^{(0)}}
		\right).
	\]
	The normalizing minor is
    $\Delta_{m,n}=\det\mathcal Q[\mathcal I_{m,n},\mathcal C_m^{(0)}]$
    from \eqref{eq:global-minor-ratios}; this explains the term $1$
    subtracted on the left.
	
	All the adjacent series converge by assumption.
	The bound \eqref{eq:global-uniform-minors} therefore
	applies to every normalized minor in this sum, with
	the constant $\mathfrak B$ defined in
    \eqref{eq:global-minor-bound-constant}, independently of the row index.
	Taking absolute values gives the asserted bound.
	In particular,
	\begin{equation}\label{eq:global-positive-margin}
		\frac{\mathcal W_{m,n}}{\Delta_{m,n}}
		\ge1-m!\mathfrak B^m\epsilon_m
		\qquad(n\ge0).
	\end{equation}
	
	There are only $p-1$ values of $m$, and each
	$\epsilon_m$ tends to zero. Thus one neighborhood of
	$\boldsymbol u=0$ can be chosen so that
    $m!\mathfrak B^m\epsilon_m<1/2$ for every $m\in\{1,\dots,p-1\}$.
	Choose this neighborhood small enough that its
	positive parameter vectors also give positive finite
	factors. For every such vector,
	\[
	\mathcal W_{m,n}>
	\frac12\Delta_{m,n}>0
	\qquad
	\bigl(m\in\{1,\dots,p-1\},\ n\ge0\bigr).
	\]
	The same choice of parameters thus works for every row index.
	More explicitly, the finitely many conditions on the continuous
	functions $\epsilon_m$ and positivity of the finite factor entries
	that are nonzero at $\boldsymbol u=0$ define one open
	neighborhood of zero. It contains a set
	$\{\boldsymbol u: \max_j|u_j|<\eta\}$ for some $\eta>0$.
	Every vector with $0<u_j<\eta$ satisfies all the displayed
	determinant inequalities, for every $m$ and every $n$.
	The uniform normalized-minor bound is what permits this same
	$\eta$ to serve all row indices.
\end{proof}

The uniform choice of parameters and the determinant formula for
extraction together construct global positive factors.
Conversely, any global positive factorization imposes strictly
positive lower bounds on the inverse-minor ratios. Their limits
must therefore be positive. The identities already established
express this condition equivalently through convergence of the
adjacent series or positivity of the scalar Jacobi fractions.

\begin{theorem}[Global positive factorization in arbitrary lower bandwidth]
\label{thm:global-lower-criterion}
For a compatible unit lower Metelmann product
$L=\widehat L_1\cdots\widehat L_p$ as in
\eqref{eq:metelmann}, the following are equivalent:
\begin{enumerate}
\item $L=L_1\cdots L_p$ with every lower bidiagonal factor unit and
every subdiagonal entry strictly positive.
\item All the adjacent series $\mathfrak S_1,\ldots,\mathfrak S_{p-1}$
in \eqref{eq:global-adjacent-series} converge.
\item The limits $\kappa_{k,\infty}=\lim_{n\to\infty}\kappa_{k,n}$
of the inverse-minor ratios \eqref{eq:global-minor-ratios} are positive
for every $k\in\{1,\dots,p-1\}$.
\item Every Jacobi fraction $\mathcal K_k$ in
\eqref{eq:global-jacobi-limit} has a strictly positive value,
for $k\in\{1,\dots,p-1\}$.
\end{enumerate}
When these conditions hold, all sufficiently small positive choices
of the $p(p-1)/2$ parameters replacing the initial zeros in
\eqref{eq:metelmann-parameters}, as in Lemma~\ref{lem:rowwise},
extend to positive factorizations of the same $L$ in every row.
For $p=1$, the factor is already positive and all lists of conditions
are empty. The result applies entrywise to bounded and unbounded
Metelmann products.
\end{theorem}
\begin{proof}
	For $p=1$, the matrix $L$ is already a positive bidiagonal
	factor, and there are no series or limits to check. Assume
	henceforth that $p\ge2$.
	
	Lemma~\ref{lem:global-inverse-minors} gives
	$\kappa_{k,\infty}=1/\mathfrak S_k$, so conditions (2) and
	(3) are equivalent. Proposition~\ref{prop:global-jacobi-fractions}
	gives $\mathcal K_k=\kappa_{k,\infty}$, proving equivalence
	with (4). It remains to prove that (1) implies (3) and that
	(2) implies (1).
	
	\emph{A positive factorization gives positive lower bounds.}
	Assume that $L=L_1\cdots L_p$
	with all the stated positivity properties.
	
	For $k=1$ and a fixed $n\ge1$, choose a leading section
	of order $d\ge\max(p+1,n+1)$. The inverse entries defining
	$\kappa_{1,n}$ agree with those of this finite section.
	Proposition~\ref{prop:inverse-parameter-bound} gives
	$[L_1]_{1,0}<\kappa_{1,d-1}$, and monotonicity gives
	$\kappa_{1,d-1}\le\kappa_{1,n}$. Hence
	$\kappa_{1,n}>[L_1]_{1,0}>0$ for every $n\ge1$.
	The lower bound is the same fixed entry of the global factor
	$L_1$ in every section. Taking the limit therefore yields
	$\kappa_{1,\infty}\ge[L_1]_{1,0}>0$.
	
	Now fix $k\in\{2,\dots,p-1\}$ and put
    $A\coloneq L_1\cdots L_k$ and $h_j\coloneq A_{j,0}$ for
    $j\in\{0,\dots,k\}$.
    Also put $V\coloneq \Sigma A^{-1}\Sigma$ and
    $R\coloneq \Sigma(L_{k+1}\cdots L_p)^{-1}\Sigma$.
	The signed inverse $\mathcal Q=\Sigma L^{-1}\Sigma$ used in the ratios
    \eqref{eq:global-minor-ratios} therefore satisfies $\mathcal Q=RV$.
	Since $A$ is a product of $k$ positive bidiagonal factors,
	$h_0=1$ and $h_j>0$ for $j\in\{1,\dots,k\}$.
	Indeed, each such entry contains a positive contribution
	obtained by taking subdiagonal entries from $j$ factors
	and diagonal entries from the remaining factors.
	
	As shown in Proposition~\ref{prop:inverse-parameter-bound},
	each matrix $\Sigma L_i^{-1}\Sigma$ has all its minors nonnegative.
	The matrices $V$ and $R$ are products of these matrices,
	so their minors are also nonnegative.
	
	The matrix $A$ has bandwidth $k$. Its first column
	therefore vanishes below row $k$, and the identity
	$V(\Sigma A\Sigma)=I$ gives
	\(\sum_{j=0}^{k}(-1)^j h_j V[:,j]=e_0\).
	Here $V[:,j]$ denotes column $j$ of $V$, and $e_0$ is
	the first coordinate vector.
	
	Let $I$ be any increasing set of $k$ row indices.
	Denote by $D_I$ the determinant whose columns are
    $V[:,0],\dots,V[:,k-2],e_0$,
	all restricted to the rows in $I$.
	If $0\notin I$, its last column is zero and $D_I=0$.
	If $0\in I$, row zero is the first selected row.
	Expansion along the last column then gives
	$(-1)^{k-1}$ times a minor of $V$, which is nonnegative.
	In both cases, $(-1)^{k-1}D_I\ge0$.
	
	Substitution of the column identity for $e_0$ gives a second
	expression for $D_I$. The terms with
	$j\in\{0,\dots,k-2\}$ repeat an existing column and
	therefore vanish. Recall that
    $\mathcal C_k^{(0)}=\{0,\dots,k-1\}$ and
    $\mathcal C_k^{(1)}=\{0,\dots,k-2,k\}$ are the two column sets
    used in \eqref{eq:global-minor-ratios}. The remaining terms give
	\[
	D_I
	=
	(-1)^{k-1}h_{k-1}\det V[I,\mathcal C_k^{(0)}]
	+
	(-1)^k h_k\det V[I,\mathcal C_k^{(1)}].
	\]
	Multiplication by $(-1)^{k-1}$ and division by
	$h_{k-1}>0$ give
	\[
	\det V[I,\mathcal C_k^{(0)}]
	-
	\frac{h_k}{h_{k-1}}\det V[I,\mathcal C_k^{(1)}]
	\ge0.
	\]
	
	To obtain the corresponding inequality for $\mathcal Q$,
	apply Cauchy--Binet to $\mathcal Q=RV$ on the consecutive rows
    $\mathcal I_{k,n}=\{n,\dots,n+k-1\}$. For every $n\ge1$,
	\begin{multline*}
		\det\mathcal Q[\mathcal I_{k,n},\mathcal C_k^{(0)}]
		-\frac{h_k}{h_{k-1}}
		\det\mathcal Q[\mathcal I_{k,n},\mathcal C_k^{(1)}]\\
		=
		\sum_{\substack{K\subseteq\{0,\dots,n+k-1\}\\|K|=k}}
		\det R[\mathcal I_{k,n},K]
		\left(
		\det V[K,\mathcal C_k^{(0)}]
		-\frac{h_k}{h_{k-1}}\det V[K,\mathcal C_k^{(1)}]
		\right)
		\ge0.
	\end{multline*}
	Every summand is nonnegative by the preceding inequality
	and the nonnegativity of the minors of $R$.
	The sum is finite because $R$ is lower triangular:
	intermediate indices larger than $n+k-1$ contribute zero.
	
	For $n\ge1$, the denominator
	$\det\mathcal Q[\mathcal I_{k,n},\mathcal C_k^{(1)}]$ is positive by
	Lemma~\ref{lem:global-inverse-minors}.
	Dividing by this determinant gives
	\begin{equation}\label{eq:global-necessary-prefix-bound}
		\kappa_{k,n}\ge\frac{h_k}{h_{k-1}}>0
		\qquad(n\ge1).
	\end{equation}
	Again, the lower bound is independent of $n$.
	Taking limits proves (3).
	
	\emph{Choosing initial parameters that work in every row.}
	Assume (2). Apply the preceding lemma on the uniform choice of
	initial parameters. Use its finite family
	$L^f=L_1^f(\boldsymbol u)\cdots L_p^f(\boldsymbol u)$,
	matrices $X_r$, and determinants $\mathcal W_{m,n}$ from
	\eqref{eq:global-W-determinants}. The lemma gives one neighborhood
	of $\boldsymbol u=0$ whose positive parameter vectors give positive
	finite factors and satisfy
	$\mathcal W_{m,n}>\frac12\Delta_{m,n}>0$ for every
	$m\in\{1,\dots,p-1\}$ and $n\ge0$.
	
	\emph{Constructing the global factors.}
	Fix any positive parameter vector in this neighborhood.
	Construct the factors one at a time.
	
	Suppose that $L_1,\dots,L_r$ have already been
	constructed as positive bidiagonal matrices and agree
	with $L_1^f,\dots,L_r^f$ on the leading section of
	order $p+1$. Suppose also that $F_r\coloneq(L_1\cdots L_r)^{-1}L$
	is unit lower triangular of bandwidth $m\coloneq p-r$, with
	positive extreme subdiagonal.
	Its leading section of order $p+1$ is then $L_{r+1}^f\cdots L_p^f$,
	because leading sections of lower triangular products
	and inverses agree with the corresponding finite
	products and inverses.
	For $r=0$, these statements hold with $F_0=L$.
	Its extreme subdiagonal is positive: in row $n\ge p$,
	the unique term contributing to $L_{n,n-p}$ is
	$\prod_{i=1}^{p}\widehat a_i(n-i+1)$, and
	$n-i+1\ge p-i+1$ is in the positive range of each factor.
	
	Assume $m\ge2$ and put $R_r\coloneq \Sigma F_r^{-1}\Sigma$.
	Inverting the definition of $F_r$ gives
    $R_r=\mathcal Q\Sigma(L_1\cdots L_r)\Sigma$.
	For $j\in\{0,\dots,m-1\}$, column $j$ of
	$\Sigma(L_1\cdots L_r)\Sigma$ vanishes below row
	$j+r\le p-1$. Its remaining entries agree with those of $X_r$, the block
    of $\Sigma^f(L_1^f\cdots L_r^f)\Sigma^f$ formed from rows $0,\dots,p-1$
    and columns $0,\dots,p-r-1$.
	Consequently,
    $[R_r]_{i,j}=\sum_{h=0}^{p-1}\mathcal Q_{i,h}[X_r]_{h,j}$
    for $i\ge0$ and $j\in\{0,\dots,m-1\}$.
	Only these first $m$ columns of $R_r$ are needed for
	the next extraction.
	
	Choose the first $m-1$ subdiagonal entries of
	$L_{r+1}^f$ as the initial data for extracting a factor
	from $F_r$. They are positive. Let $M_r$ be the
	$m\times(m-1)$ matrix with diagonal entries $1$ and
	subdiagonal entries equal to the negatives of these
	initial data, as in
	Lemma~\ref{lem:global-extraction-determinant}.
	
	To pass from the product of the first $r$ finite
	factors to that of the first $r+1$ factors, multiply
	on the right by $L_{r+1}^f$. After multiplication by
	$\Sigma^f$ on both sides, the corresponding right factor
	is $\Sigma^fL_{r+1}^f\Sigma^f$. Its first $m-1$ columns consist
	of $M_r$ followed by zero rows. Thus $X_rM_r=X_{r+1}$.
	It follows that the test determinants
    \eqref{eq:extraction-test-determinants}, with $F=F_r$ and $M=M_r$, are
	\[
	\begin{aligned}
		\det\bigl(
		R_r[\mathcal I_{m-1,n},\{0,\dots,m-1\}]M_r
		\bigr)&=
		\det\bigl(
		\mathcal Q[\mathcal I_{m-1,n},\{0,\dots,p-1\}]X_{r+1}
		\bigr)\\
		&=\mathcal W_{m-1,n}>0.
	\end{aligned}
	\]
	The last equality uses
	$r+1=p-(m-1)$ in
	\eqref{eq:global-W-determinants}.
	
	Lemma~\ref{lem:global-extraction-determinant} therefore
	constructs a factor $L_{r+1}$ with positive
	subdiagonal entries in every row. It leaves $F_{r+1}=L_{r+1}^{-1}F_r$
	with bandwidth $m-1$ and positive extreme subdiagonal.
	These are exactly the hypotheses on the remaining matrix
	needed at the next application of the determinant lemma.
	
	The new factor must also agree with $L_{r+1}^f$ on the chosen
	finite section. On that
	section, $F_r$ agrees with
	$L_{r+1}^f\cdots L_p^f$. The converse statement in
	Proposition~\ref{prop:successive-recurrences} shows
	that $L_{r+1}^f$ satisfies the extraction recurrence
	for these same coefficients. Its initial entries
	are exactly those used to construct $L_{r+1}$.
	The recurrence has a nonzero leading coefficient
	in every row, so it determines the auxiliary
	sequence uniquely from its initial data.
	The resulting bidiagonal entries therefore agree
	throughout the finite section. Removing this common
	factor also makes the leading section of $F_{r+1}$
	equal to
	$L_{r+2}^f\cdots L_p^f$.
	This proves all the assertions needed for the next
	induction step.
	
	After $p-1$ extractions, the remaining matrix is
	unit lower triangular of bandwidth one, with
	positive subdiagonal. Take this remaining matrix
	as $L_p$. This gives $L=L_1\cdots L_p$
	with every factor positive, proving (1).
	
	The parameter vector used in this construction was
	an arbitrary positive vector in the neighborhood
	chosen above. Hence every sufficiently small
	positive choice of the $p(p-1)/2$ initial parameters
	extends to a global positive factorization, as
	claimed. Each construction and each product identity is
	defined on finite leading sections and agrees on their overlaps.
	Thus it defines every entry of the global factors, including
	when their entries are unbounded.
\end{proof}

Apply the preceding theorem separately to $L$ and $U^{\mathsf T}$.
Uniqueness of normalized LDU fixes $L$ and $U$ as the triangular
products in every PBF of $T$, so the global criterion has independent
lower and upper conditions.

\begin{corollary}[The criterion for arbitrary lower and upper bandwidths]
\label{cor:global-two-sided-criterion}
Let $T=LDU$ be a semi-infinite BTP matrix of bandwidth $(p,q)$,
where $p,q\in\mathbb N$, with its normalized compatible triangular
decomposition.
From $L$, form the inverse-minor ratios $\kappa^-_{k,n}$ of
\eqref{eq:global-minor-ratios} and put
$\kappa^-_{k,\infty}\coloneq \lim_{n\to\infty}\kappa^-_{k,n}$.
Define $\kappa^+_{\ell,\infty}$ in the same way from $U^{\mathsf T}$.
Let $\mathcal K^-_k$ and $\mathcal K^+_\ell$ denote the corresponding
Jacobi fractions of Proposition~\ref{prop:global-jacobi-fractions}.
Then
\begin{equation}\label{eq:global-two-sided-criterion}
 T\in\PBF
 \quad\text{if and only if}\quad
 \begin{cases}
 \mathcal K^-_k=\kappa^-_{k,\infty}>0,&k\in\{1,\dots,p-1\},\\
 \mathcal K^+_\ell=\kappa^+_{\ell,\infty}>0,&\ell\in\{1,\dots,q-1\}.
 \end{cases}
\end{equation}
The $p+q-2$ scalar fractions always converge to nonnegative limits;
the criterion requires all those limits to be strictly positive.
Equivalently, the adjacent Metelmann series
\eqref{eq:global-adjacent-series}, formed from $L$ and $U^{\mathsf T}$,
all converge. A side of bandwidth one contributes no conditions.
\end{corollary}
\begin{proof}
The Metelmann factors give compatible unit triangular products
$L$ and $U$, with positive diagonal $D$. Uniqueness of normalized
LDU makes these products the triangular products of every PBF of $T$.
On a finite leading section, suppose
$LDU=\widetilde L\widetilde D\widetilde U$, where the lower and
upper factors are unit triangular and both diagonal factors are
invertible. Then
$\widetilde L^{-1}L=\widetilde D\widetilde U U^{-1}D^{-1}$.
The left side is unit lower triangular and the right side is upper
triangular, so both are diagonal. The diagonal of the left side is
one; comparison on the right gives $\widetilde D=D$.
It follows that $\widetilde L=L$ and $\widetilde U=U$.
For the semi-infinite matrices, every fixed entry belongs to a
finite leading section, and triangular products and inverses
restrict to that section. The same uniqueness therefore holds
entrywise throughout the matrices.

A global PBF of $T$ consequently supplies a positive bidiagonal
factorization of $L$ and, after transposition and reversal of order,
one of $U^{\mathsf T}$. Applying
Theorem~\ref{thm:global-lower-criterion} to each product proves
necessity of all the stated limits. Conversely, if those limits
are positive, apply that theorem independently to $L$ and
$U^{\mathsf T}$. Transpose the second factorization and place the
given positive diagonal $D$ between the two products. Their product
is $T$, so this is the required global $(p,q)$ PBF.
The lower list has $p-1$ conditions and the upper list has $q-1$;
when a bandwidth equals one, that factor is already positive and
its list is empty. The convergence and the equivalent series
conditions follow from the identities in
Lemma~\ref{lem:global-inverse-minors} and
Proposition~\ref{prop:global-jacobi-fractions}.
\end{proof}

The global criterion decides whether some positive initial entries
work in every row. The next proposition decides whether a particular
positive factorization of the first section extends to the whole
matrix and gives formulas for its later entries. The second part of the
proposition determines when bounded positive factors exist.

Let $L^f$ be the leading section of $L$ of order $p+1$, and choose
any positive factorization $L^f=L_1^f\cdots L_p^f$.
For $r\in\{1,\dots,p-1\}$, form $X_r$ from rows $0,\dots,p-1$
and columns $0,\dots,p-r-1$ of
$\Sigma^f(L_1^f\cdots L_r^f)\Sigma^f$, where
$\Sigma^f=\diag[1,-1,\dots,(-1)^p]$.
This is the matrix used in the proof of
Theorem~\ref{thm:global-lower-criterion}. For
$m\in\{1,\dots,p-1\}$ and $n\ge0$, the determinants in
\eqref{eq:global-W-determinants} are
\[
 \mathcal W_{m,n}
 =\det\bigl(
   \mathcal Q[\mathcal I_{m,n},\{0,\dots,p-1\}]X_{p-m}
   \bigr).
\]
The selected block of $\mathcal Q=\Sigma L^{-1}\Sigma$ has size
$m\times p$, and $X_{p-m}$ has size $p\times m$, so their product
has size $m\times m$. Thus these determinants depend only on $L$
and the chosen finite factors.

\begin{proposition}[Global extensions and bounded factors]
\label{rem:global-construction-bounds}
\leavevmode
\begin{enumerate}
\item The chosen finite factors are the leading sections of a global
positive factorization of $L$ if and only if $\mathcal W_{m,n}>0$ for every
$m\in\{1,\dots,p-1\}$ and $n\ge1$.
Set $\mathcal W_{p,n}\coloneq \Delta_{p,n}$ and $\mathcal W_{0,n}\coloneq 1$.
For $m\in\{1,\dots,p\}$, factor $i=p-m+1$ has subdiagonal entries
\begin{equation}\label{eq:global-factor-reconstruction}
 a_i(n)=
 \frac{\mathcal W_{m,n-m+1}\mathcal W_{m-1,n-m}}
 {\mathcal W_{m,n-m}\mathcal W_{m-1,n-m+1}},\qquad n\ge m.
\end{equation}
Its first $m-1$ entries are the chosen initial data.
\item A semi-infinite BTP matrix $T=LDU$ admits a PBF with bounded
factors if and only if \eqref{eq:global-two-sided-criterion} holds
and the first subdiagonal of $L$, the first superdiagonal of $U$,
and the diagonal $D$ are bounded. If $T$ is bounded, then $D$ is
automatically bounded.
\end{enumerate}
\end{proposition}
\begin{proof}
For (1), suppose first that all the stated determinants are positive.
The induction in Theorem~\ref{thm:global-lower-criterion} applies
to this particular finite factorization: at each extraction its
test determinants are the corresponding $\mathcal W_{m,n}$.
Lemma~\ref{lem:global-extraction-determinant} therefore produces
positive entries in every row, and uniqueness of the recurrence
preserves the prescribed finite section. This proves sufficiency.

For necessity, suppose a global extension exists. Put
$F_r\coloneq L_{r+1}\cdots L_p$ for $0\le r\le p$, with $F_p=I$.
The identities for the first columns used in the same induction give
\[
 \mathcal W_{m,n}
 =\det\bigl((\Sigma F_{p-m}^{-1}\Sigma)
       [\mathcal I_{m,n},\mathcal C_m^{(0)}]\bigr),
 \qquad 1\le m\le p.
\]
This includes $m=p$, since $F_0=L$ and
$\mathcal W_{p,n}=\Delta_{p,n}$.
Fix $m$ and let $v_t$ be row $t$ of the first $m$ columns of
$\Sigma F_{p-m}^{-1}\Sigma$.
Multiplication by $\Sigma F_{p-m}\Sigma$ gives,
for $t\ge m$,
$v_t=\sum_{s=1}^{m}(-1)^{s+1}[F_{p-m}]_{t,t-s}v_{t-s}$.
In the determinant with rows $v_{n+1},\dots,v_{n+m}$,
substitute this expression for the last row. The terms $s<m$
repeat one of its other rows. Only the term with $s=m$ remains.
Moving its row $v_n$ to the first position requires $m-1$
exchanges, and the resulting sign $(-1)^{m-1}$ cancels the
coefficient sign $(-1)^{m+1}$. Hence
$\mathcal W_{m,n+1}=[F_{p-m}]_{n+m,n}\mathcal W_{m,n}$.
The initial determinant is one, because the signed inverse is
unit lower triangular. The extreme subdiagonal of $F_{p-m}$ is
positive: it is the product of one positive entry from each of
its $m$ bidiagonal factors. The recurrence therefore gives
$\mathcal W_{m,n}>0$ for all $n\ge0$, proving necessity.
This calculation also applies to $m=1$, when the determinant
has a single row and no row exchanges are needed.

To obtain the reconstruction formula, put $i\coloneq p-m+1$ and compare
the extreme subdiagonals in $F_{p-m}=L_iF_{p-m+1}$:
$[F_{p-m}]_{n,n-m}=a_i(n)[F_{p-m+1}]_{n-1,n-m}$ for $n\ge m$.
Taking $n-m$ as the starting row index in the determinant
recurrences for bandwidths $m$ and $m-1$ gives
\[
 [F_{p-m}]_{n,n-m}
 =\frac{\mathcal W_{m,n-m+1}}{\mathcal W_{m,n-m}},\qquad
 [F_{p-m+1}]_{n-1,n-m}
 =\frac{\mathcal W_{m-1,n-m+1}}{\mathcal W_{m-1,n-m}}.
\]
When $m=1$, the second entry is the diagonal entry of $F_p=I$,
and the second ratio equals one by the definition
$\mathcal W_{0,n}=1$. Every denominator in these expressions is
positive. Dividing the first ratio by the second proves
\eqref{eq:global-factor-reconstruction}. The extraction preserves
the first $m-1$ entries supplied by the finite factorization.

For (2), comparison on the first subdiagonal gives
$L_{n,n-1}=\sum_{i=1}^{p}a_i(n)$.
All summands are positive, so a uniform upper bound for
$L_{n,n-1}$ bounds every sequence $a_i(n)$. A unit lower
bidiagonal matrix with bounded entries $a_i(n)$ has operator norm
at most $1+\sup_n a_i(n)$, since its strictly lower part is a
weighted shift of norm $\sup_n a_i(n)$.
Conversely, boundedness of that matrix bounds each of its entries,
and the finite sum of these bounded sequences bounds $L_{n,n-1}$.
Transposition proves the corresponding equivalence on the upper
side. The diagonal factor is a bounded operator exactly when its
diagonal sequence is bounded. Combining these facts with the
global existence criterion proves (2).

Finally, if $T$ is bounded, its nonnegative normalized factors give
$T_{n,n}=\sum_{h=0}^{n}L_{n,h}d_hU_{h,n}\ge d_n>0$.
The term $h=n$ is $d_n$ because both triangular factors are unit.
Also $|T_{n,n}|=|\langle e_n,Te_n\rangle|\le\|T\|_2$.
Thus $D$ is bounded, as claimed.
\end{proof}

\begin{example}[A bounded matrix with an unbounded normalized factor]
\label{ex:bounded-matrix-unbounded-factor}
Boundedness of $T$ alone does not bound its normalized factors:
take $p=q=1$, $L\coloneq \mathcal L(a)$ with $a(n)\coloneq n$,
$d_n\coloneq (n+1)^{-2}$, and $U\coloneq I+\Lambda$.
The decomposition $T\coloneq LDU$ is a PBF, so $T$ is BTP by
Proposition~\ref{prop:positive-product}. It is bounded because,
for $n\ge1$, its three entries in row $n$ are
$1/n$, $1/n+(n+1)^{-2}$, and $(n+1)^{-2}$.
Its unique normalized lower factor $L$ is unbounded.
\end{example}

For lower bandwidth two, the full parameter condition reduces to
one interval, including its upper endpoint.
\begin{proposition}[Bandwidth two and the endpoint parameter]
\label{rem:bandwidth-two-endpoint}
Let $L$ be a compatible lower Metelmann product of bandwidth two,
and let $\kappa_{1,\infty}$ be the inverse-minor limit defined in
Lemma~\ref{lem:global-inverse-minors}(2). The first entries in positive
factorizations $L=\mathcal L(a)\mathcal L(b)$ are exactly
$a(1)\in(0,\kappa_{1,\infty}]$. Each such first entry determines
the factorization uniquely.
\end{proposition}
\begin{proof}
Suppose first that $u\in(0,\kappa_{1,\infty}]$ and put
$G_n\coloneq \mathcal Q_{n,0}-u\mathcal Q_{n,1}$, with $G_0=1$.
For every finite $n\ge1$, strict decrease gives
$\kappa_{1,n}>\kappa_{1,n+1}\ge\kappa_{1,\infty}\ge u$.
Since $\mathcal Q_{n,1}>0$, it follows that
$G_n=\mathcal Q_{n,1}(\kappa_{1,n}-u)>0$.
For $m=2$, the matrix of initial entries in
Lemma~\ref{lem:global-extraction-determinant} is the column
$[1,-u]^{\mathsf T}$, and its test determinant in row $n$
is exactly $G_n$. The lemma therefore constructs a positive
factor $\mathcal L(a)$ with $a(1)=u$ and a remaining unit
bidiagonal factor $\mathcal L(b)$ with positive subdiagonal.

The signed inverse of this remaining factor has first column
$(G_n)_{n\ge0}$, by
$\Sigma\mathcal L(b)^{-1}\Sigma=\mathcal Q(\Sigma\mathcal L(a)\Sigma)$.
The inverse formula for a bidiagonal matrix gives
$G_n=\prod_{h=1}^{n}b(h)$, so
$b(n)=G_n/G_{n-1}>0$.
Comparison on the second subdiagonal gives
$a(n)=L_{n,n-2}/b(n-1)>0$ for $n\ge2$.
The first-subdiagonal equation at row one fixes
$b(1)=L_{1,0}-u$. Thereafter the second-subdiagonal equation
fixes $a(n)$ and the first-subdiagonal equation fixes
$b(n)=L_{n,n-1}-a(n)$. Hence the chosen first entry determines
the whole factorization uniquely.

Conversely, any positive factorization satisfies
$0<a(1)<\kappa_{1,n}$ on every finite leading section by
Proposition~\ref{prop:inverse-parameter-bound}. The same positive
number $a(1)$ occurs in every section, so taking the infimum over
$n$ gives $a(1)\le\kappa_{1,\infty}$. In particular, if this
limit is zero, no positive first entry is possible and the stated
interval is empty. If the limit is positive, the first part of
the proof includes its upper endpoint.
\end{proof}

\noindent\emph{Connection with chain sequences.}
For $L=\mathcal L(a)\mathcal L(b)$, put $A_n\coloneq L_{n,n-1}$
for $n\ge1$ and $B_n\coloneq L_{n,n-2}$ for $n\ge2$.
Then $A_n=a(n)+b(n)$ and $B_{n+1}=a(n+1)b(n)$. Consequently
\[
 g_n\coloneq\frac{a(n+1)}{A_{n+1}}\quad(n\ge0),\qquad
 c_n\coloneq\frac{B_{n+1}}{A_nA_{n+1}}
       =(1-g_{n-1})g_n\quad(n\ge1).
\]
This is Wall's chain sequence with parameter sequence $(g_n)_{n\ge0}$
\cite[\S\S19--20]{Wall1948}, in the notation of
Chihara \cite[\S2]{Chihara1962}. Conversely, each parameter sequence
with $0<g_0<1$ and $0<g_n<1$ for $n\ge1$ recovers the positive
factors by $a(n+1)=A_{n+1}g_n$ and $b(n+1)=A_{n+1}(1-g_n)$.
The Metelmann choice $\widehat a_1(1)=0$ yields the minimal
parameters $m_n=\widehat a_1(n+1)/A_{n+1}$, with $m_0=0$.
Since $A_1=\widehat a_2(1)$, the series in
\eqref{eq:global-adjacent-series} takes the form
\[
 A_1\mathfrak S_1
 =1+\sum_{j=1}^{\infty}\prod_{i=1}^j\frac{m_i}{1-m_i}.
\]
Proposition~\ref{rem:bandwidth-two-endpoint} thus gives Wall's maximal
initial parameter $M_0=\kappa_{1,\infty}/A_1=1/(A_1\mathfrak S_1)$,
with $M_0=0$ when the series diverges. A positive factorization exists
if and only if $M_0>0$, and its normalized initial entry ranges over
$(0,M_0]$.

The criterion can be checked by bounding the tails of the adjacent
series. Their terms are determined directly by two consecutive
Metelmann factors: from \eqref{eq:global-adjacent-series},
\begin{equation}\label{eq:adjacent-term-ratio}
 \vartheta_k(n)\coloneq\frac{r_k(n+1)}{r_k(n)}
 =\frac{\widehat a_{p-k}(n+1)}{\widehat a_{p-k+1}(n+1)},
 \qquad n\ge k,\quad k\in\{1,\dots,p-1\}.
\end{equation}
For an integer $N\ge k$, write
$s_{k,N}\coloneq\sum_{n=k}^{N-1}r_k(n)$ for the known partial sum.
The next inequality bounds every remaining term through a telescoping sum.

\begin{proposition}[Certified bounds for the inverse-minor limits]
\label{prop:certified-limit-bounds}
Fix $k\in\{1,\dots,p-1\}$ and an integer $N>k$.
Suppose a positive sequence $\mu(n)$ satisfies
\begin{equation}\label{eq:tail-certificate}
 \mu(n)-\vartheta_k(n)\mu(n+1)\ge1\qquad(n\ge N).
\end{equation}
Then $\mathfrak S_k$ converges and
\begin{equation}\label{eq:certified-limit-interval}
 \frac{1}{s_{k,N}+\mu(N)r_k(N)}
 \le\kappa_{k,\infty}<\kappa_{k,N-k}=\frac{1}{s_{k,N}}.
\end{equation}
Convergence and the lower bound also hold for $N=k$, with
$s_{k,k}=0$.
\end{proposition}
\begin{proof}
Multiplying \eqref{eq:tail-certificate} by $r_k(n)>0$ gives
$r_k(n)\le\mu(n)r_k(n)-\mu(n+1)r_k(n+1)$.
On summing from $N$ to $M\ge N$, the intermediate terms cancel:
\[
 \sum_{n=N}^{M}r_k(n)
 \le\mu(N)r_k(N)-\mu(M+1)r_k(M+1)
 \le\mu(N)r_k(N).
\]
The partial sums increase and are bounded, so the tail converges.
Since every term is positive,
$s_{k,N}<\mathfrak S_k\le s_{k,N}+\mu(N)r_k(N)$.
Taking reciprocals gives the lower bound for $N\ge k$ and the
strict upper bound for $N>k$.
Lemma~\ref{lem:global-inverse-minors}(2), with $n=N-k$, identifies
$1/s_{k,N}$ with the convergent $\kappa_{k,N-k}$, giving
\eqref{eq:certified-limit-interval}.
\end{proof}

These bounds also give concrete initial parameters for a global PBF.
Apply Proposition~\ref{prop:certified-limit-bounds} to every adjacent
series and use its upper bounds for $\mathfrak S_k$ in
\eqref{eq:global-minor-bound-constant}; denote the resulting bound
for $\mathfrak B$ by $\overline{\mathfrak B}$.
Choose positive initial entries and construct the factors of the section
of order $p+1$ using \eqref{eq:row-system}. Check that its factors are positive
and compute the finite errors $\epsilon_m$ from
\eqref{eq:global-initial-minor-error}.
If $m!\overline{\mathfrak B}^{\,m}\epsilon_m<1$ for every
$m\in\{1,\dots,p-1\}$, then \eqref{eq:global-positive-margin}
and the extraction argument extend these factors positively to every row.
For example, choosing all the initial entries equal to $2^{-h}$ and
increasing the integer $h$ eventually passes these finite tests:
the finite construction is continuous at the Metelmann parameters
and every $\epsilon_m$ tends to zero there.

For rational ratios \eqref{eq:adjacent-term-ratio}, convergence is
decided by the leading terms at infinity.

\begin{corollary}[A criterion for rational Metelmann ratios]
\label{cor:rational-metelmann-ratios}
Let $L=\widehat L_1\cdots\widehat L_p$ be the semi-infinite lower
Metelmann product in \eqref{eq:metelmann}. Suppose that, for each
$k\in\{1,\dots,p-1\}$, the ratio $\vartheta_k(n)$ in
\eqref{eq:adjacent-term-ratio} agrees with a rational function for
all sufficiently large integers $n$. Put
$\lambda_k\coloneq\lim_{n\to\infty}\vartheta_k(n)\in[0,+\infty]$.
When $\lambda_k=1$, define $\sigma_k$ by the expansion
\begin{equation}\label{eq:rational-ratio-asymptotics}
 \vartheta_k(n)=1-\frac{\sigma_k}{n}
                  +\mathrm O(n^{-2}),\qquad n\to\infty.
\end{equation}
Then $L$ admits a global positive factorization if and only if,
for every $k\in\{1,\dots,p-1\}$, either $\lambda_k<1$, or
$\lambda_k=1$ and $\sigma_k>1$.
\end{corollary}
\begin{proof}
Fix $k$. All the terms $r_k(n)$ are positive, and removing
finitely many of them does not affect convergence.
If $\lambda_k<1$, choose $\theta$ with
$\lambda_k<\theta<1$. For some $N$,
$\vartheta_k(n)\le\theta$ whenever $n\ge N$.
Iteration of $r_k(n+1)=\vartheta_k(n)r_k(n)$ then gives
$r_k(N+j)\le r_k(N)\theta^j$ for every $j\ge0$.
The tail is bounded by a convergent geometric series.
If $\lambda_k>1$, including $\lambda_k=+\infty$, choose
$c>1$ with $\vartheta_k(n)\ge c$ for all sufficiently large
$n$. Then $r_k(N+j)\ge r_k(N)c^j$, so the terms do not tend
to zero and the series diverges.

It remains to consider $\lambda_k=1$.
The rational function is defined and positive at every sufficiently
large integer. Its expansion at infinity is
\eqref{eq:rational-ratio-asymptotics}, so taking real logarithms gives
$\log\vartheta_k(n)=-\sigma_k/n+\mathrm O(n^{-2})$ as $n\to\infty$.
Set $b_n\coloneq \log r_k(n)+\sigma_k\log n$ for sufficiently large
positive $n$. Since
$\log(1+1/n)=1/n+\mathrm O(n^{-2})$ as $n\to\infty$,
$b_{n+1}-b_n=\log\vartheta_k(n)+\sigma_k\log(1+1/n)=\mathrm O(n^{-2})$ as $n\to\infty$.
Thus $\sum_n|b_{n+1}-b_n|$ converges. Writing $b_n$ as its
initial value plus these consecutive differences shows that it
has a finite real limit $b_\infty$.
Exponentiation gives
$n^{\sigma_k}r_k(n)\to\mathrm e^{b_\infty}>0$ as $n\to\infty$.
For all sufficiently large $n$, the term $r_k(n)$ is therefore
bounded above and below by positive constant multiples of
$n^{-\sigma_k}$. The positive series $\mathfrak S_k$ converges
exactly when $\sigma_k>1$, including divergence at $\sigma_k=1$.
Applying this test to each adjacent pair and then using
Theorem~\ref{thm:global-lower-criterion} proves the equivalence.
\end{proof}

The same calculation applies to $U^{\mathsf T}$ and gives the
$p+q-2$ tests in Corollary~\ref{cor:global-two-sided-criterion}.
Exact tail bounds can be found from polynomials $P_k,Q_k$ with
rational coefficients and an index $N_0$ such that
$\vartheta_k(n)=P_k(n)/Q_k(n)$ for every integer $n\ge N_0$.
Choose their common sign so that $Q_k$ has positive leading
coefficient. A polynomial with positive leading coefficient has
positive coefficients in powers of $n-N$ for all sufficiently
large $N$: those coefficients are its derivatives at $N$, divided
by the corresponding factorials. Testing these finitely many
rational coefficients while increasing the integer
$N\ge\max(N_0,k+1)$ therefore gives a finite procedure.

If $\lambda_k<1$, choose a rational $\theta$ with
$\lambda_k<\theta<1$. Both $Q_k$ and $\theta Q_k-P_k$ have
positive leading coefficients. Apply the coefficient test to
these two polynomials. For the resulting $N$, they are positive
on $[N,\infty)$, so $\vartheta_k(n)\le\theta$ there and
$\mu(n)\coloneq 1/(1-\theta)$ satisfies \eqref{eq:tail-certificate}.
If $\lambda_k=1$ and $\sigma_k>1$, take
$\mu(n)\coloneq 2n/(\sigma_k-1)$. Formula
\eqref{eq:rational-ratio-asymptotics} gives
$\mu(n)-\vartheta_k(n)\mu(n+1)\to2$ as $n\to\infty$.
The polynomial
$Q_k(n)(\mu(n)-1)-P_k(n)\mu(n+1)$ consequently has positive
leading coefficient. Apply the same test to this polynomial and
$Q_k$ to obtain an $N\ge\max(N_0,k+1)$ for which
\eqref{eq:tail-certificate} holds for every integer $n\ge N$.

For eventually constant Metelmann parameters, the ratios in
\eqref{eq:adjacent-term-ratio} are eventually constant as well.
Corollary~\ref{cor:rational-metelmann-ratios} then gives a strict
ordering of the constants. When the entries are constant from
their first permitted index, the factors themselves have explicit formulas.

\begin{proposition}[Eventually constant Metelmann parameters]
\label{prop:constant-metelmann-criterion}
Let $L=\widehat L_1\cdots\widehat L_p$ be the semi-infinite lower
Metelmann product in \eqref{eq:metelmann}. Suppose the positive
entries of its factors, $\widehat a_i(n)=[\widehat L_i]_{n,n-1}$,
are eventually constant, with value $\xi_i>0$, for $i\in\{1,\dots,p\}$. Then
\begin{equation}\label{eq:constant-metelmann-order}
 L\text{ admits a global positive factorization}
 \quad\text{if and only if}\quad \xi_1<\xi_2<\cdots<\xi_p.
\end{equation}
If each entry is already constant from its first permitted index,
then, for $k\in\{1,\dots,p-1\}$ and $n\ge1$, the inverse-minor
ratios \eqref{eq:global-minor-ratios} satisfy, with
$x\coloneq \xi_{p-k}$ and $y\coloneq \xi_{p-k+1}$,
\begin{equation}\label{eq:constant-metelmann-limits}
 \kappa_{k,n}=
 \frac{y^n}{\sum_{j=0}^{n-1}x^j y^{n-1-j}},\qquad
 \kappa_{k,\infty}=\max(y-x,0).
\end{equation}
In this latter case, writing $\mathcal L(b_j)$ for the unit lower
bidiagonal matrix with subdiagonal entries $b_j(n)$, an explicit
factorization is
\begin{equation}\label{eq:constant-metelmann-explicit}
 L=\mathcal L(b_p)\cdots\mathcal L(b_1),\qquad
 b_j(n)\coloneq 
 \begin{cases}
 \xi_j-\xi_{j-n},&n\in\{1,\dots,j-1\},\\
 \xi_j,&n\ge j.
 \end{cases}
\end{equation}
Its parameters are positive when the strict order holds.
\end{proposition}
\begin{proof}
For each $k$, the eventual ratio of consecutive series terms is
$\vartheta_k(n)=\xi_{p-k}/\xi_{p-k+1}$.
If this ratio is one, the expansion
\eqref{eq:rational-ratio-asymptotics} has $\sigma_k=0$.
Corollary~\ref{cor:rational-metelmann-ratios} therefore says that
all the series converge exactly when
$\xi_{p-k}<\xi_{p-k+1}$ for every $1\le k\le p-1$.
These are precisely the adjacent inequalities in
\eqref{eq:constant-metelmann-order}.

When the entries are constant from their first permitted row,
$r_k(s)=y^{-1}(x/y)^{s-k}$ for $s\ge k$.
For $n\ge1$, summing the first $n$ terms gives
\[
 \kappa_{k,n}^{-1}
 =\frac1y\sum_{j=0}^{n-1}(x/y)^j,
 \qquad
 \kappa_{k,n}
 =\frac{y^n}{\sum_{j=0}^{n-1}x^j y^{n-1-j}}.
\]
If $x<y$, the infinite geometric sum is $1/(y-x)$ and its
reciprocal is $y-x$. If $x\ge y$, the partial sums diverge
and the reciprocal tends to zero. For $x=y$, the finite formula
is $y/n$. This proves \eqref{eq:constant-metelmann-limits} in
all cases.

The explicit product can be verified entry by entry.
Write $L^{[p]}$ for the Metelmann product with parameters
$\xi_1,\dots,\xi_p$. For a list of variables $X$, denote by
$e_r(X)$ the sum of products selected by all $r$-element subsets
of the positions in that list, with $e_0(X)\coloneq 1$ and $e_r(X)\coloneq 0$ if $r<0$ or if $r$
exceeds its length. These are the elementary symmetric polynomials.
In entry $(n,n-r)$ of a product of $p$ lower bidiagonal matrices,
a term selects subdiagonal entries from exactly $r$ factors,
say $i_1<\cdots<i_r$, and diagonal entries from the others.
Its value in the Metelmann product is
$\prod_{s=1}^{r}\widehat a_{i_s}(n-s+1)$.
For $r=0$, this is the empty product $1$ on the diagonal.
For $r\ge1$, the term is nonzero exactly when
$n-s+1>p-i_s$ for every $s$.
For $n<p$, these inequalities select all the $r$-element subsets
of $\{p-n+1,\dots,p\}$: an increasing subset of that set has
$i_s\ge p-n+s$, and the converse follows already from its
first index. For $n\ge p$, every subset is permitted.
Consequently,
$[L^{[p]}]_{n,n-r}=e_r(\xi_{p-\min(n,p)+1},\dots,\xi_p)$ for $0\le r\le\min(n,p)$.

For $1\le n<p$, put $X\coloneq (\xi_{p-n+1},\dots,\xi_{p-1})$.
In the Metelmann product $L^{[p-1]}$ formed from
$\xi_1,\dots,\xi_{p-1}$, row $n$ uses the variables
$(X,\xi_{p-n})$, whereas row $n-1$ uses $X$.
The bidiagonal multiplication in
$\mathcal L(b_p)L^{[p-1]}$ therefore gives
\[
 \begin{aligned}
 [\mathcal L(b_p)L^{[p-1]}]_{n,n-r}
 &=e_r(X,\xi_{p-n})+(\xi_p-\xi_{p-n})e_{r-1}(X)\\
 &=e_r(X)+\xi_p e_{r-1}(X)
 =e_r(X,\xi_p).
 \end{aligned}
\]
Here the identity $e_r(X,t)=e_r(X)+t e_{r-1}(X)$ follows by
separating the subsets that contain $t$ from those that do not.
For $n\ge p$, rows $n$ and $n-1$ of $L^{[p-1]}$ both use all
$p-1$ variables, and $b_p(n)=\xi_p$. The same identity then
gives the row of $L^{[p]}$. Row zero is the common unit diagonal
row, and entries outside the bandwidth are zero in both products.
Thus $L^{[p]}=\mathcal L(b_p)L^{[p-1]}$.
Starting with the single factor for $p=1$ and repeating this identity
proves \eqref{eq:constant-metelmann-explicit}. Under the strict
ordering, every difference $\xi_j-\xi_{j-n}$ that occurs is
positive; all the remaining entries equal the positive number
$\xi_j$. Hence the displayed factorization is positive.
\end{proof}

When $p\ge2$ and $\xi_1<\cdots<\xi_p$, the first entry in
\eqref{eq:constant-metelmann-explicit} is
$[\mathcal L(b_p)]_{1,0}=\xi_p-\xi_{p-1}=\kappa_{1,\infty}$.
Thus this factorization attains the limiting bound obtained from
Proposition~\ref{prop:inverse-parameter-bound} on the leading
sections. For $p=2$, it gives the endpoint in
Proposition~\ref{rem:bandwidth-two-endpoint}.

Taking every permitted Metelmann entry equal to $1$ gives, for each
$p\ge2$, a bounded BTP matrix
$T_p\coloneq (\widehat L_1\cdots\widehat L_p)(I+\Lambda)$ whose finite
leading sections admit PBFs and whose global factorization is
obstructed by $\xi_1=\cdots=\xi_p=1$ in
\eqref{eq:constant-metelmann-order}; at $p=2$ this is $T_*$.

\begin{example}[Independent obstructions in bandwidth four]
\label{ex:bandwidth-four-obstructions}
For four Metelmann factors whose permitted entries are constant,
the following triples of limits follow from
\eqref{eq:constant-metelmann-limits}:
\[
 \begin{NiceArray}{c|ccc}
 (\xi_1,\xi_2,\xi_3,\xi_4)&
 \kappa_{1,\infty}&\kappa_{2,\infty}&\kappa_{3,\infty}\\ \hline
 (2,1,3,4)&1&2&0\\
 (1,3,2,4)&2&0&2\\
 (1,2,4,3)&0&2&1
 \end{NiceArray}
\]
Each condition can therefore fail while the other two hold.
For the first example,
$\kappa_{3,n}=1/(2^n-1)\to0$ as $n\to\infty$.
In each case $T\coloneq L(I+\Lambda)$, where $\Lambda$ is the upper shift
introduced in Section~\ref{sec:limit}, is bounded and BTP of
bandwidth $(4,1)$ by applying Proposition~\ref{prop:metelmann-input}
to its leading sections, and has
no global PBF by Corollary~\ref{cor:global-two-sided-criterion}.

\end{example}

\begin{example}[Equal limiting parameters and certified positive factors]
\label{ex:equal-tails-certification}
For any $p\ge2$, take the Metelmann parameters
$\widehat a_i(n)\coloneq (n+2i-p-1)/n$ for $n\ge p-i+1$ and zero
at the preceding indices, with $i\in\{1,\dots,p\}$.
All their positive tails tend to $1$ as $n\to\infty$.
The adjacent terms and their tails are
\begin{equation}\label{eq:equal-tails-exact-series}
 \begin{aligned}
 r_k(n)&=\frac{k(p-k)}{(n+p-2k)(n+p-2k+1)},\qquad n\ge k,\\
 \sum_{n=N}^{\infty}r_k(n)&=\frac{k(p-k)}{N+p-2k},\qquad N\ge k.
 \end{aligned}
\end{equation}
Indeed, the first formula has initial value $k/(p-k+1)$ and
successive ratio $(n+p-2k)/(n+p-2k+2)$, as required by
\eqref{eq:global-adjacent-series} and \eqref{eq:adjacent-term-ratio}.
The second follows by writing each term as the difference
$k(p-k)((n+p-2k)^{-1}-(n+p-2k+1)^{-1})$.
Taking $N=k$ in \eqref{eq:equal-tails-exact-series} gives
$\mathfrak S_k=k$ and $\kappa_{k,\infty}=1/k$ for
$k\in\{1,\dots,p-1\}$, and a global PBF exists.
Applying the same construction on the upper side gives examples
for every pair $(p,q)$.

For $p=4$, these sums give $\mathfrak B=6$ in
\eqref{eq:global-minor-bound-constant}.
Set all six free initial entries equal to $u\coloneq 1/4096$.
The finite construction \eqref{eq:row-system} gives positive
factors through row $4$. The errors
\eqref{eq:global-initial-minor-error} are
$\epsilon_1=3u+3u^2+u^3$,
$\epsilon_2=2u+4u^2+2u^3+u^4$, and
$\epsilon_3=u+u^2+u^3$.
Direct substitution gives $m!6^m\epsilon_m<1/3$ for
$m\in\{1,2,3\}$. Therefore
\eqref{eq:global-positive-margin} yields
$\mathcal W_{m,n}/\Delta_{m,n}>2/3$ for every $n\ge0$ and
$m\in\{1,2,3\}$, certifying the global positive extension of
these initial entries.
\end{example}

\subsection{Approximation through an initial block}
\label{subsec:boundary-approximation}

Corollary~\ref{cor:global-two-sided-criterion} decides whether the
original matrix has a global PBF. Filling the initial zeros of its
Metelmann factors in \eqref{eq:metelmann} produces a family of
matrices with PBFs for arbitrary finite bandwidths $(p,q)$.
The next theorem gives its finite support, rank, monotonicity,
and operator-norm bounds.

For $r\in\mathbb N$, let $\iota_r:\CC^r\to\ell^2(\NN)$ be the
inclusion in coordinates $\{0,\dots,r-1\}$, and let
$P_r\coloneq \iota_r\iota_r^{\mathsf T}$ be the corresponding orthogonal
projection. The transpose notation is valid because $\iota_r$ has
real entries.

\begin{theorem}[Approximation by perturbing an initial block]
\label{thm:boundary-approximation}
Let $T$ be a bounded semi-infinite BTP matrix of bandwidth $(p,q)$,
where $p,q\in\mathbb N$, and put $m_0\coloneq p+q-1$.
There is a family of bounded matrices $T_\varepsilon$, for
$\varepsilon>0$, of the same bandwidth and admitting a PBF, such that
$\Delta_\varepsilon\coloneq T_\varepsilon-T$ satisfies
\begin{equation}\label{eq:boundary-support-rank}
 \Delta_\varepsilon=P_{m_0}\Delta_\varepsilon P_{m_0},\qquad
 \operatorname{rank}\Delta_\varepsilon\le p+q-2.
\end{equation}
The family is monotone entrywise:
\begin{equation}\label{eq:boundary-order}
 T\le T_\varepsilon\le T_\delta
 \qquad(0<\varepsilon\le\delta).
\end{equation}
Moreover,
\begin{equation}\label{eq:boundary-norm-bound}
 \|T_\varepsilon-T\|_2
 \le\varepsilon\|T_1-T\|_2
 \qquad(0<\varepsilon\le1).
\end{equation}
In particular, $T_\varepsilon\to T$ in operator norm as
$\varepsilon\downarrow0$, and every entry outside the initial
$m_0\times m_0$ block is unchanged.
Here $T_\varepsilon$ is a bounded operator on $\ell^2(\NN)$, and
its PBF is the entrywise product specified in \eqref{eq:strict-definition}.
Its factors are bounded whenever the Metelmann factors of $T$ and
the positive diagonal factor $D$ in \eqref{eq:metelmann} are bounded.
\end{theorem}
\begin{proof}
Use the compatible semi-infinite Metelmann factorization
\eqref{eq:metelmann}. For $1\le i\le p$ and $1\le j\le q$, put
\[
 Z_i^-\coloneq \sum_{n=1}^{p-i}e_ne_{n-1}^{\mathsf T},\qquad
 Z_j^+\coloneq \sum_{n=1}^{q-j}e_{n-1}e_n^{\mathsf T},
\]
where an empty sum is zero, and define
\begin{equation}\label{eq:boundary-filled-factors}
 L_i(\varepsilon)\coloneq \widehat L_i+\varepsilon Z_i^-,\qquad
 U_j(\varepsilon)\coloneq \widehat U_j+\varepsilon Z_j^+.
\end{equation}
Set $L_\varepsilon\coloneq L_1(\varepsilon)\cdots L_p(\varepsilon)$,
$U_\varepsilon\coloneq U_q(\varepsilon)\cdots U_1(\varepsilon)$,
and $T_\varepsilon\coloneq L_\varepsilon D U_\varepsilon$.
At $\varepsilon=0$ these products recover $L,U,T$.
Every initial zero subdiagonal or superdiagonal entry is replaced
by $\varepsilon$, and every previously positive allowed entry is
unchanged. Thus the entrywise product is a PBF for
$\varepsilon>0$. There are finitely many banded factors, so every
entry is computed by a finite sum.

\emph{Support of the lower and upper changes.}
Expand $L_\varepsilon-L$ as a sum of products containing at
least one changed factor entry. Consider one such entry in
$L_i$, in position $(n,n-1)$ with $n\le p-i$.
In any term contributing to an entry $(a,b)$ of the full product,
the $i-1$ factors to the left of $L_i$ can increase the row
index from $n$ by at most $i-1$. Therefore
$a\le n+i-1\le p-1$. The chosen subdiagonal entry forces
$a>b$, and the remaining factors are lower triangular.
Every contribution to the difference consequently has
$1\le a\le p-1$ and $0\le b<a$.
This reasoning applies even when a term contains several changed
entries, since any one of them gives the stated bound.
Thus $L_\varepsilon-L$ has finite support, all its nonzero rows
lie in $\{1,\dots,p-1\}$, all its nonzero columns lie in
$\{0,\dots,p-2\}$, and its rank is at most $p-1$.
Apply the same argument to
$U_\varepsilon^{\mathsf T}
 =U_1(\varepsilon)^{\mathsf T}\cdots U_q(\varepsilon)^{\mathsf T}$.
After transposition, a nonzero entry $(a,b)$ of
$U_\varepsilon-U$ satisfies
$1\le b\le q-1$ and $0\le a<b$; its rank is at most $q-1$.
The respective difference is zero when $p=1$ or $q=1$.

\emph{Support and rank of the full change.}
Subtract the products using the exact identity
\begin{equation}\label{eq:boundary-product-difference}
 \Delta_\varepsilon
 =(L_\varepsilon-L)D U_\varepsilon
       +LD(U_\varepsilon-U).
\end{equation}
In the first term, the left factor restricts the row to at most
$p-1$ and the intermediate diagonal index to at most $p-2$.
The upper bandwidth $q$ then restricts the final column to at
most $p+q-2$. In the second term, the right factor restricts
the column to at most $q-1$ and the intermediate diagonal
index to at most $q-2$; the lower bandwidth $p$ restricts the
final row to at most $p+q-2$. Each nonzero term is therefore
supported in rows and columns $0,\dots,p+q-2$, proving
$\Delta_\varepsilon=P_{m_0}\Delta_\varepsilon P_{m_0}$.
When a side has bandwidth one, its zero difference contributes
no term, so these support bounds still apply.

All multiplications involving either difference reduce to
finite matrices: only finitely many entries of $D$ and of the
banded factors are used. The finite-dimensional rank inequalities
therefore bound the two summands in
\eqref{eq:boundary-product-difference} by $p-1$ and $q-1$,
respectively. Their sum has rank at most $p+q-2$.
Every finitely supported matrix is a bounded operator on
$\ell^2$, since it acts through a finite-dimensional coordinate
space. Hence $\Delta_\varepsilon$ is bounded, and so is
$T_\varepsilon=T+\Delta_\varepsilon$.

\emph{Order and norm estimate.}
Expand the full product as a polynomial in $\varepsilon$.
Each coefficient is a sum of products of matrices with
nonnegative entries. Removing its constant term gives
$\Delta_\varepsilon=\sum_{k=1}^{p+q-2}\varepsilon^k\Delta^{[k]}$, with $\Delta^{[k]}\ge0$ entrywise.
There are at most $(p-1)+(q-1)$ factors with a nonzero change,
because $Z_p^-=Z_q^+=0$; this proves the degree bound.
The support argument for each term also shows that every
coefficient matrix has finite support. If $p=q=1$, this sum is
empty and $T_\varepsilon=T$.
For $0<\varepsilon\le\delta$, each scalar power satisfies
$\varepsilon^k\le\delta^k$, proving \eqref{eq:boundary-order}.
The product has bandwidth at most $(p,q)$, and the entrywise
inequality $T_\varepsilon\ge T$ preserves the positive extreme
band entries of $T$. Thus its bandwidth is the same as that of $T$.

For $0<\varepsilon\le1$,
$\varepsilon^k\le\varepsilon$ whenever $k\ge1$, so
$0\le\Delta_\varepsilon\le\varepsilon\Delta_1$ entrywise.
For a vector $x\in\ell^2(\NN)$, let $|x|$ be its sequence of
componentwise absolute values. The triangle inequality in each
coordinate gives
\[
 |\Delta_\varepsilon x|\le\varepsilon\Delta_1|x|,
 \qquad
 \|\Delta_\varepsilon x\|_{\ell^2}
 \le\varepsilon\|\Delta_1\|_2\,\||x|\|_{\ell^2}
 =\varepsilon\|\Delta_1\|_2\,\|x\|_{\ell^2}.
\]
Taking the supremum over unit vectors proves
\eqref{eq:boundary-norm-bound} and the asserted norm convergence.
Finally, each $Z_i^-$ and $Z_j^+$ is finitely supported and bounded.
If the original Metelmann factors and $D$ are bounded, then adding
these bounded corrections leaves all factors of the displayed
PBF bounded, as asserted.
\end{proof}

\begin{remark}[The tetradiagonal initial-row perturbation]
\label{rem:two-norm-approximations}
In the construction \eqref{eq:boundary-filled-factors}, for
$(p,q)=(2,1)$ only the entry $(1,0)$ of $\widehat L_1$ is changed. Both $\widehat L_2$ and $L$ are unit lower triangular,
so their row zero is $e_0^{\mathsf T}$. Consequently,
$e_1e_0^{\mathsf T}\widehat L_2=e_1e_0^{\mathsf T}L=e_1e_0^{\mathsf T}$,
and hence
\[
 \begin{aligned}
 L_\varepsilon
 &=(\widehat L_1+\varepsilon e_1e_0^{\mathsf T})\widehat L_2=L+\varepsilon e_1e_0^{\mathsf T}
   =(I+\varepsilon e_1e_0^{\mathsf T})L,\\
 T_\varepsilon&=(I+\varepsilon e_1e_0^{\mathsf T})T.
 \end{aligned}
\]
This is the initial-row perturbation used in \cite{BFMTetra}.
Theorem~\ref{thm:boundary-approximation} extends the construction
to arbitrary lower and upper bandwidths, with explicit support and
rank bounds, entrywise monotonicity, and an operator-norm estimate.

For $T_*$ in \eqref{eq:counter-matrix},
$T_{*,\varepsilon}-T_*=\varepsilon e_1(e_0+e_1)^{\mathsf T}$, so
$\|T_{*,\varepsilon}-T_*\|_2=\sqrt2\,\varepsilon$.
The perturbation is confined to the initial $2\times2$ block.
\end{remark}

\section{Stability of positive bidiagonal factorizations}
\label{chap:stability}

The criterion in Section~\ref{sec:finite} determines whether a matrix
admits some global PBF. Stability concerns a prescribed factorization:
when the matrix changes, can nearby positive factors be chosen
analytically through those prescribed factors? The topology matters.
Positivity of the existence limits alone does not control changes
uniformly over all rows, as Theorem~\ref{stab:thm:distant} shows.

For the bounded problem, fix a bounded real matrix $T$ of bandwidth $(p,q)$
and one PBF. Perturbations $\widetilde T-T$ are measured in the
operator norm $\|\cdot\|_2$ on $\ell^2(\NN)$ and remain in the same
band. Factor sequences are measured by
$\|a\|_\infty\coloneq\sup_{n\ge1}|a(n)|$, with the maximum norm
for a finite tuple of sequences. Write $\LL(a)$ for the unit lower
bidiagonal matrix with entry $a(n)$ at $(n,n-1)$; for a constant
$c$, $\LL(c)=I+c\Lambda^{\mathsf T}$. A real analytic factor map
has convergent power-series expansions in these norms and equals
the prescribed factors at $T$. It is locally Lipschitz, with one
bound for the changes in all factor entries. Its positivity also
implies BTP by Proposition~\ref{prop:positive-product}.

\begin{example}[Perturbing a product of two factors]
\label{stab:ex:reconstruction-rows}
For $L=\LL(a_1)\LL(a_2)$, a nearby unit lower matrix $\widetilde L$
of bandwidth at most two has candidate factors determined by
\[
 \widetilde a_1(1)=u,\qquad
 \widetilde a_2(1)=\widetilde L_{1,0}-u,
\]
and, as long as the denominators are positive, by
\[
 \widetilde a_1(n)=\frac{\widetilde L_{n,n-2}}{\widetilde a_2(n-1)},
 \qquad
 \widetilde a_2(n)=\widetilde L_{n,n-1}-\widetilde a_1(n),
 \qquad n\ge2.
\]
These are the equations obtained by matching the two subdiagonals.
Keeping the initial entry fixed means taking $u=a_1(1)$ for every
$\widetilde L$. Allowing it to vary means choosing $u=u(\widetilde L)$,
with $u(L)=a_1(1)$. Either choice must keep the resulting sequences
positive and uniformly close to the original factors.
\end{example}

For $p$ lower factors the free initial entries are $a_i(n)$ with
$i\in\{1,\dots,p-1\}$ and $n\in\{1,\dots,p-i\}$;
Lemma~\ref{lem:rowwise} determines all later entries. The upper
product has the corresponding initial entries. Subsections
\ref{stab:subsec:fixed} and \ref{stab:subsec:variable} characterize
analytic dependence with these entries fixed and variable,
respectively. The proofs combine first-order recurrences, their
geometric product bounds, and the analytic inverse function theorem
\cite{BravermanKarabash,HuyMinh2001,Deimling}.

Subsection~\ref{stab:subsec:consequences} then connects these
criteria with interiority and describes all sufficiently close
factorizations of a fixed matrix. The two-factor example explains
why nearby initial entries need not give uniformly nearby factor
sequences. These results determine the local geometry used later
for Toeplitz normalization.

Finally, Subsection~\ref{stab:subsec:relative} extends the
fixed-entry criterion to arbitrary positive factor growth by
measuring changes relative to the original entries. It recovers
the bounded criterion when the relative and absolute norms are
equivalent. The products in this extension are defined entrywise;
no common-domain realization of unbounded operators is assumed.

\subsection{Continuity of the existence limits and loss of PBF}
The endpoints in Corollary~\ref{cor:global-two-sided-criterion} are
limits of the decreasing finite ratios
\eqref{eq:global-minor-ratios}. Thus
$\kappa_{k,\infty}=\inf_{n\ge1}\kappa_{k,n}$ on either triangular
side. Each finite ratio depends on only finitely many matrix entries.
This observation gives the following continuity property.

\begin{proposition}[Upper semicontinuity of the endpoints]\label{stab:prop:semicontinuity}
Fix finite bandwidths $(p,q)$. If BTP matrices $T_j$ of bandwidth
$(p,q)$ converge entrywise to a BTP matrix $T$ of this bandwidth, then
\[
 \limsup_{j\to\infty}\kappa^-_{k,\infty}(T_j)
 \le\kappa^-_{k,\infty}(T),\qquad
 \limsup_{j\to\infty}\kappa^+_{\ell,\infty}(T_j)
 \le\kappa^+_{\ell,\infty}(T)
\]
for $k\in\{1,\dots,p-1\}$ and $\ell\in\{1,\dots,q-1\}$.
In other words, these endpoints are upper semicontinuous in the
topology of entrywise convergence. If all endpoints of all $T_j$
are bounded below by one number $c>0$, then $T$ has a global PBF.
\end{proposition}
\begin{proof}
Fix one side, one order $k$, and one finite row index $n$.
The ratio $\kappa_{k,n}$ uses finitely many entries of a triangular
inverse. The compatible normalized LDU decomposition computes these
entries from a sufficiently large leading section of $T$.
In that section Gaussian elimination uses nonzero pivots, and the
triangular inverse is computed by finitely many additions,
multiplications and divisions by those pivots. The denominator minor
in the ratio is positive by Lemma~\ref{lem:global-inverse-minors}.
All these denominators remain nonzero near the given finite section.
Consequently $\kappa_{k,n}(T_j)\to\kappa_{k,n}(T)$ as
$j\to\infty$.

For each $j$,
$\kappa_{k,\infty}(T_j)\le\kappa_{k,n}(T_j)$.
Taking the upper limit in $j$ gives
\[
 \limsup_{j\to\infty}\kappa_{k,\infty}(T_j)
 \le\kappa_{k,n}(T).
\]
This holds for every $n\ge1$. Taking the infimum over $n$ on the
right gives the asserted inequality. The argument applies to $L$
and to $U^{\mathsf T}$, and hence to every lower and upper endpoint.
If the endpoints of $T_j$ are at least $c$, the inequality gives
$\kappa^-_{k,\infty}(T)\ge c$ and
$\kappa^+_{\ell,\infty}(T)\ge c$.
Corollary~\ref{cor:global-two-sided-criterion} now gives a global
PBF of $T$. A side of bandwidth one contributes no endpoints.
\end{proof}

In the next theorem the limiting matrix has positive endpoints,
whereas every approximating matrix has a zero endpoint. The change
is confined to one entry whose row index tends to infinity.

\begin{theorem}[Loss of PBF under a small distant perturbation]\label{stab:thm:distant}
Fix $c>0$, set $L_0\coloneq(I+c\Lambda^{\mathsf T})^2$, and, for $N\ge1$, put
\begin{equation}\label{stab:eq:distant}
 L_N\coloneq L_0+\frac{c^2}{N}e_{N+1}e_{N-1}^{\mathsf T}.
\end{equation}
Each $L_N$ is a compatible lower Metelmann product of bandwidth two,
but has no factorization into two unit lower bidiagonal matrices with
all subdiagonal entries strictly positive. Moreover,
$\|L_N-L_0\|_2=c^2/N$, whereas
$\kappa_{1,\infty}(L_0)=c$ and $\kappa_{1,\infty}(L_N)=0$.
For any $d>0$, the matrices $T_N\coloneq L_N(I+d\Lambda)$ are BTP of
bandwidth $(2,1)$ and satisfy
\begin{equation}\label{stab:eq:distant-full}
 \|T_N-T_0\|_2=\frac{c^2\sqrt{1+d^2}}{N}.
\end{equation}
Here $T_0$ has a bounded PBF and none of the $T_N$, $N\ge1$, has a
PBF. If $0<c,d<1$, $T_0$ and all sufficiently large $T_N$ are
invertible bounded operators; every displayed factor of $T_0$ has
a bounded inverse.
\end{theorem}
\begin{proof}
The Metelmann factors of $L_0$ have subdiagonal entries
$\widehat a(n)=c(n-1)/n$ and $\widehat b(n)=c(n+1)/n$.
Indeed, $\widehat a(1)=0$, all the other permitted entries are
positive, and multiplication of the two bidiagonal matrices gives $\widehat a(n)+\widehat b(n)=2c$ for $n\ge1$ and $\widehat a(n)\widehat b(n-1)=c^2$ for $n\ge2$.
These are exactly the first and second subdiagonals of $L_0$.

For $L_N$, retain these parameters through row $N$, and set both
parameters equal to $c$ from row $N+1$ onward. Their sum remains
$2c$ in every row. The second subdiagonal is unchanged through row
$N$, and in row $N+1$ it becomes
$c\widehat b(N)=c^2(N+1)/N$. From row $N+2$ onward it is again
$c^2$. This proves \eqref{stab:eq:distant} and supplies the required
Metelmann representation, including at $N=1$.

For two Metelmann factors, the terms in
\eqref{eq:global-adjacent-series} are
\[
 r_1(n)=\frac1{\widehat b(1)}
       \prod_{t=2}^{n}\frac{\widehat a(t)}{\widehat b(t)},
 \qquad n\ge1.
\]
For $L_0$, the product of $(t-1)/(t+1)$ telescopes, giving
$r_1(n)=1/(cn(n+1))$. Its sum is $1/c$.
For $L_N$, the terms through $N$ are the same, and each later
ratio of parameters is one. Thus $r_1(n)=1/(cN(N+1))$ for $n>N$.
This positive constant tail makes the series diverge. By
Theorem~\ref{thm:global-lower-criterion}, $L_N$ has no positive
factorization. Summing through row $n$ also gives the finite bounds
\[
 \kappa_{1,n}(L_N)=
 \begin{cases}
 c(n+1)/n,&1\le n\le N,\\
 cN(N+1)/(N^2+n-N),&n\ge N.
 \end{cases}
\]
The formulas agree at $n=N$, and their limits give the two claimed
endpoints.

Multiplication by $I+d\Lambda$ completes the Metelmann representation
of $T_N$. Proposition~\ref{prop:metelmann-input}, applied to each
leading section, gives BTP. In the normalized LDU decomposition
of $T_N$, the lower factor is $L_N$, the diagonal is $I$, and the
upper factor is $I+d\Lambda$. Uniqueness of that decomposition
forces the lower product of any PBF to equal $L_N$. The obstruction
just proved therefore excludes a PBF of $T_N$.

The operator $e_{N+1}e_{N-1}^{\mathsf T}$ has norm one. After
multiplication by the upper factor, the difference becomes
$(c^2/N)e_{N+1}(e_{N-1}+d e_N)^{\mathsf T}$, whose norm is
$(c^2/N)\sqrt{1+d^2}$, proving \eqref{stab:eq:distant-full}.
Finally, for $c,d<1$, geometric series give bounded inverses of
$I+c\Lambda^{\mathsf T}$ and $I+d\Lambda$, and hence of $T_0$.
For large $N$, $\|T_0^{-1}(T_N-T_0)\|_2<1$; a further geometric
series inverts $I+T_0^{-1}(T_N-T_0)$ and therefore $T_N$.
\end{proof}

\begin{remark}[Two topological conclusions]
Already in bandwidth $(2,1)$, a matrix with a PBF can be approached
in operator norm by BTP matrices without a PBF, as in
Theorem~\ref{stab:thm:distant}. Conversely,
Theorem~\ref{thm:boundary-approximation} shows that a bounded BTP
matrix without a PBF is a limit of bounded matrices with PBFs. Thus the PBF
class is neither open nor closed relative to the bounded BTP class.
The first phenomenon occurs even when the limiting operator is
invertible and its given factor entries have uniform positive
lower and upper bounds.

The differences $T_N-T_0$ have rank one. A finite-rank perturbation
preserves the Fredholm essential spectrum, defined as the set of
$z$ for which $zI-T$ fails to have closed range with finite-dimensional
kernel and cokernel. Thus these matrices have the same essential
spectrum and different answers to the global PBF existence question.
\end{remark}

\subsection{Bounded factors with fixed initial entries}
\label{stab:subsec:fixed}

First keep the initial entries at their prescribed values.
The aim is to control every factor change uniformly over all rows.
The lower-product criterion applies to both triangular sides of $LDU$.

\subsubsection{The criterion for a lower product}
Fix $L_0=\LL(a_1)\cdots\LL(a_p)$ with $p\ge1$, and assume
$\alpha\le a_i(n)\le A$ for all $i\in\{1,\dots,p\}$ and $n\ge1$,
where $0<\alpha\le A<\infty$. For a nearby lower matrix, prescribe
\begin{equation}\label{stab:eq:general-initial}
 \widetilde a_i(n)=a_i(n),\qquad
 i\in\{1,\dots,p-1\},\quad n\in\{1,\dots,p-i\}.
\end{equation}
These are the $p(p-1)/2$ initial entries of Lemma~\ref{lem:rowwise}.
All other factor entries are to be determined from the perturbed matrix. A variation of a
factor sequence must therefore be zero at its prescribed positions.

The equations for a fixed row determine its unknown entries from
the preceding rows. In the equations for small variations, however,
an error from a preceding row may be multiplied repeatedly as the
row index grows. The following ratios measure that multiplication:
\begin{equation}\label{stab:eq:adjacent-stability-ratio}
 \gamma_i(n)\coloneq\frac{a_i(n)}{a_{i+1}(n-1)},
 \qquad i\in\{1,\dots,p-1\},\quad n\ge2.
\end{equation}
The condition requires constants $C_i\ge1$ and $0<\rho_i<1$ for which
\begin{equation}\label{stab:eq:adjacent-stability}
 \prod_{t=n}^{n+N-1}\gamma_i(t)\le C_i\rho_i^N
 \qquad(n\ge2,\ N\ge1).
\end{equation}
The same constants work for every starting row $n$ and every length
$N$. Individual ratios may exceed one; the condition controls their
products over sufficiently long intervals.

\begin{theorem}[Analytic dependence with fixed initial entries in lower bandwidth]\label{stab:thm:general-lower-stability}
For the specified factorization $L_0$ and initial entries
\eqref{stab:eq:general-initial}, the following are equivalent:
\begin{enumerate}
\item Every neighboring pair satisfies \eqref{stab:eq:adjacent-stability}.
\item On an operator-norm neighborhood of $L_0$ among bounded real
unit lower matrices of bandwidth at most $p$, there is a real analytic
choice of positive factors
$L=\LL(\widetilde a_1(L))\cdots\LL(\widetilde a_p(L))$,
equal to the given factors at $L_0$, with the initial entries
\eqref{stab:eq:general-initial} fixed.
\end{enumerate}
This choice is the local inverse of multiplication with those
initial entries fixed. On a smaller neighborhood the factor entries
are at least $\alpha/2$. There is a finite constant $C$ such that,
for every pair of matrices $L$ and $\widetilde L$ in this smaller
neighborhood,
\begin{equation}\label{stab:eq:general-lower-lipschitz}
 \max_i\|\widetilde a_i(L)-\widetilde a_i(\widetilde L)\|_\infty
 \le C\|L-\widetilde L\|_2.
\end{equation}
Taking $\widetilde L=L_0$ gives the corresponding bound for the
change from the original factor sequences $a_i$.
For $p=1$, the list of conditions is empty and the only factor is
the matrix itself.
\end{theorem}

The proof first characterizes bounded solvability of a first-order
recurrence, then derives such recurrences from one factor extraction
and applies the calculation successively to the $p$ factors. This identifies exactly when the derivative of
multiplication has a bounded inverse. The analytic inverse function
theorem then gives the factors as analytic functions of the matrix.

\medskip\noindent\emph{A first-order recurrence on bounded sequences.}\par\nobreak\smallskip\noindent
Fix a starting index $n_0$. Given a bounded sequence $f$, consider
$v_{n_0}=f_{n_0}$ and
$v_n=f_n+\gamma_n v_{n-1}$ for $n>n_0$, with bounded positive
weights $\gamma_n$. Each term is uniquely determined, and successive
substitution gives
\[
 v_n=\sum_{s=n_0}^{n}
       \left(\prod_{t=s+1}^{n}\gamma_t\right)f_s.
\]
The term with $s=n$ has coefficient one. Thus the product of weights
from $s+1$ to $n$ is exactly the coefficient with which the datum
$f_s$ contributes to the solution in row $n$.

The recurrence determines every entry, but the resulting sequence
can be unbounded even when all the data are small. For example,
if $\gamma_n=1$ and $f_n=\varepsilon>0$ for every permitted index,
then $v_n=(n-n_0+1)\varepsilon$ for $n\ge n_0$.
Each row adds the same small contribution, and their accumulation
makes the solution unbounded. Uniform control of the variations of
factor entries in every row requires a finite constant $M$, depending
only on the weights, such that every bounded sequence of data gives $\sup_{n\ge n_0}|v_n|\le M\sup_{n\ge n_0}|f_n|$.
The constant must be independent of the row index.

The products in the solution formula describe how earlier data
reach later rows. Uniform geometric decay of these products
controls the sum of all their contributions and gives the required
bound. The next lemma proves that this decay is also necessary.
Both directions are used in Theorem~\ref{stab:thm:general-lower-stability}:
the bound supplies a bounded inverse for the derivative when
constructing perturbed factors, while the converse shows that
analytic dependence with fixed initial entries forces the decay
condition. The following operator formulation states this equivalence.

On the space of bounded sequences indexed by $n\ge n_0$, define
$(Wv)_{n_0}\coloneq0$ and $(Wv)_n\coloneq\gamma_nv_{n-1}$ for $n>n_0$.
This bounded operator is a weighted shift: it moves a sequence one
place and multiplies by the indicated weight. The recurrence is
$(I-W)v=f$. A bounded inverse of $I-W$ gives precisely the bound just described.

\begin{lemma}[Bounded inverse and products over blocks]\label{stab:lem:weighted-shift}
The operator $I-W$ is a bounded bijection on these bounded sequences
if and only if there are $C\ge1$ and $0<\rho<1$ such that
\begin{equation}\label{stab:eq:shift-products}
 \prod_{t=n}^{n+N-1}\gamma_t\le C\rho^N
 \qquad(n\ge n_0+1,\ N\ge1).
\end{equation}
In that case $\|(I-W)^{-1}\|\le1+C\rho/(1-\rho)$.
Equivalently, there are an integer $h\ge1$ and $0<\theta<1$ such
that every product of $h$ consecutive weights is at most $\theta$.
\end{lemma}
\begin{proof}
For $N\ge1$, the nonzero entries of $W^Nv$ are
\[
 (W^Nv)_n=\left(\prod_{t=n-N+1}^{n}\gamma_t\right)v_{n-N},
 \qquad n\ge n_0+N.
\]
Consequently $\|W^N\|$ is the supremum of all products of $N$
consecutive weights. Under \eqref{stab:eq:shift-products}, the
series $I+W+W^2+\cdots$ converges in operator norm. Multiplication
of its partial sum by $I-W$ gives $I-W^{N+1}$, which tends to $I$.
The series is therefore the inverse, and summing its norm bounds
proves the stated estimate.

One fixed block length suffices: if
$\|W^h\|\le\theta<1$, write $N=kh+r$, where
$r\in\{0,\dots,h-1\}$. Then
$\|W^N\|\le\theta^k\|W^r\|$.
With $\rho\coloneq\theta^{1/h}$, the finitely many quantities
$\rho^{-r}\|W^r\|$ give a constant $C$ valid for every $N$.
Conversely, \eqref{stab:eq:shift-products} gives
$\|W^h\|<1$ once $h$ is large enough.

Suppose now that $I-W$ has a bounded inverse, and put
$M\coloneq\|(I-W)^{-1}\|$. Apply the inverse to $f_n=1$ for all $n$.
The explicit solution above has nonnegative terms, so its boundedness
gives
\begin{equation}\label{stab:eq:shift-row-sum}
 \sum_{r=0}^{n-n_0}\prod_{t=n-r+1}^{n}\gamma_t\le M,
 \qquad n\ge n_0.
\end{equation}
In particular $M\ge1$, and every consecutive product of weights
is at most $M$: it occurs as one term in the sum ending at its
last index.

Fix $h\ge1$ and $n\ge n_0+h$. For
$r\in\{0,\dots,h\}$, write
$P_r\coloneq\prod_{t=n-r+1}^{n}\gamma_t$, with $P_0=1$.
The first $h-r$ factors of $P_h$ also form a consecutive product,
so $P_h\le MP_r$. Thus $P_r\ge P_h/M$ for every one of the
$h+1$ choices of $r$. Equation~\eqref{stab:eq:shift-row-sum} now
implies
\[
 \frac{h+1}{M}P_h\le\sum_{r=0}^{h}P_r\le M.
\]
This bound holds for every final row $n$. Hence
$\|W^h\|\le M^2/(h+1)$. Choose $h$ with $h+1>2M^2$.
Then $\|W^h\|<1/2$, and the preceding block argument proves
\eqref{stab:eq:shift-products}.
\end{proof}

The lemma is a positive-coefficient instance of the relation between
bounded solutions for bounded data and exponential stability of
difference equations; see \cite{BravermanKarabash}. The proof also
shows that the row-sum bound \eqref{stab:eq:shift-row-sum} alone
implies the geometric product bound. This implication controls
the diagonal variation of an LDU decomposition below.

\medskip\noindent\emph{The equations for a variation of one extracted factor.}\par\nobreak\smallskip\noindent
For two given factors, a prescribed change in their product
determines the following equations for the factor variations.
Write $a\coloneq a_1$, $b\coloneq a_2$, and
$F\coloneq\LL(a)\LL(b)$. Its leading $4\times4$ block is
\[
 \begin{bNiceMatrix}
 1&0&0&0\\
 a(1)+b(1)&1&0&0\\
 a(2)b(1)&a(2)+b(2)&1&0\\
 0&a(3)b(2)&a(3)+b(3)&1
 \end{bNiceMatrix}.
\]
Thus $F_{n,n-1}=a(n)+b(n)$ for $n\ge1$, whereas
$F_{n,n-2}=a(n)b(n-1)$ for $n\ge2$. The second identity
involves the preceding row of the second factor.

Choose a bounded strictly lower matrix $E$ of bandwidth at most
two, and consider $F+\varepsilon E$ for small real $\varepsilon$.
The diagonal remains equal to one. Define its prescribed
perturbation entries by $E_1(n)\coloneq E_{n,n-1}$ for $n\ge1$ and $E_2(n)\coloneq E_{n,n-2}$ for $n\ge2$.
The leading $4\times4$ block of $E$ is
\[
 \begin{bNiceMatrix}
 0&0&0&0\\
 E_1(1)&0&0&0\\
 E_2(2)&E_1(2)&0&0\\
 0&E_2(3)&E_1(3)&0
 \end{bNiceMatrix}.
\]
The first subdiagonal therefore changes by $\varepsilon E_1(n)$
and the second by $\varepsilon E_2(n)$. These changes are known;
the factor entries remain to be determined.

The required sequences $a_\varepsilon,b_\varepsilon$ satisfy $F+\varepsilon E=\LL(a_\varepsilon)\LL(b_\varepsilon)$, with $a_0=a$ and $b_0=b$,
with $a_\varepsilon(1)=a(1)$ because the initial entry is fixed.
Comparing the two subdiagonals gives the exact equations
\[
\begin{aligned}
 a_\varepsilon(n)+b_\varepsilon(n)
   &=a(n)+b(n)+\varepsilon E_1(n),&&n\ge1,\\
 a_\varepsilon(n)b_\varepsilon(n-1)
   &=a(n)b(n-1)+\varepsilon E_2(n),&&n\ge2.
\end{aligned}
\]
The second equation contains a product of unknown entries.
Its first-order expansion determines the corresponding linear
equations for the factor variations.

For this calculation, substitute $a(n)+\varepsilon h_n$ and
$b(n)+\varepsilon k_n$ for the factor entries. The unknown
sequences $h,k$ describe the changes per unit of $\varepsilon$.
Keeping $a(1)$ fixed requires $h_1=0$. On the first subdiagonal,
the product has entry $a(n)+b(n)+\varepsilon(h_n+k_n)$.
On the second subdiagonal, the exact expansion is
\[
 (a(n)+\varepsilon h_n)(b(n-1)+\varepsilon k_{n-1})
 =a(n)b(n-1)
 +\varepsilon\bigl(b(n-1)h_n+a(n)k_{n-1}\bigr)
 +\varepsilon^2h_nk_{n-1}.
\]
Matching the coefficients of $\varepsilon$ with the prescribed
matrix perturbation gives
\[
\begin{aligned}
 h_n+k_n&=E_1(n),&&n\ge1,\\
 b(n-1)h_n+a(n)k_{n-1}&=E_2(n),&&n\ge2.
\end{aligned}
\]
These are linear equations for $h,k$. They express the derivative
of multiplication at the original factors. The quadratic term
in the expansion remains to be accounted for when obtaining an
exact factorization.

The first equation gives $k_n=E_1(n)-h_n$, and hence
$k_1=E_1(1)$. Substituting this expression at index $n-1$
into the second equation yields $b(n-1)h_n+a(n)\bigl(E_1(n-1)-h_{n-1}\bigr)=E_2(n)$.
Division by the positive entry $b(n-1)$ gives
\[
 h_n=\frac{a(n)}{b(n-1)}h_{n-1}
      +\frac{E_2(n)-a(n)E_1(n-1)}{b(n-1)},\qquad n\ge2.
\]
The coefficient of $h_{n-1}$ is
$\gamma_1(n)=a(n)/b(n-1)$. To display how the known data enter
successive rows, put $\widehat E(n)\coloneq\frac{E_2(n)-a(n)E_1(n-1)}{b(n-1)}$ for $n\ge2$.
This sequence is determined by $E$ and the original factors.
Starting with $h_1=0$, the recurrence gives
\[
\begin{aligned}
 h_2&=\widehat E(2),\\
 h_3&=\gamma_1(3)\widehat E(2)+\widehat E(3),\\
 h_4&=\gamma_1(4)\gamma_1(3)\widehat E(2)
       +\gamma_1(4)\widehat E(3)+\widehat E(4).
\end{aligned}
\]
For instance, the contribution $\widehat E(2)$ is multiplied
by $\gamma_1(3)$ when it reaches row $3$, and by the product
$\gamma_1(4)\gamma_1(3)$ when it reaches row $4$. Each row also
adds its own contribution. This is why the stability condition
controls consecutive products of the adjacent ratios.

The bounds $\alpha\le a(n),b(n)\le A$ give $\|\widehat E\|_\infty\le\frac{\|E_2\|_\infty+A\|E_1\|_\infty}{\alpha}$.
Lemma~\ref{stab:lem:weighted-shift} shows that geometric decay
of the consecutive products of $\gamma_1$ controls the sum of
the contributions from all earlier rows. It therefore bounds
$h$ uniformly in terms of the prescribed data. The identity
$k_n=E_1(n)-h_n$ then gives the same type of bound for $k$.
The initial value and the recurrence also give uniqueness.

This calculation solves the linear equation for the derivative.
The analytic inverse function theorem, applied below, passes
from this bounded inverse of the derivative to actual factors
of $F+\varepsilon E$, including the higher-order changes.
For the resulting analytic factor sequences, $h_n$ and $k_n$
are their derivatives at $\varepsilon=0$. The factors remain
positive when their changes are smaller than the positive lower
bound on the original entries.

The same calculation applies to an arbitrary number
of remaining factors.

The next lemma uses the following norms and completeness properties.

\medskip\noindent\emph{Norms for the variations.}\par\nobreak\smallskip\noindent
Fix a positive integer $m$.
For a strictly lower matrix $E$ with bounded entries and bandwidth
at most $m$, define
\[
 \|E\|_{\mathrm{coef}}\coloneq
 \sup_{\substack{n\ge1\\1\le s\le\min(m,n)}}|E_{n,n-s}|.
\]
This is the supremum of the absolute values of its entries.
This norm is equivalent to the operator norm on $\ell^2(\NN)$.
Since each coordinate vector has norm one, $|E_{i,j}|=|\langle e_i,Ee_j\rangle|\le\|E\|_2$ for $i,j\in\NN$.
Taking the supremum gives $\|E\|_{\mathrm{coef}}\le\|E\|_2$.

For the converse inequality, let $E^{(s)}$ retain only the $s$th
subdiagonal of $E$, so that $E=\sum_{s=1}^m E^{(s)}$.
For $x\in\ell^2(\NN)$, its action is $(E^{(s)}x)_n=E_{n,n-s}x_{n-s}$ for $n\ge s$,
and the coordinates with $n<s$ are zero. Thus $E^{(s)}$ shifts a
sequence by $s$ positions and multiplies its coordinates by the
corresponding matrix entries. In particular,
\[
 \|E^{(s)}x\|_2^2
 =\sum_{n=s}^{\infty}|E_{n,n-s}|^2|x_{n-s}|^2
 \le\|E\|_{\mathrm{coef}}^2\|x\|_2^2.
\]
The triangle inequality now gives
\[
 \|E\|_{\mathrm{coef}}\le\|E\|_2
 \le\sum_{s=1}^m\|E^{(s)}\|_2
 \le m\|E\|_{\mathrm{coef}}.
\]
Hence controlling all entries uniformly is equivalent to controlling
the operator norm when the bandwidth is fixed.

To justify completeness, specify a matrix by its $m$ subdiagonal
sequences. Each is an arbitrary element of $\ell^\infty$, starting
at its first permitted row, and the coefficient norm is the maximum
of their supremum norms. A Cauchy sequence of matrices in this norm
therefore converges uniformly on each subdiagonal, because
$\ell^\infty$ is complete. The limit has the same bandwidth and
bounded entries, so the matrix space is complete in coefficient
norm. The two inequalities above give completeness in operator norm
as well.

\medskip\noindent\emph{Notation for the next lemma.}\par\nobreak\smallskip\noindent
Consider one stage of the factorization, of lower bandwidth $m\ge2$:
\begin{equation}\label{stab:eq:one-stage}
 F=\LL(a)G,\qquad G=\LL(b_1)\cdots\LL(b_{m-1}).
\end{equation}
All factor sequences are bounded above and below by positive
constants. Write $g_s(n)\coloneq G_{n,n-s}$, with $g_0(n)=1$.
At this stage the entries $a(1),\dots,a(m-1)$ are fixed.
Let $E$ be the prescribed bounded real strictly lower matrix of bandwidth
at most $m$ that gives the first-order change in $F$. The unknowns are
corresponding bounded variations $H$ of $\LL(a)$ and $R$ of $G$.
Here $h$ ranges over bounded real sequences indexed by $n\ge m$,
extended by zero for $n<m$ because the initial entries of $a$
are fixed. The only possibly nonzero entries of $H$ are
$H_{n,n-1}\coloneq h_n$ for $n\ge1$. The matrix $R$ ranges
over bounded real strictly lower matrices of bandwidth at most
$m-1$: the diagonal of $G$ remains equal to one and its
bandwidth is unchanged. Measure the pairs $(h,R)$ by
$\max\{\|h\|_\infty,\|R\|_{\mathrm{coef}}\}$ and the
target matrices $E$ by $\|\cdot\|_{\mathrm{coef}}$.
The identity $(\LL(a)+tH)(G+tR)=F+t(HG+\LL(a)R)+t^2HR$
shows that a prescribed first-order change $E$ in $F$ must satisfy
\begin{equation}\label{stab:eq:variable-derivative}
 E=HG+\LL(a)R.
\end{equation}
Thus the derivative is the linear map from the unknown
variations $(h,R)$ to the matrix variation $E$.

The norm argument above applies to the variations $R$, whose
bandwidth is at most $m-1$.
The allowed variations $h$ form the closed subspace of
$\ell^\infty$ defined by $h_n=0$ for $n<m$: a uniform limit
preserves these zero values. This space is consequently complete
in the supremum norm. These are the complete normed spaces used
for the derivative below and for the analytic inverse function
theorem.

On this space of sequences $h$, define the following operators
for $n\ge m$, extending their outputs by zero for $n<m$:
\begin{equation}\label{stab:eq:derivative-first-order}
 (\mathscr D_jh)_n\coloneq b_j(n-j)h_n-a(n)h_{n-1},
 \qquad j\in\{1,\dots,m-1\}.
\end{equation}
\begin{lemma}[First-order factors of the derivative]\label{stab:lem:variable-extraction}
The linear map $(h,R)\longmapsto HG+\LL(a)R$
from these pairs to the bounded real strictly lower matrices of
bandwidth at most $m$ is a bounded bijection if and only if every
$\mathscr D_j$ is a bounded bijection on the indicated sequence
space. For a prescribed target matrix $E$, eliminating $R$ from
$HG+\LL(a)R=E$ gives the equation
\begin{equation}\label{stab:eq:derivative-product}
 \mathscr D_{m-1}\cdots\mathscr D_1h=f,\qquad
 f_n\coloneq\sum_{r=0}^{m-1}(-1)^r\Pi_r(n)E_{m-r}(n-r),\quad n\ge m,
\end{equation}
where $E_s(n)\coloneq E_{n,n-s}$,
$\Pi_0(n)\coloneq1$, and $\Pi_r(n)\coloneq\prod_{t=0}^{r-1}a(n-t)$.
\end{lemma}
\begin{proof}
Elimination of $R$ identifies a scalar equation whose bounded
inverse gives the assertion.

For $1\le s\le\min(m,n)$, comparison of the $s$th subdiagonal
in \eqref{stab:eq:variable-derivative} gives $E_s(n)=h_ng_{s-1}(n-1)+R_s(n)+a(n)R_{s-1}(n-1)$.
Here $R_s(n)\coloneq R_{n,n-s}$, and $R_0=R_m=0$ because $R$ is strictly
lower and has bandwidth at most $m-1$. The equation for $s=1$
is $R_1(n)=E_1(n)-h_n$. In the equation for $s=2$, substitute
this expression at row $n-1$. Continuing in increasing $s$ gives
\begin{equation}\label{stab:eq:variable-remainder}
 R_s(n)=\sum_{r=0}^{s-1}(-1)^r\Pi_r(n)
 \bigl(E_{s-r}(n-r)-h_{n-r}g_{s-r-1}(n-r-1)\bigr)
\end{equation}
for $1\le s\le\min(m-1,n)$.
Every substitution introduces one factor $-a$ and reduces both the
row and the subdiagonal index by one. This explains the sign and
the product $\Pi_r(n)$ in the formula.

For example, the first two remainder equations are
\[
 R_1(n)=E_1(n)-h_n,\qquad
 R_2(n)=E_2(n)-h_ng_1(n-1)-a(n)E_1(n-1)+a(n)h_{n-1}.
\]
The second identity is used only when the second subdiagonal is
present in $R$. These equations illustrate that elimination uses
only a finite number of previous rows, bounded by the bandwidth.

For $n\ge m$ there is also an equation on the $m$th subdiagonal.
Its term $R_m(n)$ is zero, so it imposes a condition on $h$.
Substituting \eqref{stab:eq:variable-remainder} for $R_{m-1}(n-1)$
and collecting the terms in $h$ gives
\begin{equation}\label{stab:eq:variable-scalar-equation}
 \sum_{r=0}^{m-1}(-1)^r\Pi_r(n)
       g_{m-r-1}(n-r-1)h_{n-r}=f_n.
\end{equation}
For rows $n<m$, only \eqref{stab:eq:variable-remainder} is needed;
all the entries of $h$ in those rows are already zero.

The operator on the left of
\eqref{stab:eq:variable-scalar-equation} factors into first-order
operators. For $m=2$, the equation is simply
$b_1(n-1)h_n-a(n)h_{n-1}=f_n$.
For $m=3$, its left side is
\[
 \begin{aligned}
 b_1(n-1)b_2(n-2)h_n
 -a(n)\bigl(b_1(n-2)+b_2(n-2)\bigr)h_{n-1}
       +a(n)a(n-1)h_{n-2}.
 \end{aligned}
\]
Applying $\mathscr D_1$ first and then $\mathscr D_2$ gives precisely
these three coefficients. The general calculation is as follows.

For the induction step, write $G=\overline G\LL(b_{m-1})$, where
$\overline G$ has $m-2$ bidiagonal factors, and set
$\overline g_s(n)\coloneq\overline G_{n,n-s}$.
Multiplication by the last bidiagonal factor gives $g_s(n)=\overline g_s(n)+\overline g_{s-1}(n)b_{m-1}(n-s+1)$.
The first term uses its unit diagonal; the second uses the
subdiagonal entry in row $n-s+1$.
For the product $\mathscr D_{m-2}\cdots\mathscr D_1$, the induction
hypothesis gives the coefficient
$(-1)^r\Pi_r(n)\overline g_{m-r-2}(n-r-1)$ of $h_{n-r}$.
Apply $\mathscr D_{m-1}$ to this expression. Its unshifted term
multiplies that coefficient by $b_{m-1}(n-m+1)$.
Its shifted term is $-a(n)$ times the expression at row $n-1$;
using the coefficient with index $r-1$ there gives
$(-1)^r\Pi_r(n)\overline g_{m-r-1}(n-r-1)$.
Adding the two contributions gives
\begin{multline*}
(-1)^r\Pi_r(n)\overline g_{m-r-1}(n-r-1)+(-1)^r\Pi_r(n)b_{m-1}(n-m+1)\overline g_{m-r-2}(n-r-1)\\=(-1)^r\Pi_r(n)g_{m-r-1}(n-r-1).
\end{multline*}
Subdiagonals outside the range of $\overline G$ are interpreted as
zero. This also accounts for the end terms $r=0$ and $r=m-1$.
The zero initial values of $h$ give the same identity in the first
rows of the sequence space. Thus
\eqref{stab:eq:derivative-product}.

To complete the proof, first consider
$\mathscr P\coloneq\mathscr D_{m-1}\cdots\mathscr D_1$.
Each $\mathscr D_j$ is injective: if its output is zero, the first
row gives $h_m=0$, and the recurrence then gives every later
$h_n=0$. If $\mathscr P$ is surjective, its leftmost factor
$\mathscr D_{m-1}$ is surjective. It is also injective and bounded on a Banach space, so the bounded
inverse theorem gives a bounded inverse. Removing it
shows that $\mathscr D_{m-2}\cdots\mathscr D_1$ is a bounded
bijection. Repeating proves that every factor is a bounded
bijection. The reverse implication follows by composing their
bounded inverses.

If $\mathscr P$ has a bounded inverse, the finite sum defining $f$
gives $\|f\|_\infty\le C\|E\|_{\mathrm{coef}}$ for a constant
independent of the row. Hence it determines a bounded $h$.
Formula~\eqref{stab:eq:variable-remainder} then determines $R$,
with a similar bound, since $m$ is fixed and all coefficients
$a(n)$ and $g_s(n)$ are uniformly bounded. These equations solve
\eqref{stab:eq:variable-derivative} on every subdiagonal.
For zero $E$, they first give $h=0$ and then $R=0$, proving uniqueness.

Conversely, suppose the derivative in
\eqref{stab:eq:variable-derivative} is a bounded bijection.
Any bounded sequence $f_n$, $n\ge m$, can be realized by setting
$E_m(n)=f_n$ and all other subdiagonals of $E$ equal to zero.
Its unique bounded preimage supplies a bounded solution of
$\mathscr P h=f$. The product $\mathscr P$ is injective because
its first-order factors are injective. It is therefore a bounded
bijection, completing the equivalence.
\end{proof}

\begin{proof}[Proof of Theorem~\ref{stab:thm:general-lower-stability}]
\emph{Sufficiency.}
At extraction $i$, the remaining lower bandwidth is $m=p-i+1$.
In the notation of \eqref{stab:eq:one-stage}, the first sequence
is $a=a_i$, and the others are $b_j=a_{i+j}$.
Division of \eqref{stab:eq:derivative-first-order} by
$b_j(n-j)$ reduces it to an operator $I-W_{i,j}$ with weights
\[
 \gamma_{i,j}(n)\coloneq\frac{a_i(n)}{a_{i+j}(n-j)}
 =\prod_{t=i}^{i+j-1}\gamma_t\bigl(n-(t-i)\bigr),
 \qquad n\ge m.
\]
The equality is a telescoping product: the denominator
$a_{t+1}(n-(t-i)-1)$ of one ratio is the numerator of the next.
Every index of an adjacent ratio on the right is at least two,
since $j\le m-1$ and $n\ge m$.

For $N$ consecutive values of $n$, interchange the two finite
products. Applying \eqref{stab:eq:adjacent-stability} to each
adjacent ratio gives
\[
 \prod_{v=n}^{n+N-1}\gamma_{i,j}(v)
 \le\left(\prod_{t=i}^{i+j-1}C_t\right)
       \left(\prod_{t=i}^{i+j-1}\rho_t\right)^N.
\]
Denote the two constants by $C_{i,j}$ and $\rho_{i,j}$,
respectively; $\rho_{i,j}<1$. Lemma~\ref{stab:lem:weighted-shift}
gives a bounded inverse of $I-W_{i,j}$. Multiplication by
$b_j(n-j)$ also has a bounded inverse, since $b_j(n-j)\ge\alpha$.
Thus each $\mathscr D_j$ is a bounded bijection, and
Lemma~\ref{stab:lem:variable-extraction} gives a bounded inverse
for the derivative of this extraction.

To pass from variations to actual factors, consider the multiplication
map $(a,G)\mapsto\LL(a)G$, with the first $m-1$ entries of $a$
fixed. It is polynomial in bounded sequences and fixed-band matrices.
Its domain and target, after subtracting their base points, are
Banach spaces in the norms described above. Its derivative at the
given factors is a bounded bijection. The analytic inverse function
theorem \cite[Chapter~4]{Deimling} therefore gives a unique nearby
pair $(\widetilde a,\widetilde G)$ for every nearby $F$, and the
pair depends analytically on $F$.

Apply this result first to $L_0$, then to its remaining product,
and continue until one bidiagonal factor remains. At stage $i$
the fixed entries are $a_i(1),\dots,a_i(p-i)$, exactly those in
\eqref{stab:eq:general-initial}. At width one the remaining
matrix is itself the last factor. There are finitely many stages;
continuity permits shrinking the neighborhood at each stage so
that the next inverse is defined. The result is one neighborhood of
$L_0$ on which all perturbed factors exist and depend analytically on the matrix.
Shrink it once more so that the maximum change of their sequences
is less than $\alpha/2$. Every factor entry is then positive.

The derivative of this analytic inverse is continuous. It is
therefore bounded on a sufficiently small ball about $L_0$.
Integrating the derivative along the line segment between two
matrices in that ball gives \eqref{stab:eq:general-lower-lipschitz}.
Equivalence of coefficient and operator norms supplies the asserted
operator-norm estimate.

\emph{Necessity.}
Suppose the analytic choice in (2) exists. Write $\Phi$ for
multiplication of all $p$ factors with the prescribed initial
entries fixed, and $\Psi$ for that choice. The identity
$\Phi(\Psi(L))=L$ gives
$\mathrm D\Phi\,\mathrm D\Psi=I$ at $L_0$.
Thus every bounded matrix variation comes from bounded variations
of the factor sequences, with a uniform bound.

These factor variations are unique. At the first extraction the
homogeneous equation has, by
\eqref{stab:eq:variable-scalar-equation}, a scalar recurrence with
nonzero coefficient of $h_n$ and zero initial values. It forces
$h=0$, after which \eqref{stab:eq:variable-remainder} forces
$R=0$. Repeat this for the remaining factors. Hence
$\mathrm D\Phi$ is injective, and $\mathrm D\Psi$ is its bounded
inverse.

Bounded inverses at each separate extraction are also necessary.
The derivative of the first extraction is injective by the same
recurrence. It is surjective because a preimage under
$\mathrm D\Phi$ of any $E$ supplies both the first-factor
variation $H$ and the remaining-product variation $R$.
The bounded inverse theorem gives a bounded inverse. By
Lemma~\ref{stab:lem:variable-extraction}, all its first-order
operators are bounded bijections.

To continue with the remaining product, prescribe an arbitrary
bounded variation $R$ of that product and set $E\coloneq\LL(a_1)R$.
The first extraction has the solution $(H,R)=(0,R)$, which is
unique. Applying the inverse of $\mathrm D\Phi$ to this $E$
therefore gives zero variation in the first factor and variations
of the other factors whose product variation is the prescribed
$R$. Their norms are bounded by a constant times $\|R\|_2$.
This supplies a bounded inverse of multiplication for the remaining
factors. Induction proves the assertion at every extraction.

In particular, at extraction $i$ the operator $\mathscr D_1$ has
a bounded inverse. Its weight is
$a_i(n)/a_{i+1}(n-1)=\gamma_i(n)$.
Lemma~\ref{stab:lem:weighted-shift} gives the geometric bound for
all intervals lying after the first free row of that extraction.
Only finitely many earlier weights are missing. An interval meeting
those earlier rows is a product of an initial part of bounded
length and a part to which the bound applies. Increasing $C_i$
absorbs the finitely many possible initial products and the
corresponding powers of $\rho_i$. This gives
\eqref{stab:eq:adjacent-stability} for all $n\ge2$.
\end{proof}

\subsubsection{The criterion for the complete matrix}
Now let the matrix have both lower and upper diagonals. Its
specified ordered factorization is
\begin{equation}\label{stab:eq:full-given}
 T_0=L_0DU_0,\qquad
 L_0=\LL(a_1)\cdots\LL(a_p),\qquad
 U_0=\LL(b_q)^{\mathsf T}\cdots\LL(b_1)^{\mathsf T},
\end{equation}
where $p,q\ge1$ are arbitrary finite integers. The factor entries
$a_i(n),b_j(n)$ and the diagonal entries $d_n\coloneq D_{n,n}$ belong to
$[\alpha,A]$, with $0<\alpha\le A<\infty$.
Fix $a_i(n)$ for $1\le n\le p-i$ and $b_j(n)$ for
$1\le n\le q-j$. These are the initial entries in
\eqref{stab:eq:general-initial} for $L_0$ and $U_0^{\mathsf T}$,
respectively.

The perturbed PBF is constructed in two steps: first determine the
normalized triangular factors $L,D,U$ of the perturbed matrix $T$.
Then determine the bidiagonal factors of $L$ and $U^{\mathsf T}$. For the second step,
Theorem~\ref{stab:thm:general-lower-stability} requires
\eqref{stab:eq:adjacent-stability} for the ratios
$a_i(n)/a_{i+1}(n-1)$ on the lower side and
$b_j(n)/b_{j+1}(n-1)$ on the upper side. For the first step an additional sequence appears:
\begin{equation}\label{stab:eq:cross-stability-weight}
 w(n)\coloneq a_p(n)b_q(n),\qquad n\ge1.
\end{equation}
Each entry is the product of the last lower and upper factor entries in
the same row. The condition requires constants $C_0\ge1$ and $0<\rho_0<1$ with
\begin{equation}\label{stab:eq:cross-stability}
 \prod_{t=n}^{n+N-1}w(t)\le C_0\rho_0^N
 \qquad(n\ge1,\ N\ge1).
\end{equation}
As on each triangular side, this bound is uniform in the starting
row. The proof below derives this sequence from the equation for
the variation of $D$.

For $p=q=1$, its origin is already visible in ordinary Gaussian
elimination. Writing $d_n(T)$ for the diagonal pivots gives
$d_n(T)=T_{n,n}-T_{n,n-1}T_{n-1,n}/d_{n-1}(T)$ for $n\ge1$.
At the specified factors, differentiation in a direction $E$ gives,
for $n\ge1$,
\[
 \delta d_n=a_1(n)b_1(n)\delta d_{n-1}
 +E_{n,n}-b_1(n)E_{n,n-1}-a_1(n)E_{n-1,n},
\]
with $\delta d_0=E_{0,0}$. The product $w(n)$ is therefore the
coefficient carrying the preceding pivot variation into the next
row. The theorem below identifies the corresponding condition in
arbitrary bandwidth.

\Needspace{8\baselineskip}
\begin{theorem}[Analytic PBF dependence with fixed initial entries in arbitrary bandwidth]
\label{stab:thm:general-full-stability}
Fix the ordered factorization \eqref{stab:eq:full-given}, with
$a_i(n),b_j(n),d_n\in[\alpha,A]$ for
$0<\alpha\le A<\infty$, and keep the initial entries specified
above fixed. The following are equivalent:
\begin{enumerate}
\item There is $\varepsilon>0$ such that every bounded real matrix
$T$ of bandwidth at most $(p,q)$ with $\|T-T_0\|_2<\varepsilon$
has an ordered PBF with these same initial entries.
On this ball, its factor sequences and diagonal can be chosen
as real analytic functions of $T$, in their supremum norms,
and at $T=T_0$ they equal the prescribed factors.
\item The adjacent-product bounds \eqref{stab:eq:adjacent-stability}
hold on both triangular sides, and \eqref{stab:eq:cross-stability}
holds.
\end{enumerate}
The factor map in (1) is the local inverse of multiplication with
the initial entries fixed. On a smaller neighborhood all its factor
entries have a common positive lower bound, and their changes are
bounded by a constant times $\|T-\widetilde T\|_2$.
\end{theorem}

\smallskip\noindent\emph{Remark.}
If $p=1$, there are no adjacent-product conditions on the lower
side; if $q=1$, there are none on the upper side. The condition
\eqref{stab:eq:cross-stability} joining the two sides is still required.

For a matrix $E$ supported in the band $(p,q)$, put
$\|E\|_{\mathrm{coef}}\coloneq\sup_{i,j}|E_{i,j}|$. As for lower matrices,
each diagonal is a weighted shift, and there are at most $p+q+1$
diagonals. Hence
\begin{equation}\label{stab:eq:equivalent-norms}
 \|E\|_{\mathrm{coef}}\le\|E\|_2
 \le(p+q+1)\|E\|_{\mathrm{coef}}.
\end{equation}
The first inequality follows from
$|E_{i,j}|=|\langle e_i,Ee_j\rangle|\le\|E\|_2$.
Thus either norm gives the same complete space of matrices on this
fixed band. Entrywise estimates and
\eqref{stab:eq:equivalent-norms} therefore give operator-norm bounds.

\medskip\noindent\emph{Bounds for the signed triangular inverses.}\par\nobreak\smallskip\noindent
The two stages of the construction require different estimates.
The lower-band theorem already controls the passage from a triangular
matrix to its bidiagonal factors. To control the preceding LDU step,
it is necessary to determine how a change in the full matrix affects its
triangular products. Solving the linearized LDU equation introduces
both triangular inverses. Their separate entries can grow with the
row index; the required estimate concerns the combinations in which
they occur together. The next two lemmas identify these combinations
and explain the role of the crossed products in
\eqref{stab:eq:cross-stability}.

Retain the notation $\mathcal Q$ for the signed lower inverse
used in the minor ratios \eqref{eq:global-minor-ratios}, now applied
to the fixed product $L_0$, and define its upper counterpart as follows.
With the same sign matrix $\Sigma=\diag[(-1)^n]_{n\ge0}$, put $\mathcal Q\coloneq\Sigma L_0^{-1}\Sigma$ and $\mathcal Q^+\coloneq\Sigma(U_0^{\mathsf T})^{-1}\Sigma$.
These inverses are defined entry by entry using their finite leading
sections. The product formula
\eqref{eq:signed-bidiagonal-inverse-entries} makes each signed
bidiagonal inverse nonnegative. Their products $\mathcal Q$ and
$\mathcal Q^+$ are therefore also nonnegative. The adjacent
conditions compare each entry with a product from the last factor.

\begin{lemma}[Bounds in terms of the last factor]\label{stab:lem:last-factor-dominance}
Under the lower adjacent-product conditions, there is $C_L<\infty$
such that
\begin{equation}\label{stab:eq:last-factor-dominance}
 \prod_{t=s+1}^{n}a_p(t)\le\mathcal Q_{n,s}
 \le C_L\prod_{t=s+1}^{n}a_p(t),\qquad n\ge s\ge0.
\end{equation}
Under the upper adjacent-product conditions, the same bounds hold
for $\mathcal Q^+$, with a constant $C_U$ and the sequence $b_q$.
\end{lemma}
\begin{proof}
First convert the shifted ratios in
\eqref{stab:eq:adjacent-stability-ratio} into ratios at the same row.
For $v\ge u\ge2$, cancellation of the intermediate terms gives
\[
 \prod_{t=u}^{v}\frac{a_i(t)}{a_{i+1}(t)}
 =\left(\prod_{t=u}^{v}\gamma_i(t)\right)
   \frac{a_{i+1}(u-1)}{a_{i+1}(v)}.
\]
The last quotient is at most $A/\alpha$.
Multiply these estimates for $i,i+1,\dots,p-1$ over the same
interval. The intervening factor sequences cancel, leaving
$a_i(t)/a_p(t)$. This gives constants $\widetilde C_i\ge1$ and
$0<\tau_i<1$ such that
\begin{equation}\label{stab:eq:same-row-products}
 \prod_{t=u}^{u+N-1}\frac{a_i(t)}{a_p(t)}
 \le\widetilde C_i\tau_i^N,\qquad u\ge1,\ N\ge1,
 \quad i\in\{1,\dots,p-1\}.
\end{equation}
Initially this argument gives $u\ge2$. For $u=1$, separate the
first ratio, which is at most $A/\alpha$, and increase the constant
to include that ratio and one extra factor $\tau_i^{-1}$.

Write $\mathcal Q=\mathcal Q_p\cdots\mathcal Q_1$, where
$(\mathcal Q_i)_{n,s}=\prod_{t=s+1}^{n}a_i(t)$ for $n\ge s$.
In the entry $(n,s)$ of this product, choose the successive
intermediate indices between $s$ and $n$. Equivalently, divide the
interval $\{s+1,\dots,n\}$ into $p$ consecutive segments of
lengths $\ell_1,\dots,\ell_p\ge0$ with total length $n-s$.
The first segment uses $a_1$, the next uses $a_2$, and so on.
Each such division gives exactly one term in the matrix product.

Divide a term by $\prod_{t=s+1}^{n}a_p(t)$. The last segment now
contributes one. For $i<p$, an empty segment also contributes one,
whereas a nonempty segment of length $\ell_i$ contributes at most
$\widetilde C_i\tau_i^{\ell_i}$ by
\eqref{stab:eq:same-row-products}. Once the first $p-1$ lengths
are chosen, the last length is determined by their sum. To obtain
an upper bound, allow all nonnegative choices of those first
$p-1$ lengths, including choices whose sum exceeds $n-s$.
All terms are nonnegative, and the resulting separate geometric
sums give
\[
 C_L\coloneq\prod_{i=1}^{p-1}
 \left(1+\sum_{\ell=1}^{\infty}\widetilde C_i\tau_i^\ell\right)
 =\prod_{i=1}^{p-1}
 \left(1+\frac{\widetilde C_i\tau_i}{1-\tau_i}\right)<\infty.
\]
This proves the upper bound. The choice
$\ell_1=\cdots=\ell_{p-1}=0$ gives the term
$\prod_{t=s+1}^{n}a_p(t)$ itself, proving the lower bound.
For $p=1$, it is the only term and $C_L=1$.
Apply the same calculation to
$U_0^{\mathsf T}=\LL(b_1)\cdots\LL(b_q)$ to obtain the upper-side
bounds.
\end{proof}

The preceding lemma bounds the entries of each signed triangular
inverse. Combining the lower and upper bounds determines
how a prescribed change in the full matrix affects its three
triangular factors.

Let $E$ be a bounded real matrix of bandwidth at most $(p,q)$,
representing a first-order change in $T_0$. Write $X,Z,Y$ for the
corresponding first-order changes in $L,D,U$. Because the diagonals
of $L$ and $U$ remain equal to one, $X$ is strictly lower and $Y$
is strictly upper, with bandwidths at most $p$ and $q$; $Z$ is
diagonal. The coefficient of $t$ in
$(L_0+tX)(D+tZ)(U_0+tY)$ is the expression in
\eqref{stab:eq:band-ldu-derivative} below. Thus that linear map
sends the three factor variations to the resulting variation of
their product. The solution $(X,Z,Y)$ for a given $E$ must be
bounded independently of the row index. The following argument uses the
coefficient norm for each matrix and the maximum of the three
coefficient norms for a triple.

\begin{lemma}[Solving the variation equation for LDU]\label{stab:lem:band-ldu-inverse}
Under both adjacent-product conditions and
\eqref{stab:eq:cross-stability}, the linear map
\begin{equation}\label{stab:eq:band-ldu-derivative}
 (X,Z,Y)\longmapsto XDU_0+L_0ZU_0+L_0DY
\end{equation}
is a bounded bijection between the spaces just described. Explicitly,
for every bounded real $E$ of bandwidth at most $(p,q)$, the equation $XDU_0+L_0ZU_0+L_0DY=E$
has a unique bounded solution with $X$ strictly lower of bandwidth
at most $p$, $Z$ diagonal, and $Y$ strictly upper of bandwidth at
most $q$. There is a constant $C<\infty$, independent of $E$, such that $\max\{\|X\|_{\mathrm{coef}},\|Z\|_{\mathrm{coef}},\|Y\|_{\mathrm{coef}}\}\le C\|E\|_{\mathrm{coef}}$.
\end{lemma}
\begin{proof}
First solve the equation entry by entry, then verify its
bandwidths and uniform bounds. These are separate steps because
the entrywise triangular inverses may themselves be unbounded.
For any fixed entry, triangularity makes all sums finite, so the
algebraic solution is available without an estimate. The remaining
question is whether these solutions, over all rows at once, give
bounded matrices in the required band. Only after that question
has been answered does the inverse function theorem apply on the
infinite-dimensional matrix space.

For a given $E$, form $V\coloneq L_0^{-1}EU_0^{-1}$ entry by entry.
In its entry $(i,j)$ the summation indices satisfy $u\le i$ and
$v\le j$, so the sum is finite. Multiplying the equation
$E=XDU_0+L_0ZU_0+L_0DY$ by the triangular inverses gives $V=L_0^{-1}XD+Z+DYU_0^{-1}$.
The three terms on the right are strictly lower, diagonal, and
strictly upper, respectively. They are therefore determined by
the corresponding three parts $V_-,V_0,V_+$ of $V$. Solving for
the variations gives
\begin{equation}\label{stab:eq:general-ldu-derivative}
 X=L_0V_-D^{-1},\qquad Z=V_0,\qquad Y=D^{-1}V_+U_0.
\end{equation}
This proves existence and uniqueness entry by entry.

The matrices $X$ and $Y$ have the required bandwidths.
Take any finite leading section, and vary its matrix along
$T_0+tE$ near $t=0$. Its leading principal determinants are
nonzero at $t=0$, so its normalized LDU factors exist for small
$t$ and depend differentiably on $t$.
Gaussian elimination preserves the band. At column $j$, entries
to be eliminated occur only through row $j+p$, and the pivot row
has no nonzero entries past column $j+q$. In a changed row $i>j$,
these operations create neither a column beyond $i+q$ nor an
entry below the permitted lower band. The lower multipliers have
indices $i\le j+p$.
The finite factors therefore have lower and upper bandwidths
at most $p$ and $q$ throughout this variation, and so do their
derivatives. Those derivatives solve the same entrywise equation
and must equal \eqref{stab:eq:general-ldu-derivative}.
The compatible leading sections give the claimed bandwidths for
the infinite matrices.

It remains to obtain bounds independent of the section size.
The bandwidth conclusion restricts the required estimates to a
fixed number of diagonals of $V$. For example, a permitted entry
of $X=L_0V_-D^{-1}$ uses $V_{k,j}$ only with $j<k\le i$ and
$i-j\le p$. Thus even if entries of $V$ far from the diagonal
are large, they do not enter this calculation. The same observation
applies to $Y$, and $Z$ uses only the diagonal of $V$.
Put $h\coloneq\max(p,q)$ and $A_0\coloneq\max(1,A)$.
For $|i-j|\le h$, the inverse-entry estimates in
Lemma~\ref{stab:lem:last-factor-dominance} imply
\begin{equation}\label{stab:eq:near-diagonal-v-sum}
 |V_{i,j}|\le C_LC_U\|E\|_{\mathrm{coef}}
 \sum_{\substack{0\le u\le i,\ 0\le v\le j\\-q\le u-v\le p}}
 \left(\prod_{t=u+1}^{i}a_p(t)\right)
 \left(\prod_{t=v+1}^{j}b_q(t)\right).
\end{equation}
The band restriction on $(u,v)$ comes from $E_{u,v}$.
To bound one summand, consider the intervals
$\{u+1,\dots,i\}$ and $\{v+1,\dots,j\}$.
On their overlap the two factors combine into $w(t)$.
The remaining entries are at most $|i-j|+|u-v|\le2h$ in number,
and their product is at most $A_0^{2h}$.

There are $p+q+1$ possible values of $u-v$.
For one fixed value, increasing $u$ and $v$ together reduces a
positive overlap length by exactly one. Thus each positive length
occurs at most once. By \eqref{stab:eq:cross-stability}, the sum
of its bounds is at most $C_0\rho_0/(1-\rho_0)$.
For empty overlap, $\max(u,v)\ge\min(i,j)$. Since
$|u-v|\le h$, it follows that
$\min(i,j)-h\le u\le i\le\min(i,j)+h$.
Thus there are at most $2h+1$ choices for a fixed difference. The whole sum in
\eqref{stab:eq:near-diagonal-v-sum} is therefore at most $(p+q+1)A_0^{2h}\left(2h+1+\frac{C_0\rho_0}{1-\rho_0}\right)$.
In particular there is a constant $K$ independent of $i,j,E$
such that $|V_{i,j}|\le K\|E\|_{\mathrm{coef}}$ whenever
$|i-j|\le h$.

Only these near-diagonal entries occur in the allowed entries of
$X$ and $Y$. For $1\le i-j\le p$, the entry of
$X=L_0V_-D^{-1}$ uses indices $j<k\le i$, so
$1\le k-j\le p$. There are at most $p+1$ summands, each entry
of $L_0$ is at most $(1+A)^p$, and $d_j^{-1}\le\alpha^{-1}$.
Hence
$\|X\|_{\mathrm{coef}}\le
 \alpha^{-1}(p+1)(1+A)^pK\|E\|_{\mathrm{coef}}$.
The same argument gives
$\|Y\|_{\mathrm{coef}}\le
 \alpha^{-1}(q+1)(1+A)^qK\|E\|_{\mathrm{coef}}$,
and $Z=V_0$ has coefficient norm at most $K\|E\|_{\mathrm{coef}}$.
Together with \eqref{stab:eq:equivalent-norms}, these are the
required bounds for the inverse. Multiplication by the fixed
bounded factors makes the original map bounded as well.
\end{proof}

\medskip\noindent\emph{Proof of the full criterion.}\par\nobreak\smallskip\noindent
\begin{proof}[Proof of Theorem~\ref{stab:thm:general-full-stability}]
\emph{Sufficiency.}
Suppose both adjacent-product conditions and
\eqref{stab:eq:cross-stability} hold. Consider multiplication
$(L,D,U)\mapsto LDU$ on unit lower matrices of bandwidth at most
$p$, diagonal matrices, and unit upper matrices of bandwidth at
most $q$. After subtracting $(L_0,D,U_0)$, the variations belong
to the three Banach spaces of
Lemma~\ref{stab:lem:band-ldu-inverse}. Multiplication is a continuous
polynomial map, and its derivative is
\eqref{stab:eq:band-ldu-derivative}.
Lemma~\ref{stab:lem:last-factor-dominance} supplies the estimates
for the two triangular inverses.
Lemma~\ref{stab:lem:band-ldu-inverse} combines them with
\eqref{stab:eq:cross-stability} and gives a bounded inverse of
this derivative on the prescribed matrix spaces.

The analytic inverse function theorem therefore gives normalized
triangular factors $L(T),D(T),U(T)$ for every matrix $T$ in one
operator-norm neighborhood of $T_0$. These factors depend
analytically on $T$ and equal $L_0,D,U_0$ at $T_0$.
Choose the neighborhood small enough that $D(T)$ remains positive
and that $L(T)$ and $U(T)^{\mathsf T}$ lie in the respective
neighborhoods supplied by
Theorem~\ref{stab:thm:general-lower-stability}.
That theorem gives the perturbed lower and upper bidiagonal factors
with the specified initial entries. Their product with $D(T)$ is
$T$, and composition of the analytic maps gives the asserted
analytic PBF.

Continuity in the maximum of the sequence norms permits a further
restriction of the neighborhood so that every entry changes by
less than $\alpha/2$. The entries then have a common positive
lower bound. As in the lower theorem, continuity of the derivative
and integration along segments on a smaller ball give the stated
Lipschitz bound.

\emph{Necessity of the adjacent conditions.}
Suppose the analytic choice in (1) exists. Let $\widetilde L$ vary
near $L_0$ among bounded unit lower matrices of bandwidth at most
$p$, and put $T\coloneq\widetilde LDU_0$. Since
$\|T-T_0\|_2\le\|\widetilde L-L_0\|_2\|DU_0\|_2$, these matrices
lie in the domain of the analytic choice when $\widetilde L$ is
sufficiently close to $L_0$.
Their normalized triangular decomposition is already
$\widetilde L,D,U_0$. Uniqueness in every finite leading section
forces the lower product of the chosen PBF to equal
$\widetilde L$. The resulting lower factor sequences are therefore
an analytic choice of factors for $\widetilde L$, equal to the given
sequences at $L_0$ and with the lower initial entries fixed.
The necessity direction of
Theorem~\ref{stab:thm:general-lower-stability} gives all the lower
adjacent-product conditions. Varying $T=L_0D\widetilde U$ and
applying that theorem to $\widetilde U^{\mathsf T}$ gives the
upper conditions.

\emph{Necessity of the condition joining the two sides.}
Vary the full matrix along $T_0+tI$. This variation stays within
the prescribed band. Differentiate its chosen factorization at
$t=0$, and denote the variations of the triangular products by
$X,Y$ and the variation of the diagonal by $Z=\diag[z_n]$.
The sequence $(z_n)$ is bounded because the derivative of the analytic factor map is a bounded map into the sequence spaces.
Its product equation is $XDU_0+L_0ZU_0+L_0DY=I$.
The entrywise calculation in the proof of
Lemma~\ref{stab:lem:band-ldu-inverse}, leading to
\eqref{stab:eq:general-ldu-derivative}, uses no decay estimate.
It therefore applies here and identifies $Z$ as the diagonal part
of $L_0^{-1}U_0^{-1}$. Consequently
\begin{equation}\label{stab:eq:diagonal-forcing}
 z_n=\sum_{s=0}^{n}\mathcal Q_{n,s}\mathcal Q^+_{n,s}
 \ge\sum_{r=0}^{n}\prod_{t=n-r+1}^{n}w(t),\qquad n\ge0.
\end{equation}
Indeed, the signs of $(L_0^{-1})_{n,s}$ and
$(U_0^{-1})_{s,n}$ are both $(-1)^{n-s}$, so they cancel in
each summand. The lower bounds in
Lemma~\ref{stab:lem:last-factor-dominance} come from the term in
which the last factor supplies all subdiagonal entries and every
other factor supplies its unit diagonal. Thus
$\mathcal Q_{n,s}\ge\prod_{t=s+1}^{n}a_p(t)$ and
$\mathcal Q^+_{n,s}\ge\prod_{t=s+1}^{n}b_q(t)$.
Multiplication and the substitution $r=n-s$ give the inequality.

The right side of \eqref{stab:eq:diagonal-forcing} is a sum of
nonnegative terms bounded by $\sup_n|z_n|$. It is precisely the
row-sum bound \eqref{stab:eq:shift-row-sum}, with starting index
$n_0=0$ and weights $w(n)$. The last part of the proof of
Lemma~\ref{stab:lem:weighted-shift} converts this bound into
\eqref{stab:eq:cross-stability}. This proves the necessity of the
additional condition.

Finally, the inverse is unique near the specified factors with
these initial entries. If two such factorizations have the same
product, uniqueness of normalized LDU first gives the same
triangular products and diagonal. The successive row equations
with fixed initial data then give the same bidiagonal factors on
each side. Thus the analytic factor map just obtained is the
local inverse of full multiplication.
\end{proof}

The following corollary records the consequence for interior points in the
ambient space of all bounded real matrices on the given band.

\begin{corollary}[Interior points of the PBF class]
\label{stab:cor:interior-pbf}
Let $p,q\ge1$ and let $T_0$ have the ordered factorization
\eqref{stab:eq:full-given}, with all entries of $a_i,b_j,D$ in
$[\alpha,A]$ for some $0<\alpha\le A<\infty$.
Suppose \eqref{stab:eq:adjacent-stability} holds on both sides and
\eqref{stab:eq:cross-stability} holds.
There is $\varepsilon>0$ such that every bounded real matrix $T$
of bandwidth at most $(p,q)$ with $\|T-T_0\|_2<\varepsilon$
admits a PBF with bounded factors whose entries on the permitted
diagonals are uniformly bounded away from zero.
Thus $T_0$ is an interior point of the PBF class in this band space.
\end{corollary}
\begin{proof}
Use the analytic choice of factors in
Theorem~\ref{stab:thm:general-full-stability}. By continuity at
$T_0$, choose a ball on which the maximum change in any factor
sequence or diagonal entry is less than $\alpha/2$.
The resulting entries belong to $[\alpha/2,A+\alpha/2]$, and
their ordered product equals $T$ for every matrix in that ball.
\end{proof}

\subsection{Bounded factors with variable initial entries}
\label{stab:subsec:variable}

For a lower product with uniformly positive bounded factor entries,
Theorem~\ref{stab:thm:general-lower-stability} already gives a real
analytic choice of nearby positive factors with the initial entries
fixed whenever the adjacent-ratio products satisfy the uniform
geometric decay bound~\eqref{stab:eq:adjacent-stability}. In that case,
allowing the initial entries to vary is unnecessary. This subsection treats
prescribed factorizations that may fail this decay condition. For two
factors, allowing the initial entry to vary also gives an analytic
choice when the reciprocal products satisfy a uniform geometric
decay bound, as
Proposition~\ref{stab:prop:free-two-factors} shows. The factors must
still depend real analytically on the perturbed matrix and coincide
with the prescribed ones at the original matrix. The norms remain
the operator norm for matrix perturbations within the prescribed band
and the supremum norms for changes in the factor sequences.
For $0<c_1<c_2$, the factorization $\LL(c_1)\LL(c_2)$ satisfies the
decay condition with its initial entry fixed, whereas the reversed
order $\LL(c_2)\LL(c_1)$ requires that entry to vary;
Corollaries~\ref{stab:cor:two-toeplitz}
and~\ref{stab:cor:free-two-toeplitz} establish this distinction for
the same lower matrix. For two factors, adjusting
one initial entry can prevent growth of the factor variations.
The subsequent arguments treat arbitrary lower bandwidth and the
complete matrix.

\subsubsection{Two factors and the choice of the initial entry}

For two lower factors, allowing the first entry to vary introduces
one scalar unknown. The resulting recurrence can be solved either
forward from that entry or backward from the later rows. The next
proposition characterizes exactly when one of these constructions
gives an analytic choice of bounded positive factors.

\begin{proposition}[Two factors with a variable initial entry]
\label{stab:prop:free-two-factors}
Let $L_0=\LL(a)\LL(b)$, where
$0<\alpha\le a(n),b(n)\le A_*<\infty$ for every $n\ge1$, and put
$\gamma(n)\coloneq a(n)/b(n-1)$ for $n\ge2$.
For every bounded real unit lower matrix of bandwidth at most two in
an operator-norm neighborhood of $L_0$, positive factors can be chosen
real analytically in the matrix, equal to the prescribed factors at
$L_0$ and with the initial entry allowed to vary, if and only if one
of the following two bounds holds with $C\ge1$ and $0<\rho<1$, uniformly
for every $n\ge2$ and $N\ge1$:
\begin{equation}\label{stab:eq:free-two-dichotomy}
\prod_{t=n}^{n+N-1}\gamma(t)\le C\rho^N,
 \qquad\text{or}\qquad
 \prod_{t=n}^{n+N-1}\gamma(t)^{-1}\le C\rho^N.
\end{equation}
In the first case the analytic choice can keep $a(1)$ fixed.
In the second case multiplication is a local analytic bijection
on the full space of the two factor sequences. The first entry of the perturbed factor satisfies
\begin{equation}\label{stab:eq:free-two-endpoint}
 \widetilde a(1;L)=\kappa^-_{1,\infty}(L).
\end{equation}
In both cases the factor map is locally Lipschitz, using the operator
norm for matrices and the supremum norms for factor sequences.
\end{proposition}
\begin{proof}
Write the two subdiagonals of $L$ as $s(n)\coloneq L_{n,n-1}$, $n\ge1$,
and $t(n)\coloneq L_{n,n-2}$, $n\ge2$. Multiplication is the polynomial map
$(a,b)\mapsto(s,t)$ given by
$s(n)=a(n)+b(n)$ and $t(n)=a(n)b(n-1)$.
Its domain and range are products of spaces of bounded sequences:
$a,b,s$ start at row one and $t$ starts at row two. In particular,
there is no second-subdiagonal equation in row one; the first
value $a(1)$ is therefore free.
The range norm is equivalent to the operator norm: each coefficient
is bounded by $\|L-I\|_2$, and the operator norm is at most the sum
of the two subdiagonal supremum norms.

Let $h,k$ denote variations of $a,b$, and let $f,g$ be the
prescribed variations of $s,t$. The derivative equations are
$h(n)+k(n)=f(n)$ and
$b(n-1)h(n)+a(n)k(n-1)=g(n)$.
Eliminating $k=f-h$ gives
\begin{equation}\label{stab:eq:free-two-linearized}
 h(n)-\gamma(n)h(n-1)=F(n),\qquad
 F(n)\coloneq\frac{g(n)}{b(n-1)}-\gamma(n)f(n-1),\qquad n\ge2.
\end{equation}
The sequence $F$ is bounded by a constant times
$\max\{\|f\|_\infty,\|g\|_\infty\}$.
Under the first bound in \eqref{stab:eq:free-two-dichotomy}, setting
$h(1)=0$ and solving forward gives a bounded inverse on that
subspace, as in Lemma~\ref{stab:lem:weighted-shift}.
The analytic inverse function theorem gives a local analytic choice of
factors with $a(1)$ fixed.

Under the second bound, a forward solution with an arbitrary initial
value generally grows. The initial value must instead be chosen so
that this growing contribution cancels. To determine it, solve
\eqref{stab:eq:free-two-linearized} from row $n+1$ to a later row $M$
and rearrange:
\[
 h(n)=\frac{h(M)}{\prod_{t=n+1}^{M}\gamma(t)}
       -\sum_{j=n+1}^{M}
        \frac{F(j)}{\prod_{t=n+1}^{j}\gamma(t)},\qquad M>n.
\]
For any bounded solution, the first term tends to zero as
$M\to\infty$, by the reciprocal-product bound. Thus boundedness
itself prescribes the initial value and every subsequent value.
This leads to the candidate solution
\begin{equation}\label{stab:eq:free-two-backward}
 h(n)=-\sum_{j=n+1}^{\infty}
       \frac{F(j)}{\prod_{t=n+1}^{j}\gamma(t)},\qquad n\ge1.
\end{equation}
For a summation index $j=n+r$, the absolute contribution is at
most $C\rho^r\|F\|_\infty$. The geometric series over $r\ge1$
therefore converges uniformly in $n$, with
$\|h\|_\infty\le C\rho(1-\rho)^{-1}\|F\|_\infty$.
To verify the equation, multiply the series for $h(n-1)$ by
$\gamma(n)$. Its term with index $j=n$ is $-F(n)$ and its
remaining terms form exactly the series for $h(n)$.
Thus $h(n)-\gamma(n)h(n-1)=F(n)$, as required in
\eqref{stab:eq:free-two-linearized}. A bounded homogeneous solution
would satisfy
$h(n)=h(1)\prod_{t=2}^{n}\gamma(t)$.
The reciprocal-product bound forces its modulus to grow unless
$h(1)=0$, so this solution is unique.
Taking $k=f-h$ gives a bounded inverse of the full derivative.
The analytic inverse function theorem now gives a local inverse
of multiplication with both sequences free. Shrinking the
neighborhood keeps their entries at least $\alpha/2$ and gives
the asserted Lipschitz bound.

It remains to identify the initial entry of the perturbed factor. The
reciprocal-product condition persists under a sufficiently small
supremum-norm change of $a,b$. Indeed, choose a block length $h_0$
with $C\rho^{h_0}<1$. Products of that fixed length remain uniformly
less than one under a small perturbation; dividing longer intervals
into such blocks and a bounded remainder proves the condition again.
Uniformity follows from the common positive lower bounds: the
reciprocal ratios $b(n-1)/a(n)$ change uniformly when the two
sequences change in supremum norm. Subtracting two products of
length $h_0$ gives $h_0$ terms, each bounded by the change of one
ratio times a fixed bound for the other $h_0-1$ ratios. The strict
inequality therefore survives on one neighborhood, independently
of the starting row.
For any resulting factorization, suppress the tildes and let
$\mathcal Q\coloneq\Sigma L^{-1}\Sigma=\mathcal Q_b\mathcal Q_a$ be
the product of the two signed bidiagonal inverses. Their entry
formulas give
\[
 \mathcal Q_{n,0}-a(1)\mathcal Q_{n,1}
       =\prod_{t=1}^{n}b(t),\qquad
 \mathcal Q_{n,1}\ge\prod_{t=2}^{n}a(t).
\]
For the first identity, column zero of $\mathcal Q_a$ minus
$a(1)$ times its column one is $e_0$: in every row $n\ge1$ the
two products of $a$ cancel. Multiplication by $\mathcal Q_b$
therefore gives its column zero. For the inequality, retain in
$(\mathcal Q_b\mathcal Q_a)_{n,1}$ just the summand with
intermediate index $n$, whose value is the displayed product of
$a$; all other summands are nonnegative. Consequently
$0<\kappa^-_{1,n}-a(1)
 \le b(n)\prod_{t=2}^{n}\gamma(t)^{-1}$.
The right-hand side tends to zero geometrically as $n\to\infty$,
which proves \eqref{stab:eq:free-two-endpoint}.

For necessity, suppose there is a local analytic choice of factors equal to the
prescribed pair at $L_0$. Its derivative sends each bounded variation
$(f,g)$ to bounded factor variations $(h,k)$ satisfying the derivative
equations above. Choose $f(n)=0$ and $g(n)=b(n-1)$ for $n\ge2$.
Then \eqref{stab:eq:free-two-linearized} has a bounded real solution
$h$ with forcing $F(n)=1$. Fix $M\ge2$ such that $|h(n)|\le M$.
The positivity of $\gamma$ gives two possibilities.
The purpose of this positive forcing is to turn a bound on one
solution into a bound on every product of weights. Once $h$ is
nonnegative, the forcing keeps it positive thereafter, and the
telescoping product below has one factor uniformly smaller than one
for each step. If $h$ stays
negative, the same argument applies to the reciprocal weights.
These estimates will be uniform in the starting row, as required
by the statement.

If $h(r)\ge0$ at some index, then $h(n)\ge1$ for every $n\ge r+1$.
For $n\ge r+2$, the equation
$\gamma(t)h(t-1)=h(t)-1$ gives
\[
 \prod_{t=n}^{n+N-1}\gamma(t)
 =\frac{h(n+N-1)}{h(n-1)}
    \prod_{t=n}^{n+N-1}\left(1-\frac1{h(t)}\right)
 \le M\left(1-\frac1M\right)^N.
\]
The finitely many earlier weights are absorbed in a larger constant,
so the forward bound in \eqref{stab:eq:free-two-dichotomy} holds
for every starting row.

Otherwise $h(n)<0$ for all $n\ge1$. Put $v(n)\coloneq-h(n)$, so
$0<v(n)\le M$ and $\gamma(n)v(n-1)=1+v(n)$ for $n\ge2$.
For every $n\ge2$ and $N\ge1$, rearranging this identity gives
\[
 \prod_{t=n}^{n+N-1}\gamma(t)^{-1}
 =\frac{1+v(n-1)}{1+v(n+N-1)}
   \prod_{t=n}^{n+N-1}\frac{v(t-1)}{1+v(t-1)}
 \le(M+1)\left(\frac{M}{M+1}\right)^N.
\]
This is the reciprocal bound. Thus one of the two uniform geometric
bounds is necessary, and the equivalence is proved.
\end{proof}

The two alternatives in \eqref{stab:eq:free-two-dichotomy} express
a scalar exponential dichotomy: homogeneous solutions either
decay forward or decay backward, with bounds uniform in the
starting row. Formula~\eqref{stab:eq:free-two-backward} determines
the variation of the initial entry from all subsequent rows.
The relation between bounded solutions for bounded forcing and
exponential dichotomy has a general theory on the half-line
\cite{HuyMinh2001}. Here positivity reduces the converse to the
single forcing $F=1$ and the two product estimates in the proof.

\subsubsection{The criterion for an arbitrary lower product}

The same choice between forward and backward solution applies at every
extraction in a longer product. The relevant ratios involve each pair
of factors, with the row shift imposed by that extraction.
For example, extracting the first factor from a product of three
requires two first-order equations: one compares the first factor
with the second, and the other compares it with the third.
Extracting the second factor then gives the remaining comparison.
Thus allowing either direction of solution introduces conditions
for all pairs, including pairs that are not adjacent in the
original product.

The converse requires keeping track of the freedom in the initial
entries. A bounded homogeneous solution may be nonzero and still
decay to zero, so the first-order operators need not be injective
on bounded sequences. In the following quotient, these decaying
solutions represent zero. For an index set $I\coloneq\{n_0,n_0+1,\ldots\}$,
let $c_0(I)$ denote the closed subspace of $\ell^\infty(I)$ consisting
of sequences that tend to zero. This subspace is closed: a
uniform limit of sequences tending to zero also tends to zero,
as follows by first approximating uniformly and then taking a
sufficiently late row. In the Banach quotient
$\ell^\infty(I)/c_0(I)$, two bounded sequences represent the same
class precisely when their difference tends to zero. Write
$[v]$ for the class of $v$. Its norm is
\[
 \|[v]\|\coloneq\inf_{r\in c_0(I)}\|v-r\|_\infty
       =\limsup_{n\to\infty}|v_n|.
\]
For the upper bound, subtract the sequence that agrees with $v$
on a finite initial segment and is zero afterwards, then let the
segment increase. For the lower bound, subtracting any sequence
tending to zero leaves the same upper limit of absolute values.
Removing finitely many initial entries does not change the
quotient: restriction and extension by zeros are inverse maps on
classes, since their composition changes only finitely many
entries. The shifts in the recurrence retain their specified indices.

The quotient is used only to prove necessity. The factors and their
analytic dependence will still be measured in the full supremum
norm. Its role in the proof is to remove a particular obstruction:
a decaying homogeneous solution prevents injectivity on bounded
sequences, although it has no effect on the limiting behavior of
those sequences. Once the scalar equation is invertible on the
quotient, it can be removed from a composition to examine the next
equation. Surjectivity alone does not justify that deduction on
the original spaces.

\begin{lemma}[A first-order equation modulo sequences tending to zero]
\label{stab:lem:free-scalar-quotient}
Let $X\coloneq\ell^\infty(\{n_0,n_0+1,\ldots\})$ and
$Y\coloneq\ell^\infty(\{n_0+1,n_0+2,\ldots\})$, and suppose
$0<\alpha\le a(n),b(n)\le A_*<\infty$ for $n\ge n_0+1$.
Define
\[
 \mathscr D:X\longrightarrow Y,\qquad
 (\mathscr Dv)_n\coloneq b(n)v_n-a(n)v_{n-1},\qquad
 \gamma(n)\coloneq\frac{a(n)}{b(n)}.
\]
The induced map on the quotients by $c_0$ is surjective if and only
if one of the bounds
\begin{equation}\label{stab:eq:scalar-quotient-dichotomy}
 \prod_{t=n}^{n+N-1}\gamma(t)\le C\rho^N,
 \qquad\text{or}\qquad
 \prod_{t=n}^{n+N-1}\gamma(t)^{-1}\le C\rho^N
\end{equation}
holds with $C\ge1$ and $0<\rho<1$, uniformly for
$n\ge n_0+1$ and $N\ge1$.
In either case the induced map is a bounded bijection, and
$\mathscr D$ has a bounded linear right inverse on the original
sequence spaces. Moreover, every bounded solution of
$\mathscr Dv\in c_0$ belongs to $c_0$.
\end{lemma}
\begin{proof}
The bounded coefficients give
$\|\mathscr Dv\|_\infty\le2A_*\|v\|_\infty$ and ensure that
$\mathscr D$ maps sequences tending to zero to sequences tending
to zero. Thus two representatives of one class have images in
the same class. The induced map is well defined and bounded in
the quotient norms.

Suppose first that this map is surjective. A preimage of the
constant class $[1]$ has a bounded real representative $v$ with
$\mathscr Dv=1+r$, where $r_n\to0$ as $n\to\infty$.
After division by $b(n)$, the recurrence becomes $v_n-\gamma(n)v_{n-1}=F_n$, where $F_n\coloneq\frac{1+r_n}{b(n)}$.
There are an index $n_1\ge n_0$ and a constant $\delta>0$ such
that $F_n\ge\delta$ for every $n>n_1$; for example, take
$\delta\coloneq1/(2A_*)$ after increasing $n_1$.
Choose $M>\delta$ with $|v_n|\le M$ for every $n$.

If $v_r\ge0$ at some $r\ge n_1$, positivity of $\gamma$ gives
$\delta\le v_n\le M$ for all $n\ge r+1$.
For $n\ge r+2$, telescoping the recurrence yields
\[
 \prod_{t=n}^{n+N-1}\gamma(t)
 =\frac{v_{n+N-1}}{v_{n-1}}
   \prod_{t=n}^{n+N-1}\left(1-\frac{F_t}{v_t}\right)
 \le\frac M\delta\left(1-\frac\delta M\right)^N.
\]
The factors in this product are nonnegative because the preceding
values of $v$ are positive.

Otherwise $v_n<0$ for every $n\ge n_1$. Put $u_n\coloneq-v_n$, so
$0<u_n\le M$ and
$\gamma(n)u_{n-1}=F_n+u_n\ge\delta+u_n$ for $n>n_1$.
For every interval starting after $n_1$,
\[
 \prod_{t=n}^{n+N-1}\gamma(t)^{-1}
 \le\frac{\delta+u_{n-1}}{\delta+u_{n+N-1}}
   \prod_{t=n}^{n+N-1}\frac{u_{t-1}}{\delta+u_{t-1}}
 \le\frac{\delta+M}{\delta}\left(\frac M{\delta+M}\right)^N.
\]
In each case, the finitely many earlier weights are absorbed in
a larger constant $C$. This gives the appropriate bound in
\eqref{stab:eq:scalar-quotient-dichotomy} for every permitted
starting row and length.

Conversely, let $f\in Y$ and write $F_n\coloneq f_n/b(n)$.
Under the forward-product bound, set $v_{n_0}=0$ and solve
\[
 v_n=\sum_{s=n_0+1}^{n}
       \left(\prod_{t=s+1}^{n}\gamma(t)\right)F_s,
 \qquad n>n_0.
\]
The term $s=n$ has coefficient one, and the other coefficients
have total at most $C\rho/(1-\rho)$. Since
$\|F\|_\infty\le\alpha^{-1}\|f\|_\infty$, the solution has norm
at most $\alpha^{-1}(1+C\rho/(1-\rho))\|f\|_\infty$.
It defines a bounded linear right inverse.
If $f_n\to0$ as $n\to\infty$, the contribution from any fixed
finite set of indices $s$ tends to zero. The remaining terms are
bounded by a constant times the supremum of $|F_s|$ on that tail.
More explicitly, given $\varepsilon>0$, choose $s_0$ such that
$|F_s|<\varepsilon$ for $s>s_0$. The corresponding part of the
sum is at most $(1+C\rho/(1-\rho))\varepsilon$, whereas the
finitely many terms with $s\le s_0$ tend to zero as
$n\to\infty$. Taking the upper limit and then letting
$\varepsilon$ decrease to zero proves $v_n\to0$ as $n\to\infty$.
Every other bounded solution differs from this one by
$v_{n_0}\prod_{t=n_0+1}^{n}\gamma(t)$, which also tends to zero.

Under the reciprocal-product bound, use instead
\[
 v_n=-\sum_{s=n+1}^{\infty}
       \frac{F_s}{\prod_{t=n+1}^{s}\gamma(t)},
 \qquad n\ge n_0.
\]
This series converges uniformly and defines a bounded linear
right inverse, as in \eqref{stab:eq:free-two-backward}.
If $f_n\to0$ as $n\to\infty$, its modulus is at most a fixed
constant times $\sup_{s>n}|F_s|$, so it tends to zero.
A nonzero homogeneous solution grows without bound by the
reciprocal-product estimate. Hence this is the only bounded
solution.

If a bounded sequence $v$ has $\mathscr Dv\in c_0$, apply the
preceding construction to this actual output. The constructed
solution is in $c_0$. Its difference from $v$ is a bounded
homogeneous solution, which tends to zero in the forward case
and is identically zero in the reciprocal case. Thus passage to
the quotient removes the possible bounded homogeneous solutions.
Both constructions give surjectivity on bounded sequences and
therefore on their quotients. They also show that
$\mathscr Dv\in c_0$ and $v\in\ell^\infty$ imply $v\in c_0$,
which gives injectivity on the quotients. These quotients are
Banach spaces, so the bounded inverse theorem completes the proof.
\end{proof}

\begin{theorem}[Variable initial entries in arbitrary lower bandwidth]
\label{stab:thm:free-separated-limits}
Let $L_0=\LL(a_1)\cdots\LL(a_p)$, where all factor entries lie
in $[\alpha,A_*]$, with $0<\alpha\le A_*<\infty$. For $1\le i<k\le p$,
put
\begin{equation}\label{stab:eq:pairwise-weight}
 \gamma_{i,k}(n)\coloneq\frac{a_i(n)}{a_k(n-k+i)},
 \qquad n\ge k-i+1.
\end{equation}
The following conditions are equivalent:
\begin{enumerate}
\item For each pair $(i,k)$, one of the alternatives
\begin{equation}\label{stab:eq:pairwise-dichotomy}
 \prod_{t=n}^{n+N-1}\gamma_{i,k}(t)\le C_{i,k}\rho_{i,k}^{N},
 \qquad\text{or}\qquad
 \prod_{t=n}^{n+N-1}\gamma_{i,k}(t)^{-1}\le C_{i,k}\rho_{i,k}^{N}
\end{equation}
holds with $C_{i,k}\ge1$ and $0<\rho_{i,k}<1$, uniformly for
$n\ge k-i+1$ and $N\ge1$. The direction may depend on the pair.
\item There is a local real analytic choice of positive factors
for every bounded real unit lower matrix of bandwidth at most $p$
near $L_0$, equal to the prescribed factors at $L_0$, with the
initial entries allowed to vary.
\item There is a local choice of positive factors for the same matrices,
equal to the prescribed factors at $L_0$ and with variable initial
entries, that is Fr\'echet differentiable at $L_0$.
\end{enumerate}
The matrices are measured in operator norm and the factor sequences
in their supremum norms. The analytic choice in (2) is locally
Lipschitz in these norms. For $p=1$ there are no pair conditions.
\end{theorem}
\begin{proof}
For $p=1$, the matrix is its only factor. Suppose $p\ge2$.
The same extraction equations apply in both directions of the
proof, but ask different questions. For sufficiency, each prescribed
matrix variation must produce at least one bounded variation of all
the factors. For necessity, the existence of such variations must
force the product estimates for each pair. The quotient lemma
separates the first-order equations in this second argument.

\emph{Sufficiency.}
Assume (1). First construct a bounded right inverse of the
derivative at one extraction $F=\LL(a)G$, with
$G=\LL(b_1)\cdots\LL(b_{m-1})$, as in
Lemma~\ref{stab:lem:variable-extraction}. All entries of the variation
$h$, including $h_1,\ldots,h_{m-1}$, are now free.
Elimination of the remainder variation gives
\eqref{stab:eq:variable-remainder} in every row and the scalar
equation \eqref{stab:eq:variable-scalar-equation} for $n\ge m$.
These identities use multiplication and subtraction only, so their
derivation applies to arbitrary initial values of $h$.

Set $X_j\coloneq\ell^\infty(\{j+1,j+2,\ldots\})$ for
$j\in\{0,\ldots,m-1\}$. The first-order operators in the scalar
equation act between successive spaces:
\[
 \mathscr D_j:X_{j-1}\longrightarrow X_j,\qquad
 (\mathscr D_jv)_n\coloneq b_j(n-j)v_n-a(n)v_{n-1},\qquad n\ge j+1.
\]
Their domains include the first entry $v_j$. At the extraction of
$a_i$, the sequence $b_j$ is $a_{i+j}$, so division by $b_j(n-j)$
gives exactly the weight $\gamma_{i,i+j}(n)$ in
\eqref{stab:eq:pairwise-weight}.
If its forward products decay, set $v_j=0$ and solve forward.
The geometric bound makes the sum of propagated forcing terms
uniformly bounded, as in Lemma~\ref{stab:lem:weighted-shift}.
If its reciprocal products decay, use the backward sum
\eqref{stab:eq:free-two-backward}, with starting index $j$.
Division by $b_j(n-j)$ is bounded because $b_j\ge\alpha$.
Thus each $\mathscr D_j$ has a bounded linear right inverse
$S_j:X_j\to X_{j-1}$.

The algebraic factorization \eqref{stab:eq:derivative-product}
maps $X_0$ to $X_{m-1}$. Hence $h=S_1\cdots S_{m-1}f$ solves
that equation with a uniform supremum-norm bound.
Formula~\eqref{stab:eq:variable-remainder} then gives a bounded
remainder variation $R$. This constructs a bounded linear right
inverse of $(h,R)\mapsto HG+\LL(a)R$.
Repeat at the extractions of $a_2,\ldots,a_{p-1}$.
The pair conditions cover every first-order equation that occurs.
Explicitly, first send a prescribed matrix variation $E$ to the
first-factor variation $h_1$ and the remaining-product variation
$R_1$. Send $R_1$ through the next extraction to obtain $h_2,R_2$,
and continue. At the last step the remaining variation is already
bidiagonal. Substitution back into these equations shows that
$(h_1,\ldots,h_p)$ maps to $E$. Each operation is linear and
bounded, and there are only $p-1$ extractions. Their composition
is a bounded linear right inverse $S$ of the derivative of full
multiplication $\Phi$.

The right inverse $S$ selects one bounded factor variation for each
matrix variation. There may be other choices, since a variation in
the kernel changes the factors without changing their product to
first order. Therefore apply the inverse function theorem only
after restricting the factor variations to those selected by $S$.
On this restricted space the derivative is exactly the identity
in the matrix variable; uniqueness of all nearby factorizations
is neither needed nor asserted here.

To obtain the nonlinear choice, write $a^0\coloneq(a_1,\ldots,a_p)$ and
identify the target with the Banach space of strictly lower matrices
of bandwidth at most $p$. The map
$E\mapsto\Phi(a^0+SE)-L_0$ is analytic and has derivative the identity
at zero. The analytic inverse function theorem gives an inverse
$E=\mathcal E(L-L_0)$ near zero. Then
$a^0+S\mathcal E(L-L_0)$ is the required choice of factors.
Continuity keeps every entry at least $\alpha/2$ on a smaller
neighborhood, where a bound on the derivative also gives the
Lipschitz estimate. This proves (1)$\Rightarrow$(2).
An analytic choice is Fr\'echet differentiable, so
(2)$\Rightarrow$(3).

\emph{Necessity.}
Assume (3), and denote the choice by $\Psi$.
The identity $\Phi(\Psi(L))=L$ and the chain rule at $L_0$ give $\mathrm D\Phi(a^0)\,\mathrm D\Psi(L_0)=I$.
Thus the full multiplication derivative is surjective on the
bounded sequence and band-matrix spaces. Its
surjectivity on the quotients by sequences tending to zero already
forces all the pair conditions.
There are two levels to this argument. Within a single extraction,
the first-order operators must be removed one at a time, in their
actual order of composition. After that extraction has been shown
to be invertible on the quotients, it can be removed from the full
multiplication derivative, allowing the next factor to be treated. This
second step is essential: the first extraction only detects pairs
containing $a_1$.

For $m\ge1$, let $\mathcal B_m$ be the Banach space of bounded
strictly lower matrices of bandwidth at most $m$, and let
$\mathcal B_m^0$ be the closed subspace for which each of the
$m$ subdiagonal sequences tends to zero. Equip these spaces with the
coefficient norm, equivalent to the operator norm, and write
$\overline{\mathcal B}_m\coloneq\mathcal B_m/\mathcal B_m^0$.
Likewise, write $\overline X_j$ for the quotient of $X_j$ by its
sequences tending to zero. All multiplication and extraction
derivatives, and all finite sums in the elimination formulas,
preserve these subspaces: their entries use finitely many shifts
with uniformly bounded coefficients. For example, $a(n)h_{n-r}$
tends to zero whenever $h_n\to0$ as $n\to\infty$, since $a$ is bounded and $r$
is fixed. Each component of a derivative is a finite sum of
such terms, with a number of summands bounded by the bandwidth.
They therefore induce
bounded maps on the indicated quotients. Surjectivity of the
full derivative implies surjectivity of its induced map by taking
a representative of each target class.

Consider first the derivative of one extraction, $\mathscr A_m:X_0\oplus\mathcal B_{m-1}\longrightarrow\mathcal B_m$, defined by $\mathscr A_m(h,R)\coloneq HG+\LL(a)R$,
where $H_{n,n-1}=h_n$ and $G=\LL(b_1)\cdots\LL(b_{m-1})$.
Let
\[
 \mathscr P_m\coloneq\mathscr D_{m-1}\cdots\mathscr D_1,
 \qquad
 f(E)_n\coloneq\sum_{r=0}^{m-1}(-1)^r\Pi_r(n)E_{m-r}(n-r),
 \qquad n\ge m,
\]
using the notation of \eqref{stab:eq:derivative-product}.
The equation $\mathscr A_m(h,R)=E$ gives
$\mathscr P_mh=f(E)$, and
\eqref{stab:eq:variable-remainder} determines $R$ from $h,E$.
In the rows $n<m$, that formula only determines the available
entries of $R$; there is no additional scalar equation on the
free initial entries of $h$.

Suppose the induced map $\overline{\mathscr A}_m$ is surjective.
For any $f\in X_{m-1}$, choose $E$ with deepest subdiagonal
$E_m(n)=f_n$ and all other subdiagonals zero. Then $f(E)=f$.
A preimage of $[E]$ satisfies
$\mathscr A_m(h,R)-E\in\mathcal B_m^0$.
The elimination formula, which preserves vanishing sequences,
therefore gives $\mathscr P_mh-f\in c_0$.
Consequently $\overline{\mathscr P}_m=\overline{\mathscr D}_{m-1}\cdots\overline{\mathscr D}_1:\overline X_0\longrightarrow\overline X_{m-1}$
is surjective.
Its range is contained in the range of its leftmost factor
$\overline{\mathscr D}_{m-1}$, so that factor is surjective.
Lemma~\ref{stab:lem:free-scalar-quotient} gives its dichotomy and
makes it invertible on the quotients. Composing with its inverse
shows that the product of the remaining factors is surjective.
For clarity, if $BC$ is surjective and $B$ is bijective, then
for any target $y$ of $C$ there is $x$ with $BCx=By$;
injectivity of $B$ gives $Cx=y$. This is the step that removes
the leftmost factor. Repeating it proves the dichotomy for every
$\mathscr D_j$, without interchanging these operators.
In particular, $\overline{\mathscr P}_m$ is invertible.

Invertibility of the entire extraction on the quotients is also
required to continue to the remaining factors.
If $\mathscr A_m(h,R)\in\mathcal B_m^0$, its eliminated equation
gives $\mathscr P_mh\in c_0$.
The invertibility just proved forces $h\in c_0$, and
\eqref{stab:eq:variable-remainder} then gives
$R\in\mathcal B_{m-1}^0$.
Thus $\overline{\mathscr A}_m$ is injective as well as surjective.
Equivalently, a preimage of $[E]$ is obtained by solving
$\overline{\mathscr P}_m[h]=[f(E)]$ and using the finite formula
for $R$: all subdiagonals except the deepest satisfy their
equations exactly, and the remaining error tends to zero.
The bounded inverse theorem on these Banach quotients makes
$\overline{\mathscr A}_m$ a bounded bijection.

Finally, write $\Phi_p$ for multiplication of all $p$ factors
and $\Phi_{p-1}$ for multiplication of the remaining factors
$(a_2,\ldots,a_p)$, with their derivatives evaluated at the
prescribed tuples. The exact derivative identity is $\mathrm D\Phi_p=\mathscr A_p\circ\bigl(I\oplus\mathrm D\Phi_{p-1}\bigr)$.
The induced full derivative is surjective, so
$\overline{\mathscr A}_p$ is surjective and hence invertible by
the preceding argument. Composing with that inverse shows that
$I\oplus\overline{\mathrm D\Phi}_{p-1}$ is surjective.
Prescribing the target $(0,[R])$ therefore gives variations of
the remaining factors whose product variation has class $[R]$.
Projection onto the second component gives surjectivity of
$\overline{\mathrm D\Phi}_{p-1}$.
At extraction $i$, the induction hypothesis is that multiplication
of the remaining tuple $(a_i,\ldots,a_p)$ has a surjective derivative
on the quotients. The preceding argument proves the pair conditions
involving $a_i$ and preserves this hypothesis for
$(a_{i+1},\ldots,a_p)$. For the last factor the multiplication map
is the identity, so the induction is complete.

At extraction $i$, the sequence $a$ is $a_i$ and the remaining
sequence $b_j$ is $a_{i+j}$. Its scalar weight is therefore $\frac{a_i(n)}{a_{i+j}(n-j)}=\gamma_{i,i+j}(n)$ for $n\ge j+1$.
Every pair $1\le i<k\le p$ occurs once by taking $j=k-i$.
The scalar lemma provides the appropriate uniform product bound
for each of these pairs, proving (1) and completing the equivalence.
\end{proof}

The requirement on every pair in the preceding theorem is essential,
even for constant factors. The next example isolates the obstruction
from a pair that is not adjacent.

\begin{example}[Adjacent pairs and a repeated outer parameter]
\label{stab:ex:nonadjacent-pair}
Consider the constant factorization
$L_0\coloneq\LL(1)\LL(2)\LL(1)$, where $\LL(c)=I+c\Lambda^{\mathsf T}$.
Its adjacent weights are $\gamma_{1,2}=1/2$ and
$\gamma_{2,3}=2$, so both adjacent pairs satisfy
\eqref{stab:eq:pairwise-dichotomy}. The remaining weight is
$\gamma_{1,3}=1$: every forward or reciprocal product is one.
If either estimate held, it would give $1\le C\rho^N$ for
every $N\ge1$, although $C\rho^N\to0$ as $N\to\infty$
when $0<\rho<1$. The necessary part of
Theorem~\ref{stab:thm:free-separated-limits} therefore excludes a
local analytic choice through this prescribed triple, even with
all initial entries allowed to vary. Conditions on adjacent pairs
alone do not suffice.

\end{example}

\subsubsection{The criterion for the complete matrix}

The normalized triangular factors $L,D,U$ are determined by the
matrix. Their analytic dependence is characterized by a bound
on explicit row sums. Combining this bound with the pair conditions
for the two triangular products gives a complete criterion for analytic dependence through
prescribed bidiagonal factors, with their initial entries allowed to
vary.

The row sum $z_n$ defined below has a direct perturbation meaning:
it is the derivative of the $n$th pivot after adding the same small
quantity to every diagonal entry of the matrix. A uniform bound on
these derivatives is therefore necessary. For sufficiency,
the same sums also bound the response to an arbitrary
perturbation of the band. The common positive lower bound on factor
entries permits comparison of nearby entries of each signed inverse;
these comparisons reduce the other derivative entries to the sums
$z_n$.

\begin{proposition}[Row sums and analytic normalized LDU]
\label{stab:prop:ldu-row-sums}
Let $p,q\ge1$ and let $T_0=L_0DU_0$ have the ordered PBF
\eqref{stab:eq:full-given}, with all bidiagonal and diagonal entries
in $[\alpha,A]$, where $0<\alpha\le A<\infty$. With
$\Sigma\coloneq\operatorname{diag}[(-1)^n]_{n\ge0}$, define
\begin{equation}\label{stab:eq:ldu-row-sums}
 \mathcal Q\coloneq\Sigma L_0^{-1}\Sigma,\quad
 \mathcal Q^+\coloneq\Sigma(U_0^{\mathsf T})^{-1}\Sigma,\quad
 z_n\coloneq\sum_{s=0}^n\mathcal Q_{n,s}\mathcal Q^+_{n,s}.
\end{equation}
The triangular inverses are defined entrywise by their finite
leading sections. The following conditions are equivalent:
\begin{enumerate}
\item $Z_*\coloneq\sup_{n\ge0}z_n<\infty$.
\item There is $\varepsilon>0$ such that every bounded real matrix
$T$ of bandwidth at most $(p,q)$ with $\|T-T_0\|_2<\varepsilon$
has a normalized factorization $T=L(T)D(T)U(T)$, where $L(T)$
is unit lower triangular of bandwidth at most $p$, $U(T)$ is
unit upper triangular of bandwidth at most $q$, and $D(T)$ is diagonal.
The three factors depend real analytically on $T$ and equal
$L_0,D,U_0$ at $T_0$.
\item The same factorizations exist on such a neighborhood, with
the three factor maps Fr\'echet-differentiable at $T_0$.
\end{enumerate}
The triangular factors are measured in their fixed-band operator
norms and the diagonal in its supremum norm. These normalized
factors are unique, and their analytic dependence is locally
Lipschitz.
\end{proposition}
\begin{proof}
The row sums control the complete derivative of
normalized LDU, including its triangular components. No bounded
operator inverse of $L_0$ or $U_0$ is assumed.
The main estimate compares the two inverse entries arising from
one entry of the perturbation with entries in a common row and
column. Because the perturbation is banded, the numbers of row
and column shifts in this comparison are uniformly bounded.
The estimate below makes the cost of each such shift explicit.

Write $\mathcal Q=\mathcal Q_p\cdots\mathcal Q_1$, where
$(\mathcal Q_i)_{n,s}=\prod_{t=s+1}^n a_i(t)$ for $n\ge s$.
Let $A_i$ be the lower weighted shift with entry $a_i(n)$ in
position $(n,n-1)$. Entrywise, $\mathcal Q_i=(I-A_i)^{-1}$.
The identities $(I-A_p)\mathcal Q=\mathcal Q_{p-1}\cdots\mathcal Q_1\ge0$ and $\mathcal Q(I-A_1)=\mathcal Q_p\cdots\mathcal Q_2\ge0$
use entrywise inequalities; an empty product is $I$. They give
\begin{equation}\label{stab:eq:inverse-neighbor-comparisons}
 \mathcal Q_{n,s}\ge\alpha\mathcal Q_{n-1,s},\qquad
 \mathcal Q_{n,s}\ge\alpha\mathcal Q_{n,s+1},\qquad n>s.
\end{equation}
Indeed, the $(n,s)$ entries of the two nonnegative products are
$\mathcal Q_{n,s}-a_p(n)\mathcal Q_{n-1,s}$ and
$\mathcal Q_{n,s}-a_1(s+1)\mathcal Q_{n,s+1}$, respectively.
The lower bound $\alpha$ on both coefficients yields
\eqref{stab:eq:inverse-neighbor-comparisons}.
For $\mathcal Q^+$ use its factorization by the signed inverses
of $\LL(b_1),\ldots,\LL(b_q)$; the same computation has
$b_q(n)$ and $b_1(s+1)$ in place of those coefficients.

For a band variation $E$, form $V\coloneq L_0^{-1}EU_0^{-1}$ entrywise.
Its entry $(i,j)$ is a finite sum over $u\le i$ and $v\le j$.
Put $h\coloneq\max(p,q)$ and $\beta\coloneq\min(1,\alpha)$. If $|i-j|\le h$
and $E_{u,v}\ne0$, set $N\coloneq\max(i,j)$ and $s\coloneq\min(u,v)$.
Raising the two row indices to $N$ in
\eqref{stab:eq:inverse-neighbor-comparisons} takes $|i-j|$ steps
in total; lowering the two column indices to $s$ takes $|u-v|$
steps. Hence
\[
 \mathcal Q_{i,u}\mathcal Q^+_{j,v}
 \le\beta^{-(|i-j|+|u-v|)}
       \mathcal Q_{N,s}\mathcal Q^+_{N,s}
 \le\beta^{-2h}\mathcal Q_{N,s}\mathcal Q^+_{N,s}.
\]
There are $p+q+1$ permitted values of $u-v$. For a fixed value,
each $s=\min(u,v)$ determines at most one pair $(u,v)$.
Consequently, under (1),
\begin{equation}\label{stab:eq:row-sums-control-variation}
 |V_{i,j}|\le
 (p+q+1)\beta^{-2h}Z_*\|E\|_{\mathrm{coef}},
 \qquad |i-j|\le h.
\end{equation}

The equation for variations $X,Z,Y$ of $L,D,U$ is $E=XDU_0+L_0ZU_0+L_0DY$.
After multiplication by the triangular inverses, its three terms
are strictly lower, diagonal, and strictly upper, respectively.
Its unique entrywise solution is therefore $X=L_0V_-D^{-1}$, $Z=V_0$, and $Y=D^{-1}V_+U_0$,
as in \eqref{stab:eq:general-ldu-derivative}. The bandwidth
argument in Lemma~\ref{stab:lem:band-ldu-inverse} uses only finite
Gaussian elimination: each finite leading section has nonzero
pivots at $T_0$, and elimination along a sufficiently small
variation preserves its lower and upper bandwidths. Its
derivatives have bandwidths at most $p$ and $q$. Entrywise
uniqueness makes these derivatives compatible as the section size
increases, proving the same bandwidths for $X,Y$.

For an allowed entry of $X$, the sum $L_0V_-D^{-1}$ uses only
indices $j<k\le i$, with $1\le i-j\le p$. Thus only the
near-diagonal entries controlled in
\eqref{stab:eq:row-sums-control-variation} occur. There are at most
$p$ summands, the entries of $L_0$ are bounded by $(1+A)^p$,
and $d_j^{-1}\le\alpha^{-1}$. The upper estimate uses at most
$q$ summands and the bound $(1+A)^q$ for $U_0$.
The same estimate directly bounds $Z=V_0$.
By \eqref{stab:eq:equivalent-norms}, the full LDU derivative
has a bounded inverse on the fixed-band Banach spaces.
The analytic inverse function theorem applied to the polynomial
map $(L,D,U)\mapsto LDU$ proves (1)$\Rightarrow$(2).
Uniqueness follows in every finite leading section. A bounded
derivative on a smaller neighborhood gives local Lipschitz
dependence. An analytic map has a Fr\'echet derivative at every
point of its domain, which proves (2)$\Rightarrow$(3).

For (3)$\Rightarrow$(1), differentiate in the direction $E=I$.
The same finite-section calculation gives
\[
 Z_{n,n}=(L_0^{-1}U_0^{-1})_{n,n}
 =\sum_{s=0}^n\mathcal Q_{n,s}\mathcal Q^+_{n,s}=z_n.
\]
The diagonal derivative is bounded, so $Z_*<\infty$.
\end{proof}

The two requirements combine as follows. The row-sum bound controls
the normalized triangular factors, which are uniquely determined
by the matrix. The pair conditions control the subsequent choice
of bidiagonal factors on each side, with their initial entries free.

\begin{proposition}[Complete criterion with variable initial entries]
\label{stab:prop:free-full-transfer}
Let $p,q\ge1$ and let $T_0=L_0DU_0$ have the ordered PBF
\eqref{stab:eq:full-given}, with all bidiagonal and diagonal entries
in $[\alpha,A]$, where $0<\alpha\le A<\infty$.
The following conditions are equivalent:
\begin{enumerate}
\item There is $\varepsilon>0$ such that every bounded real matrix
$T$ of bandwidth at most $(p,q)$ with $\|T-T_0\|_2<\varepsilon$
has an ordered PBF with product $T$.
On this ball, the factor sequences and diagonal can be chosen as
real analytic functions of $T$, in their supremum norms, and at
$T=T_0$ they equal the prescribed factors. Their initial entries
may depend on $T$.
\item The pair conditions \eqref{stab:eq:pairwise-dichotomy} hold
for $(a_1,\ldots,a_p)$ and $(b_1,\ldots,b_q)$, and
the row sums \eqref{stab:eq:ldu-row-sums} satisfy $\sup_nz_n<\infty$.
\end{enumerate}
The same existence statement with factor maps required to be
Fr\'echet-differentiable at $T_0$ is equivalent to (1) and (2).
The analytic maps are locally Lipschitz in the supremum norms
of the factors and diagonal.

In particular, the pair conditions together with analytic
normalized LDU suffice. The row-sum condition holds if some PBF
of $T_0$ satisfies the fixed-entry criterion of
Theorem~\ref{stab:thm:general-full-stability}. It also holds when
the prescribed factor sequences have positive limits
$a_i^\infty,b_j^\infty$ satisfying
\begin{equation}\label{stab:eq:separated-ldu-limits}
 \left(\max_{1\le i\le p}a_i^\infty\right)
 \left(\max_{1\le j\le q}b_j^\infty\right)<1.
\end{equation}
\end{proposition}
\begin{proof}
Suppose first that a local PBF choice is Fr\'echet-differentiable
at $T_0$. Group its factors into $L(T),D(T),U(T)$. Grouping is
polynomial on the factor-sequence Banach spaces, so these maps
are Fr\'echet-differentiable at $T_0$.
Proposition~\ref{stab:prop:ldu-row-sums} gives $\sup_nz_n<\infty$.

To obtain the lower pair conditions, let $L$ vary over all bounded
unit lower matrices of bandwidth at most $p$ near $L_0$, and set
$S(L)\coloneq LDU_0$. This affine map preserves the band $(p,q)$ and takes
a sufficiently small neighborhood of $L_0$ into the domain of the
given choice. Its normalized triangular factors are exactly
$L,D,U_0$: the pivots are the fixed positive entries of $D$, and
normalized LDU is unique in each finite leading section.
The lower bidiagonal factors returned by the choice consequently
multiply to $L$. They give a local choice through
$(a_1,\ldots,a_p)$ that is Fr\'echet-differentiable at $L_0$.
The necessary part of
Theorem~\ref{stab:thm:free-separated-limits} gives every lower
pair condition. Applying the same argument to $S(U)\coloneq L_0DU$ and
transposing gives every upper pair condition, in the order
$U_0^{\mathsf T}=\LL(b_1)\cdots\LL(b_q)$.

Conversely, under (2),
Proposition~\ref{stab:prop:ldu-row-sums} gives analytic normalized
factors $L(T),D(T),U(T)$. Apply
Theorem~\ref{stab:thm:free-separated-limits} near the prescribed
factorizations of $L_0$ and $U_0^{\mathsf T}$, and compose its
two analytic choices with $L(T)$ and $U(T)^{\mathsf T}$.
Together with $D(T)$, their product is $T$. On a smaller
neighborhood all entries remain uniformly positive and the
derivatives are bounded. This proves (1), local Lipschitz
dependence, and the assertion about differentiable choices.

If some PBF satisfies the fixed-entry criterion, grouping its
analytic perturbed factors gives analytic normalized LDU.
Uniqueness identifies the same triangular products at $T_0$,
and Proposition~\ref{stab:prop:ldu-row-sums} gives the row-sum bound.

Finally, suppose \eqref{stab:eq:separated-ldu-limits} holds.
Choose positive $x_0<x$ and $y_0<y$ such that
$\max_i a_i^\infty<x_0$, $\max_j b_j^\infty<y_0$, and $xy<1$.
Beyond a fixed row, every lower entry is at most $x_0$ and every
upper entry at most $y_0$. A product of $r$ consecutive entries
of one lower factor is at most $C x_0^r$, uniformly in its
starting row; the finite initial prefix is included in $C$.
In the product of the $p$ signed lower inverses, an entry with
$n-s=r$ is a sum over the $\binom{r+p-1}{p-1}$ decompositions
of $r$ into $p$ nonnegative segment lengths. Each segment has
the preceding bound, so the whole entry is at most
$C^p\binom{r+p-1}{p-1}x_0^r$. The quantity
$\binom{r+p-1}{p-1}(x_0/x)^r$ is uniformly bounded in
$r\ge0$: it tends to zero because $0<x_0/x<1$ and the binomial
coefficient is a polynomial in $r$ of degree $p-1$. Increasing
the constant includes the finitely many small values of $r$.
This proves the bound by $C_Lx^r$. Applying the same calculation
to the $q$ upper sequences with $y_0<y$ gives $\mathcal Q_{n,s}\le C_Lx^{n-s}$ and $\mathcal Q^+_{n,s}\le C_Uy^{n-s}$ for $n\ge s\ge0$.
Consequently
\[
 z_n\le C_LC_U\sum_{s=0}^n(xy)^{n-s}
 \le\frac{C_LC_U}{1-xy}.
\]
This proves the final sufficient condition for analytic LDU.
\end{proof}

The sequence $w(n)=a_p(n)b_q(n)$ in
\eqref{stab:eq:cross-stability-weight} belongs to the prescribed
bidiagonal factorization. Analytic dependence of the normalized LDU factors is a property of $T_0$.

\subsection{Consequences and families of nearby factorizations}
\label{stab:subsec:consequences}

The preceding criteria settle when a prescribed PBF can be continued
analytically to nearby matrices, first with fixed initial entries
and then with variable ones. Their consequences describe interiority
and the freedom in the choice of factors. The examples retain the same
operator and supremum norms; relative norms will be introduced in
Subsection~\ref{stab:subsec:relative}.

\subsubsection{Two constant factors and endpoint initial entries}
Theorem~\ref{stab:thm:distant} starts from two equal constant factors.
For two distinct constants, the choice of initial entry determines
which factorizations satisfy the criterion with that entry fixed. The
next corollary describes the whole family. Some of its factorizations
can be varied analytically with the matrix while keeping their initial
entry fixed; others require that entry to vary as well.

\begin{corollary}[Two constant Toeplitz factors]
\label{stab:cor:two-toeplitz}
Let $0<c_1\le c_2$ and $L_0\coloneq\LL(c_1)\LL(c_2)$.
The positive factorizations $L_0=\LL(a_1)\LL(a_2)$ are uniquely
parametrized by $u=a_1(1)\in(0,c_2]$. Moreover,
$\kappa^-_{1,\infty}(L_0)=c_2$, so $u=c_2$ attains the limiting inverse-entry bound.
\begin{enumerate}
\item If $c_1<c_2$, every choice $u\in(0,c_2)$ satisfies
\eqref{stab:eq:adjacent-stability}. In particular, $L_0$ is an
interior point, in operator norm among bounded real unit lower
matrices of bandwidth at most two, of the matrices admitting
a positive bidiagonal factorization.
\item If $c_1<c_2$ and $u=c_2$, the factorization is
$\LL(c_2)\LL(c_1)$ and fails \eqref{stab:eq:adjacent-stability}.
\item If $c_1=c_2$, every positive factorization fails
\eqref{stab:eq:adjacent-stability}, and $L_0$ is an operator-norm
limit of bounded unit lower matrices of bandwidth two without
a positive bidiagonal factorization.
\end{enumerate}
\end{corollary}
\begin{proof}
Put $S\coloneq c_1+c_2$. Matching the two subdiagonals gives
\begin{equation}\label{stab:eq:two-toeplitz-recurrence}
 a_1(n+1)=\varphi(a_1(n)),\qquad
 \varphi(x)\coloneq\frac{c_1c_2}{S-x},\qquad
 a_2(n)=S-a_1(n),\qquad n\ge1.
\end{equation}
Thus $u$ determines all entries, and positivity is equivalent to
all iterates lying in $(0,S)$. The map $\varphi$ is increasing
there, fixes $c_1,c_2$, and satisfies
$\varphi(x)-x=(x-c_1)(x-c_2)/(S-x)$.
If $c_1<c_2$, the iterates starting at $0<u<c_2$ converge
monotonically to $c_1$, with $u=c_1$ giving a constant sequence;
$u=c_2$ gives the other constant sequence. If $c_2<u<S$ and
all iterates stayed in $(0,S)$, they would increase and satisfy
$a_1(n)<S-c_1c_2/S$, since $\varphi(a_1(n))<S$.
Their limit would then be a fixed point greater than $c_2$, a
contradiction. Initial values $u\ge S$ already violate positivity.
The same argument for $c_1=c_2=c$ gives the interval $(0,c]$
and convergence to $c$.

For $c_1<c_2$ and $0<u<c_2$, both sequences have positive
finite bounds, and $\gamma_1(n)\to c_1/c_2<1$ as $n\to\infty$.
Choose $\rho\in(c_1/c_2,1)$ and $n_0\ge2$ such that
$\gamma_1(n)\le\rho$ for $n\ge n_0$.
The finitely many earlier ratios give the uniform constant
$C\coloneq\prod_{t=2}^{n_0-1}\max\{1,\gamma_1(t)/\rho\}$ in
\eqref{stab:eq:adjacent-stability}; hence
Theorem~\ref{stab:thm:general-lower-stability} proves (1).
At $u=c_2$, the ratio is $c_2/c_1>1$, proving (2).
When $c_1=c_2$, the ratio tends to one. A geometric product bound
with constants $C,\rho$ would contradict
$\gamma_1(n)\ge\sigma$ on a tail, with $\sigma\in(\rho,1)$. Theorem~\ref{stab:thm:distant}
then gives the approximating matrices in (3).

Finally, the signed inverse $\mathcal Q=\Sigma L_0^{-1}\Sigma$
is lower Toeplitz with entries $\mathcal Q_{n,s}=h_{n-s}$, where
$h_r\coloneq\sum_{i=0}^{r}c_1^ic_2^{r-i}$.
This is the convolution of the two geometric inverse sequences.
Consequently $\kappa^-_{1,n}=h_n/h_{n-1}$ and
\begin{equation}\label{stab:eq:two-toeplitz-endpoint}
 0<\kappa^-_{1,n}-c_2
 =\frac{c_1^n}{h_{n-1}}\le\frac{c_1}{n},\qquad n\ge1,
\end{equation}
since $h_{n-1}\ge n c_1^{n-1}$. This identifies the endpoint.
\end{proof}

Allowing the initial entry to vary changes the conclusion at the
endpoint. Combining the classification just obtained with
Proposition~\ref{stab:prop:free-two-factors} gives the following
criterion through every factorization of this same lower matrix.

\begin{corollary}[Analytic dependence at every two-factor Toeplitz factorization]
\label{stab:cor:free-two-toeplitz}
Let $0<c_1\le c_2$ and $L_0\coloneq\LL(c_1)\LL(c_2)$.
The inequality $c_1<c_2$ is equivalent to the following statement.
For every prescribed positive factorization
$L_0=\LL(a_1)\LL(a_2)$, there is $\varepsilon>0$ such that every
bounded real unit lower matrix $L$ of bandwidth at most two with
$\|L-L_0\|_2<\varepsilon$ has positive factors
$L=\LL(\widetilde a_1(L))\LL(\widetilde a_2(L))$.
The two factor sequences can be chosen as real analytic functions
of $L$, in their supremum norms, with $\widetilde a_i(L_0)=a_i$
for $i\in\{1,2\}$ and the initial entry allowed to vary.
The neighborhood and the choice may depend on the prescribed pair.
When $c_1<c_2$, each sufficiently close matrix has exactly one
factorization in a sufficiently small supremum-norm neighborhood
of the endpoint pair $\LL(c_2)\LL(c_1)$, with initial entry
$u(L)=\kappa^-_{1,\infty}(L)$.
\end{corollary}
\begin{proof}
Corollary~\ref{stab:cor:two-toeplitz} lists all factorizations.
For $0<u<c_2$ the adjacent products decay, and for $u=c_2$
their reciprocals decay because $\gamma(n)=c_2/c_1>1$.
Proposition~\ref{stab:prop:free-two-factors} applies in each case.
If $c_1=c_2$, Theorem~\ref{stab:thm:distant} gives arbitrarily
nearby matrices without any PBF. A choice of positive factors throughout a neighborhood would factor
those matrices, proving the converse.
\end{proof}

\begin{remark}[Initial coordinates and the uniform topology]
\label{stab:rem:initial-versus-uniform}
The interval of initial entries in
Corollary~\ref{stab:cor:two-toeplitz} and the associated family of
factor sequences carry different topologies. Suppose $c_1<c_2$,
and write $a_i(n;u)$ for the factors with $a_1(1;u)=u$.
For every $0<u<c_2$ their limits are $(c_1,c_2)$, whereas
$u=c_2$ gives the constant pair $(c_2,c_1)$. Therefore
\[
 \max_{i\in\{1,2\}}\sup_{n\ge1}
       |a_i(n;u)-a_i(n;c_2)|\ge c_2-c_1,
 \qquad 0<u<c_2.
\]
The endpoint factorization is isolated in the supremum norm even
though its initial entry is a limit point of the interval $(0,c_2]$.
On each compact subinterval of $(0,c_2)$, the identity
\[
 \frac{a_1(n;u)-c_1}{a_1(n;u)-c_2}
   =\left(\frac{c_1}{c_2}\right)^{n-1}
       \frac{u-c_1}{u-c_2},\qquad n\ge1,
\]
gives continuous dependence in the supremum norm: solving for
$a_1(n;u)$ gives denominators uniformly bounded away from zero
there, and a geometric bound for the tail. Thus the factor family
has two connected components in the uniform topology, the isolated
endpoint and the family $0<u<c_2$. Its initial-coordinate domain
is a single connected interval. Evaluation of the initial entry
is continuous, but the inverse map recovering the factor sequences
is discontinuous at $u=c_2$.

When $c_1=c_2=c$, put $\varepsilon\coloneq c-u\in[0,c)$. The same
recurrence gives the exact formulas
\[
 a_1(n;u)=c-\frac{c\varepsilon}{c+(n-1)\varepsilon},\qquad
 a_2(n;u)=c+\frac{c\varepsilon}{c+(n-1)\varepsilon}.
\]
The supremum distance from the constant pair is $\varepsilon$.
For $\varepsilon>0$, the deviations satisfy
$|a_i(n;u)-c|\sim c/n$ as $n\to\infty$, and the adjacent ratios
tend to one. The endpoint is now approached in the uniform
topology, while the geometric product condition still fails.
This example separates the geometry in initial coordinates from
the local dimensions of uniformly close factor sequences studied
in Theorem~\ref{stab:thm:local-factor-family}.
\end{remark}

\begin{figure}[htbp]
\centering
\begin{tikzpicture}[x=1cm,y=1cm,font=\small,
  stable/.style={blue!65!black,line width=1.4pt},
  unstable/.style={red!65!black,line width=1.4pt,densely dashed}]
 \node[anchor=west] at (0,1.25)
   {(a) Distinct constants: $L_0=\LL(1)\LL(2)$};
 \draw[stable] (0,0) -- (8,0);
 \draw[stable,fill=white] (0,0) circle[radius=2.7pt];
 \fill[blue!65!black] (4,0) circle[radius=2.2pt];
 \fill[teal!70!black] (8,0) circle[radius=2.7pt];
 \node[blue!65!black] at (3.5,.50) {Analytic variation with fixed $u$};
 \node[teal!70!black,anchor=south,align=center,text width=2.4cm] at (8,.25)
   {Analytic variation with variable $u$};
 \node[below=5pt] at (0,0) {$0$};
 \node[below=5pt,align=center] at (4,0)
   {$u=1$\\[4pt]$\LL(1)\LL(2)$};
 \node[below=5pt,align=center] at (8,0)
   {$u=2$\\[4pt]$\LL(2)\LL(1)$};
 \node[anchor=west] at (0,-1.65)
   {(b) Equal constants: $L_0=\LL(1)^2$};
 \draw[unstable] (0,-2.85) -- (4,-2.85);
 \draw[red!65!black,line width=1.4pt,fill=white]
   (0,-2.85) circle[radius=2.7pt];
 \fill[red!65!black] (4,-2.85) circle[radius=2.7pt];
 \node[red!65!black] at (2,-2.35)
   {Nearby matrices may have no PBF};
 \node[below=5pt] at (0,-2.85) {$0$};
 \node[below=5pt,align=center] at (4,-2.85)
   {$u=1$\\[4pt]$\LL(1)^2$};
\end{tikzpicture}
\caption{Choice of the initial entry $u=a_1(1)$ for two fixed
lower Toeplitz matrices. The colored intervals are exactly the
admissible sets; the open circles exclude $u=0$. In (a), every
interior initial value permits an analytic choice of factors with $u$
fixed; the endpoint permits such a choice with $u$ varying. In (b),
every neighborhood contains matrices without any PBF, so even varying
$u$ cannot give positive factors throughout a neighborhood. Corollaries~\ref{stab:cor:two-toeplitz}
and~\ref{stab:cor:free-two-toeplitz} prove these assertions.
The marked products have constant subdiagonal entries.
The horizontal axis records $u$, not the uniform distance between
factor sequences; the upper endpoint is isolated in that norm.}
\label{fig:stable-factor-choice}
\end{figure}
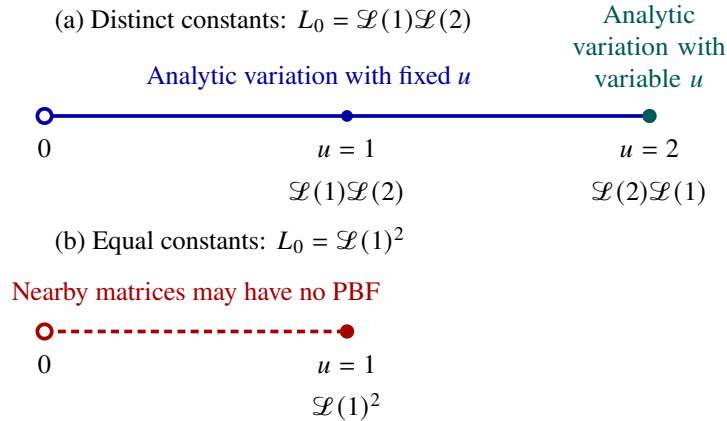

Figure~\ref{fig:stable-factor-choice} shows the admissible initial
entries in two instances of Corollary~\ref{stab:cor:two-toeplitz}.
Each point represents one factorization of the matrix fixed in its row.

\medskip\noindent\emph{An endpoint initial entry in the full matrix.}\quad
The next proposition explains how an endpoint of the allowable
initial data prevents choosing positive factors throughout a neighborhood
with those data fixed. The proof uses the bounded-factor setting
of \eqref{stab:eq:full-given} and changes only one initial matrix entry.

\begin{proposition}[An endpoint prevents keeping the initial entry fixed]
\label{stab:prop:endpoint-unstable}
Suppose $p\ge2$ and the specified PBF of $T_0=L_0DU_0$ has
$a_1(1)=\kappa^-_{1,\infty}>0$, where
$\kappa^-_{1,n}=\mathcal Q_{n,0}/\mathcal Q_{n,1}$ and
$\mathcal Q=\Sigma L_0^{-1}\Sigma$.
Every operator-norm neighborhood of $T_0$ in the fixed band contains
a matrix having no PBF with these same initial entries.
The corresponding assertion holds on the upper side when $q\ge2$
and $b_1(1)=\kappa^+_{1,\infty}$.
\end{proposition}
\begin{proof}
Let $E_{1,0}\coloneq e_1e_0^{\mathsf T}$ and, for $\varepsilon>0$, put
$L_\varepsilon\coloneq(I-\varepsilon E_{1,0})L_0$.
Row zero of the unit lower matrix $L_0$ is $e_0^{\mathsf T}$.
Thus $E_{1,0}L_0=E_{1,0}$ and
$L_\varepsilon=L_0-\varepsilon E_{1,0}$.
Also $E_{1,0}^2=0$, so
$(I-\varepsilon E_{1,0})^{-1}=I+\varepsilon E_{1,0}$.
Conjugating the inverse by $\Sigma$, and using
$\Sigma E_{1,0}\Sigma=-E_{1,0}$, gives $\mathcal Q_\varepsilon\coloneq\Sigma L_\varepsilon^{-1}\Sigma=\mathcal Q(I-\varepsilon E_{1,0})$.
Consequently column zero becomes
$\mathcal Q_{\cdot,0}-\varepsilon\mathcal Q_{\cdot,1}$, while
column one is unchanged.

Suppose $L_\varepsilon$ had a positive factorization whose first
factor had first subdiagonal entry $u=a_1(1)$. Removing this factor
would multiply $L_\varepsilon^{-1}$ on the right by that factor.
After conjugation by $\Sigma$, its effect on column zero is to
subtract $u$ times column one. The remaining positive product
would therefore have signed inverse with column zero $\mathcal Q_{\cdot,0}-(\varepsilon+u)\mathcal Q_{\cdot,1}$.
There is at least one remaining factor because $p\ge2$.
Its signed inverse is strictly positive below the diagonal, and
multiplication by the other nonnegative signed inverses preserves
that positivity. The displayed column must therefore be positive
in every row $n\ge1$. Since $\mathcal Q_{n,1}>0$, this implies
$u+\varepsilon<\kappa^-_{1,n}$ for all such $n$.
Taking the infimum gives
$u+\varepsilon\le\kappa^-_{1,\infty}=u$, a contradiction.

Finally set $T_\varepsilon\coloneq L_\varepsilon DU_0$.
Its change is $-\varepsilon E_{1,0}DU_0$, supported in the allowed
band, and its norm tends to zero with $\varepsilon$.
Its normalized triangular factors are $L_\varepsilon,D,U_0$.
Uniqueness of that decomposition would force any PBF with the
same initial entries to have the excluded lower factorization.
This proves the assertion. Apply the same argument to
$T_0^{\mathsf T}=U_0^{\mathsf T}DL_0^{\mathsf T}$ for the upper side.
\end{proof}

The perturbation in this proof has an explicit positive
factorization when the initial entry is allowed to change.
Row zero of $L_1$ is $e_0^{\mathsf T}$, so
$(I-\varepsilon E_{1,0})L_1=L_1-\varepsilon E_{1,0}$ and $L_\varepsilon=(L_1-\varepsilon E_{1,0})L_2\cdots L_p$.
For $0<\varepsilon<a_1(1)$ all these factors remain positive.
Their only changed parameter is $a_1(1)-\varepsilon$.
Thus the proposition concerns the prescribed initial entries.
For two distinct constant Toeplitz factors,
Corollary~\ref{stab:cor:free-two-toeplitz} gives an analytic choice of factors
on an entire neighborhood even at the endpoint; the endpoint itself
varies analytically with the matrix.

\subsubsection{Factor limits and interior points}
Distinct positive limits, or distinct geometric means in the periodic
case, make the pair conditions explicit. They can occur in any order
when the initial entries are allowed to vary.

\begin{corollary}[Distinct limits and distinct periodic means]
\label{stab:cor:free-limits-periodic}
The equivalent conditions of Theorem~\ref{stab:thm:free-separated-limits}
hold in either of the following cases:
\begin{enumerate}
\item $a_i(n)\to c_i>0$ as $n\to\infty$, with the $p$ limits
pairwise distinct, in any order.
\item The factor sequences have a common period $h_{\mathrm{per}}$,
and the geometric means
\[
 g_i\coloneq\left(\prod_{n=1}^{h_{\mathrm{per}}}a_i(n)\right)^{1/h_{\mathrm{per}}},
 \qquad i\in\{1,\dots,p\},
\]
are pairwise distinct, in any order.
\end{enumerate}
\end{corollary}
\begin{proof}
In (1), $\gamma_{i,k}(n)\to c_i/c_k\ne1$ as $n\to\infty$.
Choose $v(n)\coloneq\gamma_{i,k}(n)$ when $c_i<c_k$, and
$v(n)\coloneq\gamma_{i,k}(n)^{-1}$ when $c_i>c_k$. Its limit is less
than one, so there are $\rho<1$ and $n_1$ with $v(n)\le\rho$
for $n\ge n_1$. The finite constant
\[
 C\coloneq\prod_{t=k-i+1}^{n_1-1}\max\{1,v(t)/\rho\}
\]
(an empty product is one) gives
$\prod_{t=n}^{n+N-1}v(t)\le C\rho^N$ for every permitted $n,N$.
Indeed, any interval uses each earlier index at most once.
This proves \eqref{stab:eq:pairwise-dichotomy} in the chosen direction.
In (2), the product of $\gamma_{i,k}$ over a period is
$(g_i/g_k)^{h_{\mathrm{per}}}$: the shift of the denominator
permutes its entries over that period. Choose the direction in
which this product is less than one, and split any interval into
complete periods and a remainder of length less than
$h_{\mathrm{per}}$. If the chosen period product is
$\vartheta<1$, write $N=kh_{\mathrm{per}}+r$, where
$0\le r<h_{\mathrm{per}}$. Let $B$ be the largest product over
all remainder lengths and starting phases, including the empty
product. The full product is at most $B\vartheta^k$. With
$\rho\coloneq\vartheta^{1/h_{\mathrm{per}}}$, this is at most
$B\rho^{-(h_{\mathrm{per}}-1)}\rho^N$, which proves the uniform
bound in \eqref{stab:eq:pairwise-dichotomy}.
\end{proof}

For the complete matrix, normalized LDU must also be controlled.
The next proposition shows that this part of the criterion holds
at every interior point of the PBF class. Its conclusion concerns
the triangular products; choosing their bidiagonal factors still
requires the pair conditions.

\begin{proposition}[Uniform factors and analytic LDU at interior points]
\label{stab:prop:interior-ldu}
Let $p,q\ge1$ and let $T$ be an operator-norm interior point of
the PBF class among bounded real matrices of bandwidth at most
$(p,q)$. Every normalized PBF of $T$ has all its bidiagonal and
diagonal entries in a common interval $[\alpha,A]$, with
$0<\alpha\le A<\infty$; the constants may depend on the chosen
PBF. For every sufficiently close bounded real matrix of bandwidth
at most $(p,q)$, the normalized factors $L,D,U$ depend real
analytically on that matrix, using operator norm for the triangular
factors and the supremum norm for the diagonal.

More precisely, if the open ball $B(T,r)$ lies in the PBF class
and $M=\|T\|_2$, its pivots satisfy $r\le d_n\le M$, and
\begin{equation}\label{stab:eq:interior-row-sum-bound}
 \sup_{n\ge0}z_n\le\frac{M^2}{\varepsilon(r-\varepsilon)}
 \quad(0<\varepsilon<r),\qquad
 \sup_{n\ge0}z_n\le\frac{4M^2}{r^2}.
\end{equation}
\end{proposition}
\begin{proof}
The interior ball has two uses. First, a perturbation of one entry
shows that pivots and extreme diagonals cannot approach zero.
Positivity then gives uniform bounds for every factor sequence.
Second, comparison of the finite inverses of $T$ and $T-\varepsilon I$
gives the required row-sum bound.
Both matrices remain in the PBF class when $0<\varepsilon<r$,
and the signs of their inverses turn this comparison into the
row-sum bound needed in Proposition~\ref{stab:prop:ldu-row-sums}.

Consider first a matrix $S$ with a ball $B(S,s)$ contained in the
PBF class. Subtracting its normalized pivot $d_n(S)$ from entry
$(n,n)$ annuls the leading principal determinant of order $n+1$.
The cofactor of that entry is the preceding leading determinant,
so the new determinant is
\[
 \prod_{j=0}^n d_j(S)-d_n(S)\prod_{j=0}^{n-1}d_j(S)=0.
\]
The empty product for $n=0$ is one. This perturbation has norm
$d_n(S)$ and stays in the band. A positive PBF has positive
leading principal determinants, so $d_n(S)\ge s$. At $S=T$,
positivity also gives $d_n\le T_{n,n}\le M$.

Every entry on the two extreme diagonals of $T$ is at least $r$:
otherwise setting that single entry to zero would produce, within
distance less than $r$, a matrix without a positive PBF.
Fix any normalized PBF $T=LDU$. Positivity gives $L_{n,n-1}d_{n-1}\le T_{n,n-1}$ and $d_{n-1}U_{n-1,n}\le T_{n-1,n}$.
The first off-diagonal of a product of unit bidiagonal factors is
the sum of their corresponding entries. Thus every individual
bidiagonal entry is at most $B\coloneq M/r$.
On the extreme diagonals, the intermediate index in $LDU$ is
unique, and hence $T_{n,n-p}=L_{n,n-p}d_{n-p}$ and $T_{n,n+q}=d_nU_{n,n+q}$.
Each extreme entry of $L$ or $U$ is therefore at least
$r/M$. It is a product of exactly one entry from each of its
bidiagonal factors. Bounding the other entries by $B$ gives a
lower bound $(r/M)B^{-(p-1)}$ for the selected lower entry, and
$(r/M)B^{-(q-1)}$ for the selected upper entry.
For the lower factors, the extreme entry is
$L_{n,n-p}=\prod_{i=1}^p a_i(n-i+1)$. Given $a_i(t)$ with
$t\ge p-i+1$, take $n=t+i-1$ to include it in this product.
For the upper factors apply the same identity to
$U^{\mathsf T}=\LL(b_1)\cdots\LL(b_q)$.
Thus these products cover a tail of every factor sequence.
The finitely many remaining entries are positive. This proves
common positive lower and finite upper bounds for all factor
entries, and boundedness of the triangular products.

To bound the row sums, fix $0<\varepsilon<r$ and put
$T_\varepsilon\coloneq T-\varepsilon I$. Its ball of radius
$r-\varepsilon$ lies in the PBF class, so
$d_n(T_\varepsilon)\ge r-\varepsilon$.
For each $n$, let $A\coloneq T^{[n]}$ and
$B\coloneq T_\varepsilon^{[n]}=A-\varepsilon I$, of order $n+1$, and set $K_A\coloneq\Sigma_nA^{-1}\Sigma_n$, $K_B\coloneq\Sigma_nB^{-1}\Sigma_n$, and $\Sigma_n\coloneq\operatorname{diag}[(-1)^j]_{j\in\{0,\dots,n\}}$.
These matrices are entrywise nonnegative: their inverse entries
are signed cofactors of totally nonnegative matrices divided by
positive determinants. The finite resolvent identity gives $K_B-K_A=\varepsilon K_BK_A\ge0$.
The identity follows by conjugating
$B^{-1}-A^{-1}=B^{-1}(A-B)A^{-1}$, where $A-B=\varepsilon I$.
Thus $K_B\ge K_A\ge0$ entrywise. Multiplication on the right
by the nonnegative matrix $K_A$ preserves the inequality, giving
$K_BK_A\ge K_A^2$. Consequently
\[
 \varepsilon(K_A^2)_{n,n}
 \le(K_B-K_A)_{n,n}
 \le(K_B)_{n,n}
 =\frac1{d_n(T_\varepsilon)}\le\frac1{r-\varepsilon}.
\]
In $A^{-1}=U_{[n]}^{-1}D_{[n]}^{-1}L_{[n]}^{-1}$, the last row
is the last row of $L_{[n]}^{-1}$ divided by $d_n(T)$, and the
last column is the last column of $U_{[n]}^{-1}$ divided by the
same pivot. Consequently $z_n=d_n(T)^2(K_A^2)_{n,n}$.
Since $d_n(T)\le M$, this proves the first bound in
\eqref{stab:eq:interior-row-sum-bound}; setting $\varepsilon=r/2$
proves the second. Proposition~\ref{stab:prop:ldu-row-sums}
gives analytic normalized LDU. Only finite-section inverses have
been used; no operator inverse of $T$ is required.
\end{proof}

The next lemma gives a necessary condition for interior membership.
It uses finite leading sections and applies even when the chosen
factors have no uniform bounds.

\begin{lemma}[A cross-product obstruction to interior membership]
\label{qt:lem:cross-obstruction}
Let $p,q\ge1$ and let $T=LDU$ be bounded, with a positive
bidiagonal factorization
$L=\LL(a_1)\cdots\LL(a_p)$ and
$U=\LL(b_q)^{\mathsf T}\cdots\LL(b_1)^{\mathsf T}$.
If for some pair $i,j$ $a_i(n)b_j(n)\longrightarrow c\ge1$ as $n\to\infty$,
then $T$ is not an operator-norm interior point of the PBF class
in the full band $(p,q)$. In fact, for every $\varepsilon>0$,
$T-\varepsilon I$ is not totally nonnegative.
\end{lemma}
\begin{proof}
Define the entrywise signed inverses $\mathcal Q,\mathcal Q^+$
and $z_n$ by \eqref{stab:eq:ldu-row-sums}. Each signed bidiagonal
inverse is nonnegative with diagonal one. A product of such
matrices dominates any one of its factors entrywise. Therefore
\[
 \mathcal Q_{n,s}\ge\prod_{t=s+1}^n a_i(t),\qquad
 \mathcal Q^+_{n,s}\ge\prod_{t=s+1}^n b_j(t),\qquad 0\le s\le n.
\]
For every fixed integer $M\ge0$ and $n\ge M$, this gives
\[
 z_n\ge\sum_{k=0}^{M}\prod_{t=n-k+1}^{n}a_i(t)b_j(t),
 \qquad
 \liminf_{n\to\infty}z_n\ge\sum_{k=0}^{M}c^k.
\]
Since $c\ge1$ and $M$ is arbitrary, $z_n\to\infty$ as $n\to\infty$.

Normalized LDU restricts to each leading section $T^{[N]}$ of
order $N+1$. Finite triangular inversion and cyclicity of the
trace give
\begin{equation}\label{qt:eq:trace}
 \operatorname{tr}\bigl((T^{[N]})^{-1}\bigr)
 =\sum_{n=0}^{N}\frac1{d_n}
          \sum_{s=0}^{n}\mathcal Q_{n,s}\mathcal Q^+_{n,s}
 =\sum_{n=0}^{N}\frac{z_n}{d_n}.
\end{equation}
Positivity gives $0<d_n\le T_{n,n}\le\|T\|_2$.
Thus the mean of the right side tends to infinity. Explicitly,
given $H>0$, choose $n_0$ so that $z_n/d_n\ge2H$ for
$n\ge n_0$, using $d_n\le\|T\|_2$. For $N$ sufficiently
large, at least half the indices $0,\ldots,N$ lie after $n_0$,
so the mean is at least $H$.
Every leading section is oscillatory, by
Proposition~\ref{prop:positive-product} and
Corollary~\ref{cor:btp-oscillatory}, and has positive eigenvalues.
The eigenvalues of the inverse are their reciprocals, so
the trace is their sum and is at most
$(N+1)/\lambda_{\min}^{[N]}$, where $\lambda_{\min}^{[N]}$
is the smallest eigenvalue. Rearranging gives
\[
 0<\lambda_{\min}^{[N]}
 \le\frac{N+1}{\operatorname{tr}((T^{[N]})^{-1})}
 \longrightarrow0\qquad(N\to\infty).
\]
For every $\varepsilon>0$, a sufficiently large leading section
of $T-\varepsilon I$ has a negative eigenvalue. A totally
nonnegative finite matrix has only nonnegative real eigenvalues,
so that section is not totally nonnegative. Since
$\|T-(T-\varepsilon I)\|_2=\varepsilon$, interior membership
is impossible.
\end{proof}

\begin{corollary}[Interior criterion for positive factor limits]
\label{stab:cor:limiting-interior}
Let $p,q\ge1$ and let a bounded matrix $T=LDU$ have the prescribed
positive bidiagonal factorization
$L=\LL(a_1)\cdots\LL(a_p)$ and
$U=\LL(b_q)^{\mathsf T}\cdots\LL(b_1)^{\mathsf T}$.
Suppose
$a_i(n)\to a_i^\infty>0$ and $b_j(n)\to b_j^\infty>0$
as $n\to\infty$, with pairwise distinct limits within each of
the two lists, in any order. The following are equivalent:
\begin{enumerate}
\item The pivots and limiting entries satisfy
\[
 \inf_{n\ge0}d_n>0,\qquad
 \left(\max_{1\le i\le p}a_i^\infty\right)
 \left(\max_{1\le j\le q}b_j^\infty\right)<1.
\]
The second inequality is \eqref{stab:eq:separated-ldu-limits}.
\item There is $\varepsilon>0$ such that every bounded real matrix
$\widetilde T$ of bandwidth at most $(p,q)$ satisfying
$\|\widetilde T-T\|_2<\varepsilon$ has an ordered PBF with
product $\widetilde T$. On this ball the factor sequences and
diagonal can be chosen as real analytic functions of $\widetilde T$
in their supremum norms. At $\widetilde T=T$ they equal the
prescribed factors, and their initial entries are allowed to vary.
\item $T$ is an interior point of the PBF class, in operator norm
among all bounded real matrices of bandwidth at most $(p,q)$.
\end{enumerate}
The diagonal sequence $d_n$ need not converge.
\end{corollary}
\begin{proof}
Positive convergence gives common positive lower and finite upper
bounds for the finitely many bidiagonal sequences. Moreover,
$0<d_n\le\|T\|_2$. Under (1), the pivots have a positive lower
bound as well. Corollary~\ref{stab:cor:free-limits-periodic} gives
all pair conditions, and
Proposition~\ref{stab:prop:free-full-transfer} gives the row-sum
bound and the analytic choice. This proves (1)$\Rightarrow$(2).
Condition (2) supplies a full-band neighborhood in the PBF class,
so it implies (3).

If (3) holds, Proposition~\ref{stab:prop:interior-ldu} gives
$\inf_nd_n>0$. If the strict product bound in (1) failed, choose
indices attaining the two largest limits. Their product would
satisfy $a_i(n)b_j(n)\to c\ge1$ as $n\to\infty$, contrary to
Lemma~\ref{qt:lem:cross-obstruction}. This proves (3)$\Rightarrow$(1).
\end{proof}

The row-sum condition is intrinsic to normalized triangular
products. The pair conditions concern their prescribed ordered
bidiagonal factors. Proposition~\ref{stab:prop:interior-ldu}
provides the first condition at every interior point; it does not
by itself produce a factorization satisfying the second.
The preceding corollary supplies both when the factor limits are
positive and distinct within each side and obey the strict
cross-product bound.

\subsubsection{All sufficiently close factorizations}
So far, an analytic choice has assigned one factorization to each
nearby matrix. The next theorem describes all factorizations sufficiently close
to the prescribed one. Their number of independent parameters is
determined by which homogeneous variations remain bounded. This is
a question about the whole factor sequences in the supremum norm,
as illustrated by Remark~\ref{stab:rem:initial-versus-uniform}.

The alternatives in \eqref{stab:eq:pairwise-dichotomy} also determine
how many independent parameters can vary while the product stays fixed.
Let $s_{\mathrm{loc}}$
be the number of pairs for which the first estimate holds.
The two estimates cannot hold for the same pair: multiplying them
would bound one by a quantity tending to zero as $N\to\infty$.

\begin{theorem}[All nearby factorizations of a fixed matrix]
\label{stab:thm:local-factor-family}
Under the hypotheses and equivalent conditions of
Theorem~\ref{stab:thm:free-separated-limits}, set
$a^0\coloneq(a_1,\ldots,a_p)$ and
$\Phi(a)\coloneq\LL(a_1)\cdots\LL(a_p)$.
There are neighborhoods $\mathcal U$ of $a^0$ in
$(\ell^\infty)^p$, $\mathcal V$ of $L_0$ in the unit lower
matrices of bandwidth at most $p$, and $W$ of zero in
$\RR^{s_{\mathrm{loc}}}$, together with a real analytic bijection
\[
 \Theta:\mathcal V\times W\longrightarrow\mathcal U,
 \qquad
 \Phi(\Theta(L,z))=L,
 \qquad \Theta(L_0,0)=a^0,
\]
whose inverse is real analytic. All factor entries in $\mathcal U$
are positive. Thus every factorization of $L_0$ sufficiently close
to the prescribed one belongs to an analytic family with exactly
$s_{\mathrm{loc}}$ independent real parameters.

If $s_{\mathrm{loc}}=0$, each matrix sufficiently close to $L_0$
has exactly one factorization sufficiently close to $a^0$.
Closeness of factor tuples is measured by the maximum of their
supremum norms, and closeness of matrices by operator norm.
\end{theorem}
\begin{proof}
First count the kernel of the derivative and then recover the
actual factors from that calculation. A direction in this kernel
is a bounded change of the factor sequences whose first-order
change in the matrix is zero. Such a direction is an infinitesimal
candidate for varying the factors while keeping the product fixed.
Counting these directions alone does not yet give a family of
actual factorizations. The final application of the inverse function
theorem will show that they do provide independent local parameters
and that every sufficiently close factorization is included.

At an extraction of lower bandwidth $m$, use the spaces $X_j$ and
operators $\mathscr D_j:X_{j-1}\to X_j$ from the proof of
Theorem~\ref{stab:thm:free-separated-limits}.
That proof gives a bounded right inverse for each $\mathscr D_j$.
After division by $b_j(n-j)$, a homogeneous solution satisfies
\[
 v_n=\gamma(n)v_{n-1}\quad(n\ge j+1),\qquad
 v_n=v_j\prod_{t=j+1}^{n}\gamma(t)\quad(n\ge j),
\]
where $\gamma(n)\coloneq a(n)/b_j(n-j)$ and the empty product is one.
If the products of $\gamma$ decay, every value of $v_j$ gives a
bounded solution, and the kernel has dimension one.
If the reciprocal products decay, the displayed product grows
without bound. Boundedness then forces $v_j=0$, and the kernel is
zero.

These dimensions add by the following elementary identity.
Let $A:X\to Y$ and $B:Y\to Z$ be bounded surjections, and let
$S_A$ be a bounded right inverse of $A$. Then $\ker(BA)=\ker A\oplus S_A(\ker B)$.
Indeed, for $x\in\ker(BA)$,
$x=(x-S_AAx)+S_AAx$, with the first summand in $\ker A$
and $Ax\in\ker B$. The intersection of the two summands is
zero because $AS_A=I$. Thus finite kernel dimensions add when
these surjections are composed.

For one extraction, the derivative is
$(h,R)\mapsto HG+\LL(a)R$, as in
\eqref{stab:eq:variable-derivative}, now with all initial entries
of $h$ free. With $E=0$, equation
\eqref{stab:eq:variable-remainder} determines $R$ uniquely from
$h$. The remaining equation is
$\mathscr D_{m-1}\cdots\mathscr D_1h=0$.
The finite formula for $R$ has uniformly bounded coefficients,
so bounded $h$ gives bounded $R$.
Consequently the kernel of this extraction derivative is
isomorphic to $\ker(\mathscr D_{m-1}\cdots\mathscr D_1)$.
The preceding identity shows that its dimension is the number
of decaying products at that extraction.

At extraction $i$, the weights are precisely
$\gamma_{i,k}$ for $k\in\{i+1,\ldots,p\}$.
The derivative of full multiplication is the derivative of the
first extraction composed with the identity on the first sequence
and the derivative of multiplication of the remaining factors.
All these maps are surjective with bounded right inverses, by
Theorem~\ref{stab:thm:free-separated-limits} and its proof.
Adding their kernel dimensions inductively counts every pair once.
Hence, writing $A_0\coloneq\mathrm D\Phi(a^0)$, $\dim\ker A_0=s_{\mathrm{loc}}$.

Choose the bounded right inverse $S$ of $A_0$ constructed in that
proof, and put $\mathcal K\coloneq\ker A_0$ and $P_{\mathcal K}\coloneq I-SA_0$.
The identities $A_0S=I$, $A_0P_{\mathcal K}=0$, and
$P_{\mathcal K}z=z$ for $z\in\mathcal K$ show that this is
a bounded projection onto $\mathcal K$. They also imply
$P_{\mathcal K}^2=P_{\mathcal K}$ and $P_{\mathcal K}S=0$.
Every variation has the unique decomposition
$h=SA_0h+P_{\mathcal K}h$, with its second summand in the
kernel. The kernel consists of those factor variations whose
first-order change of the product is zero.
The extra coordinate $P_{\mathcal K}h$ records the factor changes
which are not determined by their product.
Consider the analytic map $h\longmapsto\bigl(\Phi(a^0+h)-L_0,P_{\mathcal K}h\bigr)$.
Its derivative at zero is $h\mapsto(A_0h,P_{\mathcal K}h)$,
with bounded inverse $(E,z)\mapsto SE+z$: applying $A_0$ to
$SE+z$ gives $E$, applying $P_{\mathcal K}$ gives $z$, and
the preceding decomposition verifies the other inverse identity.
The analytic inverse function theorem
\cite[Chapter~4]{Deimling} gives an analytic bijection between
open neighborhoods of zero. Choose a product of two smaller open
neighborhoods in its range and take its inverse image as the
factor neighborhood. This gives the product form stated in the
theorem. Identifying the finite-dimensional space
$\mathcal K$ with $\RR^{s_{\mathrm{loc}}}$ gives $\Theta$.
This parametrization includes every sufficiently close factor tuple,
because the inverse function theorem gives a bijection on a
neighborhood in the full factor space.
Taking the matrix coordinate equal to $L_0$ leaves exactly the
free coordinate $z$.
Finally shrink the neighborhood until all factor entries change
by less than $\alpha/2$. Their positivity follows from the
uniform lower bound at $a^0$.
\end{proof}

\begin{corollary}[Local dimension from limits or periodic means]
\label{stab:cor:local-dimension-limits}
In case (1) of Corollary~\ref{stab:cor:free-limits-periodic}, the
number of independent parameters in
Theorem~\ref{stab:thm:local-factor-family} is $s_{\mathrm{loc}}=\#\{(i,k):1\le i<k\le p,\ c_i<c_k\}$.
In case (2), the same formula holds with the geometric means
$g_i$ in place of $c_i$.
In particular, an increasing order has local dimension
$p(p-1)/2$, while a decreasing order has local dimension zero.
\end{corollary}
\begin{proof}
With distinct limits, $\gamma_{i,k}(n)\to c_i/c_k$ as
$n\to\infty$, so its products decay precisely when $c_i<c_k$.
For periodic factors, the product over a period is
$(g_i/g_k)^{h_{\mathrm{per}}}$, so it decays precisely when
$g_i<g_k$. The product estimates are those already proved in
Corollary~\ref{stab:cor:free-limits-periodic}.
\end{proof}

The finite-dimensional family in Section~\ref{sec:finite} uses
$p(p-1)/2$ initial entries. In the semi-infinite problem, a change
which excites a growing homogeneous solution cannot keep every row
uniformly close to the original factors. The count
$s_{\mathrm{loc}}$ gives the number of parameters compatible with
this boundedness requirement.
Local uniqueness concerns the supremum norm over all rows.
Remark~\ref{stab:rem:initial-versus-uniform} explains why an
endpoint factorization can be isolated in this norm while its
initial entries are approached by those of other positive
factorizations.

The uniqueness of normalized LDU also separates the freedom on the
two sides of a fixed matrix. Its diagonal cannot vary within a
factorization of that matrix; only the bidiagonal factors of its
fixed lower and upper products can vary.

\begin{proposition}[All nearby PBFs in the full band]
\label{stab:prop:local-full-family}
Let $T_0=L_0DU_0$ have the ordered PBF
\eqref{stab:eq:full-given}. Assume that the entries of the
bidiagonal factors lie in a common interval $[\alpha,A]$,
with $0<\alpha\le A<\infty$, and that the pair conditions
\eqref{stab:eq:pairwise-dichotomy} hold for both
$(a_1,\ldots,a_p)$ and $(b_1,\ldots,b_q)$.
Let $s_-$ and $s_+$ be their respective parameter counts from
Theorem~\ref{stab:thm:local-factor-family}.
Closeness of factors is measured by the maximum of the supremum
norms of the differences between corresponding entry sequences.
In (2), closeness of matrices is measured in operator norm.
\begin{enumerate}
\item With the matrix $T_0$ fixed, all PBFs sufficiently close to
the prescribed one form an analytic family of dimension $s_-+s_+$.
Its diagonal is fixed, and only the bidiagonal factors vary.
This assertion does not require analytic dependence of LDU.
\item Suppose in addition that
$0<d_*\le d_n\le d^*<\infty$ for all $n\ge0$, and that every
sufficiently close bounded real matrix $T$ of bandwidth at most
$(p,q)$ has normalized factors $L(T),D(T),U(T)$ depending real
analytically on $T$, in operator norm for the triangular factors
and the supremum norm for the diagonal.
Then all sufficiently close PBFs are parametrized analytically by
$(T,z_-,z_+)$, where $T$ ranges over such a neighborhood and
$z_-,z_+$ range over neighborhoods of zero in
$\RR^{s_-},\RR^{s_+}$, respectively. Each parameter tuple gives
one PBF with product $T$, and every sufficiently close PBF occurs
exactly once. At $(T_0,0,0)$ the factors are the prescribed ones.
\end{enumerate}
In (2), fixing the matrix leaves precisely the parameters $z_-,z_+$.
When $s_-+s_+=0$, this parametrization gives exactly one PBF
sufficiently close to the prescribed one for each nearby matrix,
with factors depending analytically on that matrix.
\end{proposition}
\begin{proof}
The normalized LDU decomposition is unique in every finite leading
section, so it is unique entry by entry in the semi-infinite matrix.
Every PBF of the fixed matrix $T_0$ must therefore have lower product
$L_0$, diagonal $D$, and upper product $U_0$.
The transposed upper product is
$U_0^{\mathsf T}=\LL(b_1)\cdots\LL(b_q)$, in that order.
Apply Theorem~\ref{stab:thm:local-factor-family} to these two lower
products. Their parameter choices are independent, and the diagonal
is fixed. This proves the first assertion.

When LDU depends analytically on $T$, insert $L(T)$ and
$U(T)^{\mathsf T}$ in those two parametrizations and retain
$D(T)$. The resulting factors multiply to $T$.
Conversely, uniqueness of LDU and the fact that the two
parametrizations include all nearby lower factorizations show that
every sufficiently close PBF is obtained in this way.
The inverse parameters are the bounded projections used in the
proof of Theorem~\ref{stab:thm:local-factor-family}, applied to the
two factor increments.
Shrink the neighborhood until the bidiagonal entries change by
less than $\alpha/2$ and
$\sup_n|d_n(T)-d_n|<d_*/2$. All entries remain positive.
The two parameter spaces are singletons when $s_-+s_+=0$,
which gives the stated uniqueness.
\end{proof}

\subsection{Relative perturbations with unrestricted factor growth}
\label{stab:subsec:relative}

With the initial entries fixed, relative norms extend the stability
criterion to factors with arbitrary positive entries. Each change is
divided by the original entry, so its size adapts to the scale of that
entry. The resulting criterion contains
Theorem~\ref{stab:thm:general-full-stability} and covers factors with
different, unrestricted growth rates.
Componentwise perturbation norms and finite LU sensitivity are
classical; see Higham \cite{Higham1994} and Bueno--Dopico
\cite{BuenoDopico2004}. The estimates below treat the infinitely
many rows simultaneously in arbitrary finite bandwidth.

Fix finite integers $p,q\ge1$ and a semi-infinite product
$T_0=LDU$, where $D=\diag[d_n]_{n\ge0}$ and
\[
 L=\LL(a_1)\cdots\LL(a_p),\qquad
 U=\LL(b_q)^{\mathsf T}\cdots\LL(b_1)^{\mathsf T}.
\]
All entries $a_i(n),b_j(n),d_n$ are positive, with unrestricted
magnitudes and growth. Products are defined entrywise. A \emph{path}
is a choice of intermediate indices in such a product; its contribution
is the product of the corresponding entries. Each bidiagonal factor
allows an index to stay fixed or change by one, so the number of paths
per entry is bounded solely by $p+q$.
Every entry of $T_0$ inside its band is positive: for $n\ge m$
one can use a term of $L_{n,m}$ followed by $d_m$ and the unit
diagonal of $U$, while for $n\le m$ one can use $d_nU_{n,m}$.
The required lower or upper entry is positive whenever $(n,m)$
is in the band.

For a real perturbation supported in this band, put
\begin{equation}\label{stab:eq:relative-norm}
 \|E\|_{T_0,\mathrm{rel}}
 \coloneq\sup_{\substack{n,m\ge0\\-q\le n-m\le p}}
       \frac{|E_{n,m}|}{(T_0)_{n,m}}.
\end{equation}
The space of perturbations with finite norm is Banach: entrywise
division by $T_0$ identifies it with the bounded arrays on the band.
For a unit triangular product $F$, the norm
$\|E\|_{F,\mathrm{rel}}$ is defined in the same way on its strict
triangular band; the diagonal variation is zero. This convention
will be used for $L$ and $U$ and for the remaining triangular
products during extraction. Relative factor variations are measured
by the maximum of
$\sup_{i,n}|\delta a_i(n)|/a_i(n)$,
$\sup_{j,n}|\delta b_j(n)|/b_j(n)$, and
$\sup_n|\delta d_n|/d_n$.

For a positive sequence $v(n)$, write $\mathcal G(v)$ for the condition
\begin{equation}\label{stab:eq:relative-geometric}
 \prod_{t=n}^{n+N-1}v(t)\le C_v\rho_v^N
 \qquad(n\ge1,\ N\ge1),\qquad C_v\ge1,\quad 0<\rho_v<1.
\end{equation}
The constants are independent of the starting row and the length.
Define
\begin{equation}\label{stab:eq:relative-weights}
 \gamma_i^-(n)\coloneq\frac{a_i(n)}{a_{i+1}(n)},\quad
 \gamma_j^+(n)\coloneq\frac{b_j(n)}{b_{j+1}(n)},\quad
 \gamma_0(n)\coloneq a_p(n)b_q(n)\frac{d_{n-1}}{d_n},
 \qquad n\ge1,
\end{equation}
where $1\le i<p$ and $1\le j<q$.
The first two ratios compare neighboring factors at the same row.
The third also compares the diagonal scales at two consecutive rows.
For $p=q=1$, put $\delta d_n=d_n z_n$ in the pivot variation
preceding Theorem~\ref{stab:thm:general-full-stability}. The coefficient
carrying $z_{n-1}$ into $z_n$ becomes
$a_1(n)b_1(n)d_{n-1}/d_n$, which explains the diagonal ratio.

\begin{theorem}[Relative criterion in arbitrary finite bandwidth]
\label{stab:thm:relative-criterion}
Fix the entries $a_i(n)$ for $1\le n\le p-i$ and $b_j(n)$ for
$1\le n\le q-j$. The following conditions are equivalent.
\begin{enumerate}
\item In a neighborhood of $T_0$ in \eqref{stab:eq:relative-norm}, every
real matrix on the same band has an ordered positive bidiagonal
factorization with these initial entries fixed. The factors equal
the prescribed factors at $T_0$ and depend real analytically on the
matrix in relative factor norms.
\item The conditions $\mathcal G(\gamma_i^-)$, $\mathcal G(\gamma_j^+)$, and
$\mathcal G(\gamma_0)$ hold for all the sequences in
\eqref{stab:eq:relative-weights}.
\end{enumerate}
The positive factors with these initial entries are unique. On a sufficiently
small neighborhood there is $C<\infty$ such that
\begin{equation}\label{stab:eq:relative-lipschitz}
 \max\left\{
 \sup_{i,n}\left|\frac{\widetilde a_i(n)-a_i(n)}{a_i(n)}\right|,
 \sup_{j,n}\left|\frac{\widetilde b_j(n)-b_j(n)}{b_j(n)}\right|,
 \sup_n\left|\frac{\widetilde d_n-d_n}{d_n}\right|
 \right\}
 \le C\|\widetilde T-T_0\|_{T_0,\mathrm{rel}}.
\end{equation}
A side of bandwidth one has no adjacent conditions.
\end{theorem}

The proof separates triangular extraction, controlled by the adjacent
ratios, from the LDU variation, controlled additionally by $\gamma_0$.
Both steps require estimates in the relative norms.

\subsubsection{Relative variations on one triangular side}

Lemma~\ref{stab:lem:weighted-shift} will be applied after dividing
each intermediate equation by a suitable positive sequence.
Its row-sum criterion also applies when the weights are initially
unrestricted: every weight occurs in a row sum, so a uniform row-sum
bound first makes the weights bounded.

\begin{lemma}[Relative triangular criterion]\label{stab:lem:relative-lower}
For $F_0=\LL(a_1)\cdots\LL(a_p)$ with arbitrary positive entries,
there is real analytic dependence of its positive factors with the
entries $a_i(n)$, $1\le n\le p-i$, fixed, in relative matrix
and factor norms,
if and only if $\mathcal G(a_i/a_{i+1})$ holds for every $i<p$.
\end{lemma}

\begin{proof}
For $p=1$ the matrix is its own factor. Suppose $p\ge2$.
Consider one extraction of bandwidth $m\ge2$, $F=\LL(a)G$, with $G=\LL(b_1)\cdots\LL(b_{m-1})$.
The first $m-1$ entries of $a$ are fixed. Put $c_0\coloneq a$ and
$c_j\coloneq b_j$ for $1\le j\le m-1$. The strict subdiagonals of
$F$ have the explicit finite expansion
\[
 F_{n,n-s}
 =\sum_{0\le j_1<\cdots<j_s\le m-1}
       \prod_{\nu=1}^s c_{j_\nu}(n-\nu+1),
 \qquad 1\le s\le\min(m,n).
\]
Choosing $j_\nu$ means using the subdiagonal of that factor;
the other factors contribute their unit diagonal. The increasing
indices reflect the order of the factors in the product. This
formula will justify the comparisons of positive terms below.
Let $H$ have subdiagonal
$h_n$, with $h_n=0$ for $n<m$, and let $R$ be a strictly lower
matrix of bandwidth at most $m-1$. The variation equation is $E=HG+\LL(a)R$.
Put $g_s(n)\coloneq G_{n,n-s}$, $g_0=1$, $E_s(n)\coloneq E_{n,n-s}$,
$\Pi_r(n)\coloneq\prod_{t=0}^{r-1}a(n-t)$, and $\Pi_0\coloneq1$.
Successive elimination of $R_s$, in increasing order of $s$, gives
\begin{equation}\label{stab:eq:relative-remainder}
 R_s(n)=\sum_{r=0}^{s-1}(-1)^r\Pi_r(n)
 \bigl(E_{s-r}(n-r)-h_{n-r}g_{s-r-1}(n-r-1)\bigr)
\end{equation}
for $1\le s\le\min(m-1,n)$.
Define $(\mathscr D_jh)_n\coloneq b_j(n-j)h_n-a(n)h_{n-1}$.
The final subdiagonal equation, where $R_m=0$, is
\begin{equation}\label{stab:eq:relative-scalar}
 \mathscr D_{m-1}\cdots\mathscr D_1h=f,\qquad
 f_n\coloneq\sum_{r=0}^{m-1}(-1)^r\Pi_r(n)E_{m-r}(n-r).
\end{equation}
These are the finite algebraic identities used in
Lemma~\ref{stab:lem:variable-extraction}. To verify the factorization,
the coefficient of $h_{n-r}$ on the left is
$(-1)^r\Pi_r(n)g_{m-r-1}(n-r-1)$. Induction on $m$ follows
by writing $G=\overline G\LL(b_{m-1})$ and using
$g_s(n)=\overline g_s(n)+\overline g_{s-1}(n)b_{m-1}(n-s+1)$.
Consequently \eqref{stab:eq:relative-scalar} holds for arbitrary positive entries.

\emph{Normalization of the first-order equations.}
Starting with the relative unknown $h/a$, normalize each successive
first-order equation by the coefficient of its current unknown.
This gives the positive sequences
\[
 w_0(n)\coloneq a(n),\qquad
 w_j(n)\coloneq a(n)\prod_{\ell=1}^{j}b_\ell(n-\ell),\qquad
 1\le j\le m-1.
\]
Here $w_j(n)$ is used for $n\ge j+1$, so every factor index
in its definition is permitted. In particular
$w_{m-1}(n)=F_{n,n-m}$ for $n\ge m$.
Division of each successive output by its own weight gives
\begin{equation}\label{stab:eq:relative-normalized-first-order}
 \frac{\mathscr D_j(w_{j-1}u)_n}{w_j(n)}
 =u_n-\theta_j(n)u_{n-1},\qquad
 \theta_j(n)\coloneq\prod_{\ell=0}^{j-1}
       \frac{b_\ell(n-\ell-1)}{b_{\ell+1}(n-\ell-1)},
 \quad b_0\coloneq a.
\end{equation}
All these equations are solved for $n\ge m$, with preceding
unknown zero at $n=m-1$. The equality follows
by cancelling $a(n)$ and the factors $b_\ell(n-\ell)$; explicitly
the coefficient is
\[
 \frac{a(n)w_{j-1}(n-1)}{w_j(n)}
 =\frac{a(n-1)}{b_j(n-j)}
     \prod_{\ell=1}^{j-1}
       \frac{b_\ell(n-\ell-1)}{b_\ell(n-\ell)}.
\]
This is exactly the product in \eqref{stab:eq:relative-normalized-first-order}.
Consecutive products of each $\theta_j$ satisfy a geometric bound
whenever all adjacent same-row ratios do: interchange the finite
products and apply their bounds at the shifted starting rows.
For the composition, set $u_0\coloneq h/a$ and
$u_j\coloneq(\mathscr D_j\cdots\mathscr D_1h)/w_j$. The identity
reads $u_j=(I-W_j)u_{j-1}$. All $W_j$ act on the rows
$n\ge m$, with zero at the preceding row. The weighted-shift
lemma gives a bounded inverse for each $I-W_j$. Solving from
$u_{m-1}=f/w_{m-1}$ back to $u_0$ therefore gives a bounded
inverse for the normalized scalar equation.
Since $u_0=h/a$, the required factor estimate follows once
the supremum norm of $f/w_{m-1}$ is bounded by a constant times
$\|E\|_{F,\mathrm{rel}}$. The remainder $R$
must also be bounded relative to $G$.

\emph{Bounds for its forcing term and the remainder.}
Each term is compared with a positive summand of its denominator,
using factor entries at the same row. This avoids comparisons between
rows, where growth is unrestricted.
The adjacent geometric bounds at length one bound each
$b_i(t)/b_{i+1}(t)$ uniformly. Multiplying at most $m-1$ such
bounds gives a uniform number $K\ge1$ such that
$b_i(t)\le Kb_j(t)$ whenever $0\le i<j\le m-1$, with $b_0=a$.
A term in a subdiagonal entry of a bidiagonal product selects
increasing factor indices, at successive rows decreasing by one.
Among $m$ factors, the $t$th selected index in a selection of $s$
indices is at most the $t$th index of the last $s$ factors.
Comparing at the very same row therefore bounds every path term by
$K^s$ times that last-factor term. The last-$s$ term of $G$ is
\[
 \prod_{\nu=1}^s b_{m-s+\nu-1}(n-\nu+1).
\]
It is a positive summand in $G_{n,n-s}$. Using this summand as
a lower bound for the denominator avoids comparing entries at
different rows, whose relative sizes are unrestricted.

This comparison also applies when the first few chosen factors are
repeated copies of $a=b_0$. In
$\Pi_r(n)F_{n-r,n-m}$ there are $m$ entries: $r$ initial copies
of $b_0$, followed by a selection of $m-r$ factors. At every
position its index is at most that of the corresponding term in
$w_{m-1}(n)$. Writing $E_{n-r,n-m}$ as $F_{n-r,n-m}$ times its relative
variation therefore gives uniformly bounded coefficients for these
relative data in $f/w_{m-1}$. Similarly, every path in
$\Pi_r(n)F_{n-r,n-s}$ and in
$\Pi_r(n)a(n-r)g_{s-r-1}(n-r-1)$ is bounded by a constant
times the last-$s$ term of $G_{n,n-s}$, which occurs in that
positive entry. There are $s$ entries at rows $n,n-1,\ldots,n-s+1$.
The selected factor at the $\nu$th row has index at most
$m-s+\nu-1$, since otherwise there would be too few later
factors to complete the selection. Repeated initial copies of
$b_0$ satisfy this inequality as well. Thus every term is at
most $K^s$ times the displayed last-$s$ term. The data enter
through $|E_{n-r,n-s}|\le\|E\|_{F,\mathrm{rel}}F_{n-r,n-s}$
and $|h_{n-r}|\le\|h/a\|_\infty a(n-r)$.
Since the number of terms depends only on $m$,
\begin{equation}\label{stab:eq:relative-extraction-bounds}
 \left\|\frac{f}{w_{m-1}}\right\|_\infty
 \le C_m\|E\|_{F,\mathrm{rel}},\qquad
 \|R\|_{G,\mathrm{rel}}
 \le C_m\left(\|E\|_{F,\mathrm{rel}}+\|h/a\|_\infty\right).
\end{equation}
For instance, $C_m\coloneq m2^mK^m$ bounds the finite sums.
Only entries on strict lower subdiagonals are included in these norms.
The recurrence and \eqref{stab:eq:relative-extraction-bounds} give a
bounded inverse for the derivative of this extraction. Repeat for
each remaining factor to obtain a bounded inverse for the full
relative multiplication derivative.

\emph{Necessity.}
The assumed analytic factor map has bounded derivative in the
relative norms. At the first extraction, group the variations of the
remaining factors into $R$. Differentiating their positive finite
product gives
$\|R\|_{G,\mathrm{rel}}\le(m-1)\max_j\|\delta b_j/b_j\|_\infty$:
each original path contribution is multiplied by the sum of its
relative factor changes, with at most $m-1$ summands. The given
derivative therefore supplies relatively bounded $h,R$ for every
prescribed $E$. Prescribe only the
last subdiagonal variation $E_m(n)=F_{n,n-m}f_n$. It has relative
norm $\|f\|_\infty$, and the normalized equation becomes $(I-W_{m-1})\cdots(I-W_1)(h/a)=f$.
The matrices $W_j$ need not yet define bounded operators on
$\ell^\infty$. Their entrywise inverses are nevertheless given
by finite forward substitution in each row. Each inverse
with zero preceding data $(I-W_j)^{-1}$ has nonnegative entries
and unit diagonal. The inverse of the product is the product of these inverses in
reverse order. All act on the same rows $n\ge m$, and multiplication
by a nonnegative matrix with unit diagonal can only increase a
nonnegative matrix entrywise. Their product therefore dominates
each individual inverse. Taking $f=1$
therefore bounds the row sums of $(I-W_1)^{-1}$ by the bound for
$h/a$. The first-order argument gives $\mathcal G(\theta_1)$,
which is $\mathcal G(a/b_1)$ after shifting the index by one.
The row-sum bound first bounds every weight of $W_1$; the proof
of Lemma~\ref{stab:lem:weighted-shift} then gives geometric
decay. The finitely many weights preceding the rows under
consideration are positive finite numbers and can be included
by increasing the constant. No boundedness of these intermediate
operators was used before obtaining their row-sum bounds.

This reasoning applies at every extraction. More explicitly,
an arbitrary relative variation $R$ of the remaining product gives
$E=\LL(a)R$, with $\|E\|_{F,\mathrm{rel}}\le\|R\|_{G,\mathrm{rel}}$
by positivity. The first-extraction equation has the formal solution
$(h,R)=(0,R)$. Its zero initial conditions and nonzero pivots make
that solution unique. The bounded inverse of the full derivative
therefore supplies bounded relative variations of all the remaining
factors. Induction gives every adjacent condition.

\emph{Analyticity and uniqueness.}
Division by the positive original output entry gives the normalized
multiplication map
\[
 \sum_\pi\lambda_\pi
       \left(\prod_{\nu\in\pi}(1+x_\nu)-1\right),\qquad
 \lambda_\pi>0,\qquad \sum_\pi\lambda_\pi=1.
\]
Here $\pi$ indexes its paths, $\lambda_\pi$ is the fraction of
the original entry contributed by that path, and $x_\nu$ is the
relative change of a selected factor entry.
There are at most $p$ relative variables in a lower path. Thus this
is a continuous polynomial map between the indicated Banach spaces,
with uniform bounds on all derivatives on bounded balls. For
example, on $\|x\|_\infty\le R$ the first derivative has norm
at most $p(1+R)^{p-1}$: differentiate each product, bound its at
most $p$ summands, and use $\sum_\pi\lambda_\pi=1$.
The same argument bounds the higher derivatives, which vanish
above degree $p$.
The analytic inverse function theorem gives the local factor map.
Relative changes smaller than $1/2$ keep all entries positive.
The finite row equations, with their positive pivots and prescribed
initial entries, give uniqueness row by row.
\end{proof}

\subsubsection{The full triangular variation}

The remaining step controls the relative changes of $L,D,U$ when
$T_0$ changes. Signed triangular inverses have nonnegative entries;
comparison of their positive contributions with the original factor
entries will give bounds independent of the row.

Assume the geometric conditions of Theorem~\ref{stab:thm:relative-criterion}.
Set
\[
 \Pi^-(n,k)\coloneq\prod_{t=k+1}^{n}a_p(t),\qquad
 \Pi^+(n,k)\coloneq\prod_{t=k+1}^{n}b_q(t),\qquad n\ge k,
\]
with empty products equal to one. Write
$\Sigma=\diag[(-1)^n]_{n\ge0}$,
$\mathcal Q\coloneq\Sigma L^{-1}\Sigma$, and
$\mathcal Q^+\coloneq\Sigma(U^{\mathsf T})^{-1}\Sigma$.
Triangular inverses are defined using their finite leading sections.

\begin{lemma}[Positive estimates for the relative LDU derivative]
\label{stab:lem:relative-ldu}
Under the geometric conditions in
Theorem~\ref{stab:thm:relative-criterion}, there is a constant
independent of the row indices such that, for every band perturbation
$E$ of finite relative norm, the unique entrywise solution of $E=XDU+LZU+LDY$
satisfies
\[
 \|X\|_{L,\mathrm{rel}}+\sup_n\frac{|Z_{n,n}|}{d_n}+\|Y\|_{U,\mathrm{rel}}
 \le C\|E\|_{T_0,\mathrm{rel}}.
\]
Here $X$ and $Y$ are strictly triangular with bandwidths at most
$p$ and $q$, and $Z$ is diagonal.
\end{lemma}

\begin{proof}
\emph{Triangular inverse bounds.}
The positive inverse product gives
\begin{equation}\label{stab:eq:relative-inverse-domination}
 \Pi^-(n,k)\le \mathcal Q_{n,k}\le C_L\Pi^-(n,k),\qquad
 \Pi^+(n,k)\le \mathcal Q^+_{n,k}\le C_U\Pi^+(n,k).
\end{equation}
To verify the first upper bound, expand $\mathcal Q=\mathcal Q_p\cdots\mathcal Q_1$,
where $(\mathcal Q_i)_{n,k}=\prod_{t=k+1}^n a_i(t)$.
Every term divides the interval into consecutive segments using
the successive sequences $a_1,\ldots,a_p$. Divide by $\Pi^-(n,k)$.
On each segment using $a_i$ with $i<p$, the quotient has a
uniformly geometric bound, since $a_i/a_p$ is the product of the
adjacent ratios. Summing independently over the first $p-1$
segment lengths gives a finite product of geometric sums. The
term using only $a_p$ gives the lower bound. For the upper
estimate, use $U^{\mathsf T}=\LL(b_1)\cdots\LL(b_q)$ and
divide its inverse expansion by $\Pi^+(n,k)$. The same sums,
now for $b_j/b_q$ with $j<q$, give the constant $C_U$.

The adjacent conditions also give a uniform pointwise bound
$a_i(t)\le K a_p(t)$ and $b_j(t)\le K b_q(t)$.
Consequently the finite positive path expansions imply
\begin{equation}\label{stab:eq:relative-triangular-comparison}
 \begin{split}
 L_{n,k}&\le C \Pi^-(n,k),\qquad
 L_{n,k}\le C L_{n,j}\Pi^-(j,k)\quad(k\le j\le n, n-k\le p),\\
 U_{k,n}&\le C \Pi^+(n,k),\qquad
 U_{k,n}\le C U_{j,n}\Pi^+(j,k)\quad(k\le j\le n, n-k\le q).
 \end{split}
\end{equation}
For the second lower inequality, split each path after its first
$n-j$ subdiagonal choices. That prefix is a term of $L_{n,j}$.
Bound the remaining choices at each row by $K a_p$ and sum over
the at most $2^p$ completions. Thus a prefix contributes at most
$2^pK^{j-k}$ times itself times $\Pi^-(j,k)$; summing over
prefixes gives the second lower inequality. For the upper one
apply this calculation to $U^{\mathsf T}$ at indices $(n,k)$
and $(n,j)$. All comparisons use entries at the same row.

Combining the first inequalities in
\eqref{stab:eq:relative-triangular-comparison} with \eqref{stab:eq:relative-inverse-domination}
gives $(\mathcal Q L)_{n,k}\le C\Pi^-(n,k)$ and $(\mathcal Q^+U^{\mathsf T})_{n,k}\le C\Pi^+(n,k)$.
For example, in $(\mathcal Q L)_{n,k}$ the intermediate index
$r$ satisfies $k\le r\le\min(n,k+p)$. Each term is at most
$C\Pi^-(n,r)\Pi^-(r,k)=C\Pi^-(n,k)$, so at most $p+1$
terms occur. The upper convolution has at most $q+1$ terms.

\emph{The central geometric sum.}
The two inverse products meet at an index $k$ with diagonal weight
$d_k$. After division by the output diagonal weight, their product
is controlled by $\mathcal G(\gamma_0)$.
Put $V\coloneq L^{-1}EU^{-1}$, entry by entry, and
$\varepsilon\coloneq\|E\|_{T_0,\mathrm{rel}}$.
Every such entry is a finite sum. Positivity gives
\[
 |V_{i,j}|\le\varepsilon[\mathcal Q T_0(\mathcal Q^+)^{\mathsf T}]_{i,j}
 \le C\varepsilon\sum_{k=0}^{\min(i,j)}d_k\Pi^-(i,k)\Pi^+(j,k).
\]
For $i\ge j$, this last sum is
$d_j\Pi^-(i,j)\sum_{k=0}^j\prod_{t=k+1}^j\gamma_0(t)$.
Indeed, the product of the diagonal ratios telescopes:
\[
 d_k\Pi^-(j,k)\Pi^+(j,k)
 =d_j\prod_{t=k+1}^j
       \left(a_p(t)b_q(t)\frac{d_{t-1}}{d_t}\right).
\]
The term $k=j$ contributes one to the normalized sum, and the
remaining terms have total at most $C_0\rho_0/(1-\rho_0)$,
where $C_0,\rho_0$ are constants for $\mathcal G(\gamma_0)$.
For $i\le j$, instead extract $d_i\Pi^+(j,i)$ and sum to $i$.
Thus
\begin{equation}\label{stab:eq:relative-v-bound}
 |V_{i,j}|\le C\varepsilon
 \begin{cases}
 d_j \Pi^-(i,j),&i\ge j,\\
 d_i \Pi^+(j,i),&i\le j.
 \end{cases}
\end{equation}
The same comparisons in the original finite product show that
\begin{equation}\label{stab:eq:relative-entry-comparison}
 d_jL_{i,j}\le (T_0)_{i,j}\le C d_jL_{i,j}\quad(i\ge j),\qquad
 d_iU_{i,j}\le (T_0)_{i,j}\le C d_iU_{i,j}\quad(i\le j)
\end{equation}
on the allowed band. For example, when $i\ge j$, expand
$(T_0)_{i,j}=\sum_{k\le j}L_{i,k}d_kU_{k,j}$, taking entries
outside the bands to be zero. For each nonzero term use
$L_{i,k}\le C L_{i,j}\Pi^-(j,k)$ and
$U_{k,j}\le C\Pi^+(j,k)$. The telescoping identity above bounds
the sum by $C d_jL_{i,j}$ times a uniformly bounded geometric sum.
The term $k=j$ supplies the lower bound.

\emph{Solving while retaining the cancellations.}
The triangular parts of $V$ determine the unique formal solution: $X=LV_-D^{-1}$, $Z=V_0$, and $Y=D^{-1}V_+U$.
To recover the bandwidth and compare directly with $L_{i,j}$,
retain the cancellations in $LV=EU^{-1}$ before taking absolute
values. The equivalent expression is
$X=[EU^{-1}-L(V_0+V_+)]D^{-1}$.
When $i>j$, both sums contributing to $X_{i,j}$ involve only
$\max(0,i-p)\le k\le j$. Thus they vanish for $i-j>p$.
In the first sum, \eqref{stab:eq:relative-entry-comparison} and
\eqref{stab:eq:relative-inverse-domination} bound a term by
$C\varepsilon d_kL_{i,k}\Pi^+(j,k)$. The second sum has the same
bound by \eqref{stab:eq:relative-v-bound}. In both cases
\[
 d_kL_{i,k}\Pi^+(j,k)
 \le C d_jL_{i,j}\prod_{t=k+1}^{j}\gamma_0(t)
 \le C' d_jL_{i,j}.
\]
There are at most $p+1$ terms. Division by $d_j$ proves
$|X_{i,j}|\le C\varepsilon L_{i,j}$.
For the upper side use
$Y=D^{-1}[L^{-1}E-(V_0+V_-)U]$. If $i<j$, intermediate
indices in either sum lie between $\max(0,j-q)$ and $i$;
there are none when $j-i>q$. A term is bounded by
$C\varepsilon d_k\Pi^-(i,k)U_{k,j}$. The comparison
$U_{k,j}\le C U_{i,j}\Pi^+(i,k)$ changes this to
$C\varepsilon d_iU_{i,j}\prod_{t=k+1}^i\gamma_0(t)$.
There are at most $q+1$ terms. Their products are uniformly
bounded, including the empty one. Division by $d_i$ gives
$|Y_{i,j}|\le C\varepsilon U_{i,j}$ and proves its bandwidth.
Finally the diagonal case of
\eqref{stab:eq:relative-v-bound} gives $|Z_{n,n}|\le C\varepsilon d_n$.
Uniqueness follows by splitting the equation for $V$ into its three
triangular parts. This proves the lemma.
\end{proof}

\subsubsection{Proof of the full criterion and the bounded specialization}

\begin{proof}[Proof of Theorem~\ref{stab:thm:relative-criterion}]
\emph{Sufficiency.}
Lemma~\ref{stab:lem:relative-ldu} gives bounded relative variations of $L,D,U$
for every relative variation of $T_0$. Lemma~\ref{stab:lem:relative-lower} then
gives bounded relative variations of the individual bidiagonal
factors on both sides, with the prescribed initial variations zero.
Thus the derivative of full multiplication has a bounded inverse.

As in the lower lemma, division of every output by its positive
base entry expresses the full multiplication map as positive weighted
averages of monomials in the relative variables. Each monomial has
degree at most $p+q+1$. The map is therefore analytic between the
relative Banach spaces. The analytic inverse function theorem applies.
On a smaller ball all relative factor changes are smaller than $1/2$,
so positivity holds. The inverse has continuous derivative, bounded
on a sufficiently small ball, which gives \eqref{stab:eq:relative-lipschitz}.
Uniqueness of normalized LDU on every finite leading section,
followed by the rowwise triangular uniqueness on each side, gives
uniqueness of the positive factors with those initial entries.

\emph{Necessity on the triangular sides.}
Vary $T=\widetilde LDU$ with $\widetilde L$ near $L$ in relative
norm. Positivity of $D,U$ gives $\|(\widetilde L-L)DU\|_{T_0,\mathrm{rel}}\le\|\widetilde L-L\|_{L,\mathrm{rel}}$.
Normalized LDU uniqueness forces the lower product of the chosen
PBF to be $\widetilde L$. It therefore gives the analytic relative
factor map required in Lemma~\ref{stab:lem:relative-lower}. This proves all
$\mathcal G(\gamma_i^-)$. Varying the upper part proves all
$\mathcal G(\gamma_j^+)$.

\emph{Necessity across the diagonal.}
Choose the diagonal direction $E=D$. It belongs to the relative
perturbation space,
since $(T_0)_{n,n}\ge d_n$. Its diagonal factor derivative satisfies
\[
 \frac{Z_{n,n}}{d_n}
 =\frac1{d_n}\sum_{k=0}^{n}\mathcal Q_{n,k}\mathcal Q^+_{n,k}d_k
 \ge\sum_{k=0}^{n}\prod_{t=k+1}^{n}\gamma_0(t).
\]
The signs of the two triangular inverses cancel. The lower bound
uses the terms contributed by their last bidiagonal factors and
holds for every positive PBF. The analytic derivative bounds the
left side uniformly. The first-order row-sum argument therefore
gives $\mathcal G(\gamma_0)$, completing the necessity proof.
\end{proof}

\begin{corollary}[Recovery of the operator-norm criterion]
\label{stab:cor:relative-bounded}
If every $a_i(n),b_j(n),d_n$ belongs to one interval
$[\alpha,A]$ with $0<\alpha\le A<\infty$, the relative theorem is
equivalent to Theorem~\ref{stab:thm:general-full-stability}.
On either triangular side, Lemma~\ref{stab:lem:relative-lower}
specializes to Theorem~\ref{stab:thm:general-lower-stability}.
\end{corollary}
\begin{proof}
Every permitted entry of $T_0$ is uniformly bounded above and
below by positive constants, since it is a fixed finite sum of
positive path products and at least one path is present. Thus the
relative matrix norm is equivalent to the coefficient norm and,
on the fixed band, to operator norm. Relative and absolute factor
supremum norms are also equivalent.

For $v\ge u\ge2$, the products in the two criteria satisfy
\[
 \prod_{n=u}^{v}\frac{a_i(n)}{a_{i+1}(n-1)}
 =\frac{a_{i+1}(v)}{a_{i+1}(u-1)}\prod_{n=u}^{v}\gamma_i^-(n),\qquad
 \prod_{n=u}^{v}\gamma_0(n)
 =\frac{d_{u-1}}{d_v}\prod_{n=u}^{v}a_p(n)b_q(n).
\]
Every endpoint quotient is between $\alpha/A$ and $A/\alpha$.
Hence each geometric condition is equivalent to its counterpart
in Theorem~\ref{stab:thm:general-full-stability}. The upper ratios obey the same identity.
The finitely many initial weights are absorbed into the constants.
\end{proof}

\begin{example}[Factors with different growth rates]
\label{stab:ex:relative-different-growth}
Take $p=q=2$ and put $a_1(n)\coloneq n$, $a_2(n)\coloneq2n^2$, $b_1(n)\coloneq n$, $b_2(n)\coloneq2n^3$, and $d_n\coloneq8^n(n!)^5$.
The criterion weights are $\gamma_1^-(n)=\frac1{2n}$, $\gamma_1^+(n)=\frac1{2n^2}$, and $\gamma_0(n)=\frac12$.
All consecutive products are at most $2^{-N}$, so
Theorem~\ref{stab:thm:relative-criterion} gives positive factors
depending analytically on every sufficiently small relative band
perturbation. Every bidiagonal factor is unbounded. Under a common
diagonal conjugation of the lower factors, the ratio
$a_2(n)/a_1(n)=2n$ is preserved. Thus these factors have different
growth rates even after that normalization.

Direct multiplication gives
\[
 \frac{(T_0)_{n,n}}{d_n}
 =1+\frac{(n+2n^2)(n+2n^3)}{8n^5}+\frac1{16n^3},\qquad n\ge2.
\]
In particular $d_n/(T_0)_{n,n}\to2/3$ as $n\to\infty$.
The contribution of the other factor paths remains comparable
to $d_n$, and the weight $\gamma_0=1/2$ controls the propagation
of their relative variations between rows.
\end{example}

\section[Toeplitz matrices and mixed-type polynomials]{Banded Toeplitz matrices and mixed-type orthogonal polynomials}
\label{sec:mixed-systems}

Global PBFs with bounded factors for BTP Toeplitz matrices follow
from shifted Metelmann factors and from the negative roots of their
polynomial symbols. For distinct roots, the construction also gives
analytic dependence of the factors
(Theorem~\ref{thm:toeplitz-stable-simple-roots}); its exchanges extend
to upper--lower products with distinct limits
(Corollary~\ref{cor:asymptotic-ul-stability}). The chosen factors specify
the polynomial initial conditions, and the bounded spectral Favard
theorem \cite{BFM2023} supplies their matrix of positive measures.
Explicit formulas determine these polynomials, the measure, and its
endpoint behavior, including repeated roots. Although the chosen factors determine the
polynomial vectors and their measure, the associated weighted linear
forms are independent of that choice. Section~\ref{subsec:spectral-kernels}
develops the Chebyshev convolutions and finite and Poisson kernels.
The density and explicit oscillatory amplitudes then give interior
zero spacing and the local sine-kernel limit.

Section~\ref{subsec:determinants-equilibrium} relates determinants
of the polynomial vectors to Schur polynomials and to the vector
equilibrium problem of Duits and Kuijlaars. Section~\ref{subsec:quasi-toeplitz}
treats reordered constant factors with a common Toeplitz tail,
their factorization stability, and the polynomials associated with
finite changes in the initial block, together with their measures
under the bounded PBF hypothesis.

The binomial and geometric families give explicit factor formulas,
verified by finite summation identities.
The geometric family is related to that of any ordered partition
of the same roots by
\eqref{eq:q-change-polynomials}--\eqref{eq:q-change-measure}.
At $q=1$, \eqref{eq:lima-rescaling} recovers Lima's system from the
binomial family; Proposition~\ref{prop:q-affine-obstruction} explains
why the geometric parameter cannot generally be absorbed into a
change of variable and positive diagonal normalization.

\subsection{Root criterion and global positive factorization}\label{sec:toeplitz}

A BTP Toeplitz matrix admits two constructions of a global PBF.
The first uses shifted Metelmann factors; the second uses the negative
roots of its polynomial symbol. Both produce an upper--lower product.

The term ``global'' means that the factors are defined on every
row of the semi-infinite matrix, with compatible entries as the
finite leading sections increase.
The exchange lemma below puts this product in the required
lower--diagonal--upper order, with positive entries in every row.

Let
\begin{equation}\label{eq:general-toeplitz}
 T\coloneq [t_{i-j}]_{i,j\in\NN},\qquad
 P(z)\coloneq \sum_{k=-q}^{p}t_k z^{k+q},\qquad a(z)\coloneq z^{-q}P(z),
\end{equation}
where $p,q\in\mathbb N$, $t_k\in\RR$, $t_k=0$ for
$k\notin\{-q,\dots,p\}$, and $t_{-q},t_p>0$.
These matrices are bounded because each is a finite linear
combination of powers of the upper and lower shifts.

More explicitly, the upper shift acts by $(\Lambda v)_n=v_{n+1}$,
whereas $(\Lambda^{\mathsf T}v)_0=0$ and
$(\Lambda^{\mathsf T}v)_n=v_{n-1}$ for $n\ge1$. Both have operator
norm one on $\ell^2$. Hence
\[
 T=t_0I+\sum_{k=1}^{p}t_k(\Lambda^{\mathsf T})^k
       +\sum_{k=1}^{q}t_{-k}\Lambda^k,
 \qquad \|T\|_2\le\sum_{k=-q}^{p}|t_k|.
\]
The Laurent polynomial $a$ records these diagonal coefficients:
the coefficient of $z^k$ is the entry on diagonal $i-j=k$.
Its ordinary polynomial multiple $P=z^qa$ will be used to locate
roots without a pole at the origin.

The following lemma is a version of \cite[Lemma~2.2]{ManasPBF2026}
with a variable positive superdiagonal; the cited result uses $v_n=1$.
The proof below specifies the normalization used here.

\begin{lemma}[Positive upper--lower exchange]\label{lem:exchange}
Let $U$ be upper bidiagonal with diagonal $u_n>0$ and
superdiagonal $v_n>0$, and let $L$ be unit lower bidiagonal with
subdiagonal $\ell_n>0$. There is a unique decomposition
$UL=L'U'$ with $L'$ unit lower bidiagonal and $U'$ upper
bidiagonal. Its entries are
\begin{align}
 d_0&=u_0+v_0\ell_1,\label{eq:swap-initial}\\
 d_n&=u_n+v_n\ell_{n+1}
       -\frac{u_n\ell_n v_{n-1}}{d_{n-1}},\qquad n\in\mathbb N,
       \label{eq:swap-diagonal}\\
 L'_{n,n-1}&=\frac{u_n\ell_n}{d_{n-1}},\qquad
 U'_{n,n}=d_n,\qquad U'_{n,n+1}=v_n.\label{eq:swap-factors}
\end{align}
All displayed parameters are positive.
\end{lemma}
\begin{proof}
Write $c_n\coloneq L'_{n,n-1}$ and $d_n\coloneq U'_{n,n}$. Direct multiplication gives
\[
 \begin{aligned}
 (UL)_{n,n-1}&=u_n\ell_n,&
 (L'U')_{n,n-1}&=c_nd_{n-1},\qquad n\in\mathbb N,\\
 (UL)_{n,n}&=u_n+v_n\ell_{n+1},&
 (L'U')_{n,n}&=d_n+c_nv_{n-1},\qquad n\in\mathbb N.
 \end{aligned}
\]
Both upper diagonals are $v_n$, and the entry $(0,0)$ gives
$d_0=u_0+v_0\ell_1$. Thus $c_n=u_n\ell_n/d_{n-1}$, and substitution
in the diagonal equation gives \eqref{eq:swap-diagonal}. These
relations determine the factors uniquely whenever their denominators
are nonzero.

There are no further entries to check: each of the two products is
tridiagonal. Equality on the lower, main and upper diagonals is
therefore equality of the whole matrix. The row-zero equation is
listed separately because no term involving row $-1$ is present.

To verify this and positivity simultaneously, prove the stronger
inequality $d_n>v_n\ell_{n+1}$. It holds at $n=0$ because $u_0>0$.
If it holds at $n-1$, then
\[
 d_n=v_n\ell_{n+1}
       +u_n\left[1-\frac{v_{n-1}\ell_n}{d_{n-1}}\right]
       >v_n\ell_{n+1}>0.
\]
This also makes every denominator and lower parameter positive.
The recursion consequently constructs all rows with positive
parameters determined by the given semi-infinite factors.
\end{proof}

\begin{theorem}[Toeplitz characterization and global PBF]
\label{thm:toeplitz}
Let $T=[t_{i-j}]_{i,j\in\NN}$ be a real semi-infinite Toeplitz
matrix with $t_k=0$ for $k\notin\{-q,\dots,p\}$, where
$p,q\in\mathbb N$ and $t_{-q},t_p>0$. Define
$P(z)\coloneq\sum_{k=-q}^{p}t_kz^{k+q}$.
Then the following conditions are equivalent:
\begin{enumerate}
\item $T$ is TN.
\item $T$ is BTP.
\item All zeros of $P$ are real and negative.
\item $T$ has a semi-infinite $(p,q)$ PBF.
\end{enumerate}
In condition (3), multiplicities are unrestricted. Equivalently,
\begin{equation}\label{eq:root-parameters}
 P(z)=t_{-q}\prod_{\nu=1}^{p+q}(1+\beta_\nu z),
 \qquad \beta_\nu>0.
\end{equation}
In condition (4), the factors can be chosen bounded, with their allowed positive
entries bounded away from zero.
\end{theorem}
\begin{proof}
Put $m\coloneq p+q$, $\omega\coloneq t_{-q}$ and
$c_r\coloneq t_{r-q}/\omega$, with $c_r=0$ outside
$\{0,\dots,m\}$. The coefficient matrix
$C\coloneq[c_{i-j}]_{i,j\in\NN}$ is unit lower triangular of
bandwidth $m$. For finite increasing index lists $I,J$,
\begin{equation}\label{eq:toeplitz-coefficient-minors}
 C[I,J]=\omega^{-1}T[I,J+q],\qquad
 T[I,J]=\omega C[I+q,J],
\end{equation}
where adding $q$ to a list means adding it to every index.
In particular, $C$ is TN if and only if $T$ is TN.

For a minor of order $r$, the two identities multiply its determinant
by $\omega^{-r}$ and $\omega^r$, respectively. Since $\omega>0$,
they preserve the sign of every minor and justify both implications. If $T$ is BTP,
then $C$ is BTP: the admissibility inequalities
$0\le i_a-j_a\le m$ become
$-q\le i_a-(j_a+q)\le p$ in the first identity.

\emph{Direct construction from BTP.}
Assume condition (2). The Metelmann factors of $C$ are obtained
by applying Proposition~\ref{prop:metelmann-input} to the leading
sections of $C(I+\Lambda)$ and using their compatibility from
Section~\ref{sec:limit}. Here $\Lambda$ is the upper shift.
This product is BTP of bandwidth $(m,1)$. Indeed, both factors are
TN. For an admissible minor with row list $I$ and column list $J$,
the list $K\coloneq (\min(i_a,j_a))_a$ is strictly increasing and satisfies
$0\le i_a-\min(i_a,j_a)\le m$ and
$0\le j_a-\min(i_a,j_a)\le1$.
Thus the Cauchy--Binet term with intermediate list $K$ is positive:
$\det C[I,K]>0$ and $\det(I+\Lambda)[K,J]=1$.
All other terms are nonnegative. Here the list of minima is strictly
increasing: increasing each of two integer entries strictly also
increases their minimum strictly. In the upper bidiagonal minor,
any admissible matching is forced to preserve order, and its
nonzero entries are all one; its determinant is therefore one.
Only finitely many intermediate lists can contribute because both
factors are banded. Thus Cauchy--Binet is a finite sum for this minor.
The normalized triangular decomposition of this product is $C\,I\,(I+\Lambda)$, so uniqueness
identifies $C$ with its lower Metelmann product. Consequently
$C=\widehat L_1\cdots\widehat L_m$. Write
$\widehat a_i(n)\coloneq(\widehat L_i)_{n,n-1}$,
where $\widehat a_i(n)=0$ for $n\le m-i$ and
$\widehat a_i(n)>0$ for $n>m-i$.

Delete the first $m$ rows and columns of each factor and reindex
from zero. Denote the resulting unit lower bidiagonal matrices by
$B_i$, so their subdiagonal entries are
$b_i(n)\coloneq\widehat a_i(n+m)>0$ for $n\ge1$.
Their product is still $C$:
\begin{equation}\label{eq:toeplitz-shifted-metelmann}
 C=B_1\cdots B_m.
\end{equation}
To verify this, consider an entry $(u+m,v+m)$ with $u,v\in\NN$.
In a product of lower triangular matrices, every contributing
intermediate index lies between $v+m$ and $u+m$, hence is at least
$m$. Deleting indices below $m$ therefore leaves the entry unchanged.
After reindexing, the product has entry $C_{u+m,v+m}=C_{u,v}$
in position $(u,v)$, since $C$ is Toeplitz. This proves
\eqref{eq:toeplitz-shifted-metelmann} and removes all the initial
zeros without changing the matrix.

The first and last subdiagonals of the original product give
$c_1=\sum_{i=1}^{m}\widehat a_i(n)$ for $n\ge1$ and $c_m=\prod_{i=1}^{m}\widehat a_i(N-i+1)$ for $N\ge m$.
Here $c_1,c_m>0$. The sum bounds every entry above by $c_1$.
To bound $\widehat a_i(t)$ below for $t\ge m$, take $N=t+i-1$
in the product. Every other factor has an index at least one and
is at most $c_1$. Hence
\begin{equation}\label{eq:toeplitz-shifted-bounds}
 \frac{c_m}{c_1^{m-1}}\le b_i(n)\le c_1,
 \qquad i\in\{1,\dots,m\},\quad n\ge1.
\end{equation}

The second identity in \eqref{eq:toeplitz-coefficient-minors} gives
$T=\omega\Lambda^q C$. For each $i\in\{1,\dots,q\}$, define
an upper bidiagonal matrix $V_i$ by
$(V_i)_{n,n}\coloneq b_i(n+q-i+1)$ and $(V_i)_{n,n+1}\coloneq 1$ for $n\in\NN$.
Entrywise multiplication shows
$\Lambda^{q-i+1}B_i=V_i\Lambda^{q-i}$.
Applying these identities successively yields
\begin{equation}\label{eq:toeplitz-metelmann-ul}
 T=\omega V_1\cdots V_q B_{q+1}\cdots B_m.
\end{equation}
All diagonal entries of the $V_i$ and all subdiagonal entries of
the $B_i$ satisfy the uniform positive bounds
\eqref{eq:toeplitz-shifted-bounds}.

\emph{Positive exchanges and diagonal normalization.}
In \eqref{eq:toeplitz-metelmann-ul}, take the lower factors in
their left-to-right order and move each across all upper factors
to its left, starting with the nearest one, using
Lemma~\ref{lem:exchange}. The relative order within each group is
preserved. There are $pq$ upper--lower pairs to exchange; afterwards
all lower factors are on the left. All parameters remain positive.
Write each resulting upper factor as a positive diagonal times a
unit upper bidiagonal. Moving these diagonals to the left by
diagonal conjugation gives \eqref{eq:strict-definition}, with
$\omega$ included in $D$.

These operations also preserve uniform bounds. Every upper factor
before diagonal normalization has superdiagonal one. If its diagonal
$u_n$ and the subdiagonal $\ell_n$ of the lower factor being exchanged
satisfy $0<u_-\le u_n\le u_+$ and
$0<\ell_-\le\ell_n\le\ell_+$, Lemma~\ref{lem:exchange} gives
$\ell_-<d_n\le u_++\ell_+$ and $u_-\ell_-/(u_++\ell_+)\le L'_{n,n-1}<u_+$.
Indeed, its proof gives $d_n>\ell_{n+1}$, while dropping the
nonnegative subtracted term in \eqref{eq:swap-diagonal} gives the
upper bound. Substitution in $L'_{n,n-1}=u_n\ell_n/d_{n-1}$
gives the remaining inequalities. Starting from
\eqref{eq:toeplitz-shifted-bounds}, these estimates apply through
all $pq$ exchanges. Diagonal normalization involves only finite
products, inverses and adjacent-entry ratios of positive sequences
bounded above and away from zero. It preserves both bounds.

For example, moving a diagonal $E\coloneq \operatorname{diag}(e_n)$ to the
left of a unit upper factor $V$ uses
$VE=E(E^{-1}VE)$. Its superdiagonal entry $v_n$ becomes
$v_ne_{n+1}/e_n$. If $e_-\le e_n\le e_+$ with $e_->0$, this
ratio lies between $e_-/e_+$ and $e_+/e_-$. The analogous estimate
for every extracted diagonal proves the asserted bounds after the
finite normalization procedure.
This proves (2)$\Rightarrow$(4), including the asserted boundedness,
directly from BTP. Conversely, Proposition~\ref{prop:positive-product}
gives (4)$\Rightarrow$(2). Also, (2)$\Rightarrow$(1), since the
remaining minors are forced to vanish by the band pattern.

\emph{The root criterion.}
A sequence whose coefficient Toeplitz matrix is TN is called a
P\'olya-frequency sequence. The finite-sequence
Aissen--Schoenberg--Whitney theorem \cite{ASW} states that a finite
nonnegative sequence has this property if and only if its generating
polynomial has only real nonpositive zeros. Here the generating
polynomial of $(c_r)$ is $P(z)/\omega$. If $T$ is TN, its entries
and hence the $c_r$ are nonnegative, so the theorem applies.
Conversely, negative roots and $P(0)=\omega>0$ give
\eqref{eq:root-parameters}, whose expansion has nonnegative
coefficients; the same theorem makes $C$, and hence $T$, TN.
The nonzero constant term excludes zero roots, and $t_p>0$ fixes
the degree at $m$. Thus (1) and (3) are equivalent, including
repeated roots.

\emph{Alternative construction from the roots.}
Assume (3). The roots give the exact upper--lower product
\begin{equation}\label{eq:positive-ul-toeplitz}
 T=t_{-q}
 \prod_{\nu=1}^{q}(\Lambda+\beta_\nu I)
 \prod_{\nu=q+1}^{p+q}(I+\beta_\nu\Lambda^{\mathsf T}).
\end{equation}
To verify it, expand the first block in powers of $\Lambda$
and the second in powers of $\Lambda^{\mathsf T}$. Every term
has all upper shifts on the left. Since
$\Lambda\Lambda^{\mathsf T}=I$,
\[
 \Lambda^r(\Lambda^{\mathsf T})^s=
 \begin{cases}
  \Lambda^{r-s},&r\ge s,\\
  (\Lambda^{\mathsf T})^{s-r},&s\ge r.
 \end{cases}
\]
Under $\Lambda\leftrightarrow z^{-1}$ and
$\Lambda^{\mathsf T}\leftrightarrow z$, the expanded product has
Laurent polynomial
\[
 t_{-q}\prod_{\nu=1}^{q}(z^{-1}+\beta_\nu)
       \prod_{\nu=q+1}^{p+q}(1+\beta_\nu z)
 =z^{-q}P(z)=\sum_{k=-q}^{p}t_kz^k.
\]
Its coefficient $t_k$ occupies the diagonal $i-j=k$, proving
\eqref{eq:positive-ul-toeplitz} entry by entry.
Apply the same exchange procedure and diagonal normalization to
this product. The bounds now start from the positive constants
$\beta_\nu$. This proves (3)$\Rightarrow$(4) and completes the
equivalence. It also gives a second construction under condition
(2), through the root criterion. This second construction records the
chosen ordered partition of the roots and will be used below.
\end{proof}

\begin{remark}[Factor order and bidiagonal entries]
The product in \eqref{eq:positive-ul-toeplitz} is in upper--lower
order, where $\Lambda\Lambda^{\mathsf T}=I$.
In the reverse order,
$\Lambda^{\mathsf T}\Lambda=I-e_0e_0^{\mathsf T}$.
The positive exchanges account for this initial-row term and give
a lower--diagonal--upper PBF whose bidiagonal entries generally vary
with the row index.
\end{remark}

\subsubsection{Simple roots and analytic dependence of factors}
The root construction also selects factors to which the stability
criterion of Section~\ref{chap:stability} applies. The relevant choice
puts the larger root parameters in the upper group and orders both
groups increasingly before the exchanges. The limits through the construction follow from the behavior of a
single exchange.

\begin{lemma}[Limits under positive exchange]
\label{lem:exchange-separated-limits}
In Lemma~\ref{lem:exchange}, suppose $v_n=1$ for every $n\ge0$,
$u_n\to u>0$ and $\ell_n\to\ell>0$ as $n\to\infty$, with $u\ne\ell$.
Then
$U'_{n,n}\to\max(u,\ell)$ and
$L'_{n,n-1}\to\min(u,\ell)$ as $n\to\infty$.
If $u<\ell$, put $\delta_n\coloneq U'_{n,n}-\ell_{n+1}$. Then
\begin{equation}\label{eq:exchange-swapped-limits}
 U'_{n,n}=\ell_{n+1}+\delta_n,\qquad
 L'_{n,n-1}=u_n-\delta_n,\qquad n\ge1,
\end{equation}
where $\delta_n>0$ and, for every $\theta\in(u/\ell,1)$,
$\delta_n=\mathrm O(\theta^n)$ as $n\to\infty$.
\end{lemma}
\begin{proof}
Write $d_n=U'_{n,n}$. The proof of Lemma~\ref{lem:exchange}
gives $d_n>\ell_{n+1}$. Set
$\delta_n\coloneq d_n-\ell_{n+1}>0$ and $z_n\coloneq 1/\delta_n$.
Equation~\eqref{eq:swap-diagonal} becomes
\begin{equation}\label{eq:exchange-reciprocal-excess}
 \delta_n=\frac{u_n\delta_{n-1}}{\ell_n+\delta_{n-1}},\qquad
 z_n=\frac{\ell_n}{u_n}z_{n-1}+\frac1{u_n},\qquad n\ge1.
\end{equation}
First suppose $u>\ell$. The limiting affine recurrence for $z_n$
has the positive fixed point $z_*\coloneq1/(u-\ell)$.
Put $e_n\coloneq z_n-z_*$ and
$f_n\coloneq (\ell_n/u_n)z_*+1/u_n-z_*$.
Then $f_n\to0$ as $n\to\infty$ and
$e_n=(\ell_n/u_n)e_{n-1}+f_n$.
Choose $\theta\in(\ell/u,1)$ and an integer $N\ge1$ such that
$\ell_n/u_n\le\theta$ for every $n>N$; this is possible because
$\ell_n/u_n\to\ell/u<\theta$ as $n\to\infty$.
Iterating from row $N$ gives
\[
 |e_n|\le\theta^{n-N}|e_N|
       +\sum_{j=N+1}^{n}\theta^{n-j}|f_j|,
 \qquad n\ge N.
\]
Given $\varepsilon>0$, increase $N$ until
$|f_j|\le(1-\theta)\varepsilon$ for $j>N$.
The sum is then at most $\varepsilon$, and the first term tends
to zero as $n\to\infty$. Hence $z_n\to z_*$,
$\delta_n\to u-\ell$, and $d_n\to u$, as $n\to\infty$.
Equation~\eqref{eq:swap-factors} gives
$L'_{n,n-1}=u_n\ell_n/d_{n-1}\to\ell$ as $n\to\infty$.

Now suppose $u<\ell$ and choose $\theta\in(u/\ell,1)$.
Choose $N_\theta\ge1$ such that $u_n/\ell_n\le\theta$ for
every $n>N_\theta$, using $u_n/\ell_n\to u/\ell<\theta$ as
$n\to\infty$. The first recurrence in
\eqref{eq:exchange-reciprocal-excess} then gives
$0<\delta_n\le\theta\delta_{n-1}$ for $n>N_\theta$.
Hence $\delta_n\le\delta_{N_\theta}\theta^{n-N_\theta}$ for
$n\ge N_\theta$, proving the stated geometric bound.
The first identity in \eqref{eq:exchange-swapped-limits} is its
definition. For the second, substitute
$d_{n-1}=\ell_n+\delta_{n-1}$ in \eqref{eq:swap-factors}:
\[
 L'_{n,n-1}
 =\frac{u_n\ell_n}{\ell_n+\delta_{n-1}}
 =u_n-\frac{u_n\delta_{n-1}}{\ell_n+\delta_{n-1}}
 =u_n-\delta_n.
\]
Thus $d_n\to\ell$ and $L'_{n,n-1}\to u$ as $n\to\infty$.
The geometric estimate measures the error relative to the input
entries in \eqref{eq:exchange-swapped-limits}; their convergence to
$u,\ell$ determines the remaining error relative to these constants.
\end{proof}

The next theorem treats analytic choices of positive factors with fixed
or variable initial entries. In both cases, the matrix varies in
operator norm within its band and the factors vary in their supremum
norms.

Thus a neighborhood includes perturbations of all entries on the
allowed diagonals, not just perturbations that preserve the Toeplitz
pattern. A real analytic choice is a map from this open set of
matrices to the bounded sequences of factor entries that has a
convergent local power-series expansion in the corresponding norms. The existence of a PBF that extends to such an analytic choice
is a property of the Toeplitz matrix itself.

\begin{theorem}[Analytic PBF dependence and distinct roots]
\label{thm:toeplitz-stable-simple-roots}
Let $T=[t_{i-j}]_{i,j\in\NN}$ be a BTP Toeplitz matrix of bandwidth
$(p,q)$, with $p,q\ge1$ and $t_{-q},t_p>0$. Write
$P(z)=\omega\prod_{\nu=1}^{m}(1+\beta_\nu z)$, where
$\omega\coloneq t_{-q}$, $m\coloneq p+q$ and
$\beta_1\ge\cdots\ge\beta_m>0$.
The following conditions are equivalent:
\begin{enumerate}
\item $\beta_1>\cdots>\beta_m$.
\item There exist a PBF of $T$, whose positive factor entries are
bounded above and bounded away from zero uniformly over all factors
and rows, and $\varepsilon>0$ with the following property.
Every bounded real matrix $\widetilde T$ of bandwidth at most $(p,q)$
satisfying $\|\widetilde T-T\|_2<\varepsilon$ admits a PBF whose
factors can be chosen to depend real analytically on $\widetilde T$,
in the supremum norms of their entry sequences.
At $\widetilde T=T$, they agree with the chosen factors of $T$.
The initial entries of both triangular sides, as specified in
Lemma~\ref{lem:rowwise}, remain fixed at their values in this PBF.
\item There exists a PBF of $T$ with positive factor entries uniformly
bounded above and away from zero, for which the same assertion holds
with the initial entries allowed to depend on $\widetilde T$.
At $\widetilde T=T$, the chosen factors still agree with this PBF of $T$.
\item $T$ is an interior point of the PBF class, in operator norm
among all bounded real matrices of bandwidth at most $(p,q)$.
\end{enumerate}
Conditions (2) and (3) each assert existence of a suitable PBF of $T$.
For a prescribed PBF, keeping its initial entries fixed and allowing
them to vary are different requirements.
When these conditions hold, the root construction with
\begin{equation}\label{eq:toeplitz-stable-root-order}
 \eta_i\coloneq \beta_{q-i+1}\quad(i\in\{1,\dots,q\}),\qquad
 \xi_j\coloneq \beta_{m-j+1}\quad(j\in\{1,\dots,p\})
\end{equation}
gives such a PBF. In the order \eqref{stab:eq:full-given}, its
entries satisfy, for $i\in\{1,\dots,p\}$ and $j\in\{1,\dots,q\}$,
\begin{equation}\label{eq:toeplitz-stable-factor-limits}
 a_i(n)\to\beta_{m-i+1},\qquad
 b_j(n)\to\beta_j^{-1},\qquad
 d_n\to\omega\prod_{\nu=1}^{q}\beta_\nu,
 \qquad n\to\infty.
\end{equation}
\end{theorem}
\begin{proof}
\emph{Construction for distinct roots.}
Assume (1), and use the order in
\eqref{eq:toeplitz-stable-root-order}. Equation~\eqref{eq:positive-ul-toeplitz}
gives the exact product
$T=\omega\prod_{i=1}^{q}(\Lambda+\eta_i I)\prod_{j=1}^{p}(I+\xi_j\Lambda^{\mathsf T})$.
Both lists $\eta_1,\dots,\eta_q$ and $\xi_1,\dots,\xi_p$ are
increasing, and every $\eta_i$ is larger than every $\xi_j$.
Take the lower factors in order, moving each to the left across
all upper factors as in the proof of Theorem~\ref{thm:toeplitz}.
At each exchange, Lemma~\ref{lem:exchange-separated-limits}
preserves the limiting upper diagonal and lower subdiagonal.

At the first exchange its two limiting values satisfy $\eta_i>\xi_j$.
The conclusion of that lemma leaves these values attached to the
same types of factors. Inductively the next exchange again has an
upper limit from the $\eta$ list and a lower limit from the $\xi$
list, so the strict inequality required for the same conclusion
holds at every step.
After the finitely many exchanges the product is
$T=\omega L_1\cdots L_p V_1\cdots V_q$, where
$[L_j]_{n,n-1}\to\xi_j$ and $[V_i]_{n,n}\to\eta_i$ as
$n\to\infty$. Each $V_i$ still has superdiagonal one by \eqref{eq:swap-factors}.

Put $v_i(n)\coloneq [V_i]_{n,n}$. Extract the diagonal from each $V_i$
and move these diagonals to the left through the preceding upper
factors. The resulting identity
$\omega V_1\cdots V_q=DU_q\cdots U_1$ has entries
\begin{equation}\label{eq:toeplitz-stable-normalization}
 \begin{aligned}
 d_n&=\omega\prod_{h=1}^{q}v_h(n),&& n\ge0,\\
 b_{q-i+1}(n)&=\frac1{v_i(n-1)}
       \prod_{h=i+1}^{q}\frac{v_h(n)}{v_h(n-1)},
       &&n\ge1,\quad i\in\{1,\dots,q\}.
 \end{aligned}
\end{equation}
Indeed, dividing $V_i$ by its diagonal gives superdiagonal
$1/v_i(n-1)$ in position $(n-1,n)$. Moving a later diagonal
with entries $v_h(n)$ to its left multiplies that entry by
$v_h(n)/v_h(n-1)$. This proves
\eqref{eq:toeplitz-stable-normalization} and hence
\eqref{eq:toeplitz-stable-factor-limits}. Every parameter sequence
is positive and tends to a positive finite limit. There are
finitely many such sequences, so all their entries have common
positive lower and finite upper bounds.

As $n\to\infty$, the ratios required by the full stability
 theorem satisfy
\begin{equation}\label{eq:toeplitz-stability-gap-limits}
 \begin{aligned}
 \frac{a_i(n)}{a_{i+1}(n-1)}
 &\longrightarrow\frac{\beta_{m-i+1}}{\beta_{m-i}}<1,
 &&i\in\{1,\dots,p-1\},\\
 \frac{b_j(n)}{b_{j+1}(n-1)}
 &\longrightarrow\frac{\beta_{j+1}}{\beta_j}<1,
 &&j\in\{1,\dots,q-1\},\\
 a_p(n)b_q(n)&\longrightarrow\frac{\beta_{q+1}}{\beta_q}<1,
 &&
 \end{aligned}
\end{equation}
For any positive sequence $w(t)$ tending to a number less than
one, choose $\rho<1$ larger than that limit and $N_0$ such that
$w(t)\le\rho$ for $t\ge N_0$. With
$C\coloneq \prod_{1\le t<N_0}\max\{1,w(t)/\rho\}$, every product of
$N$ consecutive terms is at most $C\rho^N$.
For a sequence starting at index two, use the product from $t=2$.
Applying this observation to the three kinds of sequence in
\eqref{eq:toeplitz-stability-gap-limits} proves
\eqref{stab:eq:adjacent-stability} on both sides and
\eqref{stab:eq:cross-stability}.
Theorem~\ref{stab:thm:general-full-stability} proves (1)$\Rightarrow$(2).
An analytic choice with fixed initial entries is also a choice with
those entries allowed to vary, proving (2)$\Rightarrow$(3).
The choice in (3) supplies a PBF throughout a neighborhood, proving
(3)$\Rightarrow$(4).

\emph{A repeated root prevents interiority.}
Suppose two root parameters equal $\beta>0$, and write
$P(z)=\omega(1+\beta z)^2R(z)$, where $R$ is the product of the
remaining factors $1+\beta_\nu z$; for $m=2$, set $R\coloneq 1$.
For $\varepsilon>0$, define
$P_\varepsilon(z)\coloneq\omega\bigl((1+\beta z)^2+\varepsilon z^2\bigr)R(z)$,
and let $T_\varepsilon$ be the Toeplitz matrix of bandwidth
$(p,q)$ with Laurent symbol $z^{-q}P_\varepsilon(z)$.
The perturbed quadratic has discriminant $-4\varepsilon<0$.

Its leading coefficient is $\beta^2+\varepsilon>0$, so it has two
nonreal conjugate roots. Multiplication by the polynomial $R$ cannot
remove either of them: no division or cancellation is involved.
Thus $P_\varepsilon$ has nonreal zeros, and
Theorem~\ref{thm:toeplitz} excludes a PBF for $T_\varepsilon$.
The constant and leading coefficients remain positive.
Since the coefficients of $R$ are nonnegative, expressing
$T_\varepsilon-T$ as a sum of shifts gives
$\|T_\varepsilon-T\|_2\le\omega\varepsilon R(1)\longrightarrow0$ as $\varepsilon\downarrow0$.
Consequently $T$ is not an interior point. This proves
(4)$\Rightarrow$(1) and completes the equivalence.
\end{proof}

The $p-1$ lower gaps, the $q-1$ upper gaps and the central gap
in \eqref{eq:toeplitz-stability-gap-limits} account for all
$p+q-1$ gaps between consecutive root parameters. The indicated
ordering gives an analytic choice of factors with fixed initial entries.
Allowing the initial entries to vary covers every ordered partition
of the roots, including partitions that initially place smaller
parameters in the upper group.

\begin{corollary}[Variable initial entries for every ordered root partition]
\label{cor:toeplitz-free-root-orders}
Under the distinct-root hypothesis of
Theorem~\ref{thm:toeplitz-stable-simple-roots}, partition
$\{\beta_1,\ldots,\beta_m\}$ into any ordered upper list of length
$q$ and any ordered lower list of length $p$.
Fix the PBF obtained by moving each lower factor to the left across
the upper factors. There is $\varepsilon>0$ such that every bounded
real matrix $\widetilde T$ of bandwidth at most $(p,q)$ with
$\|\widetilde T-T\|_2<\varepsilon$ has positive factors whose product
is $\widetilde T$. These factors can be chosen to depend real
analytically on $\widetilde T$ in the supremum norms of their entry
sequences, and they equal the fixed factors at $\widetilde T=T$.
Their initial entries may vary with $\widetilde T$.
This includes the factorization attaining the normalization minima
in Theorem~\ref{thm:toeplitz-normalization-minima}.
\end{corollary}
\begin{proof}
At each exchange, Lemma~\ref{lem:exchange-separated-limits}
places the smaller of the two limits in the lower factor and the
larger in the upper factor. Move the first lower factor successively
across the $q$ upper factors. Its limit after these exchanges is
the minimum of its initial limit and the $q$ upper limits: at each
step it retains the smaller value and passes the larger to the
upper factor. The upper list therefore retains the $q$ largest
of these $q+1$ values.

Repeat this procedure with the next lower factor. If the upper
list contains the $q$ largest values encountered so far, inserting
one further value and removing the minimum again leaves the $q$
largest values encountered. Induction through the $p$ lower factors
shows that the final upper limits are
$\beta_1,\ldots,\beta_q$, in some order, and the final lower
limits are $\beta_{q+1},\ldots,\beta_m$, in some order.
Each parameter has been retained exactly once. In particular,
the limits within either group are pairwise distinct.

Normalize the upper factors by
\eqref{eq:toeplitz-stable-normalization}. The lower limits are
unchanged, and the normalized upper limits are the reciprocals of
the final upper diagonal limits, in reversed order. All sequences
are positive and converge to positive constants, so their entries
have common positive lower and finite upper bounds.
Corollary~\ref{stab:cor:free-limits-periodic} gives the pair
conditions on both sides. Moreover,
$\left(\max_i a_i^\infty\right)\left(\max_j b_j^\infty\right)=\beta_{q+1}/\beta_q<1$.
Proposition~\ref{stab:prop:free-full-transfer} gives analytic
normalized LDU factors and the asserted analytic choice of bidiagonal
factors.
\end{proof}

The cross weight in the criterion with fixed initial entries depends on the chosen
PBF. For the ordering in \eqref{eq:toeplitz-stable-root-order}, its
limit is $\beta_{q+1}/\beta_q$; for the ordering that minimizes the initial-condition entries
it is $\beta_m/\beta_1$.
The normalized triangular factors $L,D,U$ are the same, by uniqueness.

The same exchanges apply to positive upper--lower products whose
entries approach distinct limits. This gives a class of
asymptotically Toeplitz matrices with the same analytic dependence of their factors.

\begin{corollary}[Upper--lower products with distinct limits]
\label{cor:asymptotic-ul-stability}
Let $p,q\ge1$, $\omega>0$, and
$T=\omega V_1\cdots V_q B_1\cdots B_p$.
Each $V_i$ is upper bidiagonal with positive diagonal $u_i(n)$,
$n\ge0$, and superdiagonal one. Each $B_j=\LL(\ell_j)$ is unit
lower bidiagonal with positive subdiagonal $\ell_j(n)$, $n\ge1$.
Suppose $u_i(n)\to\eta_i>0$ and $\ell_j(n)\to\xi_j>0$ as
$n\to\infty$, and all $p+q$ limits are pairwise distinct.
The PBF produced by the positive exchanges has positive entries
bounded above and bounded away from zero uniformly over all factors
and rows. There is $\varepsilon>0$ such that every bounded real
matrix $\widetilde T$ of bandwidth at most $(p,q)$ with
$\|\widetilde T-T\|_2<\varepsilon$ has a PBF whose factors depend
real analytically on $\widetilde T$ in their supremum norms and
equal the exchanged factors at $\widetilde T=T$. Their initial
entries may vary with $\widetilde T$. In particular, $T$ is an
interior point of the PBF class in this space of banded matrices.
If the limits are ordered as
\begin{equation}\label{eq:asymptotic-ul-limit-order}
 0<\xi_1<\cdots<\xi_p<\eta_1<\cdots<\eta_q,
\end{equation}
the analytic choice can keep the initial entries fixed.
The diagonals of $T$ converge to those of the Toeplitz matrix with symbol
\begin{equation}\label{eq:asymptotic-ul-symbol}
 a_\infty(z)\coloneq \omega z^{-q}
       \prod_{i=1}^{q}(1+\eta_i z)
       \prod_{j=1}^{p}(1+\xi_j z).
\end{equation}
\end{corollary}
\begin{proof}
Apply Lemma~\ref{lem:exchange-separated-limits} at each of the
$pq$ exchanges. The argument in
Corollary~\ref{cor:toeplitz-free-root-orders} uses only the limiting
parameters. It therefore places the $q$ largest limits in the upper
factors and the remaining $p$ in the lower factors.
Normalize by \eqref{eq:toeplitz-stable-normalization}, and denote
the final upper diagonal limits by $\delta_1,\dots,\delta_q$.
The unit bidiagonal factor sequences have positive distinct limits
within each side, and the largest lower limit divided by the
smallest upper diagonal limit is strictly less than one. Moreover,
$d_n\longrightarrow\omega\prod_{j=1}^{q}\delta_j>0$ as $n\to\infty$.
Each pivot is positive. Convergence and the finitely many initial
entries therefore give common positive lower and finite upper
bounds for both the pivots and all bidiagonal entries. Thus
Corollary~\ref{stab:cor:free-limits-periodic} and
Proposition~\ref{stab:prop:free-full-transfer} give an analytic choice of positive factors throughout a neighborhood,
which proves the interior-point assertion.

Under \eqref{eq:asymptotic-ul-limit-order}, every exchange preserves
its two limits. In the order \eqref{stab:eq:full-given},
$a_j(n)\to\xi_j$, $b_i(n)\to\eta_{q-i+1}^{-1}$ and
$d_n\to\omega\prod_{h=1}^q\eta_h$ as $n\to\infty$.
The lower adjacent ratios tend to $\xi_j/\xi_{j+1}<1$, the upper
adjacent ratios to $\eta_{q-i}/\eta_{q-i+1}<1$, and the cross weight
to $\xi_p/\eta_1<1$, as $n\to\infty$.
Absorbing the finitely many initial terms into the constants gives
all the geometric product bounds of
Theorem~\ref{stab:thm:general-full-stability}, which permits fixed
initial entries.

Finally, each entry of the original product is a finite sum of
products of factor entries at a bounded distance from its row index.
Passing to the limit along each diagonal gives the product of the
constant upper and lower factors with parameters $\eta_i,\xi_j$.
Equation~\eqref{eq:positive-ul-toeplitz} identifies its Laurent symbol
as \eqref{eq:asymptotic-ul-symbol}.
\end{proof}

\begin{example}[The root criterion and finite sections]
The negative-root condition characterizes simultaneous BTP of all
leading sections of a semi-infinite Toeplitz matrix.
At a fixed size, Theorem~\ref{thm:finite} characterizes BTP by
existence of a PBF. For example,
\[
 \begin{bmatrix}3/2&1\\1&3/2\end{bmatrix}
\]
is totally positive and admits a PBF, while
$1+(3/2)z+z^2$ has nonreal conjugate roots.
All four entries of the displayed matrix are positive and its
only minor of order two is $9/4-1=5/4>0$. The polynomial has
discriminant $9/4-4=-7/4<0$, so it fails the root criterion for
simultaneous positivity of all semi-infinite Toeplitz minors even
though this particular finite section is totally positive.
\end{example}

The constant factors in the root construction also give matrices
with the same Toeplitz tail and different initial entries. The next
proposition locates these differences and supplies positive
factorizations for every ordering. Fix positive parameters
$\omega$, $\eta_i$ for $i\in\{1,\dots,q\}$ and $\xi_j$ for
$j\in\{1,\dots,p\}$. With $\Lambda$ the upper shift, put
\begin{equation}\label{eq:constant-root-factors}
 \mathcal U_i^{(0)}\coloneq\Lambda+\eta_i I,\qquad
 \mathcal L_j^{(0)}\coloneq I+\xi_j\Lambda^{\mathsf T}.
\end{equation}
Let $T_0\coloneq\omega\mathcal U_1^{(0)}\cdots
\mathcal U_q^{(0)}\mathcal L_1^{(0)}\cdots\mathcal L_p^{(0)}$.
By \eqref{eq:positive-ul-toeplitz}, this is the Toeplitz matrix
with symbol $\omega z^{-q}\prod_{i=1}^{q}(1+\eta_i z)
\prod_{j=1}^{p}(1+\xi_j z)$.
For an ordering $w$ of these $p+q$ factors, denote their product,
multiplied by $\omega$, by $T_w$.

\begin{proposition}[Positive products with a common Toeplitz tail]
\label{prop:constant-factor-orders}
Every $T_w$ is BTP of bandwidth $(p,q)$ and admits a PBF with
bounded factors and allowed positive entries bounded away from zero.
Its difference from $T_0$ satisfies
\begin{equation}\label{eq:constant-order-boundary}
 (T_w-T_0)_{a,b}=0\quad\text{if }a\ge p\text{ or }b\ge q,
 \qquad \operatorname{rank}(T_w-T_0)\le\min(p,q).
\end{equation}
Moreover, $T_w\le T_0$ entrywise. For a fixed ordering $w$, the
ordered PBF obtained by adjacent positive upper--lower exchanges is
independent of the sequence of exchanges.
\end{proposition}
\begin{proof}
The shift identities $\Lambda\Lambda^{\mathsf T}=I$ and
$\Lambda^{\mathsf T}\Lambda=I-e_0e_0^{\mathsf T}$ give
\begin{equation}\label{eq:constant-factor-commutator}
 \mathcal U_i^{(0)}\mathcal L_j^{(0)}
 -\mathcal L_j^{(0)}\mathcal U_i^{(0)}
 =\xi_j e_0e_0^{\mathsf T},
\end{equation}
where $e_0$ is the first coordinate column vector.
Move each constant upper factor to the left across the lower
factors preceding it. One such change replaces a product
$\omega A\mathcal L_j^{(0)}\mathcal U_i^{(0)}B$ by
$\omega A\mathcal U_i^{(0)}\mathcal L_j^{(0)}B$ and adds
$\omega\xi_j(Ae_0)(e_0^{\mathsf T}B)$.
Here $A$ and $B$ are the products before and after the exchanged
pair, respectively; an empty product is $I$.

Starting with $e_0$, multiplication by an upper bidiagonal factor
does not increase the largest row index in its support, and a lower
factor increases it by at most one. Thus $Ae_0$ is supported in
rows $\{0,\dots,p-1\}$: the exchanged pair already contains one
of the $p$ lower factors. Similarly, $e_0^{\mathsf T}B$ is
supported in columns $\{0,\dots,q-1\}$. All entries of both vectors
are nonnegative. Summing these changes proves the support assertion
and $T_w\le T_0$, because the final product is $T_0$; constant
factors of the same triangular type commute. A matrix supported
in these $p$ rows and $q$ columns has rank at most $\min(p,q)$.

To factor the chosen $T_w$, apply Lemma~\ref{lem:exchange} to each
upper--lower pair until all lower factors precede all upper factors.
Each of these refactorizations preserves $T_w$ and keeps every
allowed entry positive. Initially the entries are positive constants.
The estimates in the proof of Theorem~\ref{thm:toeplitz} preserve
upper bounds and positive lower bounds at each exchange and at the
final diagonal normalization. There are at most $pq$ exchanges,
so they give the asserted bounded PBF. Its extreme diagonals are
positive, and Proposition~\ref{prop:positive-product} gives BTP.

To verify independence of the exchange sequence, retain the labels
of the factors as they move. Every exchange decreases by one the
number of pairs in which an upper factor lies to the left of a
lower factor. Each such initially inverted pair crosses once;
factors of the same type
never cross. Two upper--lower pairs available at the same stage
cannot share a factor, so their exchanges act on disjoint pairs
and commute. Induction on the number of remaining inverted pairs
now compares any two sequences: if their first exchanges differ,
perform both to reach the same intermediate list and apply the
induction hypothesis to the shorter sequences. Thus the final
ordered bidiagonal factors, and their prescribed diagonal
normalization, are the same.
\end{proof}

\begin{remark}[Cyclic products and Christoffel transformations]
\label{rem:constant-products-christoffel}
Moving the first factor of a product to the end gives its cyclic
transform. More generally, writing a product as $AB$ and replacing
it by $BA$ gives the intertwining identities
$(AB)A=A(BA)$ and $B(AB)=(BA)B$; this operation is called a Darboux
transformation. Proposition~\ref{prop:constant-factor-orders}
shows that all cyclic products of the constant factors remain
within the same finite modification of a Toeplitz matrix.

The row and column forms of permuting Christoffel transformations
can be written as follows. Let
$\dd\Psi$ be a $q\times p$ matrix of measures and define the
$r\times r$ matrix $\mathfrak C_r(x)$ by
$[f_1,\dots,f_r]\mathfrak C_r(x)=[f_2,\dots,f_r,xf_1]$ for every row
$[f_1,\dots,f_r]$; for $r=1$, $\mathfrak C_1(x)=[x]$.
The permuting Christoffel operations are
\begin{equation}\label{eq:rectangular-christoffel-shifts}
 \dd\Psi\longmapsto\dd\Psi \mathfrak C_p(x),\qquad
 \dd\Psi\longmapsto \mathfrak C_q(x)^{\mathsf T}\dd\Psi.
\end{equation}
They move the first column or row to the last position and multiply
it by $x$. Since $\mathfrak C_r(x)^r=xI_r$, $p$ column steps or $q$ row steps
give $x\,\dd\Psi$. When the finite moment systems defining the
transformed polynomials have nonzero determinants, their polynomial
families can be constructed. The relation with cyclic factor transfers, including
the degree conditions on the initial polynomial data and the
constant changes of normalization, is developed in
\cite[Section~6]{ManasPBF2026}.

\end{remark}

\begin{example}[The uniform Christoffel chain]
\label{ex:uniform-christoffel-chain}
For $p=2$, $q=1$, $\omega=1$ and $\eta_1=\xi_1=\xi_2=h=4/27$, write
$\mathcal U\coloneq \Lambda+hI$ and $\mathcal L\coloneq I+h\Lambda^{\mathsf T}$.
The products $\mathcal L^2\mathcal U$, $\mathcal L\mathcal U\mathcal L$
and $\mathcal U\mathcal L^2=T_0$ are the monic recurrence matrices
of \cite[Sections~4.1--4.3]{Uniform}, respectively. Their differences
from $T_0$ are $-2h e_0e_0^{\mathsf T}-h^2e_1e_0^{\mathsf T}$,
$-h e_0e_0^{\mathsf T}$ and zero.
These differences can also be read directly from the commutator.
With $E\coloneq e_0e_0^{\mathsf T}$, one has
$\mathcal U\mathcal L-\mathcal L\mathcal U=hE$,
$E\mathcal L=E$, and
$\mathcal L E=E+h e_1e_0^{\mathsf T}$. Consequently
\[
 \mathcal U\mathcal L^2-\mathcal L\mathcal U\mathcal L=hE,
 \qquad
 \mathcal U\mathcal L^2-\mathcal L^2\mathcal U
 =hE\mathcal L+h\mathcal L E
 =2hE+h^2e_1e_0^{\mathsf T}.
\]
The column operation in
\eqref{eq:rectangular-christoffel-shifts} becomes
$(W_1,W_2)\mapsto(W_2,xW_1)$; normalizing the masses gives the
semi-stochastic Christoffel chain of \cite[Theorem~5.1(ii)]{Uniform}.
\end{example}
\subsection{The general Toeplitz class: polynomials and matrix measure}
\label{sec:general-polynomials}

The negative-root Toeplitz class admits explicit polynomial and measure
formulas. An ordered partition of the roots selects a PBF
and the polynomial initial conditions. The measure formulas in
Theorems~\ref{thm:general-resolvent}--\ref{thm:general-measure}
apply to any bounded PBF of the same Toeplitz matrix.
Here an ordered partition means a choice of
which $q$ positive root parameters are assigned to the upper block
and which $p$ to the lower block, together with an order within each
block. Root multiplicities are arbitrary.

Repeated values are retained with their
multiplicities in every list and product below.

Put $m\coloneq p+q$ and $\omega\coloneq t_{-q}>0$.
Choose and order the positive root parameters so that
\begin{equation}\label{eq:ordered-root-split}
 P(z)=\omega\prod_{i=1}^{q}(1+\eta_i z)
                  \prod_{j=1}^{p}(1+\xi_j z),
 \qquad \eta_i,\xi_j>0.
\end{equation}
The initial upper--lower factorization is
\begin{equation}\label{eq:ordered-raw-factorization}
 T=\omega\,\mathcal U^{(0)}_1\cdots\mathcal U^{(0)}_q
                 \mathcal L^{(0)}_1\cdots\mathcal L^{(0)}_p,
 \quad
 \mathcal U^{(0)}_i\coloneq\Lambda+\eta_i I,
 \quad \mathcal L^{(0)}_j\coloneq I+\xi_j\Lambda^{\mathsf T}.
\end{equation}
Apply Lemma~\ref{lem:exchange} first to
$\mathcal U_q^{(0)}\mathcal L_1^{(0)}$.
Move the resulting lower factor past any upper factors still to its
left, starting with the nearest one. It is now the leftmost factor.
For each $j\in\{2,\dots,p\}$, move $\mathcal L_j^{(0)}$ past the
$q$ updated upper factors, starting with the rightmost one.
This places it immediately after the first $j-1$ lower factors.
After these $p$ rounds, all lower factors precede all upper factors.
Write each upper factor as a positive diagonal times a unit upper
bidiagonal matrix, and move these diagonals to the left of the upper
block, as in the proof of Theorem~\ref{thm:toeplitz}. Their product,
together with $\omega$, is the diagonal factor $D$.

The initial-condition matrices are specified first, before their closed
formulas. Write the resulting ordered PBF as
$T=L_1\cdots L_pD U_q\cdots U_1$, as in
\eqref{eq:strict-definition}. Thus $U_1$ is the rightmost upper
factor and $L_1$ is the leftmost lower factor.

Let $e_r$ denote the coordinate column vector with a $1$ in position
$r$ and zeros elsewhere, with indices starting at zero. Define the
unnormalized matrices by
\begin{align}
	\widetilde H_{i,j}
	&\coloneq e_0^{\mathsf T}U_i\cdots U_1e_j,
	&&i,j\in\{0,\dots,q-1\},
	\label{eq:raw-H-definition}\\
	\widetilde K_{i,j}
	&\coloneq e_i^{\mathsf T}L_1\cdots L_je_0,
	&&i,j\in\{0,\dots,p-1\}.
	\label{eq:raw-K-definition}
\end{align}
All products are matrix products, and empty products are understood
as the identity matrix. Thus row $i$ of $\widetilde H$ consists of
the first $q$ entries of row zero of $U_i\cdots U_1$, the product
of the $i$ rightmost upper factors in their displayed order.
Similarly, column $j$ of $\widetilde K$ consists of the first $p$
entries of column zero of $L_1\cdots L_j$, the product of the $j$
leftmost lower factors. In particular, row zero of $\widetilde H$
is $(1,0,\ldots,0)$, and column zero of $\widetilde K$ is
$(1,0,\ldots,0)^{\mathsf T}$.
The first matrix is lower triangular,
the second upper triangular, and their diagonal entries are positive products of allowed bidiagonal entries. Set
\begin{equation}\label{eq:normalized-boundary-definition}
 H_{i,j}\coloneq\frac{\widetilde H_{i,j}}{\widetilde H_{i,i}},\qquad
 K_{i,j}\coloneq\frac{\widetilde K_{i,j}}{\widetilde K_{j,j}}.
\end{equation}
Thus $H$ is normalized by rows and $K$ by columns. These definitions
record the order of the partial products, even when full products
commute. Equivalently, the column of $H^{-1}$ indexed by $b-1$ is
obtained from the first $q$ entries of the unique formal sequence $v$
satisfying $v_0=\cdots=v_{b-2}=0$, $v_{b-1}=1$, and
$U_b\cdots U_1v=0$. These equations are interpreted entrywise.
The characterization is proved below.

The partial products defining $H,K$ also occur in the necessity proof
of the global criterion. This gives a direct connection between the
adjacent series of Section~\ref{sec:finite} and the initial conditions
used here. The following bounds apply to every global PBF, with the
same row and column normalization as in
\eqref{eq:normalized-boundary-definition}.

\begin{proposition}[Bounds on the initial-condition matrices]
\label{prop:normalization-series-bounds}
Let $T=LDU$ be semi-infinite BTP of bandwidth $(p,q)$, and let
$T=L_1\cdots L_pDU_q\cdots U_1$ be a global PBF.
Define $H,K$ by
\eqref{eq:raw-H-definition}--\eqref{eq:normalized-boundary-definition}.
Let $\mathfrak S_k^-$ and $\mathfrak S_\ell^+$ be the adjacent
series \eqref{eq:global-adjacent-series} computed from the
Metelmann factors of $L$ and $U^{\mathsf T}$, respectively.
Then
\begin{equation}\label{eq:normalization-series-bounds}
 K_{k-1,k}\ge\mathfrak S_k^-\quad(k\in\{1,\dots,p-1\}),
 \qquad
 H_{\ell,\ell-1}\ge\mathfrak S_\ell^+
 \quad(\ell\in\{1,\dots,q-1\}).
\end{equation}
\end{proposition}
\begin{proof}
Fix $k\in\{1,\dots,p-1\}$ and put
$h_j\coloneq e_j^{\mathsf T}L_1\cdots L_ke_0$ for
$j\in\{k-1,k\}$. Both numbers are positive.
Column normalization in \eqref{eq:normalized-boundary-definition}
gives $K_{k-1,k}=h_{k-1}/h_k$.
The inequality \eqref{eq:global-necessary-prefix-bound}, proved
for this same partial product, gives
$1/\kappa^-_{k,n}\le h_{k-1}/h_k$ for every $n\ge1$.
By Lemma~\ref{lem:global-inverse-minors}, its left-hand side
increases to $\mathfrak S_k^-$. This proves the lower bound.
For the inequality involving $H$, apply the same argument to $T^{\mathsf T}$.
Its first $\ell$ lower factors are
$U_1^{\mathsf T},\dots,U_\ell^{\mathsf T}$, so their partial
product is $(U_\ell\cdots U_1)^{\mathsf T}$. Its column
normalization is the transpose of the row normalization defining
$H$, giving the asserted entry $H_{\ell,\ell-1}$.

In particular, the order of the factors after transposition is
$U_1^{\mathsf T}\cdots U_\ell^{\mathsf T}$; reversing this order
would use a different partial product and would not give the
entry of $H$ in the statement.
\end{proof}

For a finite multiset $X=(x_1,\ldots,x_r)$, define the complete
homogeneous symmetric polynomial of degree $\ell\ge0$ by
\begin{equation}\label{eq:complete-homogeneous-definition}
 \hc_\ell(X)\coloneq
 \sum_{\substack{\nu_1,\ldots,\nu_r\in\NN\\
                  \nu_1+\cdots+\nu_r=\ell}}
 x_1^{\nu_1}\cdots x_r^{\nu_r}.
\end{equation}
In particular, $\hc_0(X)=1$; set $\hc_\ell(X)\coloneq0$ for $\ell<0$.
If $s\in\NN$, then $X^{\langle s\rangle}$ denotes the multiset in
which each member of $X$ is repeated $s$ times.

A multiset retains repetitions as separate entries; for example,
$\hc_1(a,a,b)=2a+b$. The generating identity follows by multiplying
the geometric series for each copy of each variable. For the empty
multiset this convention gives $\hc_0=1$ and $\hc_\ell=0$ for
$\ell>0$. Thus
\begin{equation}\label{eq:h-reciprocal}
 \prod_{x\in X}(1+xz)^{-s}
 =\sum_{\ell\in\NN}(-1)^\ell
\hc_\ell\bigl(X^{\langle s\rangle}\bigr)z^\ell.
\end{equation}

The homogeneous equations for successive partial products determine
the initial data in Lemma~\ref{lem:ordered-boundary-flag}.
These matrices are independent of the spectral parameter $x$.
Theorem~\ref{thm:general-toeplitz-polynomials} extends them to polynomial
families, and Corollary~\ref{cor:change-normalization} relates the
families and measures obtained from different PBFs.

\begin{lemma}[The initial-condition matrices and their inverses]
\label{lem:ordered-boundary-flag}
For the ordered-root PBF constructed from
\eqref{eq:ordered-raw-factorization}, let $H\in\RR^{q\times q}$ be
the unit lower triangular matrix and $K\in\RR^{p\times p}$ the
unit upper triangular matrix defined by
\eqref{eq:raw-H-definition}--\eqref{eq:normalized-boundary-definition}.
For
$n\in\{0,\dots,q-1\}$ and $b\in\{1,\dots,q\}$, put $r\coloneq n-b+1$. Then
\begin{equation}\label{eq:H-general-inverse}
 (H^{-1})_{n,b-1}=
 \begin{cases}
 0,&n<b-1,\\
 (-1)^r\hc_r(\eta_{q-b+1},\ldots,\eta_q,
                    \xi_1,\ldots,\xi_p),&n\ge b-1.
 \end{cases}
\end{equation}
Define
\begin{equation}\label{eq:reversed-root-data}
 P^\vee(z)\coloneq z^mP(1/z),\qquad
 \omega^\vee\coloneq t_p=\omega\prod_{i=1}^q\eta_i
                              \prod_{j=1}^p\xi_j,
\end{equation}
and $\xi_i^\vee\coloneq \xi_{p+1-i}^{-1}$ for $i\in\{1,\dots,p\}$,
with $\eta_j^\vee\coloneq \eta_{q+1-j}^{-1}$ for $j\in\{1,\dots,q\}$.
Then
\begin{equation}\label{eq:P-reversed-factorization}
 P^\vee(z)=\omega^\vee
 \prod_{i=1}^{p}(1+\xi_i^\vee z)
 \prod_{j=1}^{q}(1+\eta_j^\vee z),
\end{equation}
and, for $a\in\{1,\dots,p\}$, $n\in\{0,\dots,p-1\}$, and $r\coloneq n-a+1$,
\begin{equation}\label{eq:K-general-inverse}
 (K^{-1})_{a-1,n}=
 \begin{cases}
 0,&n<a-1,\\
 (-1)^r\hc_r(\xi_{p-a+1}^\vee,\ldots,\xi_p^\vee,
                    \eta_1^\vee,\ldots,\eta_q^\vee),&n\ge a-1.
 \end{cases}
\end{equation}
In particular, $H^{-1}$ is lower triangular and $K^{-1}$ is upper
triangular, both with diagonal entries equal to one.
\end{lemma}
\begin{proof}
\emph{A sequence satisfying the homogeneous equation for a partial product.}
For $v=\begin{bNiceMatrix}v_0&v_1&\Cdots\end{bNiceMatrix}^{\mathsf T}$,
write $\mathcal Gv(z)\coloneq\sum_{n\in\NN}v_nz^n$.

At this stage the series is formal: no convergence or summability
of $v$ is needed. Multiplication by a banded matrix uses finitely
many entries in each row. The upper-shift formula subtracts $v_0$
because $\sum_{n\ge0}v_{n+1}z^n=(\mathcal Gv-v_0)/z$, while the
lower shift inserts an initial zero and multiplies the series by $z$.
The elementary actions are
\begin{align}
 \mathcal G((I+\xi\Lambda^{\mathsf T})v)(z)&=(1+\xi z)\mathcal Gv(z),
 \label{eq:GF-lower-action}\\
 \mathcal G((\Lambda+\eta I)v)(z)
 &=\frac{(1+\eta z)\mathcal Gv(z)-v_0}{z}.
 \label{eq:GF-upper-action}
\end{align}
For $b\in\{1,\dots,q\}$, let $v^{(b)}$ have generating function
\begin{equation}\label{eq:zero-flag-vector}
 \mathcal Gv^{(b)}(z)\coloneq 
 \frac{z^{b-1}}
 {\prod_{i=q-b+1}^{q}(1+\eta_i z)
  \prod_{j=1}^{p}(1+\xi_jz)}.
\end{equation}
Its first nonzero entry has index $b-1$ and equals one. Applying
the $p$ initial lower factors $\mathcal L_j^{(0)}$ from
\eqref{eq:ordered-raw-factorization} cancels the second product in the denominator,
leaving $z^{b-1}/\prod_{i=q-b+1}^{q}(1+\eta_i z)$.
After the rightmost $t$ upper factors have acted, for
$t\in\{0,\dots,b-1\}$, the generating function is
\[
 \frac{z^{b-1-t}}{\prod_{i=q-b+1}^{q-t}(1+\eta_i z)}.
\]
For $t<b-1$ its constant coefficient is zero, so
\eqref{eq:GF-upper-action} cancels one denominator factor and lowers
the power of $z$ by one. At $t=b-1$ the function is
$(1+\eta_{q-b+1}z)^{-1}$, whose constant coefficient is one; the
last upper factor therefore sends it to zero. Consequently
$\mathcal U^{(0)}_{q-b+1}\cdots\mathcal U^{(0)}_q\mathcal L^{(0)}_1\cdots\mathcal L^{(0)}_p v^{(b)}=0$.

\emph{Preservation of this equation under the exchanges.}
Fix $b$, and let $V_j^{(b)}$ be the product of the rightmost $b$
upper factors after the first $j$ lower factors have been moved
entirely to the left. Induction gives
$V_j^{(b)}\mathcal L_{j+1}^{(0)}\cdots\mathcal L_p^{(0)}v^{(b)}=0$ for $j\in\{0,\dots,p\}$.
The case $j=0$ is the homogeneous equation just established.
At step $j$, the factor $\mathcal L_j^{(0)}$ first crosses the
rightmost $b$ upper factors. The $b$ exchange identities combine to
give
$V_{j-1}^{(b)}\mathcal L_j^{(0)} =\widetilde{\mathcal L}_j V_j^{(b)}$.
The lower bidiagonal factor $\widetilde{\mathcal L}_j$ is formally
invertible: its diagonal entries are nonzero, and each row of its
inverse has only finitely many entries. It can therefore be
cancelled from that identity, proving the induction step.

In detail, if $y$ is the sequence on its right and
$\widetilde{\mathcal L}_j y=0$, row zero gives $y_0=0$.
Row $n$ then expresses $y_n$ as a multiple of $y_{n-1}$, so
induction gives $y_n=0$ for every $n$. This cancellation is valid
for formal sequences even if the triangular inverse is not a
bounded operator on $\ell^2$.
The remaining exchanges move it past the other upper factors and
preserve the rightmost $b$ upper factors. At $j=p$ this yields
$V_p^{(b)}v^{(b)}=0$.

Diagonal normalization can be performed from right to left. For
each $b$ it writes $V_p^{(b)}=E_bU_b\cdots U_1$, with $E_b$
an invertible positive diagonal matrix. Indeed, after normalizing
the rightmost $b-1$ factors, one absorbs their extracted diagonal
into the next upper factor and normalizes that product. Hence
$U_b\cdots U_1v^{(b)}=0$ for the final unit upper factors.

\emph{Identification with the inverse initial-condition matrix.}
An upper product of $b$ bidiagonal factors has positive $b$th
superdiagonal. In its equation at row $n$, the coefficient of
$v_{n+b}$ is therefore nonzero. Prescribing
$v_0=\cdots=v_{b-2}=0$ and $v_{b-1}=1$ determines all remaining
entries recursively.

The row-$n$ equation has only the entries $v_n,\ldots,v_{n+b}$.
Once the first $b$ entries have been fixed, it solves for $v_b$
at row zero, then for $v_{b+1}$ at row one, and so on. Thus the
uniqueness concerns all formal sequences, not an additional decay
condition at infinity. The sequence $v^{(b)}$ consequently satisfies
the required homogeneous equation and initial entries uniquely.

Let $w^{(b)}\coloneq [v_0^{(b)},\ldots,v_{q-1}^{(b)}]^{\mathsf T}$.
By \eqref{eq:raw-H-definition}, row $i$ of $\widetilde H$ is
the zeroth row of $U_i\cdots U_1$, which uses only vector entries
through index $i$. For $i<b-1$ its
product with $w^{(b)}$ vanishes. At $i=b-1$ the product is
$\widetilde H_{b-1,b-1}$, since the only nonzero contributing initial
entry is $v_{b-1}^{(b)}=1$. For $i\ge b$, the product vanishes
because $U_i\cdots U_1$ has $U_b\cdots U_1$ as its rightmost
subproduct. Row normalization by
\eqref{eq:normalized-boundary-definition} gives $Hw^{(b)}=e_{b-1}$. This
proves that $w^{(b)}$ is the column of $H^{-1}$ indexed by $b-1$.
Expanding the denominator in \eqref{eq:zero-flag-vector} by
\eqref{eq:h-reciprocal}, and accounting for the factor $z^{b-1}$,
gives precisely \eqref{eq:H-general-inverse}.

\emph{The left initial conditions directly from row sequences.}
For a row sequence $w=[w_0,w_1,\ldots]$, use the same generating
function convention $\mathcal Gw(z)\coloneq \sum_{n\in\NN}w_nz^n$.
Right multiplication gives
\[
 \begin{aligned}
 \mathcal G\bigl(w(\Lambda+\eta I)\bigr)(z)
   &=(\eta+z)\mathcal Gw(z),&
 \mathcal G\bigl(w(I+\xi\Lambda^{\mathsf T})\bigr)(z)
   &=\frac{(z+\xi)\mathcal Gw(z)-\xi w_0}{z}.
 \end{aligned}
\]
Fix $a\in\{1,\dots,p\}$ and define $w^{(a)}$ by
\begin{equation}\label{eq:left-zero-flag-vector}
 \mathcal Gw^{(a)}(z)\coloneq 
 \frac{z^{a-1}}
 {\prod_{i=1}^{q}(1+\eta_i^{-1}z)
  \prod_{j=1}^{a}(1+\xi_j^{-1}z)}.
\end{equation}
Its first nonzero entry is $w^{(a)}_{a-1}=1$.
Multiplication by the initial upper factors $\mathcal U_i^{(0)}$
in \eqref{eq:ordered-raw-factorization} cancels the first product
in the denominator, with scalar multiplier
$\prod_{i=1}^q\eta_i$.

This multiplier comes from
$\eta_i+z=\eta_i(1+\eta_i^{-1}z)$ for each upper factor.
For a lower factor the corresponding identity is
$z+\xi_j=\xi_j(1+\xi_j^{-1}z)$, which gives the successive
multipliers $\xi_j$ in the next display. After the first $t$ lower factors, for
$t\in\{0,\dots,a-1\}$, the generating function is
\[
 \left(\prod_{i=1}^q\eta_i\right)
 \left(\prod_{j=1}^t\xi_j\right)
 \frac{z^{a-1-t}}{\prod_{j=t+1}^{a}(1+\xi_j^{-1}z)}.
\]
At the last step the constant-term subtraction sends this sequence
to zero. Therefore
\[
 w^{(a)}\mathcal U^{(0)}_1\cdots\mathcal U^{(0)}_q
             \mathcal L^{(0)}_1\cdots\mathcal L^{(0)}_a=0.
\]
After the first $a$ rounds of the prescribed exchanges, their product
identity is
\[
 \mathcal U^{(0)}_1\cdots\mathcal U^{(0)}_q
 \mathcal L^{(0)}_1\cdots\mathcal L^{(0)}_a
 =L_1\cdots L_a\,\mathcal V^{(a)},
\]
where $\mathcal V^{(a)}$ is the current upper product, with positive
diagonal. It is injective for right multiplication on formal row
sequences: in $h\mathcal V^{(a)}=0$, column zero gives $h_0=0$,
and successive columns then give $h_n=0$ for every $n\ge1$.
It follows that $w^{(a)}L_1\cdots L_a=0$.
Later exchanges leave these first $a$ lower factors unchanged.
Diagonal normalization of the final upper block extracts its diagonal
only as far as the position between the lower and upper blocks, and
therefore also leaves them unchanged. This proves the homogeneous
equation for the actual lower factors used to define $K$.

Let $\widehat w^{(a)}\coloneq [w^{(a)}_0,\ldots,w^{(a)}_{p-1}]$.
By \eqref{eq:raw-K-definition}, column $j$ of $\widetilde K$
is column zero of $L_1\cdots L_j$, with support in rows indexed
by $\{0,\dots,j\}$. Its pairing with
$\widehat w^{(a)}$ is zero for $j<a-1$, equals
$\widetilde K_{a-1,a-1}$ for $j=a-1$, and is zero for $j\ge a$
because $w^{(a)}L_1\cdots L_a=0$. Column normalization by
\eqref{eq:normalized-boundary-definition} consequently gives
$\widehat w^{(a)}K=e_{a-1}^{\mathsf T}$.
Thus $\widehat w^{(a)}$ is the row of $K^{-1}$ indexed by $a-1$.
By the reciprocal-parameter definitions following
\eqref{eq:reversed-root-data}, the multisets
$(\xi_1^{-1},\ldots,\xi_a^{-1})$ and
$(\xi_{p-a+1}^\vee,\ldots,\xi_p^\vee)$ agree, as do the reciprocal
upper parameters and $(\eta_1^\vee,\ldots,\eta_q^\vee)$.
Expanding \eqref{eq:left-zero-flag-vector} by
\eqref{eq:h-reciprocal} proves \eqref{eq:K-general-inverse} directly.
The factorization of $P^\vee$ in
\eqref{eq:P-reversed-factorization} follows by reversing the factors
of $P$ and extracting their positive leading coefficients.
\end{proof}

For Toeplitz matrices, the series in
Proposition~\ref{prop:normalization-series-bounds} can be evaluated
from the roots. Their values will also select initial conditions
that attain all the bounds simultaneously. Let $s_\lambda(x)$ denote
the Schur polynomial associated with a partition $\lambda$.
Concretely, draw left-aligned rows of lengths
$\lambda_1\ge\lambda_2\ge\cdots$ and fill their boxes with
indices in $\{1,\dots,m\}$, weakly increasing along rows and
strictly increasing down columns. Then $s_\lambda(x)$ is the sum,
over these fillings, of the products of the variables indexed by
their boxes. The notation $(n^r)$ means $r$ rows of length $n$;
$(n^{r-1},n-1)$ has its last row shortened by one box.

\begin{theorem}[Exact Toeplitz limits and minimal initial conditions]
\label{thm:toeplitz-normalization-minima}
Let $T$ be the BTP Toeplitz matrix in
Theorem~\ref{thm:toeplitz}, with root parameters ordered as
$\beta_1\ge\cdots\ge\beta_m>0$, where $m\coloneq p+q$.
Use the lower and upper ratios of
Corollary~\ref{cor:global-two-sided-criterion} and the adjacent
series of Proposition~\ref{prop:normalization-series-bounds}.
\begin{enumerate}
\item For $k\in\{1,\dots,p-1\}$, $n\ge1$, and $r\coloneq q+k$,
\begin{equation}\label{eq:toeplitz-schur-ratio}
 \frac1{\kappa^-_{k,n}}
 =\frac{s_{(n^{r-1},n-1)}(\beta_1,\dots,\beta_m)}
        {s_{(n^r)}(\beta_1,\dots,\beta_m)}.
\end{equation}
For $k\in\{1,\dots,p-1\}$ and $\ell\in\{1,\dots,q-1\}$,
 the exact values of the adjacent series are
\begin{equation}\label{eq:toeplitz-series-root-values}
 \mathfrak S_k^-=\sum_{\nu=1}^{q+k}\beta_\nu^{-1},
 \qquad
 \mathfrak S_\ell^+=\sum_{\nu=q-\ell+1}^{m}\beta_\nu.
\end{equation}
\item The construction following \eqref{eq:ordered-raw-factorization},
with $\eta_i\coloneq \beta_i$ for $i\in\{1,\dots,q\}$ and
$\xi_j\coloneq \beta_{q+j}$ for $j\in\{1,\dots,p\}$,
attains every bound in \eqref{eq:normalization-series-bounds}.
Thus these entries of $H,K$ are simultaneously minimal over all
global PBFs of $T$. Within the PBFs obtained by that ordered-root
construction, the factorization attaining all these minima is unique.
\end{enumerate}
\end{theorem}
\begin{proof}
\emph{The finite identity.}
Let $\mathsf e_j(\beta)$ be the elementary symmetric polynomial of degree
$j$, with $\mathsf e_0=1$ and $\mathsf e_j\coloneq 0$ outside $\{0,\dots,m\}$.
The symbol \eqref{eq:root-parameters} gives
$T_{i,j}=\omega \mathsf e_{q+i-j}(\beta)$, where $\omega=t_{-q}$.
Fix $k,n$ as in (1) and work in the leading section of order
$n+k$ of the normalized lower factor $L$.
The complementary-minor identity for
$\mathcal Q=\Sigma L^{-1}\Sigma$ gives
\begin{align*}
 \Delta_{k,n}
 &=\det L[\{k,\dots,n+k-1\},\{0,\dots,n-1\}],\\
 \det\mathcal Q[\mathcal I_{k,n},\mathcal C_k^{(1)}]
 &=\det L[\{k-1,k+1,\dots,n+k-1\},\{0,\dots,n-1\}].
\end{align*}
Here the row and column sets on the left are those of
\eqref{eq:global-minor-ratios}; when $n=1$, the second row set
on the right consists only of $k-1$.
The signs of complementary minors are canceled by $\Sigma$, and
$\det L=1$ on this section.

Each minor on the right uses the first $n$ columns. In $T=LDU$,
upper triangularity makes Cauchy--Binet reduce to these same
intermediate columns. Its corresponding minor of $T$ is therefore
the minor of $L$ multiplied by $d_0\cdots d_{n-1}$.
These factors cancel in the quotient defining $1/\kappa^-_{k,n}$,
as do the powers $\omega^n$.

The reduction to one term can be checked directly. A possibly
nonzero minor of the upper triangular matrix $DU$ with column list
$0,\ldots,n-1$ must have intermediate row indices $h_j\le j$.
An increasing list of $n$ nonnegative integers also has $h_j\ge j$,
so the only choice is $h_j=j$ for every $j$. Its determinant is
$d_0\cdots d_{n-1}$.
The denominator becomes $\det[\mathsf e_{r+i-j}]_{i,j\in\{0,\dots,n-1\}}$;
the numerator has its first row replaced by
$[\mathsf e_{r-1-j}]_{j\in\{0,\dots,n-1\}}$.
Reverse both row and column order. The dual Jacobi--Trudi formula
identifies these determinants with the denominator and numerator
in \eqref{eq:toeplitz-schur-ratio}, respectively; see
\cite[Section~3]{MaximenkoMoctezuma2017} for this use of Schur
polynomials for Toeplitz minors.

In the dual Jacobi--Trudi formula the entries are
$\mathsf e_{\lambda'_i-i+j}$, where $\lambda'$ is the partition
of column lengths of the diagram of $\lambda$. For $(n^r)$ these
lengths are $r$ repeated $n$ times. For $(n^{r-1},n-1)$ they are
$r$ repeated $n-1$ times followed by $r-1$. Reversing both orders
places the altered row in exactly that last position. The two
reversals have the same sign, so their signs multiply to one.

\emph{The limit, including repeated roots.}
Fix $1\le r<m$ and $x_1\ge\cdots\ge x_m>0$, and write
$\mathscr D_n(x)\coloneq s_{(n^r)}(x)$ and
$\mathscr N_n(x)\coloneq s_{(n^{r-1},n-1)}(x)$.
The diagonal derivative identity for Schur polynomials
\cite[Remark~2.4]{GrinbergEtAl2024} gives
\begin{equation}\label{eq:rectangular-schur-derivative}
 \sum_{i=1}^{m}\partial_{x_i}\mathscr D_n(x)
 =(n+m-r)\mathscr N_n(x).
\end{equation}
For completeness, this also follows from the determinants just
used, since $\sum_{h=1}^{m}\partial_{x_h}\mathsf e_{r+i-j}
=(m-r+1-i+j)\mathsf e_{r+i-j-1}$.
Split the coefficient as $(m-r+1-i)+j$ and use multilinearity
first in rows and then in columns. In the row contribution,
every replaced row except the first repeats its preceding row.
In the column contribution, every replaced column except the last
repeats its following column. The two surviving determinants are
$\mathscr N_n$, with coefficients $m-r+1$ and $n-1$.
For the last-column determinant, transposition and reversal of
both orders give the first-row determinant. Their sum proves
\eqref{eq:rectangular-schur-derivative}.

Every column in a filling of the rectangle contains $r$ distinct
indices. Its monomial is therefore at most $(x_1\cdots x_r)^n$.
The filling with every entry in row $i$ equal to $i$ attains this
value. A weakly increasing row is determined by the counts of
its $m$ possible entries, each between $0$ and $n$. There are at
most $(n+1)^{rm}$ fillings, and hence
\begin{equation}\label{eq:rectangular-schur-growth}
 (x_1\cdots x_r)^n\le\mathscr D_n(x)
 \le(n+1)^{rm}(x_1\cdots x_r)^n.
\end{equation}
Set $x_i(t)\coloneq x_i\mathrm e^{t/x_i}$ and
$f_n(t)\coloneq n^{-1}\log\mathscr D_n(x_1(t),\dots,x_m(t))$.
For all $t$ in a sufficiently small interval about zero, the
order of the variables is preserved; equal variables remain equal.
Applying \eqref{eq:rectangular-schur-growth} shows that, uniformly
on that interval as $n\to\infty$,
$f_n(t)\to\sum_{i=1}^{r}\log x_i+t\sum_{i=1}^{r}x_i^{-1}$.
Each $f_n$ is convex, since it is the logarithm of a positive sum
of exponentials of linear functions, divided by $n$.

This argument does not require distinct $x_i$. If $x_i=x_j$, their
perturbations $x_i\mathrm e^{t/x_i}$ and $x_j\mathrm e^{t/x_j}$
coincide for every $t$. If two values differ, continuity preserves
their order on a sufficiently small interval; the intersection of
these finitely many intervals works for all pairs. The error in
$n^{-1}\log\mathscr D_n$ from the displayed growth estimate is
at most $rm\log(n+1)/n$, independently of $t$ there.
For any fixed sufficiently small $h>0$, convexity gives
$(f_n(0)-f_n(-h))/h\le f_n'(0)\le(f_n(h)-f_n(0))/h$.
Both bounds tend to $\sum_{i=1}^{r}x_i^{-1}$.
The chain rule and \eqref{eq:rectangular-schur-derivative} now give
\[
 \frac{\mathscr N_n(x)}{\mathscr D_n(x)}
 =\frac{n}{n+m-r}f_n'(0)
 \longrightarrow\sum_{i=1}^{r}x_i^{-1},\qquad n\to\infty.
\]
With $x_i=\beta_i$ and $r=q+k$, this proves the lower formula in
\eqref{eq:toeplitz-series-root-values}, by
Lemma~\ref{lem:global-inverse-minors}.
Apply it to $T^{\mathsf T}$, whose bandwidth is $(q,p)$ and whose
root parameters in decreasing order are
$\beta_m^{-1},\dots,\beta_1^{-1}$. It gives
$\mathfrak S_\ell^+=\sum_{\nu=q-\ell+1}^{m}\beta_\nu$.

\emph{Attainment and the ordered-root construction.}
The first off-diagonal of a unit triangular inverse is the
negative of the first off-diagonal of the original matrix.
Taking degree one in \eqref{eq:H-general-inverse} and
\eqref{eq:K-general-inverse} therefore gives
\begin{equation}\label{eq:ordered-first-normalization-entries}
 K_{k-1,k}=\sum_{i=1}^{q}\eta_i^{-1}
                +\sum_{j=1}^{k}\xi_j^{-1},
 \qquad
 H_{\ell,\ell-1}=\sum_{j=1}^{p}\xi_j
                +\sum_{i=q-\ell+1}^{q}\eta_i.
\end{equation}
The ordered choice in (2) gives precisely
\eqref{eq:toeplitz-series-root-values}. Proposition~\ref{prop:normalization-series-bounds}
proves its minimality over all global PBFs.

To identify the equality cases among ordered-root constructions,
the lower sum in \eqref{eq:ordered-first-normalization-entries}
must use the $q+k$ largest parameters, counted with multiplicity.
The upper sum must use the $p+\ell$ smallest.

For the first assertion, replacing a selected parameter by a larger
unselected one decreases its reciprocal strictly unless the two
parameters are equal. Thus a minimizing sum of $q+k$ reciprocals
uses precisely the largest values, allowing equal values at the
cut. The second assertion follows from the same replacement
argument for the parameters themselves, selecting the smallest
$p+\ell$ values. Equality at a tie changes labels only, not the
selected values or their multiplicities.
Differences between successive selected multisets, together with
the full multiset at $k=p$ and $\ell=q$, give
$\xi_j=\beta_{q+j}$ for $j\in\{2,\dots,p\}$ and
$\eta_i=\beta_i$ for $i\in\{1,\dots,q-1\}$.
The remaining pair satisfies
$\{\eta_q,\xi_1\}=\{\beta_q,\beta_{q+1}\}$.
Both orders of this pair give the same factorization: since
$\Lambda\Lambda^{\mathsf T}=I$,
\[
 (\Lambda+aI)(I+b\Lambda^{\mathsf T})
 =\Lambda+(a+b)I+ab\Lambda^{\mathsf T}
 =(\Lambda+bI)(I+a\Lambda^{\mathsf T}).
\]
This is the first pair exchanged in the prescribed procedure.
Lemma~\ref{lem:exchange} gives a unique lower--upper decomposition
of its product, so the factors after that step coincide. Every
subsequent exchange and diagonal normalization is then identical.
Repeated parameters allow exchanges of equal values, which also
leave the factors unchanged. This proves (2).
\end{proof}

The minimum just obtained has a stability consequence. For this
ordered-root PBF, $K_{0,1}=1/a_1(1)=\mathfrak S^-_1$ when $p\ge2$,
and $H_{1,0}=1/b_1(1)=\mathfrak S^+_1$ when $q\ge2$.
Proposition~\ref{stab:prop:endpoint-unstable} therefore excludes a choice of PBFs throughout a neighborhood
with those initial entries fixed whenever one side
has width at least two. When all root parameters are distinct,
Corollary~\ref{cor:toeplitz-free-root-orders} gives an analytic choice of factors through this minimizing PBF by
allowing its initial entries to vary. The ordering in
Theorem~\ref{thm:toeplitz-stable-simple-roots} instead permits those
entries to remain fixed.

The parameter count also distinguishes these two
orders. The minimizing PBF is locally unique when all factor
entries are allowed to vary, while the order with fixed initial
entries belongs to a family of maximal local dimension.

\begin{corollary}[Local uniqueness at the normalization minimum]
\label{cor:toeplitz-minimizer-isolated}
Assume that the root parameters in
Theorem~\ref{thm:toeplitz-normalization-minima} are pairwise distinct.
The PBF attaining the normalization minima in that theorem is the
only PBF of $T$ in some supremum-norm neighborhood of its factor
tuple. For every bounded real matrix of bandwidth at most $(p,q)$
sufficiently close to $T$ in operator norm, there is exactly one
PBF sufficiently close to this minimizing tuple in the supremum
norm of its factor entries. Its factors depend real analytically
on the matrix.

The ordering \eqref{eq:toeplitz-stable-root-order}, which permits
fixed initial entries, instead belongs to a local family of PBFs
of the same $T$ with dimension
$p(p-1)/2+q(q-1)/2$.
\end{corollary}
\begin{proof}
For the minimizing order, the upper parameters before the exchanges are
$\eta_j=\beta_j$ and the lower parameters are
$\xi_i=\beta_{q+i}$. Every upper parameter exceeds every lower
one, so the exchanges preserve these limits by
Lemma~\ref{lem:exchange-separated-limits}.
Equation~\eqref{eq:toeplitz-stable-normalization} gives
$a_i(n)\longrightarrow\beta_{q+i}$ and $b_j(n)\longrightarrow\beta_{q-j+1}^{-1}$ as $n\to\infty$.
Both lists of limiting factor entries are strictly decreasing.
Corollary~\ref{stab:cor:local-dimension-limits} therefore gives
$s_-=s_+=0$.
The diagonal entries converge to
$\omega\prod_{\nu=1}^{q}\beta_\nu>0$ as $n\to\infty$.
They are positive in every row, so they have uniform positive lower
and finite upper bounds.

To see the lower bound, choose a row beyond which each entry exceeds
half of its positive limit, then take the minimum of those bounds
and the finitely many earlier entries. The upper bound follows in
the same way from convergence and a finite maximum.
Analytic dependence of normalized LDU was established in
Corollary~\ref{cor:toeplitz-free-root-orders}.
Proposition~\ref{stab:prop:local-full-family} now proves local
uniqueness and analytic dependence near the minimizing tuple.

For the order \eqref{eq:toeplitz-stable-root-order}, the limits
in \eqref{eq:toeplitz-stable-factor-limits} are increasing on each
side. Every pair contributes one parameter, giving the stated
maximal dimension.

Indeed, on the lower side all $p(p-1)/2$ pairs satisfy the increasing
limit inequality in Corollary~\ref{stab:cor:local-dimension-limits};
on the upper side all $q(q-1)/2$ pairs do. In normalized LDU the
diagonal factor is fixed by $T$, and the two triangular products
vary independently as specified in
Proposition~\ref{stab:prop:local-full-family}. Their dimensions add.
\end{proof}

This local uniqueness concerns all PBFs sufficiently close in the
supremum norm on their factor entries. It does not assert global
uniqueness among all PBFs of $T$ or among all factorizations attaining
the normalization minima.

The matrices $H,K$ normalized in
\eqref{eq:normalized-boundary-definition} now provide the initial
conditions for the polynomials. Arrange the left vectors $A_n(x)\in\RR[x]^p$ as
columns and the right vectors $B_n(x)\in\RR[x]^{1\times q}$ as rows:
\begin{equation}\label{eq:polynomial-arrays}
	A(x)\coloneq 
	\begin{bNiceMatrix}[margin=2pt]
		A_0(x) & A_1(x) & \Cdots[shorten-end=5pt]
	\end{bNiceMatrix},
	\qquad
	B(x)\coloneq 
	\begin{bNiceMatrix}
		B_0(x)\\
		B_1(x)\\
		\Vdots
	\end{bNiceMatrix}.
\end{equation}
Thus $A(x)$ has size $p\times\NN$ and $B(x)$ has size
$\NN\times q$. Denote the $a$th component of $A_n$ by
$A_n^{(a)}$ and the $b$th component of $B_n$ by $B_n^{(b)}$,
for $a\in\{1,\dots,p\}$ and $b\in\{1,\dots,q\}$. The equations $A(x)T=xA(x)$ and
$TB(x)=xB(x)$ are supplemented by
\begin{equation}\label{eq:general-initial-polynomials}
	\begin{bNiceMatrix}
		A_0 & \Cdots & A_{p-1}
	\end{bNiceMatrix}
	=K^{-1},
	\qquad
	\begin{bNiceMatrix}
		B_0\\
		\Vdots\\
		B_{q-1}
	\end{bNiceMatrix}
	=H^{-1}.
\end{equation}
For the right recurrence, row $n$ determines $B_{n+q}(x)$ from
earlier rows because its coefficient is $t_{-q}>0$. The left
recurrence similarly determines $A_{n+p}(x)$ because its coefficient
is $t_p>0$. The first $q$ right rows and first $p$ left columns in
\eqref{eq:general-initial-polynomials} therefore specify unique
polynomial families. The following theorem solves both recurrences
by generating functions.

\begin{theorem}[Generating functions for the full Toeplitz class]
\label{thm:general-toeplitz-polynomials}
The solutions of $A(x)T=xA(x)$ and $TB(x)=xB(x)$ with
\eqref{eq:general-initial-polynomials} are determined by
\begin{align}
 \sum_{n\in\NN}B_n^{(b)}(x)z^n
 &=\frac{\omega z^{b-1}\prod_{i=1}^{q-b}(1+\eta_i z)}
 {P(z)-xz^q},&&b\in\{1,\dots,q\},
 \label{eq:general-right-GF}\\
 \sum_{n\in\NN}A_n^{(a)}(x)z^n
 &=\frac{\omega^\vee z^{a-1}
                  \prod_{i=1}^{p-a}(1+\xi_i^\vee z)}
 {P^\vee(z)-xz^p},&&a\in\{1,\dots,p\}.
 \label{eq:general-left-GF}
\end{align}
Here $P,\omega,\eta_i$ are as in \eqref{eq:ordered-root-split};
$P^\vee,\omega^\vee$ are defined in
\eqref{eq:reversed-root-data}, and $\xi_i^\vee=\xi_{p+1-i}^{-1}$
as in Lemma~\ref{lem:ordered-boundary-flag}.
\end{theorem}
\begin{proof}
Denote the first quotient by $F_b(z;x)$. Its denominator has
nonzero constant term $\omega$, so it defines a formal power series
in $z$ with polynomial coefficients in $x$. Indeed,
\[
 \frac1{P(z)-xz^q}=\sum_{k\ge0}\frac{x^kz^{qk}}{P(z)^{k+1}},
\]
and only finitely many indices $k$ can contribute to a fixed
coefficient of $z$. This is an identity of formal series and does
not require a bound on $x$. At $x=0$ the quotient $F_b$ is

\[
 \frac{z^{b-1}}
 {\prod_{i=q-b+1}^{q}(1+\eta_i z)
  \prod_{j=1}^{p}(1+\xi_jz)},
\]
whose first $q$ coefficients form the column of $H^{-1}$ indexed by
$b-1$, by Lemma~\ref{lem:ordered-boundary-flag}. The dependence on $x$
first enters at a power of $z$ at least $q$, so these initial
coefficients are unchanged for arbitrary $x$.

The identity
$(P(z)-xz^q)F_b(z;x) =\omega z^{b-1}\prod_{i=1}^{q-b}(1+\eta_i z)$
has a right-hand side of degree $q-1$. For $n\in\NN$, its coefficient
of $z^{n+q}$ is therefore zero. Since
$P(z)=\sum_{k=-q}^{p}t_kz^{k+q}$, the corresponding coefficient on
the left gives
$\sum_{k=-q}^{p}t_kB_{n-k}^{(b)}(x)=xB_n^{(b)}(x)$,
where $B_j^{(b)}\coloneq 0$ for $j<0$.
This is the component form of $TB=xB$. The term $t_{-q}B_{n+q}^{(b)}$
has nonzero coefficient, so the recurrence and the first $q$ values
determine every later coefficient uniquely.

For the second quotient, the constant term of $P^\vee$ is
$\omega^\vee=t_p$. Its value at $x=0$ is
\[
 \frac{z^{a-1}}
 {\prod_{i=p-a+1}^{p}(1+\xi_i^\vee z)
  \prod_{j=1}^{q}(1+\eta_j^\vee z)}.
\]
The first $p$ coefficients are the row of $K^{-1}$ indexed by $a-1$.
They are independent of $x$. Now
$P^\vee(z)=\sum_{k=-q}^{p}t_kz^{p-k}$, and the numerator of the
second quotient has degree $p-1$. Comparison of the coefficient of
$z^{n+p}$ consequently gives
$\sum_{k=-q}^{p}t_kA_{n+k}^{(a)}(x)=xA_n^{(a)}(x)$,
where $A_j^{(a)}\coloneq 0$ for $j<0$.
This is $AT=xA$, with $t_pA_{n+p}^{(a)}$ determining the next value.
Thus both displayed generating functions have the required initial
conditions and solve their respective recurrences.
\end{proof}

The generating functions also give finite formulas for each polynomial
coefficient, using the following lists of variables.
For $b\in\{1,\dots,q\}$ and $k\in\NN$, let $\mathcal X_{b,k}$ be the list
containing $k$ copies of each of $\eta_1,\ldots,\eta_{q-b}$ and
$k+1$ copies of each of $\eta_{q-b+1},\ldots,\eta_q,\xi_1,\ldots,\xi_p$.
Explicitly,
\begin{equation}\label{eq:Xbk}
\begin{multlined}[t][.9\textwidth]
 \mathcal X_{b,k}\coloneq\bigl(
 \underbrace{\eta_1,\ldots,\eta_1}_{k},\ldots,
 \underbrace{\eta_{q-b},\ldots,\eta_{q-b}}_{k},
 \underbrace{\eta_{q-b+1},\ldots,\eta_{q-b+1}}_{k+1},\ldots,
 \underbrace{\eta_q,\ldots,\eta_q}_{k+1},\\
 \underbrace{\xi_1,\ldots,\xi_1}_{k+1},\ldots,
 \underbrace{\xi_p,\ldots,\xi_p}_{k+1}\bigr).
\end{multlined}
\end{equation}
Blocks of multiplicity zero and empty index ranges are omitted.
For $a\in\{1,\dots,p\}$, define $\mathcal Z_{a,k}$ analogously:
\begin{equation}\label{eq:Zak}
\begin{multlined}[t][.9\textwidth]
 \mathcal Z_{a,k}\coloneq\bigl(
 \underbrace{\xi_1^\vee,\ldots,\xi_1^\vee}_{k},\ldots,
 \underbrace{\xi_{p-a}^\vee,\ldots,\xi_{p-a}^\vee}_{k},
 \underbrace{\xi_{p-a+1}^\vee,\ldots,\xi_{p-a+1}^\vee}_{k+1},\ldots,
 \underbrace{\xi_p^\vee,\ldots,\xi_p^\vee}_{k+1},\\
 \underbrace{\eta_1^\vee,\ldots,\eta_1^\vee}_{k+1},\ldots,
 \underbrace{\eta_q^\vee,\ldots,\eta_q^\vee}_{k+1}\bigr).
\end{multlined}
\end{equation}
Thus the first $p-a$ variables have multiplicity $k$, and all the
remaining variables have multiplicity $k+1$.
The parameters $\eta_i,\xi_j$ are those of
\eqref{eq:ordered-root-split}; the reciprocal parameters
$\xi_i^\vee,\eta_j^\vee$ are defined in
Lemma~\ref{lem:ordered-boundary-flag} after \eqref{eq:reversed-root-data}.
The coefficient formulas use the complete homogeneous polynomials
$\hc_\ell$ of \eqref{eq:complete-homogeneous-definition}.

\begin{corollary}[Finite coefficient formulas]
\label{cor:general-toeplitz-coefficients}
For the lists $\mathcal X_{b,k}$ and $\mathcal Z_{a,k}$ in
\eqref{eq:Xbk} and \eqref{eq:Zak}, the polynomial families of
Theorem~\ref{thm:general-toeplitz-polynomials} have the following
coefficients. For $b\in\{1,\dots,q\}$, $B_n^{(b)}=0$ for $n<b-1$, while otherwise
\begin{equation}\label{eq:general-B-coefficients}
 B_n^{(b)}(x)=
 \sum_{k=0}^{\lfloor(n-b+1)/q\rfloor}
 (-1)^{N_{n,b,k}}\omega^{-k}
 \hc_{N_{n,b,k}}(\mathcal X_{b,k})x^k,
 \quad N_{n,b,k}\coloneq n-b+1-qk.
\end{equation}
For $a\in\{1,\dots,p\}$, $A_n^{(a)}=0$ for $n<a-1$, while otherwise
\begin{equation}\label{eq:general-A-coefficients}
 A_n^{(a)}(x)=
 \sum_{k=0}^{\lfloor(n-a+1)/p\rfloor}
 (-1)^{M_{n,a,k}}(\omega^\vee)^{-k}
 \hc_{M_{n,a,k}}(\mathcal Z_{a,k})x^k,
 \quad M_{n,a,k}\coloneq n-a+1-pk.
\end{equation}
Every nonzero component has degree given by the corresponding upper
summation limit.
\end{corollary}
\begin{proof}
Expand the first generating function as a geometric series in
$xz^q/P(z)$. For a fixed $k\in\NN$, its contribution is
\[
 \frac{\omega^{-k}x^kz^{qk+b-1}}
 {\prod_{i=1}^{q-b}(1+\eta_i z)^k
  \prod_{i=q-b+1}^{q}(1+\eta_i z)^{k+1}
  \prod_{j=1}^{p}(1+\xi_j z)^{k+1}}.
\]
The three products explain the multiplicities in $\mathcal X_{b,k}$.
By \eqref{eq:h-reciprocal}, the coefficient still needed to obtain
$z^n$ is
$(-1)^{N_{n,b,k}}\hc_{N_{n,b,k}}(\mathcal X_{b,k})$, where
$N_{n,b,k}=n-b+1-qk$. It vanishes when $N_{n,b,k}<0$, which gives
exactly the upper limit in \eqref{eq:general-B-coefficients}.

For the left generating function the contribution of $x^k$ is
\[
 \frac{(\omega^\vee)^{-k}x^kz^{pk+a-1}}
 {\prod_{i=1}^{p-a}(1+\xi_i^\vee z)^k
  \prod_{i=p-a+1}^{p}(1+\xi_i^\vee z)^{k+1}
  \prod_{j=1}^{q}(1+\eta_j^\vee z)^{k+1}}.
\]
Its remaining coefficient has index $M_{n,a,k}=n-a+1-pk$ and equals
$(-1)^{M_{n,a,k}}\hc_{M_{n,a,k}}(\mathcal Z_{a,k})$.
This proves \eqref{eq:general-A-coefficients}, including its range.
At the largest allowed $k$, the complete homogeneous polynomial is
strictly positive: all its variables are positive, and its value
at index zero is one.

The variable list is nonempty: on the right it contains the $p$
lower parameters with multiplicity $k+1$, and on the left it
contains the $q$ reciprocal upper parameters with multiplicity
$k+1$. For positive degree, taking the whole exponent on any one
of these variables gives a positive summand. No cancellation
between distinct $k$ can remove this coefficient, since distinct
$k$ correspond to distinct powers of $x$. The highest displayed power of $x$ therefore
has a nonzero coefficient in each nonzero polynomial.
\end{proof}

The same generating functions express the polynomials in terms of
the roots of their denominators. This is the partial-fraction method
for rational generating functions; see
\cite[Chapter~IV]{FlajoletSedgewick2009}. Here the numerators record
the initial conditions of each component, so they must be retained
when passing from the recurrence to a root formula.
Put $m\coloneq p+q$ and define
\begin{equation}\label{eq:polynomial-root-data}
 \begin{aligned}
 Q(z;x)&\coloneq P(z)-xz^q,
 &Q^\vee(z;x)&\coloneq P^\vee(z)-xz^p=z^mQ(1/z;x),\\
 \mathsf N_b^B(z)&\coloneq
 \omega z^{b-1}\prod_{i=1}^{q-b}(1+\eta_i z),
 &\mathsf N_a^A(z)&\coloneq
 \omega^\vee z^{a-1}\prod_{i=1}^{p-a}(1+\xi_i^\vee z).
 \end{aligned}
\end{equation}
The indices range over $b\in\{1,\dots,q\}$ and
$a\in\{1,\dots,p\}$, with the ordered parameters and reciprocals
of \eqref{eq:ordered-root-split} and \eqref{eq:reversed-root-data}.
For each fixed $x\in\CC$, both denominators have degree $m$ and
nonzero constant term. Their roots are therefore nonzero.
These roots depend on $x$; at $x=0$ they are the roots of $P$ and
$P^\vee$, respectively.

\begin{proposition}[Root formulas for the polynomial components]
\label{prop:polynomial-root-formulas}
Use the data in \eqref{eq:polynomial-root-data}. If the roots
$z_1,\dots,z_m$ of $Q(\,\cdot\,;x)$ are simple, then, for every
$n\in\NN$,
\begin{align}
 B_n^{(b)}(x)&=-\sum_{j=1}^m
 \frac{\mathsf N_b^B(z_j)}{Q_z(z_j;x)z_j^{n+1}},
 \qquad b\in\{1,\dots,q\},
 \label{eq:B-root-formula}\\
 A_n^{(a)}(x)&=\sum_{j=1}^m
 \frac{\mathsf N_a^A(z_j^{-1})z_j^{n+m-1}}{Q_z(z_j;x)},
 \qquad a\in\{1,\dots,p\}.
 \label{eq:A-root-formula}
\end{align}
Here $Q_z$ denotes differentiation with respect to $z$.

For multiple roots, either component is obtained as follows.
Take $(\mathsf N,\mathsf D)\coloneq (\mathsf N_b^B,Q)$ for a right
component, or $(\mathsf N_a^A,Q^\vee)$ for a left component,
with $x$ fixed. Let $\rho_1,\dots,\rho_J$ be the distinct roots
of $\mathsf D$, of multiplicities $\mu_1,\dots,\mu_J$.
Define $\mathsf D_j(z)\coloneq
\mathsf D(z)/(z-\rho_j)^{\mu_j}$ and, for
$r\in\{1,\dots,\mu_j\}$,
\begin{equation}\label{eq:root-principal-coefficients}
 \gamma_{j,r}\coloneq\frac{1}{(\mu_j-r)!}
 \left.
 \frac{\mathrm d^{\mu_j-r}}{\mathrm dz^{\mu_j-r}}
 \frac{\mathsf N(z)}{\mathsf D_j(z)}\right|_{z=\rho_j}.
\end{equation}
Then the corresponding polynomial is
\begin{equation}\label{eq:confluent-polynomial-root-formula}
 [z^n]\frac{\mathsf N(z)}{\mathsf D(z)}
 =\sum_{j=1}^J\sum_{r=1}^{\mu_j}
 (-1)^r\gamma_{j,r}
 \binom{n+r-1}{r-1}\rho_j^{-n-r}.
\end{equation}
\end{proposition}
\begin{proof}
The numerator degrees are $q-1$ and $p-1$, both less than $m$.
For simple roots, the partial-fraction decomposition of the right
generating function is therefore
\[
 \frac{\mathsf N_b^B(z)}{Q(z;x)}
 =\sum_{j=1}^m\frac{\mathsf N_b^B(z_j)}{Q_z(z_j;x)}
                  \frac{1}{z-z_j}.
\]
At $z=0$, the coefficient of $z^n$ in $(z-z_j)^{-1}$ is
$-z_j^{-n-1}$. Comparing with \eqref{eq:general-right-GF}
gives \eqref{eq:B-root-formula}, including the initial indices
where some components vanish.

The roots of $Q^\vee$ are $z_j^{-1}$. Differentiating
$Q^\vee(w;x)=w^mQ(1/w;x)$ at $w=z_j^{-1}$ gives
$(Q^\vee)_z(z_j^{-1};x)=-z_j^{2-m}Q_z(z_j;x)$;
the term containing $Q(z_j;x)$ vanishes. Applying the same
coefficient calculation to \eqref{eq:general-left-GF} thus gives
\[
 -\frac{\mathsf N_a^A(z_j^{-1})}
 {(Q^\vee)_z(z_j^{-1};x)(z_j^{-1})^{n+1}}
 =\frac{\mathsf N_a^A(z_j^{-1})z_j^{n+m-1}}{Q_z(z_j;x)}
\]
for each root. Summing proves \eqref{eq:A-root-formula}.

At a root $\rho_j$ of multiplicity $\mu_j$, the quotient
$\mathsf N/\mathsf D_j$ is analytic because
$\mathsf D_j(\rho_j)\ne0$. Its Taylor coefficients of orders
$0,\dots,\mu_j-1$, divided by $(z-\rho_j)^{\mu_j}$,
give the principal part
$\sum_{r=1}^{\mu_j}\gamma_{j,r}(z-\rho_j)^{-r}$.
Subtract all these principal parts from $\mathsf N/\mathsf D$.
The difference is a rational function without finite poles and
tends to zero at infinity, so it is zero. Finally,
\[
 [z^n](z-\rho_j)^{-r}
 =(-1)^r\rho_j^{-r}[z^n](1-z/\rho_j)^{-r}
 =(-1)^r\binom{n+r-1}{r-1}\rho_j^{-n-r},
\]
which proves \eqref{eq:confluent-polynomial-root-formula}.
\end{proof}

A root of multiplicity $\mu$ thus contributes a polynomial in $n$
of degree at most $\mu-1$, multiplied by its reciprocal to the
power $n$. Formula \eqref{eq:confluent-polynomial-root-formula}
also covers cancellation between numerator and denominator:
the corresponding principal coefficients simply vanish.
When simple roots merge as $x$ varies, the sum of their contributions
has this finite value, since its generating coefficient is the
polynomial already determined by
Theorem~\ref{thm:general-toeplitz-polynomials}.

The individual simple-root fractions need not have finite limits:
their denominators $Q_z(z_j;x)$ may tend to zero. The cancellation
is justified by summing before taking the limit. Equivalently, a
small contour surrounding the colliding roots expresses their
combined contribution as a residue sum; the integrand remains
analytic on that contour as the roots merge, and its integral has
the limit given by the principal part in
\eqref{eq:root-principal-coefficients}.

The root formulas also give differential equations when the root
parameters of the symbol are arbitrary. Successive differentiation
and reduction modulo a polynomial equation give the classical
differential-resolvent construction
\cite[Section~2.1]{BostanChyzakLecerfSalvySchost2007}. Keeping an exponent
as a parameter is treated in \cite[Sections~3--4]{Nahay2004Powersum}.
Applied to the numerators in \eqref{eq:polynomial-root-data}, this
construction retains each component's normalization.

Consider a quotient $N(z)/Q(z;x)$, where
$Q(z;x)=P(z)-xz^q$, $P$ has degree $m=p+q$, $P(0)\ne0$, and
$\deg N<m$. The right components use
$N=\mathsf N_b^B$ in \eqref{eq:polynomial-root-data}. The left components
use $P^\vee$, $p$ in place of $q$, and $N=\mathsf N_a^A$.
Thus the same construction applies to both arrays.

Let $\nu$ be an auxiliary parameter. All polynomial remainders in the
following construction are taken modulo $Q$ and have degree less than
$m$ in $z$.

Here reduction modulo $Q$ means replacing a polynomial by the
remainder after division by $Q$. Two polynomials represent the
same element if their difference is a multiple of $Q$. Thus their
values agree at every root of $Q$. Arithmetic on these remainders
is performed by ordinary addition and multiplication followed by
this reduction. Thus each remainder is represented by its $m$ coefficients
in the basis $1,z,\dots,z^{m-1}$. Work over the field $\mathbb R(\nu,x)$ and let
$\mathsf C_Q(z;x)$ be the inverse of $Q_z$ modulo $Q$:
\begin{equation}\label{eq:general-ode-modular-inverse}
 \operatorname{rem}(\mathsf C_Q Q_z,Q)=1,\qquad \deg_z\mathsf C_Q<m.
\end{equation}
The coefficients of $\mathsf C_Q$ are rational functions of $x$; they are obtained
by the Euclidean algorithm, or by solving a linear system of order $m$.

To be explicit, the Euclidean algorithm finds polynomials $A,B$
with coefficients in $\mathbb R(x)$ such that $AQ+BQ_z=1$;
then $\mathsf C_Q$ is the remainder of $B$. The proof below shows
that $Q$ and $Q_z$ have no common factor over this field. The
formula is therefore well defined even though particular values of
$x$ can make these two polynomials share a root.
Define an operation on polynomial remainders by
\begin{equation}\label{eq:general-ode-connection}
 \mathscr D_\nu f\coloneq
 \operatorname{rem}\!\left(
 \partial_x f+z^q\mathsf C_Q\,\partial_z f+
 \bigl((q-\nu-1)z^{q-1}\mathsf C_Q-z^q Q_{zz}\mathsf C_Q^2\bigr)f,Q\right).
\end{equation}
Starting with $f_0\coloneq N$, calculate $f_{j+1}\coloneq \mathscr D_\nu f_j$.
Write $\mathbf c_j$ for the column of the $m$ coefficients of $f_j$
in the ordered basis $1,z,\ldots,z^{m-1}$.

\begin{proposition}[Differential equations for an arbitrary polynomial symbol]
\label{prop:general-toeplitz-differential}
Suppose that $p,q\ge1$, $P$ has degree $m=p+q$ and nonzero constant
term, and $N\ne0$ has degree less than $m$. The construction
\eqref{eq:general-ode-modular-inverse}--\eqref{eq:general-ode-connection}
gives polynomials $\mathsf d_0(\nu,x),\ldots,\mathsf d_r(\nu,x)$, with $r\le m$,
such that, for every $n\ge0$,
\begin{equation}\label{eq:general-toeplitz-differential}
 \sum_{j=0}^{r}\mathsf d_j(n,x)\frac{\mathrm d^j}{\mathrm dx^j}
 \left([z^n]\frac{N(z)}{P(z)-xz^q}\right)=0.
\end{equation}
The coefficients can be chosen so that, for each fixed $n$, they
are not all identically zero as polynomials in $x$. The order may
decrease after $\nu=n$ is substituted. In particular, every right
and left component in Theorem~\ref{thm:general-toeplitz-polynomials}
has an equation of order at most $p+q$, with coefficients calculated
from its generating function.

Here is an explicit elimination rule. Let $r$ be the first index for
which $\mathbf c_0,\ldots,\mathbf c_r$ are linearly dependent over
$\mathbb R(\nu,x)$. Choose $r$ rows for which the restriction $R$ of
$[\mathbf c_0\ \cdots\ \mathbf c_{r-1}]$ is nonsingular, and let
$R_j$ be obtained by replacing its column indexed by $j$ with the same
rows of $\mathbf c_r$, for $j\in\{0,\ldots,r-1\}$. A coefficient
list before polynomial normalization is
\begin{equation}\label{eq:general-ode-determinants}
 (\mathsf d_0,\ldots,\mathsf d_r)\coloneq (-\det R_0,\ldots,-\det R_{r-1},\det R).
\end{equation}
Clear the denominators in this list and divide out the common
polynomial factor of its entries in $\mathbb R[\nu,x]$.
\end{proposition}

\begin{proof}
First, $Q_z$ is invertible modulo $Q$. Indeed, a common zero of
$Q$ and $Q_z$ must be nonzero and must satisfy
$zP'(z)-qP(z)=0$. This is a nonzero polynomial independent of $x$:
its constant coefficient is $-qP(0)\ne0$ and its leading coefficient
is $p$ times the leading coefficient of $P$.
Its zeros give only finitely many possible values of
$x=P(z)/z^q$.

The discriminant of a polynomial is a polynomial expression in its
coefficients that vanishes exactly when two roots coincide, provided
the leading coefficient is nonzero. Here the leading coefficient
of $Q$ is the fixed nonzero leading coefficient of $P$. The finite
exceptional set just described is therefore exactly the possible
zero set of this discriminant. Thus the discriminant of $Q$, as a polynomial in $x$,
is not identically zero, and \eqref{eq:general-ode-modular-inverse}
is well defined over $\mathbb R(x)$. This argument allows repeated
roots of $P$.

Choose a complex value of $x$ away from that finite discriminant set
and a sufficiently small simply connected neighborhood. The roots
$z_1(x),\ldots,z_m(x)$ of $Q$ are distinct, analytic and nonzero
there. Choose an analytic logarithm of each root on this neighborhood.
Implicit differentiation of $Q(z_i(x);x)=0$ gives
$z_i'=z_i^q/Q_z(z_i;x)$. Put
$w_i(\nu,x)\coloneq z_i(x)^{-\nu-1}/Q_z(z_i(x);x)$.
The product and chain rules give
\begin{equation}\label{eq:general-ode-root-weight}
 \frac{\partial_x w_i}{w_i}
 =\frac{(q-\nu-1)z_i^{q-1}}{Q_z(z_i;x)}
  -\frac{z_i^q Q_{zz}(z_i;x)}{Q_z(z_i;x)^2}.
\end{equation}
For clarity, the second denominator in the weight depends on $x$
both through the root and explicitly through $Q$:
$\partial_x[Q_z(z_i;x)]=Q_{zz}(z_i;x)z_i'-qz_i^{q-1}$.
Substitution of this expression gives exactly
\eqref{eq:general-ode-root-weight}.

For any polynomial remainder $f$ with coefficients depending on $x$,
the same calculation now gives
\begin{equation}\label{eq:general-ode-differentiate-weight}
 \partial_x\bigl(w_i f(z_i;x)\bigr)
 =w_i(\mathscr D_\nu f)(z_i;x).
\end{equation}
The modular inverse $\mathsf C_Q$ has value $1/Q_z(z_i;x)$ at each root,
and reduction modulo $Q$ preserves that value. Consequently the
recursion defining $f_j$ computes the successive derivatives of
$w_i N(z_i)$, for every root separately.

For the derivative of $f(z_i(x);x)$, the two contributions are
$\partial_xf(z_i;x)$ and $z_i'\partial_zf(z_i;x)$. Adding the
logarithmic derivative of $w_i$ produces the three terms in
\eqref{eq:general-ode-connection}. Evaluation at a root is also
independent of the representative: a multiple of $Q$ vanishes
identically along $z_i(x)$, so its total derivative along that
root vanishes as well.

There are $m$ coordinates in the quotient by $Q$. Hence a first
linear dependence occurs at some $r\le m$. The preceding $r$
columns are independent, so the matrix $R$ in the statement exists.
Since this is the first dependence, $\mathbf c_r$ is a linear
combination of $\mathbf c_0,\ldots,\mathbf c_{r-1}$ on all $m$
rows. Its coefficients are uniquely determined by the chosen $r$
rows because $R$ is nonsingular. Cramer's rule computes those
coefficients and therefore gives the identity of full remainders
$(\det R)f_r-\sum_{j=0}^{r-1}(\det R_j)f_j=0$.
Multiply its value at each root by $w_i$ and use
\eqref{eq:general-ode-differentiate-weight}. The differential operator
with coefficients \eqref{eq:general-ode-determinants} annihilates
each function $w_iN(z_i)$, and therefore their negative sum.

Clearing denominators gives polynomial coefficients in $\nu$ and
$x$. Dividing their common factor leaves the same identity wherever
that factor is nonzero. Both sides then extend analytically, so the
identity also holds at its zeros. This includes values of $\nu$
excluded by intermediate divisions: with the chosen logarithms,
$z_i^{-\nu-1}$ is an entire function of $\nu$.
After this normalization, the coefficients cannot all vanish
identically in $x$ at a fixed value $\nu=n$. If they did, the factor
$\nu-n$ would divide every coefficient in $\mathbb R[\nu,x]$,
contrary to the removal of the common factor.

More explicitly, regard each coefficient as a polynomial in $x$
whose coefficients are polynomials in $\nu$. If specialization at
$\nu=n$ gives the zero polynomial in $x$, each of these coefficient
polynomials vanishes at $n$ and is divisible by $\nu-n$. This is
the common factor that normalization has excluded. Some leading
coefficients may vanish at $n$, which is why only a decrease in
order, rather than a fixed order for every $n$, is asserted.

Finally, when $\nu=n\ge0$, the partial-fraction formula gives
\[
 -\sum_{i=1}^{m}w_i(n,x)N(z_i(x))
 =[z^n]\frac{N(z)}{P(z)-xz^q}.
\]
This coefficient is a polynomial in $x$, since $P(0)\ne0$.
Equation~\eqref{eq:general-toeplitz-differential} holds on the open
set used above. After polynomial normalization, its left-hand side
is a polynomial in $x$, so it vanishes identically. In particular
the equation remains valid at values of $x$ where roots of $Q$
coalesce, including $x=0$ when $P$ has repeated roots.
\end{proof}

The polynomial index enters through $q-\nu-1$ in
\eqref{eq:general-ode-connection}; the desired polynomials
are not needed as input. The initial numerator enters through $f_0=N$, so different components
may give different scalar equations.
The matrix form of \eqref{eq:general-ode-connection} is common to
all right components, with a corresponding matrix form for the
left components. The scalar elimination may introduce additional zeros of the
coefficient of the highest derivative. Away from the discriminant,
the roots, the weights $w_i$ and their products with $N(z_i)$ are
analytic on a neighborhood after branches have been chosen.
Their scalar equation may nevertheless have a vanishing highest
coefficient because of the determinants used to eliminate the
other coordinates. Such a singularity is called apparent: it comes
from that scalar representation, while the local root expressions
remain analytic. This explanation uses the root expressions
before their sum is identified with a polynomial.

There is also one explicit first-order system for all numerators
associated with the same denominator. Put
$g_h(x)\coloneq [z^h]Q(z;x)^{-1}$ for $h\ge0$, and $g_h\coloneq 0$ for $h<0$,
and define
\begin{equation}\label{eq:general-ode-common-system}
 \begin{aligned}
 \mathbf v_n&\coloneq [g_n\ g_{n-1}\ \cdots\ g_{n-m+1}]^{\mathsf T},
 &\mathbf v_n'&=M_n(x)\mathbf v_n,&
 [M_\nu(x)]_{r,s}&\coloneq [z^s]\mathscr D_\nu z^r,
 &&0\le r,s\le m-1.
 \end{aligned}
\end{equation}
Indeed, the root formula applied to $z^r/Q$ identifies the $r$th
entry of $\mathbf v_n$ with $-\sum_i w_i(n,x)z_i^r$;
\eqref{eq:general-ode-differentiate-weight} gives the indicated row
of the system. If $N(z)=\sum_{r=0}^{m-1}N_rz^r$, its generating
coefficient is $\sum_rN_rg_{n-r}$. In particular, all the right
components are constant linear combinations of the same vector
$\mathbf v_n$. The entries of $M_\nu$ are rational in $x$ and affine
in $\nu$, and their poles are among the discriminant values of $Q$.

The affine dependence follows because the only occurrence of
$\nu$ in \eqref{eq:general-ode-connection} is the multiplier
$-\nu z^{q-1}\mathsf C_Q f$. Polynomial reduction introduces no
new dependence on $\nu$. Poles in $x$ come from the modular
inverse; the linear system defining it is nonsingular wherever
$Q$ has distinct roots.
The coefficient recursion used in the scalar elimination is then
$\mathbf c_{j+1}=\partial_x\mathbf c_j+M_\nu^{\mathsf T}\mathbf c_j$.

\begin{example}[A quadratic symbol]
Let $p=q=1$, $P(z)=c_0+c_1z+c_2z^2$, with $c_0c_2\ne0$,
and $N=c_0$. Write $t\coloneq x-c_1$ and
$D\coloneq t^2-4c_0c_2$. In the basis $1,z$, direct reduction modulo
$c_2z^2-tz+c_0$ gives
\[
 \mathscr D_\nu 1
 =\frac{\nu t-2(\nu+1)c_2z}{D},\qquad
 \mathscr D_\nu z
 =\frac{2\nu c_0-(\nu+1)tz}{D}.
\]
For a scalar rational function $h(x)$ and a remainder $f$,
$\mathscr D_\nu(hf)=h'f+h\mathscr D_\nu f$.
Applying this rule to the displayed expression for
$\mathscr D_\nu1$, using $t'=1$ and $D'=2t$, gives
$\mathscr D_\nu^2 1=\nu(\nu+2)/D-(3t/D)\mathscr D_\nu1$.
Thus $D\mathscr D_\nu^2 1+3t\mathscr D_\nu1-
\nu(\nu+2)=0$. Since the numerator $N=c_0$ is constant in $x$,
multiplication by $c_0$ leaves this equation unchanged. The
resulting equation is
\begin{equation}\label{eq:general-ode-quadratic-example}
 \bigl((x-c_1)^2-4c_0c_2\bigr)(B_n^{(1)})''
 +3(x-c_1)(B_n^{(1)})'-n(n+2)B_n^{(1)}=0.
\end{equation}
This is the shifted and scaled Chebyshev equation of the second
kind.
\end{example}

Different ordered partitions of the roots generally give different
PBFs, initial-condition matrices, and polynomial bases; the bases are
related by constant triangular matrices with diagonal entries equal
to one. When all root
parameters coincide, changing their order or partition leaves these
matrices unchanged.

The next step identifies the measure for these polynomial families. A $q\times p$ matrix of positive measures means an
array $[\Psi_{b,a}]_{b\in\{1,\dots,q\},\,a\in\{1,\dots,p\}}$ in which each entry is a finite nonnegative
Borel measure. Positivity here is entrywise positivity, not
positive semidefiniteness of a square matrix; the array need not
be square. Integrals of the matrix measure mean the integrals of
its entries. Its support is the union of the supports of these
finitely many scalar measures. The spectral Favard construction of
\cite[Sections~7--9]{BFM2023} produces such measures from positive
bidiagonal factors; the unbounded theory is developed in
\cite{BFMUnbounded}. The following corollary proves the bounded
construction directly for the ordered boundary products $H,K$ in
\eqref{eq:normalized-boundary-definition}. This fixes both the
factor order and the moments used in the explicit root calculations
that follow. Compact support gives uniqueness from the moments.

For such a measure, the Stieltjes transform is
$W(\lambda)\coloneq \int(\lambda-x)^{-1}\,\dd\Psi(x)$, entry by entry.
Where a continuous density exists, it is recovered as
$-\pi^{-1}\im W(x+\mathrm{i}0)$.

The notation $x+\mathrm i0$ denotes the limit from
$x+\mathrm i\epsilon$ as $\epsilon\downarrow0$. The sign follows
from
$\operatorname{Im}(x+\mathrm i\epsilon-y)^{-1}
=-\epsilon/((x-y)^2+\epsilon^2)$; after multiplication by
$-1/\pi$ this is the positive Poisson kernel, whose integral tends
to the density at a point of continuity. Thus computing the boundary values
of this transform will convert the polynomial construction into
an explicit measure.

\begin{corollary}[Positive measures for the prescribed normalization]
\label{cor:general-favard}
Let $T$ be a bounded semi-infinite banded matrix with a PBF
$T=L_1\cdots L_pDU_q\cdots U_1$ whose factors are bounded.
Define $H,K$ from its partial products by
\eqref{eq:raw-H-definition}--\eqref{eq:normalized-boundary-definition}.
Let $B_n$ be the row polynomial solutions of $TB=xB$, with the
matrix of first $q$ rows equal to $H^{-1}$, and let $A_n$ be the
column polynomial solutions of $AT=xA$, with the matrix of first
$p$ columns equal to $K^{-1}$.
There is a unique $q\times p$ matrix $\Psi$ of entrywise positive,
compactly supported finite measures such that
\begin{equation}\label{eq:general-moments}
 \int_{\RR} x^\ell\,\dd\Psi(x)=HE_qT^\ell E_p^{\mathsf T}K,
 \qquad \ell\in\NN,
\end{equation}
where $E_r\coloneq [I_r\ 0]$ selects the first $r$ coordinates. It satisfies
\begin{equation}\label{eq:general-biorthogonality}
 \int_{\RR} B_n(x)\,\dd\Psi(x)\,A_l(x)=\delta_{n,l}.
\end{equation}
\end{corollary}
\begin{proof}
Positive measures for the ordered boundary products defining $H,K$
can be constructed directly. The finite-dimensional ingredients are the
classical oscillation, interlacing, and eigenvector sign theorems;
see \cite[Theorems~1.4--1.6]{BFM2023}.

Fix $N\ge\max\{1,p-1,q-1\}$ and denote leading sections by a
superscript $[N]$. Triangularity gives
$T^{[N]}=L_1^{[N]}\cdots L_p^{[N]}D^{[N]}
U_q^{[N]}\cdots U_1^{[N]}$: in a selected entry of the lower--diagonal--upper
product, every intermediate index is at most $N$. Put
\[
 F_{j,N}\coloneq L_1^{[N]}\cdots L_j^{[N]},\qquad
 G_{i,N}\coloneq U_i^{[N]}\cdots U_1^{[N]},
 \qquad F_{0,N}=G_{0,N}=I_{N+1}.
\]
The matrices $F_{j,N}^{-1}T^{[N]}F_{j,N}$ and
$G_{i,N}T^{[N]}G_{i,N}^{-1}$ are cyclic rearrangements of these
finite factors. Each is totally nonnegative and invertible, and
has positive first subdiagonal and first superdiagonal entries.
Indeed, each such entry includes a positive contribution from an
appropriate bidiagonal factor, with diagonal entries selected from
all the other factors. The Gantmacher--Krein criterion makes these
matrices, and $T^{[N]}$, oscillatory. Their eigenvalues agree by
similarity.

Write these eigenvalues as
$\lambda_1^{[N]}<\cdots<\lambda_{N+1}^{[N]}$, and choose real
right and left eigenvectors $u_k,w_k$ of $T^{[N]}$, with $w_k$ a
row, such that $(u_k)_N=(w_k)_N=1$. The eigenvector sign theorem
gives nonzero endpoint entries and exactly $N+1-k$ sign changes.
It follows that their first entries have sign
$\sigma_{k,N}\coloneq(-1)^{N+1-k}$.
The vectors $G_{i,N}u_k$ and $w_kF_{j,N}$ are eigenvectors of the
corresponding cyclic matrices. Their last entries remain one,
because $G_{i,N}$ is unit upper triangular and $F_{j,N}$ is unit
lower triangular. The same eigenvector sign theorem therefore
gives sign $\sigma_{k,N}$ for their first entries as well.

Let $E_{r,N}\coloneq[I_r\ 0]\in\RR^{r\times(N+1)}$ and define
\[
 \alpha_k\coloneq HE_{q,N}u_k,\qquad
 \beta_k^{\mathsf T}\coloneq w_kE_{p,N}^{\mathsf T}K.
\]
The definitions of the boundary products give, for
$0\le i<q$ and $0\le j<p$,
\[
 (\alpha_k)_i=
 \frac{e_0^{\mathsf T}G_{i,N}u_k}{\widetilde H_{i,i}},\qquad
 (\beta_k)_j=
 \frac{w_kF_{j,N}e_0}{\widetilde K_{j,j}}.
\]
There is no boundary correction in these identities: the first
row of $G_{i,N}$ is supported in $\{0,\ldots,i\}$ and the first
column of $F_{j,N}$ in $\{0,\ldots,j\}$. Their positive
denominators and the preceding sign rule show that every component
of $\alpha_k$ and $\beta_k$ has sign $\sigma_{k,N}$.

Put $\chi_{N+1}(\lambda)\coloneq
\det(\lambda I_{N+1}-T^{[N]})$ and
$\chi_N(\lambda)\coloneq\det(\lambda I_N-T^{[N-1]})$.
The spectral projector at $\lambda_k^{[N]}$ is
$\Pi_k=c_k u_kw_k$, where, by its last diagonal entry,
\[
 c_k=(\Pi_k)_{N,N}
 =\frac{\chi_N(\lambda_k^{[N]})}
        {\chi_{N+1}'(\lambda_k^{[N]})}>0.
\]
The strict inequality follows from interlacing of the two monic
characteristic polynomials: numerator and denominator have the
same sign. Consequently the finite matrix measure
\[
 \dd\Psi^{[N]}\coloneq
 \sum_{k=1}^{N+1}c_k\alpha_k\beta_k^{\mathsf T}
                   \delta_{\lambda_k^{[N]}}
\]
has strictly positive entries at every atom, and the spectral
decomposition gives
\[
 \int x^\ell\,\dd\Psi^{[N]}(x)
 =HE_{q,N}(T^{[N]})^\ell E_{p,N}^{\mathsf T}K.
\]
All atoms belong to $[0,\|T\|_2]$, since $T^{[N]}$ is a
compression of $T$. The total mass matrix is the fixed matrix
$HE_qE_p^{\mathsf T}K$. Compactness of finite positive measures
on this interval, applied to the finitely many entries with a
common subsequence, yields an entrywise positive limit $\Psi$.
For every fixed $\ell$, bandedness makes the displayed moment
equal to $HE_qT^\ell E_p^{\mathsf T}K$ once $N$ is sufficiently
large: the finitely many contributing paths then remain inside
the section. Weak convergence on the common compact interval
therefore proves \eqref{eq:general-moments}.

These moments also determine the polynomial normalization.
Write the finite coefficient expansions as
$B_n(x)=\sum_s B_{n,s}x^s$ and
$A_l(x)=\sum_t A_{l,t}x^t$,
where $B_{n,s}$ are row vectors and $A_{l,t}$ are column vectors.
All sums below are finite because $A_l$ and $B_n$ are polynomial
vectors. Each power of the banded matrix $T$ has finite bandwidth,
so selecting finitely many of its rows or columns gives vectors
with finite support. Define the row and column vectors
$\boldsymbol r_n\coloneq\sum_s B_{n,s}HE_qT^s$ and
$\boldsymbol c_l\coloneq\sum_t T^tE_p^{\mathsf T}KA_{l,t}$.
The initial polynomial matrices $H^{-1}$ and $K^{-1}$ give
$\boldsymbol r_n=e_n^{\mathsf T}$ for $n\in\{0,\dots,q-1\}$
and $\boldsymbol c_l=e_l$ for $l\in\{0,\dots,p-1\}$.
The equations $TB=xB$ and $AT=xA$ imply
$\boldsymbol r_nT=\sum_j T_{n,j}\boldsymbol r_j$ and
$T\boldsymbol c_l=\sum_j T_{j,l}\boldsymbol c_j$.
In the first equation, the nonzero coefficient $T_{n,n+q}$
determines $\boldsymbol r_{n+q}$ from the preceding rows.
In the second, $T_{l+p,l}\ne0$ determines
$\boldsymbol c_{l+p}$ from the preceding columns.
The coordinate vectors satisfy the same equations and initial data.
Induction therefore gives
$\boldsymbol r_n=e_n^{\mathsf T}$ and $\boldsymbol c_l=e_l$
for all $n,l\in\NN$.
For each pair of coefficient indices $s,t$, the corresponding
integral is
$B_{n,s}HE_qT^{s+t}E_p^{\mathsf T}KA_{l,t}$.
Factoring $T^{s+t}=T^sT^t$ and summing these finitely many terms
gives, by the definitions of $\boldsymbol r_n$ and $\boldsymbol c_l$,

\[
 \int_{\RR} B_n(x)\,\dd\Psi(x)\,A_l(x)
 =\boldsymbol r_n\boldsymbol c_l=e_n^{\mathsf T}e_l=\delta_{n,l}.
\]

Finally, two compactly supported measures with the prescribed moments
have a difference that annihilates every polynomial. Polynomial density
in the continuous functions on a compact interval containing both
supports makes that difference zero, entry by entry. This proves
uniqueness and also applies to finite signed measures.
\end{proof}

The bounds in Proposition~\ref{prop:normalization-series-bounds}
also constrain the total masses of these measures. In
\eqref{eq:general-moments} with $\ell=0$, row zero of $H$ is
$e_0^{\mathsf T}$ and column zero of $K$ is $e_0$.
Thus $\int_{\RR}\dd\Psi_{1,a}=K_{0,a-1}$ and
$\int_{\RR}\dd\Psi_{b,1}=H_{b-1,0}$, with
$a\in\{1,\dots,p\}$ and $b\in\{1,\dots,q\}$.
In particular,
\begin{equation}\label{eq:favard-mass-series-bounds}
 \int_{\RR}\dd\Psi_{1,2}\ge\mathfrak S_1^-
 \quad(p\ge2),\qquad
 \int_{\RR}\dd\Psi_{2,1}\ge\mathfrak S_1^+
 \quad(q\ge2).
\end{equation}
For Toeplitz matrices, the bounded PBF in
Theorem~\ref{thm:toeplitz-normalization-minima} attains these
minimum masses in the unit-diagonal normalization
\eqref{eq:normalized-boundary-definition}.

Corollary~\ref{cor:general-favard} applies to the ordered-root PBF
constructed from \eqref{eq:ordered-raw-factorization} and to the
geometric-root PBF in \eqref{eq:q-PBF}. More generally, fix a
bounded banded matrix $T$ with two bounded PBFs. Their different
partial products may change $H,K$ in
\eqref{eq:raw-H-definition}--\eqref{eq:normalized-boundary-definition}.
The recurrence matrix remains $T$, but the polynomial initial
conditions and the measure normalization may change. The following
corollary gives their exact relation using
\eqref{eq:general-moments} and uniqueness of the compactly supported
measures.

\begin{corollary}[Change of factorization and normalization]
\label{cor:change-normalization}
Suppose two bounded PBFs of the same bounded matrix $T$ give initial
matrices $(H,K)$ and $(H',K')$ as in
Corollary~\ref{cor:general-favard}. Denote their measures by
$\Psi,\Psi'$ and their polynomial families by $(B_n,A_n)$ and
$(B'_n,A'_n)$. Then, as matrices of measures and polynomial vectors,
\begin{align}
 \dd\Psi'&=H'H^{-1}\,\dd\Psi\,K^{-1}K',
 \label{eq:change-measure}\\
 B'_n&=B_nHH'^{-1},\qquad
 A'_n=K'^{-1}KA_n,\qquad n\in\NN.
 \label{eq:change-polynomials}
\end{align}
\end{corollary}
\begin{proof}
For every $\ell\in\NN$, equation~\eqref{eq:general-moments} gives
\[
 \int x^\ell H'H^{-1}\,\dd\Psi(x)\,K^{-1}K'
 =H'E_qT^\ell E_p^{\mathsf T}K'.
\]
These are the moments of $\Psi'$. Uniqueness for compactly
supported signed measures therefore gives equality
with $\dd\Psi'$. The first $q$ rows of the proposed right family
form $H^{-1}HH'^{-1}=H'^{-1}$, and the first $p$ columns of the
proposed left family form $K'^{-1}KK^{-1}=K'^{-1}$.
Multiplication by these constant matrices commutes with the
recurrences: for example,
$T(BHH'^{-1})=(TB)HH'^{-1}=xBHH'^{-1}$.
The left identity is checked with the constant matrix on the left.
Thus the transformations in \eqref{eq:change-polynomials} have
the required initial data and recurrences. The nonzero extreme
diagonals determine each subsequent vector uniquely.
\end{proof}

For these changes of basis, positivity of the transformed measure
follows from the second PBF and Corollary~\ref{cor:general-favard}.

The distinction between existence of measures and maximal polynomial
degrees is part of the normality theory in \cite{BFM2024}. For the
present Toeplitz families, the degrees are verified directly by
Corollary~\ref{cor:general-toeplitz-coefficients}.

Here the relevant degree condition means that every component
attains the maximum degree allowed by its position in the step-line
ordering. The degree statement is needed to turn the first $n$
polynomial vectors into a basis of the first $n$ coordinate
monomials, rather than only a spanning list with possible linear
dependences. The mixed relations
then follow by the basis argument underlying
\cite[Corollary~9.2]{BFM2023}.

\begin{corollary}[Orthogonality of mixed type for every ordered partition of roots]
\label{cor:general-step}
For $r\in\mathbb N$ and $n\in\NN$, define $\nu_r(n)\in\NN^r$ by writing $n=rk+s$,
$s\in\{0,\dots,r-1\}$, and setting its first $s$ entries equal to $k+1$ and
the remaining entries equal to $k$. Then
\begin{equation}\label{eq:general-step-degrees}
 \deg B_n^{(b)}=\nu_q(n+1)_b-1,
 \qquad
 \deg A_n^{(a)}=\nu_p(n+1)_a-1,
\end{equation}
where degree $-1$ denotes the zero polynomial. The polynomial families of
Theorem~\ref{thm:general-toeplitz-polynomials} satisfy, for
$a\in\{1,\dots,p\}$ and $b\in\{1,\dots,q\}$,
\begin{align}
 \sum_{b=1}^q\int_{\RR} x^j B_n^{(b)}(x)\,\dd\Psi_{b,a}(x)&=0,
 &&j\in\{0,\dots,\nu_p(n)_a-1\},
 \label{eq:general-mixed-right}\\
 \sum_{a=1}^p\int_{\RR} x^j A_n^{(a)}(x)\,\dd\Psi_{b,a}(x)&=0,
 &&j\in\{0,\dots,\nu_q(n)_b-1\}.
 \label{eq:general-mixed-left}
\end{align}
Their common biorthogonal normalization is
\eqref{eq:general-biorthogonality}.
\end{corollary}
\begin{proof}
Corollary~\ref{cor:general-toeplitz-coefficients} gives the stated
degrees and nonzero highest coefficients. Use $e_0,\ldots,e_{p-1}$ for the coordinate columns in $\RR^p$.
Order the column monomials as
$e_0,\ldots,e_{p-1}, xe_0,\ldots,xe_{p-1}, x^2e_0,\ldots$.
The vector $A_j$ contains only the first $j+1$ monomials in this
list, with a nonzero coefficient of the last one. The coefficient
matrix of $A_0,\ldots,A_{n-1}$ in the first $n$ monomials is
therefore triangular with nonzero diagonal. These vectors form a
basis of
$\operatorname{span}\{x^je_{a-1}:
a\in\{1,\dots,p\},\ j\in\{0,\dots,\nu_p(n)_a-1\}\}$.
Biorthogonality of $B_n$ against $A_0,\ldots,A_{n-1}$ is thus
equivalent to \eqref{eq:general-mixed-right}.

For each allowed $x^je_{a-1}$, expand that coordinate monomial in
this basis. Pairing the expansion with $B_n$ gives a linear
combination of $\delta_{n,l}=0$ for $l<n$, hence zero. Conversely,
if all the coordinate-monomial pairings vanish, they vanish for
all their linear combinations, including those $A_l$. The same
reasoning proves both directions of the asserted equivalence.
For the other relation, use the row coordinate vectors
$e_{b-1}^{\mathsf T}\in\RR^{1\times q}$, ordered by increasing
degree and then by $b$. The same triangular argument shows that
$B_0,\ldots,B_{n-1}$ span the first $n$ row monomials, proving
\eqref{eq:general-mixed-left}. The pairing at equal indices is one
by \eqref{eq:general-biorthogonality}.
\end{proof}

The resolvent at $\lambda$ is the bounded inverse
$(\lambda I-T)^{-1}$, when this inverse exists. For
$|\lambda|>\|T\|_2$ it is given by the norm-convergent Neumann
series $\sum_{\ell\ge0}\lambda^{-\ell-1}T^\ell$.
Its entries are computed by separating the roots of the polynomial
$P(z)-\lambda z^q$, where $P$ is the symbol in
\eqref{eq:general-toeplitz}. For sufficiently large $|\lambda|$,
let $S=S(\lambda)$ be its multiset of $q$ roots tending to zero and
let $\mathcal Z=\mathcal Z(\lambda)$ be its multiset of $p$ roots
tending to infinity. All roots are counted with multiplicity. To justify this separation,
fix any radius $r>0$. For sufficiently large $|\lambda|$,
$|\lambda|r^q>\max_{|z|=r}|P(z)|$. Rouch\'e's theorem then gives
exactly $q$ roots inside that circle, counted with multiplicity.
Applying this to arbitrarily small fixed radii shows that these
$q$ roots tend to zero. Applying it also to arbitrarily large
fixed radii shows that the other $p$ roots tend to infinity.
The unordered groups do not require a globally single-valued labeling
of the individual roots. Thus

\begin{equation}\label{eq:general-kernel}
 Q(z;\lambda)=P(z)-\lambda z^q
 =t_p\prod_{u\in S}(z-u)\prod_{v\in\mathcal Z}(z-v).
\end{equation}
The notation $\mathcal Z^{-1}$ means the multiset of reciprocals
$v^{-1}$ for $v\in\mathcal Z$. The polynomials $\hc_\ell$ from
\eqref{eq:complete-homogeneous-definition} are evaluated on these root
multisets. Put
\begin{equation}\label{eq:general-gamma}
 \mathfrak g(\lambda)\coloneq\frac{(-1)^{q+1}}{\omega}\prod_{u\in S}u.
\end{equation}

\begin{theorem}[Finite root formula for the resolvent]
\label{thm:general-resolvent}
With the root multisets $S,\mathcal Z$ in \eqref{eq:general-kernel}
and $\mathfrak g$ in \eqref{eq:general-gamma}, for $|\lambda|$ sufficiently large,
\begin{equation}\label{eq:general-R}
 [ (\lambda I-T)^{-1}]_{i,j}
 =\mathfrak g(\lambda)\sum_{k=0}^{\min(i,j)}
 \hc_{j-k}(S)\hc_{i-k}(\mathcal Z^{-1}).
\end{equation}
For any bounded PBF of this Toeplitz matrix, let $H,K$ be its
initial-condition matrices \eqref{eq:normalized-boundary-definition}.
The Stieltjes transform of the measure in
Corollary~\ref{cor:general-favard} is
\begin{equation}\label{eq:general-W}
 W(\lambda)\coloneq \int_{\RR}\frac{\dd\Psi(x)}{\lambda-x}
 =H[ (\lambda I-T)^{-1}_{i,j}]_{
       \substack{i\in\{0,\dots,q-1\}\\j\in\{0,\dots,p-1\}}}K.
\end{equation}
\end{theorem}
\begin{proof}
For column $j$, let $F_j(z)\coloneq \sum_{i\in\NN}[(\lambda I-T)^{-1}]_{i,j}z^i$.
For $|\lambda|>\|T\|_2$, each resolvent column belongs to $\ell^2$,
so $F_j$ is analytic for $|z|<1$.

Indeed, for $0\le r<1$, Cauchy--Schwarz bounds its absolute sum by
$\|[(\lambda I-T)^{-1}]_{\cdot,j}\|_2(1-r^2)^{-1/2}$ on
$|z|\le r$. The resulting locally uniform convergence justifies
multiplication by the polynomial $P$ and coefficient comparison. Increasing $|\lambda|$ if necessary
puts all roots in $S$ in this disk and all roots in $\mathcal Z$
outside it. Multiply the column equation
$(\lambda I-T)[(\lambda I-T)^{-1}]_{\cdot,j}=e_j$ by $z^{n+q}$
and sum over $n\in\NN$. Write $r_h\coloneq [(\lambda I-T)^{-1}]_{h,j}$ and set $r_h\coloneq 0$ for
$h<0$. For $n\ge0$, the column equation reads
$\lambda r_n-\sum_{k=-q}^{p}t_kr_{n-k}=\delta_{n,j}$.
For $k\ge0$, summing $r_{n-k}z^{n+q}$ begins at $n=k$ and
produces $z^{k+q}F_j$. For $k<0$, the entries
$r_0,\ldots,r_{-k-1}$ must be removed from that full product.
The resulting polynomial is
\[
 N_j(z)\coloneq \sum_{k=-q}^{-1}t_k
             \sum_{h=0}^{-k-1}r_hz^{h+k+q},
 \qquad \deg N_j<q.
\]
Thus the summed equation is
$Q(z;\lambda)F_j(z)=N_j(z)-z^{q+j}$.
Analyticity of $F_j$ requires the numerator on the right to vanish
at every root in $S$, with its multiplicity. Hence $N_j$ is exactly
the remainder when $z^{q+j}$ is divided by
$\prod_{u\in S}(z-u)$.

This conclusion includes multiple roots: if $u$ has multiplicity
$\mu$, an analytic quotient at $u$ requires the numerator and
its first $\mu-1$ derivatives to vanish there. Hence the whole
factor $(z-u)^\mu$ divides the numerator. The product of these
factors has degree $q$, and polynomial division has a unique
remainder of degree less than $q$.

To compute the quotient, use the complete homogeneous polynomials
of \eqref{eq:complete-homogeneous-definition} and expand at infinity:
\[
 \frac{z^{q+j}}{\prod_{u\in S}(z-u)}
 =z^j\prod_{u\in S}(1-u/z)^{-1}
 =\sum_{r\in\NN}\hc_r(S)z^{j-r}.
\]
Its polynomial part is $\sum_{k=0}^{j}\hc_{j-k}(S)z^k$.

Terms with $r\le j$ have nonnegative powers of $z$ and form the
quotient of polynomial division; terms with $r>j$ tend to zero
at infinity and are the remainder divided by the degree-$q$
divisor. Therefore $z^{q+j}-N_j$ is that divisor times the
stated polynomial part. This also fixes the minus sign in the
next expression for $F_j$.
Cancellation of the factor corresponding to $S$ therefore gives
\[
 F_j(z)=-\frac{\sum_{k=0}^{j}\hc_{j-k}(S)z^k}
 {t_p\prod_{v\in\mathcal Z}(z-v)}.
\]
Since the constant term of $Q$ is $\omega$, its roots satisfy
$\prod_{u\in S}u\prod_{v\in\mathcal Z}v=(-1)^m\omega/t_p$.
Consequently
\[
 -\frac1{t_p\prod_{v\in\mathcal Z}(z-v)}
 =\mathfrak g(\lambda)\prod_{v\in\mathcal Z}(1-z/v)^{-1}
 =\mathfrak g(\lambda)\sum_{r\in\NN}\hc_r(\mathcal Z^{-1})z^r.
\]
The coefficient of $z^i$ in its product with the quotient above
is the finite sum in \eqref{eq:general-R}.

Finally, the Neumann series for $(\lambda I-T)^{-1}$ and the
geometric expansion of $(\lambda-x)^{-1}$ converge for sufficiently
large $|\lambda|$. Their coefficients are, respectively, $T^\ell$
and the moments of $\Psi$. After selecting the first $q$ rows and
$p$ columns and multiplying by $H,K$, equation
\eqref{eq:general-moments} identifies those coefficients term by term.
This proves \eqref{eq:general-W}.
\end{proof}

For $|\lambda|$ sufficiently large, write $R\coloneq(\lambda I-T)^{-1}$.
The entry $(T^n)_{i,0}$ is the total weight of the $n$-step paths from
$0$ to $i$ that stay in $\NN$, with weight $t_k$ for a step of size $k$.
Thus $\lambda R_{i,0}$ is their length generating function evaluated at
$1/\lambda$. The formulas of Banderier and Flajolet
\cite[Theorem~2, equations~(19)--(20)]{BanderierFlajolet2002}
give this first column in terms of the kernel roots. Application to
$T^{\mathsf T}$ gives the first row. Deleting row and column zero from
$\lambda I-T$ gives the same Toeplitz matrix after reindexing, so the
Schur complement gives
\[
 R_{i,j}=R_{i-1,j-1}+\frac{R_{i,0}R_{0,j}}{R_{0,0}},\qquad i,j\ge1.
\]
Iteration recovers \eqref{eq:general-R}. Its initial rectangular block,
with the factors $H,K$ in \eqref{eq:general-W}, gives the transform whose
jump determines the Favard density in Theorem~\ref{thm:general-measure}.

The measure support considered next is also the limiting set of
eigenvalues of the finite Toeplitz sections. The classical results
of Schmidt--Spitzer and Ullman, together with Hirschman's limiting
eigenvalue distribution, are developed in
\cite[Chapter~11]{BottcherGrudsky2005}.
Hirschman's measure \cite{Hirschman1967} is the weak limit obtained
by giving each eigenvalue of the section of order $N$ mass $1/N$.
Duits and Kuijlaars \cite{DuitsKuijlaars2008} characterize it as
one component of a vector equilibrium problem.
This measure describes the distribution of eigenvalues and is
determined by $T$. The matrix measure $\Psi$ represents the moments
\eqref{eq:general-moments} and depends on the initial matrices $H,K$.
Example~\ref{ex:scalar-chebyshev} compares their densities explicitly.
For a nonsymmetric Toeplitz operator, the limiting eigenvalue set
can differ from the full spectrum on $\ell^2$.

The spectrum of the infinite operator is the set of $\lambda$ for
which the bounded inverse on $\ell^2$ fails to exist. In contrast,
the limiting eigenvalue set is formed from limits of eigenvalues of
increasing finite leading sections. The arguments below continue
entries of the Stieltjes transform as analytic functions; they do
not claim that the whole bounded operator resolvent exists at every
point reached by that continuation.
The root expressions for the Stieltjes transform are continued from
large $|\lambda|$, where the operator resolvent formula has been proved.

Retain the Laurent symbol $a(z)=z^{-q}P(z)$ from
\eqref{eq:general-toeplitz}. Order the root parameters in
\eqref{eq:root-parameters} as $\beta_1\ge\cdots\ge\beta_m>0$ and define
\begin{equation}\label{eq:F-critical}
 \mathcal F(z)\coloneq \sum_{\nu=1}^{m}\frac{\beta_\nu z}{1+\beta_\nu z}.
\end{equation}

\begin{theorem}[Root separation and a support bound]
\label{thm:general-support}
There is a unique $z_+>0$ with $\mathcal F(z_+)=q$. Put
$\lambda_+\coloneq a(z_+)$. If $\beta_q>\beta_{q+1}$, let
$z_-\in(-\beta_{q+1}^{-1},-\beta_q^{-1})$ be the unique point
with $\mathcal F(z_-)=q$, and put $\lambda_-\coloneq a(z_-)$. If
$\beta_q=\beta_{q+1}$, put $\lambda_-\coloneq 0$.
The limiting set of eigenvalues of the finite Toeplitz sections is
$[\lambda_-,\lambda_+]$. For the measure associated with any bounded
PBF of this Toeplitz matrix,
\begin{equation}\label{eq:general-support-inclusion}
 \supp\Psi\subseteq[\lambda_-,\lambda_+].
\end{equation}
For each $x\in(\lambda_-,\lambda_+)$, order the roots of
$Q(z;x)=P(z)-xz^q$ from \eqref{eq:general-kernel} by increasing
modulus. The roots in positions $q$ and $q+1$ form a nonreal conjugate pair
$\zeta(x),\overline{\zeta(x)}$, with $\im\zeta(x)<0$, and they are the
only roots on their common modulus circle. There are exactly $q-1$
roots of smaller modulus and $p-1$ of larger modulus.
\end{theorem}
\begin{proof}
\emph{The critical points.}
On the positive real axis,
$\mathcal F'(z)=\sum_{\nu=1}^m\beta_\nu/(1+\beta_\nu z)^2>0$,
and $\mathcal F$ increases from $0$ to $m$. Since $0<q<m$, it takes
the value $q$ exactly once. On the negative interval specified in
the theorem, its derivative is again positive and its limits at the
left and right endpoints are $-\infty$ and $+\infty$. This gives
the stated existence and uniqueness of $z_-$. Logarithmic
differentiation of $a=z^{-q}P$ gives
\begin{equation}\label{eq:a-critical}
 \frac{a'(z)}{a(z)}=\frac{\mathcal F(z)-q}{z}.
\end{equation}
For $z>0$, this derivative is negative before $z_+$ and positive
after $z_+$. Hence
$\lambda_+=a(z_+)=\min_{z>0}a(z)$.

The minimum is attained in the interior because $a(z)$ tends to
$+\infty$ both as $z\downarrow0$ and as $z\to+\infty$:
the respective dominant terms are $\omega z^{-q}$ and $t_pz^p$.
On the strict negative central gap, exactly $q$ factors
$1+\beta_\nu z$ are negative and all others are positive. The
endpoint poles of $\mathcal F$ therefore occur on either side
of the unique stationary point without a zero of $P$ in between.

The sign of $P(z_-)$ is $(-1)^q$, the same as that of $z_-^q$;
hence $\lambda_->0$ when the gap is strict. Also,
$a''(z_\pm)=a(z_\pm)\mathcal F'(z_\pm)/z_\pm$, so these critical
points are quadratic, with $a''(z_-)<0<a''(z_+)$.

\emph{The interval determined by equal root moduli.}
Let $r_+\coloneq z_+$. On $|z|=r_+$,
$|P(z)|\le P(r_+)=\lambda_+r_+^q$, with equality only at $z=r_+$.

For each factor,
$|1+\beta_\nu r_+\mathrm e^{\mathrm i\theta}|
\le1+\beta_\nu r_+$, with equality only when
$\theta=0$ modulo $2\pi$. Multiplication gives the stated strict
inequality everywhere else on the circle.
For $x>\lambda_+$, the term $-xz^q$ strictly dominates $P$ on
this circle. Rouch\'e's theorem therefore gives exactly $q$ zeros
of $Q(z;x)$ inside it and $p$ outside it. In particular, the roots
in positions $q$ and $q+1$, ordered by increasing modulus, have
different moduli.
When the central gap is strict, let $r_-\coloneq -z_-$. Then
$\beta_q^{-1}<r_-<\beta_{q+1}^{-1}$ and, factor by factor,
$|P(r_-\mathrm{e}^{\mathrm{i}\theta})|\ge |P(-r_-)|=\lambda_-r_-^q$,
with equality only at $\theta=\pi$.

This time each factor has its minimum modulus at the negative
point: its squared modulus is
$1+\beta_\nu^2r_-^2+2\beta_\nu r_-\cos\theta$.
None of these factors vanishes on the circle because $r_-$ lies
strictly between two consecutive root radii. Consequently the
product inequality is strict at every other point. For $0\le x<\lambda_-$,
the polynomial $P$ strictly dominates $xz^q$. It has exactly $q$
zeros inside the circle, namely $-\beta_1^{-1},\ldots,-\beta_q^{-1}$.
Rouch\'e again separates the first $q$ roots from the remaining $p$.
At $x=\lambda_+$ the two groups meet at the double root $z_+$;
at $x=\lambda_-$ they meet at the double root $z_-$ when the gap
is strict. Every other root lies strictly inside or outside the relevant
circle, by the strict modulus inequality away from that real point.

The double root splits into one root on each side of the circle
on the separating side of the endpoint. At $z_+$ this follows from
$a''(z_+)>0$ and $\lambda>\lambda_+$; at $z_-$ it follows from
$a''(z_-)<0$ and $\lambda<\lambda_-$. Since the separated side
has exactly $q$ roots inside, the double root occupies positions
$q$ and $q+1$ at the endpoint. This verifies that each displayed
endpoint belongs to the equal-modulus set.
If $\beta_q=\beta_{q+1}$, the $q$th
and $(q+1)$st roots of $P=Q(\,\cdot\,;0)$ already have equal modulus.

Let $\Sigma_{\mathrm{lim}}$ be the set of complex $\lambda$ for which those two
ordered moduli are equal. The Schmidt--Spitzer theorem
\cite{SchmidtSpitzer} identifies $\Sigma_{\mathrm{lim}}$ with the limiting set of
eigenvalues of the finite Toeplitz sections. Those sections are
BTP by Theorem~\ref{thm:toeplitz}, and
Corollary~\ref{cor:btp-oscillatory} makes their eigenvalues positive
and simple; thus $\Sigma_{\mathrm{lim}}\subseteq[0,\infty)$. Ullman's theorem
\cite{Ullman} states that this set is connected. The preceding
root counts give $\Sigma_{\mathrm{lim}}\subseteq[\lambda_-,\lambda_+]$ and put
both endpoints in $\Sigma_{\mathrm{lim}}$. A connected subset
of the real line containing both endpoints contains the interval
between them. Hence $\Sigma_{\mathrm{lim}}=[\lambda_-,\lambda_+]$.

\emph{The roots inside the interval and the measure support.}
Fix $x\in(\lambda_-,\lambda_+)$. Since $x\in\Sigma_{\mathrm{lim}}$, the roots in
positions $q,q+1$ have a common modulus $r$. On any circle $|z|=r$,
\[
 |P(r\mathrm{e}^{\mathrm{i}\theta})|^2=\omega^2\prod_{\nu=1}^{m}
 (1+\beta_\nu^2r^2+2\beta_\nu r\cos\theta),
\]
which is strictly increasing in $\cos\theta$.

For $-1<\cos\theta\le1$, every factor is positive and strictly
increasing. At $\cos\theta=-1$ a factor may vanish, but its value
becomes positive immediately as $\cos\theta$ increases. Thus the
whole product is strictly increasing on the entire interval
$[-1,1]$. Equal moduli of $Q$-roots on a fixed circle therefore
force equal values of $\cos\theta$, allowing only a conjugate pair. At a root of $Q(z;x)$
on this circle its value must be $x^2r^{2q}$. Thus at most two
distinct roots can have modulus $r$, and, when nonreal, they are
conjugates. Every nonreal root is simple for $x>0$:
a multiple root would satisfy $\mathcal F(z)=q$, whereas
$\im\mathcal F(z)=\im z\sum_{\nu=1}^{m} \frac{\beta_\nu}{|1+\beta_\nu z|^2}$
is nonzero for nonreal $z$.

There is no positive root, since $a(z)\ge\lambda_+>x$ for $z>0$.
Two distinct negative roots cannot have the same modulus. Thus the
two middle roots are either a nonreal conjugate pair or a repeated
negative root $c$. In the latter case, $c$ is not a zero of $P$
because $x>0$, and \eqref{eq:a-critical} gives $\mathcal F(c)=q$.
Such a point must lie in a strict gap between consecutive distinct
roots of $P$: on $(-\infty,-\beta_m^{-1})$ each summand in
$\mathcal F$ is greater than $1$, so $\mathcal F>m>q$; on
$(-\beta_1^{-1},0)$ each summand is negative, so $\mathcal F<0<q$.
In each strict gap, $\mathcal F'>0$ and its endpoint limits are
$-\infty,+\infty$, so there is exactly one such critical point.

Let $k$ count the zeros of $P$ in $(c,0)$, with multiplicity.
Since $c$ is a root of $Q(z;x)$ for the fixed $x>0$,
$x_*\coloneq a(c)=x>0$.
On the circle $|z|=-c$, the product formula gives
$|P(z)|\ge |P(c)|=x_*|c|^q$,
with equality only at $z=c$. For $0<y<x_*$ sufficiently close to
$x_*$, Rouch\'e's theorem therefore counts exactly $k$ roots of
$Q(z;y)$ inside this circle. Moreover,
$a''(c)=a(c)\mathcal F'(c)/c<0$, so the two roots near $c$ are real
and have the form
\[
 c\pm\sqrt{\frac{2(x_*-y)}{-a''(c)}}
   +\mathrm O(x_*-y),\qquad y\uparrow x_*.
\]
One lies inside the circle and the other outside. No other root
is on the circle at $y=x_*$. Thus the colliding roots occupy
positions $k,k+1$ in increasing modulus. They are the middle roots
only when $k=q$, which gives $x_*=\lambda_-$, outside the open
interval under consideration. If $\beta_q=\beta_{q+1}$, there is
no strict gap with $k=q$, so no negative collision can involve the
middle roots at any $x>0$. This excludes the remaining real case.
The two middle roots are therefore nonreal and conjugate throughout
the open interval. They are simple and are the only roots on their
modulus circle, so exactly $q-1$ roots have smaller modulus and
$p-1$ have larger modulus.

The root formula admits analytic continuation, including through
collisions inside a group. At any $\lambda_0$ outside
$\Sigma_{\mathrm{lim}}$, choose a circle whose radius lies strictly
between the $q$th and $(q+1)$st root moduli. By continuity of the
roots as multisets, this same circle separates the two groups for
$\lambda$ in a neighborhood of $\lambda_0$. For every nonnegative
integer $h$, the power sum of the roots inside it is
\[
 \sum_{u\in S}u^h=
 \frac{1}{2\pi\mathrm i}\int_{|z|=r}
       z^h\frac{Q_z(z;\lambda)}{Q(z;\lambda)}\,\mathrm dz.
\]
The integrand is analytic in $\lambda$ on the contour, so these
power sums are analytic even when enclosed roots coincide.
Newton's identities then make the elementary symmetric functions,
and hence the complete homogeneous polynomials, analytic in
$\lambda$. The same reasoning applies to the outside group,
using the reciprocal polynomial to work with its reciprocal roots.
The product defining $\mathfrak g$ is another elementary symmetric
function of the inside group. Thus every finite expression in
\eqref{eq:general-R} is analytic off $\Sigma_{\mathrm{lim}}$.

These locally defined expressions agree on overlaps because each
uses the same group of the $q$ smallest root moduli. They therefore
define a single analytic function there without a choice of labels
for individual roots. The matrix obtained by the finite selection
and multiplication in \eqref{eq:general-W} agrees with $W$ for
large $|\lambda|$. Both are analytic in each open half-plane:
$W$ is the Stieltjes transform of measures on the real line, and
the root formula has no group collision off
$\Sigma_{\mathrm{lim}}\subset\mathbb R$. The identity theorem
extends their equality to both half-planes.

For real $\lambda$ outside this interval, each separated group is
invariant under conjugation. Its symmetric functions, and hence
the continued matrix formula, are real. On a closed real interval
in the complement, the formula extends analytically through the
axis. Its imaginary boundary values are zero and locally uniform;
Stieltjes inversion gives zero measure on that interval. Boundedness
of this analytic extension also excludes an atom at any of its
points. Covering the complement by such intervals gives
\eqref{eq:general-support-inclusion}. This establishes the support bound.
The interval cannot collapse to a single point. Indeed, the first
entry of the measure has total mass one and first moment $t_0$,
by \eqref{eq:general-moments} and the first row and column of $H,K$.
Its variance is
$(T^2)_{0,0}-t_0^2=\sum_{j=1}^{\min(p,q)}t_{-j}t_j>0$.
Every summand is positive because
\eqref{eq:root-parameters} has positive coefficients and $p,q\ge1$.
A probability measure supported at one point has zero variance,
so $\lambda_-<\lambda_+$. The next theorem determines
the density and proves equality of the support with this interval.
\end{proof}

The orientation of the boundary roots follows from the quadratic
critical point at $z_+$. For real $\lambda>\lambda_+$ sufficiently
close to $\lambda_+$, the smaller of the two roots near $z_+$ belongs
to $S$. Its local expansion is
\[
 z_{\rm s}(\lambda)=z_+
 -\sqrt{\frac{2}{a''(z_+)}}\sqrt{\lambda-\lambda_+}
 +\mathrm O(\lambda-\lambda_+),\qquad \lambda\to\lambda_+,
\]
where the square root is positive for real $\lambda>\lambda_+$.
On approaching $x<\lambda_+$ from the upper half-plane, that square
root tends to $\mathrm i\sqrt{\lambda_+-x}$.
Consequently the root in $S$ has negative imaginary part just below
$\lambda_+$.

Here the square root is continued from its positive value on the
positive real axis through the upper half-plane. Thus a negative
real argument has boundary value with argument $\pi/2$, rather
than $-\pi/2$. The preceding minus sign in $z_{\rm s}$ then
selects the root with negative imaginary part. By Theorem~\ref{thm:general-support}, the middle pair
remains simple, nonreal, and separated from the other roots in the
open interval. The choice therefore remains $\zeta(x)$ throughout
that interval.

A change to the conjugate root would require a collision of the
pair or passage through the real axis. Both are excluded in the
open interval by Theorem~\ref{thm:general-support}. Continuation
along overlapping neighborhoods of its points therefore preserves
the selection. For $x\in(\lambda_-,\lambda_+)$, define the root
multisets
\begin{equation}\label{eq:general-density-root-sets}
 \begin{aligned}
 \mathcal S_0&\coloneq\{u:Q(u;x)=0,\ |u|<|\zeta(x)|\},&
 \mathcal Z_0&\coloneq\{v:Q(v;x)=0,\ |v|>|\zeta(x)|\}.
 \end{aligned}
\end{equation}
Roots are counted with multiplicity: $\mathcal S_0$ has $q-1$
members and $\mathcal Z_0$ has $p-1$. The boundary values of
$S(\lambda),\mathcal Z(\lambda)$ from the upper half-plane, and the
two enlarged multisets used for the density, are
\begin{equation}\label{eq:general-root-boundary-sets}
 \begin{aligned}
 S_+&=\mathcal S_0\cup\{\zeta\},&
 \mathcal Z_+&=\mathcal Z_0\cup\{\bar\zeta\},&
 \widehat{\mathcal S}&\coloneq\mathcal S_0\cup\{\zeta,\bar\zeta\},&
 \widehat{\mathcal Z}&\coloneq\mathcal Z_0\cup\{\zeta,\bar\zeta\}.
 \end{aligned}
\end{equation}
The dependence of all these multisets on $x$ is suppressed in the
notation; unions include multiplicities. Using $\hc_\ell$ from
\eqref{eq:complete-homogeneous-definition}, define
\begin{equation}\label{eq:general-density-data}
\begin{gathered}
 \rho(x)\coloneq \frac{-\im\zeta(x)}{\pi\omega}
             \prod_{u\in \mathcal S_0}|u|,\\
 f(x)\coloneq H[\hc_i(\widehat{\mathcal Z}^{-1})]_{i\in\{0,\dots,q-1\}},\qquad
 g(x)^{\mathsf T}\coloneq [\hc_j(\widehat{\mathcal S})]_{j\in\{0,\dots,p-1\}}K.
\end{gathered}
\end{equation}

\begin{theorem}[Algebraic rank-one measure]
\label{thm:general-measure}
For any bounded PBF of the Toeplitz matrix in
\eqref{eq:general-toeplitz}, let $H,K$ be its initial-condition
matrices \eqref{eq:normalized-boundary-definition}. Use the endpoints
$\lambda_-,\lambda_+$ and root $\zeta$ of
Theorem~\ref{thm:general-support}, and form $\rho,f,g$ by
\eqref{eq:general-density-data} with the multisets in
\eqref{eq:general-density-root-sets}--\eqref{eq:general-root-boundary-sets}.
The measure $\Psi$ of Corollary~\ref{cor:general-favard} is purely
absolutely continuous, has
\begin{equation}\label{eq:general-support-interval}
 \supp\Psi=[\lambda_-,\lambda_+],
\end{equation}
and satisfies
\begin{equation}\label{eq:general-weight}
 \dd\Psi(x)=\rho(x)f(x)g(x)^{\mathsf T}\dd x,
 \qquad \lambda_-<x<\lambda_+.
\end{equation}
Here $\rho>0$, while every entry of $f$ and $g$ is nonnegative. The
density has rank exactly one, and its upper-left entry is $\rho$.
\end{theorem}
\begin{proof}
\emph{The jump of the Stieltjes transform.}
First sum \eqref{eq:general-R} over both indices. Reindexing its
finite inner sum by $i=k+r$ and $j=k+\ell$ separates the three sums
and gives, for formal variables $s,t$,
\[
 \sum_{i,j\in\NN}[(\lambda I-T)^{-1}]_{i,j}s^it^j
 =\frac{\mathfrak g(\lambda)}{1-st}
  \prod_{u\in S}(1-ut)^{-1}
  \prod_{v\in\mathcal Z}(1-s/v)^{-1}.
\]
Across the interval, the smaller and larger root multisets
$\mathcal S_0,\mathcal Z_0$ of \eqref{eq:general-density-root-sets}
are unchanged, while $\zeta$ and $\bar\zeta$ exchange groups as in
\eqref{eq:general-root-boundary-sets}. The elementary identity
\begin{equation}\label{eq:general-jump-identity}
 \frac{\zeta}{(1-\zeta t)(1-s/\bar\zeta)}
 -\frac{\bar\zeta}{(1-\bar\zeta t)(1-s/\zeta)}
 =\frac{(\zeta-\bar\zeta)(1-st)}
 {(1-\zeta t)(1-\bar\zeta t)
  (1-s/\zeta)(1-s/\bar\zeta)}
\end{equation}
cancels the factor $(1-st)^{-1}$. The remaining generating function
is a product of a function of $s$ and a function of $t$.

In detail, put
$\chi(x)\coloneq (-1)^{q+1}\omega^{-1}\prod_{u\in\mathcal S_0}u$.
The sign computation below shows that this scalar is positive.
The difference between the upper and lower boundary values of the
double generating function is
\[
 \chi(x)(\zeta-\bar\zeta)
 \prod_{u\in\widehat{\mathcal S}}(1-ut)^{-1}
 \prod_{v\in\widehat{\mathcal Z}}(1-s/v)^{-1}.
\]
Each coefficient of $s^it^j$ is therefore $\chi(x)(\zeta-\bar\zeta)$
times the product of the $i$th coefficient in the second product
and the $j$th coefficient in the first. No infinite series needs
to be evaluated on the support: the computation can be made in
formal variables $s,t$, coefficient by coefficient. Its coefficients
give a matrix of rank one, with column
$[\hc_i(\widehat{\mathcal Z}^{-1})]$ and row $[\hc_j(\widehat{\mathcal S})]$.

To determine the scalar factor and its sign, observe that $\mathcal S_0$ is
invariant under conjugation. For $z>0$ and
$x\in(\lambda_-,\lambda_+)$, the minimum computed in the proof of
Theorem~\ref{thm:general-support} gives $a(z)\ge\lambda_+>x$. Since $Q(0;x)=\omega>0$, every real root of
$Q(z;x)$ is therefore negative. Conjugate pairs in $\mathcal S_0$ have positive
products. As $\mathcal S_0$ contains $q-1$ roots,
$\prod_{u\in\mathcal S_0}u=(-1)^{q-1}\prod_{u\in\mathcal S_0}|u|$.
The scalar factor $\mathfrak g$ in \eqref{eq:general-gamma} consequently has
upper boundary value
$\mathfrak g_+(x)=\frac{\zeta(x)}{\omega}\prod_{u\in\mathcal S_0}|u|$.
The lower boundary value is the conjugate of the upper one. Hence
the density contributed by the scalar jump is
$-(2\pi\mathrm{i})^{-1}(\mathfrak g_+-\mathfrak g_-)
=-\pi^{-1}\im\mathfrak g_+=\rho>0$.
Premultiplication by $H$ and postmultiplication by $K$ turn the two
coefficient vectors into $f$ and $g^{\mathsf T}$ as defined in
\eqref{eq:general-density-data}.
Entrywise Stieltjes inversion now gives \eqref{eq:general-weight}
on the open interval.

\emph{Positivity, rank, and support.}
Since $H$ is lower triangular and $K$ is upper triangular, both with
diagonal entries equal to one, their first row and first column,
respectively, select the constant
complete homogeneous polynomial $\hc_0=1$. Thus the first components
of $f$ and $g$ are one. Corollary~\ref{cor:general-favard}
constructs $\Psi$ as a limit of positive finite spectral measures,
so its density entries are nonnegative almost everywhere. Those
entries are continuous in the open interval, and are consequently
nonnegative at every point there.

Continuity follows from the simple separated central pair and the
symmetric dependence on the remaining two groups, established in
the support proof. If a density entry were negative at a point,
it would remain negative on a small interval by continuity, giving
negative measure to that interval and contradicting this
positivity. The first column of the density is
$\rho f$ and its first row is $\rho g^{\mathsf T}$; since $\rho>0$,
this proves $f,g\ge0$. Since $f_1=g_1=1$, both vectors are nonzero,
so the rank is exactly one. The first density entry
is strictly positive on the whole open interval. Every open subinterval
therefore has positive measure. Taking the closure and combining this
with \eqref{eq:general-support-inclusion} proves
\eqref{eq:general-support-interval}.

\emph{Absolute continuity on the full support.}
At a coincidence of roots inside one of the groups, the complete
homogeneous polynomials in \eqref{eq:general-R} are symmetric in
those roots, so the formula is invariant under permutation of their
individual branches. At a coincidence between the two groups, both limiting
formulas remain finite. Indeed, near any finite $\lambda_0$ the
coefficients of $Q$ are bounded and its leading and constant
coefficients are fixed and nonzero. Its roots are therefore uniformly
bounded and bounded away from zero.

For the first bound apply the usual polynomial root bound after
dividing $Q$ by its fixed leading coefficient. For the lower bound,
apply the same bound to $z^mQ(1/z;\lambda)$, whose leading
coefficient is the fixed nonzero constant term of $Q$.
Both bounds are uniform on a small closed neighborhood of
$\lambda_0$, since all coefficients are bounded there. The finite sums in
\eqref{eq:general-R}, and hence every entry of $W$, remain bounded.

Let $[a,b]$ be a closed interval avoiding the finitely many algebraic
branch values. The branch of each entry of $W$ obtained from the
upper half-plane extends analytically to a neighborhood of this
interval. Thus $W(x+\mathrm i\epsilon)$ converges uniformly to
its boundary value on $[a,b]$ as $\epsilon\downarrow0$.
The Stieltjes inversion formula on subintervals then identifies
the whole restriction of the measure with the density
$-\pi^{-1}\im W(x+\mathrm i0)$; no singular part remains there.
The exceptional branch values are finite because they are among
the zeros of the nonzero discriminant discussed above. An algebraic
entry of $W$ has analytic branches locally away from that set;
there the integrated inversion formula already accounts for its
entire measure. Any singular remainder is therefore supported on
a finite set. A finite measure supported on a finite set is the
sum of its point masses, so only atoms remain to be excluded.
For each scalar entry and each real $x_0$, an atom is detected by
$\Psi_{b,a}(\{x_0\}) =\lim_{\epsilon\downarrow0} [-\epsilon\im W_{b,a}(x_0+\mathrm{i}\epsilon)]$.
Boundedness of $W$ makes this limit zero.

For a positive scalar measure the displayed atom formula follows
by dominated convergence from
\[
 -\epsilon\operatorname{Im}W_{b,a}(x_0+\mathrm i\epsilon)
 =\int\frac{\epsilon^2}{(x_0-x)^2+\epsilon^2}
                  \,\mathrm d\Psi_{b,a}(x).
\]
The integrand is between zero and one and tends to the indicator
of $\{x_0\}$. Thus the limit is exactly the point mass, including
when $x_0$ is an endpoint where several roots meet. The measure is therefore
purely absolutely continuous on its full support.
\end{proof}

\begin{remark}[Computing the endpoints and choosing the density root]
When $\beta_q>\beta_{q+1}$, the strict monotonicity of the function
$\mathcal F$ defined in \eqref{eq:F-critical} permits bisection inside
the negative interval in
Theorem~\ref{thm:general-support}. The resulting critical points
and endpoints satisfy
$-\beta_{q+1}^{-1}<z_-<-\beta_q^{-1}<0<z_+$ and
$0\le\lambda_-<\lambda_+$.
The case $\beta_q=\beta_{q+1}$ instead uses $\lambda_-=0$ directly.

If the roots are stored in increasing modulus with zero-based indices,
the central conjugate pair occupies indices $q-1,q$, and $\zeta$
is its member with negative imaginary part.
The smaller and larger sets contain exactly $q-1$ and $p-1$ roots.
For $q=1$ the product over the smaller roots in \eqref{eq:general-density-data}
is empty and equals one, giving
$\rho(x)=-\im\zeta(x)/(\pi\omega)>0$.
\end{remark}

\begin{proposition}[Density near the endpoints]
\label{prop:general-edges}
Use the endpoint data $\lambda_\pm,z_\pm$ of
Theorem~\ref{thm:general-support}, with $z_-$ defined when
$\beta_q>\beta_{q+1}$, the Laurent symbol
$a(z)=z^{-q}P(z)$ in \eqref{eq:general-toeplitz}, and the scalar
density $\rho$ in \eqref{eq:general-density-data}.
Let $\mathcal S_{0,+}$ be the limit of the smaller-root multiset
$\mathcal S_0$ in \eqref{eq:general-density-root-sets} as
$x\uparrow\lambda_+$. For every positive root-parameter multiset,
including coincident parameters,
\begin{equation}\label{eq:general-upper-edge}
 \rho(x)\sim C_+\sqrt{\lambda_+-x},\qquad
 C_+\coloneq\frac{\prod_{u\in\mathcal S_{0,+}}|u|}{\pi\omega}
          \sqrt{\frac{2}{a''(z_+)}},
 \qquad x\uparrow\lambda_+.
\end{equation}
If $\beta_q>\beta_{q+1}$, let $\mathcal S_{0,-}$ denote the
corresponding root set at the lower endpoint. Then
\begin{equation}\label{eq:general-lower-edge}
 \rho(x)\sim C_-\sqrt{x-\lambda_-},\qquad
 C_-\coloneq\frac{\prod_{u\in\mathcal S_{0,-}}|u|}{\pi\omega}
          \sqrt{\frac{2}{-a''(z_-)}},
 \qquad x\downarrow\lambda_-.
\end{equation}
Denote the endpoint limits of the vectors in
\eqref{eq:general-density-data} by $f_+,g_+$ and, when
$\beta_q>\beta_{q+1}$, by $f_-,g_-$. The corresponding matrix limits
of the density are
\[
 \begin{aligned}
 \frac{1}{\sqrt{\lambda_+-x}}\,\frac{\dd\Psi}{\dd x}(x)
 &\longrightarrow C_+f_+g_+^{\mathsf T}
 &&(x\uparrow\lambda_+),\\
 \frac{1}{\sqrt{x-\lambda_-}}\,\frac{\dd\Psi}{\dd x}(x)
 &\longrightarrow C_-f_-g_-^{\mathsf T}
 &&(x\downarrow\lambda_-,\ \beta_q>\beta_{q+1}).
 \end{aligned}
\]
An entrywise asymptotic equivalent follows whenever the corresponding
entry of $f_\pm g_\pm^{\mathsf T}$ is nonzero. If it is zero, the
limit gives a smaller order and does not specify the first nonzero term.
If instead the root
$-\beta_q^{-1}$ belongs to a central block of multiplicity $d\ge2$,
then $\lambda_-=0$ and $\rho(x)\asymp x^{1/d}$ as $x\downarrow0$;
here $\asymp$ means that the ratio is bounded above and below by
positive constants near the stated limit. Individual matrix
entries can vanish to higher orders.
\end{proposition}
\begin{proof}
At a critical point, differentiation of \eqref{eq:a-critical} gives
$a''(z_*)=a(z_*)\frac{\mathcal F'(z_*)}{z_*}$.
The upper critical point exists for every root multiset, and
$a''(z_+)>0$. When the central gap is strict, the lower critical
point also exists and $a''(z_-)<0$. At each available critical
point, Taylor expansion gives
$a(z)=\lambda_\pm+\tfrac12 a''(z_\pm)(z-z_\pm)^2
+\mathrm O((z-z_\pm)^3)$ as $z\to z_\pm$.
Choose in each case the root with negative imaginary part. Inverting
these expansions yields
\[
 \begin{aligned}
 \zeta(x)&=z_--\mathrm{i}\sqrt{\frac{2}{-a''(z_-)}}
           \sqrt{x-\lambda_-}+\mathrm O(x-\lambda_-),
           &&x\downarrow\lambda_-,\\
 \zeta(x)&=z_+-\mathrm{i}\sqrt{\frac{2}{a''(z_+)}}
           \sqrt{\lambda_+-x}+\mathrm O(\lambda_+-x),
           &&x\uparrow\lambda_+.
 \end{aligned}
\]
At the upper endpoint the remaining roots stay separated from
$z_+$ and have finite nonzero limits. The same holds at the lower
endpoint when the central gap is strict. Substitution into
\eqref{eq:general-density-data} proves the corresponding scalar
equivalents; the limits of $f$ and $g$ give the matrix behavior.
In particular, the statement at the upper endpoint is unaffected by
coincident roots of $P$, whose collisions occur at $\lambda=0$.

For the multiple-root case write $z_0\coloneq -\beta_q^{-1}$ and
$a(z)=C(z-z_0)^d[1+\mathrm O(z-z_0)]$ as $z\to z_0$, with real
$C\ne0$. Let $s$ be the number of root parameters strictly larger
than $\beta_q$, and put $k\coloneq q-s\in\{1,\dots,d-1\}$. To determine the sign of $C$, note that the $s$ factors coming
from larger parameters are negative at $z_0$, the factors from
smaller parameters are positive, and $z_0^{-q}$ has sign $(-1)^q$.
The repeated block contributes the positive coefficient
$\beta_q^d$ of $(z-z_0)^d$. Thus the sign of
$C$ is $(-1)^{s-q}=(-1)^k$.

Put $r_0\coloneq -z_0>0$ and $t\coloneq |C|^{-1/d}x^{1/d}$.
After writing $z=z_0+t\alpha$ and dividing $a(z)=x$ by $x$,
the leading equation is $\alpha^d=(-1)^k$.
Its $d$ roots are distinct, so the analytic implicit function
theorem applied to this scaled equation gives displacements
$z-z_0=t\mathrm e^{\mathrm i\theta}+\mathrm O(t^2)$ as $t\downarrow0$, with $d\theta\equiv k\pi\pmod{2\pi}$.
Their moduli satisfy
$|z_0+t\mathrm e^{\mathrm i\theta}+\mathrm O(t^2)|
=r_0-t\cos\theta+\mathrm O(t^2)$ as $t\downarrow0$.
Thus the increasing order of root moduli is, to first order,
the decreasing order of $\cos\theta$.
There are exactly $k-1$ allowed angles strictly between
$-\pi k/d$ and $\pi k/d$: they are
$\theta=(k-2j)\pi/d$ for $j\in\{1,\dots,k-1\}$.
The next two angles are $\pm\pi k/d$, with equal cosine;
all remaining allowed angles have smaller cosine. The strict
inequalities between these cosine values persist for small $t$.

Outside this block, precisely $s$ roots have smaller limiting
modulus, because they correspond to the strictly larger root
parameters. Within the block the roots in positions $k,k+1$
therefore occupy positions $q,q+1$ in the whole list. This is the
conjugate pair with leading arguments $\pm\pi k/d$, and its
negative-imaginary member has argument $-\pi k/d$.
In particular,
$-\im\zeta(x)\sim |C|^{-1/d}\sin(\pi k/d)x^{1/d}$
as $x\downarrow0$.
Here $\sin(\pi k/d)>0$. As $x\downarrow0$, the other factors in
$\rho$ have positive finite limits. Hence
$\rho(x)\asymp x^{1/d}$ as $x\downarrow0$. Some components
of $f$ or $g$ can tend to zero as $x\downarrow0$; their products with $\rho$ can
therefore have a higher order of vanishing.
\end{proof}

\Needspace{12\baselineskip}
\begin{corollary}[Factorization stability and endpoints of the orthogonality support]
\label{cor:reconstruction-support-edges}
Let $T$ be the BTP Toeplitz matrix with symbol
$a(z)=\omega z^{-q}\prod_{\nu=1}^{p+q}(1+\beta_\nu z)$, where
$p,q\ge1$, $\omega>0$, and
$\beta_1\ge\cdots\ge\beta_{p+q}>0$.
Let $[\lambda_-,\lambda_+]$ be the support in
Theorem~\ref{thm:general-measure} and let $\rho$ be the scalar
density defined by \eqref{eq:general-density-data}.
Interior and boundary below refer to the PBF class in operator
norm among bounded real matrices of bandwidth at most $(p,q)$.
\begin{enumerate}
\item If $T$ is an interior point, then $\lambda_->0$ and $\rho$
has the square-root equivalents
\eqref{eq:general-upper-edge} and \eqref{eq:general-lower-edge}
at the two endpoints.
\item If $\beta_q=\beta_{q+1}$ and this parameter has multiplicity
$d$, then $T$ is a boundary point, $\lambda_-=0$, and
$\rho(x)\asymp x^{1/d}$ as $x\downarrow0$.
\item If $\beta_q>\beta_{q+1}$ and some parameters coincide,
then $T$ is a boundary point with $\lambda_->0$ and the same
square-root equivalents at both endpoints.
\end{enumerate}
\end{corollary}
\begin{proof}
Theorem~\ref{thm:toeplitz-stable-simple-roots} identifies the
interior points in this Toeplitz class with the matrices whose
root parameters are all distinct. Such parameters have a strict
central gap $\beta_q>\beta_{q+1}$. Theorem~\ref{thm:general-support}
then gives $\lambda_->0$, and
Proposition~\ref{prop:general-edges} gives the two square-root
equivalents. This proves (1).

In (2) and (3), repeated parameters exclude membership in the
interior by the same stability theorem. Theorem~\ref{thm:toeplitz}
supplies a bounded PBF for $T$, so $T$ belongs to the PBF class
itself and is consequently a boundary point. A central repetition
gives $\lambda_-=0$ and the exponent $1/d$ by
Theorem~\ref{thm:general-support} and
Proposition~\ref{prop:general-edges}, proving (2). When the central
gap remains strict, those results instead give a positive lower
endpoint and square-root behavior at both ends, proving (3).
\end{proof}

The central gap determines whether the orthogonality support reaches
zero, while an analytically varying PBF requires separation of every
pair of root parameters. For instance, take $p=3$, $q=2$, $\omega=1$
and $(\beta_1,\ldots,\beta_5)=(2,2,1,1,1)$, so
$a(z)=(1+2z)^2(1+z)^3/z^2$.
Since $\mathcal F(-2/3)=2$ and $-2/3\in(-1,-1/2)$,
Theorem~\ref{thm:general-support} gives $z_-=-2/3$ and
$\lambda_-=a(-2/3)=1/108>0$.
The gap $\beta_2>\beta_3$ is strict, but the repeated parameters
make $T$ a boundary point of the PBF class. Its density $\rho$
still has the square-root equivalents at both endpoints.

\begin{example}[Scalar Chebyshev polynomials for arbitrary roots]
\label{ex:scalar-chebyshev}
The scalar case makes the Chebyshev connection explicit for arbitrary
positive root parameters. For $p=q=1$ and
$P(z)=\omega(1+\eta z)(1+\xi z)$, the initial matrices are $H=K=[1]$. Set
\begin{equation}\label{eq:scalar-chebyshev-variable}
 y\coloneq \frac{x-\omega(\eta+\xi)}{2\omega\sqrt{\eta\xi}},
 \qquad \lambda_\pm=\omega(\sqrt\eta\pm\sqrt\xi)^2.
\end{equation}
Then
\begin{equation}\label{eq:scalar-chebyshev-polynomials}
 B_n(x)=(\eta\xi)^{n/2}U_n(y),\qquad
 A_n(x)=(\eta\xi)^{-n/2}U_n(y),
\end{equation}
where $U_n$ is the Chebyshev polynomial of the second kind, and
\begin{equation}\label{eq:scalar-chebyshev-measure}
 \dd\Psi(x)
 =\frac{\sqrt{(\lambda_+-x)(x-\lambda_-)}}
 {2\pi\omega^2\eta\xi}\,\dd x
 =\frac{2}{\pi}\sqrt{1-y^2}\,\dd y,
 \qquad -1<y<1.
\end{equation}
Indeed, the diagonal coefficients are
$t_{-1}=\omega$, $t_0=\omega(\eta+\xi)$ and
$t_1=\omega\eta\xi$. The right recurrence is
\[
 \omega B_{n+1}+\omega(\eta+\xi)B_n
       +\omega\eta\xi B_{n-1}=xB_n,
 \qquad B_{-1}=0,\quad B_0=1.
\]
Substituting $B_n=(\eta\xi)^{n/2}U_n(y)$ and dividing by
$\omega(\eta\xi)^{(n+1)/2}$ gives
$U_{n+1}=2yU_n-U_{n-1}$, with $U_0=1$ and $U_1=2y$.
The left recurrence has the two off-diagonal coefficients reversed,
which gives the reciprocal scale $(\eta\xi)^{-n/2}$.
The root with negative imaginary part is
\[
 \zeta(x)=\frac{x-\omega(\eta+\xi)
       -\mathrm i\sqrt{(\lambda_+-x)(x-\lambda_-)}}
                  {2\omega\eta\xi}.
\]
Since $q=1$, the product in $\rho$ is empty and $f=g=1$.
Hence $\rho=-\operatorname{Im}\zeta/(\pi\omega)$ gives the
first density formula. The Jacobian
$\dd x=2\omega\sqrt{\eta\xi}\,\dd y$ gives the last expression.
For $\omega=\eta=\xi=1$, both families reduce to
$U_n((x-2)/2)$. In this case $P(z)=(1+z)^2$ and the two densities
on $(0,4)$ are $\sqrt{x(4-x)}/(2\pi)$ for the orthogonality measure
and $1/(\pi\sqrt{x(4-x)})$ for the limiting eigenvalue distribution
\cite{DuitsKuijlaars2008}.
For this last comparison, the finite section of order $N$ has
eigenvalues $2+2\cos(j\pi/(N+1))$, $j\in\{1,\dots,N\}$.
They follow by inserting the vectors with entries
$\sin((n+1)j\pi/(N+1))$ into the tridiagonal recurrence and
checking the zero boundary values at $n=-1,N$.
The normalized counting measure therefore converges to the image
of $\mathrm d\theta/\pi$ under $x=2+2\cos\theta$,
$0<\theta<\pi$. Since
$|\mathrm dx/\mathrm d\theta|=2\sin\theta=\sqrt{x(4-x)}$,
its density is $1/(\pi\sqrt{x(4-x)})$. This computation also
shows directly that its total mass is one.
\end{example}

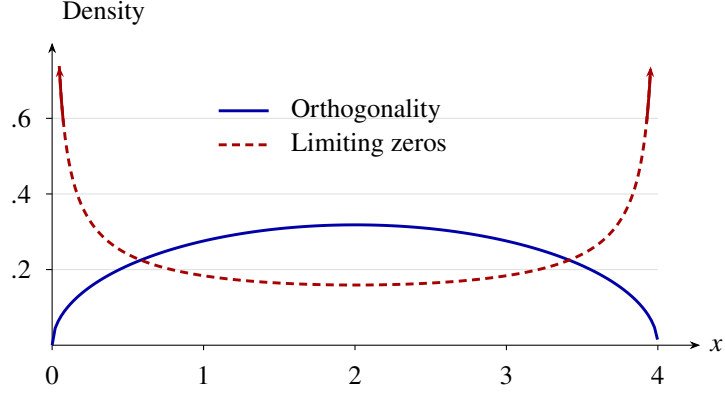
\begin{figure}[htbp]
\centering
\begin{tikzpicture}[x=2cm,y=5cm,font=\small,
 orth/.style={blue!65!black,line width=1.1pt},
 zeros/.style={red!65!black,line width=1.1pt,densely dashed}]
 \draw[gray!25] (0,.2) -- (4,.2);
 \draw[gray!25] (0,.4) -- (4,.4);
 \draw[gray!25] (0,.6) -- (4,.6);
 \draw[-{Stealth[length=4pt]}] (0,0) -- (4.28,0)
   node[right] {$x$};
 \draw[-{Stealth[length=4pt]}] (0,0) -- (0,.80);
 \node[anchor=south west] at (0,.82) {Density};
 \foreach \x in {0,1,2,3,4}{
   \draw (\x,0) -- (\x,-.013);
   \node[below=4pt] at (\x,0) {$\x$};
 }
 \foreach \y in {.2,.4,.6}{
   \draw (0,\y) -- (-.04,\y);
   \node[left=4pt] at (0,\y) {$\y$};
 }
 \draw[orth,domain=0:4,samples=201]
   plot (\x,{sqrt(\x*(4-\x))/(2*pi)});
 \draw[zeros,domain=.047:3.953,samples=301]
   plot (\x,{1/(pi*sqrt(\x*(4-\x)))});
 \draw[red!65!black,line width=1.1pt,-{Stealth[length=4pt]},
   domain=.075:.047,samples=12]
   plot (\x,{1/(pi*sqrt(\x*(4-\x)))});
 \draw[red!65!black,line width=1.1pt,-{Stealth[length=4pt]},
   domain=3.925:3.953,samples=12]
   plot (\x,{1/(pi*sqrt(\x*(4-\x)))});
 \draw[orth] (1.10,.62) -- (1.43,.62);
 \node[anchor=west] at (1.51,.62) {Orthogonality};
 \draw[zeros] (1.10,.53) -- (1.43,.53);
 \node[anchor=west] at (1.51,.53) {Limiting zeros};
\end{tikzpicture}
\caption{The scalar case $P(z)=(1+z)^2$ in
Example~\ref{ex:scalar-chebyshev}. The orthogonality density
$\sqrt{x(4-x)}/(2\pi)$ (solid) vanishes at the endpoints.
The limiting zero density $1/(\pi\sqrt{x(4-x)})$ (dashed)
diverges there, as indicated by the arrows. Both have total mass
one and support $[0,4]$.}
\label{fig:orthogonality-zero-densities}
\end{figure}

Figure~\ref{fig:orthogonality-zero-densities} compares these two
probability densities on their common support.

\subsection{Chebyshev formulas, zeros, and spectral kernels}
\label{subsec:spectral-kernels}

The mixed-type orthogonality relations involve weighted linear
forms in the polynomial components. The rank-one density of
Theorem~\ref{thm:general-measure} gives two such forms, one from
each polynomial family. These forms depend only on $T$ and admit
exact Chebyshev convolutions for arbitrary bandwidths. Separation of
the central roots determines their interior zeros and controls the
sums of their products, yielding a sine-kernel limit.

Throughout this subsection, $T$ is the Toeplitz matrix in
\eqref{eq:general-toeplitz}, with $p,q\ge1$ and positive root
parameters in \eqref{eq:root-parameters}, allowing repetitions.
Use the support $\mathcal I=[\lambda_-,\lambda_+]$ and the central root
$\zeta$ of Theorem~\ref{thm:general-support}. Its polar coordinates are
\begin{equation}\label{eq:spectral-polar-data}
 \zeta(x)=r(x)\mathrm e^{-\mathrm i\theta(x)},\qquad
 r(x)>0,\quad 0<\theta(x)<\pi,\qquad x\in \mathcal I^\circ.
\end{equation}
The multisets $\mathcal S_0,\mathcal Z_0$ and
$\widehat{\mathcal S},\widehat{\mathcal Z}$ are those of
\eqref{eq:general-density-root-sets}--\eqref{eq:general-root-boundary-sets}.
They depend on $x$. The first two contain the $q-1$ roots of modulus
less than $r$ and the $p-1$ roots of modulus greater than $r$,
respectively. The enlarged multisets also contain $\zeta,\bar\zeta$.

\subsubsection{Linear forms determined by the matrix}

The density $\dd\Psi(x)=\rho(x)f(x)g(x)^{\mathsf T}\dd x$
contains a column $f(x)$ of length $q$ and a column $g(x)$ of
length $p$. The row $B_n(x)$ can therefore be multiplied by
$f(x)$, and the column $A_n(x)$ by $g(x)^{\mathsf T}$ on the left.
Each product is a scalar weighted sum of polynomial components.
For each fixed $x$, it is linear in those components; its weights
are functions of $x$. Such products are called weighted linear forms.
They are generally algebraic functions of $x$, even though the
components of $A_n$ and $B_n$ are polynomials.

\begin{theorem}[Linear forms independent of the factorization]
\label{prop:spectral-functions}
For any bounded PBF of $T$, let $A_n,B_n$ have the initial conditions
\eqref{eq:general-initial-polynomials}, and let $\rho,f,g$ be given by
\eqref{eq:general-density-data}. The weighted linear forms
\begin{equation}\label{eq:spectral-functions}
 \mathcal A_n(x)\coloneq g(x)^{\mathsf T}A_n(x),\qquad
 \mathcal B_n(x)\coloneq B_n(x)f(x),\qquad n\in\NN,
\end{equation}
satisfy
\begin{equation}\label{eq:spectral-functions-roots}
 \mathcal A_n(x)=\hc_n(\widehat{\mathcal S}),\qquad
 \mathcal B_n(x)=\hc_n(\widehat{\mathcal Z}^{-1}),
 \qquad x\in \mathcal I^\circ,
\end{equation}
where $\hc_n$ is defined in \eqref{eq:complete-homogeneous-definition}.
They are independent of the chosen PBF and obey
\begin{equation}\label{eq:spectral-biorthogonality}
 \int_{\mathcal I}\mathcal B_n(x)\mathcal A_\ell(x)\rho(x)\dd x
 =\delta_{n,\ell},\qquad n,\ell\in\NN.
\end{equation}
\end{theorem}
\begin{proof}
A change of PBF leaves the linear forms unchanged by cancellation
of the normalization matrices. Let $(H',K')$ be the initial matrices for a second
bounded PBF of the same $T$. Formula~\eqref{eq:general-density-data}
gives $f'=H'H^{-1}f$ and
${g'}^{\mathsf T}=g^{\mathsf T}K^{-1}K'$, while
Corollary~\ref{cor:change-normalization} gives
$B'_n=B_nHH'^{-1}$ and $A'_n=K'^{-1}KA_n$. Hence
\[
 \begin{aligned}
 {g'}^{\mathsf T}A'_n
   &=g^{\mathsf T}K^{-1}K'K'^{-1}KA_n=g^{\mathsf T}A_n,\\
 B'_nf'&=B_nHH'^{-1}H'H^{-1}f=B_nf.
 \end{aligned}
\]
The weight $\rho$ also remains unchanged, since its root formula
in \eqref{eq:general-density-data} depends only on $T$.

To compute the linear forms, fix $x\in \mathcal I^\circ$.
Multiplying $AT=xA$ by $g(x)^{\mathsf T}$ gives the scalar left
recurrence for $\mathcal A_n$.
Its first $p$ entries are
$[\hc_0(\widehat{\mathcal S})\ \cdots\
\hc_{p-1}(\widehat{\mathcal S})]KK^{-1}$, by
\eqref{eq:general-initial-polynomials} and
\eqref{eq:general-density-data}.
The sequence on the right of the first formula in
\eqref{eq:spectral-functions-roots} has these entries and satisfies
the same recurrence. Indeed, putting
$Q^\vee(t;x)=t^{p+q}Q(t^{-1};x)$ gives
$Q^\vee(t;x)\prod_{u\in\widehat{\mathcal S}}(1-ut)^{-1}=t_p\prod_{v\in\mathcal Z_0}(1-vt)$.
The right-hand side has degree $p-1$. Comparing coefficients of
$t^{n+p}$ therefore gives
$\sum_{j=-q}^{p}t_j\hc_{n+j}(\widehat{\mathcal S})
=x\hc_n(\widehat{\mathcal S})$, with negative-index terms zero.
Since the coefficient of the next term is $t_p>0$, the first $p$
terms determine the whole sequence.

The first $q$ entries for the right recurrence form the vector
$H^{-1}f=[\hc_i(\widehat{\mathcal Z}^{-1})]_{i\in\{0,\dots,q-1\}}$.
The identity
$Q(t;x)\prod_{v\in\widehat{\mathcal Z}}(1-t/v)^{-1}=\omega\prod_{u\in\mathcal S_0}(1-t/u)$
has a right-hand side of degree $q-1$. Its coefficient of $t^{n+q}$
gives the right recurrence. The coefficient $\omega=t_{-q}>0$
determines each next term. This proves both formulas.
Finally, the rank-one density gives the identity
$B_n\dd\Psi A_\ell=(B_nf)(g^{\mathsf T}A_\ell)\rho\dd x
=\mathcal B_n\mathcal A_\ell\rho\dd x$.
Thus the biorthogonality integral with the matrix measure in
Corollary~\ref{cor:general-favard} is exactly the scalar integral
\eqref{eq:spectral-biorthogonality}.
\end{proof}

The root formulas for these linear forms now give their Chebyshev
convolutions. The coefficients below depend on $x$ through the
remaining roots.

\begin{proposition}[Exact Chebyshev convolutions]
\label{prop:spectral-chebyshev-convolutions}
For $n\in\NN$ and $x\in \mathcal I^\circ$, with $r,\theta$ from
\eqref{eq:spectral-polar-data},
\begin{align}
 \mathcal A_n(x)
 &=\sum_{j=0}^{n}\hc_j(\mathcal S_0)\,
       r^{n-j}U_{n-j}(\cos\theta),\label{eq:spectral-chebyshev-left}\\
 \mathcal B_n(x)
 &=\sum_{j=0}^{n}\hc_j(\mathcal Z_0^{-1})\,
       r^{-(n-j)}U_{n-j}(\cos\theta).\label{eq:spectral-chebyshev-right}
\end{align}
Here $U_k(\cos\theta)=\sin((k+1)\theta)/\sin\theta$.
For an empty multiset, $\hc_0=1$ and $\hc_j=0$ for $j>0$.
\end{proposition}
\begin{proof}
The Chebyshev generating function \cite[Section~18.12]{NISTDLMF}
gives
\[
 \frac{1}{(1-\zeta t)(1-\bar\zeta t)}
 =\frac{1}{1-2r\cos\theta\,t+r^2t^2}
 =\sum_{n\ge0}r^nU_n(\cos\theta)t^n.
\]
Multiply by $\prod_{u\in\mathcal S_0}(1-ut)^{-1}$ and compare
coefficients to obtain \eqref{eq:spectral-chebyshev-left}.
The reciprocal pair $\zeta^{-1},\bar\zeta^{-1}$ has modulus $r^{-1}$
and the same normalized real part $\cos\theta$. Multiplication by
$\prod_{v\in\mathcal Z_0}(1-t/v)^{-1}$ proves
\eqref{eq:spectral-chebyshev-right}.
\end{proof}

For $p=q=1$ both additional root groups are empty:
$\mathcal A_n=r^nU_n(\cos\theta)$ and
$\mathcal B_n=r^{-n}U_n(\cos\theta)$, recovering
Example~\ref{ex:scalar-chebyshev}.

\subsubsection{Finite kernels and a Poisson kernel}

Define the $p\times q$ matrix kernel and its scalar contraction by
\begin{equation}\label{eq:spectral-kernel-definition}
 \begin{aligned}
 \mathscr K_N(x,y)&\coloneq\sum_{n=0}^{N-1}A_n(x)B_n(y),&
 k_N(x,y)&\coloneq g(x)^{\mathsf T}\mathscr K_N(x,y)f(y)
             =\sum_{n=0}^{N-1}\mathcal A_n(x)\mathcal B_n(y).
 \end{aligned}
\end{equation}
This scalar kernel acts as the identity on the span of
$\mathcal A_0,\ldots,\mathcal A_{N-1}$ when integration uses the
functions $\mathcal B_n$ and the weight $\rho$. More explicitly,
expanding the two finite sums gives
\[
 \int_{\mathcal I}k_N(x,t)k_N(t,y)\rho(t)\dd t
 =\sum_{n,l=0}^{N-1}\mathcal A_n(x)\mathcal B_l(y)
   \int_{\mathcal I}\mathcal B_n(t)\mathcal A_l(t)\rho(t)\dd t.
\]
By \eqref{eq:spectral-biorthogonality}, the inner integral is
$\delta_{n,l}$. Thus the sum reduces to $k_N(x,y)$, proving the
reproduction identity
\begin{equation}\label{eq:spectral-kernel-reproduction}
 \int_{\mathcal I} k_N(x,t)k_N(t,y)\rho(t)\dd t=k_N(x,y).
\end{equation}
Christoffel--Darboux formulas for mixed orthogonality are developed in
\cite{DaemsKuijlaars2007,AriznabarretaManas2015}. The next formula
specializes the recurrence-matrix calculation to the Toeplitz band
and expresses the scalar kernels in terms of roots.

\begin{proposition}[Boundary formula and root sums]
\label{prop:spectral-kernel-root-sums}
For $N\ge\max\{p,q\}$,
\begin{align}
 (x-y)\mathscr K_N(x,y)
 =\sum_{d=1}^{p}t_d\sum_{j=N-d}^{N-1}A_{j+d}(x)B_j(y)
 -\sum_{d=1}^{q}t_{-d}\sum_{j=N-d}^{N-1}A_j(x)B_{j+d}(y).
 \label{eq:spectral-CD}
\end{align}
For $x,y\in \mathcal I^\circ$, put $X\coloneq\widehat{\mathcal S}(x)$ and
$Y\coloneq\widehat{\mathcal Z}(y)^{-1}$.
When each multiset has distinct elements, define
$c_u(X)\coloneq\prod_{w\in X\setminus\{u\}}(1-w/u)^{-1}$.
Then
\begin{equation}\label{eq:spectral-finite-root}
 k_N(x,y)=\sum_{u\in X}\sum_{v\in Y}c_u(X)c_v(Y)G_N(uv).
\end{equation}
The Poisson kernel is
\begin{equation}\label{eq:spectral-Poisson}
 \mathscr P_t(x,y)\coloneq
 \sum_{n\ge0}t^n\mathcal A_n(x)\mathcal B_n(y)
 =\sum_{u\in X}\sum_{v\in Y}\frac{c_u(X)c_v(Y)}{1-tuv},
 \qquad |t|\frac{r(x)}{r(y)}<1,
\end{equation}
where $G_N(w)\coloneq\sum_{n=0}^{N-1}w^n$, so $G_N(1)=N$ and
$G_N(w)=(1-w^N)/(1-w)$ for $w\ne1$.
Coincident elements are treated by continuity or by partial fractions
with repeated poles.
\end{proposition}
\begin{proof}
Let $\mathsf P_N$ be the diagonal projection onto the first $N$ coordinates. The recurrences give
$(x-y)\mathscr K_N=A(x)(T\mathsf P_N-\mathsf P_NT)B(y)$.
An entry of the commutator vanishes if both indices are less than
$N$ or both are at least $N$. The lower entries crossing this cut
have positions $(j+d,j)$ with $N-d\le j<N$ and coefficient $t_d$.
The upper ones have positions $(j,j+d)$ and coefficient $-t_{-d}$.
These are exactly the two sums in \eqref{eq:spectral-CD}.

Partial fractions in \eqref{eq:complete-homogeneous-definition} give
$\hc_n(X)=\sum_{u\in X}c_u(X)u^n$.
Multiplying the corresponding expressions for $X,Y$ and summing the
finite geometric progression proves \eqref{eq:spectral-finite-root}.
In the domain of \eqref{eq:spectral-Poisson},
$|tuv|\le |t|r(x)/r(y)<1$ for every pair, so the infinite geometric
series converge absolutely. This gives the Poisson formula.
For each fixed $n$, the coefficient $\hc_n$ is a polynomial in
the root variables. Hence the finite sum on the left of
\eqref{eq:spectral-finite-root} has a finite, well-defined limit
when roots coincide, even though separate partial-fraction
coefficients may diverge. For the infinite sum, fix a point in
the strict domain $|t|r(x)/r(y)<1$. Choose numbers $R_X,R_Y$
slightly larger than the largest root moduli in $X,Y$ so that
$|t|R_XR_Y<1$. Cauchy's coefficient estimate, applied on the
circles of radii $R_X^{-1},R_Y^{-1}$ in the generating variables,
gives $|\hc_n(X)\hc_n(Y)|\le C(R_XR_Y)^n$, uniformly for
nearby root multisets. This summable bound permits passage to the
limit in the Poisson series. It proves the asserted extension at
coincident roots without taking limits of its individual singular
partial-fraction coefficients.
\end{proof}

On the diagonal, \eqref{eq:spectral-Poisson} converges for every
$|t|<1$. In the scalar case $r$ is constant, and the formula becomes
\[
 \sum_{n\ge0}t^nU_n(\cos\theta)U_n(\cos\varphi)
 =\frac{1-t^2}
 {(1-2t\cos(\theta-\varphi)+t^2)
  (1-2t\cos(\theta+\varphi)+t^2)}.
\]

\subsubsection{Interior oscillations and zeros}

The scaled roots $\mathcal S_0/r$ and $r\mathcal Z_0^{-1}$ have
modulus less than one. On each compact interval inside $\mathcal I^\circ$
they stay uniformly inside the unit disk. This gives geometric
estimates for the coefficients in the two Chebyshev convolutions.
Define the nonzero amplitudes
\begin{equation}\label{eq:spectral-amplitudes}
 C_A(x)\coloneq\prod_{u\in\mathcal S_0}(1-u/\zeta)^{-1},\qquad
 C_B(x)\coloneq\prod_{v\in\mathcal Z_0}(1-\zeta/v)^{-1}.
\end{equation}
On a compact interval choose continuous arguments
$\delta_A\coloneq\arg C_A$ and $\delta_B\coloneq\arg C_B$.

\Needspace{8\baselineskip}
\begin{theorem}[Uniform oscillatory formulas and interior zeros]
\label{thm:spectral-oscillations}
Let $J$ be a compact interval in $\mathcal I^\circ$.
\begin{enumerate}
\item There exists $0<\vartheta<1$ such that, uniformly for $x\in J$,
\begin{align}
 r^{-n}\mathcal A_n
 &=\frac{|C_A|}{\sin\theta}
   \sin((n+1)\theta-\delta_A)+\mathrm O(\vartheta^n),
   &&n\to\infty,\label{eq:spectral-left-asymptotic}\\
 r^n\mathcal B_n
 &=\frac{|C_B|}{\sin\theta}
   \sin((n+1)\theta+\delta_B)+\mathrm O(\vartheta^n),
   &&n\to\infty.\label{eq:spectral-right-asymptotic}
\end{align}
The derivatives of the remainders are $\mathrm O(n\vartheta^n)$ as
$n\to\infty$. With $\mathcal F$ from \eqref{eq:F-critical}, the
phase and its associated density satisfy
\begin{equation}\label{eq:spectral-phase-density}
 \varpi(x)\coloneq-\frac{\theta'(x)}{\pi}
 =-\frac{\im\mathcal F(\zeta(x))}
        {\pi x|\mathcal F(\zeta(x))-q|^2}>0.
\end{equation}
\item Let $J_0$ be a compact interval in $J^\circ$. Each solution
in $J_0$ of $(n+1)\theta-\delta_A=k\pi$, with $k\in\mathbb Z$,
has, for all sufficiently large $n$,
exactly one zero of $\mathcal A_n$ within distance
$\mathrm O(\vartheta^n/n)$ as $n\to\infty$. The corresponding statement for
$\mathcal B_n$ uses $(n+1)\theta+\delta_B=k\pi$.
These zeros are simple.

\item For $[a,b]\subset \mathcal I^\circ$, let $N_n^A([a,b])$ and
$N_n^B([a,b])$ count the real zeros of these scalar functions.
Then
\begin{equation}\label{eq:spectral-zero-count}
 \frac{N_n^A([a,b])}{n},\ \frac{N_n^B([a,b])}{n}
 \longrightarrow \int_a^b\varpi(x)\dd x
 =\frac{\theta(a)-\theta(b)}{\pi},\qquad n\to\infty.
\end{equation}
For consecutive zeros converging to $x_0\in \mathcal I^\circ$, the spacing
is $(n\varpi(x_0))^{-1}(1+\mathrm o(1))$ as $n\to\infty$.
\end{enumerate}
\end{theorem}
\begin{proof}
Divide \eqref{eq:spectral-chebyshev-left} by $r^n$ and set
$b_j\coloneq\hc_j(\mathcal S_0/r)$. Their generating function has its
poles uniformly outside a circle of radius greater than one for
$x\in J$. Cauchy's coefficient estimate therefore gives
$|b_j|\le C\vartheta_0^j$ with $\vartheta_0<1$.
In the trigonometric expression for the convolution, extend the
sum from $j\le n$ to all $j\ge0$. The added terms have total size
$\mathrm O(\vartheta_0^n)$, uniformly on $J$ as $n\to\infty$. Thus
\[
 r^{-n}\mathcal A_n
 =\frac{1}{\sin\theta}\im\left(
   \mathrm e^{\mathrm i(n+1)\theta}
       \sum_{j\ge0}b_j\mathrm e^{-\mathrm i j\theta}\right)
   +\mathrm O(\vartheta_0^n),\qquad n\to\infty.
\]
The sum is $\prod_{u\in\mathcal S_0}(1-u/\bar\zeta)^{-1}
=\overline{C_A}$, since $\mathcal S_0$ is invariant under
conjugation. This proves \eqref{eq:spectral-left-asymptotic}.
For the right convolution use $b_j\coloneq\hc_j(r\mathcal Z_0^{-1})$.
The corresponding sum is $C_B$, proving
\eqref{eq:spectral-right-asymptotic}.

The analyticity needed here concerns the whole root group,
not an arbitrary labeling of its individual roots. Locally near
each point of $J$, surround one group by a contour containing no
other root. Its power sums are contour integrals of
$z^kQ_z(z;x)/Q(z;x)$ and are analytic in $x$. Newton's identities
then show that the coefficients of its monic polynomial are
analytic, even if roots inside the contour coincide. These local
polynomials agree on overlaps because they describe the same
separated group. The positive function $r$ also has a local
analytic extension, obtained from the square root of the product
of the two central roots. Thus the generating functions for
$b_j$ are analytic in $x$ on fixed smaller neighborhoods.
Cauchy's estimate first in the generating variable and then in
$x$ bounds $b_j$ and its derivative geometrically, with uniform
constants after taking a finite cover of $J$. Differentiating the
trigonometric terms contributes at most a factor proportional to
$n+j$. Summing their tails gives the derivative bound
$\mathrm O(n\vartheta^n)$ as $n\to\infty$ for a choice
$\vartheta_0<\vartheta<1$.

Differentiation of $a(\zeta(x))=x$, using \eqref{eq:a-critical},
gives
\begin{equation}\label{eq:spectral-root-log-derivative}
 \frac{\zeta'}{\zeta}=\frac{1}{x(\mathcal F(\zeta)-q)}.
\end{equation}
The imaginary part of the left-hand side is $-\theta'$.
Moreover,
$\im\mathcal F(\zeta)
=\im\zeta\sum_\nu\beta_\nu/|1+\beta_\nu\zeta|^2<0$.
Taking imaginary parts proves \eqref{eq:spectral-phase-density}.

To locate the zeros, divide \eqref{eq:spectral-left-asymptotic}
by the positive amplitude $|C_A|/\sin\theta$.
The resulting expression is the sine of a phase plus a remainder
of size $\mathrm O(\vartheta^n)$, with derivative
$\mathrm O(n\vartheta^n)$ as $n\to\infty$.
The phase derivative $(n+1)\theta'-\delta_A'$ is negative
and has modulus at least $cn$ on $J$, for all sufficiently large $n$.
Write this phase as $\phi_n$ and the normalized remainder as
$\varepsilon_n$. Around a solution $x_{k,n}$ of
$\phi_n(x_{k,n})=k\pi$, consider the interval on which
$|\phi_n-k\pi|\le\pi/4$. For solutions in a fixed compact
subinterval of $J^\circ$, this whole interval lies in $J$ when
$n$ is large, because $|\phi_n'|\ge cn$. At its endpoints the
sine has opposite signs and absolute value $1/\sqrt2$, whereas
$\varepsilon_n$ tends uniformly to zero. A zero therefore exists.
On that interval, $|\cos\phi_n|\ge1/\sqrt2$ and its sign is
constant. The bound $|\varepsilon_n'|\le Cn\vartheta^n$ implies
that $(\sin\phi_n+\varepsilon_n)'$ has constant nonzero sign
for sufficiently large $n$. The zero is consequently unique and
simple. At the zero, $|\sin\phi_n|\le C\vartheta^n$, so
$|\phi_n-k\pi|\le C'\vartheta^n$. The mean value theorem and
$|\phi_n'|\ge cn$ give the distance
$\mathrm O(\vartheta^n/n)$ from $x_{k,n}$ as $n\to\infty$.
Outside these phase intervals, $|\sin\phi_n|\ge1/\sqrt2$,
so there are no further zeros in the compact subinterval.

The number of multiples of $\pi$ between the phase values at $a,b$
is $(n+1)(\theta(a)-\theta(b))/\pi+\mathrm O(1)$ as $n\to\infty$.
The zeros near the two endpoints alter this count by a bounded
number. This proves \eqref{eq:spectral-zero-count}.
Two consecutive zeros have phase difference
$\pi+\mathrm O(\vartheta^n)$ as $n\to\infty$; the mean value theorem
gives their spacing. The right-hand sequence follows by the same
argument with its specified phase.
\end{proof}

The density $\varpi$ in \eqref{eq:spectral-phase-density} is
Hirschman's limiting eigenvalue density \cite{Hirschman1967}.
In its expression in terms of the boundary values of the roots
\cite[equation~(1.12)]{DuitsKuijlaars2008}, the roots in
$\mathcal S_0$ have the same boundary values on the two sides of
$\mathcal I^\circ$, while $\zeta,\bar\zeta$ exchange positions.
The remaining contribution is
$(2\pi\mathrm i)^{-1}
(\zeta'/\zeta-\bar\zeta'/\bar\zeta)=\varpi$.
Thus the same density governs the interior zeros of
$\mathcal A_n,\mathcal B_n$ and the limiting eigenvalues of the
finite Toeplitz sections.

\subsubsection{Diagonal limits and the sine kernel}

The next identity connects the amplitudes
\eqref{eq:spectral-amplitudes}, the weight $\rho$ in
\eqref{eq:general-density-data}, and the derivative of the central
root. It supplies the coefficients needed to sum the oscillations.

\begin{proposition}[Amplitude identity and diagonal limits]
\label{prop:spectral-amplitude-identity}
For $x\in \mathcal I^\circ$,
\begin{equation}\label{eq:spectral-amplitude-identity}
 \frac{\rho C_AC_B}{2\sin^2\theta}
 =-\frac{\mathrm i}{\pi}\frac{\zeta'}{\zeta}
 =\varpi-\frac{\mathrm i}{\pi}(\log r)'.
\end{equation}
Uniformly on each compact subset of $\mathcal I^\circ$,
\begin{align}
 \frac{\rho(x)}{N}k_N(x,x)
 &=\varpi(x)+\mathrm O(N^{-1}),&&N\to\infty,
 \label{eq:spectral-diagonal-limit}\\
 (1-t)\rho(x)\mathscr P_t(x,x)
 &=\varpi(x)+\mathrm O(1-t),&&t\uparrow1.
 \label{eq:spectral-Abel-limit}
\end{align}
\end{proposition}
\begin{proof}
Write $s_0\coloneq\prod_{u\in\mathcal S_0}|u|$.
As in the proof of Theorem~\ref{thm:general-measure},
$\prod_{u\in\mathcal S_0}u=(-1)^{q-1}s_0$.
The product of all roots of $Q$ is $(-1)^{p+q}\omega/t_p$.
To keep track of the signs, write the amplitudes as
\[
 C_A=\frac{\zeta^{q-1}}{\prod_{u\in\mathcal S_0}(\zeta-u)},
 \qquad
 C_B=\frac{(-1)^{p-1}\prod_{v\in\mathcal Z_0}v}
                {\prod_{v\in\mathcal Z_0}(\zeta-v)}.
\]
Since $\zeta\bar\zeta=r^2$, the product of all roots also gives
\[
 \prod_{v\in\mathcal Z_0}v
   =\frac{(-1)^{p+1}\omega}{t_pr^2s_0},\qquad
 Q_z(\zeta;x)=t_p(\zeta-\bar\zeta)
       \prod_{u\in\mathcal S_0}(\zeta-u)
       \prod_{v\in\mathcal Z_0}(\zeta-v).
\]
Multiplying the two amplitude formulas and substituting these
identities cancels the signs and gives
\[
 C_AC_B
 =\frac{\omega(\zeta-\bar\zeta)\zeta^{q-1}}
        {r^2s_0Q_z(\zeta;x)}.
\]
Now $\zeta-\bar\zeta=-2\mathrm i r\sin\theta$,
$\rho=r\sin\theta\,s_0/(\pi\omega)$, and
$\zeta'=\zeta^q/Q_z(\zeta;x)$.
Substitution gives the first equality in
\eqref{eq:spectral-amplitude-identity}. The second follows from
$\zeta'/\zeta=(\log r)'-\mathrm i\theta'$ and
$\varpi=-\theta'/\pi$.

Put $a_A\coloneq|C_A|/\sin\theta$, $a_B\coloneq|C_B|/\sin\theta$, and
$\phi\coloneq\delta_A+\delta_B$. Multiplying the two formulas in
Theorem~\ref{thm:spectral-oscillations} gives
\[
 \mathcal A_n\mathcal B_n
 =\frac{a_Aa_B}{2}
   \bigl(\cos\phi-\cos(2(n+1)\theta+\delta_B-\delta_A)\bigr)
   +\mathrm O(\vartheta^n),\qquad n\to\infty.
\]
The real part of \eqref{eq:spectral-amplitude-identity} says
$\rho a_Aa_B\cos\phi/2=\varpi$.
On a compact interior interval, $\sin\theta$ is bounded away from
zero. Hence the finite sum of the oscillating cosines is uniformly
bounded, as is the sum of the geometric remainders.
Division by $N$ proves \eqref{eq:spectral-diagonal-limit}.
The Abel sum of the oscillation is the real part of a bounded
multiple of $(1-t\mathrm e^{2\mathrm i\theta})^{-1}$.
Multiplication by $1-t$ proves \eqref{eq:spectral-Abel-limit}.
\end{proof}

For type-II multiple orthogonal polynomials,
\cite{SwiderskiVanAssche2022} relates the weak limit of the
normalized diagonal kernel measure to the limiting zero measure.
Here \eqref{eq:spectral-diagonal-limit} gives an explicit interior
limit with a uniform error bound for the scalar functions of the mixed
Toeplitz system. Summing the explicit oscillatory asymptotics and
using both parts of the amplitude identity
\eqref{eq:spectral-amplitude-identity} gives the following
sine-kernel limit, including its normalization and uniform error.

Bleher and Kuijlaars \cite{BleherKuijlaars2004} obtained the bulk
sine kernel for the Gaussian external-source ensemble through multiple
Hermite polynomials. For mixed Gaussian weights, Daems, Kuijlaars,
and Veys \cite[Theorem~2.2]{DaemsKuijlaarsVeys2008} obtained the
same limit for a conjugated correlation kernel of nonintersecting
Brownian paths with two starting and two ending points.
For the mixed Toeplitz system below, the modulus $r$ of the central
symbol root gives the conjugating factor explicitly. The amplitude
identity gives the local normalization and the uniform error for
arbitrary $p,q$, including repeated root parameters.

\begin{theorem}[Interior sine-kernel limit]
\label{thm:spectral-sine-kernel}
Let $J$ be a compact interval in $\mathcal I^\circ$.
For $x_0\in J$, put $\varpi_0\coloneq\varpi(x_0)>0$,
$x_N\coloneq x_0+u/(N\varpi_0)$ and
$y_N\coloneq x_0+v/(N\varpi_0)$ for real $u,v$. Then
\begin{equation}\label{eq:spectral-sine-kernel}
 \frac{\rho(x_0)}{N\varpi_0}
 \left(\frac{r(y_N)}{r(x_N)}\right)^{\!N}k_N(x_N,y_N)
 =
 \frac{\sin\pi(u-v)}{\pi(u-v)}
 +\mathrm O(N^{-1}),\qquad N\to\infty.
\end{equation}
The sine quotient equals $1$ at $u=v$. For every $M>0$, the error
is uniform for $x_0\in J$ and $|u|,|v|\le M$.
This holds for all $p,q\ge1$ and all positive root parameters
in \eqref{eq:root-parameters}, including repetitions; the error
constant may depend on $J,M$ and these fixed symbol data.
\end{theorem}
\begin{proof}
Choose a compact interval $J_*\subset\mathcal I^\circ$ with
$J\subset J_*^\circ$. On $J_*$, both $\varpi$ and $\sin\theta$
have positive lower bounds, and the derivatives of the amplitudes,
phases and $\log r$ used below are bounded. Consequently, for
$x_0\in J$ and $|u|,|v|\le M$, the points $x_N,y_N$ lie in $J_*$
for all sufficiently large $N$, uniformly in these parameters.
The oscillatory formulas and all the following estimates are
uniform on these sets.

Write $\Delta\coloneq u-v$, $h_0\coloneq(\log r)'(x_0)$ and
$c_0\coloneq h_0/\varpi_0$.
Taylor expansion gives, uniformly for $0\le n<N$,
\begin{align*}
 \left(\frac{r(x_N)}{r(y_N)}\right)^n
 &=\mathrm e^{c_0\Delta n/N}(1+\mathrm O(N^{-1})),\\
 (n+1)(\theta(x_N)-\theta(y_N))
 &=-\pi\Delta n/N+\mathrm O(N^{-1}),
\end{align*}
where all limits in this proof are taken as $N\to\infty$.
Multiply \eqref{eq:spectral-left-asymptotic} at $x_N$ by
\eqref{eq:spectral-right-asymptotic} at $y_N$, and use
$2\sin a\sin b=\cos(a-b)-\cos(a+b)$.
The term with the sum of the phases has geometric ratio approaching
$\mathrm e^{2\mathrm i\theta(x_0)}\ne1$. Including the real factor
$(r(x_N)/r(y_N))^n$, its finite geometric sum stays uniformly bounded:
the ratio stays away from one and its $N$th power is bounded.
Indeed, $|1-\mathrm e^{2\mathrm i\theta(x_0)}|=2\sin\theta(x_0)$
has a positive lower bound for $x_0\in J$.
After division by $N$ its contribution is $\mathrm O(N^{-1})$.

The amplitudes and the additional phases differ from their values
at $x_0$ by $\mathrm O(N^{-1})$. The geometric errors in
Theorem~\ref{thm:spectral-oscillations} have a bounded sum.
The remaining, difference-phase term divided by $N$ is therefore
a Riemann sum with error $\mathrm O(N^{-1})$.
Set $a_A=|C_A|/\sin\theta$, $a_B=|C_B|/\sin\theta$ and
$\phi=\delta_A+\delta_B$ as in the preceding proof.
The real and imaginary parts of
\eqref{eq:spectral-amplitude-identity}, evaluated at $x_0$, are
\[
 \frac{\rho(x_0)a_A(x_0)a_B(x_0)}{2\varpi_0}\cos\phi(x_0)=1,
 \qquad
 \frac{\rho(x_0)a_A(x_0)a_B(x_0)}{2\varpi_0}\sin\phi(x_0)
 =-\frac{c_0}{\pi}.
\]
They give
\[
 \frac{\rho(x_0)}{N\varpi_0}k_N(x_N,y_N)
 =\int_0^1\mathrm e^{c_0\Delta t}
   \left(\cos(\pi\Delta t)+\frac{c_0}{\pi}\sin(\pi\Delta t)\right)
       \dd t+\mathrm O(N^{-1}).
\]
For $\Delta\ne0$, the integrand is the derivative of
$\mathrm e^{c_0\Delta t}\sin(\pi\Delta t)/(\pi\Delta)$.
The integral is consequently
$\mathrm e^{c_0\Delta}\sin(\pi\Delta)/(\pi\Delta)$, with value
$1$ at $\Delta=0$ by continuity.
Finally,
$(r(y_N)/r(x_N))^N
=\mathrm e^{-c_0\Delta}(1+\mathrm O(N^{-1}))$.
This proves \eqref{eq:spectral-sine-kernel}, including uniformity
across the diagonal.
\end{proof}

The geometric-sum bound involves
$(\min_{x\in J}\sin\theta(x))^{-1}$, which grows without bound
as $J$ approaches an endpoint of the support. The estimates above
therefore give uniformity on compact interior intervals.

The factor $(r(y)/r(x))^N$ is a ratio of one function evaluated
at $y$ and at $x$, so it preserves the reproduction identity
\eqref{eq:spectral-kernel-reproduction}. It equals one when $r$
is constant. The PBF independence of the functions
\eqref{eq:spectral-functions} and kernels
\eqref{eq:spectral-kernel-definition} makes these results applicable
to both the binomial family and the explicit
geometric-root PBF below.

\begin{example}[Four distinct roots and an explicit sine limit]
\label{ex:spectral-quartic-kernel}
Let $p=q=2$ and
$P(z)\coloneq(1+3z+z^2)(1+4z+z^2)
=1+7z+14z^2+7z^3+z^4$.
The Toeplitz matrix has diagonal $14$, first off-diagonals $7$,
and second off-diagonals $1$. Its support is $[2,30]$.
Indeed, with $w\coloneq z+z^{-1}$ the equation $Q(z;x)=0$ becomes
$w^2+7w+12-x=0$. Set
\[
 w_\pm(x)\coloneq\frac{-7\pm\sqrt{1+4x}}2,\quad
 \cos\theta(x)=\frac{w_+(x)}2,\quad
 s(x)\coloneq\frac{w_-(x)+\sqrt{w_-(x)^2-4}}2.
\]
For $2<x<30$, $-2<w_+<2$, $w_-<-2$ and $-1<s<0$.
The central roots are $\mathrm e^{\pm\mathrm i\theta}$ and the
remaining roots are $s,s^{-1}$. Thus $r=1$ and
\begin{equation}\label{eq:spectral-quartic-functions}
 \mathcal A_n(x)=\mathcal B_n(x)
 =\sum_{j=0}^{n}s(x)^jU_{n-j}(\cos\theta(x)).
\end{equation}
The weight and limiting density are
\[
 \rho(x)=\frac{|s(x)|\sin\theta(x)}{\pi},\qquad
 \varpi(x)=\frac{1}{2\pi\sin\theta(x)\sqrt{1+4x}}.
\]
Here the two scalar sequences coincide, so
\eqref{eq:spectral-biorthogonality} is an orthogonality relation
for the algebraic functions \eqref{eq:spectral-quartic-functions}.
At $x_0=12$, $s=(-7+3\sqrt5)/2$, $\theta=\pi/2$,
$\rho=(7-3\sqrt5)/(2\pi)$ and $\varpi_0=1/(14\pi)$.
In particular,
$\mathcal A_n(12)=\mathcal B_n(12)
=\sum_{j=0}^{\lfloor n/2\rfloor}(-1)^js^{n-2j}$, and
\[
 \frac{49-21\sqrt5}{N}
 k_N\left(12+\frac{14\pi u}{N},12+\frac{14\pi v}{N}\right)
 \longrightarrow\frac{\sin\pi(u-v)}{\pi(u-v)},\qquad N\to\infty.
\]
\end{example}

\subsection{Determinantal polynomials and vector equilibrium}
\label{subsec:determinants-equilibrium}

The determinants of consecutive left and right polynomial vectors
recover the characteristic polynomial of a finite section. Replacing
one vector gives the determinantal polynomials of
\cite[Sections~2--4]{BFM2023}, which describe its left and right
eigenvectors. For Toeplitz matrices these determinants are Schur
polynomials. This identifies three of them with the shifted Toeplitz
determinants of Duits and Kuijlaars and supplies their limiting zero
measures. The BTP condition constrains the geometry of their vector
equilibrium problem, and the finite eigenvectors recover the scalar
functions of Theorem~\ref{prop:spectral-functions}.

Retain the assumptions and notation of
Section~\ref{subsec:spectral-kernels}. Let $Z(x)$ be the multiset
of all $m=p+q$ roots of $Q(z;x)=P(z)-xz^q$.
Write its elements as $z_1(x),\ldots,z_m(x)$, counting
multiplicities. Every root is nonzero.
To distinguish the characteristic polynomial from the symbol
$P$, put $\chi_0\coloneq1$ and
$\chi_N(x)\coloneq\det(xI_N-T^{[N-1]})$ for $N\ge1$.
Define the nonzero constants
\begin{equation}\label{eq:determinantal-scales}
 c_N^-\coloneq(-1)^{(p-1)N}t_p^N,\qquad
 c_N^+\coloneq(-1)^{(q-1)N}t_{-q}^N,\qquad N\in\NN.
\end{equation}
For $0\le n\le N$, define
\begin{equation}\label{eq:determinantal-QR}
 \begin{aligned}
 \mathsf Q_{n,N}(x)&\coloneq\det
 [A_n(x)\ A_{N+1}(x)\ \cdots\ A_{N+p-1}(x)],\\
 \mathsf R_{n,N}(x)&\coloneq\det
 [B_n(x)^{\mathsf T}\ B_{N+1}(x)^{\mathsf T}\ \cdots\
 B_{N+q-1}(x)^{\mathsf T}].
 \end{aligned}
\end{equation}
For $p=1$ or $q=1$, the corresponding determinant contains just
its first vector. These are the polynomials denoted by $Q_{n,N}$
and $R_{n,N}$ in \cite[Definition~3.1]{BFM2023}.
The Schur notation used below is fixed by
\begin{equation}\label{eq:determinantal-Schur-definition}
 s_\lambda(Z)=\det[\hc_{\lambda_i-i+j}(Z)]_{i,j\in\{1,\dots,d\}},
 \qquad \hc_0=1,\quad \hc_j=0\quad(j<0),
\end{equation}
where $d$ is any length obtained by appending zeros to the partition
$\lambda$. The notation $(N^r,n)$ means $r$ parts equal to $N$
followed by $n$; zero parts may be omitted.

\subsubsection{Schur formulas and finite eigenvectors}

\begin{proposition}[Determinants, cofactors, and their kernel]
\label{prop:determinantal-Schur-cofactors}
For $N\in\NN$ and $0\le n\le N$,
\begin{align}
 \mathsf Q_{n,N}(x)&=s_{(N^{p-1},n)}(Z(x)),&
 \mathsf R_{n,N}(x)&=s_{(N^{q-1},n)}(Z(x)^{-1}),
 \label{eq:determinantal-Schur-QR}\\
 \chi_N(x)&=c_N^-s_{(N^p)}(Z(x))
          =c_N^+s_{(N^q)}(Z(x)^{-1}).
 \label{eq:determinantal-Schur-chi}
\end{align}
In particular, these determinants are independent of the chosen PBF.
If $M_N(x)\coloneq xI_{N+1}-T^{[N]}$, then
\begin{equation}\label{eq:determinantal-adjugate}
 c_N^-\mathsf Q_{n,N}(x)=[\operatorname{adj}M_N(x)]_{N,n},\qquad
 c_N^+\mathsf R_{n,N}(x)=[\operatorname{adj}M_N(x)]_{n,N}.
\end{equation}
The polynomial kernel
\begin{equation}\label{eq:determinantal-kernel}
 \mathcal D_N(x,y)\coloneq c_N^-c_N^+
       \sum_{n=0}^{N}\mathsf Q_{n,N}(x)\mathsf R_{n,N}(y)
 =\frac{\chi_{N+1}(x)\chi_N(y)-\chi_N(x)\chi_{N+1}(y)}{x-y}
\end{equation}
is symmetric. On the diagonal it equals
$\chi_{N+1}'\chi_N-\chi_N'\chi_{N+1}$ and is strictly positive
for real $x$.
\end{proposition}
\begin{proof}
The sequences $[\hc_{n-a}(Z)]_{a\in\{0,\dots,p-1\}}$ satisfy the left
recurrence. To see the boundary condition as well as the recurrence,
their generating functions are
$t^a\prod_{z\in Z}(1-zt)^{-1}$, and
$Q^\vee(t;x)=t_p\prod_{z\in Z}(1-zt)$.
Multiplying gives $t_pt^a$, of degree at most $p-1$.
Its coefficients of degrees $n+p$, $n\ge0$, vanish, which is
exactly the left recurrence with negative-index entries zero.
The first $p$ columns form an upper unit triangular matrix.
They are independent of $x$, since the first coefficient of
$Q^\vee$ involving $x$ has degree $p$.
The initial matrix $K^{-1}$ in
\eqref{eq:general-initial-polynomials} is also upper unit triangular.
Thus the actual vectors $A_n$ are obtained by a constant change of
basis of determinant one. The determinants in
\eqref{eq:determinantal-QR} are unchanged.

For the right vectors use $[\hc_{n-b}(Z^{-1})]_{b\in\{0,\dots,q-1\}}$ and
$Q(t;x)=t_{-q}\prod_{z\in Z}(1-t/z)$.
The same argument applies to the lower unit triangular initial
matrix $H^{-1}$. For the left determinant, before this constant
change of basis, entry $(a,j)$ is $\hc_{n_j-a}(Z)$, where
$a\in\{0,\ldots,p-1\}$ and the column indices are
$n_1\coloneq n$, $n_j\coloneq N+j-1$ for $2\le j\le p$.
Reverse both rows and columns and then transpose. The two
reversals have the same sign, so they leave the determinant
unchanged. In one-based indices the resulting entry is
$\hc_{N-i+j}(Z)$ for $i<p$ and $\hc_{n-p+j}(Z)$ for $i=p$.
These are exactly the entries of
\eqref{eq:determinantal-Schur-definition} for the partition
$(N^{p-1},n)$. Using $q$ and $Z^{-1}$ proves the right formula
in the same way, including the cases of one-component vectors. The characteristic-polynomial
identity of \cite[Theorem~2.7]{BFM2023}, with the extreme diagonals
constant, now gives \eqref{eq:determinantal-Schur-chi}.
This is consistent with the general Toeplitz cofactor formulas in
\cite{MaximenkoMoctezuma2017}.

For completeness, the cofactor identification follows directly
from the finite boundary equations. Put
$\boldsymbol q_N\coloneq[\mathsf Q_{0,N}\ \cdots\ \mathsf Q_{N,N}]$
and $\boldsymbol r_N\coloneq[\mathsf R_{0,N}\ \cdots\
\mathsf R_{N,N}]^{\mathsf T}$.
Extend the determinant definitions to the next indices. They vanish
at $N+1,\ldots,N+p-1$ on the left and at
$N+1,\ldots,N+q-1$ on the right, because a vector is repeated.
The recurrences consequently give
$c_N^-\boldsymbol q_NM_N=\chi_{N+1}e_N^{\mathsf T}$ and $c_N^+M_N\boldsymbol r_N=\chi_{N+1}e_N$.
For example, in the last component of the left equation the
recurrence leaves only $t_p\mathsf Q_{N+p,N}$. Moving its
first column to the last position gives
\[
 \mathsf Q_{N+p,N}
 =(-1)^{p-1}\det[A_{N+1}\ \cdots\ A_{N+p}]
 =(-1)^{p-1}\frac{\chi_{N+1}}{c_{N+1}^-}.
\]
Since $c_{N+1}^-=(-1)^{p-1}t_pc_N^-$, multiplying by $c_N^-$
makes this component exactly $\chi_{N+1}$. Every earlier
component vanishes by the recurrence and the repeated-column
identities. The right equation uses the identical calculation
with $q,t_{-q},c_N^+$ and the vectors $B_j^{\mathsf T}$. For $\chi_{N+1}\ne0$, multiplication by $M_N^{-1}$ proves
\eqref{eq:determinantal-adjugate}; polynomial continuation proves
it for every $x$.

The resolvent identity
$M_N(x)^{-1}M_N(y)^{-1}
=(M_N(y)^{-1}-M_N(x)^{-1})/(x-y)$, in entry $(N,N)$, proves
\eqref{eq:determinantal-kernel}. This is the identity of
\cite[Proposition~4.1]{BFM2023} in the present notation.
Its diagonal expression follows by taking $y\to x$.
The characteristic polynomials of consecutive oscillatory sections
have simple interlacing real zeros
\cite[Proposition~7.5]{BFM2023}. Let $\lambda_1<\cdots<\lambda_{N+1}$ be the zeros of
$\chi_{N+1}$. Partial fractions give
\[
 \frac{\chi_N(x)}{\chi_{N+1}(x)}
  =\sum_{j=1}^{N+1}\frac{w_j}{x-\lambda_j},\qquad
 w_j\coloneq\frac{\chi_N(\lambda_j)}{\chi_{N+1}'(\lambda_j)}>0.
\]
The positivity follows directly from interlacing: both numerator
and denominator in $w_j$ have sign $(-1)^{N+1-j}$, since their
polynomials are monic. Differentiating the finite sum shows that
\[
 \frac{\mathcal D_N(x,x)}{\chi_{N+1}(x)^2}
 =-\left(\frac{\chi_N}{\chi_{N+1}}\right)'(x)
 =\sum_{j=1}^{N+1}\frac{w_j}{(x-\lambda_j)^2}>0
\]
away from these zeros.
At a zero of $\chi_{N+1}$ the value is
$\chi_{N+1}'\chi_N>0$, again by interlacing.
\end{proof}

At a zero of $\chi_{N+1}$ the boundary equations in this proof
give left and right eigenvectors of $T^{[N]}$; their last entries
are nonzero since $\chi_N$ and $\chi_{N+1}$ have no common zero.
For example, if $p=2$, then
$\mathsf Q_{n,N}=\hc_N(Z)\hc_n(Z)-\hc_{N+1}(Z)\hc_{n-1}(Z)$
and $\chi_N=(-t_2)^N(\hc_N(Z)^2-\hc_{N+1}(Z)\hc_{N-1}(Z))$.
For $q=2$ the same formulas for $\mathsf R_{n,N}$ use $Z^{-1}$.

\subsubsection{The equilibrium problem on alternating semiaxes}

For $N\ge1$ and $k\in\{-q+1,\ldots,p-1\}$, order the roots in $Z(x)$ by
increasing modulus and set
\begin{equation}\label{eq:DK-curves-determinants}
 \begin{aligned}
 \Gamma_k&\coloneq \{x\in\CC:|z_{q+k}(x)|=|z_{q+k+1}(x)|\},\\
\mathfrak D_{k,N}(x)
&\coloneq \det T_N\big(\bigl(a(z)-x\bigr)z^{-k}\big)
=\det\bigl[t_{i-j+k}-x\,\delta_{j,i+k}\bigr]_{i,j\in\{0,\dots,N-1\}}.
 \end{aligned}
\end{equation}
For a Laurent polynomial $b$, $T_N(b)\coloneq[b_{i-j}]_{i,j\in\{0,\dots,N-1\}}$.
These are the curves and shifted determinants of
\cite[equations~(2.7) and~(2.15)]{DuitsKuijlaars2008}.
Their construction applies to general complex banded Toeplitz
symbols. The positive root parameters in \eqref{eq:root-parameters} impose the
following geometry on all of their curves.

\begin{proposition}[Real equilibrium supports with alternating signs]
\label{prop:DK-alternating-supports}
For the BTP symbol \eqref{eq:general-toeplitz},
\begin{equation}\label{eq:DK-alternating-supports}
 \Gamma_k\subseteq\{x\in\RR:(-1)^kx\ge0\},
 \qquad k\in\{-q+1,\ldots,p-1\}.
\end{equation}
Moreover, $\Gamma_0=[\lambda_-,\lambda_+]$, with the endpoints
in Theorem~\ref{thm:general-support}.
\end{proposition}
\begin{proof}
For every $r>0$ the function
\begin{equation}\label{eq:DK-circle-modulus}
 |P(r\mathrm e^{\mathrm i\theta})|^2
 =\omega^2\prod_{\nu=1}^{m}
       (1+\beta_\nu^2r^2+2\beta_\nu r\cos\theta)
\end{equation}
is strictly increasing in $\cos\theta\in[-1,1]$.
Each factor is positive on $(-1,1]$ and strictly increasing;
at $-1$ it is nonnegative. Thus two distinct roots of
$P(z)-xz^q$ on the same circle must be conjugate, and there are
exactly two roots on that circle. Since $a$ has real coefficients,
the common value $a(z)=a(\bar z)=x$ is real.

A repeated root at $x\ne0$ satisfies $\mathcal F(z)=q$, with
$\mathcal F$ defined in \eqref{eq:F-critical}. For nonreal $z$,
\[
 \im\mathcal F(z)=\im z\sum_{\nu=1}^{m}
          \frac{\beta_\nu}{|1+\beta_\nu z|^2}\ne0.
\]
Thus a repeated root is real, and its image $x=a(z)$ is real.
At $x=0$ all roots are negative real numbers. This proves
$\Gamma_k\subseteq\RR$, including the multiple-root cases.

Suppose first that the equal-modulus roots for $x\in\Gamma_k$
are distinct. They form a nonreal conjugate pair of modulus $r$.
There are exactly $q+k-1$ roots of smaller modulus.
Equation~\eqref{eq:DK-circle-modulus} gives
$|P(-r)|<|x|r^q<P(r)$.
Hence $Q(r;x)>0$ and
$\sgn Q(-r;x)=(-1)^{q+1}\sgn x$.
The number of real roots in $(-r,r)$ has the parity of $q+k-1$,
because the other roots inside the circle occur in conjugate pairs.
Counting multiplicities, the change in sign from $-r$ to $r$
therefore gives
$\sgn Q(-r;x)=(-1)^{q+k-1}\sgn Q(r;x)$.
Combining the two signs yields $\sgn x=(-1)^k$.

The remaining collisions can be checked directly. Let $c<0$ be a
repeated root at $x_*\coloneq a(c)\ne0$, and let $j$ count the roots of
$P$ in $(c,0)$, with multiplicities. On $|z|=-c$, the minimum
of $|P(z)|$ is attained only at $c$.
For $y=(1-\varepsilon)x_*$, $0<\varepsilon\ll1$, Rouch\'e's
theorem counts exactly $j$ roots of $P(z)-yz^q$ inside this circle.
Also $\mathcal F'(c)>0$ and
$a''(c)=x_*\mathcal F'(c)/c$, so the colliding roots are real
for these $y$, one on each side of $c$. Their positions at the
collision are consequently $j,j+1$ in increasing modulus.
Thus $j=q+k$, while
$\sgn a(c)=(-1)^{j-q}=(-1)^k$.
At a positive collision, \eqref{eq:DK-circle-modulus} instead
has its maximum at $c$. For $y=(1+\varepsilon)a(c)$, Rouch\'e
counts $q$ roots inside $|z|=c$. Since $a''(c)>0$, the two nearby
roots for these $y$ are real, one inside and one outside the circle.
Thus the collision is in positions
$q,q+1$ and belongs to $\Gamma_0$ with positive image.
Collisions at $x=0$ satisfy \eqref{eq:DK-alternating-supports}
for every $k$. The identification of $\Gamma_0$ is
Theorem~\ref{thm:general-support}.
\end{proof}

Let $\mu_k$ denote the equilibrium measures of
\cite[Theorem~2.3]{DuitsKuijlaars2008}. Their supports are
$\Gamma_k$ and their masses are $(q+k)/q$ for $k\le0$ and
$(p-k)/p$ for $k\ge0$. With
$I(\nu,\sigma)\coloneq\iint\log|x-y|^{-1}\dd\nu(x)\dd\sigma(y)$
and $I(\nu)\coloneq I(\nu,\nu)$, they uniquely minimize
\begin{equation}\label{eq:DK-vector-energy}
 \mathcal E(\boldsymbol\nu)
 \coloneq\sum_{k=-q+1}^{p-1}I(\nu_k)
  -\sum_{k=-q+1}^{p-2}I(\nu_k,\nu_{k+1})
\end{equation}
over their admissible class with these supports and masses.
Following the physical interpretation in
\cite[equations~(2.13)--(2.14) and the following discussion]{DuitsKuijlaars2008},
each positive measure describes a distribution of charged particles.
Particles in the same distribution repel one another; particles in
consecutive distributions attract one another with half that strength.
Nonconsecutive distributions do not interact directly.
The factor $1/2$ agrees with \eqref{eq:DK-vector-energy}: the
self-interaction counts each pair twice, whereas each mixed interaction
appears only once. The signs of the semiaxes specify particle positions,
not signs of electrical charges.
Proposition~\ref{prop:DK-alternating-supports} places consecutive
components on opposite semiaxes. The central component has mass
one and occupies the compact positive interval
$[\lambda_-,\lambda_+]$.
For $U^\nu(x)\coloneq\int\log|x-y|^{-1}\dd\nu(y)$, its equilibrium
condition is $2U^{\mu_0}-U^{\mu_{-1}}-U^{\mu_1}=\ell_0$
on that interval, with an absent neighbor interpreted as zero.
Thus the neighboring measures supply the external field seen by
the central one.

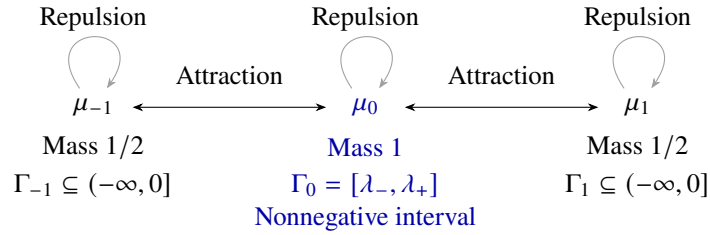
\begin{figure}[htbp]
\centering
\begin{tikzpicture}[x=1cm,y=1cm,font=\small]
 \node at (0,0) {$\mu_{-1}$};
 \node[blue!65!black] at (3.6,0) {$\mu_0$};
 \node at (7.2,0) {$\mu_1$};
 \draw[{Stealth[length=4pt]}-{Stealth[length=4pt]}]
   (.52,0) -- (3.08,0) node[midway,above=5pt] {Attraction};
 \draw[{Stealth[length=4pt]}-{Stealth[length=4pt]}]
   (4.12,0) -- (6.68,0) node[midway,above=5pt] {Attraction};
 \foreach \x in {0,3.6,7.2}{
   \draw[-{Stealth[length=4pt]},gray!75]
     (\x-.18,.25) .. controls (\x-.8,1.15) and (\x+.8,1.15)
     .. (\x+.18,.25);
   \node at (\x,1.17) {Repulsion};
 }
 \node at (0,-.53) {Mass $1/2$};
 \node[blue!65!black] at (3.6,-.53) {Mass $1$};
 \node at (7.2,-.53) {Mass $1/2$};
 \node at (0,-1.03) {$\Gamma_{-1}\subseteq(-\infty,0]$};
 \node[blue!65!black] at (3.6,-1.03)
   {$\Gamma_0=[\lambda_-,\lambda_+]$};
 \node at (7.2,-1.03) {$\Gamma_1\subseteq(-\infty,0]$};
 \node[blue!65!black] at (3.6,-1.48) {Nonnegative interval};
\end{tikzpicture}
\caption{The equilibrium problem \eqref{eq:DK-vector-energy}
for $p=q=2$. The diagram orders the measures by index, not by
their positions on the real line. Each distribution has repulsive
logarithmic interaction with itself and attractive interaction of
half that strength with its neighbors. The two distributions
$\mu_{-1}$ and $\mu_1$ do not interact directly.
The BTP condition places the two neighboring supports on the
nonpositive semiaxis and the central support on a nonnegative
interval, by Proposition~\ref{prop:DK-alternating-supports}.
The neighboring measures generate the external field
$-U^{\mu_{-1}}-U^{\mu_1}$ seen by $\mu_0$.}
\label{fig:equilibrium-chain}
\end{figure}

Figure~\ref{fig:equilibrium-chain} displays this interaction and
the support restrictions when $p=q=2$.

\begin{corollary}[Shifted determinants and limiting zero measures]
\label{cor:DK-determinantal-identification}
For $N\ge1$,
\begin{equation}\label{eq:DK-determinantal-identification}
 \begin{aligned}
 \mathfrak D_{0,N}&=(-1)^N\chi_N,&
 \mathfrak D_{1,N}&=c_N^-\mathsf Q_{0,N} &&(p\ge2),&
 \mathfrak D_{-1,N}&=c_N^+\mathsf R_{0,N} &&(q\ge2).
 \end{aligned}
\end{equation}
Giving each zero mass $1/N$, with multiplicity, the zero measures
of $\chi_N$, $\mathsf Q_{0,N}$, and $\mathsf R_{0,N}$ converge
weakly as $N\to\infty$ to $\mu_0$, $\mu_1$, and $\mu_{-1}$,
respectively, whenever the corresponding index is present.
\end{corollary}
\begin{proof}
The first equality is the definition of $\chi_N$.
In \eqref{eq:determinantal-adjugate}, the cofactor for
$\mathsf Q_{0,N}$ deletes row $0$ and column $N$ of $M_N$.
The remaining rows are $1,\ldots,N$ and columns $0,\ldots,N-1$.
Its determinant is $(-1)^N\mathfrak D_{1,N}$, and the cofactor
sign is also $(-1)^N$. The signs cancel. Deleting row $N$ and
column $0$ gives the reciprocal shift and proves the formula for
$\mathsf R_{0,N}$.
All multiplying constants are nonzero and independent of $x$.
The weak limits therefore follow from
\cite[Theorem~2.6]{DuitsKuijlaars2008}: all coefficients in the
filled band are positive, so the nonzero exponents of $a$ have
greatest common divisor one, as required there. In particular, the last
two limiting measures have masses $(p-1)/p$ and $(q-1)/q$;
the normalization used here is $1/N$, as in that theorem.
\end{proof}

\begin{example}[The cubic binomial equilibrium]
\label{ex:DK-cubic-equilibrium}
For $a(z)=4(1+z)^3/(27z)$, $p=2$ and $q=1$.
Duits and Kuijlaars compute
$\Gamma_0=[0,1]$ and $\Gamma_1=(-\infty,0]$
in \cite[Section~6, Example~1]{DuitsKuijlaars2008}.
The equilibrium masses are $1$ and $1/2$.
Corollary~\ref{cor:DK-determinantal-identification} identifies
their two sequences of determinants with $\chi_N$ and a constant
multiple of $\mathsf Q_{0,N}$. The sign pattern is the $p=2,q=1$
instance of Proposition~\ref{prop:DK-alternating-supports}.
\end{example}

The limiting zero measure of the characteristic polynomials is unchanged
when the diagonals converge to the Toeplitz coefficients.
The proof applies the classical method of moments: diagonal
convergence makes normalized traces asymptotically equal, while
BTP gives real positive eigenvalues and boundedness places them
in a common compact interval.

\begin{proposition}[Characteristic zeros for convergent Toeplitz diagonals]
\label{prop:asymptotic-toeplitz-zero-measure}
Let $T=[t_{i-j}]_{i,j\in\NN}$ be the BTP Toeplitz matrix above,
of bandwidth $(p,q)$, with limiting characteristic zero measure
$\mu_0$. Let $\widetilde T$ be a bounded BTP matrix of the same
bandwidth such that
$\widetilde T_{n,n-k}\longrightarrow t_k$ as $n\to\infty$ for $-q\le k\le p$,
where entries with negative indices are omitted. Define
\[
 \widetilde\chi_N(x)\coloneq\det(xI_N-\widetilde T^{[N-1]}),\qquad
 \widetilde\nu_N\coloneq\frac1N
       \sum_{\widetilde\chi_N(\lambda)=0}\delta_\lambda,
 \qquad N\ge1,
\]
counting zeros with multiplicity. Then
$\widetilde\nu_N\Longrightarrow\mu_0$ weakly as $N\to\infty$.
Equivalently, for every continuous function $f$ on a compact real
interval containing all these zeros,
\[
 \lim_{N\to\infty}\frac1N
       \sum_{\widetilde\chi_N(\lambda)=0}f(\lambda)
       =\int_{\mathcal I}f(x)\varpi(x)\dd x.
\]
In particular, this conclusion holds after any finite change of
the entries that preserves boundedness and BTP.
\end{proposition}
\begin{proof}
Put $A_N\coloneq\widetilde T^{[N-1]}$, $B_N\coloneq T^{[N-1]}$, and
$E_N\coloneq A_N-B_N$. Let $\|M\|_{\mathrm F}^2\coloneq\sum_{i,j}|M_{i,j}|^2$
be the Frobenius norm. There are only $p+q+1$ diagonals. On each one the squared
differences form a bounded sequence tending to zero. To see why
its average tends to zero, fix $\varepsilon>0$ and choose an
index after which every squared difference is at most
$\varepsilon$. In the average of the first $N$ terms, the
finitely many preceding terms contribute a quantity tending to
zero as $N\to\infty$, and the remaining terms contribute at
most $\varepsilon$. Summing over the finitely many diagonals gives
$N^{-1}\|E_N\|_{\mathrm F}^2\longrightarrow0$ as $N\to\infty$.
Choose $M\ge\max\{1,\|\widetilde T\|_2,\|T\|_2\}$.
The finite sections have operator norms at most $M$. For each
fixed integer $k\ge1$, the telescoping identity and the
Frobenius Cauchy--Schwarz inequality yield
\[
 \begin{aligned}
 A_N^k-B_N^k
   &=\sum_{j=0}^{k-1}A_N^jE_NB_N^{k-1-j},\\
 \frac1N\left|\operatorname{tr}(A_N^k-B_N^k)\right|
   &\le kM^{k-1}\frac{\|E_N\|_{\mathrm F}}{\sqrt N}
      \longrightarrow0,\qquad N\to\infty.
 \end{aligned}
\]
For the trace estimate, cyclically move $E_N$ to the first
position and use
$\|B_N^{k-1-j}A_N^j\|_{\mathrm F}\le\sqrt N M^{k-1}$.

By Corollary~\ref{cor:btp-oscillatory}, both finite sections are
oscillatory. Their eigenvalues are positive and real, and the
operator-norm bound places them all in $[0,M]$. The traces are
the sums of the corresponding powers of the eigenvalues. Hence
the empirical characteristic zero measures for $A_N$ and $B_N$
have asymptotically equal moments of every nonnegative integer
order. Corollary~\ref{cor:DK-determinantal-identification} gives
the weak limit $\mu_0$ for $B_N$. For a continuous $f$ on $[0,M]$, choose a polynomial $h$
with $\sup_{[0,M]}|f-h|<\varepsilon$. All empirical measures
have mass one, as does $\mu_0$. Their integrals of $f-h$ thus
have absolute value at most $\varepsilon$. The convergence of
the moments gives convergence of their integrals of $h$.
Consequently the limiting discrepancy for $f$ is at most
$2\varepsilon$. Letting $\varepsilon\downarrow0$ proves the
same weak limit for $A_N$. The density of $\mu_0$ is $\varpi$,
as identified after Theorem~\ref{thm:spectral-oscillations}.
\end{proof}

\subsubsection{Limits of finite eigenvectors and their kernel}

The shifted determinants describe global zero distributions. Their
ratios also recover the scalar functions that occur in the Favard
integral. The following formulas and estimates hold on a compact interval
$J\subset \mathcal I^\circ$, using the root groups and phase already
defined in Section~\ref{subsec:spectral-kernels}.

\begin{lemma}[Uniform remainder for separated eigenvalues]
\label{lem:separated-power-remainder}
Let $M(x)$ be a finite matrix analytic near a compact real interval
$J$. Suppose its distinguished eigenvalues $\mu_1(x),\ldots,\mu_s(x)$
are simple and have the same positive modulus $r(x)$ on $J$.
Assume every other eigenvalue has modulus at most $\vartheta_0r(x)$,
where $\vartheta_0<1$, and that $r$ has a nonvanishing analytic extension
near each point of $J$. Let $\Pi_j$ denote the corresponding spectral
projections. For some $C<\infty$ and $\vartheta\in(\vartheta_0,1)$,
\begin{equation}\label{eq:separated-power-remainder}
 \frac{M(x)^N}{r(x)^N}
 =\sum_{j=1}^{s}\left(\frac{\mu_j(x)}{r(x)}\right)^N\Pi_j(x)
       +E_N(x),\qquad
 \sup_{x\in J}\bigl(\|E_N(x)\|+\|E_N'(x)\|\bigr)\le C\vartheta^N
\end{equation}
for every $N\ge0$. Eigenvalues in the remaining group may coincide.
\end{lemma}
\begin{proof}
Fix $x_0\in J$ and choose $\vartheta\in(\vartheta_0,1)$.
In a sufficiently small complex neighborhood of $x_0$, the circle
$|w|=\vartheta$ encloses all the remaining eigenvalues of $M(x)/r(x)$
and excludes the distinguished ones. On its product with a smaller
closed neighborhood, $wI-M(x)/r(x)$ is invertible and its inverse
is uniformly bounded. The part of the matrix power contributed
by the enclosed eigenvalues is exactly
\begin{equation}\label{eq:separated-power-contour}
 E_N(x)=\frac{1}{2\pi\mathrm i}\int_{|w|=\vartheta}
          w^N\bigl(wI-M(x)/r(x)\bigr)^{-1}\,\mathrm dw.
\end{equation}
This formula follows by applying Cauchy's formula to each Jordan
block, and so includes multiple eigenvalues. It gives
$\|E_N(x)\|\le C_0\vartheta^N$ throughout that complex neighborhood.
Cauchy's estimate in the variable $x$ gives the same exponential
bound for $E_N'$ on a smaller neighborhood. The excluded simple
eigenvalues contribute the displayed terms with $\Pi_j$.
A finite cover of $J$ makes the constants uniform.
\end{proof}

\begin{theorem}[Boundary ratios and characteristic-polynomial oscillations]
\label{thm:determinantal-asymptotics}
Let $J$ be a compact interval in $\mathcal I^\circ$.
For each fixed $n\in\NN$, uniformly for $x\in J$,
\begin{equation}\label{eq:determinantal-boundary-limits}
 \frac{\mathsf Q_{n,N}(x)}{\mathsf Q_{0,N}(x)}
 =\mathcal A_n(x)+\mathrm O(\vartheta^N),\qquad
 \frac{\mathsf R_{n,N}(x)}{\mathsf R_{0,N}(x)}
 =\mathcal B_n(x)+\mathrm O(\vartheta^N),\qquad N\to\infty,
\end{equation}
for some $0<\vartheta<1$; the denominators are nonzero for all
sufficiently large $N$.
Put $V(x)\coloneq\mathcal Z_0(x)\cup\{\zeta(x)\}$ and
\begin{equation}\label{eq:determinantal-leading-data}
 R_*(x)\coloneq t_pr(x)\prod_{v\in\mathcal Z_0(x)}|v|,\qquad
 C_*(x)\coloneq\prod_{v\in V(x)}\prod_{u\in Z(x)\setminus V(x)}
                   (1-u/v)^{-1}.
\end{equation}
Both $R_*>0$ and $C_*\ne0$. On $J$ choose a continuous
argument $\delta_* \coloneq\arg C_*$. Then
\begin{equation}\label{eq:determinantal-chi-asymptotic}
 \frac{\chi_N(x)}{R_*(x)^N}
 =2|C_*(x)|\cos(N\theta(x)-\delta_*(x))
       +\mathrm O(\vartheta^N),\qquad N\to\infty.
\end{equation}
The derivative of the remainder in
\eqref{eq:determinantal-chi-asymptotic} is
$\mathrm O(N\vartheta^N)$ as $N\to\infty$, uniformly on $J$.
\end{theorem}
\begin{proof}
First suppose that the roots are distinct. Expansion of the alternant formula
for a Schur polynomial along its first $d$ rows gives
\begin{equation}\label{eq:determinantal-Schur-expansion}
 s_{(N^d,n)}(Z)
 =\sum_{\substack{V\subset Z\\|V|=d}}
   \left(\prod_{v\in V}v\right)^N
   \prod_{\substack{v\in V\\u\in Z\setminus V}}(1-u/v)^{-1}
   \hc_n(Z\setminus V),\qquad 0\le n\le N.
\end{equation}
Indeed, the first $d$ exponents in the alternant are
$N+m-1,\ldots,N+m-d$. Their determinant supplies the power
$N+m-d$ of each selected root and its Vandermonde determinant.
The remaining determinant, divided by its Vandermonde determinant,
is $\hc_n$ on the complementary roots. Dividing by the full
Vandermonde leaves exactly the cross factors displayed above.
For $d=0$ the formula is simply $s_{(n)}=\hc_n$.
For $n=0$, this is the Schur form of Widom's formula given in
\cite[Theorem~6]{Alexandersson2012}; the last part $n$ supplies the
additional complete homogeneous polynomial in the present formula.

Take $d=p-1$. The unique selection with largest product modulus
is $\mathcal Z_0$, separated from its complement by
Theorem~\ref{thm:general-support}. Its complement is $\widehat{\mathcal S}$. Denote this
largest product by $\mu\coloneq\prod_{v\in\mathcal Z_0}v$ and its
coefficient in the expansion by $c$. The cross factors defining
$c$ are nonzero because the two groups have different moduli.
For fixed $n$, the finite expansion therefore reads
\[
 \mathsf Q_{n,N}
 =\mu^N\bigl(c\hc_n(\widehat{\mathcal S})
                       +\mathrm O(\vartheta^N)\bigr),\qquad
 \mathsf Q_{0,N}
 =\mu^N\bigl(c+\mathrm O(\vartheta^N)\bigr),
 \quad N\to\infty.
\]
On a compact set without root coincidences, $|c|$ has a positive
minimum, so the second expression is nonzero for all sufficiently
large $N$. Dividing gives $\hc_n(\widehat{\mathcal S})=
\mathcal A_n$ with a geometric error. The treatment below makes
these estimates uniform also at coincidences within a group.
For $\mathsf R_{n,N}$ apply the same expansion to $Z^{-1}$
with $d=q-1$. The unique largest group is
$\mathcal S_0^{-1}$, with complement
$\widehat{\mathcal Z}^{-1}$. This proves both limits, with a
geometric error because the product-modulus gaps are uniform on $J$.

For $s_{(N^p)}$, use \eqref{eq:determinantal-Schur-expansion}
with $d=p$ and $n=0$. Exactly two selections have largest
product modulus: $\mathcal Z_0\cup\{\zeta\}$ and its conjugate.
Every real root of $Q(z;x)$ is negative for $x\in \mathcal I^\circ$,
as proved in Theorem~\ref{thm:general-measure}. Hence
$\prod_{v\in\mathcal Z_0}v=(-1)^{p-1}
\prod_{v\in\mathcal Z_0}|v|$.
Multiplication by $c_N^-$ in
\eqref{eq:determinantal-Schur-chi} gives the leading terms
$R_*^NC_*\mathrm e^{-\mathrm iN\theta}$ and its conjugate.
Their sum is \eqref{eq:determinantal-chi-asymptotic}.
This is the two-term specialization of Widom's expansion used in
\cite[equations~(2.22)--(2.25)]{DuitsKuijlaars2008}.

Uniformity when roots inside a separated group coincide follows
from matrix powers, whose entries remain analytic through such a
coincidence. Let $\mathsf e_j(Z)$ denote the
elementary symmetric polynomial of degree $j$ in all the roots,
and let $C(x)$ be the $m\times m$ companion matrix with first row
$[\mathsf e_1(Z),-\mathsf e_2(Z),\ldots,(-1)^{m+1}\mathsf e_m(Z)]$, ones on the first
subdiagonal and zeros elsewhere. These symmetric coefficients are
polynomial in $x$. Put $v\coloneq[1\ 0\ \cdots\ 0]^{\mathsf T}$.
The generating function for $\hc_N(Z)$ gives
\[
 C(x)^Nv=[\hc_N(Z)\ \hc_{N-1}(Z)\ \cdots\
                 \hc_{N-m+1}(Z)]^{\mathsf T},\qquad N\ge0,
\]
where $\hc_j=0$ for $j<0$.

For $d\ge1$, the first $d$ rows of the Jacobi--Trudi determinant
for $s_{(N^d,n)}$ have entries $\hc_{N-i+j}(Z)$.
They are the first $d$ coordinates of the vectors
$C^NC^{j-1}v$, $1\le j\le d+1$.
The last row has entries $\hc_{n-d-1+j}(Z)$ and is independent
of $N$. Expanding along that row expresses the Schur polynomial
as a linear functional, analytic in $x$, of the matrix
$(\bigwedge^d C)^N$. Here $\bigwedge^d C$ is the matrix of all
$d\times d$ minors of $C$; the identity
$\bigwedge^d(C^N)=(\bigwedge^d C)^N$ is Cauchy--Binet.
The same representation with $d$ rows applies to $s_{(N^d)}$.
For $d=0$, the boundary formula is simply $\hc_n$.

The eigenvalues of $\bigwedge^d C$ are the products of $d$ roots,
with multiplicities; this follows by triangularizing $C$.
For the boundary ratios, the unique leading product contains
the entire large-root group. For the characteristic polynomial,
the two leading products contain that entire group and one of
the two distinct central roots. Thus these leading eigenvalues
are simple even when roots within the other groups coincide.
Their coefficients are the nonzero cross products already computed
from \eqref{eq:determinantal-Schur-expansion}: those products contain
no differences within a group. Continuity gives the same coefficients
at a coincidence.

Theorem~\ref{thm:general-support} gives a uniform modulus gap on $J$.
Apply Lemma~\ref{lem:separated-power-remainder} to these compound
matrices, and to the companion matrix for $Z^{-1}$ on the right.
For a boundary ratio one can first divide the compound matrix by
its unique leading eigenvalue and take $r=1$. The resulting leading
coefficient is nonzero, so the denominator stays nonzero for large
$N$ and division preserves the exponential estimate.
For the characteristic polynomial, include the scalar
$(-1)^{p-1}t_p$ in the compound matrix. Its two leading eigenvalues
are $R_*\mathrm e^{\mp\mathrm i\theta}$.
The function $R_*$ has a local analytic extension: it is the
positive square root of their product on the real interval.
The lemma gives uniform exponential bounds for the normalized
remainder and its derivative. In particular it gives the stated
$\mathrm O(N\vartheta^N)$ derivative bound. This proves all assertions
also at coincidences inside the separated root groups.
\end{proof}

B\"ottcher, Grudsky, and Ram\'irez de Arellano \cite{BottcherGrudskyRamirez2006}
proved convergence of normalized finite-section eigenvectors under generic
hypotheses. Bogoya, Gasca, and Grudsky \cite{BogoyaGascaGrudsky2026}
derived eigenvector asymptotics for non-Hermitian tetradiagonal Toeplitz
matrices in the single-arc case. For the BTP symbols of arbitrary
bandwidth considered here, \eqref{eq:determinantal-boundary-limits}
identifies the left and right cofactor ratios, at each fixed index $n$,
with approximations to the Favard forms $\mathcal A_n$ and $\mathcal B_n$.
The geometric error bound holds throughout compact interior intervals,
including points between consecutive eigenvalues.

Bogoya and Grudsky \cite{BogoyaGrudsky2025} obtained individual
eigenvalue expansions for non-Hermitian banded Toeplitz matrices
near simple nondegenerate points of the limiting set. For the BTP
symbols here, the two central roots determine the phase on the
whole open support interval. The Schur expansion gives the phase
shift $\delta_*$ explicitly and controls the remainder and its
derivative on compact interior intervals, also when roots coincide
within a separated group.

The characteristic-polynomial formula locates its interior zeros
near $N\theta-\delta_*=(j+1/2)\pi$ with error
$\mathrm O(\vartheta^N/N)$ as $N\to\infty$.
For $\phi_N\coloneq N\theta-\delta_*$, equation
\eqref{eq:spectral-phase-density} gives
$\phi_N'(x)=-N\pi\varpi(x)-\delta_*'(x)$.
The function $\varpi$ has a positive minimum on $J$, and
$\delta_*'$ is bounded there. Thus $\phi_N'$ is negative, with
magnitude bounded above and below by positive multiples of $N$,
uniformly on $J$ for all sufficiently large $N$.
The sign-change and derivative argument of
Theorem~\ref{thm:spectral-oscillations} proves uniqueness of each
nearby zero and spacing $(N\varpi(x))^{-1}(1+\mathrm o(1))$.
Here the limit is $N\to\infty$ within the fixed compact interval $J$.
The resulting zero density is the measure $\mu_0$ already
identified by Hirschman and Duits--Kuijlaars.

\begin{corollary}[Local limit of the determinantal kernel]
\label{cor:determinantal-sine-kernel}
Let $x_0\in \mathcal I^\circ$, $\varpi_0\coloneq\varpi(x_0)$, and put
$x_N\coloneq x_0+u/(N\varpi_0)$, $y_N\coloneq x_0+v/(N\varpi_0)$.
For the kernel in \eqref{eq:determinantal-kernel},
\begin{equation}\label{eq:determinantal-sine-kernel}
 \frac{\mathcal D_N(x_N,y_N)}
 {\sqrt{\mathcal D_N(x_N,x_N)\mathcal D_N(y_N,y_N)}}
 \longrightarrow\frac{\sin\pi(u-v)}{\pi(u-v)},\qquad N\to\infty,
\end{equation}
uniformly for $u,v$ in compact subsets of $\RR$, with value one
when $u=v$.
\end{corollary}
\begin{proof}
Choose a compact interval $J\subset\mathcal I^\circ$ with $x_0$ in
its interior. Write $R\coloneq R_*$, $b\coloneq2|C_*|$, and
$\phi_N\coloneq N\theta-\delta_*$. These amplitudes and phases
are smooth on $J$, and $R,b$ are positive there. By
Theorem~\ref{thm:determinantal-asymptotics},
\[
 \frac{\chi_N}{R^Nb}=\cos\phi_N+\varepsilon_N,
 \qquad
 \frac{\chi_{N+1}}{R^{N+1}b}
       =\cos(\phi_N+\theta)+\eta_N,
\]
where, uniformly on $J$ as $N\to\infty$,
$|\varepsilon_N|+|\eta_N|=\mathrm O(\vartheta^N)$ and
$|\varepsilon_N'|+|\eta_N'|=\mathrm O(N\vartheta^N)$.
Also $\phi_N'=\mathrm O(N)$ uniformly on $J$ as $N\to\infty$.

Put $a\coloneq R\cos\theta$ and $h\coloneq R\sin\theta$.
Substitution in \eqref{eq:determinantal-kernel}, followed by the
angle-addition formula, gives, for $x\ne y$ in $J$,
\[
 \begin{aligned}
 \frac{\mathcal D_N(x,y)}{R(x)^NR(y)^Nb(x)b(y)}
 ={}&h(x)\frac{\sin(\phi_N(y)-\phi_N(x))}{x-y}\\
 &+\frac{a(x)-a(y)}{x-y}\cos\phi_N(x)\cos\phi_N(y)\\
 &-\frac{h(x)-h(y)}{x-y}\cos\phi_N(x)\sin\phi_N(y)
   +\frac{E_N(x,y)}{x-y}.
 \end{aligned}
\]
Here $E_N$ collects the terms containing at least one of
$\varepsilon_N,\eta_N$. It is antisymmetric, so $E_N(y,y)=0$,
and the preceding derivative bounds give
$\partial_xE_N=\mathrm O(N\vartheta^N)$ uniformly on $J\times J$
as $N\to\infty$. Thus
\[
 \frac{E_N(x,y)}{x-y}
 =\int_0^1\partial_1E_N(y+t(x-y),y)\,\mathrm dt
 =\mathrm O(N\vartheta^N),\qquad N\to\infty.
\]
This identity also defines its continuous value at $x=y$.
The other two divided differences are bounded by
$\sup_J|a'|$ and $\sup_J|h'|$. Hence the displayed kernel equals
its first term plus $\mathrm O(1)$, uniformly across the diagonal
as $N\to\infty$.

For $u,v$ in a fixed compact set,
$x_N,y_N=x_0+\mathrm O(N^{-1})$ as $N\to\infty$.
Using $\theta'=-\pi\varpi$ and integrating $\phi_N'$ between these
two points gives the stronger phase estimate
\[
 \phi_N(y_N)-\phi_N(x_N)
 =\pi(u-v)+\mathrm O\bigl(|u-v|/N\bigr),\qquad N\to\infty.
\]
Indeed, throughout this interval
$-\phi_N'=N\pi\varpi_0+\mathrm O(1)$ as $N\to\infty$.
Since $|\sin s-\sin t|\le|s-t|$, division by
$x_N-y_N=(u-v)/(N\varpi_0)$ now yields
\[
 \frac{\sin(\phi_N(y_N)-\phi_N(x_N))}{x_N-y_N}
 =N\pi\varpi_0\frac{\sin\pi(u-v)}{\pi(u-v)}
       +\mathrm O(1),\qquad N\to\infty,
\]
uniformly also at $u=v$ by continuous extension.
The smooth amplitudes vary by $\mathrm O(N^{-1})$ on this scale
as $N\to\infty$;
their multiplication by a term of order $N$ therefore gives
\[
 \frac{\mathcal D_N(x_N,y_N)}{R(x_N)^NR(y_N)^N}
 =R(x_0)b(x_0)^2
 \left(N\pi\varpi_0\sin\theta(x_0)
       \frac{\sin\pi(u-v)}{\pi(u-v)}+\mathrm O(1)\right)
\]
as $N\to\infty$, uniformly for bounded $u,v$.
Taking $y=x$ in the same divided-difference identity gives
\[
 \mathcal D_N(x,x)=R(x)^{2N+1}b(x)^2
       \bigl(N\pi\varpi(x)\sin\theta(x)+\mathrm O(1)\bigr),
 \qquad N\to\infty,
\]
uniformly near $x_0$. Its leading coefficient is positive.
Consequently the square root of the two diagonal values, divided
by $R(x_N)^NR(y_N)^N$, equals
$N R(x_0)b(x_0)^2\pi\varpi_0\sin\theta(x_0)
 (1+\mathrm O(N^{-1}))$ as $N\to\infty$.
Division proves \eqref{eq:determinantal-sine-kernel}, with an error
$\mathrm O(N^{-1})$ as $N\to\infty$, uniform for bounded $u,v$,
including $u=v$.
\end{proof}

The matrix measure $\Psi$ is the Favard orthogonality measure;
the measures $\mu_k$ solve \eqref{eq:DK-vector-energy}.
Their connection at $k=0$ is explicit:
\eqref{eq:spectral-diagonal-limit} reads
$\rho(x)k_N(x,x)/N\to\dd\mu_0/\dd x$ in the interior.
The limit is taken as $N\to\infty$.
The boundary ratios \eqref{eq:determinantal-boundary-limits}
recover the functions used in this kernel from finite-section
eigenvectors.
The determinantal kernel of \cite{BFM2023} has the local limit
\eqref{eq:determinantal-sine-kernel}.

\subsection[Products with a Toeplitz tail]{Products with a Toeplitz tail: stability, polynomials, and measures}
\label{subsec:quasi-toeplitz}
\noeqref{qt:eq:trace,qt:eq:finite-B,qt:eq:finite-A,qt:eq:root-B,qt:eq:root-A,qt:eq:zero-mass}

Proposition~\ref{prop:constant-factor-orders} constructs BTP matrices
that have the same Toeplitz tail and differ only in a finite initial
block. The following criterion determines which of these products
are interior points of the PBF class. Explicit formulas then give
their polynomial families and matrix measures. The stability criterion compares the largest lower limit with
the smallest upper diagonal limit. The polynomial formulas modify the numerators
of the Toeplitz generating functions, and the measure formula uses a
finite matrix inverse. All statements allow arbitrary $p,q\ge1$.

For clarity, write
$\mathcal L(c)\coloneq I+c\Lambda^{\mathsf T}$ and $\mathcal U(c)\coloneq\Lambda+cI$ for $c>0$,
where $\Lambda$ is the upper shift. Fix positive $\omega$, lower
parameters $\xi_1,\ldots,\xi_p$, and upper parameters
$\eta_1,\ldots,\eta_q$, and put $m\coloneq p+q$. The product with all upper factors first
and its polynomial symbol are
\begin{align*}
 T_0&\coloneq\omega\mathcal U(\eta_1)\cdots\mathcal U(\eta_q)
                 \mathcal L(\xi_1)\cdots\mathcal L(\xi_p),\\
 P(z)&\coloneq\omega\prod_{i=1}^q(1+\eta_i z)
                 \prod_{j=1}^p(1+\xi_j z).
\end{align*}
Let $T_w$ be the product, multiplied by $\omega$, in any ordering
$w$ of these $m$ factors. With $E_r$ the restriction to the first
$r$ coordinates, \eqref{eq:constant-order-boundary} reads
\begin{equation}\label{qt:eq:boundary}
 T_w=T_0+E_p^{\mathsf T}C_wE_q,\qquad C_w\in\RR^{p\times q}.
\end{equation}
Each $T_w$ has a PBF with bounded entries bounded away from zero;
moreover $C_w\le0$ entrywise. The three products in
Example~\ref{ex:uniform-christoffel-chain} are the specialization
already studied in \cite[Sections~4.1--4.3]{Uniform}.

\subsubsection{Limiting groups and interior points}
Suppose first that all $m$ parameters are distinct. By Lemma~\ref{lem:exchange-separated-limits}, the limits under
one adjacent positive exchange obey
\begin{equation}\label{qt:eq:compare}
 (\mathcal U(u),\LL(l))\longmapsto
 (\LL(\min(u,l)),\mathcal U(\max(u,l))).
\end{equation}
Fix a sequence of adjacent upper--lower exchanges that moves all
lower factors to the left, and apply this rule at each step.
For the fixed initial ordering $w$, the resulting ordered factors
are independent of the exchange sequence by
Proposition~\ref{prop:constant-factor-orders}; hence so are their
limits below.
Let $\alpha_1,\ldots,\alpha_p$ be their parameters and let
$\delta_1,\ldots,\delta_q$ be the parameters of the upper factors
in their left-to-right order. The actual exchanged lower entries tend
to the $\alpha_i$. The upper diagonals tend to the $\delta_j$ and
their superdiagonals remain one. After moving their positive diagonal
parts to the middle, the unit factor $U_{q-j+1}$ has upper-entry
limit $1/\delta_j$, for $j\in\{1,\ldots,q\}$,
and the middle diagonal tends to $\omega\prod_j\delta_j$, as
$n\to\infty$. Indeed, the diagonal normalizations in
\eqref{eq:toeplitz-stable-normalization} introduce ratios of
shifted diagonal entries with the same positive limit; those
ratios tend to one. The final two groups therefore partition the
original root parameters.

The general interior criterion in
Corollary~\ref{stab:cor:limiting-interior} applies to the limiting
groups produced by the exchanges. The following theorem also
excludes repeated root parameters by perturbing the common tail.

\begin{theorem}[Interior criterion for reordered constant factors]\label{qt:thm:classification}
Let $T_w$ be a product of the positive constant factors in
\eqref{qt:eq:boundary}, with $p,q\ge1$. The following conditions
are equivalent.
\begin{enumerate}
\item All root parameters are distinct, and the comparison procedure
\eqref{qt:eq:compare} ends with
\begin{equation}\label{qt:eq:gap}
 \max_i\alpha_i<\min_j\delta_j.
\end{equation}
\item There is $\varepsilon>0$ such that every bounded real matrix
$\widetilde T$ of bandwidth at most $(p,q)$ satisfying
$\|\widetilde T-T_w\|_2<\varepsilon$ has a PBF whose factors can
be chosen to depend real analytically on $\widetilde T$.
At $\widetilde T=T_w$, they equal the factors produced by the
positive exchanges. Their initial entries may vary with $\widetilde T$.
\item $T_w$ is an interior point of the PBF class in operator norm
among all bounded real matrices of bandwidth at most $(p,q)$.
\end{enumerate}
The factor sequences in (2) are measured in their supremum norms.
\end{theorem}
\begin{proof}
If (1) holds, the exchanged lower factors have distinct positive
limits, as do the unit upper factors. Condition \eqref{qt:eq:gap} says
that the product of the largest lower and largest upper limits is
less than one. Thus \eqref{stab:eq:separated-ldu-limits} holds.
The diagonal normalization also gives
$d_n\longrightarrow\omega\prod_{j=1}^{q}\delta_j>0$ as $n\to\infty$.
All factor entries and pivots are positive in every row. Their
positive finite limits give common positive lower and finite upper
bounds, including the finitely many initial rows.
Corollary~\ref{stab:cor:free-limits-periodic} supplies the pair
conditions on both sides, and
Proposition~\ref{stab:prop:free-full-transfer} first gives analytic
normalized LDU factors and then an analytic choice of their
bidiagonal factors. This is the map in (2). Continuity keeps its factor entries positive in a
smaller neighborhood, proving (2)$\Rightarrow$(3).

Suppose (3) holds. First there can be no repeated root parameter.
If $\beta$ is repeated, write
$P(z)=\omega(1+\beta z)^2R(z)$ and replace it by
$P_\varepsilon(z)\coloneq\omega\bigl[(1+\beta z)^2+\varepsilon z^2\bigr]R(z)$ for $\varepsilon>0$.
The quadratic factor has discriminant $-4\varepsilon$ and nonreal
zeros. Theorem~\ref{thm:toeplitz} implies that the Toeplitz matrix with
this polynomial symbol is not totally nonnegative. It has a finite
negative minor. Translate the row and column lists of this minor
far enough to avoid the finite block in \eqref{qt:eq:boundary}.
The same negative minor remains after adding $E_p^{\mathsf T}C_wE_q$.
These matrices tend to $T_w$ in operator norm as $\varepsilon\to0$,
contradicting (3).

All parameters are consequently distinct and the comparison procedure
applies. If \eqref{qt:eq:gap} fails, choose $i,j$ with
$\alpha_i>\delta_j$. The corresponding unit factor entries satisfy
$a_i(n)b_{q-j+1}(n)\to\alpha_i/\delta_j>1$ as $n\to\infty$.
Lemma~\ref{qt:lem:cross-obstruction} contradicts (3). This proves (1).
\end{proof}

\begin{corollary}[A separated initial partition]
If all parameters are distinct and every original upper parameter
$\eta_i$ exceeds every original lower parameter $\xi_j$, then
every ordering of the constant factors satisfies the three conditions
of Theorem~\ref{qt:thm:classification}.
\end{corollary}
\begin{proof}
Initially every upper parameter is larger than every lower
parameter. Suppose this is still true before an exchange. In the
pair being exchanged, the minimum is therefore the lower
parameter and the maximum is the upper parameter. Rule
\eqref{qt:eq:compare} leaves both parameter groups unchanged.
Induction over the finite number of exchanges preserves their
strict separation. The final lower maximum is consequently less
than the final upper minimum, which is condition (1) of
Theorem~\ref{qt:thm:classification}.
\end{proof}

\begin{remark}[The effect of the initial block]
For a general initial partition, distinct roots alone do not suffice.
For example $(I+b\Lambda^{\mathsf T})(\Lambda+aI)$, with $b>a>0$, has a constant
PBF but fails the gap condition. Reversing the two factors gives the
exact Toeplitz matrix, whose positive exchange places $a$ below and
$b$ above. Thus the two matrices have the same distinct tail roots
but different interior behavior. This distinction concerns the
matrix, not the choice to fix or vary its initial factor entries.
\end{remark}

\subsubsection{Generating functions and formulas in roots}
Let $T=T_0+E_p^{\mathsf T}CE_q$, where $C\in\RR^{p\times q}$.
The boundary block leaves both positive extreme diagonals of $T_0$
unchanged. Write $B_n(x)$ as a row with $q$ components and $A_n(x)$ as a column
with $p$ components, satisfying $TB=xB$ and $AT=xA$.
Let their initial matrices be $B_{[0,q)}=H^{-1}$ and
$A_{[0,p)}=K^{-1}$, with $H,K$ nonsingular. For a PBF, use its
normalized Favard matrices $H,K$. Put
$v_r(z)\coloneq[1,z,\ldots,z^{r-1}]$ and $P^\vee(z)\coloneq z^mP(1/z)$.
The notation $[F(z)]_{<r}$ means componentwise deletion of powers
$z^j$ with $j\ge r$.

\begin{proposition}[Polynomial formulas for a finite boundary block]\label{qt:prop:generating}
Let $T=T_0+E_p^{\mathsf T}CE_q$, with $T_0$ and $P$ as above,
$C\in\RR^{p\times q}$, and let $H,K$ be nonsingular initial
matrices of sizes $q,p$. The corresponding solutions of $TB=xB$
and $AT=xA$ satisfy, as formal power series in $z$,
\begin{align}
 \sum_{n\ge0}B_n(x)z^n
 &=\frac{\bigl([P(z)v_q(z)]_{<q}-z^qv_p(z)C\bigr)H^{-1}}
 {P(z)-xz^q},\label{qt:eq:Bgen}\\*
 \sum_{n\ge0}A_n(x)z^n
 &=\frac{K^{-1}\bigl([P^\vee(z)v_p(z)^{\mathsf T}]_{<p}
             -z^pCv_q(z)^{\mathsf T}\bigr)}
 {P^\vee(z)-xz^p}.\label{qt:eq:Agen}
\end{align}
Both numerators are independent of $x$ and have degree at most
$m-1$. The boundary block changes none of the extreme diagonals,
so the given initial matrices determine both polynomial families
uniquely by their recurrences.
\end{proposition}
\begin{proof}
Multiply the equation in row $n$ of $TB=xB$ by $z^{n+q}$ and sum.
For the Toeplitz part, changing the index from $n$ to $j$ in
$\sum_{j\ge0}t_{n-j}B_j$ gives $P(z)\sum_jB_jz^j$, apart from
the terms corresponding to $n<0$. These omitted terms have degree
less than $q$ and are exactly
$[P(z)\sum_{j=0}^{q-1}B_jz^j]_{<q}$.
The boundary correction contributes
$z^qv_p(z)CB_{[0,q)}=z^qv_p(z)CH^{-1}$.
Moving this contribution to the other side proves \eqref{qt:eq:Bgen}.
The denominators have nonzero constant terms, so the formal divisions
are well defined. Transpose the recurrence and interchange $p,q$;
the symbol polynomial becomes $P^\vee$ and the boundary block becomes
$C^{\mathsf T}$. Transposing the resulting identity gives
\eqref{qt:eq:Agen}. The degree bounds follow directly from the dimensions
of the boundary block.
\end{proof}

Write $N_B(z)$ and $N_A(z)$ for the full row and column numerators
in \eqref{qt:eq:Bgen} and \eqref{qt:eq:Agen}, including the initial
matrices. Expanding first in the polynomial variable $x$ gives
finite coefficient formulas, valid for every $n\in\NN$:
\begin{align}
 B_n(x)&=\sum_{k=0}^{\lfloor n/q\rfloor}
 x^k[z^{n-qk}]\frac{N_B(z)}{P(z)^{k+1}},\label{qt:eq:finite-B}\\*
 A_n(x)&=\sum_{k=0}^{\lfloor n/p\rfloor}
 x^k[z^{n-pk}]\frac{N_A(z)}{P^\vee(z)^{k+1}}.\label{qt:eq:finite-A}
\end{align}
Here $[z^r]$ extracts the Taylor coefficient at zero. These are
finite polynomial formulas because a term with $qk>n$, or $pk>n$,
has zero coefficient. The denominators have nonzero constant terms,
so each indicated coefficient requires only finitely many arithmetic
operations on the coefficients of $P,N_B,N_A$.

The same functions give finite sums in the spectral roots. Put
$Q(z;x)\coloneq P(z)-xz^q$ and $Q^\vee(z;x)\coloneq P^\vee(z)-xz^p$.
When the corresponding roots are simple, partial fractions give
\begin{align}
 B_n(x)&=-\sum_{Q(\zeta;x)=0}
 \frac{N_B(\zeta)}{Q_z(\zeta;x)\zeta^{n+1}},\label{qt:eq:root-B}\\*
 A_n(x)&=-\sum_{Q^\vee(\zeta;x)=0}
 \frac{N_A(\zeta)}{Q_z^\vee(\zeta;x)\zeta^{n+1}}.\label{qt:eq:root-A}
\end{align}
Indeed, the residue at a simple root is $N(\zeta)/Q_z(\zeta;x)$,
and $[z^n](z-\zeta)^{-1}=-\zeta^{-n-1}$. There is no polynomial
part, since each numerator has degree below $m$. At repeated roots,
the higher-order partial fractions give the derivative formulas of
Proposition~\ref{prop:polynomial-root-formulas}. Thus the boundary
block changes the numerators while preserving the algebraic equations
that generate the two polynomial sequences.

\subsubsection{Densities and point masses from a finite inverse}
The finite resolvent identity used here has classical antecedents.
For scalar Jacobi matrices, Webb and Olver
\cite[Theorems~4.14 and~4.22]{WebbOlver2021} give explicit measures
for finite perturbations of Toeplitz matrices. Bini, Iannazzo, Meini, Meng,
and Robol \cite{BiniIannazzoMeiniMengRobol2023} reduce the operator
eigenvalue problem for finitely representable quasi-Toeplitz matrices
to finite nonlinear equations. The following formulas apply the
finite resolvent identity to the rectangular Favard measures for
arbitrary bandwidths.

Assume $T$ has a bounded PBF, and take the corresponding $H,K$.
For large $|\lambda|$, put
$R_0(\lambda)\coloneq(\lambda I-T_0)^{-1}$ and $G(\lambda)\coloneq E_qR_0(\lambda)E_p^{\mathsf T}$.
Theorem~\ref{thm:general-resolvent} gives every entry of $G$
in roots of $P(z)-\lambda z^q$. Use the multisets $S(\lambda)$
and $\mathcal Z(\lambda)$ from \eqref{eq:general-kernel}, containing
the $q$ small roots and $p$ large roots at infinity, respectively.
With the scalar factor in \eqref{eq:general-gamma}, the formula is
\[
 G_{i,j}(\lambda)=\mathfrak g(\lambda)
 \sum_{k=0}^{\min(i,j)}
 \hc_{j-k}(S(\lambda))\hc_{i-k}(\mathcal Z(\lambda)^{-1}),
\]
for $i\in\{0,\ldots,q-1\}$ and $j\in\{0,\ldots,p-1\}$.
Here $\hc_r$ is defined in \eqref{eq:complete-homogeneous-definition}.

\begin{proposition}[The measure under a finite boundary change]\label{qt:prop:measure}
Let $T=T_0+E_p^{\mathsf T}CE_q$ have a bounded PBF, and let
$H,K$ be its normalized initial matrices in
\eqref{eq:normalized-boundary-definition}. The Stieltjes transform
of its $q\times p$ Favard measure is
\begin{equation}\label{qt:eq:W}
 W_C(\lambda)=H(I_q-G(\lambda)C)^{-1}G(\lambda)K
             =HG(\lambda)(I_p-CG(\lambda))^{-1}K.
\end{equation}
The identity holds initially for large $|\lambda|$ and then by
analytic continuation of the Stieltjes transforms. Let
$\mathcal I=[\lambda_-,\lambda_+]$ be the Toeplitz Favard support interval.
Put $G^\pm(x)\coloneq\lim_{y\downarrow0}G(x\pm\mathrm i y)$.
Except at finitely many points where these values are singular or
the displayed inverses do not exist, the density on $\mathcal I$ is
\begin{equation}\label{qt:eq:density}
 \frac{\dd\Psi_C}{\dd x}
 =H(I_q-G^+C)^{-1}
   \frac{G^--G^+}{2\pi\mathrm i}
   (I_p-CG^-)^{-1}K.
\end{equation}
It is algebraic and has rank one almost everywhere in the interior
of $\mathcal I$. There is no singular continuous part. Outside $\mathcal I$ there are
at most finitely many atoms; they can occur only at uncanceled poles
of \eqref{qt:eq:W}, among the real zeros of
$\det(I_q-GC)=\det(I_p-CG)$. Their mass matrices are the residues.
Endpoint or other exceptional atoms can be recovered from
$\lim_{y\downarrow0}\mathrm i y W_C(x+\mathrm i y)$.
\end{proposition}
\begin{proof}
For sufficiently large $|\lambda|$, both operator resolvents
exist by the convergent Neumann series for bounded operators.
Write $R_C\coloneq(\lambda I-T)^{-1}$ and $G_C\coloneq E_qR_CE_p^{\mathsf T}$.
The resolvent identity gives
$R_C=R_0+R_0E_p^{\mathsf T}CE_qR_C$.
After restricting its rows and columns, this is
$G_C=G+GCG_C$, whence
$G_C=(I_q-GC)^{-1}G$. Multiplication by $H,K$ identifies the Stieltjes transform:
the Neumann expansion and the moment identity give
\[
 HG_C(\lambda)K
 =\sum_{l=0}^{\infty}
   \frac{HE_qT^lE_p^{\mathsf T}K}{\lambda^{l+1}}
 =\sum_{l=0}^{\infty}\frac{\int x^l\dd\Psi_C(x)}{\lambda^{l+1}}
 =\int\frac{\dd\Psi_C(x)}{\lambda-x}.
\]
All series converge for $|\lambda|$ larger than both the operator
norm and the compact support bound supplied by
Corollary~\ref{cor:general-favard}. This proves \eqref{qt:eq:W}; the rectangular alternative follows by
multiplying $(I_q-GC)G=G(I_p-CG)$.

For any matrices $G_1,G_2$ for which the displayed inverses exist,
\[
 (I_q-G_1C)^{-1}G_1-(I_q-G_2C)^{-1}G_2
 =(I_q-G_1C)^{-1}(G_1-G_2)(I_p-CG_2)^{-1}.
\]
Apply this to the upper and lower boundary values. Stieltjes
inversion proves \eqref{qt:eq:density}. Its middle factor has rank
one and is nonzero in the interior of $\mathcal I$, by Theorem~\ref{thm:general-measure} and multiplication by
the inverses of the Toeplitz initial matrices; all other factors are invertible at regular points.

Every entry of $G$ is algebraic. The two determinants are algebraic
and tend to one at infinity. They are therefore nonzero algebraic functions, with only finitely
many zeros and singularities. To see the finiteness of the zeros,
choose the irreducible polynomial of such a function over
$\CC(\lambda)$. After clearing denominators, its constant
coefficient is a nonzero polynomial: otherwise irreducibility
would force the function to be identically zero. Every zero is
among the zeros of this polynomial. The polynomial discriminant
and leading coefficient similarly contain all branch points and
poles. These sets are finite. On each real
interval avoiding these points and the algebraic branch points,
Stieltjes inversion gives a regular density on $\mathcal I$ and zero density
outside it. The remaining singular measure is supported on a finite
set, so it consists of atoms. Positivity and finiteness of each
Favard measure entry ensure that its real poles off $\mathcal I$ are simple,
with residues equal to the masses. The analytic continuation used here is that of the finite
Stieltjes matrix. Its measure support is the one recovered by
Stieltjes inversion, as in Section~\ref{sec:general-polynomials}.
\end{proof}

The finite determinant also identifies exactly when an atom occurs
inside the Toeplitz interval. Under the hypotheses of
Proposition~\ref{qt:prop:measure}, define
\begin{equation}\label{qt:eq:boundary-determinant}
 \begin{aligned}
 \Delta_C(\lambda)
   &\coloneq\det(I_q-G(\lambda)C)=\det(I_p-CG(\lambda)),\\
 d_C(x)&\coloneq\det(I_q-G^+(x)C),\qquad x\in\mathcal I^\circ.
 \end{aligned}
\end{equation}
Here the boundary value is taken in the continuation of the finite
Stieltjes data specified in Proposition~\ref{qt:prop:measure}.

\begin{lemma}[Interior determinant zeros and atoms]
\label{qt:lem:boundary-zero-atom}
Under the hypotheses of Proposition~\ref{qt:prop:measure}, for every
$x_0\in\mathcal I^\circ$,
\begin{equation}\label{qt:eq:boundary-zero-atom}
 d_C(x_0)=0\quad\text{if and only if}\quad
 \Psi_C(\{x_0\})\ne0.
\end{equation}
At such a point the nonremovable poles of the entries of $W_C$
are simple, and their residues form the nonzero mass matrix
$\Psi_C(\{x_0\})$.
\end{lemma}
\begin{proof}
The middle root pair is simple and nonreal at $x_0$, and the two
root groups in the upper boundary value are separated. In a small
complex neighborhood of $x_0$, choose disjoint contours enclosing
these groups. Their symmetric functions are holomorphic there,
including when roots within a group coincide. The root formula for
$G$ therefore gives a holomorphic matrix $G_+$ in this neighborhood
whose restriction to the upper half-plane is the prescribed
continuation. Put
$F_+\coloneq(I_q-G_+C)^{-1}G_+$.
This is meromorphic. Its determinant denominator is not identically
zero: the same algebraic branch tends to one at infinity.

Wherever the inverse exists,
$(I_q-G_+C)^{-1}=I_q+F_+C$.
If $d_C(x_0)=0$ and every entry of $F_+$ were removable at $x_0$,
the right-hand side would extend holomorphically. Multiplying by
$I_q-G_+C$ and evaluating at $x_0$ would give a right inverse of a
singular square matrix. Thus some entry of $F_+$ has a genuine pole.
Since $H,K$ are invertible, the same is true of $W_C=HF_+K$.

Every entry of $W_C$ is the Stieltjes transform of a finite positive
measure. Consequently it is $\mathrm O(y^{-1})$ as $y\downarrow0$
at $x_0+\mathrm i y$, and dominated convergence gives
$\lim_{y\downarrow0}\mathrm i y W_C(x_0+\mathrm i y)=\Psi_C(\{x_0\})$.
Indeed, the scalar kernel
$\mathrm i y/(x_0+\mathrm i y-t)$ has modulus at most one and
converges to the indicator of $\{x_0\}$. A meromorphic pole is
therefore at most simple, and a nonremovable pole has a nonzero
residue equal to the corresponding mass. Conversely, if
$d_C(x_0)\ne0$, the local formula for $W_C$ is holomorphic at
$x_0$, so this limit vanishes. This proves
\eqref{qt:eq:boundary-zero-atom} and the residue assertion.
\end{proof}

\begin{proposition}[Classes with no interior determinant zeros]
\label{qt:prop:interior-nonvanishing}
Under the hypotheses of Proposition~\ref{qt:prop:measure}, each of
the following conditions implies
$d_C(x)\ne0$ for every $x\in\mathcal I^\circ$:
\begin{enumerate}
\item $\min(p,q)=1$;
\item the boundary change affects only the $(0,0)$ entry, so that
      $T_C=T_0+c e_0e_0^{\mathsf T}$ for some $c\in\RR$;
\item $T_C$ is a reordering of the positive constant bidiagonal
      factors in \eqref{qt:eq:boundary}, and $p+q\le5$.
\end{enumerate}
In these cases the Favard measure has no atoms in
$\mathcal I^\circ$, and $\inf_{x\in J}|d_C(x)|>0$ for every compact
interval $J\subset\mathcal I^\circ$.
\end{proposition}
\begin{proof}
For (2), the determinant lemma and the Toeplitz jump formula give
$d_C(x)=1-c\mathfrak g_+(x)$ and $\im\mathfrak g_+(x)=-\pi\rho(x)<0$, where $\mathfrak g_+=(R_0^+)_{0,0}$.
When $c=0$ the determinant equals one; otherwise its imaginary part
is nonzero. This proves (2).

For (1), first suppose $q=1$. The small-root group consists only
of $\zeta$, and \eqref{eq:general-R} gives in the local holomorphic continuation from the upper half-plane
\[
 G(\lambda)=\omega^{-1}
       [\zeta(\lambda),\zeta(\lambda)^2,\ldots,
                                      \zeta(\lambda)^p],\qquad
 F_C(\lambda)\coloneq(I_1-G(\lambda)C)^{-1}G(\lambda)
                 =\frac{G(\lambda)}{\Delta_C(\lambda)}.
\]
If also $p=1$, then $C$ is a $1\times1$ matrix and
$T_C=T_0+C_{0,0}e_0e_0^{\mathsf T}$, so (2) applies.
Suppose $p\ge2$ and $d_C(x_0)=0$. The first entry of $G$ is nonzero at $x_0$, so the
first entry of $F_C$ has a nonremovable pole there. Since
$F_C=H^{-1}W_CK^{-1}$ and the initial matrices are real, each entry
of $F_C$ is the Stieltjes transform of a finite real signed measure.
The bound and limit used in the preceding lemma show that this
pole is simple and has a real nonzero residue $r_0$. But
$(F_C)_{0,1}=\zeta(F_C)_{0,0}$, so the second entry has residue
$\zeta(x_0)r_0$, which is nonreal. This contradicts the reality of
its signed measure. Thus $d_C$ never vanishes. Transposing treats
$p=1$: the transposed system again has a bounded PBF, its finite
block is $G^{\mathsf T}$, and Sylvester's determinant identity
identifies its determinant with $\Delta_C$.

For (3), cases with $\min(p,q)=1$ are already covered. The remaining
pairs are $(2,2),(2,3),(3,2)$. First take $p=2$ and
$q\in\{2,3\}$. Place the factors in their given order on a circle.
The two gaps between the lower factors contain altogether $q$
upper factors, so one gap contains zero or one. If it contains
zero, a cyclic rotation places all the upper factors before the
two adjacent lower factors and gives $T_0$. If it contains one,
a cyclic rotation gives
$W=\omega U_*\mathcal L(b)\mathcal U(c)\mathcal L(d)$,
where $U_*$ contains the other $q-1$ upper factors. Let
$\alpha>0$ be the product of their diagonal parameters. The identities
\[
 \mathcal L(b)\mathcal U(c)
   =\mathcal U(c)\mathcal L(b)-b e_0e_0^{\mathsf T},
 \qquad U_*e_0=\alpha e_0,
 \qquad e_0^{\mathsf T}\mathcal L(d)=e_0^{\mathsf T}
\]
give $W=T_0-\kappa e_0e_0^{\mathsf T}$ with
$\kappa=\omega\alpha b>0$. Thus in both cases there is a cyclic
rotation $W=T_0-\kappa e_0e_0^{\mathsf T}$ with $\kappa\ge0$.
The factors of the same triangular type commute, so the product
with all upper factors first is exactly the common $T_0$.

Regularity transfers from $W$ to the original product $T_w$ through
the cyclic intertwining identities.
Fix $x_0\in\mathcal I^\circ$ and write the factors in the order
of $W$ as $F_1,\ldots,F_m$, where $m=p+q$. Set
\[
 \begin{aligned}
 V_j&\coloneq\omega F_{j+1}\cdots F_mF_1\cdots F_j,
             &&0\le j<m,\qquad V_m\coloneq V_0,\\
 A_j&\coloneq F_{j+1},\qquad
 B_j\coloneq\omega F_{j+2}\cdots F_mF_1\cdots F_j,
             &&0\le j<m.
 \end{aligned}
\]
Empty products are the identity. Then $V_0=W$, some $V_r$ equals
$T_w$, and $V_j=A_jB_j$, $V_{j+1}=B_jA_j$. All these matrices
are bounded and have finite bandwidth. Choose $M$ larger than
their operator norms and that of $T_0$. On $|\lambda|>M$, write
$R_V=(\lambda I-V)^{-1}$. The resolvent identities give
\[
 \begin{aligned}
 (R_W)_{i,n}
   &=(R_0)_{i,n}
     -\frac{\kappa(R_0)_{i,0}(R_0)_{0,n}}
                  {1+\kappa\mathfrak g},\\
 (R_{V_{j+1}})_{i,n}
   &=\lambda^{-1}\left(\delta_{i,n}
       +\sum_{k,\ell\ge0}(B_j)_{i,k}(R_{V_j})_{k,\ell}
                                      (A_j)_{\ell,n}\right).
 \end{aligned}
\]
The second identity follows from
$(\lambda I-BA)^{-1}=\lambda^{-1}[I+B(\lambda I-AB)^{-1}A]$.
For fixed $i,n$, only finitely many summands occur: row $i$ of
$B_j$ and column $n$ of $A_j$ have finite support.

Continue these scalar entries from the common exterior; all subsequent
statements of holomorphy refer to those continued entries.
Every entry of $R_0$ given by \eqref{eq:general-R} is holomorphic
in the upper half-plane. At $x_0$, the separated root groups in
\eqref{eq:general-root-boundary-sets} give holomorphic symmetric
functions on a disk $D$ about $x_0$; hence every entry has a
holomorphic continuation to that same disk. Coincidences inside
a group do not affect its symmetric functions. Moreover,
$\mathfrak g=(R_0)_{0,0}$ is the Stieltjes transform of the
nonzero positive measure with density $\rho$, so
$\im\mathfrak g(\lambda)<0$ in the upper half-plane and
$\im\mathfrak g_+(x_0)=-\pi\rho(x_0)<0$.
Thus $1+\kappa\mathfrak g$ has no zero in the upper half-plane
or at $x_0$. Shrink $D$ so that this denominator is nonzero
throughout $D$ and $0\notin D$; this is possible because
$x_0>0$. The first resolvent identity makes every continued
entry of $R_W$ holomorphic on the upper half-plane together with
$D$. Successive use of the second identity preserves this
property: each entry is a finite sum of holomorphic functions
multiplied by $\lambda^{-1}$. Induction from $V_0=W$ through
$V_r=T_w$ therefore proves the same regularity for every
continued entry of $R_{T_w}$.

For $|\lambda|>M$ sufficiently large, the moment identity
identifies $HE_qR_{T_w}E_p^{\mathsf T}K$ with the Stieltjes
matrix $W_C$. Both are holomorphic in the upper half-plane,
so the identity theorem identifies them there. The continued
finite block just constructed is therefore a holomorphic
extension of $W_C$ to $D$. In particular its mass at $x_0$
vanishes, and Lemma~\ref{qt:lem:boundary-zero-atom} gives
$d_C(x_0)\ne0$. Since $x_0$ was arbitrary, this proves the
claim for $(2,2)$ and $(2,3)$.

Finally,
$\mathcal L(c)^{\mathsf T}=c\mathcal U(c^{-1})$ and
$\mathcal U(c)^{\mathsf T}=c\mathcal L(c^{-1})$.
Transposition therefore takes $(3,2)$ to the proved case $(2,3)$,
with a positive overall scalar and reciprocal parameters. Its
Toeplitz reference is $T_0^{\mathsf T}$, its finite block is
$G^{\mathsf T}$, and its boundary change is $C^{\mathsf T}$.
The finite Toeplitz sections have the same eigenvalues after
transposition, so their limiting interval is the same $\mathcal I$.
Sylvester's identity gives
$\det(I_p-G^{\mathsf T}C^{\mathsf T})
 =\det(I_p-CG)=\det(I_q-GC)$,
so nonvanishing transfers back to the original determinant.

The absence of interior atoms follows from
Lemma~\ref{qt:lem:boundary-zero-atom}. The positive lower bound
on each compact interval follows from continuity and nonvanishing.
\end{proof}

\begin{example}[A common tail with different stability and a mass at zero]
\label{ex:quasi-toeplitz-two-two}
Set $\omega=1$ and
$T_a\coloneq\mathcal U(1)\LL(3)\mathcal U(4)\LL(2)$, with $P(z)\coloneq(1+z)(1+2z)(1+3z)(1+4z)$.
The product with all upper factors first is $T_0\coloneq\mathcal U(1)\mathcal U(4)\LL(2)\LL(3)$,
and the commutator gives $T_a=T_0-3e_0e_0^{\mathsf T}$.
The limit exchanges end at
$(\LL(1),\LL(2),\mathcal U(3),\mathcal U(4))$, so Theorem
\ref{qt:thm:classification} proves analyticity and interior membership.

For the actual positive exchanges, first interchange $\mathcal U(1)\LL(3)$
and $\mathcal U(4)\LL(2)$, then interchange the remaining middle
upper--lower pair. The first lower entry is $3/4$. The rightmost raw
upper factor has diagonal entry $6$ at row zero; its unit
normalization has superdiagonal $1/6$ there. The Favard definitions
therefore give
\[
 H=\begin{bmatrix}1&0\\6&1\end{bmatrix},\qquad
 K=\begin{bmatrix}1&4/3\\0&1\end{bmatrix}.
\]
Proposition~\ref{qt:prop:generating} becomes
\begin{align*}
 \sum_{n\ge0}B_n(x)z^n
 &=\frac{[1+4z+3z^2,\ z]}
 {1+10z+(35-x)z^2+50z^3+24z^4},\\
 \sum_{n\ge0}A_n(x)z^n
 &=\frac{[24+18z+3z^2,\ 24z]^{\mathsf T}}
 {24+50z+(35-x)z^2+10z^3+z^4}.
\end{align*}
For example, the coefficients of $z^2$ give
$B_2(x)=[x+28,\,-10]$ and $A_2(x)=[x/24+13/9,\,-25/12]^{\mathsf T}$.
Let $s$ be the sum of the two small roots of $P(z)-\lambda z^2$,
$t$ the sum of the reciprocals of the two large roots, and
$\mathfrak g$ the negative product of the two small roots. Then
\[
 G=\mathfrak g\begin{bmatrix}1&s\\t&1+st\end{bmatrix},\qquad
 W_C=\frac{\mathfrak g}{1+3\mathfrak g}
 \begin{bmatrix}1&s+4/3\\t+6&(t+6)(s+4/3)\end{bmatrix}
 +\mathfrak g\begin{bmatrix}0&0\\0&1\end{bmatrix}.
\]
This determines the matrix of measures by its jump and its poles.

For comparison, $T_b\coloneq\LL(2)\LL(3)\mathcal U(1)\mathcal U(4)$ has the
same $P$, but
\[
 C_b=\begin{bmatrix}-31&-5\\-30&-6\end{bmatrix}.
\]
Its limiting groups are $\{2,3\}$ below and $\{1,4\}$ above.
It has a bounded positive factorization, but is not interior and
admits no analytic choice of PBF on a full-band neighborhood.
The difference is also visible in the measure. The two middle root
parameters are $3>2$, so the common Toeplitz interval has
$\lambda_->0$. The constant PBF of $T_b$ in the displayed order gives
\[
 H_b=\begin{bmatrix}1&0\\4&1\end{bmatrix},\qquad
 K_b=\begin{bmatrix}1&1/2\\0&1\end{bmatrix}.
\]
At $\lambda=0$ the two small roots are $-1/4,-1/3$ and the other
two are $-1/2,-1$. Differentiating
$P(\zeta(\lambda))-\lambda\zeta(\lambda)^2=0$ gives
$\zeta'(0)=\zeta(0)^2/P'(\zeta(0))$. Substitution in the root
formula yields
\[
 G(0)=\begin{bmatrix}-1/12&7/144\\1/4&-11/48\end{bmatrix},\qquad
 G'(0)=\begin{bmatrix}-5/72&59/864\\13/36&-157/432\end{bmatrix}.
\]
For $D(\lambda)\coloneq I_2-G(\lambda)C_b$, multiplication gives
\[
 D(0)=\frac18\begin{bmatrix}-1&-1\\7&7\end{bmatrix},\qquad
 D'(0)=\frac1{48}\begin{bmatrix}-5&3\\14&-18\end{bmatrix}.
\]
The equal columns give $\det D(0)=0$. Differentiating the two
products in the determinant gives
$(\det D)'(0)=(-35+18-21+14)/384=-1/16$.
The simple-pole formula therefore proves the exact mass identity
\begin{equation}\label{qt:eq:zero-mass}
 \Psi_b(\{0\})
 =\frac{H_b\operatorname{adj}D(0)G(0)K_b}
 { (\det D)'(0)}
 =\begin{bmatrix}2/3&1/9\\2&1/3\end{bmatrix}.
\end{equation}
The matrix $T_a$ has boundary denominator
$1+3\mathfrak g(0)=3/4$, so its measure has no atom at zero.
Thus the same four distinct tail roots allow a boundary change that
creates a point mass at zero and destroys interior membership.

\end{example}

\subsubsection{Interior oscillations after a boundary change}
\label{qt:subsec:interior-asymptotics}
\noeqref{qt:eq:boundary-scalar-data,qt:eq:boundary-spectral-functions,qt:eq:real-rank-one-density,qt:eq:boundary-biorthogonality,qt:eq:boundary-kernel,qt:eq:boundary-reproduction,qt:eq:left-oscillations,qt:eq:right-oscillations,qt:eq:right-transport,qt:eq:left-transport,qt:eq:boundary-sine-kernel,qt:eq:boundary-diagonal-limit,qt:eq:transported-amplitudes}

B\"ottcher, Fukshansky, Garcia, and Maharaj
\cite{BottcherFukshanskyGarciaMaharaj2015} derived exact and asymptotic
determinant formulas for finite Toeplitz sections with fixed corner blocks.
For scalar Jacobi matrices, Breuer, Last, and Simon
\cite{BreuerLastSimon2014} proved persistence of the local sine kernel
under a change in the first diagonal entry, with the same local spacing
and a modified spectral weight. For the semi-infinite banded system
below, a bounded PBF and nonvanishing of the initial-block determinant
$d_C$ on an interior interval give the modified density
$\rho_C=\rho/|d_C|^2$. The argument of $d_C$ supplies the common phase
shift of the two weighted forms. Their amplitudes determine the
normalized mixed kernel in Theorem~\ref{qt:thm:boundary-sine}.

Let $T_C=T_0+E_p^{\mathsf T}CE_q$, $C\in\RR^{p\times q}$, where $T_0$
is the BTP Toeplitz matrix in \eqref{eq:general-toeplitz}, with
$p,q\ge1$. Assume that $T_C$ has a PBF with bounded factors.
Let $H_C,K_C$, $A_n^C,B_n^C$, and $\Psi_C$ be its Favard data from
Corollary~\ref{cor:general-favard}. Use the interval $\mathcal I$,
root $\zeta=r\mathrm e^{-\mathrm i\theta}$, root groups, amplitudes
$C_A,C_B$, and densities $\rho,\varpi$ of
Section~\ref{subsec:spectral-kernels} for $T_0$.
Thus these quantities describe the common Toeplitz tail.

Write $R_0^\pm$ for the entrywise boundary values of the root formula
\eqref{eq:general-R}, and put $G^\pm\coloneq E_qR_0^\pm E_p^{\mathsf T}$.
Define the column and row
\[
 f_0\coloneq[\hc_i(\widehat{\mathcal Z}^{-1})]_{i\in\{0,\dots,q-1\}},\qquad
 g_0^{\mathsf T}\coloneq[\hc_j(\widehat{\mathcal S})]_{j\in\{0,\dots,p-1\}}.
\]
Use $d_C$ from \eqref{qt:eq:boundary-determinant}. For
$x\in\mathcal I^\circ$, set
\begin{equation}\label{qt:eq:boundary-scalar-data}
 \begin{aligned}
 f_C&\coloneq\operatorname{adj}(I_q-G^+C)f_0,&
 g_C^{\mathsf T}&\coloneq g_0^{\mathsf T}\operatorname{adj}(I_p-CG^-).
 \end{aligned}
\end{equation}
At points where $d_C\ne0$, define
\begin{equation}\label{qt:eq:boundary-spectral-functions}
 \rho_C\coloneq\frac{\rho}{|d_C|^2},\qquad
 \mathcal A_n^C\coloneq g_C^{\mathsf T}K_CA_n^C,\qquad
 \mathcal B_n^C\coloneq B_n^CH_Cf_C,
 \qquad n\in\NN.
\end{equation}
These are weighted linear forms in the components of $A_n^C$
and $B_n^C$, with weight vectors $K_C^{\mathsf T}g_C$ and
$H_Cf_C$, respectively. They are scalar algebraic functions, as
in Theorem~\ref{prop:spectral-functions}.
The exceptional interior points form a finite set by
Proposition~\ref{qt:prop:measure}. That proposition also gives the
full decomposition of $\Psi_C$ into its density on $\mathcal I$
and finitely many atoms, with no singular continuous part.
Proposition~\ref{qt:prop:interior-nonvanishing} makes the condition
$d_C\ne0$ automatic throughout $\mathcal I^\circ$ in its three
stated classes; in those cases every compact interior interval
is admissible in the results below.

\begin{lemma}[Linear forms and the atomic contribution]
\label{qt:lem:scalar-biorthogonality}
The vectors $f_C,g_C$ are real, and the density of $\Psi_C$ at
regular interior points is
\begin{equation}\label{qt:eq:real-rank-one-density}
 \frac{\dd\Psi_C}{\dd x}=\rho_CH_Cf_Cg_C^{\mathsf T}K_C.
\end{equation}
Let $x_1,\ldots,x_s$ be all its atoms, including any endpoint
or exceptional interior atoms, and put $M_j\coloneq\Psi_C(\{x_j\})$.
Then, for $n,l\in\NN$,
\begin{equation}\label{qt:eq:boundary-biorthogonality}
 \int_{\mathcal I}\mathcal B_n^C(x)\mathcal A_l^C(x)
                    \rho_C(x)\dd x
 +\sum_{j=1}^{s}B_n^C(x_j)M_jA_l^C(x_j)=\delta_{n,l}.
\end{equation}
The linear forms on the continuous part and $\rho_C$ are independent of
the chosen PBF of $T_C$.
\end{lemma}
\begin{proof}
The Toeplitz jump formula gives
$G^- -G^+=2\pi\mathrm i\rho f_0g_0^{\mathsf T}$.
For a square matrix $M$ and column and row vectors $u,v^{\mathsf T}$,
\[
 \operatorname{adj}(M+tuv^{\mathsf T})u=\operatorname{adj}(M)u,
 \qquad
 v^{\mathsf T}\operatorname{adj}(M+tuv^{\mathsf T})
 =v^{\mathsf T}\operatorname{adj}(M).
\]
For invertible matrices these identities follow by combining the
rank-one inverse formula with the determinant formula. They hold
for all matrices because both sides are polynomial in their entries.
Apply the first identity to $I_q-G^+C$, with $u=f_0$, and the
second to $I_p-CG^-$, with $v^{\mathsf T}=g_0^{\mathsf T}$.
The adjugate-vector products are unchanged when $G^+$ and $G^-$
are interchanged. These two boundary values are conjugate, while
$f_0,g_0$ and $C$ are real. Hence $f_C,g_C$ are real.

In \eqref{qt:eq:density}, replace the two inverses by their
adjugates divided by their determinants. The determinants are
$d_C$ and $\bar d_C$, respectively. This proves
\eqref{qt:eq:real-rank-one-density}.
Proposition~\ref{qt:prop:measure} supplies the full decomposition
of $\Psi_C$ into this density and finitely many atoms, with no
singular continuous part. Substitution in Favard biorthogonality
\eqref{eq:general-biorthogonality} proves
\eqref{qt:eq:boundary-biorthogonality}. Its integrand is integrable
because \eqref{qt:eq:real-rank-one-density} identifies it with
$B_n^C(\dd\Psi_C/\dd x)A_l^C$.

Finally, $B_n^CH_C$ and $K_CA_n^C$ are the recurrence solutions
with identity initial matrices. The unchanged positive extreme
diagonals determine these solutions uniquely. Thus they depend
only on $T_C$. The vectors $f_C,g_C$ and the weight $\rho_C$
depend only on $C,T_0$, proving the last assertion.
\end{proof}

Define the matrix kernel and its scalar contraction by
\begin{equation}\label{qt:eq:boundary-kernel}
 \begin{aligned}
 \mathscr K_N^C(x,y)&\coloneq\sum_{n=0}^{N-1}A_n^C(x)B_n^C(y),\\
 k_N^C(x,y)&\coloneq g_C(x)^{\mathsf T}K_C\mathscr K_N^C(x,y)H_Cf_C(y)
           =\sum_{n=0}^{N-1}\mathcal A_n^C(x)\mathcal B_n^C(y).
 \end{aligned}
\end{equation}
If atoms are present, their terms in
\eqref{qt:eq:boundary-biorthogonality} also occur in kernel
reproduction. More explicitly, expansion of the finite sums gives
\begin{equation}\label{qt:eq:boundary-reproduction}
 \begin{aligned}
 k_N^C(x,y)
 &=\int_{\mathcal I}k_N^C(x,t)k_N^C(t,y)\rho_C(t)\dd t\\
 &\quad+g_C(x)^{\mathsf T}K_C
   \left(\sum_{j=1}^{s}\mathscr K_N^C(x,x_j)M_j
                          \mathscr K_N^C(x_j,y)\right)H_Cf_C(y).
 \end{aligned}
\end{equation}
The following limits concern $x,y$ in the regular continuous
interval and use this same kernel.

\begin{theorem}[Boundary phase and interior zeros]
\label{qt:thm:boundary-oscillations}
Let $T_C=T_0+E_p^{\mathsf T}CE_q$ admit a PBF with bounded
factors, where $T_0$ is the BTP Toeplitz matrix in
\eqref{eq:general-toeplitz}. Use the weighted linear forms in
\eqref{qt:eq:boundary-spectral-functions}. Let $J$ be a compact
interval in $\mathcal I^\circ$ on which $d_C\ne0$. Choose continuous arguments
$\delta_A\coloneq\arg C_A$, $\delta_B\coloneq\arg C_B$, and
$\delta_C\coloneq\arg d_C$. There exists $\vartheta\in(0,1)$ such that,
uniformly on $J$,
\begin{align}
 r^{-n}\mathcal A_n^C
 &=\frac{|d_C||C_A|}{\sin\theta}
     \sin((n+1)\theta-\delta_A+\delta_C)
       +\mathrm O(\vartheta^n),
       &&n\to\infty,\label{qt:eq:left-oscillations}\\
 r^n\mathcal B_n^C
 &=\frac{|d_C||C_B|}{\sin\theta}
     \sin((n+1)\theta+\delta_B+\delta_C)
       +\mathrm O(\vartheta^n),
       &&n\to\infty.\label{qt:eq:right-oscillations}
\end{align}
The derivatives of the remainders are $\mathrm O(n\vartheta^n)$
as $n\to\infty$. On every compact interval inside $J^\circ$,
for all sufficiently large $n$, each zero of either displayed sine
has exactly one corresponding real zero of its scalar function
within distance $\mathrm O(\vartheta^n/n)$ as $n\to\infty$;
these zeros are simple.
For $[a,b]\subset J^\circ$, the number of real zeros of each scalar
function in $[a,b]$, divided by $n$, tends to
$\int_a^b\varpi(x)\dd x$ as $n\to\infty$.
Consecutive zeros converging to $x_0\in J^\circ$ have spacing
$(n\varpi(x_0))^{-1}(1+\mathrm o(1))$ as $n\to\infty$.
\end{theorem}
\begin{proof}
Let $a_n^0\coloneq\hc_n(\widehat{\mathcal S})$ and
$b_n^0\coloneq\hc_n(\widehat{\mathcal Z}^{-1})$ be the Toeplitz scalar
functions in \eqref{eq:spectral-functions-roots}.
The first step is to establish the exact identities
\begin{align}
 \mathcal B_n^C
 &=d_Cb_n^0+\sum_{j=0}^{p-1}[R_0^+]_{n,j}(Cf_C)_j,
       \label{qt:eq:right-transport}\\
 \mathcal A_n^C
 &=\bar d_Ca_n^0+
   \sum_{i=0}^{q-1}(g_C^{\mathsf T}C)_i[R_0^-]_{i,n}.
       \label{qt:eq:left-transport}
\end{align}
The right-hand side of \eqref{qt:eq:right-transport} has first
$q$ entries $d_Cf_0+G^+Cf_C=f_C$, because
$(I_q-G^+C)f_C=d_Cf_0$.
Applying $T_0-xI$ to this sequence gives $-E_p^{\mathsf T}Cf_C$.
The boundary block cancels that expression, since the sequence
has first $q$ entries $f_C$. It therefore solves $T_Cb=xb$ with
the same initial entries as $\mathcal B_n^C$. Uniqueness of the
right recurrence proves \eqref{qt:eq:right-transport}.
All products here are finite in each row, so the algebraic
boundary values suffice. The left identity follows by the same
argument, using
$g_C^{\mathsf T}(I_p-CG^-)=\bar d_Cg_0^{\mathsf T}$.

For fixed $j$ and $n\ge j$, \eqref{eq:general-R} gives
\[
 [R_0^+]_{n,j}=\mathfrak g_+
   \sum_{k=0}^{j}\hc_{j-k}(\mathcal S_0\cup\{\zeta\})
                \hc_{n-k}(\mathcal Z_0^{-1}\cup\{\bar\zeta^{-1}\}).
\]
After multiplication by $r^n$, the term associated with
$\bar\zeta^{-1}$ is proportional to $\mathrm e^{-\mathrm i n\theta}$;
all the other terms decrease geometrically. There is no term
proportional to $\mathrm e^{\mathrm i n\theta}$ in this correction.
Thus the coefficient of $\mathrm e^{\mathrm i n\theta}$ in
$r^n\mathcal B_n^C$ is $d_C$ times its Toeplitz coefficient.
Since the scalar function is real, the coefficient of
$\mathrm e^{-\mathrm i n\theta}$ is its conjugate.
Formula~\eqref{eq:spectral-right-asymptotic} now gives
\eqref{qt:eq:right-oscillations}.

For the left sequence, the lower boundary value has root groups
$\mathcal S_0\cup\{\bar\zeta\}$ and
$\mathcal Z_0\cup\{\zeta\}$. Hence, for fixed $i$ and $n\ge i$,
\[
 [R_0^-]_{i,n}=\mathfrak g_-
   \sum_{k=0}^{i}\hc_{n-k}(\mathcal S_0\cup\{\bar\zeta\})
                \hc_{i-k}(\mathcal Z_0^{-1}\cup\{\zeta^{-1}\}).
\]
After division by $r^n$, the correction has only the term
$\mathrm e^{\mathrm i n\theta}$ and geometrically decreasing terms.
The coefficient of $\mathrm e^{-\mathrm i n\theta}$ therefore
remains $\bar d_C$ times its Toeplitz value. Reality determines
the other coefficient, giving \eqref{qt:eq:left-oscillations}.
This explains why the same additional phase $\delta_C$ occurs
on both sides.

The remaining roots, after scaling, have modulus uniformly less
than one on $J$. Their generating functions are analytic on a
uniformly larger disk. Cauchy coefficient estimates give the
geometric remainders and their derivatives, as in the proof of
Theorem~\ref{thm:spectral-oscillations}. These estimates include
coincidences inside either separated group, because the symmetric
generating functions remain analytic. The factors $d_C,d_C^{-1}$
and their derivatives are bounded on $J$. Both leading amplitudes
are nonzero there.

The derivatives of the new phases are
$-(n+1)\pi\varpi-\delta_A'+\delta_C'$ and
$-(n+1)\pi\varpi+\delta_B'+\delta_C'$.
Their leading negative terms have size proportional to $n$,
while the other terms are bounded. Near each multiple of $\pi$
the sine therefore changes sign with derivative of that size.
The remainder displaces the zero by $\mathrm O(\vartheta^n/n)$,
and its derivative bound makes the zero unique and simple, as
$n\to\infty$. Counting the multiples of $\pi$ gives the limiting
zero density; applying the mean value theorem to consecutive
phases gives the spacing.
\end{proof}

\begin{theorem}[Interior kernel limit with a boundary block]
\label{qt:thm:boundary-sine}
Under the hypotheses of Theorem~\ref{qt:thm:boundary-oscillations},
fix $x_0\in J^\circ$ and put $\varpi_0\coloneq\varpi(x_0)$,
$x_N\coloneq x_0+u/(N\varpi_0)$, and $y_N\coloneq x_0+v/(N\varpi_0)$.
For the kernel \eqref{qt:eq:boundary-kernel},
\begin{equation}\label{qt:eq:boundary-sine-kernel}
 \frac{\rho_C(x_0)}{N\varpi_0}
 \left(\frac{r(y_N)}{r(x_N)}\right)^N k_N^C(x_N,y_N)
 =\frac{\sin\pi(u-v)}{\pi(u-v)}+\mathrm O(N^{-1}),
 \qquad N\to\infty,
\end{equation}
uniformly for $u,v$ in compact subsets of $\RR$.
The sine quotient is interpreted as one at $u=v$. In particular,
\begin{equation}\label{qt:eq:boundary-diagonal-limit}
 \frac{\rho_C(x)}{N}k_N^C(x,x)
 =\varpi(x)+\mathrm O(N^{-1}),\qquad N\to\infty,
\end{equation}
uniformly on $J$.
\end{theorem}
\begin{proof}
In the convention of \eqref{eq:spectral-amplitudes}, the new
amplitudes are $C_A^C\coloneq\bar d_CC_A$ and $C_B^C\coloneq d_CC_B$.
Their product satisfies the exact identity
\begin{equation}\label{qt:eq:transported-amplitudes}
 \frac{\rho_CC_A^CC_B^C}{2\sin^2\theta}
 =\frac{\rho C_AC_B}{2\sin^2\theta}
 =-\frac{\mathrm i}{\pi}\frac{\zeta'}{\zeta}
 =\varpi-\frac{\mathrm i}{\pi}(\log r)'.
\end{equation}
The last two equalities are
\eqref{eq:spectral-amplitude-identity}. Thus the boundary factor
cancels from the weighted product of the amplitudes.

Apply the summation argument of
Theorem~\ref{thm:spectral-sine-kernel} to
\eqref{qt:eq:left-oscillations}--\eqref{qt:eq:right-oscillations}.
For completeness, put $c_0\coloneq(\log r)'(x_0)/\varpi_0$ and
$\Delta\coloneq u-v$. The difference
$\delta_C(x_N)-\delta_C(y_N)$ and the relative changes of the
amplitudes are $\mathrm O(N^{-1})$ as $N\to\infty$.
The sum of the geometric remainders is bounded, and the
series involving the sum of the phases is bounded because
$\mathrm e^{2\mathrm i\theta(x_0)}\ne1$.
After division by $N$, the remaining sum involving the difference of the phases is a
Riemann sum with error $\mathrm O(N^{-1})$ as $N\to\infty$.
The real and imaginary parts of \eqref{qt:eq:transported-amplitudes}
give its limiting integral as
\[
 \int_0^1\mathrm e^{c_0\Delta t}
 \left(\cos(\pi\Delta t)+\frac{c_0}{\pi}\sin(\pi\Delta t)\right)\dd t
 =\mathrm e^{c_0\Delta}\frac{\sin(\pi\Delta)}{\pi\Delta}.
\]
The equality at $\Delta=0$ follows by continuity.
The factor $(r(y_N)/r(x_N))^N$ equals
$\mathrm e^{-c_0\Delta}(1+\mathrm O(N^{-1}))$ as $N\to\infty$.
This proves \eqref{qt:eq:boundary-sine-kernel}, uniformly across
the diagonal. Taking $u=v=0$ gives
\eqref{qt:eq:boundary-diagonal-limit}; the estimates are uniform
when $x_0$ ranges over $J$.
\end{proof}

These formulas apply to both matrices in
Example~\ref{ex:quasi-toeplitz-two-two}, on compact interior
intervals where the corresponding $d_C$ does not vanish.
For $T_a$, $d_C=1+3\mathfrak g_+$; for $T_b$ it is
$\det(I_2-G^+C_b)$. The atom of $T_b$ at zero contributes to
\eqref{qt:eq:boundary-biorthogonality}, while the interior kernel
has the same limit as for $T_a$ and the common Toeplitz tail.

In the scalar Jacobi case, boundary-symbol formulas for the
measure and its atoms are given by Geronimo and Iliev
\cite[Theorem~2.1]{GeronimoIliev2017} and by Webb and Olver
\cite[Theorem~4.22]{WebbOlver2021}. The squared denominator and
additional phase above describe the same effect for the two scalar
functions selected by the rectangular Favard measure.

\subsection{The binomial family: explicit polynomials and matrix measure}\label{sec:binomial}

Take $P(z)=(1+z)^{p+q}$ in \eqref{eq:general-toeplitz}. Thus all root parameters in
\eqref{eq:root-parameters} equal one, and the Toeplitz entries
are binomial coefficients. In this case the individual
bidiagonal factors and the positive diagonal have closed formulas.
Proposition~\ref{prop:pbf} verifies the factorization by a binomial
convolution. Its partial products determine the initial-condition
matrices, and Theorem~\ref{thm:polynomials} gives finite sums for
the polynomial families. Theorem~\ref{thm:measure} computes the
matrix density. A trigonometric parametrization makes the support
explicit, and Theorem~\ref{thm:hardedge} and
Proposition~\ref{prop:softedge} describe its two endpoints.
The scalar case $p=q=1$ is given in
Example~\ref{ex:scalar-chebyshev}.

Fix $p,q\in\mathbb N$, and put $m\coloneq p+q$. Throughout, a binomial coefficient with an
integer lower argument outside its usual range is zero. Set
\begin{equation}\label{eq:T}
 T\coloneq\left[\binom{m}{q+i-j}\right]_{i,j\in\NN},\qquad
 a(z)\coloneq\frac{(1+z)^m}{z^q},\qquad
 c\coloneq\frac{m^m}{p^p q^q}.
\end{equation}
Thus $T$ has $p$ subdiagonals and $q$ superdiagonals. This is the
repeated-root case of Theorem~\ref{thm:toeplitz}; its individual PBF
factors have the following closed formulas.

For $r\ge1$, let $\mathsf L_r$ be the unit lower bidiagonal matrix with
\begin{equation}\label{eq:Lr}
 (\mathsf L_r)_{n,n-1}\coloneq\frac{n}{n+r}\quad(n\in\mathbb N).
\end{equation}
Let $D$ be the positive diagonal matrix
\begin{equation}\label{eq:D}
 D_{n,n}\coloneq\frac{(n+m)!\,n!}{(n+p)!\,(n+q)!},\qquad n\in\NN.
\end{equation}

\begin{proposition}[An explicit PBF for the binomial matrix]\label{prop:pbf}
Every finite leading section of the matrix $T$ in \eqref{eq:T} is BTP.
With $\mathsf L_r$ and $D$ defined by \eqref{eq:Lr} and \eqref{eq:D},
a bounded PBF of $T$ is
\begin{equation}\label{eq:pbf}
 T=\mathsf L_q^{\,p}D(\mathsf L_p^{\mathsf T})^{q}.
\end{equation}
\end{proposition}
\begin{proof}
Write $\mathsf L_q=I+N_q$, where $(N_q)_{n,n-1}\coloneq n/(n+q)$ and all other
entries vanish. The only nonzero entry of $N_q^r$ in row $i$ is
at column $i-r$, and it equals
\[
 \prod_{h=0}^{r-1}\frac{i-h}{i-h+q}
 =\frac{i!\,(i-r+q)!}{(i-r)!\,(i+q)!},
 \qquad r\in\{0,\dots,i\}.
\]
Expanding $(I+N_q)^p$, and applying the same calculation before
transposing the upper factor, gives
\begin{align*}
 (\mathsf L_q^p)_{i,k}&=\binom p{i-k}
 \frac{i!\,(k+q)!}{k!\,(i+q)!},&
 ((\mathsf L_p^{\mathsf T})^q)_{k,j}&=\binom q{j-k}
 \frac{j!\,(k+p)!}{k!\,(j+p)!}.
\end{align*}
After multiplying by $D_{k,k}$, the factors involving $k+p$ and
$k+q$ cancel. The entry $(i,j)$ of the proposed product becomes
\[
 \frac{m!\,i!\,j!}{(i+q)!\,(j+p)!}
 \sum_{k\in\NN}\binom p{i-k}\binom q{j-k}\binom{m+k}{k}.
\]
It remains to use the identity
\begin{equation}\label{eq:convolution}
 \sum_{k\in\NN}\binom p{i-k}\binom q{j-k}\binom{m+k}{k}
 =\binom{i+q}{j}\binom{j+p}{i}.
\end{equation}
For the left-hand side, put $i=k+r$, $j=k+\ell$ and sum over
$i,j\in\NN$. The resulting generating function is
\[
 \left(\sum_{r=0}^p\binom pr s^r\right)
 \left(\sum_{\ell=0}^q\binom q\ell t^\ell\right)
 \left(\sum_{k\in\NN}\binom{m+k}{k}(st)^k\right)
 =\frac{(1+s)^p(1+t)^q}{(1-st)^{m+1}}.
\]
The right-hand side has the same generating function. One way to
check this last assertion using only Pascal's identity is to put
$c_{p,q}(i,j)\coloneq\binom{i+q}{j}\binom{j+p}{i}$, with value zero when
either index is negative. Pascal's identity gives
$c_{p+1,q}(i,j)-c_{p+1,q}(i-1,j-1)=c_{p,q}(i,j)+c_{p,q}(i-1,j)$.
For the generating functions this means multiplication by
$(1+s)/(1-st)$ when $p$ is increased by one. Interchanging
$(p,i,s)$ with $(q,j,t)$ gives multiplication by $(1+t)/(1-st)$
when $q$ is increased by one. The initial case is
$c_{0,0}(i,j)=\delta_{i,j}$, with generating function $(1-st)^{-1}$.
Induction on $p+q$ proves the assertion and hence
\eqref{eq:convolution}.

Substitution in the product entry above gives
$m!/[(q+i-j)!(p+j-i)!]=\binom m{q+i-j}$ whenever the indices lie
in the band, and gives zero otherwise. This proves the factorization.
The bidiagonal entries are between zero and one, and $D_{n,n}$ tends to
one as $n\to\infty$. The finite leading sections factor by the corresponding finite
factors: every summation index is at most $\min(i,j)$ in the
lower-diagonal-upper product. The implication $\PBF\Rightarrow\BTP$
then follows from Proposition~\ref{prop:positive-product}.
\end{proof}

For $(p,q)=(2,1)$, the lower product in \eqref{eq:pbf} is
$L=\mathsf L_1^2$. Its signed inverse
$\mathcal Q=\Sigma L^{-1}\Sigma$, with $\Sigma=\diag[(-1)^n]_{n\in\NN}$,
has $\mathcal Q_{n,0}=1$ and $\mathcal Q_{n,1}=2n/(n+1)$ for $n\ge1$,
so \eqref{eq:global-minor-ratios} gives
$\kappa_{1,n}=(n+1)/(2n)\to1/2$ as $n\to\infty$.
Thus $(\mathsf L_1)_{1,0}=1/2$ attains the endpoint of the admissible
interval $(0,1/2]$ in Proposition~\ref{rem:bandwidth-two-endpoint}.

\begin{example}[Lima's constant-coefficient system]
\label{ex:lima-system}
For $p=r$ and $q=1$, the binomial matrix in \eqref{eq:T} recovers
the constant-coefficient hypergeometric recurrence of
\cite[Theorem~8.5(a), equation~(8.25)]{Lima2023} by an explicit
algebraic diagonal rescaling. Put
$h_r\coloneq r^r/(r+1)^{r+1}$ and
$S_r\coloneq\operatorname{diag}[h_r^n]_{n\in\NN}$.
Write $B_n\coloneq B_n^{(1)}$ for the right polynomial family of
Theorem~\ref{thm:general-toeplitz-polynomials}, specialized to
$P(z)=(1+z)^{r+1}$ and $q=1$. The monic recurrence matrix
$\mathcal T_r$ and polynomials $P_n$ in Lima's normalization
are related to this family by the entrywise identities
\begin{equation}\label{eq:lima-rescaling}
 \mathcal T_r=h_rS_rTS_r^{-1},\qquad
 P_n(x)=h_r^n B_n(x/h_r).
\end{equation}
Indeed, the entry on subdiagonal $k$ becomes
$\binom{r+1}{k+1}h_r^{k+1}$, and the superdiagonal remains one.
For $r=2$, $h_2=4/27$ and the resulting monic recurrence matrix
is the one in \cite[Section~4.3]{Uniform}, with parameters
$(4/3,5/3,2,5/2)$, or with the first two interchanged.
Lima also gives hypergeometric polynomial and Meijer--$G$ measure
formulas for general $r$. Transposition covers $p=1,q=r$.
The present construction provides a rectangular mixed-type
normalization for arbitrary positive bandwidths $p,q$.
The geometric-root family in Section~\ref{sec:q-binomial} contains
this constant-coefficient polynomial system at $\tau=1$, through
the same rescaling \eqref{eq:lima-rescaling}.
\end{example}

The binomial factorization makes the normalization of the initial conditions
explicit. This gives finite sums for both polynomial families and
identifies the corresponding mixed-type orthogonality.

Use the arrangement in \eqref{eq:polynomial-arrays}: the left
polynomial vectors are the columns of $A(x)\in\RR[x]^{p\times\NN}$,
and the right vectors are the rows of $B(x)\in\RR[x]^{\NN\times q}$.
Their components are
$A_n(x)=[A_n^{(1)}(x),\ldots,A_n^{(p)}(x)]^{\mathsf T}$ and $B_n(x)=[B_n^{(1)}(x),\ldots,B_n^{(q)}(x)]$.
The row and column normalization in
\eqref{eq:normalized-boundary-definition} gives the following
initial-condition matrices for \eqref{eq:pbf}, with indices starting at zero:
\begin{equation}\label{eq:boundary}
 \begin{aligned}
 K_{i,j}&=\begin{cases}\binom{j+q}{i+q},&i\le j,\\0,&i>j,
 \end{cases}\quad &&i,j\in\{0,\dots,p-1\},\\
 H_{i,j}&=\begin{cases}\binom{i+p}{j+p},&j\le i,\\0,&j>i,
 \end{cases}\quad &&i,j\in\{0,\dots,q-1\}.
 \end{aligned}
\end{equation}
Thus $K\in\RR^{p\times p}$ is upper triangular and
$H\in\RR^{q\times q}$ is lower triangular, both with diagonal
entries equal to one. Set
\begin{equation}\label{eq:initial}
 [A_0\ \cdots\ A_{p-1}]=K^{-1},\qquad
 [B_0^{\mathsf T}\ \cdots\ B_{q-1}^{\mathsf T}]^{\mathsf T}=H^{-1}.
\end{equation}
Inverses are obtained by multiplying the entry in position $(i,j)$
by $(-1)^{i-j}$.

\begin{theorem}[Finite formulas for the binomial polynomial families]\label{thm:polynomials}
For the binomial matrix $T$ in \eqref{eq:T}, the solutions of
$AT=xA$ and $TB=xB$ with the initial matrices in
\eqref{eq:boundary}--\eqref{eq:initial} are,
for $n\in\NN$, $a\in\{1,\dots,p\}$ and $b\in\{1,\dots,q\}$,
\begin{align}
 B_n^{(b)}(x)
 &=\sum_{k=0}^{\lfloor(n-b+1)/q\rfloor}
 (-1)^{n-b+1-qk}
 \binom{n+p(k+1)}{n-b+1-qk}x^k,
 \label{eq:B}\\
 A_n^{(a)}(x)
 &=\sum_{k=0}^{\lfloor(n-a+1)/p\rfloor}
 (-1)^{n-a+1-pk}
 \binom{n+q(k+1)}{n-a+1-pk}x^k,
 \label{eq:A}
\end{align}
A sum with negative upper limit is zero. Every nonzero component has
the degree indicated by the upper summation limit.
\end{theorem}
\begin{proof}
For the right components the generating functions are
\begin{equation}\label{eq:GF}
 \sum_{n\in\NN}B_n^{(b)}(x)z^n
 =\frac{z^{b-1}(1+z)^{q-b}}{(1+z)^m-xz^q}.
\end{equation}
Expanding the denominator first in $xz^q$ gives
$\sum_{k\in\NN}x^k z^{qk+b-1}(1+z)^{-mk-p-b}$.
For a term with fixed $k$, set $r\coloneq n-b+1-qk$ and $M\coloneq mk+p+b$.
Here $[z^r]F(z)$ denotes the coefficient of $z^r$ in the power
series $F(z)$. The negative-binomial identity
$[z^r](1+z)^{-M}=(-1)^r\binom{M+r-1}{r}$,
valid for $r\in\NN$,
gives $M+r-1=n+p(k+1)$. Only $r\ge0$ contributes, so the resulting
sum is exactly \eqref{eq:B}, with the stated upper limit.
For $n<q$, the coefficients are the entries of $H^{-1}$.
Multiplication of \eqref{eq:GF} by its denominator gives the right
recurrence, including the boundary convention $B_n=0$ for $n<0$.
For the left components, the reversed symbol is also $(1+z)^m$,
and the left generating function in
Theorem~\ref{thm:general-toeplitz-polynomials} is
\[
 \sum_{n\in\NN}A_n^{(a)}(x)z^n
 =\frac{z^{a-1}(1+z)^{p-a}}{(1+z)^m-xz^p}
 =\sum_{k\in\NN}x^kz^{pk+a-1}(1+z)^{-mk-q-a}.
\]
Now $r\coloneq n-a+1-pk$ and $M\coloneq mk+q+a$, so $M+r-1=n+q(k+1)$.
The same coefficient identity gives \eqref{eq:A}. Its first $p$
coefficients are the rows of $K^{-1}$, and multiplication by the
denominator gives $AT=xA$.
The last binomial coefficient in each nonempty sum is positive, so its
leading coefficient is nonzero.
\end{proof}

The normalization is the one associated with \eqref{eq:pbf}. Indeed,
the unnormalized matrix $\widetilde K$ in \eqref{eq:raw-K-definition},
formed here from the first columns of
$I,\mathsf L_q,\ldots,\mathsf L_q^{p-1}$, has entries
\[
 \widetilde K_{i,j}=\frac{\binom ji}{\binom{i+q}{q}}.
\]
Dividing each column by its diagonal entry, as in
\eqref{eq:normalized-boundary-definition}, gives $K$ in
\eqref{eq:boundary}. The transposed construction from the first rows of
$I,\mathsf L_p^{\mathsf T},\ldots,(\mathsf L_p^{\mathsf T})^{q-1}$ gives $H$ after
row normalization. These matrices are triangular totally positive in
the sense of positivity of all structurally admissible minors.

Apply Corollary~\ref{cor:general-favard} to the bounded PBF
\eqref{eq:pbf} with the initial matrices \eqref{eq:boundary}.
This gives the unique $q\times p$ matrix $\Psi$ of entrywise positive,
compactly supported finite measures with moments
\eqref{eq:general-moments} and biorthogonality
\eqref{eq:general-biorthogonality}. The mixed-type relations follow
from Corollary~\ref{cor:general-step}. Since all root parameters
equal one, Theorem~\ref{thm:general-support} gives
$\lambda_-=0$, $z_+=q/p$, and $\lambda_+=a(q/p)=c$; hence
$\supp\Psi\subseteq[0,c]$.

The finite coefficients also give a differential equation for each
binomial component. Its order is $m=p+q$. The derivation below
keeps the component index and the initial normalization explicit.
It also identifies the equations with generalized hypergeometric
equations, extending to both sides of the mixed-type system the
description for $q=1$ in \cite{Lima2023}.
Write $\mathscr E_x\coloneq x\,\mathrm d/\mathrm dx$ for the
Euler differential operator, so that $\mathscr E_x x^k=kx^k$.

\begin{proposition}[Differential equations for the binomial components]
\label{prop:binomial-differential-equations}
Let $p,q\ge1$, $m\coloneq p+q$, and let $B_n^{(b)},A_n^{(a)}$ be
the polynomials of Theorem~\ref{thm:polynomials}.
For $b\in\{1,\dots,q\}$ and $n\ge b-1$,
\begin{equation}\label{eq:binomial-B-differential}
 \Biggl[
 \prod_{j=0}^{m-1}(m\mathscr E_x+b-q+j)
 -x\prod_{j=1}^{p}(p\mathscr E_x+n+p+j)
       \prod_{j=0}^{q-1}(q\mathscr E_x-n+b-1+j)
 \Biggr]B_n^{(b)}(x)=0.
\end{equation}
For $a\in\{1,\dots,p\}$ and $n\ge a-1$,
\begin{equation}\label{eq:binomial-A-differential}
 \Biggl[
 \prod_{j=0}^{m-1}(m\mathscr E_x+a-p+j)
-x\prod_{j=1}^{q}(q\mathscr E_x+n+q+j)
       \prod_{j=0}^{p-1}(p\mathscr E_x-n+a-1+j)
 \Biggr]A_n^{(a)}(x)=0.
\end{equation}
In each equation the products are polynomials in $\mathscr E_x$;
the multiplication by $x$ in the second term is performed after
that product of operators. The equation determines the formal
power series solution uniquely from its constant coefficient.
\end{proposition}
\begin{proof}
Fix $b$ and $n\ge b-1$, and put $N\coloneq n-b+1$ and
$d\coloneq\lfloor N/q\rfloor$. By \eqref{eq:B},
$B_n^{(b)}(x)=\sum_{k=0}^d c_kx^k$, where
\[
 c_k\coloneq(-1)^{N-qk}
 \frac{(n+p+pk)!}{(N-qk)!\,(p+b-1+mk)!}.
\]
For $1\le k\le d$, dividing consecutive coefficients gives
\begin{equation}\label{eq:binomial-coefficient-ratio}
 \frac{c_k}{c_{k-1}}=(-1)^q
 \frac{\prod_{j=1}^{p}(n+pk+j)
                    \prod_{j=1}^{q}(N-qk+j)}
      {\prod_{j=0}^{m-1}(mk+b-q+j)}.
\end{equation}
The denominator comes from the $m$ consecutive factors between
$(p+b-1+m(k-1))!$ and $(p+b-1+mk)!$.
The two numerator products similarly come from the change of
$p$ and $q$ in the other factorials.

To translate this identity into an equation, define the scalar
polynomials
$D_b(s)\coloneq\prod_{j=0}^{m-1}(ms+b-q+j)$ and $R_{n,b}(s)\coloneq\prod_{j=1}^{p}(ps+n+p+j)\prod_{j=0}^{q-1}(qs-N+j)$.
The $q$ factors in the second product at $s=k-1$ are the
negatives of $N-qk+1,\dots,N-qk+q$, in reverse order.
Thus \eqref{eq:binomial-coefficient-ratio} says exactly that
$D_b(k)c_k=R_{n,b}(k-1)c_{k-1}$.
Because $\mathscr E_x$ multiplies $x^k$ by $k$, this is the
vanishing of the coefficient of $x^k$ in
$[D_b(\mathscr E_x)-xR_{n,b}(\mathscr E_x)]B_n^{(b)}$.

There are two boundary coefficients to check. For $k=0$,
$D_b(0)=0$: its factor with $j=q-b$ is zero.
For $k=d+1$, write $N=qd+r$ with $0\le r<q$.
The second product of $R_{n,b}(d)$ has a zero factor at $j=r$,
so the coefficient of $x^{d+1}$ also vanishes. All higher
coefficients are zero. This proves
\eqref{eq:binomial-B-differential} at every degree.

For every integer $k\ge1$, all factors of $D_b(k)$ are positive,
since its smallest factor is at least $p+b>0$.
Consequently the coefficient recurrence determines each $c_k$
from $c_{k-1}$, starting with
$c_0=(-1)^N\binom{n+p}{N}$; the zero of $R_{n,b}(d)$ makes
the series terminate. This proves the uniqueness assertion.
Formula \eqref{eq:A} is \eqref{eq:B} with $p,q$ interchanged
and $b$ replaced by $a$. Performing those substitutions in
the coefficient calculation gives
\eqref{eq:binomial-A-differential} and its uniqueness assertion.
\end{proof}

The coefficient of the highest derivative in either displayed
operator is $x^m(m^m-p^p q^q x)$. Indeed, the leading coefficient
of $\mathscr E_x^m$ is $x^m$, and the two products have leading
coefficients $m^m$ and $p^p q^q$ as polynomials in $\mathscr E_x$.
Every coefficient of the differential equation has a factor $x$,
because $D_b(0)=0$, and likewise on the left.
After cancelling that common factor, the leading coefficient is
$p^p q^q x^{m-1}(c-x)$, where
$c=m^m/(p^p q^q)$ is the upper endpoint in \eqref{eq:T}.
The finite singular points of these equations are thus precisely
the endpoints $0,c$ of the binomial support.

\begin{remark}[Hypergeometric form and the scalar equation]
\label{rem:binomial-hypergeometric-equations}
For $N\coloneq n-b+1\ge0$, form the following lists of upper and lower
parameters:
\[
 \begin{aligned}
 \boldsymbol\alpha_{n,b}
 &\coloneq\left(\frac{n+p+i}{p}\right)_{i\in\{1,\dots,p\}}
       \mathbin{\cup}\left(\frac{-N+j}{q}\right)_{j\in\{0,\dots,q-1\}},&
 \boldsymbol\beta_b
 &\coloneq\left(1+\frac{b-q+j}{m}\right)_
          {\substack{0\le j\le m-1\\j\ne q-b}}.
 \end{aligned}
\]
Here the union concatenates lists, preserving multiplicities.
The lower list omits the parameter $1$ corresponding to the
factor $\mathscr E_x$ in $D_b(\mathscr E_x)$.
With the usual generalized hypergeometric notation
\cite[Sections~16.2 and 16.8]{NISTDLMF},
\begin{equation}\label{eq:binomial-hypergeometric}
 B_n^{(b)}(x)=(-1)^N\binom{n+p}{N}\,
 {}_mF_{m-1}\left(
 \begin{matrix}\boldsymbol\alpha_{n,b}\\\boldsymbol\beta_b\end{matrix}
 ;\frac{x}{c}\right).
\end{equation}
For clarity, the coefficient of $z^k$ in
${}_mF_{m-1}(\boldsymbol\alpha;\boldsymbol\beta;z)$ is
$\prod_i(\alpha_i)_k/(k!\prod_j(\beta_j)_k)$, with
$(u)_0\coloneq1$ and $(u)_k\coloneq u(u+1)\cdots(u+k-1)$.
All lower parameters above are positive. If $N=qd+r$,
$0\le r<q$, the upper parameter with $j=r$ is $-d$,
which terminates the series. The ratio of its consecutive
coefficients, including the factor $1/c=p^p q^q/m^m$,
is \eqref{eq:binomial-coefficient-ratio}; its constant coefficient
also agrees. This proves \eqref{eq:binomial-hypergeometric}.
The left formula follows by the same interchange $p\leftrightarrow q$,
$b\leftrightarrow a$.

For $q=1$, the rescaling
$P_n^{\mathrm{Lima}}(x)=c^{-n}B_n^{(1)}(cx)$ in
\eqref{eq:lima-rescaling} turns
\eqref{eq:binomial-hypergeometric} into
\cite[Theorem~8.5, equation~(8.26)]{Lima2023}.
The differential equations are correspondingly instances of the
generalized hypergeometric equation; compare
\cite[Proposition~6.4]{Lima2023} and
\cite[equation~16.8.3]{NISTDLMF}.
For $p=q=1$, both components coincide and, after cancelling $x$,
\eqref{eq:binomial-B-differential} becomes
\[
 x(4-x)y''(x)+(6-3x)y'(x)+n(n+2)y(x)=0.
\]
This is the Chebyshev equation for
$y(x)=U_n((x-2)/2)$, in the normalization of
Example~\ref{ex:scalar-chebyshev}.
\end{remark}

For $|\lambda|$ sufficiently large, write
\begin{equation}\label{eq:kernel}
 Q(z;\lambda)=(1+z)^m-\lambda z^q
 =\prod_{u\in S}(z-u)\prod_{v\in\mathcal Z}(z-v).
\end{equation}
Here $S=S(\lambda)$ consists of the $q$ roots tending to zero as
$|\lambda|\to\infty$, and $\mathcal Z=\mathcal Z(\lambda)$ consists
of the other $p$ roots, which tend to infinity.
Theorem~\ref{thm:general-resolvent}, with $\omega=t_p=1$, gives
the Stieltjes matrix $W$ through \eqref{eq:general-W}.
The measure extends $W$ analytically to $\CC\setminus[0,c]$,
where it is algebraic by the finite root formula
\eqref{eq:general-R}. The following lemma parametrizes the central
conjugate pair trigonometrically and fixes its boundary orientation.

\begin{lemma}[Trigonometric parametrization of the middle roots]\label{lem:roots}
For $0<x<c$, the roots of the kernel $Q(z;x)$ in \eqref{eq:kernel},
with $c$ as in \eqref{eq:T}, contain a unique middle conjugate
pair $\zeta,\bar\zeta$, with $\im\zeta<0$, having exactly $q-1$
roots of smaller modulus and $p-1$ of larger modulus.
In this lemma and its proof, $\theta$ denotes the polar angle in
\eqref{eq:spectral-polar-data} divided by $m$. The pair is
parametrized by
\begin{align}
 \zeta(\theta)&=\frac{\sin(q\theta)}{\sin(p\theta)}\mathrm{e}^{-\mathrm{i}m\theta},
 \label{eq:alpha}\\
 x(\theta)&=\frac{\sin^m(m\theta)}
 {\sin^p(p\theta)\sin^q(q\theta)},
 \qquad 0<\theta<\frac\pi m.\label{eq:trig}
\end{align}
The function $x(\theta)$ decreases strictly from $c$ to zero. The
small-root set continued from the upper half-plane has boundary value
$S_+=\mathcal S_0\cup\{\zeta\}$, where $\mathcal S_0$ is the set of $q-1$ smaller roots.
The large-root set has boundary value $\mathcal Z_+=\mathcal Z_0\cup\{\bar\zeta\}$,
where $\mathcal Z_0$ is the set of $p-1$ larger roots.
\end{lemma}
\begin{proof}
The addition formula for sine gives
$\sin(p\theta)+\sin(q\theta)\mathrm e^{-\mathrm i m\theta}=\sin(m\theta)\mathrm e^{-\mathrm i q\theta}$.
Dividing by $\sin(p\theta)$ yields
$1+\zeta(\theta)=\sin(m\theta)\mathrm e^{-\mathrm i q\theta}/
\sin(p\theta)$. In $(1+\zeta)^m/\zeta^q$, the two complex
phases are both $\mathrm e^{-\mathrm i mq\theta}$ and cancel.
The remaining positive quotient is exactly \eqref{eq:trig}. The only finite critical values of $a(z)$ are
$0$ and $c$: away from $z=-1,0$,
$a'(z)=a(z)\frac{pz-q}{z(1+z)}$.
The parametrized $\zeta$ is nonreal in the open parameter interval,
and $\zeta'(\theta)\ne0$ because its argument has derivative $-m$.
Thus $x'(\theta)=a'(\zeta)\zeta'\ne0$. The endpoint limits
give strict decrease.

For this symbol, Theorem~\ref{thm:general-support} gives
$\lambda_-=0$, $z_+=q/p$ and $\lambda_+=c$, together with the
separation of the middle conjugate pair on $(0,c)$.
As $\theta\downarrow0$, the displayed $\zeta(\theta)$ tends to
$q/p$ with negative imaginary part. It is therefore the lower member
of the middle pair near $c$. Simplicity and separation of this pair
extend the identification throughout $(0,c)$. The boundary orientation in \eqref{eq:general-root-boundary-sets}
places it in $S_+$, the limit of the small-root group from the
upper half-plane.
\end{proof}

Using the smaller and larger root groups $\mathcal S_0,\mathcal Z_0$
and the middle pair from Lemma~\ref{lem:roots}, define
\begin{equation}\label{eq:binomial-density-root-sets}
 \widehat{\mathcal S}\coloneq\mathcal S_0\cup\{\zeta,\bar\zeta\},\qquad
 \widehat{\mathcal Z}\coloneq\mathcal Z_0\cup\{\zeta,\bar\zeta\}.
\end{equation}
These conjugation-invariant multisets have $q+1$ and $p+1$ elements,
respectively. Define
\begin{align}
 \rho(x)&\coloneq\frac{-\im\zeta(x)}{\pi}\prod_{u\in \mathcal S_0}|u|,
 \label{eq:rho}\\
 f_b(x)&\coloneq\hc_{b-1}\bigl((1+v^{-1})_{v\in\widehat{\mathcal Z}}\bigr),
 &&b\in\{1,\dots,q\},\label{eq:f}\\
 g_a(x)&\coloneq\hc_{a-1}\bigl((1+u)_{u\in\widehat{\mathcal S}}\bigr),
 &&a\in\{1,\dots,p\}.\label{eq:g}
\end{align}
Empty products equal one. In particular $f_1=g_1=1$.

\begin{theorem}[Algebraic matrix weight]\label{thm:measure}
The binomial measure obtained from
Corollary~\ref{cor:general-favard} is purely absolutely continuous. With $\rho,f_b,g_a$ defined by \eqref{eq:rho}--\eqref{eq:g},
its density is
\begin{equation}\label{eq:weight}
 \dd\Psi(x)=\rho(x)
 \begin{bmatrix}f_1(x)\\\vdots\\f_q(x)\end{bmatrix}
 \begin{bmatrix}g_1(x)&\cdots&g_p(x)\end{bmatrix}\dd x,
 \qquad 0<x<c.
\end{equation}
The scalar density $\rho$ is positive in the open interval;
$f_b,g_a$ are real and nonnegative there. The matrix density has rank
exactly one at every interior point, and $\supp\Psi=[0,c]$.
\end{theorem}
\begin{proof}
Apply Theorem~\ref{thm:general-measure} to the bounded PBF
\eqref{eq:pbf}. Its initial matrices are \eqref{eq:boundary}.
The scalar density in \eqref{eq:general-density-data} is exactly
\eqref{eq:rho}. It remains to evaluate the two matrix-vector products
that define $f,g$ for these binomial initial matrices.

The normalization by $H,K$ uses the identity
\begin{equation}\label{eq:h-shift}
 \hc_k(1+y_1,\ldots,1+y_r)
 =\sum_{i=0}^k\binom{r+k-1}{r+i-1}\hc_i(y_1,\ldots,y_r).
\end{equation}
To prove this identity, use the generating function
\[
 \begin{aligned}
 \prod_{j=1}^{r}(1-(1+y_j)t)^{-1}
 &=(1-t)^{-r}\prod_{j=1}^{r}
       \left(1-\frac{y_jt}{1-t}\right)^{-1}\\
 &=\sum_{i\in\NN}\hc_i(y_1,\ldots,y_r)t^i(1-t)^{-r-i}.
 \end{aligned}
\]
The coefficient of $t^k$ on the last line is
$\sum_{i=0}^k\binom{r+k-1}{k-i}\hc_i(y_1,\ldots,y_r)$,
which is \eqref{eq:h-shift}.
For the column density vector, $\widehat{\mathcal Z}$ has $p+1$ roots and
the coefficient $\binom{p+i}{p+j}$ is $H_{i,j}$. Thus multiplication
by $H$ replaces each variable $v^{-1}$ by $1+v^{-1}$.
For the row density vector, $\widehat{\mathcal S}$ has $q+1$ roots and
$K_{i,j}=\binom{q+j}{q+i}$ performs the replacement $u\mapsto1+u$.
This gives \eqref{eq:f}--\eqref{eq:weight} from the general density formula.

Absolute continuity, positivity, rank one, and exact support
$[0,c]$ now follow from Theorem~\ref{thm:general-measure}.
\end{proof}

At the origin, the coalescing roots produce different powers in the
entries of the binomial density. The following theorem determines
these powers and their leading constants. At the upper endpoint, a double-root
expansion gives the square-root behavior.

\begin{theorem}[Behavior of all density entries at zero]\label{thm:hardedge}
For the density factors $\rho,f_b,g_a$ in
\eqref{eq:rho}--\eqref{eq:g}, put $\phi\coloneq\pi/m$ and define positive constants
$F_b\coloneq\prod_{j=1}^{b-1}\frac{\sin((p+j)\phi)}{\sin(j\phi)}$ and
$G_a\coloneq\prod_{j=1}^{a-1}\frac{\sin((q+j)\phi)}{\sin(j\phi)}$.
Then
\begin{align}
 \rho(x)&=\frac{\sin(q\phi)}{\pi}x^{1/m}
 [1+\mathrm O(x^{1/m})],&&x\downarrow0,\label{eq:rho-hard}\\
 f_b(x)&=F_bx^{(b-1)/m}[1+\mathrm O(x^{1/m})],&&x\downarrow0,\label{eq:f-hard}\\
 g_a(x)&=G_ax^{(a-1)/m}[1+\mathrm O(x^{1/m})],&&x\downarrow0.\label{eq:g-hard}
\end{align}
Consequently every entry has a nonzero leading term:
\begin{equation}\label{eq:entry-hard}
 \frac{\dd\Psi_{b,a}}{\dd x}(x)
 =\frac{\sin(q\phi)}{\pi}F_bG_a\,
 x^{(a+b-1)/m}[1+\mathrm O(x^{1/m})],\qquad x\downarrow0.
\end{equation}
\end{theorem}
\begin{proof}
Set $\upsilon\coloneq x^{1/m}$ and write $z=-1+\upsilon w$.
The root equation becomes
$w^m=(-1+\upsilon w)^q$.
At $\upsilon=0$ its roots are simple roots of $(-1)^q$; the implicit
function theorem gives analytic expansions in $\upsilon$.
The middle root has
$1+\zeta=\upsilon \mathrm{e}^{-\mathrm{i}q\phi}+\mathrm O(\upsilon^2)$ as
$\upsilon\downarrow0$.
The product in \eqref{eq:rho} tends to one, proving
\eqref{eq:rho-hard}.

Recall the root multisets $\widehat{\mathcal S},\widehat{\mathcal Z}$
from \eqref{eq:binomial-density-root-sets}. The rescaled quantities
$(1+u)/\upsilon$, for $u\in\widehat{\mathcal S}$, tend, as $\upsilon\downarrow0$, to the $q+1$ points
$\mathrm{e}^{\mathrm{i}(-q+2j)\phi}$, $j\in\{0,\dots,q\}$. This follows either from
\eqref{eq:alpha} or from ordering
$|-1+\upsilon w|=1-\upsilon\operatorname{Re}w
+\mathrm O(\upsilon^2)$ as $\upsilon\downarrow0$.
Likewise $(1+v^{-1})/\upsilon$, for $v\in\widehat{\mathcal Z}$, tend,
as $\upsilon\downarrow0$, to
$\mathrm{e}^{\mathrm{i}(-p+2j)\phi}$, $j\in\{0,\dots,p\}$.
The following coefficient identity evaluates these limits.
For a parameter $t$, put
$H_r(z)\coloneq\prod_{j=0}^r(1-t^jz)^{-1}$ and
$h_k\coloneq\hc_k(1,t,\ldots,t^r)$. The identity
$(1-t^{r+1}z)H_r(tz)=(1-z)H_r(z)$ implies
$(1-t^k)h_k=(1-t^{r+k})h_{k-1}$ for $k\ge1$.
Starting with $h_0=1$, it follows that
\[
 h_k=\prod_{j=1}^k\frac{1-t^{r+j}}{1-t^j}
\]
whenever the denominators are nonzero. Set
$t=\mathrm e^{2\mathrm i\phi}$ and multiply all variables by
$\mathrm e^{-\mathrm i r\phi}$. Homogeneity of degree $k$
contributes $\mathrm e^{-\mathrm i rk\phi}$, and the phase
from the $k$ displayed quotients is its inverse. Hence
\[
 \hc_k(\mathrm{e}^{-\mathrm{i}r\phi},\mathrm{e}^{\mathrm{i}(-r+2)\phi},\ldots,\mathrm{e}^{\mathrm{i}r\phi})
 =\prod_{j=1}^k\frac{\sin((r+j)\phi)}{\sin(j\phi)}.
\]
The relevant degrees satisfy $k<m$, so none of these denominators
vanishes. Each rescaled root equals its limiting value plus
$\mathrm O(\upsilon)$ as $\upsilon\downarrow0$. Since $\hc_k$
is a polynomial of degree $k$, substitution gives its limiting
value plus $\mathrm O(\upsilon)$, with the overall factor
$\upsilon^k$. Taking $(r,k)=(p,b-1)$ and $(q,a-1)$ gives
\eqref{eq:f-hard} and \eqref{eq:g-hard}. In the ranges used here,
all numerator and denominator angles lie strictly between zero and
$\pi$, so the constants are positive. Multiplication proves
\eqref{eq:entry-hard}.
\end{proof}

\begin{proposition}[Upper endpoint]\label{prop:softedge}
Let $c$ be as in \eqref{eq:T}, and let $\mathcal S_{0,c}$ be the
limit as $x\uparrow c$ of the smaller-root multiset $\mathcal S_0$
in Lemma~\ref{lem:roots}. For the scalar density $\rho$ in
\eqref{eq:rho},
\begin{equation}\label{eq:softedge}
 \rho(x)=\frac1\pi\left(\prod_{u\in \mathcal S_{0,c}}|u|\right)
 \sqrt{\frac{2mq}{c p^3}}\sqrt{c-x}\,[1+\mathrm O(\sqrt{c-x})],\quad x\uparrow c.
\end{equation}
The functions $f_b,g_a$ in \eqref{eq:f}--\eqref{eq:g} remain bounded
at $c$, and every entry of the density in \eqref{eq:weight} is
$\mathrm O(\sqrt{c-x})$ as $x\uparrow c$.
\end{proposition}
\begin{proof}
For the Laurent symbol $a(z)=(1+z)^m/z^q$ in \eqref{eq:T},
$a'(z_*)=0$ at $z_*=q/p$ and $a(z_*)=c$.
The logarithmic derivative is $m/(1+z)-q/z$. Its derivative
at $z_*$ equals
$-m/(1+z_*)^2+q/z_*^2=p^3/(mq)$,
so $a''(z_*)=cp^3/(mq)>0$. Taylor expansion now gives
$x-c=\tfrac12a''(z_*)(\zeta-z_*)^2+\mathrm O((\zeta-z_*)^3)$ as $\zeta\to z_*$.
Writing $\zeta-z_*=\sqrt{c-x}\,w$ and dividing by $c-x$
shows that $w$ tends to a root of
$-1=\tfrac12a''(z_*)w^2$. The choice $\im\zeta<0$ selects
$w=-\mathrm i\sqrt{2/a''(z_*)}$. These roots are simple, so
the implicit function theorem gives an error
$\mathrm O(\sqrt{c-x})$ in $w$ as $x\uparrow c$. Thus
$-\im\zeta(x)=\sqrt{2/a''(z_*)}\sqrt{c-x}
+\mathrm O(c-x)$ as $x\uparrow c$.
The derivative is
$a'(z)=(1+z)^{m-1}(pz-q)/z^{q+1}$.
Its other zero is $z=-1$, where $a(-1)=0\ne c$. Thus only
$z_*=q/p$ can be a multiple root of $P(z)-cz^q$; all the other
roots are simple. They are nonzero because the constant term of
this polynomial is one. Their product of moduli tends to the
positive value displayed in
\eqref{eq:softedge}. Substituting into \eqref{eq:rho} proves
that formula. The remaining assertions follow from the finite
expressions \eqref{eq:f}--\eqref{eq:g}.
\end{proof}

For $p=q=2$, the general density has the following expression in radicals.

\begin{example}[A two-by-two binomial weight]
\label{ex:two-by-two-binomial}
For $p=q=2$, the matrix has diagonal $6$, adjacent diagonals $4$, and
outer diagonals $1$. With $s\coloneq\sqrt{x}$, put
\[
 r(x)\coloneq\frac{s+\sqrt{s(s+4)}}2,\qquad
 w(x)\coloneq\frac{\sqrt{s(4-s)}}{2\pi[1+r(x)]}.
\]
The general measure reduces to
\begin{equation}\label{eq:22}
 \dd\Psi(x)=w(x)
 \begin{bmatrix}1&r(x)\\r(x)&r(x)^2\end{bmatrix}\dd x,
 \qquad 0<x<16.
\end{equation}
Indeed, with $s=\sqrt{x}$ the polynomial in the root equation factors as
$[(1+z)^2-sz][(1+z)^2+sz]$.
The first quadratic supplies the central unit-modulus pair with
$\zeta+\bar\zeta=s-2$. The other two roots are
$-1/(1+r)$ and $-(1+r)$, because $r^2=s(1+r)$.
Substitution into \eqref{eq:rho}--\eqref{eq:g} gives
$\rho=w$ and $f_2=g_2=r$, proving \eqref{eq:22} directly.
\end{example}

\subsection{The family with geometrically spaced roots}
\label{sec:q-binomial}

Geometrically spaced root parameters give a family of Toeplitz
matrices $T_\tau$ within the class of Section~\ref{sec:general-polynomials}.
Gaussian binomial identities verify an explicit PBF. Its ordered
partial products give the initial-condition matrices and hence both
polynomial generating functions. The matrix measure then follows from the
general algebraic root formula. Any ordered partition of these
root parameters also gives polynomial families by
Theorem~\ref{thm:general-toeplitz-polynomials}. The constant changes
\eqref{eq:q-change-polynomials}--\eqref{eq:q-change-measure} below
relate those families and their measures to the present ones.
At $\tau=1$, the formulas recover the binomial family with scale
$\alpha$.

The parameters $\alpha>0$ and $\tau>0$ specify the scale and ratio
of the root parameters. The base usually denoted by $q$ in basic
hypergeometric identities is written $\tau$ here; the integer $q$
continues to denote the upper bandwidth. These identities evaluate
the factors and generating functions; the positive measures are
supplied by the PBF and the spectral Favard theorem.
For $N\in\NN$ and $r\in\{0,\dots,N\}$,
define the finite product and Gaussian binomial coefficient by
\begin{equation}\label{eq:gaussian-notation}
 (a;\tau)_N\coloneq\prod_{\nu=0}^{N-1}(1-a\tau^\nu),\qquad
 \gbinom{N}{r}\coloneq
 \frac{(\tau;\tau)_N}{(\tau;\tau)_r(\tau;\tau)_{N-r}}.
\end{equation}
Set $\gbinom{N}{r}\coloneq0$ for integer $r\notin\{0,\dots,N\}$. At $\tau=1$ all
formulas are interpreted by continuous extension.

\begin{lemma}[$\tau$-Vandermonde--Pfaff convolution]
\label{lem:q-vandermonde-pfaff}
For $p,q,i,j\in\NN$ and $m\coloneq p+q$,
\begin{equation}\label{eq:q-vandermonde-pfaff}
 \sum_{k\in\NN}\tau^{(i-k)(j-k)}
 \gbinom{p}{i-k}\gbinom{q}{j-k}\gbinom{m+k}{k}
 =\gbinom{i+q}{j}\gbinom{j+p}{i}.
\end{equation}
\end{lemma}
\begin{proof}
Suppose $i\le j$ and put $d\coloneq j-i$. If $d>q$, both sides vanish.
Assume $d\in\{0,\dots,q\}$, and set $r=i-k$. After factoring out
the $r=0$ term, the left-hand side is
\[
 \gbinom{q}{d}\gbinom{m+i}{i}
 {}_3\phi_2\!\left[
 \begin{matrix}\tau^{-p},\tau^{-(q-d)},\tau^{-i}\\
 \tau^{d+1},\tau^{-(m+i)}\end{matrix};\tau,\tau\right].
\]
Here is an explicit check of the conversion to this series.
The term with index $r$ before factoring is
\[
 T_r\coloneq\tau^{r(d+r)}\gbinom{p}{r}\gbinom{q}{d+r}
                         \gbinom{m+i-r}{i-r}.
\]
For consecutive nonzero terms its ratio is
\[
 \frac{T_{r+1}}{T_r}
 =\tau^{d+2r+1}
  \frac{(1-\tau^{p-r})(1-\tau^{q-d-r})(1-\tau^{i-r})}
       {(1-\tau^{r+1})(1-\tau^{d+r+1})(1-\tau^{m+i-r})}.
\]
Replacing each factor $1-\tau^{a-r}$ by
$-\tau^{a-r}(1-\tau^{-a+r})$ gives the consecutive-term
ratio of the displayed ${}_3\phi_2$ series. Its first term is
one, and $T_0=\gbinom{q}{d}\gbinom{m+i}{i}$, so the conversion
holds term by term, including its termination.
The series notation in this calculation means the terminating sum
\[
 {}_3\phi_2\!\left[\begin{matrix}a,b,\tau^{-N}\\c,d\end{matrix};\tau,z\right]
 \coloneq\sum_{r=0}^{N}
 \frac{(a;\tau)_r(b;\tau)_r(\tau^{-N};\tau)_r}
 {(c;\tau)_r(d;\tau)_r(\tau;\tau)_r}z^r,
\]
when its displayed denominators do not vanish.
The terminating $\tau$-Pfaff--Saalsch\"utz identity used here is
\[
 {}_3\phi_2\!\left[
 \begin{matrix}a,b,\tau^{-N}\\c,ab\tau^{1-N}/c\end{matrix};\tau,\tau\right]
 =\frac{(c/a;\tau)_N(c/b;\tau)_N}
 {(c;\tau)_N(c/(ab);\tau)_N}
\]
\cite{GasperRahman}. With
$a=\tau^{-p}$, $b=\tau^{-(q-d)}$, $c=\tau^{d+1}$, and $N=i$,
the second denominator parameter is
$ab\tau^{1-N}/c=\tau^{-(m+i)}$, and the sum terminates at or before
$r=i$. Thus the identity applies to the displayed series and gives
\[
 \gbinom{q}{d}\gbinom{m+i}{i}
 \frac{(\tau^{p+d+1};\tau)_i(\tau^{q+1};\tau)_i}
 {(\tau^{d+1};\tau)_i(\tau^{m+1};\tau)_i}
 =\gbinom{i+q}{i+d}\gbinom{i+d+p}{i}.
\]
This is the right-hand side. The case $i>j$ follows by interchanging
$(i,p)$ with $(j,q)$, and $\tau=1$ by continuity.
\end{proof}

\begin{proposition}[Factorization and polynomials for geometrically spaced roots]
\label{prop:q-binomial-family}
Use the finite products and Gaussian coefficients defined in
\eqref{eq:gaussian-notation}. Let $p,q\in\mathbb N$, $m\coloneq p+q$,
$\alpha>0$, $\tau>0$, and
\begin{equation}\label{eq:q-symbol}
 P_\tau(z)\coloneq\prod_{\nu=0}^{m-1}(1+\alpha\tau^\nu z)
 =(-\alpha z;\tau)_m.
\end{equation}
Let $T_\tau$ be the $(p,q)$-banded Toeplitz matrix with symbol
$z^{-q}P_\tau(z)$.
\begin{enumerate}
\item Its entries are
\begin{equation}\label{eq:q-binomial-entries}
 (T_\tau)_{i,j}=\alpha^h\tau^{\binom h2}\gbinom{m}{h},
 \qquad h\coloneq q+i-j,
\end{equation}
where the right-hand side is zero for $h\notin\{0,\dots,m\}$.

For factor indices $r\in\{0,\dots,p-1\}$ and
$s\in\{0,\dots,q-1\}$, define unit bidiagonal matrices by
\begin{align}
 (L_r)_{n,n-1}
 &\coloneq\alpha\tau^{q+r}\frac{1-\tau^n}{1-\tau^{n+q}},
 &&n\in\mathbb N,\label{eq:q-La}\\
 (U_s)_{n,n+1}
 &\coloneq\alpha^{-1}\tau^{-s}
   \frac{1-\tau^{n+1}}{1-\tau^{n+p+1}},
 &&n\in\NN,\label{eq:q-Ub}
\end{align}
and let
\begin{equation}\label{eq:q-D}
 D_{n,n}\coloneq\alpha^q\tau^{\binom q2}
 \frac{(\tau^{n+p+1};\tau)_q}{(\tau^{n+1};\tau)_q}.
\end{equation}
Then a bounded PBF is
\begin{equation}\label{eq:q-PBF}
 T_\tau=L_0L_1\cdots L_{p-1}\,D\,U_{q-1}U_{q-2}\cdots U_0.
\end{equation}
The lower factors are ordered by increasing index from left to right;
the upper factors by increasing index from right to left, so $U_0$
is the rightmost factor. Factors within either full block commute.
Their partial products, which determine the initial-condition
matrices, depend on this specified order.

\item The initial-condition matrices are defined by the partial products in
\eqref{eq:raw-H-definition}--\eqref{eq:normalized-boundary-definition},
with factor indices now starting at zero. For
$i\in\{0,\dots,q-1\}$ and $j\in\{0,\dots,i\}$, the entries of $H$ and $H^{-1}$ are
\begin{align}
 H_{i,j}&=\alpha^{i-j}\tau^{\binom{i-j}{2}}\gbinom{p+i}{i-j},
 \label{eq:q-H}\\
 (H^{-1})_{i,j}&=(-\alpha)^{i-j}\gbinom{p+i}{i-j}.
 \label{eq:q-Hinv}
\end{align}
For $j\in\{0,\dots,p-1\}$ and $i\in\{0,\dots,j\}$, the entries of $K$ and $K^{-1}$ are
\begin{align}
 K_{i,j}&=\alpha^{-(j-i)}
 \tau^{-(j-i)(q+i)-\binom{j-i}{2}}\gbinom{q+j}{j-i},
 \label{eq:q-K}\\
 (K^{-1})_{i,j}&=(-\alpha^{-1})^{j-i}
 \tau^{-(j-i)(q+i)}\gbinom{q+j}{j-i}.
 \label{eq:q-Kinv}
\end{align}
All omitted triangular entries are zero.

\item Let $A_n,B_n$ be the left and right polynomial vectors arranged as
in \eqref{eq:polynomial-arrays}, with the initial condition
\eqref{eq:general-initial-polynomials} formed from
\eqref{eq:q-Hinv} and \eqref{eq:q-Kinv}. For $b\in\{1,\dots,q\}$,
the right components satisfy
\begin{align}
 \sum_{n\in\NN}B_n^{(b)}(x)z^n
 &=\frac{z^{b-1}(-\alpha\tau^{p+b}z;\tau)_{q-b}}
 {P_\tau(z)-xz^q},\label{eq:q-B-GF}\\
 B_n^{(b)}(x)
 &=\sum_{k=0}^{\lfloor(n-b+1)/q\rfloor}
 x^k[z^{n-b+1-qk}]
 \frac1{(-\alpha z;\tau)_{p+b}(-\alpha z;\tau)_m^k}.
 \label{eq:q-B-coeff}
\end{align}
Here $[z^j]F(z)$ denotes the coefficient of $z^j$ in $F(z)$.
For the left components, use the reversed-symbol convention in
\eqref{eq:reversed-root-data}: put $t_p=\alpha^m\tau^{\binom m2}$ and
$P_\tau^\vee(z)\coloneq z^mP_\tau(1/z)=t_p(-\alpha^{-1}z;\tau^{-1})_m$.
For $a\in\{1,\dots,p\}$, the left components satisfy
\begin{equation}\label{eq:q-A-GF}
 \sum_{n\in\NN}A_n^{(a)}(x)z^n
 =\frac{t_pz^{a-1}
 (-\alpha^{-1}\tau^{-(q+a)}z;\tau^{-1})_{p-a}}
 {P_\tau^\vee(z)-xz^p}.
\end{equation}
\item At $\tau=1$, denote the present matrix and polynomial families by
$T_\alpha$, $B_{n,\alpha}^{(b)}$ and $A_{n,\alpha}^{(a)}$, and use
the subscript $\mathrm{bin}$ for the normalized binomial family in
\eqref{eq:T}, \eqref{eq:B} and \eqref{eq:A}. Then
\begin{equation}\label{eq:tau-one-matrix-scale}
 T_\alpha=\alpha^qD_\alpha T_{\rm bin}D_\alpha^{-1},
 \qquad D_\alpha\coloneq\diag[\alpha^n]_{n\in\NN}.
\end{equation}
This identity holds entrywise in the semi-infinite setting.
On every finite leading section it is a diagonal similarity
together with multiplication by $\alpha^q$.
The corresponding polynomial formulas are
\begin{equation}\label{eq:tau-one-polynomial-scale}
 \begin{aligned}
 B_{n,\alpha}^{(b)}(x)
 &=\alpha^{n-b+1}B_{n,\rm bin}^{(b)}(x/\alpha^q),
 &&b\in\{1,\dots,q\},\\
 A_{n,\alpha}^{(a)}(x)
 &=\alpha^{a-1-n}A_{n,\rm bin}^{(a)}(x/\alpha^q),
 &&a\in\{1,\dots,p\},
 \end{aligned}
 \qquad n\in\NN.
\end{equation}
\end{enumerate}
\end{proposition}
\begin{proof}
The finite $\tau$-binomial theorem gives
$(-\alpha z;\tau)_m=\sum_{h=0}^m \alpha^h\tau^{\binom h2}\gbinom{m}{h}z^h$,
which proves \eqref{eq:q-binomial-entries}. Define weighted shifts
$\mathcal N,\mathcal M$ by
$\mathcal N_{n,n-1}\coloneq\frac{1-\tau^n}{1-\tau^{n+q}}$ and
$\mathcal M_{n,n+1}\coloneq\frac{1-\tau^{n+1}}{1-\tau^{n+p+1}}$.
Then $L_r=I+\alpha\tau^{q+r}\mathcal N$ and
$U_s=I+\alpha^{-1}\tau^{-s}\mathcal M$, proving the commutation
claim. Set $\mathcal L\coloneq L_0\cdots L_{p-1}$ and
$\mathcal U\coloneq U_{q-1}\cdots U_0$. Another application of the finite
$\tau$-binomial theorem gives
\[
 \mathcal L=\sum_{r=0}^{p}\alpha^r\tau^{qr+\binom r2}
   \gbinom{p}{r}\mathcal N^r,\qquad
 \mathcal U=\sum_{s=0}^{q}\alpha^{-s}\tau^{-qs+\binom{s+1}{2}}
   \gbinom{q}{s}\mathcal M^s.
\]
The powers of the weighted shifts have entries
\[
 (\mathcal N^{i-k})_{i,k}
 =\frac{(\tau^{k+1};\tau)_{i-k}}
        {(\tau^{k+q+1};\tau)_{i-k}},\qquad
 (\mathcal M^{j-k})_{k,j}
 =\frac{(\tau^{k+1};\tau)_{j-k}}
        {(\tau^{k+p+1};\tau)_{j-k}},
\]
for $k\in\{0,\dots,\min(i,j)\}$. Each identity is the product of the entries
along the unique path of that shift. For $r=i-k$ and $s=j-k$,
substitution of these ratios and of $D_{k,k}$ cancels the factors
with endpoints $k+p$ and $k+q$. The corresponding product
$\mathcal L_{i,k}D_{k,k}\mathcal U_{k,j}$ is
\[
 \alpha^h\tau^{\binom h2}
 \frac{(\tau;\tau)_m(\tau;\tau)_i(\tau;\tau)_j}
 {(\tau;\tau)_{i+q}(\tau;\tau)_{j+p}}
 \tau^{(i-k)(j-k)}
 \gbinom{p}{i-k}\gbinom{q}{j-k}\gbinom{m+k}{k}.
\]
Summing in $k$ and using Lemma~\ref{lem:q-vandermonde-pfaff}
reduces the last factors to $\gbinom{m}{h}$, proving
\eqref{eq:q-PBF}. For positivity, the numerator and denominator of
each ratio in \eqref{eq:q-La}--\eqref{eq:q-D} have the same sign
when $\tau\ne1$. At $\tau=1$, a ratio
$(1-\tau^n)/(1-\tau^{n+r})$ has the positive limit $n/(n+r)$.
For fixed $\tau>0$ and $r\in\mathbb N$, this ratio tends, as
$n\to\infty$, to $1$ for $\tau\le1$ and to $\tau^{-r}$ for $\tau>1$.
The diagonal entries are finite products of such ratios or their
inverses, times a positive constant. Each sequence of allowed entries
therefore converges to a finite positive limit. Its finitely many
initial entries are also positive, so the entire sequence is bounded
above and away from zero. This proves the claimed bounded PBF.

For the initial-condition matrices, the products in
\eqref{eq:raw-H-definition}--\eqref{eq:normalized-boundary-definition},
with the present zero-based factor indices, are
$\widetilde K_{i,j}=e_i^{\mathsf T}L_0\cdots L_{j-1}e_0$ and
$\widetilde H_{i,j}=e_0^{\mathsf T}U_{i-1}\cdots U_0e_j$.
For $\widetilde K_{i,j}$, the finite binomial expansion of the first
$j$ lower factors selects the coefficient of $\mathcal N^i$ and the
entry $(\mathcal N^i)_{i,0}$. For $\widetilde H_{i,j}$, the corresponding
expansion of the last $i$ upper factors selects $\mathcal M^j$ and
$(\mathcal M^j)_{0,j}$. The power formulas above give
\[
 \widetilde K_{i,j}=\alpha^i\tau^{qi+\binom i2}
 \gbinom{j}{i}\frac{(\tau;\tau)_i}{(\tau^{q+1};\tau)_i},\qquad
 \widetilde H_{i,j}=\alpha^{-j}\tau^{-ij+\binom{j+1}{2}}
 \gbinom{i}{j}\frac{(\tau;\tau)_j}{(\tau^{p+1};\tau)_j}.
\]
Divide the first expression by its value at $i=j$ to normalize each
column, and the second by its value at $j=i$ to normalize each row.
The powers of $\alpha$ become $\alpha^{-(j-i)}$ and $\alpha^{i-j}$.
For the powers of $\tau$, the simplifications are
$\binom i2-\binom j2-q(j-i)=-(j-i)(q+i)-\binom{j-i}{2}$ and $\binom i2-ij+\binom{j+1}{2}=\binom{i-j}{2}$.
The Pochhammer ratios combine with the Gaussian coefficients to give
\eqref{eq:q-H}--\eqref{eq:q-K}.

To verify the inverses, multiply the proposed triangular matrices.
For $d\in\{0,\dots,N\}$ and $r\in\{0,\dots,d\}$, use
$\gbinom{N}{r}\gbinom{N-r}{d-r}=\gbinom{N}{d}\gbinom{d}{r}$,
which follows by cancelling their defining Pochhammer products.
For every allowed $d\in\NN$, the resulting entries are
\[
 \begin{aligned}
 (HH^{-1})_{k+d,k}
 &=(-\alpha)^d\gbinom{p+k+d}{d}(1;\tau)_d,\\
 (KK^{-1})_{i,i+d}
 &=\alpha^{-d}\tau^{-d(q+i)-\binom d2}
   \gbinom{q+i+d}{d}(1;\tau)_d.
 \end{aligned}
\]
Here the finite sums have been evaluated using
$\sum_{r=0}^{d}(-1)^r\tau^{\binom r2}\gbinom{d}{r} =(1;\tau)_d$.
For $d=0$ every factor is one; for $d>0$ the product $(1;\tau)_d$
contains $1-\tau^0=0$. Thus both products are identity matrices,
which verifies the asserted inverses.

Finally,
\[
 \frac{z^{b-1}(-\alpha\tau^{p+b}z;\tau)_{q-b}}
 {P_\tau(z)}=\frac{z^{b-1}}{(-\alpha z;\tau)_{p+b}}.
\]
The reciprocal $\tau$-binomial theorem gives the coefficient of
$z^n$ in this expression as
$(-\alpha)^{n-b+1}\gbinom{p+n}{n-b+1}$ when $n\ge b-1$, and zero
otherwise. For $n<q$ this is the entry $(n,b-1)$ of $H^{-1}$.
Adding $-xz^q$ to the denominator preserves those initial
coefficients, and multiplication by that denominator gives the right
recurrence by comparison of the coefficients of $z^{n+q}$.
This proves \eqref{eq:q-B-GF}. Expanding geometrically in
$xz^q/P_\tau$, with $P_\tau$ defined in \eqref{eq:q-symbol}, leaves, at power $x^k$, the reciprocal product
$[(-\alpha z;\tau)_{p+b}(-\alpha z;\tau)_m^k]^{-1}$ and the factor
$z^{b-1+qk}$. Their coefficient of $z^n$ is
\eqref{eq:q-B-coeff}.

For the left family, setting $x=0$ gives
$z^{a-1}/(-\alpha^{-1}z;\tau^{-1})_{q+a}$.
Its coefficient at $z^n$, with $r\coloneq n-a+1\ge0$, is
\[
 (-\alpha^{-1})^r\genfrac{[}{]}{0pt}{}{q+n}{r}_{\!\tau^{-1}}
 =(-\alpha^{-1})^r\tau^{-r(q+a-1)}\gbinom{q+n}{r},
\]
where the equality uses
$\genfrac{[}{]}{0pt}{}{N}{r}_{\!\tau^{-1}}
=\tau^{-r(N-r)}\gbinom{N}{r}$.
For $n<p$ these are exactly the entries $(a-1,n)$ of $K^{-1}$.
Since the numerator in \eqref{eq:q-A-GF} has degree $p-1$,
multiplication by $P_\tau^\vee-xz^p$ and comparison at $z^{n+p}$
gives the left recurrence. Its nonzero extreme coefficient $t_p$
ensures uniqueness.

Finally, at $\tau=1$ the entries are
$\alpha^{q+i-j}\binom{m}{q+i-j}$, which proves
\eqref{eq:tau-one-matrix-scale}. Substituting $w=\alpha z$ in
\eqref{eq:q-B-GF}, and $w=z/\alpha$ in \eqref{eq:q-A-GF},
reduces the two generating functions to those of the binomial family
with spectral variable $x/\alpha^q$. Comparing coefficients gives
\eqref{eq:tau-one-polynomial-scale}.
\end{proof}

To make the relation with the general construction explicit, choose
any ordered partition of the root parameters of $P_\tau$ and denote
its initial matrices, polynomial vectors, and measure by
$H_{\mathrm{ord}},K_{\mathrm{ord}}$,
$B_{n,\mathrm{ord}},A_{n,\mathrm{ord}}$, and $\Psi_{\mathrm{ord}}$.
Keep $H,K$ from \eqref{eq:q-H}--\eqref{eq:q-K} for the present
geometric family, and denote its measure by $\Psi$. Then
Corollary~\ref{cor:change-normalization} gives, for every $n\in\NN$,
\begin{equation}\label{eq:q-change-polynomials}
 B_n(x)=B_{n,\mathrm{ord}}(x)H_{\mathrm{ord}}H^{-1},\qquad
 A_n(x)=K^{-1}K_{\mathrm{ord}}A_{n,\mathrm{ord}}(x),
\end{equation}
and
\begin{equation}\label{eq:q-change-measure}
 \dd\Psi(x)=HH_{\mathrm{ord}}^{-1}\,
 \dd\Psi_{\mathrm{ord}}(x)\,K_{\mathrm{ord}}^{-1}K.
\end{equation}
Indeed, the first two changes give initial matrices $H^{-1}$ and
$K^{-1}$ and preserve the recurrences for the same $T_\tau$, so
uniqueness determines every polynomial vector. The measure identity
follows from their moments and compact support, as in the proof of
that corollary.

The value $\tau=1$ gives the binomial family by
\eqref{eq:tau-one-matrix-scale}--\eqref{eq:tau-one-polynomial-scale}.
The next proposition determines whether other values of $\tau$ can
be obtained from that family by translating and scaling the spectral
variable and rescaling the polynomial vectors and their components. Let
$T_{\rm bin}$ be the matrix in \eqref{eq:T}, with the same
bandwidths, and put $D_s\coloneq\diag[s^n]_{n\in\NN}$ for $s>0$.
The matrix identities below are understood entrywise, or on each
finite leading section.

\begin{proposition}[When the geometric family reduces to the binomial family]
\label{prop:q-affine-obstruction}
For $p+q\ge3$, an identity
\begin{equation}\label{eq:q-affine-equivalence}
 T_\tau=aD_sT_{\rm bin}D_s^{-1}+bI,
 \qquad a,s>0,\quad b\in\RR,
\end{equation}
holds if and only if $\tau=1$. In that case,
$a=\alpha^q$, $s=\alpha$, and $b=0$.
For $p=q=1$, the identity holds for every $\tau>0$, with
$a=s=\alpha\sqrt{\tau}$ and $b=\alpha(\sqrt{\tau}-1)^2$.
Allowing any positive diagonal matrix in place of $D_s$ gives
the same possibilities.
\end{proposition}
\begin{proof}
Put $m\coloneq p+q$ and write
$P_\tau(z)=\sum_{r=0}^m c_rz^r=\prod_{\nu=0}^{m-1}(1+\beta_\nu z)$,
where $\beta_\nu\coloneq\alpha\tau^\nu$. The coefficient $c_r$ lies on
diagonal $r-q$ of $T_\tau$. The term $bI$ changes only the main
diagonal, so \eqref{eq:q-affine-equivalence} requires
\begin{equation}\label{eq:q-affine-coefficients}
 c_r=a s^{r-q}\binom mr,
 \qquad r\in\{0,\dots,m\}\setminus\{q\}.
\end{equation}
Both extreme coefficients occur in this list because $p,q\ge1$.
The equality $c_0=1$ gives $a=s^q$, and
$c_m=\prod_{\nu=0}^{m-1}\beta_\nu$ then gives
$s=(\prod_{\nu=0}^{m-1}\beta_\nu)^{1/m}
=\alpha\tau^{(m-1)/2}$.

If $q\ge2$, the coefficient $c_1$ also lies outside the main
diagonal. Equation~\eqref{eq:q-affine-coefficients} requires
$\sum_{\nu=0}^{m-1}\beta_\nu=ms$. Thus the arithmetic mean of
the positive $\beta_\nu$ equals their geometric mean. Equality
in the arithmetic--geometric mean inequality forces all these
parameters to coincide, hence $\tau=1$.

If $q=1$ and $m\ge3$, use $c_{m-1}$ instead. It lies outside
the main diagonal, and
$c_{m-1}=s^m\sum_{\nu=0}^{m-1}\beta_\nu^{-1}$.
Equation~\eqref{eq:q-affine-coefficients} therefore requires
$\sum_{\nu=0}^{m-1}\beta_\nu^{-1}=m/s$.
The same inequality for the reciprocals again forces all the
parameters to coincide. At $\tau=1$, the main diagonal gives
$b=0$, and \eqref{eq:tau-one-matrix-scale} proves the converse.

For $p=q=1$, the two off-diagonals instead give
$a=s=\alpha\sqrt{\tau}$. The diagonal of $T_\tau$ is
$\alpha(1+\tau)$, whereas that of $aT_{\rm bin}$ is $2a$.
Their difference is the stated value of $b$.

Finally, suppose a positive diagonal matrix with entries $d_n$
replaces $D_s$. The first superdiagonal of its conjugation of
$T_{\rm bin}$ is a fixed positive coefficient times $d_n/d_{n+1}$.
Since $T_\tau$ also has a nonzero constant first superdiagonal,
these ratios must be constant. Thus $d_n=d_0s^n$ for some $s>0$,
and the scalar factor $d_0$ cancels in the conjugation.
\end{proof}

For $q=1$, $\alpha=1$, and $\tau=1$, the rescaling
\eqref{eq:lima-rescaling} therefore recovers Lima's
constant-coefficient polynomial system. For $p\ge2$ and
$\tau\ne1$, Proposition~\ref{prop:q-affine-obstruction} shows
that the geometric family cannot be obtained from that system by
an affine change of variable and positive diagonal normalization.
The scalar case $p=q=1$ remains an affine rescaling of the
Chebyshev polynomials described in
\eqref{eq:scalar-chebyshev-variable}--\eqref{eq:scalar-chebyshev-measure}.

\begin{corollary}[Explicit algebraic measure for the deformation]
\label{cor:q-measure}
For the PBF \eqref{eq:q-PBF}, use $H,K$ from \eqref{eq:q-H} and
\eqref{eq:q-K}, and let $E_r=[I_r\ 0]$ select the first $r$ coordinates.
The unique measure $\Psi_\tau$ supplied by
Corollary~\ref{cor:general-favard} satisfies
$\int_{\RR} x^\ell\dd\Psi_\tau(x)=HE_qT_\tau^\ell E_p^{\mathsf T}K$ and
$\int_{\RR} B_n\dd\Psi_\tau A_l=\delta_{n,l}$.
Here $A_n,B_n$ are the polynomial families in
Proposition~\ref{prop:q-binomial-family}.
Its endpoints are given by Theorem~\ref{thm:general-support}; its
exact support and rank-one density are given by Theorem~\ref{thm:general-measure} with
the root parameters $\alpha\tau^\nu$, $\nu\in\{0,\dots,m-1\}$,
arranged as $\beta_1\ge\cdots\ge\beta_m$ as in
Theorem~\ref{thm:general-support}. The measure is therefore
given by the algebraic root formula \eqref{eq:general-density-data}.
It is purely absolutely continuous on one interval.
For $\tau\ne1$ the scalar density $\rho$ in
\eqref{eq:general-density-data} vanishes as a square root at both
endpoints; at $\tau=1$ the lower endpoint is zero and has the
binomial behavior of Theorem~\ref{thm:hardedge}, with the scaling
\eqref{eq:tau-one-matrix-scale}.
\end{corollary}
\begin{proof}
Proposition~\ref{prop:q-binomial-family} shows that
\eqref{eq:q-PBF} is a PBF of $T_\tau$ with bounded factors and
computes its initial matrices \eqref{eq:q-H}--\eqref{eq:q-K}.
Corollary~\ref{cor:general-favard} therefore gives the unique
positive measure with the two identities, in precisely this
normalization. The matrix $T_\tau$ belongs to the Toeplitz class of
Theorem~\ref{thm:toeplitz}, with root parameters $\alpha\tau^\nu$ and
$\omega=1$. Theorems~\ref{thm:general-support}
and~\ref{thm:general-measure} hold for every bounded PBF of such a
matrix; applied to \eqref{eq:q-PBF}, they give the endpoints, the
exact support, absolute continuity, and the rank-one density
\eqref{eq:general-density-data} formed with the present $H,K$.
For $\tau\ne1$ the root parameters are distinct, so
$\beta_q>\beta_{q+1}$ and Proposition~\ref{prop:general-edges} gives
square-root vanishing of $\rho$ at both endpoints. At $\tau=1$ the
scaling \eqref{eq:tau-one-matrix-scale} reduces the family to the binomial
one, whose behavior at zero is given by Theorem~\ref{thm:hardedge}.
\end{proof}

Finite non-symmetric band matrices with Gaussian binomial
coefficients, including determinant, inverse, and LU formulas, were
studied by Ar{\i}kan and K{\i}l{\i}\c{c}
\cite{ArikanKilic2017}.

\section{Conclusions and outlook}

The results connect the choice of positive bidiagonal factors with
both their stability and the spectral information determined by the
matrix. The first question is whether a positive factorization exists. In finite dimension, BTP and PBF are equivalent. For semi-infinite matrices, the
difficulty is to choose initial entries that yield positive factors
in every row simultaneously. The global criterion identifies this
possibility precisely with convergence of the associated Metelmann series. Their terms are
positive, so convergence gives finite, strictly positive sums, from
which suitable initial entries can be chosen. This existence
criterion requires no boundedness assumption on the matrix or
its factors.

Existence alone does not ensure stability: arbitrarily small perturbations in distant rows can destroy a PBF. The stability criteria determine when every sufficiently small perturbation within the same band admits nearby positive factors that continue the prescribed PBF analytically. With fixed initial entries, relative entrywise bounds accommodate arbitrary positive growth, including different unbounded growth rates. With variable initial entries and uniformly positive bounded factors, the criterion combines conditions on every pair on each triangular side with joint row sums formed from the signed inverses of $L$ and $U^{\mathsf T}$. These sums control normalized LDU and are bounded at every operator-norm interior point of the PBF class. In this bounded setting, closeness is uniform over all rows; nearby initial entries alone need not give nearby factor sequences.

For BTP Toeplitz matrices, the relation between stability and spectral normalization becomes explicit through the roots of the symbol. The matrix lies in the interior of the PBF class precisely when these roots are distinct. In that case, the factorization attaining the minimum normalization bounds is locally isolated among the factorizations of the same matrix and has a unique nearby analytic continuation as the matrix varies. The roots also yield explicit formulas for the vector polynomials and their matrix of orthogonality measures. Although the polynomial vectors and the measure depend on the PBF, the changes in their normalizations cancel in the associated weighted linear forms. These forms therefore recover information intrinsic to the matrix and connect the Favard representation with finite-section eigenvectors, the spacing of zeros, local limits of the kernels, and the vector equilibrium problem. This description remains useful when the matrix is modified in a finite boundary block: such changes can introduce point masses into the measure while preserving the interior asymptotics on intervals where the boundary determinant does not vanish. Beyond finite boundary changes, convergence of the diagonals of a BTP matrix to those of a BTP Toeplitz matrix is sufficient for the characteristic polynomials to have the same limiting zero distribution as the latter.

Two directions remain. First, does every norm-interior matrix admit
some PBF satisfying the pairwise conditions? Its normalized LDU already
depends analytically; the remaining issue is the choice of bidiagonal
factors. The global geometry of all admissible initial entries, and its
relation to the uniform topology of factors, also remains to be determined
beyond two factors.

Second, does the boundary determinant remain nonzero throughout the
interior for every admissible finite block in arbitrary bandwidth?
For diagonals approaching those of a BTP Toeplitz matrix, the global zero
distribution is settled here. It remains to determine additional tail
hypotheses yielding strong polynomial asymptotics, local kernel limits,
and explicit Favard measures.

\section*{Funding}

The author was supported by the research project PID2024-155133NB-I00,
\emph{Ortogonalidad, aproximaci\'on e integrabilidad: aplicaciones en procesos
estoc\'asticos cl\'asicos y cu\'anticos}.

\section*{Competing interests}

The author declares no competing interests.

\section*{Code availability}

The \href{https://github.com/ManuelManas/btp-pbf-verification/releases/tag/final}{verification package} \cite{BTPPBFVerification} contains Lean formalizations of Lemmas~\ref{lem:rowwise}, \ref{lem:global-inverse-minors} and \ref{lem:global-extraction-determinant}, and Propositions~\ref{prop:successive-recurrences} and \ref{prop:inverse-parameter-bound}. It also formalizes the equivalence of conditions~(1)--(3) in Theorem~\ref{thm:global-lower-criterion}, including the extension of sufficiently small positive initial entries, and the inverse-minor limit criterion in Corollary~\ref{cor:global-two-sided-criterion}. These semi-infinite results take compatible Metelmann factors as input. For Theorem~\ref{thm:toeplitz}, the formalization proves the BTP--PBF equivalence and the existence of bounded positive factors; the full characterization involving total nonnegativity and negative roots assumes the classical root criterion explicitly. 

For stability, the general criteria in Theorems~\ref{stab:thm:general-lower-stability}, \ref{stab:thm:free-separated-limits}, \ref{stab:thm:general-full-stability} and~\ref{stab:thm:relative-criterion}, Propositions~\ref{stab:prop:ldu-row-sums} and~\ref{stab:prop:free-full-transfer}, and Lemma~\ref{stab:lem:relative-lower} are formalized on infinite sequence spaces, including positivity, uniqueness for fixed initial entries and local Lipschitz dependence. The bounded full-band criteria are also proved in the operator norm on $\ell^2$. Supporting scalar estimates, the abstract contraction argument and recorded computational checks are included. The supplementary verification table lists the result numbers and page numbers in the accompanying manuscript, together with the Lean declarations, hypotheses and scope. The Jacobi-fraction formulation and spectral results are not formalized.

\section*{Declaration of AI use}

During the preparation of this manuscript, the author used ChatGPT
and Codex (OpenAI) and Claude (Anthropic) as aids in exposition,
mathematical exploration, and the discussion and review of mathematical
arguments. These tools also assisted with computational checks,
the development of formal verification code, and the review of
notation and references.


\begin{thebibliography}{99}
\bibitem{ASW}
M. Aissen, I. J. Schoenberg, and A. M. Whitney,
\emph{On the generating functions of totally positive sequences I},
Journal d'Analyse Math\'ematique \textbf{2} (1952), 93--103.
\href{https://doi.org/10.1007/BF02786970}{doi:10.1007/BF02786970}.

\bibitem{Alexandersson2012}
P.~Alexandersson,
\emph{Schur polynomials, banded Toeplitz matrices and Widom's formula},
The Electronic Journal of Combinatorics \textbf{19} (2012), no.~4, Paper P22, 13~pp.
\href{https://doi.org/10.37236/2651}{doi:10.37236/2651}.

\bibitem{AptekarevKalyaginSaff2009}
A. I. Aptekarev, V. A. Kalyagin, and E. B. Saff,
\emph{Higher-order three-term recurrences and asymptotics of multiple
orthogonal polynomials}, Constructive Approximation \textbf{30}
(2009), no.~2, 175--223.
\href{https://doi.org/10.1007/s00365-008-9032-0}
{doi:10.1007/s00365-008-9032-0}.

\bibitem{ArikanKilic2017}
T. Ar{\i}kan and E. K{\i}l{\i}\c{c},
\emph{A class of non-symmetric band determinants with the Gaussian
$q$-binomial coefficients}, Quaestiones Mathematicae \textbf{40}
(2017), no.~5, 645--660.
\href{https://doi.org/10.2989/16073606.2017.1306596}
{doi:10.2989/16073606.2017.1306596}.

\bibitem{AriznabarretaManas2015}
G. Ariznabarreta and M. Ma\~nas,
\emph{A Jacobi type Christoffel--Darboux formula for multiple orthogonal
polynomials of mixed type}, Linear Algebra and its Applications
\textbf{468} (2015), 154--170.
\href{https://doi.org/10.1016/j.laa.2014.04.030}{doi:10.1016/j.laa.2014.04.030}.

\bibitem{BanderierFlajolet2002}
C. Banderier and P. Flajolet,
\emph{Basic analytic combinatorics of directed lattice paths},
Theoretical Computer Science \textbf{281} (2002), 37--80.
\href{https://doi.org/10.1016/S0304-3975(02)00007-5}{doi:10.1016/S0304-3975(02)00007-5}.

\bibitem{BarriosGeronimoLopez}
D. Barrios Rolan\'{\i}a, J. S. Geronimo, and G. L\'opez Lagomasino,
\emph{High-order recurrence relations, Hermite--Pad\'e approximation
and Nikishin systems}, Sbornik: Mathematics \textbf{209} (2018),
385--420.
\href{https://doi.org/10.1070/SM8724}{doi:10.1070/SM8724}.

\bibitem{Bauer1956}
F. L. Bauer,
\emph{Beitr\"age zur Entwicklung numerischer Verfahren f\"ur
programmgesteuerte Rechenanlagen. II. Direkte Faktorisierung eines Polynoms},
Sitzungsberichte der Bayerischen Akademie der Wissenschaften,
Mathematisch-Naturwissenschaftliche Klasse (1956), 163--203.
\url{https://publikationen.badw.de/de/003896006}.

\bibitem{BerensteinFominZelevinsky1996}
A. Berenstein, S. Fomin, and A. Zelevinsky,
\emph{Parametrizations of canonical bases and totally positive matrices},
Advances in Mathematics \textbf{122} (1996), no.~1, 49--149.
\href{https://doi.org/10.1006/aima.1996.0057}{doi:10.1006/aima.1996.0057}.

\bibitem{BleherKuijlaars2004}
P. M. Bleher and A. B. J. Kuijlaars,
\emph{Large $n$ limit of Gaussian random matrices with external source, part I},
Communications in Mathematical Physics \textbf{252} (2004), 43--76.
\href{https://doi.org/10.1007/s00220-004-1196-2}{doi:10.1007/s00220-004-1196-2}.

\bibitem{BogoyaGascaGrudsky2026}
M. Bogoya, J. Gasca, and S. Grudsky,
\emph{On the eigenvectors of non-Hermitian tetradiagonal Toeplitz matrices},
Numerical Linear Algebra with Applications \textbf{33} (2026), no.~2, e70074.
\href{https://doi.org/10.1002/nla.70074}{doi:10.1002/nla.70074}.

\bibitem{BogoyaGrudsky2025}
M. Bogoya and S. M. Grudsky,
\emph{Eigenvalues of non-Hermitian banded Toeplitz matrices approaching simple points of the limiting set},
Computational Mathematics and Mathematical Physics \textbf{65} (2025), 1453--1471.
\href{https://doi.org/10.1134/S0965542525700745}{doi:10.1134/S0965542525700745}.

\bibitem{BostanChyzakLecerfSalvySchost2007}
A. Bostan, F. Chyzak, B. Salvy, G. Lecerf, and \'{E}. Schost,
\emph{Differential equations for algebraic functions},
Proceedings of the 2007 International Symposium on Symbolic and
Algebraic Computation (ISSAC '07), ACM, New York, 2007, 25--32.
\href{https://doi.org/10.1145/1277548.1277553}{doi:10.1145/1277548.1277553}.

\bibitem{BiniIannazzoMeiniMengRobol2023}
D. A. Bini, B. Iannazzo, B. Meini, J. Meng, and L. Robol,
\emph{Computing eigenvalues of semi-infinite quasi-Toeplitz matrices},
Numerical Algorithms \textbf{92} (2023), 89--118.
\href{https://doi.org/10.1007/s11075-022-01381-0}
{doi:10.1007/s11075-022-01381-0}.

\bibitem{BottcherFukshanskyGarciaMaharaj2015}
A. B\"ottcher, L. Fukshansky, S. R. Garcia, and H. Maharaj,
\emph{Toeplitz determinants with perturbations in the corners},
Journal of Functional Analysis \textbf{268} (2015), 171--193.
\href{https://doi.org/10.1016/j.jfa.2014.10.023}{doi:10.1016/j.jfa.2014.10.023}.

\bibitem{BottcherGrudsky2005}
A. B\"ottcher and S. M. Grudsky,
\emph{Spectral Properties of Banded Toeplitz Matrices},
Society for Industrial and Applied Mathematics, Philadelphia, 2005.
\href{https://doi.org/10.1137/1.9780898717853}
{doi:10.1137/1.9780898717853}.

\bibitem{BottcherGrudskyRamirez2006}
A. B\"ottcher, S. Grudsky, and E. Ram\'irez de Arellano,
\emph{On the asymptotic behavior of the eigenvectors of large banded Toeplitz matrices},
Mathematische Nachrichten \textbf{279} (2006), 121--129.
\href{https://doi.org/10.1002/mana.200310349}{doi:10.1002/mana.200310349}.

\bibitem{BottcherSilbermann2006}
A. B\"ottcher and B. Silbermann,
\emph{Analysis of Toeplitz Operators}, second ed., Springer Monographs
in Mathematics, Springer, Berlin, 2006.
\href{https://doi.org/10.1007/3-540-32436-4}{doi:10.1007/3-540-32436-4}.

\bibitem{Uniform}
A. Branquinho, J. E. F. D\'{\i}az, A. Foulqui\'e-Moreno,
and M. Ma\~nas,
\emph{Uniform multiple orthogonal polynomials and Markov chains},
Numerical Algorithms \textbf{102} (2026), 1743--1772.
\href{https://doi.org/10.1007/s11075-025-02195-6}
{doi:10.1007/s11075-025-02195-6}.

\bibitem{BFMTetra}
A. Branquinho, A. Foulqui\'e-Moreno, and M. Ma\~nas,
\emph{Positive bidiagonal factorization of tetradiagonal Hessenberg
matrices}, Linear Algebra and its Applications \textbf{677} (2023),
132--160.
\href{https://doi.org/10.1016/j.laa.2023.08.001}
{doi:10.1016/j.laa.2023.08.001}.

\bibitem{BFM2023}
A. Branquinho, A. Foulqui\'e-Moreno, and M. Ma\~nas,
\emph{Spectral theory for bounded banded matrices with positive
bidiagonal factorization and mixed multiple orthogonal polynomials},
Advances in Mathematics \textbf{434} (2023), 109313.
\href{https://doi.org/10.1016/j.aim.2023.109313}{doi:10.1016/j.aim.2023.109313}.

\bibitem{BFM2024}
A. Branquinho, A. Foulqui\'e-Moreno, and M. Ma\~nas,
\emph{Banded totally positive matrices and normality for mixed-type multiple
orthogonal polynomials}, Studies in Applied Mathematics (2026), e70303.
\href{https://doi.org/10.1111/sapm.70303}{doi:10.1111/sapm.70303}.
Numbered results are cited according to
\href{https://arxiv.org/abs/2404.13965v1}{arXiv:2404.13965v1}.

\bibitem{BFMBidiagonal}
A. Branquinho, A. Foulqui\'e-Moreno, and M. Ma\~nas,
\emph{Bidiagonal factorization of banded recursion matrices for
mixed-type multiple orthogonal polynomials},
Linear Algebra and its Applications \textbf{752} (2027), 486--535.
\href{https://doi.org/10.1016/j.laa.2026.09.018}{doi:10.1016/j.laa.2026.09.018}.

\bibitem{BFMUnbounded}
A. Branquinho, A. Foulqui\'e-Moreno, and M. Ma\~nas,
\emph{Unbounded banded matrices, shifted positive bidiagonal
factorizations, and mixed-type multiple orthogonality},
Journal of Spectral Theory, to appear.
\href{https://arxiv.org/abs/2601.12453v3}{arXiv:2601.12453v3}.

\bibitem{BravermanKarabash}
E. Braverman and I. M. Karabash,
\emph{Bohl--Perron-type stability theorems for linear difference equations
with infinite delay}, Journal of Difference Equations and Applications
\textbf{18} (2012), no.~5, 909--939.
\href{https://doi.org/10.1080/10236198.2010.531276}{doi:10.1080/10236198.2010.531276}.

\bibitem{BuenoDopico2004}
M. I. Bueno and F. M. Dopico,
\emph{Stability and sensitivity of tridiagonal LU factorization
without pivoting}, BIT Numerical Mathematics \textbf{44} (2004), no.~4, 651--673.
\href{https://doi.org/10.1007/s10543-004-6025-7}
{doi:10.1007/s10543-004-6025-7}.

\bibitem{BreuerLastSimon2014}
J. Breuer, Y. Last, and B. Simon,
\emph{Stability of asymptotics of Christoffel--Darboux kernels},
Communications in Mathematical Physics \textbf{330} (2014), 1155--1178.
\href{https://doi.org/10.1007/s00220-014-1913-4}{doi:10.1007/s00220-014-1913-4}.

\bibitem{Chihara1962}
T. S. Chihara,
\emph{Chain sequences and orthogonal polynomials},
Transactions of the American Mathematical Society \textbf{104} (1962), 1--16.
\href{https://doi.org/10.1090/S0002-9947-1962-0138933-7}{doi:10.1090/S0002-9947-1962-0138933-7}.

\bibitem{Cryer1973}
C. W. Cryer,
\emph{The LU-factorization of totally positive matrices},
Linear Algebra and its Applications \textbf{7} (1973), no.~1, 83--92.
\href{https://doi.org/10.1016/0024-3795(73)90039-6}
{doi:10.1016/0024-3795(73)90039-6}.

\bibitem{DaemsKuijlaars2007}
E. Daems and A. B. J. Kuijlaars,
\emph{Multiple orthogonal polynomials of mixed type and non-intersecting
Brownian motions}, Journal of Approximation Theory
\textbf{146} (2007), 91--114.
\href{https://arxiv.org/abs/math/0511470}{arXiv:math/0511470}.

\bibitem{DaemsKuijlaarsVeys2008}
E. Daems, A. B. J. Kuijlaars, and W. Veys,
\emph{Asymptotics of non-intersecting Brownian motions and a $4\times4$ Riemann--Hilbert problem},
Journal of Approximation Theory \textbf{153} (2008), 225--256.
\href{https://doi.org/10.1016/j.jat.2007.10.003}{doi:10.1016/j.jat.2007.10.003}.

\bibitem{DavisHershMiller}
J. F. Davis, P. Hersh, and E. Miller,
\emph{Fibers of maps to totally nonnegative spaces},
preprint, version of December 1, 2024.
\href{https://arxiv.org/abs/1903.01420v2}{arXiv:1903.01420v2}.

\bibitem{deBoorPinkus}
C. de Boor and A. Pinkus,
\emph{The approximation of a totally positive band matrix by a
strictly banded totally positive one}, Linear Algebra and its
Applications \textbf{42} (1982), 81--98.
\href{https://doi.org/10.1016/0024-3795(82)90139-2}
{doi:10.1016/0024-3795(82)90139-2}.

\bibitem{Deimling}
K. Deimling, \emph{Nonlinear Functional Analysis},
Springer, 1985. \href{https://doi.org/10.1007/978-3-662-00547-7}{doi:10.1007/978-3-662-00547-7}.

\bibitem{DelvauxLopez}
S. Delvaux, A. L\'opez Garc\'{\i}a, and G. L\'opez Lagomasino,
\emph{A family of Nikishin systems with periodic recurrence
coefficients}, Sbornik: Mathematics \textbf{204} (2013), 43--74.
\href{https://doi.org/10.1070/SM2013v204n01ABEH004291}
{doi:10.1070/SM2013v204n01ABEH004291}.

\bibitem{DuitsKuijlaars2008}
M. Duits and A. B. J. Kuijlaars,
\emph{An equilibrium problem for the limiting eigenvalue distribution
of banded Toeplitz matrices},
SIAM Journal on Matrix Analysis and Applications \textbf{30} (2008),
no.~1, 173--196.
\href{https://doi.org/10.1137/070687141}{doi:10.1137/070687141}.

\bibitem{ErbeMathsenPeterson2001}
L. H. Erbe, R. M. Mathsen, and A. Peterson,
\emph{Factoring linear differential operators on measure chains},
Journal of Inequalities and Applications \textbf{6} (2001), no.~3,
287--303.
\href{https://doi.org/10.1155/S1025583401000170}{doi:10.1155/S1025583401000170}.

\bibitem{FallatJohnson2011}
S. M. Fallat and C. R. Johnson,
\emph{Totally Nonnegative Matrices},
Princeton Series in Applied Mathematics,
Princeton University Press, Princeton, NJ, 2011.

\bibitem{FlajoletSedgewick2009}
P. Flajolet and R. Sedgewick,
\emph{Analytic Combinatorics},
Cambridge University Press, Cambridge, 2009.
\href{https://doi.org/10.1017/CBO9780511801655}{doi:10.1017/CBO9780511801655}.

\bibitem{FominZelevinsky2000}
S. Fomin and A. Zelevinsky,
\emph{Total positivity: Tests and parametrizations},
The Mathematical Intelligencer \textbf{22} (2000), no.~1, 23--33.
\href{https://doi.org/10.1007/BF03024444}{doi:10.1007/BF03024444}.

\bibitem{GantmacherKrein}
F. R. Gantmacher and M. G. Krein,
\emph{Oscillation Matrices and Kernels and Small Vibrations of
Mechanical Systems}, revised ed., AMS Chelsea Publishing 345,
American Mathematical Society, Providence, RI, 2002.
\href{https://doi.org/10.1090/chel/345}{doi:10.1090/chel/345}.

\bibitem{GascaMicchelliPena1992}
M. Gasca, C. A. Micchelli, and J. M. Pe\~na,
\emph{Almost strictly totally positive matrices},
Numerical Algorithms \textbf{2} (1992), no.~2, 225--236.
\href{https://doi.org/10.1007/BF02145387}{doi:10.1007/BF02145387}.

\bibitem{GascaPena1995}
M. Gasca and J. M. Pe\~na,
\emph{On the characterization of almost strictly totally positive
matrices}, Advances in Computational Mathematics \textbf{3}
(1995), no.~3, 239--250.
\href{https://doi.org/10.1007/BF02432001}{doi:10.1007/BF02432001}.

\bibitem{GasperRahman}
G. Gasper and M. Rahman,
\emph{Basic Hypergeometric Series}, second ed., Encyclopedia of
Mathematics and its Applications 96, Cambridge University Press,
Cambridge, 2004.
\href{https://doi.org/10.1017/CBO9780511526251}
{doi:10.1017/CBO9780511526251}.

\bibitem{GrinbergEtAl2024}
D. Grinberg, N. Korniichuk, K. Molokanov, and S. Khomych,
\emph{The diagonal derivative of a skew Schur polynomial},
\href{https://arxiv.org/abs/2402.14217v1}{arXiv:2402.14217v1} (2024).

\bibitem{GeronimoIliev2017}
J. S. Geronimo and P. Iliev,
\emph{Bernstein--Szeg\H{o} measures, Banach algebras, and scattering theory},
Transactions of the American Mathematical Society \textbf{369} (2017),
no.~8, 5581--5600.
\href{https://doi.org/10.1090/tran/6841}{doi:10.1090/tran/6841}.

\bibitem{HarrisKruegerTrench1998}
B. Harris, R. Krueger, and W. F. Trench,
\emph{Trench's canonical form for a disconjugate $n$th-order linear
difference equation}, Panamerican Mathematical Journal
\textbf{8} (1998), 55--71.
\href{https://ramanujan.math.trinity.edu/wtrench/research/papers/TRENCH_RP_90.PDF}
{Author's version}.

\bibitem{Hartman1978}
P. Hartman,
\emph{Difference equations: disconjugacy, principal solutions,
Green's functions, complete monotonicity},
Transactions of the American Mathematical Society \textbf{246} (1978),
1--30.
\href{https://doi.org/10.1090/S0002-9947-1978-0515528-6}{doi:10.1090/S0002-9947-1978-0515528-6}.

\bibitem{Higham1994}
N. J. Higham,
\emph{A survey of componentwise perturbation theory in numerical
linear algebra}, in W.~Gautschi (ed.), \emph{Mathematics of Computation
1943--1993: A Half-Century of Computational Mathematics},
Proceedings of Symposia in Applied Mathematics \textbf{48},
American Mathematical Society, Providence, RI, 1994, pp.~49--77.
\href{https://eprints.maths.manchester.ac.uk/651/}{Author's version}.

\bibitem{Higham2024}
N. J. Higham,
\emph{The power of bidiagonal matrices},
Electronic Journal of Linear Algebra \textbf{40} (2024), 453--474.
\href{https://doi.org/10.13001/ela.2024.8297}
{doi:10.13001/ela.2024.8297}.

\bibitem{Hirschman1967}
I. I. Hirschman, Jr.,
\emph{The spectra of certain Toeplitz matrices},
Illinois Journal of Mathematics \textbf{11} (1967), no.~1, 145--159.
\href{https://doi.org/10.1215/ijm/1256054792}
{doi:10.1215/ijm/1256054792}.

\bibitem{HuyMinh2001}
N. T. Huy and N. V. Minh,
\emph{Exponential dichotomy of difference equations and applications
to evolution equations on the half-line},
Computers \& Mathematics with Applications \textbf{42} (2001),
nos.~3--5, 301--311.
\href{https://doi.org/10.1016/S0898-1221(01)00155-9}
{doi:10.1016/S0898-1221(01)00155-9}.

\bibitem{Lima2023}
H. Lima,
\emph{Multiple orthogonal polynomials associated with branched
continued fractions for ratios of hypergeometric series},
Advances in Applied Mathematics \textbf{147} (2023), 102505.
\href{https://doi.org/10.1016/j.aam.2023.102505}{doi:10.1016/j.aam.2023.102505}.
\href{https://arxiv.org/abs/2208.03539v2}{arXiv:2208.03539v2}.
Numbering refers to this arXiv version.

\bibitem{LimaSymmetric2026}
H. Lima,
\emph{Bidiagonal matrix factorisations associated with symmetric
multiple orthogonal polynomials and lattice paths},
Numerical Algorithms \textbf{102} (2026), no.~3, 1527--1552.
\href{https://doi.org/10.1007/s11075-025-02093-x}
{doi:10.1007/s11075-025-02093-x}.

\bibitem{Lusztig1994}
G. Lusztig,
\emph{Total positivity in reductive groups},
in \emph{Lie Theory and Geometry}, Progress in Mathematics,
vol.~123, Birkh\"auser, Boston, 1994, 531--568.
\href{https://doi.org/10.1007/978-1-4612-0261-5_20}
{doi:10.1007/978-1-4612-0261-5\_20}.

\bibitem{BTPPBFVerification}
M. Ma\~nas,
\emph{BTP/PBF verification package},
verification archive (2026).
\url{https://github.com/ManuelManas/btp-pbf-verification/releases/tag/final}.

\bibitem{ManasPBF2026}
M. Ma\~nas,
\emph{Positive bidiagonal factorizations: Combinatorial structure,
spectral theory, and matrix continued fractions}, preprint (2026).
\href{https://arxiv.org/abs/2609.26066v2}{arXiv:2609.26066v2}.

\bibitem{MaximenkoMoctezuma2017}
E. A. Maximenko and M. A. Moctezuma-Salazar,
\emph{Cofactors and eigenvectors of banded Toeplitz matrices:
Trench formulas via skew Schur polynomials},
Operators and Matrices \textbf{11} (2017), no.~4, 1149--1169.
\href{https://doi.org/10.7153/oam-2017-11-79}{doi:10.7153/oam-2017-11-79}.

\bibitem{Metelmann}
K. Metelmann,
\emph{Ein Kriterium f\"ur den Nachweis der Totalnichtnegativit\"at
von Bandmatrizen}, Linear Algebra and its Applications
\textbf{7} (1973), 163--171.
\href{https://doi.org/10.1016/0024-3795(73)90050-5}
{doi:10.1016/0024-3795(73)90050-5}.

\bibitem{Nahay2004Powersum}
J. M. Nahay,
\emph{Powersum formula for differential resolvents},
International Journal of Mathematics and Mathematical Sciences
\textbf{2004} (2004), no.~7, 365--371.
\href{https://doi.org/10.1155/S0161171204210602}{doi:10.1155/S0161171204210602}.

\bibitem{NISTDLMF}
National Institute of Standards and Technology,
\emph{Digital Library of Mathematical Functions},
Sections~\href{https://dlmf.nist.gov/16.2}{16.2},
\href{https://dlmf.nist.gov/16.8}{16.8}, and
\href{https://dlmf.nist.gov/18.12}{18.12} (accessed 17 October 2026).

\bibitem{PetreolleSokalZhu}
M. P\'etr\'eolle, A. D. Sokal, and B.-X. Zhu,
\emph{Lattice paths and branched continued fractions: An infinite
sequence of generalizations of the Stieltjes--Rogers and
Thron--Rogers polynomials, with coefficientwise Hankel-total positivity},
Memoirs of the American Mathematical Society \textbf{291} (2023),
no.~1450.
\href{https://doi.org/10.1090/memo/1450}{doi:10.1090/memo/1450}.

\bibitem{SchmidtSpitzer}
P. Schmidt and F. Spitzer,
\emph{The Toeplitz matrices of an arbitrary Laurent polynomial},
Mathematica Scandinavica \textbf{8} (1960), 15--38.
\href{https://doi.org/10.7146/math.scand.a-10588}
{doi:10.7146/math.scand.a-10588}.

\bibitem{Sokal2024}
A. D. Sokal,
\emph{Multiple orthogonal polynomials, $d$-orthogonal polynomials,
production matrices, and branched continued fractions},
Transactions of the American Mathematical Society, Series B
\textbf{11} (2024), 762--797.
\href{https://doi.org/10.1090/btran/133}{doi:10.1090/btran/133}.

\bibitem{SwiderskiVanAssche2022}
G. \'Swiderski and W. Van Assche,
\emph{Christoffel functions for multiple orthogonal polynomials},
Journal of Approximation Theory \textbf{283} (2022), 105820.
\href{https://doi.org/10.1016/j.jat.2022.105820}{doi:10.1016/j.jat.2022.105820}.

\bibitem{Trench1974}
W. F. Trench,
\emph{Canonical forms and principal systems for general disconjugate
equations},
Transactions of the American Mathematical Society \textbf{189} (1974),
319--327.
\href{https://doi.org/10.1090/S0002-9947-1974-0330632-X}{doi:10.1090/S0002-9947-1974-0330632-X}.

\bibitem{Ullman}
J. L. Ullman,
\emph{A problem of Schmidt and Spitzer},
Bulletin of the American Mathematical Society \textbf{73} (1967),
883--885.
\href{https://doi.org/10.1090/S0002-9904-1967-11826-3}
{doi:10.1090/S0002-9904-1967-11826-3}.

\bibitem{Wall1948}
H. S. Wall,
\emph{Analytic Theory of Continued Fractions},
D. Van Nostrand, New York, 1948.

\bibitem{WebbOlver2021}
M. Webb and S. Olver,
\emph{Spectra of Jacobi operators via connection coefficient matrices},
Communications in Mathematical Physics \textbf{382} (2021), 657--707.
\href{https://doi.org/10.1007/s00220-021-03939-w}
{doi:10.1007/s00220-021-03939-w}.

\bibitem{YoulaKazanjian1978}
D. C. Youla and N. N. Kazanjian,
\emph{Bauer-type factorization of positive matrices and the theory of
matrix polynomials orthogonal on the unit circle},
IEEE Transactions on Circuits and Systems \textbf{25} (1978),
no.~2, 57--69.
\href{https://doi.org/10.1109/TCS.1978.1084443}{doi:10.1109/TCS.1978.1084443}.

\end{thebibliography}
\end{document}